\documentclass[10pt, a4paper]{article}
\usepackage{graphicx}
\usepackage{color}
\usepackage{amsmath}
\usepackage{amssymb}
\usepackage{amscd}
\usepackage{bbm}
\usepackage{tikz}
\usepackage{framed}
\usepackage[normalem]{ulem}
\usetikzlibrary{positioning,arrows.meta,bending,calc}
\usepackage{hyperref}
\hypersetup{
    bookmarks=true,         
    colorlinks=true,       
    linkcolor=red,          
    citecolor=green,        
    filecolor=magenta,      
    urlcolor=cyan           
}
\usepackage[most]{tcolorbox}
\usepackage[utf8]{inputenc}
\usepackage{amsthm}
\usepackage{bbm} 
\usepackage{glossaries}
\usepackage[linesnumbered,lined,boxed,commentsnumbered]{algorithm2e}

\usepackage{marginnote}

\newtheorem{Theorem}{Theorem}
\newtheorem{Assumption}{Assumption}

\newtheorem{Lemma}{Lemma}
\newtheorem{Proposition}{Proposition}
\newtheorem{Corollary}{Corollary}
\newtheorem{Remark}{Remark}
\newtheorem{Example}{Example}

\newenvironment{Propositionrewrite}[1]
    {\renewcommand{\theProposition}{\ref{#1} (rewrite)}%
    \addtocounter{Proposition}{-1}%
    \begin{Proposition}}
    {\end{Proposition}}
\newenvironment{Examplecontinue}[1]
    {\renewcommand{\theExample}{\ref{#1} (continue)}%
    \addtocounter{Example}{-1}%
    \begin{Example}}
    {\end{Example}}
    
\newtheorem*{Propositionextension}{Proposition \extensionnumber\space (extension)}
\newtheorem*{Lemmaextension}{Lemma \extensionnumber\space (extension)}
\makeatletter
\newenvironment{Propositionextend}[1]
    {\def\extensionnumber{\ref{#1}}%
    \begin{Propositionextension}\phantomsection
    \protected@edef\@currentlabel{\getrefnumber{#1} (extension)}}
    {\end{Propositionextension}}
\newenvironment{Lemmaextend}[1]
    {\def\extensionnumber{\ref{#1}}%
    \begin{Lemmaextension}\phantomsection
    \protected@edef\@currentlabel{\getrefnumber{#1} (extension)}}
    {\end{Lemmaextension}}
\makeatother

\newcounter{BigConst}                     

\newcounter{SmallConst}                     

\newcounter{gamma}                     

\newcounter{kappa}                     

\newcounter{delta}                     

\newcounter{theta}                     

\newcounter{L}                     

\newcounter{eps}                     

\newcommand{\inr}[1]{\bigl< #1 \bigr>}

\newcommand{\norm}[1]{\left\|#1\right\|}%
\newcommand{\vertiii}[1]{{\left\vert\kern-0.25ex\left\vert\kern-0.25ex\left\vert #1 
    \right\vert\kern-0.25ex\right\vert\kern-0.25ex\right\vert}}

\newcommand{\vertiiii}[1]{{\left\vert\kern-0.25ex\left\vert\kern-0.25ex\left\vert\kern-0.25ex\left\vert #1 
    \right\vert\kern-0.25ex\right\vert\kern-0.25ex\right\vert\kern-0.25ex\right\vert}}

\newcommand{\beginproof}{{\bf Proof. {\hspace{0.2cm}}}}
\def \endproof
{{\mbox{}\nolinebreak\hfill\rule{2mm}{2mm}\par\medbreak}}

\DeclareMathOperator*{\argmax}{argmax}
\DeclareMathOperator*{\argmin}{argmin}
\DeclareMathOperator*{\supp}{supp}

\DeclareMathOperator{\diag}{diag}
\DeclareMathOperator{\Span}{span}

\DeclareMathOperator*{\dist}{dist}

\DeclareMathOperator*{\diam}{diam}

\DeclareMathOperator*{\sign}{sign}

\newcommand{\Proj}{{\mathrm{Proj}}}
\def\ds1{\textrm{1\kern-0.25emI}} 

\newcommand{\D}{\ensuremath{\mathbf{D}}}

\newcommand{\1}{\ensuremath{\mathbbm{1}}}
\DeclareMathOperator{\Var}{\mathrm{Var}}

\newcommand \N{\mathbb{N}}
\newcommand \cA{{\cal A}}
\newcommand \cB{{\cal B}}
\newcommand \cC{{\cal C}}

\newcommand \cF{{\cal F}}
\newcommand \cG{{\cal G}}

\newcommand \cL{{\cal L}}
\newcommand \cM{{\cal M}}
\newcommand \cN{{\cal N}}

\newcommand \cQ{{\cal Q}}

\newcommand \cS{{\cal S}}
\newcommand \cT{{\cal T}}

\newcommand \cX{{\cal X}}
\newcommand \cY{{\cal Y}}

\newcommand \bE{{\mathbb E}}

\newcommand \bG{{\mathbb G}}

\newcommand \bN{{\mathbb N}}

\newcommand \bP{{\mathbb P}}

\newcommand \bR{{\mathbb R}}

\newcommand \bX{{\mathbb X}}

\newcommand{\bbeta}{{\boldsymbol{\beta}}}

\newcommand{\bu}{\boldsymbol{u}}

\newcommand{\va}{{\boldsymbol{a}}}
\newcommand{\vb}{{\boldsymbol{b}}}

\newcommand{\ve}{{\boldsymbol{e}}}

\newcommand{\vg}{{\boldsymbol{g}}}
\newcommand{\vh}{{\boldsymbol{h}}}

\newcommand{\vs}{{\boldsymbol{s}}}
\newcommand{\vt}{{\boldsymbol{t}}}
\newcommand{\vu}{{\boldsymbol{u}}}
\newcommand{\vv}{{\boldsymbol{v}}}
\newcommand{\vw}{{\boldsymbol{w}}}
\newcommand{\vx}{{\boldsymbol{x}}}
\newcommand{\vy}{{\boldsymbol{y}}}
\newcommand{\vz}{{\boldsymbol{z}}}

\newcommand{\vA}{{\boldsymbol{A}}}

\newcommand{\vG}{{\boldsymbol{G}}}

\newcommand{\vbeta}{{\boldsymbol{\beta}}}

\newcommand{\vdelta}{{\boldsymbol{\delta}}}

\newcommand{\vveps}{{\boldsymbol{\varepsilon}}}
\newcommand{\vzeta}{{\boldsymbol{\zeta}}}

\newcommand{\vlambda}{{\boldsymbol{\lambda}}}
\newcommand{\vmu}{{\boldsymbol{\mu}}}
\newcommand{\vnu}{{\boldsymbol{\nu}}}

\newcommand{\vomega}{{\boldsymbol{\omega}}}
\newcommand{\vzero}{\boldsymbol{0}}

\newcommand{\bxi}{{\boldsymbol{\xi}}}

\DeclareMathOperator*{\Tr}{Tr}

\newcommand{\mysymbol}[3]{%
\newglossaryentry{#1}{%
      name={\ensuremath{#2}},%
      text={\ensuremath{#2}},%
      description={#3},%
      sort={#1}%
    }%
\expandafter\newcommand\expandafter{\csname smb#1\endcsname}{\gls{#1}}%
\expandafter\newcommand\expandafter{\csname #1\endcsname}{\ensuremath{#2}}%
}

\newcommand{\op}{\text{op}}

\usepackage[english]{babel}
\usepackage{authblk}
\usepackage{enumitem}
\begin{document}

\title{Sharp convergence rate of the estimation error of LASSO via the Feature Space Decomposition method}
\author{Guillaume Lecu\'e\thanks{ESSEC Business School, 3 avenue Bernard Hirsch, 95021 Cergy-Pontoise, France. E-mail: \texttt{lecue@essec.edu}}\,\,\, and\,\,\, Zong Shang\thanks{School of Mathematics, Georgia Institute of Technology, Atlanta, GA, 30332-0160, USA. E-mail: \texttt{zong.shang@ensae.fr}}\,\,\,%
                \thanks{Department of Statistics, CREST-ENSAE, Institut Polytechnique de Paris, 5 Avenue Henry Le Chatlier, Palaiseau, 91120, France.}
        }
\date{\today}
\maketitle
\begin{abstract}
    In this paper, we use the Feature Space Decomposition (FSD) method to establish an explicit, non-asymptotic, instance-optimal convergence rate for the estimation error of the LASSO in the Gaussian linear regression problem for any signal, any diagonal covariance and any regularization parameter above $\sigma_\xi/\sqrt{N}$ where $\sigma_\xi$ is the noise variance and $N$ the number of data. 
    
\end{abstract}

This paper is the sixth one in the series on the Feature Space Decomposition following \cite{lecueGeometricalViewpointBenign2024}, \cite{gavrilopoulosGeometricalAnalysisKernel2025}, \cite{lecueSharpConvergenceRates2025}, \cite{shangGeometryStatisticalFeature2026},  and \cite{adamczakFeatureSpaceDecompositioninpreperation}. The position of this paper within the FSD series is as follows: by studying LASSO, we demonstrate how the FSD method can be used to study general nonlinear (Regularized) Empirical Risk Minimization, and offer new perspectives and tools in statistical learning theory.

\section{Introduction}

We consider a Gaussian linear regression problem. Let $p > N$ be two positive integers, let $\Sigma=\operatorname{diag}(\sigma_1,\ldots,\sigma_p)$ with $\sigma_j>0$ for any $j\in[p]$, let $\mathbb{P}_X=\mathcal N(\boldsymbol0,\Sigma)$ be the centered Gaussian distribution with covariance $\Sigma$ on $\mathbb{R}^p$, and let $X$ be a random vector with distribution $\mathbb{P}_X$, called the design vector; let $\bxi = (\xi_i)_{i=1}^N$ be a centered Gaussian random vector on $\mathbb{R}^N$ with covariance matrix $\sigma_\xi^2 I_N$, where $\sigma_\xi^2>0$ is called the variance of noise. Assume that $X$ and $\bxi$ are independent. Let $X_1,\dots,X_N$ be $N$ independent copies of $X$. Denote by $\bX: \vv \in \mathbb{R}^p \mapsto (\langle X_i, \vv \rangle)_{i=1}^N \in \mathbb{R}^N$ the design matrix. Assume that $\bbeta^\star \in \mathbb{R}^p$ is a nonzero vector, called the signal, and define $\vy = \bX \bbeta^\star + \bxi$, called the response vector. Let $\ve_1,\dots,\ve_p$ be an canonical basis of $\mathbb{R}^p$, and the $\ell_1^p$ norm $\|\cdot\|_1$. Let $\lambda \ge 0$ be a real number, and define the LASSO with tuning parameter $\lambda\geq 0$ as the following nonlinear operator from $\mathbb{R}^N$ to $\mathbb{R}^p$: when $\lambda>0$,
\begin{align*}
    &\hat\bbeta=\hat\bbeta[\vy],\mbox{ where } \hat\bbeta[\cdot]:\vy\in\bR^N\mapsto \hat\bbeta[\vy]\in \operatorname*{argmin}_{\boldsymbol\beta\in\mathbb R^p}\left(\frac{1}{2N\lambda}\|\vy-\mathbb X\boldsymbol\beta\|_2^2+\|\boldsymbol\beta\|_1\right),
\end{align*}and when $\lambda=0$, $\hat\bbeta[\vy]\in\argmin(\|\bbeta\|_1:\, \bX\bbeta=\vy)$ is known as the basis pursuit \cite{donoho_compressed_2006,foucartMathematicalIntroductionCompressive2013}. One of the most important statistical properties of $\hat\bbeta$ is its excess generalization risk with respect to the squared-loss that coincides with the  estimation error with respect to $\bbeta^\star$, namely the random variable $\|\Sigma^{1/2}(\hat\bbeta - \bbeta^\star)\|_2^2$. In this paper, we use the Feature Space Decomposition (FSD) method developed in a series of works by \cite{lecueGeometricalViewpointBenign2024,gavrilopoulosGeometricalAnalysisKernel2025,lecueSharpConvergenceRates2025,shangGeometryStatisticalFeature2026,adamczakFeatureSpaceDecompositioninpreperation} to establish non-asymptotic, matching upper and lower bounds for the estimation error.

\subsection{Main results}

To present the main results of this paper, we first define some quantities. A discussion of these quantities is deferred to Section~\ref{sec:discussion}. Let \(\Phi_{\mathrm{G}}(t)=\mathbb P\{G\geq t\}\) for \(G\sim \cN(0,1)\) being a standard Gaussian random variable, and let \(\Phi_{\mathrm{G}}^{-1}\) denote its inverse.  We use $\sign(x)$ to denote the sign of a real number $x$; and use the convention \(\operatorname{sign}(0)=0\). Let $\beta_j^\star=\langle\bbeta^\star,\ve_j\rangle$ be the $j$-th coordinate of $\bbeta^\star$. Let $[p]=\{1,2,\cdots,p\}$ and let
\begin{equation}\label{eq:set_of_selected_features_LASSO}
    J_\star:=J_\star(\lambda)=\{j\in[p]:\sigma_j|\beta_j^\star|>\lambda\},\, 
\end{equation} $V_{J_\star}:=\operatorname{span}(\boldsymbol e_j:j\in J_\star)$,  $V_{J_\star^c}:=\operatorname{span}(\boldsymbol e_j:j\in J_\star^c)$ and $J_\star^c : =J_\star^c(\lambda)=[p]\backslash J_\star(\lambda).$

We adopt the convention that for any $\bbeta \in \mathbb{R}^p$, $\bbeta_{J_\star} = \sum_{j \in J_\star} \langle \bbeta, \ve_j \rangle \ve_j$, and $\bbeta_{J_\star^c} = \sum_{j \in J_\star^c} \langle \bbeta, \ve_j \rangle \ve_j$. For any $J\subset[p]$, denote $\Sigma_{J} = \sum_{j\in J}\sigma_j\ve_j\otimes\ve_j$. For a finite set $A$, we denote its cardinality by $|A|$. In addition, we define the following two quantities: for any $0<r<|J_\star^c|$, define $z_r$ be the solution of
\begin{align*}
\sum_{j\in J_\star^c}2\Phi_G\bigg(\frac{z_r}{\sqrt{\sigma_j}}\bigg)=r,\quad\mbox{ and let }\quad d_r=\sum_{j\in J_\star^c}\mathbb E\bigg(|G|-\frac{z_r}{\sqrt{\sigma_j}}\bigg)_+^2,\mbox{ where }(\cdot)_+:x\in\bR\mapsto \max\{x,0\}.
\end{align*}For $\lambda>0$ and $|J_\star^c|>N$, define $q=q(\lambda)$ to be the unique solution of the following equation (also known as the calibration equation, \cite{bayatiLASSORiskGaussian2015})
\begin{align}
&\frac{z_q(N-q)}{\sqrt N(N-d_q)^{\frac{1}{2}}}=\frac{\sqrt N\lambda}{\sigma_\xi},
\label{eq:common-q-calibration}
\end{align}whose existence and uniqueness is given by Lemma~\ref{lem:common-calibration} in the Appendix.
We define
\begin{align}
&\mathrm{Var}_{J_\star^c}=\frac{\sigma_\xi^2d_q}{N-d_q} = \frac{\sigma_\xi^2\sum_{j\in J_\star^c}\mathbb E\big(|G|-\frac{z_q}{\sqrt{\sigma_j}}\big)_+^2}{N-\sum_{j\in J_\star^c}\mathbb E\big(|G|-\frac{z_q}{\sqrt{\sigma_j}}\big)_+^2},
\label{eq:common-free-variance}
\end{align}
called the free subspace variance\footnote{In fact, in the regime of interest in this paper, namely in Theorem~\ref{theo:upper_LASSO} and Theorem~\ref{theo:lower_LASSO}, $\Var_{J_\star^c}$ is always equivalent to $\sigma_\xi^2\frac{d_q}{N}$, as demonstrated in the anisotropic case by \eqref{eq:estimate_d_q_Var_anisotropic}, and in the isotropic case since $d_q<0.31N$ in \eqref{eq:common-calibration-identities}. We retain the original form of $\Var_{J_\star^c}$ in order to draw an analogy with mean-field asymptotic theory \cite{montanariMeanFieldAsymptotics2018}.}. Let
\begin{align}\label{eq:def_r_star_without_constant}
    r^2(V_{J_\star},V_{J_\star^c}) =  \sigma_\xi^2\frac{|J_\star|}{N} + \lambda^2\sum_{j\in J_\star}\sigma_j^{-1} + \|\Sigma_{J_\star^c}^{1/2}\bbeta_{J_\star^c}^\star\|_2^2 + \mathrm{Var}_{J_\star^c}.
\end{align}
  
Let $z_N^{\mathrm{all}}>0$ be the unique solution of $\sum_{j=1}^p2\Phi_G(\frac{z_N^{\mathrm{all}}}{\sqrt{\sigma_j}})=N$. We now introduce two sufficients assumptions: one for the upper bound in Theorem~\ref{theo:upper_LASSO} and one for the lower bound result in Theorem~\ref{theo:lower_LASSO}.

\begin{Assumption}\label{assumption:LASSO_upper}    
    For a sufficiently small absolute constant $0<c<\frac{1}{100}$, assume that $1\leq |J_\star|\leq cN$ and $r^2(V_{J_\star},V_{J_\star^c})\leq c\sigma_\xi^2$. Moreover,
    \begin{enumerate}
        \item if $\Sigma=I_p$, assume $N\leq c|J_\star^c|$, and let $\lambda$ satisfy $\frac{\sigma_\xi}{\sqrt{N}}\leq\lambda\leq 2^{-9}\sigma_\xi\sqrt{\frac{\log(e|J_\star^c|/N)}{N}}$ or $\lambda\geq 2^{8}\sigma_\xi\sqrt{\frac{\log(e|J_\star^c|/N)}{N}}$. When  $\frac{\sigma_\xi}{\sqrt N}\leq\lambda\leq2^{-9}\sigma_\xi\sqrt{\frac{\log(\frac{e|J_\star^c|}{N})}{N}},$ we further assume that $|J_\star|\log^2\big(e|J_\star^c|/N\big)\leq cN.$

        \item Alternatively, for a general diagonal covariance matrix $\Sigma$, assume $|J_\star^c|>N$ and $\lambda\ge\frac{8\sigma_\xi z_N^{\mathrm{all}}}{\sqrt N}$.
    \end{enumerate}
\end{Assumption}

Assumption~\ref{assumption:LASSO_upper} covers anistropic designs and a wide range of regularization parameters: any $\lambda \ge \frac{\sigma_\xi}{\sqrt N}$ — except on the interval $(2^{-9}\sigma_\xi\sqrt{\frac{\log(e|J_\star^c|/N)}{N}}, 2^8\sigma_\xi\sqrt{\frac{\log(e|J_\star^c|/N)}{N}})$ whose endpoints are equal up to multiplicative constants. It also granted nothing on the signal, not even sparsity. As a conqequence, it allows to obtain convergence results for the LASSO in scenarii beyond the setups for which the LASSO was designed (see Section~\ref{sec:examples}). Note also that, in the isotropic case, Theorem~\ref{theo:lower_LASSO} below shows that if $r^2(V_{J_\star}, V_{J_\star^c}) > c\sigma_\xi^2$, then $\hat\bbeta$ is inconsistent. Let us now introduce sufficient conditions for the lower bound result from Theorem~\ref{theo:lower_LASSO}. 

\begin{Assumption}\label{assumption:LASSO_lower}
    Fix a sufficiently small absolute constant $0<c<\frac{1}{100}$, assume that $1\leq |J_\star|\leq cN$. Moreover,
    \begin{enumerate}
        \item if $\Sigma=I_p$, then assume $\lambda\geq\frac{\sigma_\xi}{\sqrt N}$ and $N\leq c|J_\star^c|$;
        \item if $\Sigma$ is a general positive diagonal, assume that $|J_\star^c|>N$, $\lambda\ge\frac{8\sigma_\xi z_N^{\mathrm{all}}}{\sqrt N}$, and $r^2(V_{J_\star},V_{J_\star^c})\le c\sigma_\xi^2$.
    \end{enumerate}
\end{Assumption}
The following two theorems are the main results of this paper; they establish high-probability matching upper and lower bounds for the estimation error of the LASSO for a wide range of regularization parameter $\lambda$. Here, $\mathbb{P}^{\otimes N} = \mathbb{P}_{\bX} \otimes \mathbb{P}_{\bxi}$, where $\mathbb{P}_{\bX}$ denotes the probability distribution of $\bX$ and $\mathbb{P}_{\bxi}$ denotes the probability distribution of $\bxi$.
\begin{Theorem}\label{theo:upper_LASSO}
    Grant Assumption~\ref{assumption:LASSO_upper}. There exist absolute constants $c,C$ such that
    \begin{align*}
        \bP^{\otimes N}\bigg\{\norm{\Sigma^{1/2}(\hat\bbeta - \bbeta^\star)}_2^2 \leq Cr^2(V_{J_\star},V_{J_\star^c})\bigg\} \geq 1-C\exp\big(-c\min\{|J_\star|,d_q\}\big).
    \end{align*}
\end{Theorem}

\begin{Theorem}\label{theo:lower_LASSO}
    Grant Assumption~\ref{assumption:LASSO_lower}. There exist absolute constants $c_1,c_2,C>0$ such that
    \begin{align*}
        \bP^{\otimes N}\bigg\{\norm{\Sigma^{1/2}(\hat\bbeta - \bbeta^\star)}_2^2 \geq c_2 r^2(V_{J_\star},V_{J_\star^c})\bigg\} \geq 1-C\exp(-c_1\min\{|J_\star|,d_q\}).
    \end{align*}
\end{Theorem}

The proofs of Theorem~\ref{theo:upper_LASSO} and Theorem~\ref{theo:lower_LASSO} in the isotropic case can be found in Section~\ref{sec:proof_upper_bound} and Section~\ref{sec:proof_lower_bound}, respectively. We present the differences in the proof for the anisotropic design case in Section~\ref{sec:diagonal-upper} and Section~\ref{sec:diagonal-lower}.

\subsection{Examples}
The first example that may come to mind is the classical setup of the LASSO when the signal is $s$-sparse, the design vector is isotropic and the regularization parameter is $\lambda\sim \sigma_\xi\sqrt{\log(ep)/N}$. In that case, we have for the minimax worst case $\sup_{\bbeta^\star\in U_s}r^2(V_{J_\star},V_{J_\star^c})\leq \sigma^2_\xi (s\log(ep)/N)$ where $U_s$ is the set of all $s$-sparse vectors in $\bR^p$. As we can see here, the classical 'LASSO convergence rate' $(s\log(ep)/N)$ is the worst case 'minimax' scenario over model $U_s$; the actual rate  $r^2(V_{J_\star},V_{J_\star^c})$ can be much smaller, of the order of $\sigma_\xi/N$ for $s$-sparse signals with small non-zero coordinates (below $\lambda$) and $\lambda\sim \sigma_\xi/N$. Only an {\it{instance-dependent}} analysis such as the one behind Theorem~\ref{theo:upper_LASSO} and Theorem~\ref{theo:lower_LASSO} can reveal the gap between the worst case  rate (ie the minimax rate over $U_s$) and the 'true', instance-dependent, rate  $r^2(V_{J_\star},V_{J_\star^c})$.

The Lasso has long been the flagship estimator of high-dimensional statistics, where sparsity is the central structural assumption. Consequently, it has historically been tied to this assumption from its inception, just as the ridge estimator has been traditionally associated with regularity conditions. However, the examples below demonstrate that the Lasso can outperform ridge regression even under regularity assumptions. This highlights another value of instance-optimal guarantees—such as Theorems~\ref{theo:upper_LASSO} and~\ref{theo:lower_LASSO}—which can reveal regimes where an estimator achieves minimax optimality over a signal class well beyond its original motivation.

The computation of the following examples may be found in Section~\ref{sec:examples}. The following example is the setting in which the eigenvalues have polynomial decay and the signal satisfies a source condition of order $s$, which is a common assumption in inverse problems, \cite{lecueSharpConvergenceRates2025}.
\begin{Example}\label{example:polynomial_decay}
    Let $(\sigma_j)_{j\geq 1}$ be a sequence such that there exist absolute constants $c,C$ and $\alpha>0$ such that for any $j\leq p$, $cj^{-\alpha}\leq \sigma_j \leq C j^{-\alpha}$; let $(\beta_j^\star)_{j\geq 1}$ be a sequence such that there exists some $s>0$ such that $\sum_{j\geq 1}\sigma_j^{1-s}(\beta_j^\star)^2<\infty$. Fix $\sigma_\xi\in\bR_+$ be a real number. Take $\bbeta^\star = (\beta_j^\star)_{j=1}^p$, and $\Sigma = \diag(\sigma_1,\cdots,\sigma_p)$.  All assertions below hold for sufficiently large integers $N$; the implicit constants may depend on $\alpha$, $s$, $c,C$, and the signal-to-noise ratio $\|\Sigma^{\frac{1-s}{2}}\boldsymbol\beta^\star\|_2/\sigma_\xi$. Then for every large enough $N$,
    \begin{align*}
        \mbox{for } \, \lambda\asymp \frac{\sigma_\xi}{\sqrt{N}}N^{-\frac{\alpha}{2(\alpha s +1)}},\qquad \bP^{\otimes N}\Big(\|\Sigma^{1/2}(\hat\bbeta - \bbeta^\star)\|_2^2 \asymp \sigma_\xi^2 N^{-\frac{\alpha s}{\alpha s+1}} \Big)\geq 1-c\exp(-N^{-\frac{1}{\alpha s+1}}).
    \end{align*}
    We claim that the LASSO with this tuning parameter  attains the same convergence rate as the asymptotic minimax rate for this class of regression problems, see, for instance, \cite{lecueSharpConvergenceRates2025}. In particular, it has no saturation effect compared to ridge regression \cite{gavrilopoulosGeometricalAnalysisKernel2025}, which achieves the slower convergence rate $\sigma_\xi^2 N^{-\frac{\alpha \tilde s}{\alpha \tilde s+1}}$ for the 'saturated' $\tilde s = \min\{s,2\}$.
\end{Example}
We further present a result for the infinite-dimensional LASSO. Here, we allow $\boldsymbol\beta^\star \notin \ell_1$, namely the model-misspecified situation.
\begin{Example}\label{example:infinite}
    Suppose $X=\sum_{j\geq 1}\sqrt{\sigma_j}G_j\in c_0$, that is, $\sqrt{\sigma_j}G_j\to 0$ almost surely, as $j\to\infty$. Suppose $\sum_{j\geq 1}\sigma_j(\beta_j^\star)^2<\infty$. Consider $\hat\bbeta:\bR^N\to\ell_1$, and consider $Y=\sum_{j\geq 1}x_j\beta_j^\star + \xi$, whose convergence is in $L^2(\bP)$. For any finite $p\in\bN_+$, let $\xi^{(p)}=\xi + \sum_{j>p}x_j\beta_j^\star$. Since $\xi^{(p)}$ is independent with $(x_1,\cdots,x_p)$, and $\xi^{(p)}\sim\cN(0,\sigma_\xi^2+\sum_{j>p}\sigma_j(\beta_j^\star)^2)$, Theorem~\ref{theo:upper_LASSO} and Theorem~\ref{theo:lower_LASSO} apply to the signal $\bbeta_{(p)}^\star=(\beta_1^\star,\cdots,\beta_p^\star)^\top\in\bR^p$, the estimator $\hat\bbeta_{(p)}$ in $\bR^p$, the design vector $X_{(p)}=(x_1,\cdots,x_p)^\top$, and the response $Y=\langle X_{(p)},\bbeta_{(p)}^\star\rangle+\xi^{(p)}=\sum_{j\geq 1}x_j\beta_j^\star+\xi$. In particular, this does not necessarily imply that $\bbeta^\star \in \ell_1$; that is, if we regard $\cF_{\ell_1} = \{\vx \mapsto \sum_j x_j\beta_j : \bbeta \in \ell_1\}$ as the model class, then $f^\star(\vx) = \sum_j x_j\beta_j^\star$ may not be in $\cF_{\ell_1}$, so that the model is misspecified in this case. Applying Theorem~\ref{theo:upper_LASSO} and Theorem~\ref{theo:lower_LASSO}, we obtain that for every $p\in\bN_+$, $p>N$, such that Assumption~\ref{assumption:LASSO_upper} and Assumption~\ref{assumption:LASSO_lower} hold, with probability at least $1-C\exp(-c\min\{|J_\star^{(p)}|,d_{q_p}\})$, $\|\Sigma_{(p)}^{1/2}(\hat\bbeta_{(p)}-\bbeta_{(p)}^\star)\|_2^2\asymp r_p^2$, where $d_{q_p}$ is defined in the $p$-dimensional case, and $\Sigma_{(p)}=\bE[X_{(p)}\otimes X_{(p)}]$; and $r_p^2$ is defined similarly, but with $\sigma_\xi^2$ replaced by $\sigma_\xi^2+\sum_{j>p}\sigma_j(\beta_j^\star)^2$. 
    %

    We make $r_p^2$ explicit. Let $\sigma_{\xi,p}^2:=\sigma_\xi^2+\sum_{j>p}\sigma_j(\beta_j^\star)^2$, $J_\star^{(p)}:=\{j\in[p]:\sigma_j|\beta_j^\star|>\lambda\}$ and $J_\star^{(p),c}:=[p]\setminus J_\star^{(p)}$. Let $z_r^{(p)}$, $d_r^{(p)}$ and $z_N^{\mathrm{all},(p)}$ be defined as $z_r$, $d_r$ and $z_N^{\mathrm{all}}$, with $J_\star^c$ replaced by $J_\star^{(p),c}$, and let $q_p$ be the solution of \eqref{eq:common-q-calibration} with $J_\star^c$ replaced by $J_\star^{(p),c}$ and $\sigma_\xi$ replaced by $\sigma_{\xi,p}$; thus $d_{q_p}=d_{q_p}^{(p)}$. Then \eqref{eq:def_r_star_without_constant} reads
    \begin{align}\label{eq:r_p_infinite_example}
        r_p^2=\sigma_{\xi,p}^2\frac{|J_\star^{(p)}|}{N}+\lambda^2\sum_{j\in J_\star^{(p)}}\frac{1}{\sigma_j}+\sum_{j\in J_\star^{(p),c}}\sigma_j(\beta_j^\star)^2+\frac{\sigma_{\xi,p}^2d_{q_p}^{(p)}}{N-d_{q_p}^{(p)}}.
    \end{align}
    The truncation thus converts the part $\sum_{j>p}\sigma_j(\beta_j^\star)^2$ of the bias into additional noise, while the bias on $[p]$ is carried only by the free coordinates $J_\star^{(p),c}$. The limit of \eqref{eq:r_p_infinite_example} as $p\to\infty$ is explicit. Let $J_\star^{\infty}:=\{j\geq1:\sigma_j|\beta_j^\star|>\lambda\}$ and $J_\star^{\infty,c}:=\bN_+\setminus J_\star^\infty$; let $z_r^\infty$, $d_r^\infty$, $z_N^{\mathrm{all},\infty}$ and $q_\infty$ be defined as $z_r$, $d_r$, $z_N^{\mathrm{all}}$ and $q$, with the sums over $J_\star^c$ and $[p]$ replaced by sums over $J_\star^{\infty,c}$ and $\bN_+$, respectively; and let
    \begin{align}\label{eq:r_infty_infinite_example}
        r_\infty^2:=\sigma_\xi^2\frac{|J_\star^{\infty}|}{N}+\lambda^2\sum_{j\in J_\star^{\infty}}\frac{1}{\sigma_j}+\sum_{j\in J_\star^{\infty,c}}\sigma_j(\beta_j^\star)^2+\frac{\sigma_\xi^2d_{q_\infty}^{\infty}}{N-d_{q_\infty}^{\infty}}.
    \end{align}
    Assume that $1\leq|J_\star^\infty|\leq cN$, $r_\infty^2<c\sigma_\xi^2$ and $\lambda>\frac{8\sigma_\xi z_N^{\mathrm{all},\infty}}{\sqrt N}$, where $c$ is the constant of Assumption~\ref{assumption:LASSO_upper}. It is shown in Section~\ref{sec:examples} that $J_\star^\infty$ is finite, that all the quantities in \eqref{eq:r_infty_infinite_example} are well defined, that $r_p^2\to r_\infty^2$ and $d_{q_p}^{(p)}\to d_{q_\infty}^\infty$, and that Assumption~\ref{assumption:LASSO_upper} holds for the $p$-dimensional problem for every sufficiently large $p$. Consequently, Theorem~\ref{theo:upper_LASSO} yields, for every sufficiently large $p$,
    \begin{align*}
        \bP^{\otimes N}\Big\{\|\Sigma_{(p)}^{1/2}(\hat\bbeta_{(p)}-\bbeta_{(p)}^\star)\|_2^2\leq Cr_p^2\Big\}\geq1-C\exp\big(-c\min\{|J_\star^\infty|,d_{q_p}^{(p)}\}\big),
    \end{align*}
    where $c,C$ are the absolute constants of Theorem~\ref{theo:upper_LASSO}. Letting $p\to\infty$, there is an event of probability at least $1-C\exp(-c\min\{|J_\star^\infty|,d_{q_\infty}^\infty\})$ on which the LASSO\footnote{LASSO can also be defined on the infinite-dimensional absolutely summable space $\ell_1$, that is, assuming $X \in c_0$, where $c_0 = \{\vx = (x_j)_{j \ge 1}, x_j \to 0\}$ becomes a separable Banach space when equipped with the $\|\cdot\|_\infty$ norm; and the feasible region of the LASSO optimization problem is replaced by $\ell_1$. In this case, $\langle X, \bbeta \rangle = \sum_{j \ge 1} x_j \beta_j$ can be understood as the dual pairing between $c_0$ and $\ell_1$ (noting that $c_0^* \cong \ell_1$).} problem over $\ell_1$,
    \begin{align*}
        \min_{\bbeta\in\ell_1}\bigg(\frac{1}{2N\lambda}\sum_{i=1}^N\big(Y_i-\langle X_i,\bbeta\rangle\big)^2+\|\bbeta\|_1\bigg),
    \end{align*}
    admits a minimizer $\hat\bbeta$, obtained as a weak-$*$ limit point of $(\hat\bbeta_{(p)})_{p}$, such that
    \begin{align*}
        \sum_{j\geq1}\sigma_j(\hat\beta_j-\beta_j^\star)^2\leq Cr_\infty^2.
    \end{align*}
    In particular, whenever this minimizer is unique, the LASSO over $\ell_1$ achieves the rate $r_\infty^2$, although $\bbeta^\star$ need not belong to $\ell_1$.
\end{Example}

As mentioned above, these two examples demonstrate that the Lasso can achieve optimal statistical performance beyond its original setting. We now turn to the discussion of the main results of this paper. We divide the discussion into two parts:
\begin{enumerate}
    \item a discussion of Theorem~\ref{theo:upper_LASSO} and Theorem~\ref{theo:lower_LASSO} concerning the theory of LASSO; and
    \item a discussion of the FSD method used in proving Theorem~\ref{theo:upper_LASSO} and Theorem~\ref{theo:lower_LASSO} from the perspective of mathematical statistical methodology.
\end{enumerate}
These two parts are placed in Section~\ref{sec:discussion} and Section~\ref{sec:FSD}, respectively.

\section{Why do we still study LASSO?}\label{sec:discussion}

LASSO, or more generally, $\ell_1$-regularization methods, are widely recognized as one of the greatest breakthroughs in mathematical statistics in recent decades, and the mathematical theory behind them has been continuously presented at the International Congress of Mathematicians from 2006 to 2018; see, e.g., \cite{fanStatisticalChallengesHigh2007,vandegeer1RegularizationHighDimensionalStatistical,wainwrightConstrainedFormsStatistical2014,candesMathematicsSparsityFew,montanariMeanFieldAsymptotics2018}. Its theory and applications have also been included in any basic textbook on (modern) mathematical statistics; see, e.g., \cite{koltchinskiiOracleInequalitiesEmpirical2011,buhlmannStatisticsHighDimensionalData2011,giraudIntroductionHighDimensionalStatistics2014,wainwrightHighDimensionalStatisticsNonAsymptotic2019,vershyninHighDimensionalProbabilityIntroduction2018,vandervaartWeakConvergenceEmpirical2023}. 
A natural question is: why does this paper still study LASSO? We claim that despite at least two decades of development, \textbf{its estimation error, even in Gaussian linear regression problem, is still not clearly and completely understood.} In this section, we organize our discussion around the following main themes:
\begin{enumerate}
    \item What are the advantages of the results of this paper over the existing ones?
    \item What new insights do the conclusions of this paper bring to the understanding of LASSO?
\end{enumerate}
Since the anisotropic theory of LASSO is relatively scarce, in this section we focus mainly on the comparison with existing theory in the isotropic case. Only in the end of this section, we will discuss the anisotropic case and answer questions such as: what does the LASSO when the signal is not sparse? What is the interplay between the LASSO and the design covariance matrix?

\subsection{Comparison with existing results}\label{sec:comparision}
Here, we first make a comparison with the existing LASSO theory. The theoretical literature on LASSO can be roughly divided into two strands: non-asymptotic theory and asymptotic theory.

\paragraph{The non-asymptotic theory of LASSO.} Early LASSO theory, even in the optimally tuned Gaussian $s$-sparse linear regression setting, was still unable to provide minimax optimal convergence rates $\sigma_\xi^2\frac{s\log(ep/s)}{N}$, and they only achieve $\sigma_\xi^2\frac{s\log(ep)}{N}$; see, e.g., the classical textbooks \cite[Theorem 7.13]{wainwrightHighDimensionalStatisticsNonAsymptotic2019}, \cite[Corollary 5.3]{giraudIntroductionHighDimensionalStatistics2014}, \cite[Lemma 6.10]{buhlmannStatisticsHighDimensionalData2011}, \cite[Theorem 10.6.1]{vershyninHighDimensionalProbabilityIntroduction2018}, \cite[pp. 456]{vandervaartWeakConvergenceEmpirical2023}. To the best of our knowledge, the earliest result showing that an optimally tuned LASSO can achieve the minimax optimal convergence rate in this problem are \cite{lecueRegularizationSmallballMethod2018, bellecSlopeMeetsLasso2018}. However, these existing non-asymptotic theory of LASSO suffers from the following deficiencies:
    \begin{enumerate}
        \item \textbf{They can only establish upper bounds on the estimation error, and cannot provide matching lower bounds.} Some works have attempted to establish estimation error lower bounds for LASSO, such as \cite{chatterjeeNewPerspectiveLeast2014,vandegeerTightBoundsLasso2018,bellecNoiseBarrierLarge2025}; however, they either match the upper bounds only in very limited settings, \cite{vandegeerTightBoundsLasso2018,bellecNoiseBarrierLarge2025}, or make the convergence rate difficult to compute, \cite{chatterjeeNewPerspectiveLeast2014}.
        
        \item \label{item:only_certain_lambda} \textbf{Conclusions can be obtained only for certain choices of the tuning parameter.} Although \cite{lecueRegularizationSmallballMethod2018,bellecSlopeMeetsLasso2018,bellecNoiseBarrierLarge2025} can obtain the minimax optimal convergence rate over $s$-sparse linear regression problems, and even a sharp phase transition, when $\lambda = \sigma_\xi \sqrt{2\frac{\log(ep/s)}{N}}$, they fail to yield meaningful results when $\lambda$ lies outside this range. This is a significant deficiency, especially when we simultaneously consider the following deficiency in item~\emph{\ref{item:only_for_sparse_vectors}}. On the other side, Theorem~\ref{theo:upper_LASSO} and \ref{theo:lower_LASSO} together reveal the influence of the regularization parameter on the statistical properties of the LASSO for a wide range of this hyper-parameter.

        \item \label{item:only_for_sparse_vectors} \textbf{These theories assume that $\bbeta^\star$ is $s$-sparse and cannot accommodate the case where $\bbeta^\star$ is only approximately sparse.} This deficiency is even more critical for LASSO used in practical applications --- in practice, many works have reported that the signal is never strictly $s$-sparse; see, e.g., \cite{devoreOptimalComputation2007,candesStableSignalRecovery2006,candesMathematicsSparsityFew,shangGeometryStatisticalFeature2026}. To the best of our knowledge, existing finite-sample LASSO theory cannot provide a rule for selecting the theoretically optimal $\lambda$ or a precise characterization of the convergence rate for Gaussian non-sparse linear regression problems. Two exceptions that we are aware of are \cite{bellecSlopeMeetsLasso2018} and \cite{yeRateMinimaxityLasso}. In \cite{bellecSlopeMeetsLasso2018}, the authors provided, for arbitrary $\bbeta^\star$, an oracle inequality for $\frac{1}{N}\|\bX(\hat\bbeta-\bbeta^\star)\|_2^2$ over $\{\bbeta : \|\bbeta\|_0 \le s\}$ at $\lambda \asymp \sigma_\xi \sqrt{\frac{\log(ep/s)}{N}}$. However, this does not imply an upper bound on the estimation error $\|\hat\bbeta-\bbeta^\star\|_2^2$, and more importantly, for a specific choice of  $\bbeta^\star$, there is no evidence that $\lambda \asymp \sigma_\xi \sqrt{\frac{\log(ep/s)}{N}}$ yields the optimal convergence rate; yet as we noted in the previous point, that is, item~\emph{\ref{item:only_certain_lambda}}, their theorem applies only to these $\lambda$, let alone its optimality. In \cite{yeRateMinimaxityLasso}, the authors did examine the upper bound of $\|\hat\bbeta-\bbeta^\star\|_2^2$, but likewise did not study arbitrary tuning parameters, nor did they provide a sharp convergence rate for the estimation error.
    \end{enumerate}
    In contrast, our Theorem~\ref{theo:upper_LASSO} and Theorem~\ref{theo:lower_LASSO} comprehensively achieve the following objectives:
    \begin{enumerate}
        \item From Theorem~\ref{theo:upper_LASSO} and Theorem~\ref{theo:lower_LASSO}, the estimation error $\|\Sigma^{1/2}(\hat\bbeta-\bbeta^\star)\|_2^2$ is (up to multiplicative constants) of the order of $r^2(V_{J_\star}, V_{J_\star^c})$. It is a quantity depending explicitly on the spectrum of $\Sigma$, the signal $\bbeta^\star$ and the regularization parameter $\lambda$.  Therefore, we say that Theorem~\ref{theo:upper_LASSO} and Theorem~\ref{theo:lower_LASSO} provide an \uline{{\it{instance-optimal}} convergence rate} for the estimation error of LASSO in the Gaussian linear regression problem.

        \item The results of Theorem~\ref{theo:upper_LASSO} and Theorem~\ref{theo:lower_LASSO} provide matching upper and lower bounds for almost every $\lambda \gtrsim \frac{\sigma_\xi}{\sqrt N}$, any signal and any diagonal covariance matrix. In particular, $\bbeta^\star$ may have arbitrary patterns beyond the classical structural sparsity assumption for which the LASSO was designed.
    \end{enumerate}

\paragraph{The asymptotic theory of LASSO.} The classical asymptotic theory of LASSO faces the same three deficiencies as the non-asymptotic theory. For instance, the limit theorems for the estimation error established in \cite{meinshausenLassotypeRecoverySparse2009,zhangSparsityBiasLasso2008,yeRateMinimaxityLasso} do not yield sharp convergence rates --- they only achieve (worst case) minimax optimality in $s$-sparse linear regression problems under specific choices of $\lambda$. This situation has been greatly improved by the mean-field asymptotics theory; we refer the reader to \cite{montanariMeanFieldAsymptotics2018} for a comprehensive survey. Such conclusions are typically established via Approximate Message Passing (AMP) or the Convex Gaussian Min-max Theorem (CGMT). In the mean-field asymptotics theory, one usually considers the regime $p, N \to \infty$ with $\frac{N}{p}$ tending to a fixed constant (called the aspect ratio), in order to obtain a precise characterization of the limiting behavior of $\hat\bbeta$; this regime is referred to as the proportional regime. Let $u\geq 0$, $S_u(\cdot)=\sign(\cdot)(|\cdot|-u)_+$ be the soft-max function. Let $G\sim\cN(0,1)$. We refer to the following system of nonlinear equations of $(\tau,z)\in\bR^2$ as the ``mean-field fixed-point equations for the LASSO'', \cite{bayatiLASSORiskGaussian2015}.
    \begin{align}\label{eq:montanari_equation}
        \begin{aligned}
            \tau^2 &= 1+\frac{1}{N}\sum_{j=1}^p\bE\Big(S_{\tau z}\big(\frac{\sqrt{N}}{\sigma_\xi}\beta_j^\star + \tau G\big) - \frac{\sqrt{N}}{\sigma_\xi}\beta_j^\star \Big)^2,\\
            \frac{\sqrt{N}\lambda}{\sigma_\xi} &= \tau z\Big[ 1-\frac{1}{N}\sum_{j=1}^p\bP\big\{ \big|\frac{\sqrt{N}}{\sigma_\xi}\beta_j^\star + \tau G\big| > \tau z \big\}\Big].
        \end{aligned}
    \end{align}In the isotropic situation, the solution of this nonlinear system exists and is unique, and is denoted by $(\tau_\star,z_\star)$. It is proved in \cite{bayatiLASSORiskGaussian2015} that $\|\hat\bbeta - \bbeta^\star\|_2^2 \to \sigma_\xi^2(\tau_\star^2-1)$ under certain limit, thereby even characterizing the exact constant of the limiting behavior of the estimation error; see \cite{bayatiLASSORiskGaussian2015,thrampoulidisPreciseErrorAnalysis2016}. 
    Here, for clarity, we refer to such conclusions as ``precise estimation error'', whereas what we establish in this paper is an ``instance-optimal convergence rate''.  However, this theory has many deficiencies as well.
    \begin{enumerate}
        \item \textbf{This approach is widely regarded as not being particularly helpful for understanding the theoretical nature of the estimation error itself.} Let us temporarily concede the statistical significance of the proportional regime itself. As Ren\'e Thom's dictum \emph{Pr\'edire n'est pas expliquer}, \cite{renePredireNestPas}, is particularly relevant here. The precise asymptotic prediction of the estimation error need not amount to an thorough understanding of the statistical behavior of LASSO. Indeed, \eqref{eq:montanari_equation} has no known closed-form solution, and apart from qualitative properties such as existence and uniqueness (proved in \cite{bayatiLASSORiskGaussian2015}), there is no satisfactory quantitative theory describing the behavior of this solution (\cite{liMinimum$ell_1$normInterpolators2021,bellecOutofsampleErrorEstimation2023}) --- hence such theory does not provide a theoretical answer to how the estimation error varies with the aspect ratio, $\bbeta^\star$, $\lambda$, and the properties of $\sigma_\xi$, but can only be solved numerically. In contrast, Theorem~\ref{theo:upper_LASSO} and Theorem~\ref{theo:lower_LASSO} provide clear interpretability, especially the meaning of each term in $r^2(V_{J_\star}, V_{J_\star^c})$, as well as the quantitative relationships with respect to $N, p, \bbeta^\star, \lambda$, and $\sigma_\xi$.

        \item Moreover, \textbf{the proportional regime itself is also a significant deficiency of this approach.} Prescribing $\frac{N}{p} \to \delta$ as a scaling trajectory does not fully reveal the scaling law, \cite[pp. 14]{wainwrightHighDimensionalStatisticsNonAsymptotic2019}, and this proportional regime is often merely an artifact of the technical tools used in the proofs, rather than a fundamental statistical phenomenon; see \cite{chengDimensionFreeRidge2022}. In contrast, Theorem~\ref{theo:upper_LASSO} and Theorem~\ref{theo:lower_LASSO} are non-asymptotic, and hence allow any form of the limit as $p, N \to \infty$ if one wants to establish asymptotic results.

        \item Although the mean-field asymptotics theory can formally accommodate arbitrary $\bbeta^\star$ (except for \cite{thrampoulidisPreciseErrorAnalysis2016}, which assumes the signal itself to be a random vector), in order to obtain its precise limit, \textbf{the mean-field asymptotics theory, often requires the assumption that as $p, N \to \infty$, $\frac{1}{p}\sum_{j=1}^p \delta_{(\sqrt{N}\beta_j^\star)}$ converges to some prior distribution}, where $\delta$ denotes the Dirac measure; see \cite[Section 1.1]{bellecDebiasingConvexRegularized2021}, and also \cite{bellecObservableAdjustmentsSingleindex2025}. In contrast, Theorem~\ref{theo:upper_LASSO} and Theorem~\ref{theo:lower_LASSO} do not require any assumption that the coordinates of  $\bbeta^\star$ converge to any distribution as $p, N \to \infty$.

        \item Finally, it is also well known that the \textbf{mean-field asymptotics theory can provide precise limits only for a subset of $\lambda$}. Early limit theories established via AMP held only pointwise in the tuning parameter; see \cite{bayatiLASSORiskGaussian2015}. Although \cite{miolaneDistributionLassoUniform2021} extended this limit theory to hold uniformly for certain $\lambda \in [\lambda_{\min}, \lambda_{\max}]$ by using CGMT, we emphasize that $\lambda_{\min}$ and $\lambda_{\max}$ in \cite{miolaneDistributionLassoUniform2021} must be of order $1$ as $N, p \to \infty$, and must stay away from $0$, \cite[Remark 3.1]{miolaneDistributionLassoUniform2021} (hence excluding the optimally tuned parameter $\sigma_\xi \sqrt{(2+\varepsilon) \frac{\log(ep/s)}{N}}$). In the same remark, Miolane and Montanari emphasize that this is the only interesting regime in the proportional setting. Indeed, in \cite[Lemma 4.1]{miolaneDistributionLassoUniform2021}, the authors construct a special $s$-sparse $\bbeta^\star$ in the proportional regime such that $\|\hat\bbeta-\bbeta^\star\|_2^2 \gtrsim \frac{s}{p}\lambda^2$; thus if $\lambda = \omega(1)$, the estimation error has a non-vanishing lower bound and so consistency fails. However, we point out that beyond the proportional regime, one must allow $\lambda \to \infty$ as $N, p \to \infty$ in order to possibly cover the optimally tuned parameter in Gaussian $s$-sparse linear regression. In contrast, our Theorem~\ref{theo:upper_LASSO} and Theorem~\ref{theo:lower_LASSO} hold for almost every $\lambda \gtrsim \sigma_\xi/\sqrt N$.
    \end{enumerate}

\subsection{Insights into LASSO from the FSD Method}\label{sec:insights_into_LASSO}

In this subsection, we introduce the new insights into the understanding of LASSO brought by the FSD method. The  FSD approach reveals a geometric understanding of the estimation mechanism behind high-dimensional estimators (for a comprehensive treatment of the FSD method, we refer the reader to the second-named author's PhD thesis, \cite{shangFeatureSpaceDecomposition2026}). In the case of the LASSO estimator, one can deduce such an understanding from the \emph{FSD's proof} of Theorem~\ref{theo:upper_LASSO} and Theorem~\ref{theo:lower_LASSO}.

\paragraph{LASSO's estimation mechanism.}

The estimation part $\hat\bbeta_{J_\star}$ of LASSO is a regularized empirical risk minimization 'RERM'\footnote{Here, we adopt the general definition of RERM; we do not require the empirical risk to be separable; see \cite{lohRegularizedMestimatorsNonconvexity2015,shangFeatureSpaceDecomposition2026}}.  However, its empirical risk is not the one inherited from $\hat\bbeta$: it behaves like a Huber function of $\norm{\boldsymbol y-\mathbb X_{J_\star}\bbeta}_2$. As a consequence, $\hat\bbeta_{J_\star}$ interpolates between a LASSO estimator and a square-root LASSO estimator \cite{Bellonietal2011}. This is revealed by the following proposition together with the following analysis. We recall that for any $J \subset [p]$, $V_J = \Span(\ve_j : j \in J)$, $\bbeta_J = P_{V_J}\bbeta$, $\bbeta_{J^c}=P_{V_{J^c}}\bbeta$, and $\bX_J = \bX P_{V_J}$, $\bX_{J^c} = \bX P_{V_{J^c}}$.

\begin{Propositionrewrite}{prop:effective-regularization-lasso}[LASSO's FSD]
    This proposition holds for any positive diagonal $\Sigma$ and any $\lambda>0$. Let
\begin{align}\tag{\ref{eq:def_cA}}
    \cA:\vmu\in\bR^N\mapsto \cA[\vmu]\in \argmin\bigg( \frac{1}{2N\lambda}\norm{\vmu - \bX_{J_\star^c}\bbeta}_2^2 + \norm{\bbeta}_1:\, \bbeta\in V_{J_\star^c} \bigg).
\end{align}
Then $\hat\bbeta = \hat\bbeta_{J_\star} + \hat\bbeta_{J_\star^c}$ where $\hat\bbeta_{J_\star^c} = \cA[\boldsymbol y-\mathbb X_{J_\star}\hat\bbeta_{J_\star}]$ and
\begin{align*}
\widehat{\boldsymbol\beta}_{J_\star}\in\operatorname*{argmin}_{\boldsymbol\beta\in V_{J_\star}}\left(P_N\ell_{\boldsymbol\beta}^{(\lambda,{J_\star})}+\|\boldsymbol\beta\|_1\right),
\end{align*}
where
\begin{align}
P_N\ell_\bullet^{(\lambda,{J_\star})}:\boldsymbol\beta\in V_{J_\star}&\mapsto \frac{1}{2N\lambda}\|\boldsymbol y-\mathbb X_{J_\star}\boldsymbol\beta\|_2^2-\frac{1}{2N\lambda}\operatorname{dist}\left(\boldsymbol y-\mathbb X_{J_\star}\boldsymbol\beta, N\lambda(\mathbb X_{J_\star^c}B_1^{J_\star^c})^\circ\right)^2 \tag{\ref{eq:implicit_reg_effective_loss}}
\end{align}and $N\lambda(\mathbb X_{J_\star^c}B_1^{J^c})^\circ = \left\{z\in\bR^N: \norm{\mathbb X_{J_\star^c}^\top z}_\infty\leq N\lambda \right\}$. 
\end{Propositionrewrite}

The estimation part $\widehat{\boldsymbol\beta}_{J_\star}$ of the LASSO with regularization parameter $\lambda$ is performed over $V_{J_\star}$ for $J_\star = \{j\in[p]: \sigma_j |\beta^\star_j|>\lambda\}$ as introduced in \eqref{eq:set_of_selected_features_LASSO}. Hence, the $j$th feature is selected as worth to estimate the signal into that direction when both the variance of the feature (i.e. $\sigma_j$) and the signal strength into that direction (i.e. $|\beta^\star_j|$) are sufficiently strong such that their product goes beyond $\lambda$: '$\sigma_j |\beta^\star_j|>\lambda$'. This is what we expect from an efficient 'model selection based estimator': keep only directions with enough variability because they are the only one that can be used to explain the variability in the output $Y$; and features associated to a strong enough signal because they are the one contributing the most in the sum $\inr{X,\bbeta}$ for the prediction of $Y$: they are the only features that can explain the observed variability in the output and that have a chance to predict it well. The other features (the one with small variability or weak signal strength) should not appear in the final estimator, they would just add useless noise to the final prediction\footnote{Here we assumed that both design vector and noise are centered; this implies that there is no intercept in our analysis. In practice, one usually keep the intercept in the list of features even though it is a feature with no variability. The discussion here goes only for non constant features}. Hence, the set $J_\star$ of selected features is fully orientated toward a better prediction of $Y$ from (a selected set coordinates of) $X$.

 A FSD decomposition of a (high-dimensional) estimator (such as the one described in Proposition~\ref{prop:effective-regularization-lasso} for the LASSO) shed some light onto the mechanism behind its signal estimation process. Indeed, the estimation part of the LASSO is the 'RERM':
\begin{equation}
    \label{eq:esti_part_LASSO}
    \hat\bbeta_{J_\star} \in \argmin_{\boldsymbol\beta \in V_{J_\star}}\left(\frac{1}{2N\lambda}\|\boldsymbol y-\mathbb X_{J_\star}\boldsymbol\beta \|_2^2-\frac{1}{2N\lambda}\operatorname{dist}\left(\boldsymbol y-\mathbb X_{J_\star}\boldsymbol\beta, N\lambda(\mathbb X_{J_\star^c}B_1^{J_\star^c})^\circ\right)^2 + \norm{\boldsymbol\beta}_1\right).
\end{equation}
We call the functional $P_N\ell_\bullet^{(\lambda,J_\star)}: V_{J_\star} \to \mathbb{R}$ defined in \eqref{eq:implicit_reg_effective_loss} the \textbf{effective} empirical risk. That is, although statisticians use the squared risk to define $\hat\bbeta$ via RERM, when actually estimating $\bbeta^\star$ on the relevant set of features indexed by $J_\star$, the estimator  uses a different goodness-of-fit term. 

There are three terms in the objective function defining $\hat\bbeta_{J_\star}$ in \eqref{eq:esti_part_LASSO}: two goodness-of-fit terms and one regularization term. We believe that the $\ell_1$ regularization term  (inherited  from the initial LASSO itself)  has no more utility at this stage, i.e. after the set $J_\star$ has been selected. We perceive it more like a {\it{nuisance regularization term}} responsible for a bias toward $0$ of the coordinates of the LASSO for features (the one in $J_\star$) where there is no need to bias toward zero (since $J_\star$ identifies coordinates with large signal strength); ideally, a better estimator would not have this 'residual regularization term'. The more interesting components in \eqref{eq:esti_part_LASSO}  are the two goodness-of-fit terms gathered in the effective empirical risk
\begin{equation}
    \label{eq:implicit_reg_LASSO}
    P_N\ell_\bullet^{(\lambda,J_\star)}: \boldsymbol\beta\in V_{J_\star} \mapsto \frac{1}{2N\lambda}\left( 
    \operatorname{dist}\left(\boldsymbol y-\mathbb X_{J_\star}\boldsymbol\beta, 0\right)^2-
    \operatorname{dist}\left(\boldsymbol y-\mathbb X_{J_\star}\boldsymbol\beta, N\lambda(\mathbb X_{J_\star^c}B_1^{J_\star^c})^\circ\right)^2 
    \right)
\end{equation}which is forcing the residual $\boldsymbol y-\mathbb X_J\boldsymbol\beta$ to be close to $0$ and far from  $N\lambda(\mathbb X_{J^c}B_1^{J^c})^\circ$ (whereas $0\in N\lambda(\mathbb X_{J^c}B_1^{J^c})^\circ$). In Theorem~\ref{theo:equiv_dist_to_polytope}, we showed that with probability at least $1-4\exp(-cN)$, for all $\bbeta\in V_{J_\star}$,
\begin{equation*}
    P_N\ell_{\boldsymbol\beta}^{(\lambda,J_\star)}  = \frac{1}{2N\lambda} \norm{\boldsymbol y-\mathbb X_{J_\star}\boldsymbol\beta}_2^2  \mbox{ if } \norm{\boldsymbol y-\mathbb X_{J_\star}\boldsymbol\beta}_2\leq R_-
\end{equation*} and
\begin{equation*}
     \frac{R_-}{2N\lambda} \norm{\boldsymbol y-\mathbb X_{J_\star}\boldsymbol\beta}_2 \leq P_N\ell_{\boldsymbol\beta}^{(\lambda,J_\star)}  \leq  \frac{R_+}{N\lambda} \norm{\boldsymbol y-\mathbb X_{J_\star}\boldsymbol\beta}_2  \mbox{ if } \norm{\boldsymbol y-\mathbb X_{J_\star}\boldsymbol\beta}_2\geq 2R_+
\end{equation*} for $R_- = \sqrt{N}\lambda/2$ and $R_+ = N\lambda/(c \sqrt{\log(e|J^c_\star|/N)})$. It shows that the effective goodness-of-fit term is sandwiched between two Huber functions of the Euclidean norm of the residual: it is the MSE close to $0$, and it grows linearly in $\|\boldsymbol y-\mathbb X_J\boldsymbol\beta\|_2$, between $(R_-/(2N\lambda))\|\boldsymbol y-\mathbb X_J\boldsymbol\beta\|_2$ and $(R_+/(N\lambda))\|\boldsymbol y-\mathbb X_J\boldsymbol\beta\|_2$, for residuals of norm at least $2R_+$. As a result $\hat\bbeta_{J_\star}$ interpolates between the LASSO and the  square-root LASSO estimators. This is the typical (unexpected) insight the FSD analysis can provide onto an estimator.


The 'none estimating part' of the LASSO:  $\hat\bbeta_{J_\star^c} = \cA[\boldsymbol y-\mathbb X_{J_\star}\hat\bbeta_{J_\star}]$ is performing a LASSO estimator on the residual $\boldsymbol y-\mathbb X_{J_\star}\hat\bbeta_{J_\star}$ of the estimation part but using only the remaining features, the one in $J_\star^c$, with either low variability or small signal strength. Hence, the predictive power and the utility of these features is low, and we do not expect to extract no more info from the residual of an efficient estimator such as $\hat\bbeta_{J_\star}$ (which has already exhausted all the interesting features) from the remaining low predictive power ones in $\mathbb X_{J_\star^c}$. This explains why we do not expect $\hat\bbeta_{J_\star^c}$ to be an estimator of $\bbeta_{J_\star^c}^\star$. This observation goes beyond the LASSO case, it is common to all high-dimensional estimators that have been analyzed so far using the FSD. It explains why any FSD based analysis of an estimator starts with  $\norm{\Sigma^{1/2}(\hat\bbeta - \bbeta^\star)}_2 \leq \norm{\Sigma^{1/2}_{J_\star}(\hat\bbeta_{J_\star} - \bbeta_{J_\star}^\star)}_2  +  \norm{\Sigma^{1/2}_{J_\star^c}\hat\bbeta_{J_\star^c}}_2 + \norm{\Sigma^{1/2}_{J_\star^c}\bbeta_{J_\star^c}^\star}_2$, which appears to be the right way to start given matching lower bounds such as Theorem~\ref{theo:lower_LASSO} where the {\it{alignment term}} $\norm{\Sigma^{1/2}_{J_\star^c}\bbeta_{J_\star^c}^\star}_2$ appears in the instance-optimal convergence rate $r^2(V_{J_\star}, V_{J_\star^c})$.

\begin{Remark}[Improving the LASSO]
\label{rem:improved_LASSO}
 The FSD decomposition of the LASSO reveals implicit mechanisms behind its statistical performances. Bypassing, we may understand its weakness and try to overcome them. In particular, we have seen that the selected set of features is based on the product $\sigma_j |\beta^\star_j|> \lambda$. There is, however, no evidence that it is the best way to proceed; maybe a better estimator would use two regularization parameters $\lambda_1$ and $\lambda_2$ to select features that comply both conditions $\sigma_j\geq \lambda_1$ and $|\beta^\star_j|\geq\lambda_2$. Another improved version of the LASSO would avoid the 'residual $\ell_1$ regularization term' appearing in \eqref{eq:esti_part_LASSO} that seems to be not any more relevant after a model selection step has been performed and may be responsible for an unwanted  bias toward $0$. 
\end{Remark}

\begin{Remark}[Comparing ridge and LASSO FSDs]
\label{rem:ridge_VS_LASSO}
    The FSD of the ridge estimator (see \cite{tsiglerBenignOverfittingRidge2023,gavrilopoulosGeometricalAnalysisKernel2025}) selects the set of all features given by $J_\star^{ridge}:=\{1,\ldots , k^*\}$ where $k^*$ is called the {\it{estimation dimension}} in \cite{shangFeatureSpaceDecomposition2026} (See Proposition~6 in the supplementary material of \cite{gavrilopoulosGeometricalAnalysisKernel2025}). Both $k^*$ and $J_\star^{ridge}$ depend only on the spectrum of $\Sigma$: it is the same set of selected features whatever the signal $\bbeta^\star$ is: ridge estimator does not adapt its estimation space $V_{J_\star}$ according to the signal strength. This is a major difference with the LASSO that selects features depending on both $\Sigma$ and $\bbeta^\star$: feature $j$ is selected by the LASSO IFF $\sigma_j |\beta^\star_j|>\lambda$. As a consequence, the LASSO learns features adapted to the prediction task of $Y$ given $X$ whereas ridge only keep features in $X$ with large explainability (i.e. variance), regardless of the outcome $Y$. Hence, when one uses ridge the hope is that the signal aligns well with the top $k^*$ singular vectors of $\Sigma$, whereas LASSO does it by itself: it has a model selection skill whose aim is to make the alignment term $\norm{\Sigma_{J_\star^c}\bbeta^\star_{J_\star^c}}_2$ small. That is the reason why we said that LASSO has the {\it{feature learning property}} in \cite{shangGeometryStatisticalFeature2026} whereas ridge does not. On the free part of the feature space $V_{J_\star^c}$, both LASSO and ridge estimators are trying to extract some more information from the residual $\boldsymbol y-\mathbb X_J\hat\bbeta_{J_\star}$ using features in $\bX_{J_\star^c}$ but fail to do so and end up estimating $0$, the 'signal' with the smallest $\ell_1$ and $\ell_2$ norm. In this part of the space, prediction capability of the features is too low to make the goodness-of-fit term small and so only the regularization term is optimized, explaining why $0$ is the target in this part of the space for both estimators. Finally, we end our comparison remark between ridge's FSD and LASSO's FSD by observing that ridge estimation part $\hat\bbeta^{(ridge)}_{J_\star^{ridge}}$ is equivalent to a ridge estimator but with a shifted regularization parameter (compared with the original one), hence, at a first glance, the effect of the projection of the ridge estimator on the estimation space $V_{J_\star^{ridge}}$ may look to be on the regularization part (see \cite{gavrilopoulosGeometricalAnalysisKernel2025}). Whereas for the LASSO, as we have seen above, the regularization term is still w.r.t. the $\ell_1$-norm with the same regularization parameter in \eqref{eq:esti_part_LASSO} but the goodness-of-fit term is the one modified by the projection onto the estimation space $V_{J_\star}$: the MSE is 'replaced' by a 'Huber'-like two regimes function of $\norm{\boldsymbol y-\mathbb X_J\boldsymbol\beta}_2$. Hence, the effects may seem different. However, Proposition~\ref{prop:effective-regularization-lasso} can be extended to any regularization term that can be decomposed in the canonical basis such as $\bbeta\to\sum_{j}|\beta_j|^p$ for $p\geq1$. For $p=2$, i.e. for the ridge estimator, the effective empirical risk is similar to the one we got in \eqref{eq:implicit_reg_LASSO} for  the LASSO, except that $B_1^{J_\star^c}$ is replaced by the Euclidean ball  $B_2^{(J_\star^{ridge})^c}$ so that the convex body onto which the residue is projected in the ridge case is $N\lambda(\mathbb X_{(J_\star^{ridge})^c}B_2^{(J_\star^{ridge})^c})^\circ = \left\{z\in\bR^N: \norm{\mathbb X_{(J_\star^{ridge})^c}^\top z}_2\leq N\lambda \right\}$. It appears that this body has a simple isomorphic description thanks to the Dvoresky-Milmann theorem and to the choice of $J_\star^{ridge}$ (which is exactly defined so that the Dvoresky-Milmann theorem can be applied to $\mathbb X_{(J_\star^{ridge})^c}$): with high probability, $N\lambda(\mathbb X_{(J_\star^{ridge})^c}B_2^{(J_\star^{ridge})^c})^\circ$  is isomorphic to $N \lambda \sqrt{{\rm Tr}(\Sigma_{(J_\star^{ridge})^c})} B_2^N$. As a consequence, on the same event, the effective empirical risk of the estimation part of ridge has been proved to be proportional to the original MSE: it is equivalent to $\lambda + {\rm Tr}(\Sigma_{(J_\star^{ridge})^c})$ (see Proposition~15 in the Supplementary material of \cite{gavrilopoulosGeometricalAnalysisKernel2025}); this can be interpreted into a shift in the regularization parameter.

\end{Remark}

\paragraph{Effective regularization from effective risk.} 
Proposition~\ref{prop:effective-regularization-lasso} itself is a reformulation of the definition of $\hat\bbeta$ via convex analysis --- calling the functional in \eqref{eq:def_effective_loss} the effective empirical risk may of course be viewed as a notational sleight of hand. What property, then, legitimizes this nomenclature? To argue that this reformulation offers a new perspective on LASSO, we need to be able to use this effective empirical risk to obtain the estimation error of $\hat\bbeta_{J_\star}$. 

In classical statistical learning theory, this task is accomplished via the one-sided Restricted Isomorphic Property (RIP). To establish this property, the first step is to examine the lower bound on the curvature of the effective empirical \emph{excess} risk at $\bbeta_{J_\star}^\star$, namely, the (local) Bernstein's property, \cite{bartlettConvexityClassificationRisk2006,bartlettEmpiricalMinimization2006,bartlettLargeMarginClassifiers2003,chinotRobustHighDimensional2021,lecueLearningLinearLoss2024}. Define $P\ell_\bullet^{(\lambda,J_\star)}:\bbeta_{J_\star}\in V_{J_\star}\mapsto \bE[P_N\ell_{\bbeta_{J_\star}}^{(\lambda,J_\star)}|\bX_{J_\star^c}]$ be the effective population risk functional. Here we emphasize that this risk itself is determined by the design matrix $\bX_{J_\star^c}$ on the free subspace $V_{J_\star^c}$, and hence is a random loss. Consequently, its population risk is a conditional expectation given $\bX_{J_\star^c}$. The following property establishes the effective regularization of this effective empirical risk, \cite{bartlettDeepLearningStatistical2021}, see also \cite[pp. 20]{shangFeatureSpaceDecomposition2026}. For any $q\geq 1$, $J\subset[p]$, $\bbeta_J\in V_J$ and $r>0$, let $S_q^J(\bbeta_J;r)=\{\vv\in V_J: \|\vv-\bbeta_J\|_q=r\}$; $B_q^J(\bbeta_J;r)=\{\vv\in V_J:\|\vv-\bbeta_J\|_q\leq r\}$.

\begin{Propositionrewrite}{prop:local-bernstein} Grant Assumption~\ref{assumption:LASSO_upper}.
There exist absolute constants $c_1,c_2,C_\beta>0$ such that with probability at least \(1-4\exp(-c_1N)\) with respect to \(\mathbb X_{J_\star^c}\), for any $r\in\bR$ such that $r \geq C_\beta\|\boldsymbol\beta_{J_\star^c}^\star\|_2$, 
\begin{align*}
    &\forall \bbeta_{J_\star}\in S_2^{J_\star}(\bbeta_{J_\star}^\star;r),\qquad P\mathcal L_{\boldsymbol\beta_{J_\star}}^{(\lambda,J_\star)}\geq c_{\mathrm{LB}}r^2\min\bigg\{\frac1\lambda,\frac{1}{\sqrt{\sigma_\xi^2+r^2}}\sqrt{\frac{N}{\log\big(\frac{e|J_\star^c|}{N}\big)}}\bigg\}.
    \end{align*}
\end{Propositionrewrite}

Proposition~\ref{prop:local-bernstein} establishes a phase transition phenomenon in the local curvature of the functional $P_N\ell_\bullet^{(\lambda,J_\star)}$ at $\bbeta_{J_\star}^\star$. When $\lambda \ge \sigma_\xi \sqrt{\frac{1}{N}\log\big(\frac{e|J_\star^c|}{N}\big)}$, i.e., in the well-regularized regime, the local curvature is governed by $r^2/\lambda$; conversely, when the regularization is too weak (in particular, this includes the case when $\lambda=0$  , and we do not use $\lambda\geq\frac{\sigma_\xi}{\sqrt N}$ in its proof), the local curvature is no longer determined by the tuning parameter prescribed by the statistician, but is instead governed by a parameter depending on $N$, $\sigma_\xi$, and $|J_\star^c|$ --- parameters that are chosen by the estimator itself  . This phenomenon we call self-regularization, \cite{bartlettDeepLearningStatistical2021}, see also \cite[pp. 20]{shangFeatureSpaceDecomposition2026} 
This phenomenon, in which the estimator autonomously decides, based on the geometric structure of the free subspace, whether to adopt self-regularization or the tuning parameter specified by the statistician, we call \textbf{effective regularization}. In the FSD method, such effective regularization is typically established via the Dvoretzky-Milman theorem (\cite{dvoretzkyTHEOREMCONVEXBODIES1959,dvoredskyResultsConvexBodies1961,milmanNewProofTheorem1971}), and Proposition~\ref{prop:local-bernstein} is no exception. Here, to establish self-regularization, i.e., the case where $\lambda$ is small, we employ a fundamental result from Asymptotic Geometric Analysis (AGA), namely Gluskin's theorem, \cite{gluskinExtremalPropertiesOrthogonal1989}, which can be viewed as a one-sided form of the Dvoretzky-Milman theorem for $\ell_2\to\ell_1$ random operators. In general, the Dvoretzky–Milman theorem needs to be used in the less-regularized regime \cite{gavrilopoulosGeometricalAnalysisKernel2025,lecueSharpConvergenceRates2025}.

\paragraph{Uniform convergence of the effective empirical excess risk.}
Another important ingredient in establishing the estimation error of the estimation part $\hat\bbeta_{J_\star}$ is a uniform lower bound on $P_N\cL_{\bbeta_{J_\star}}^{(\lambda,J_\star)}$ over $\bbeta_{J_\star} \in B_1^{J_\star}(\bbeta_{J_\star}^\star;\rho_\star) \cap S_2^{J_\star}(\bbeta_{J_\star}^\star;r_\star)$ for some $\rho_\star>r_\star>0$, where we define $P_N\cL_{\bbeta_{J_\star}}^{(\lambda,J_\star)} = P_N\ell_{\bbeta_{J_\star}}^{(\lambda,J_\star)} - P_N\ell_{\bbeta_{J_\star}^\star}^{(\lambda,J_\star)}$ as the effective empirical excess risk of $\bbeta_{J_\star}$. This conclusion is established by the following two corollary and proposition, respectively.

\begin{Corollary}\label{coro:RIP}
    Grant Assumption~\ref{assumption:LASSO_upper}. Assume $\lambda\geq\sigma_\xi\sqrt{\log(e|J_\star^c|/N)/N}$, $C_\beta\|\boldsymbol\beta_{J_\star^c}^\star\|_2\leq r\leq\sigma_\xi$, and $\rho=2r\sqrt{|J_\star|}$. Then, with $\bP^{\otimes N}$-probability at least $1-C\exp(-c|J_\star|)$, for all $\boldsymbol\beta_{J_\star} \in B_1^{J_\star}(\boldsymbol\beta_{J_\star}^\star;\rho)\cap S_2^{J_\star}(\boldsymbol\beta_{J_\star}^\star;r)$,
\begin{align*}
P_N\mathcal L_{\boldsymbol\beta_{J_\star}}^{(\lambda,J_\star)}\ge c_2\frac{r^2}{\lambda} - C\frac{\sigma_\xi}{\lambda}r\sqrt{\frac{|J_\star|}{N}}.
\end{align*}
\end{Corollary}

\begin{Propositionrewrite}{prop:RIP-small-lambda}Grant Assumption~\ref{assumption:LASSO_upper}, and $\frac{\sigma_\xi}{\sqrt N}\le\lambda\le \sigma_\xi\sqrt{\frac{1}{N}\log(e|J_\star^c|/N)}$. Let $r,\rho$ satisfy $r\leq \rho$ and \eqref{eq:cond_r_rho_prop5} in the full statement of Proposition~\ref{prop:RIP-small-lambda}. Then, there exists a constant $c_{\mathrm{LB}}$, such that with $\bP^{\otimes N}$-probability at least $1-C\exp(-c|J_\star|)$, for all $\boldsymbol\beta_{J_\star}\in B_1^{J_\star}(\boldsymbol\beta_{J_\star}^\star;\rho)\cap S_2^{J_\star}(\boldsymbol\beta_{J_\star}^\star;r)$,
\begin{align}
&P_N\mathcal L_{\boldsymbol\beta_{J_\star}}^{(\lambda,J_\star)}\ge\frac{c_{\mathrm{LB}}}{2}\frac{r^2}{\sqrt{\sigma_\xi^2+r^2}}\sqrt{\frac{N}{\log(e|J_\star^c|/N)}}.\tag{\ref{eq:RIP-small-lambda-lower}}
\end{align}
\end{Propositionrewrite}

Corollary~\ref{coro:RIP} and Proposition~\ref{prop:RIP-small-lambda} establish the one-sided RIP for the effective empirical excess risk, which provides the final ingredient for studying $\hat\bbeta_{J_\star}$. Here, the most important observation --- we emphasize again --- is that we do not establish the one-sided RIP for the empirical excess risk of the squared risk, but rather for this \textbf{effective} empirical excess risk. In Section~\ref{sec:empirical_process_RIP}, we demonstrate the implications of this observation for more general mathematical statistical theory.

\paragraph{\texorpdfstring{$\Var_{J_\star^c}$}{free variance}: mean-field fixed point in the pure-noise case.} At the end of this section, let us discuss the meaning of $\Var_{J_\star^c}$. We point out that $\Var_{J_\star^c}$ is precisely the solution to \eqref{eq:montanari_equation} in the case $\bbeta^\star = \vzero$ and when $p$ is replaced by $|J_\star^c|$. Indeed, at this time, \eqref{eq:montanari_equation} is readily reduced to
\begin{align}\tag{\ref{eq:montanari_equation}}
    \tau^2 &= 1+ \frac{|J_\star^c|}{N}\tau^2\bE\big(|G|-z\big)_+^2,\quad \mbox{ and }\quad \frac{\sqrt{N}\lambda}{\sigma_\xi} = \tau z\big(1-2\frac{|J_\star^c|}{N}\Phi_G(z)\big),
\end{align}where we have used the fact that $S_{\tau z}(\tau G) = \tau S_z(G)$, see, for instance, \eqref{eq:scaling_identity_softmax}; together with the definition of $\Phi_G$. Recall the definition of $z_q=\Phi_G^{-1}(\frac{q}{2|J_\star^c|})$ and $d_q=|J_\star^c|\bE(|G|-z_q)_+^2$, then \eqref{eq:montanari_equation} is equivalent to \eqref{eq:common-q-calibration} and \eqref{eq:common-free-variance}, and consequently, $\Var_{J_\star^c} = \sigma_\xi^2\frac{d_q}{N-d_q} = \sigma_\xi^2(\tau_\star^2-1)$. Namely, $\Var_{J_\star^c}$ is precisely the quantity $\|\hat\bbeta[\bxi]\|_2^2$ characterized by the fixed point of the mean-field fixed-point equations; see the discussion following \eqref{eq:montanari_equation}. This explains the meaning of the variance term $\Var_{J_\star^c}$ on the free subspace: it ``is'' $\|\hat\bbeta_{J_\star^c}[\bxi]\|_2^2$, namely the squared $\ell_2$ norm of the LASSO estimator with $\bX_{J_\star^c}$ as the design matrix, $\bxi$ as the response vector, and $\lambda$ as the tuning parameter. This completes the final piece of the puzzle for understanding Theorem~\ref{theo:upper_LASSO} and Theorem~\ref{theo:lower_LASSO}.

In \cite[pp. 30]{shangFeatureSpaceDecomposition2026}, the second-named author proposed a research direction: for RERM, what is the correct mathematical tool for understanding the energy of the free part of a nonlinear estimator, namely $\|\hat\bbeta_{J_\star^c}\|_2^2$, and suggested that the pursuit of this question might promote the development of new results in AGA. In the ridge regression problem, this tool is the upper side of the Dvoretzky-Milman theorem, \cite{gavrilopoulosGeometricalAnalysisKernel2025}; in the minimum $\ell_q$-norm interpolant estimator, this tool is the newly developed average-case version of the upper side of the Dvoretzky-Milman theorem in \cite{adamczakFeatureSpaceDecompositioninpreperation}. In this paper, the tool we use to study $\|\hat\bbeta_{J_\star^c}\|_2^2$ is CGMT. We do not claim that CGMT is the correct tool for the free-part energy of general nonlinear RERM. Although using CGMT does yield the correct answer --- in the present paper we do not aim to identify the correct AGA tool --- we believe that behind this proof there must necessarily exist a new, as yet unidentified, AGA result. From this perspective, our understanding of $V_{J_\star^c}$ remains incomplete; perhaps we first need to understand the correspondence between the zero-signal case of mean-field asymptotics theory and AGA.

\section{What can FSD offer to mathematical statistics?}\label{sec:FSD}

Next, let us turn to the methodological discussion of this paper from the perspective of mathematical statistics: what does the FSD method used to prove Theorem~\ref{theo:upper_LASSO} and Theorem~\ref{theo:lower_LASSO} bring? In this section, we discuss and revisit the following three points: 1) empirical processes methods, 2) the bias-variance decomposition and 3) the minimax paradigm. To highlight the main points, in this section, we still restrict our discussion to the isotropic design case, unless otherwise stated.

\subsection{Fundamental Proof Strategy: Is Empirical Process Necessary for (R)ERM?}\label{sec:empirical_process_RIP}

The first question we discuss is whether empirical process theory is necessarily required in the study of the estimation error of (R)ERM. Let $\cF \subset L^2(\mathbb{P}_X)$ be a separable class of real-valued measurable functions. We call the random process $(|P_N - P|(f))_{f \in \cF}$ the empirical process indexed by $\cF$, with its marginal given by $|P_N - P|(f) = |\frac{1}{N}\sum_{i=1}^N f(X_i) - \mathbb{E}[f(X)]|$, see, for instance, the standard reference \cite[pp.~129]{vandervaartWeakConvergenceEmpirical2023}.

Empirical process theory is a key tool for studying the estimation error of (R)ERM, dating back at least to the work of Vapnik and Chervonenkis in 1971, \cite{Vapnik1971}, which related the population excess risk of ERM to upper bounds on the empirical process over the entire statistical model --- a method known as the uniform convergence argument. This approach has been refined over the past three decades. In \cite{koltchinskiiRademacherProcessesBounding2000}, the authors pointed out that the upper bound on the population excess risk is determined not by the empirical process indexed by the entire statistical model, but rather by a localized subset near the oracle; this view was further developed in subsequent years, \cite{koltchinskiiRademacherPenaltiesStructural2001,koltchinskiiEmpiricalMarginDistributions2002,bartlettLocalRademacherComplexities2005,koltchinskiiLocalRademacherComplexities2006,koltchinskii2004IMSMedallion2006,blanchardDiscussion2004IMS2006}, giving rise to local Rademacher complexity or other type of local measure of complexity such as Gaussian mean width \cite{lecueLearningSubgaussianClasses2016}  for the control the generalization error of (R)ERM. This has given rise to many classic textbooks on empirical processes, for example, \cite{vandegeerEmpiricalProcessesMEstimation2006,koltchinskiiOracleInequalitiesEmpirical2011,vandervaartWeakConvergenceEmpirical2023}.

As the title of this subsection suggests, in this subsection we challenge this belief. To this end, we first briefly review how classical statistical learning theory uses empirical processes to study the estimation error of (R)ERM, why it faces inherent difficulties, and then we introduce the new insights brought by the FSD method. Since the empirical process method is well known to provide only upper bounds on the estimation error, in this subsection we restrict our consideration to comparisons of upper bounds.

\subsubsection{Empirical processes methods in classical statistical learning theory}\label{sec:empirical_process_in_classical_stats}

In this section, we take LASSO as an example to briefly introduce the localization argument and fixed points as upper bounds for its estimation error. In this subsubsection, we proceed quickly; we refer the reader to \cite{lecueRegularizationSmallballMethod2018} for the relevant details.

By definition, the regularized empirical excess risk at the LASSO is non-positive: $\frac{1}{2N\lambda}\|\vy-\bX\hat\bbeta\|_2^2 - \frac{1}{2N\lambda}\|\vy-\bX\bbeta^\star\|_2^2 + (\|\hat\bbeta\|_1 - \|\bbeta^\star\|_1) \le 0$. Therefore, any $\bbeta \in \mathbb{R}^p$ satisfying $\frac{1}{2N\lambda}\|\vy-\bX\bbeta\|_2^2 - \frac{1}{2N\lambda}\|\vy-\bX\bbeta^\star\|_2^2 + (\|\bbeta\|_1 - \|\bbeta^\star\|_1) > 0$ cannot be $\hat\bbeta$. Hence, to prove that $\|\hat\bbeta-\bbeta^\star\|_2^2 \le r^2$ for some $r>0$ holds with high probability, it suffices to show that with high probability, every $\bbeta \notin B_2(\bbeta^\star; r)$ has positive regularized empirical excess risk; then, conversely, we must have $\hat\bbeta \in B_2(\bbeta^\star; r)$. This is called the localization technique, one of the basic techniques in statistical learning theory, \cite{koltchinskiiOracleInequalitiesEmpirical2011}.

We use, without explanation, the homogeneity argument and the sparsity equation (see, for instance, \cite{lecueRegularizationSmallballMethod2018}) to reduce the problem to two empirical processes on the following index sets. Let $$\partial K = \partial(B_1(\bbeta^\star;\rho)\cap B_2(\bbeta^\star;r)) = (B_1(\bbeta^\star;\rho)\cap S_2(\bbeta^\star;r))\sqcup (S_1(\bbeta^\star;\rho)\cap B_2(\bbeta^\star;r)),$$ where $\rho$ is some parameter determined by the so-called sparsity equation, \cite{lecueRegularizationSmallballMethod2018}, such that for any $\bbeta\in S_1(\bbeta^\star;\rho)\cap B_2(\bbeta^\star;r)$, there holds $\|\bbeta\|_1-\|\bbeta^\star\|_1\geq \frac{4}{5}\rho$, then
\begin{align*}
    &\frac{1}{N}\inf_{\bbeta\in \partial K}\Big(\|\vy-\bX\bbeta\|_2^2 - \|\vy-\bX\bbeta^\star\|_2^2\Big) \geq \inf_{\bbeta \in \partial K}\Big(\frac{1}{N}\sum_{i=1}^N\langle X_i,\bbeta-\bbeta^\star\rangle^2\Big) - 2\sup_{\bbeta \in \partial K}\Big(\frac{1}{N}\sum_{i=1}^N\xi_i\langle X_i,\bbeta-\bbeta^\star\rangle\Big),
\end{align*}
which corresponds to the two empirical processes: quadratic process and the multiplier process \cite{lecueLearningSubgaussianClasses2016}. The fixed points defined by these two random processes --- the quadratic fixed point $r_Q$ and the multiplier fixed point $r_M$ --- are characterizing the estimation error of many estimators, including LASSO, e.g., \cite{lecueLearningSubgaussianClasses2016,mendelsonLocalVsGlobal2017,alquierEstimationBoundsSharp2019}. Let $0<\delta_Q,\delta_M<\frac{1}{4}$, $\theta>0$, and let
\begin{align*}
    r_Q(\rho,\delta_Q) &= \inf\bigg(r>0:\bP\bigg\{\inf_{\bbeta\in B_1(\bbeta^\star;\rho)\cap S_2(\bbeta^\star;r)}\frac{1}{N}\sum_{i=1}^N\langle X_i,\bbeta-\bbeta^\star\rangle^2 \geq \theta r^2 \bigg\} \geq 1-\delta_Q\bigg),\mbox{ and}\\
    r_M(\rho,\delta_M) &= \inf\bigg(r>0:\bP\bigg\{\sup_{\bbeta\in B_1(\bbeta^\star;\rho)\cap B_2(\bbeta^\star;r)}\bigg|\frac{1}{N}\sum_{i=1}^N\xi_i\langle X_i,\bbeta-\bbeta^\star\rangle\bigg| \leq \frac{1}{4}\theta r^2 \bigg\} \geq 1-\delta_M\bigg).
\end{align*}Therefore, to prove that $\frac{1}{2N\lambda}\|\vy-\bX\bbeta\|_2^2 - \frac{1}{2N\lambda}\|\vy-\bX\bbeta^\star\|_2^2 + (\|\bbeta\|_1 - \|\bbeta^\star\|_1) > 0$ for any $\bbeta\in\partial K$, it suffices to take $r>\max\{r_Q(\rho,\delta_Q),r_M(\rho,\delta_M)\}$ such that the following two cases hold simultaneously:
\begin{enumerate}
    \item \label{item:almost_sparse_main} when $\bbeta\in B_1(\bbeta^\star;\rho)\cap S_2(\bbeta^\star;r)$, $\theta r^2 > \frac{1}{2}\theta r^2 + \lambda \rho$; and
    \item \label{item:well_spread_main} when $\bbeta\in S_1(\bbeta^\star;\rho)\cap B_2(\bbeta^\star;r)$, $\frac{4}{5}\lambda\rho>\frac{1}{2}\theta r^2$.
\end{enumerate}
Therefore, the smallest $r$ satisfying the above two equations and exceeding $\max\{r_Q(\rho,\delta_Q), r_M(\rho,\delta_M)\}$ serves as a $1-\delta_Q-\delta_M$ high-probability upper bound for the estimation error (up to squaring). This is the general method in classical statistical learning theory for obtaining upper bounds on the estimation error of (R)ERM via empirical processes.

\subsubsection{Bottleneck of empirical process methods}\label{sec:bottleneck_empirical_process}

The classical method described above has the following significant obstacles.
\begin{enumerate}
    \item An upper bound can only be given for $\lambda$ belonging to a very narrow range. This is because the system of inequalities $\theta r^2 > \frac{1}{2}\theta r^2 + \lambda \rho$ and $\frac{4}{5}\lambda\rho > \frac{1}{2}\theta r^2$ simultaneously constrains $\lambda$ to be neither too large nor too small --- this is precisely the reason, mentioned in Section~\ref{sec:comparision}, item~\emph{\ref{item:only_certain_lambda}}, why the existing non-asymptotic theory of LASSO holds only for certain choices of $\lambda$.

    \item Most crucially, for a considerable number of estimators, estimating $r_Q(\rho,\delta_Q)$ and $r_M(\rho,\delta_M)$ is itself very difficult.
\end{enumerate}

To control the empirical process, one needs to impose  strong tail assumptions on $X$. To control the upper bound of the empirical process on a general index set such as $B_1(\bbeta^\star;\rho)\cap B_2(\bbeta^\star;r)$, the most effective method is generic chaining, \cite{talagrandUpperLowerBounds2021}, which applies well to stochastic processes whose increments have uniformly bounded sub-Gaussian norms on the index set. Although there have been efforts to avoid the ``strong coupling'' of generic chaining and sub-Gaussian increments, these methods either currently show no promise of being extendable to more general index sets such as \cite{latalaBoundingSupremaCanonical2023,zhivotovskiyDimensionfreeBoundsSums2024}, or still require  strong assumptions on the marginal distribution of the design vector such as \cite{mendelsonUpperBoundsProduct2016,mendelsonMultiplierProcessesWeak2017}. Beyond the sub-gaussian increments case, chaining methods give rise to complexity measure of local subsets that are usually hard to handle.

Let us temporarily assume that the empirical process satisfies this sub-Gaussian increment property. This further leads to two scenarios: whether or not the equivalence between the $\psi_2(\mathbb P_X)$ norm and the $L^2(\mathbb{P}_X)$ norm over $\cF$ holds uniformly on the index set \cite{lecueLearningSubgaussianClasses2016}.
\begin{enumerate}[label = (2\Alph*.)]
    \item \label{item:estimate_gamma_2} If the equivalence between the $\psi_2(\mathbb P_X)$ norm and the $L^2(\mathbb{P}_X)$ norm holds (ie $\cF$ is a {\it{sub-gaussian class}} as defined in \cite{lecueLearningSubgaussianClasses2016}), then the quadratic process and the multiplier process (after assuming certain conditions on $\xi$) can be expressed in terms of Talagrand's $\gamma_2$-functional, \cite[Definition 2.7.3]{talagrandUpperLowerBounds2021}, see also \cite{bednorzConcentrationChainingMethod2014,dirksenTailBoundsGeneric2015,mendelsonUpperBoundsProduct2016,mendelsonMultiplierProcessesWeak2017}. However, \textbf{estimating the $\gamma_2$ functional is in general a well-known challenge} in particular for {\it{local}} subsets of the model. When the index set is the unit ball of some metric space, it is closely related to the metric properties of that metric space, such as 2-convexity and 1-unconditionality, making sharp estimates highly nontrivial; see, e.g., van Handel's interpolation and entropy contraction principle, \cite{vanhandelChainingInterpolationConvexity2018,vanhandelChainingInterpolationConvexity2018a}, see also \cite[Chapter 4]{talagrandUpperLowerBounds2021}. In particular, in the anisotropic case, classical methods such as \cite{lecueRegularizationSmallballMethod2018} require estimating $\gamma_2(\rho\Sigma^{1/2}B_1^p\cap rB_2^p)$ — a well-known challenging problem for non diagonal covariance matrices $\Sigma$. We emphasize that while the $\gamma_2$ functional can certainly be estimated crudely by Dudley's entropy, as in \cite{vandervaartWeakConvergenceEmpirical2023,vandegeerEmpiricalProcessesMEstimation2006}, even when the index set is an ellipsoid, Dudley's entropy integral introduces an extra logarithmic loss (\cite[Exercise 2.13.3]{talagrandUpperLowerBounds2021}), and hence cannot yield optimal (up to multiplicative constant) upper bound on the generalization error in general.

    \item \label{item:other_complexity_parameters} If the uniform equivalence between the $\psi_2(\mathbb P_X)$ norm and the $L^2(\mathbb{P}_X)$ norm over the entire model $\cF$ does not hold, there are two further cases: either a finite-order norm equivalence holds, or one assumes uniform boundedness in the $L^\infty(\mathbb{P}_X)$ norm, namely the bounded setup, \cite{koltchinskiiLocalRademacherComplexities2006,alquierEstimationBoundsSharp2019}. 
    In this case, one typically needs to use Gin\'e-Zinn's symmetrization inequality, \cite[Section 11.8]{talagrandUpperLowerBounds2021}, together with Talagrand's contraction principle, \cite[Lemma 6.4.4]{talagrandUpperLowerBounds2021}, to convert the original quadratic process and multiplier process into the problem of upper bounding a localized Rademacher complexity. However, first, the contraction principle itself may already introduce slackness in $r_Q(\rho,\delta_Q)$; more critically, there is currently no general good method for estimating the localized Rademacher complexity other than using Dudley's entropy, which in general also introduces slackness.
\end{enumerate}

\subsubsection{The low-dimensional estimation space \texorpdfstring{$V_{J_\star}$}{VJstar}}
The two obstacles of Section~\ref{sec:bottleneck_empirical_process} are not removed by a new probabilistic tool. They disappear because of the structure revealed by the FSD: once the features in $J_\star$ have been selected, estimation takes place on the low-dimensional subspace $V_{J_\star}$, where no genuine regularization is needed. We explain this in three steps below; the classical road-map of Section~\ref{sec:empirical_process_in_classical_stats} (localization, decomposition into a multiplier and a quadratic component, generic chaining) is followed throughout, and only the objects to which it is applied change.

\paragraph{The processes are not empirical processes, but they have sub-Gaussian increments.} By Proposition~\ref{prop:effective-regularization-lasso}, $\hat\bbeta_{J_\star}$ minimizes the effective empirical risk $P_N\ell_\bullet^{(\lambda,J_\star)}$ plus $\|\cdot\|_1$ over $V_{J_\star}$. The relevant stochastic process is therefore the effective empirical excess risk $(P_N\cL_{\bbeta_{J_\star}}^{(\lambda,J_\star)})_{\bbeta_{J_\star}\in V_{J_\star}}$, and not the empirical excess risk of the squared loss. Through the projection onto the random polytope $(\bX_{J_\star^c}B_1^{J_\star^c})^\circ$ in \eqref{eq:implicit_reg_effective_loss}, this process depends on the whole sample at once and is not a sum of $N$ independent terms; the empirical process machinery (symmetrization, contraction, localized Rademacher complexities) does not apply to it directly. As in the proof of Proposition~\ref{prop:RIP}, we split it, conditionally on $\bX_{J_\star^c}$, as $(P_N-P)\cL_{\bbeta_{J_\star}}^{(\lambda,J_\star)}=P_N\cM_{\bbeta_{J_\star}}+(P_N-P)\cQ_{\bbeta_{J_\star}}$, see \eqref{eq:modified-process} and \eqref{eq:original-modified-decomposition}. The multiplier component $P_N\cM_{\bbeta_{J_\star}}=\langle-\vlambda^\star[\vzeta],\bX_{J_\star}(\bbeta_{J_\star}-\bbeta_{J_\star}^\star)\rangle$ is linear in the Gaussian matrix $\bX_{J_\star}$, with a multiplier $\vlambda^\star[\vzeta]$ independent of $\bX_{J_\star}$. The quadratic component $P_N\cQ_{\bbeta_{J_\star}}$ is nonnegative, and the $\frac{1}{N\lambda}$-Lipschitz continuity of $\vlambda^\star[\cdot]$, together with Gaussian concentration for Lipschitz functions, shows that its increments are sub-Gaussian, with sub-Gaussian norm proportional to $\|\bbeta_{J_\star}-\bbeta'_{J_\star}\|_2$; see \eqref{eq:noise-gradient} and \eqref{eq:design-gradient}. This is the situation commonly met with dependent data, for instance in time series, where the empirical process structure is lost and most of the work consists in proving a sub-Gaussian increment property, after which generic chaining \cite[Theorem 2.7.13]{talagrandUpperLowerBounds2021}, \cite[Theorem 3.2]{dirksenTailBoundsGeneric2015} applies. What is specific to the present setting is that such a process arises from i.i.d.\ data, as a consequence of the FSD.

\paragraph{No $\ell_1$-localization is needed on $V_{J_\star}$.} Since $|J_\star|\leq cN$, estimation on $V_{J_\star}$ is a low-dimensional problem, for which the expected rate is the parametric rate $\sigma_\xi^2|J_\star|/N$ of least squares in dimension $|J_\star|$. Such a rate does not require the $\ell_1$-ball coming from the regularization: in Section~\ref{sec:proof_upper_bound} we take $\rho_\star=2r_\star\sqrt{|J_\star|}$, so that $B_1^{J_\star}(\bbeta_{J_\star}^\star;\rho_\star)\cap B_2^{J_\star}(\bbeta_{J_\star}^\star;r_\star)=B_2^{J_\star}(\bbeta_{J_\star}^\star;r_\star)$, since $\|\vv\|_1\leq\sqrt{|J_\star|}\|\vv\|_2$ for $\vv\in V_{J_\star}$. The localization set is a Euclidean sphere, and there is no need to separate the almost sparse and the well-spread directions as in items~\emph{\ref{item:almost_sparse_main}} and~\emph{\ref{item:well_spread_main}} of Section~\ref{sec:empirical_process_in_classical_stats}. As a consequence, the constraint on $\lambda$ coming from the simultaneous validity of these two items disappears; the range of $\lambda$ in Theorem~\ref{theo:upper_LASSO} is instead dictated by the geometry of the free subspace $V_{J_\star^c}$, through Proposition~\ref{prop:local-bernstein}. The $\ell_1$-norm is still present in the objective defining $\hat\bbeta_{J_\star}$, where it produces the shrinkage bias $\lambda^2|J_\star|=\lambda^2\sum_{j\in J_\star}\sigma_j^{-1}$ of $r^2(V_{J_\star},V_{J_\star^c})$, but it plays no role in the localization nor in the complexity of the problem.

\paragraph{The complexity is that of a Euclidean ball.} Because the increments are sub-Gaussian with respect to $\ell_2$ and the index set is $B_2^{J_\star}(\bbeta_{J_\star}^\star;r)$, the only complexity parameter that enters is
\begin{equation*}
    \gamma_2\big(B_2^{J_\star}(\vzero;r),\|\cdot\|_2\big)\asymp\bE\sup_{\vv\in B_2^{J_\star}(\vzero;r)}\langle\vg,\vv\rangle=r\,\bE\|\vg\|_2\asymp r\sqrt{|J_\star|},
\end{equation*}
where $\vg$ is a standard Gaussian vector in $V_{J_\star}$, by the majorizing measure theorem \cite[Theorem 2.10.1]{talagrandUpperLowerBounds2021}. Combining this bound with the curvature of order $r^2/\lambda$ from Proposition~\ref{prop:local-bernstein} and with the sub-Gaussian norm of order $\frac{\sqrt{\sigma_\xi^2+r^2}}{\lambda\sqrt N}\|\bbeta_{J_\star}-\bbeta_{J_\star}'\|_2$ of the increments (Proposition~\ref{prop:RIP}), and with the bound $\|\bbeta_{J_\star}\|_1-\|\bbeta_{J_\star}^\star\|_1\geq-\sqrt{|J_\star|}\,\|\bbeta_{J_\star}-\bbeta_{J_\star}^\star\|_2$ on the excess regularization, the localization argument requires $\frac{r^2}{\lambda}\gtrsim\frac{\sigma_\xi}{\lambda}r\sqrt{\frac{|J_\star|}{N}}+r\sqrt{|J_\star|}$ for $r\leq\sigma_\xi$, that is, $r^2\gtrsim\sigma_\xi^2\frac{|J_\star|}{N}+\lambda^2|J_\star|$, which are the first two terms of $r^2(V_{J_\star},V_{J_\star^c})$ in \eqref{eq:def_r_star_without_constant}. In contrast, the classical analysis of Section~\ref{sec:empirical_process_in_classical_stats} requires sharp estimates of $\gamma_2(\rho B_1^p\cap rB_2^p,\|\cdot\|_2)$ in the isotropic case, and of $\gamma_2(\rho\Sigma^{1/2}B_1^p\cap rB_2^p,\|\cdot\|_2)$ in the anisotropic case, that is, of the complexity of an interpolation body between two norms; this is the difficulty described in item~\emph{\ref{item:estimate_gamma_2}} of Section~\ref{sec:bottleneck_empirical_process}. The same simplification occurs for the minimum $\ell_q$-norm interpolant in \cite{adamczakFeatureSpaceDecompositioninpreperation}, where the Gaussian mean width of $\Sigma_{J_\star}^{1/2}B_q^{J_\star}(\vzero;\rho)\cap B_2^{J_\star}(\vzero;r)$ is replaced by the bound $r\sqrt{|J_\star|}$. The simplicity of the index set would persist under heavier tails: if the increments were only sub-exponential with respect to $\|\cdot\|_\infty$, the relevant functional would be $\gamma_1(B_2^{J_\star}(\vzero;r),\|\cdot\|_\infty)\asymp r\sqrt{|J_\star|}$, which again requires no entropy computation, and the corresponding chaining bound \cite[Theorem 3.2]{dirksenTailBoundsGeneric2015} would avoid the difficulty of item~\emph{\ref{item:other_complexity_parameters}} of Section~\ref{sec:bottleneck_empirical_process}.\footnote{The upper bound follows from Lata\l{}a's majorizing measure theorem for canonical processes with exponential variables, \cite{latalaSudakovMinorationPrinciple1997}, together with the bound on the Euclidean norm of a vector with independent symmetric exponential coordinates. For the lower bound, let $J=J_\star$, $S=\{|J|^{-1/2}\vveps:\vveps\in\{\pm1\}^J\}\subset B_2^J$, let $(\cT_s)_{s\geq0}$ be an admissible sequence, and let $s_0$ satisfy $2^{s_0}\leq\frac12|J|<2^{s_0+1}$. Distinct points of $S$ are at $\ell_\infty$-distance $2|J|^{-1/2}$, so every open ball $B_\infty(\va,|J|^{-1/2})$ contains at most one point of $S$. Since $|\cT_{s_0}|\leq2^{2^{s_0}}<2^{|J|}=|S|$, there exists $\vt\in S$ with $d_\infty(\vt,\cT_{s_0})\geq|J|^{-1/2}$, and hence $\sup_{\vu\in B_2^J}\sum_{s\geq0}2^sd_\infty(\vu,\cT_s)\geq2^{s_0}|J|^{-1/2}\geq\frac14\sqrt{|J|}$. Taking the infimum over admissible sequences gives $\gamma_1(B_2^J,\|\cdot\|_\infty)\geq\frac14\sqrt{|J|}$.}

In summary, the difficult geometry has not disappeared: it has been moved to the free subspace $V_{J_\star^c}$, where it is handled by Gluskin's theorem (Proposition~\ref{prop:local-bernstein}) and by the analysis of $\Var_{J_\star^c}$. On the estimation subspace $V_{J_\star}$, the analysis reduces to a sub-Gaussian increments process indexed by a Euclidean ball of dimension $|J_\star|\leq cN$, whose complexity is explicit.

\subsection{Redefine the Bias-Variance Decomposition via FSD}\label{sub:redefine_the_bias_variance_decomposition_via_fsd}

In this subsection, we discuss how the FSD method leads us to rethink the bias-variance decomposition of the estimation error. It is well known that the bias-variance decomposition is often used to study upper and lower bounds on the expectation of the estimation error. The classical bias-variance decomposition of the estimation error of a given estimator  $\hat\bbeta$ of $\bbeta^\star$ is
\begin{align*}
    \bE_{\bX,\bxi}\|\langle X,\hat\bbeta - \bbeta^\star\rangle\|_{L^2(\bP_X)}^2 = \|\langle X,\bE_{\bX,\bxi}\hat\bbeta-\bbeta^\star\rangle\|_{L^2(\bP_X)}^2 + \bE_{\bX,\bxi}\|\langle X,\hat\bbeta - \bE_{\bX,\bxi}\hat\bbeta\rangle\|_{L^2(\bP_X)}^2,
\end{align*}where the first term on the right-hand-side is the bias of $\hat\bbeta$ while the second its variance.

The advantage of the bias-variance decomposition defined in this manner is that it is essentially a reformulation of the definition of the variance of a random variable $\bE\zeta^2=(\bE\zeta)^2+\Var(\zeta)$ --- and hence provides an exact decomposition. However, its drawback is also evident: it is difficult to assign statistical meaning to the bias and variance in terms of $\hat\bbeta$, $\bbeta^\star$, and $\sigma_\xi$, especially when $\hat\bbeta$ is an RERM; \textbf{characterizing the effect of the tuning parameter $\lambda$ on the bias and variance is a well-known challenge}. Even if we revert to the fixed design setting and define the fixed-design bias and variance as $\|\bX(\mathbb{E}_{\bxi}\hat\bbeta-\bbeta^\star)\|_2^2$ and $\mathbb{E}_{\bxi}\|\bX(\hat\bbeta-\mathbb{E}_{\bxi}\hat\bbeta)\|_2^2$, the monotonicity of the fixed-design bias and variance with respect to the tuning parameter is, to the best of our knowledge, only known when $\hat\bbeta$ is some special linear operator, such as when $\hat\bbeta$ is ridge regression; see \cite[Exercise 5.5.7]{giraudIntroductionHighDimensionalStatistics2014}. When $\hat\bbeta: \mathbb{R}^N \to \mathbb{R}^p$ is a nonlinear operator, it is generally difficult to characterize the monotonicity of the bias and variance with respect to $\lambda$; see \cite{bellecNoiseBarrierLarge2025}.

It is well known that understanding the monotonicity of bias and variance with respect to the tuning parameter is often a prerequisite for adaptive selection of the tuning parameter, as in Lepskii's principle, \cite{lepskiiProblemAdaptiveEstimation1991}. Adaptive (data-dependent) tuning parameter selection is one of the most important topics in mathematical statistics, \cite{lepskiTheoryAdaptiveEstimation2023}. To address this issue, \cite{bellecNoiseBarrierLarge2025} proposed using the large signal bias and the noise barrier as substitutes for the bias-variance decomposition, and established the monotonicity of these two quantities. Although this approach has enjoyed some success, there is currently no evidence that the large signal bias and the noise barrier themselves genuinely constitute a sharp/optimal decomposition of the estimation error --- even up to multiplicative constants. 

In this subsection, we show how to use FSD to redefine the bias-variance decomposition, which we call the FSD bias-variance decomposition, so that the FSD-bias and FSD-variance naturally exhibit monotonicity with respect to the tuning parameter. This should be viewed as yet another instance of FSD as a theoretical framework; see \cite[Section 1.5.4]{shangFeatureSpaceDecomposition2026}. We define
\begin{align*}
    \mathrm{BIAS}^{\mathrm{FSD}} &= \mathrm{BIAS}_{J_\star} + \mathrm{BIAS}_{J_\star^c},\qquad \mbox{ where }\mathrm{BIAS}_{J_\star} = \lambda^2\sum_{j\in J_\star}\sigma_j^{-1},\, \mathrm{BIAS}_{J_\star^c} = \|\Sigma_{J_\star^c}^{1/2}\bbeta_{J_\star^c}^\star\|_2^2,\\
    \mathrm{VAR}^{\mathrm{FSD}} &= \mathrm{VAR}_{J_\star} + \mathrm{VAR}_{J_\star^c},\qquad\mbox{ where }\,\mathrm{VAR}_{J_\star} = \sigma_\xi^2\frac{|J_\star|}{N},\mbox{ and } \mathrm{VAR}_{J_\star^c}\mbox{ is defined in }\eqref{eq:common-free-variance}.
\end{align*}
The advantage of this definition is that $\|\Sigma^{1/2}(\hat\bbeta-\bbeta^\star)\|_2^2 \asymp \mathrm{BIAS}^{\mathrm{FSD}} + \mathrm{VAR}^{\mathrm{FSD}}$ holds with high probability thanks to Theorem~\ref{theo:upper_LASSO} and Theorem~\ref{theo:lower_LASSO}, thereby preserving the form of the bias-variance decomposition up to multiplicative constant. Moreover, both $\mathrm{BIAS}^{\mathrm{FSD}}$ and $\mathrm{VAR}^{\mathrm{FSD}}$ are decomposed into components on the estimation subspace $V_{J_\star}$ and the free subspace $V_{J_\star^c}$, respectively, and hence admit a clear interpretability, as we have already explained in Section~\ref{sec:insights_into_LASSO}. The most important is the following fact. Here, for any $\lambda>0$, let $q(\lambda_-), z_{q(\lambda_-)}$ and $d_{q(\lambda_-)}$ be the left limit at $\lambda$; and $J_\star(\lambda_-)=\{j: \sigma_j|\beta_j^\star|\geq \lambda\}$.
\begin{Proposition}\label{prop:monotonicity_bias_variance}
    Suppose there exists $\lambda_0\geq0$ such that $|J_\star^c(\lambda_0)|\geq 100N$. Then for $\lambda\in(\lambda_0,\infty)$, $\mathrm{BIAS}^{\mathrm{FSD}}$ is non-decreasing for every positive diagonal $\Sigma$. Moreover, $\lambda\in\{\lambda_0\leq \lambda<\infty\}\mapsto \mathrm{VAR}^{\mathrm{FSD}}$ is non-increasing, if and only if for any $\lambda>\lambda_0$ such that $\lambda = \sigma_j|\beta_j^\star|$ for some $j\in[p]$, there holds
    \begin{align*}
        d_{q(\lambda)} - d_{q(\lambda_-)} \leq \Big|\Big\{j\in[p]:\, \sigma_j|\beta_j^\star|=\lambda\Big\}\Big|\Big(1-\frac{d_{q(\lambda)}}{N}\Big)\Big(1-\frac{d_{q(\lambda_-)}}{N}\Big).
    \end{align*}
\end{Proposition}
\begin{proof}[Proof of Proposition~\ref{prop:monotonicity_bias_variance}]
    It is straightforward to check that $\mathrm{BIAS}^{\mathrm{FSD}} = \sum_{j=1}^p\min\{\frac{\lambda^2}{\sigma_j},\sigma_j(\beta_j^\star)^2\}$, indicating that $\mathrm{BIAS}^{\mathrm{FSD}}$ is increasing. For $\mathrm{VAR}^{\mathrm{FSD}}$, the monotonicity of $\mathrm{VAR}^{\mathrm{FSD}}$ is a consequence of Lemma~\ref{lem:monotonicity_total_fsd_variance} together with the observation that $\lambda\mapsto|J_\star^c(\lambda)|$ is non-decreasing.
\end{proof}
Therefore, under the assumptions of Proposition~\ref{prop:monotonicity_bias_variance}, we claim that by redefining the bias-variance decomposition via the FSD method (which we call the \uline{FSD bias-variance decomposition}), $\mathrm{BIAS}^{\mathrm{FSD}}$ and $\mathrm{VAR}^{\mathrm{FSD}}$ recover monotonicity, thereby restoring the statistical advantage of the bias-variance decomposition. In fact, this is a process opposite to the classical bias-variance decomposition --- the classical approach first provides the decomposition via the definition of variance, and then studies its statistical meaning; the FSD bias-variance decomposition, on the other hand, first establishes the statistical meaning (namely Theorem~\ref{theo:upper_LASSO} and Theorem~\ref{theo:lower_LASSO}), and then defines the bias and variance in retrospect. Remarkably, this approach nevertheless preserves the statistical advantage of the bias-variance decomposition. For general spectral methods, such as ridge regression, gradient descent, and gradient flow, we established similar monotonicity results for the FSD bias-variance decomposition in \cite{lecueSharpConvergenceRates2025}. Let us instance it for the LASSO, the proof of the following example may be found in Section~\ref{sec:examples}.
\begin{Examplecontinue}{example:polynomial_decay}
Consider the setup of Example~\ref{example:polynomial_decay}. Let $\lambda_0\asymp \frac{\sigma_\xi}{\sqrt{N}}$.  When $N$ is large enough, on $\{\lambda\in\bR:\lambda\geq \lambda_0\}$, $\lambda\mapsto \mathrm{BIAS}^{\mathrm{FSD}}$ is non-decreasing, while $\lambda\mapsto\mathrm{VAR}^{\mathrm{FSD}}$ is non-increasing.
\end{Examplecontinue}
Finally, let us emphasize that $\mathrm{BIAS}^{\mathrm{FSD}}$ and $\mathrm{VAR}^{\mathrm{FSD}}$ themselves still depend on the unobservable $\Sigma$, $\bbeta^\star$ and $\sigma_\xi$, and hence cannot serve as an adaptive tuning parameter selection principle in practice. Converting them, via their empirical versions, into a truly implementable, data-dependent tuning parameter selection principle is an interesting direction for future research, \cite{blanchardLepskiiPrincipleSupervised2019}.

\subsection{Minimax Lower Bounds and the {\it{instance-optimal}} results from FSD}\label{sec:minimax_vs_problem_specific}

In this subsection, we discuss insights that the FSD method brings to the discussion of minimax lower bounds versus {\it{instance-optimal}} estimation error lower bounds. The minimax lower bound is one of the most fundamental paradigms in mathematical statistics \cite{tsybakovIntroductionNonparametricEstimation2009}. Whether an estimator achieves minimax optimality on a given class model is a standard question in mathematical statistics; see, e.g., the classic textbook \cite[Chapter 2]{tsybakovIntroductionNonparametricEstimation2009}, \cite[Chapter 3]{giraudIntroductionHighDimensionalStatistics2014}, \cite[Chapter 6]{gineMathematicalFoundationsInfiniteDimensional2015}, \cite[Chapter 15]{wainwrightHighDimensionalStatisticsNonAsymptotic2019}. 

Let us first recall the motivation for minimax optimality paradigm. We take the linear regression problem as an example, and let $\cB \subset \mathbb{R}^p$ be a bounded set. 
\begin{enumerate}
    \item The primary motivation for minimax optimality is that \textbf{statisticians typically do not know the exact $\bbeta^\star$, but only know a set $\cB$ such that $\bbeta^\star \in \cB$}, for example, $\cB$ being some ellipsoid; see \cite[Chapter 2]{tsybakovIntroductionNonparametricEstimation2009}. This viewpoint can be traced back to one of the origins of minimax optimality, namely the work of A. Wald, \cite[pp. 173]{waldStatisticalDecisionFunctions1949}, according to which $\bbeta^\star$ is chosen by nature, and the statistician cannot know nature's choice; when there is insufficient prior information about the preference of this choice, investigating the worst-case estimation error $\sup\{\|\hat\bbeta-\bbeta^\star\|_2^2 : \bbeta^\star \in \cB\}$ becomes a natural approach; see \cite[Chapter 1, Section II, pp. 18]{bergerStatisticalDecisionTheory1985}.

    \item This worst-case behavior itself does not necessarily entail an outer infimum over all measurable functions, and the motivation for this outer infimum is \textbf{to establish a negative theory} --- namely, to set a stopping rule such that once a statistician has constructed an estimator $\hat\bbeta$ satisfying $\inf_{\tilde\bbeta}(\sup(\bE_{\bbeta^\star}\|\tilde\bbeta-\bbeta^\star\|_2^2 : \bbeta^\star \in \cB) : \tilde\bbeta\mbox{ measurable}) \gtrsim \bE_{\bbeta^\star}\|\hat\bbeta-\bbeta^\star\|_2^2$ for all $\bbeta^\star\in\cB$, one can cease further investigation of the problem class from the perspective of estimation error --- because no estimator can have a uniformly smaller estimation error on that class $\cB$, see, for instance, \cite[pp. 55]{giraudIntroductionHighDimensionalStatistics2014}.

    \item The minimax paradigm has some other motivations in modern mathematical statistics, which we briefly mention but do not discuss here. In recent years, the minimax paradigm has evolved to incorporate computational constraints into the outer infimum, for example by restricting to computationally tractable algorithms, thereby seeking a statistical-computational trade-off, where the object of comparison is replaced by all ``computable'' estimators; see the surveys \cite{wainwrightConstrainedFormsStatistical2014,carpentierStatisticalComputationalChallenges2026}. In addition, there is also the viewpoint of using the minimax paradigm as a tool to characterize when a statistical problem is unsolvable; see \cite[Section 3.4.3]{giraudIntroductionHighDimensionalStatistics2014}.
\end{enumerate}

Of course, since its introduction, discussions of this paradigm have never ceased, e.g, \cite[Section 5.1]{bergerStatisticalDecisionTheory1985} and \cite{JohnstoneICM}. In the following, we introduce one of the main criticism. \textbf{The minimax lower bound is pessimistic for evaluating the statistical properties of an estimator}. In \cite{bellecNoiseBarrierLarge2025}, the author points out that the minimax lower bound only proves the existence of some $\bbeta_{\mathrm{worst}}^\star$ within the class of problems under consideration, such that no estimator can have a smaller estimation error. However, this lower bound provides no information about the estimation error of the estimator for other $\bbeta^\star$ in the same class, in particular, for the 'true one' behind the data generating process (DGP). Hence, the minimax approach studies a scenario that is likely not the one we want to study given for our dataset: it solves another problem for another set of data generated by another DGP.  On the contrary, {\it{instance-optimal}} results  attempt to characterize matching upper and lower bounds on the estimation error of an estimator for {\underline{any}} $\bbeta^\star$ and more generally any DGP. Characterizing the estimation error itself is one of the main motivations for the recent development of mean-field asymptotics theory; see \cite[pp. 2993]{montanariMeanFieldAsymptotics2018}. On top of that, the minimax approach may also drive the choice of hyperparameters only against the worst case situation whereas an {\it{instance-optimal}} approach adapts the choice of hyperparameters to the true DGP.

Here, we claim that by using the FSD method to obtain an instance-optimal convergence rate, we can achieve a compromise between the two viewpoints. Namely, combining Theorem~\ref{theo:lower_LASSO} and Theorem~\ref{theo:upper_LASSO}, the statistician can claim to have obtained a correct understanding of the estimation error for estimating $\bbeta^\star$ using $\hat\bbeta$, rather than merely a pessimistic lower bound. If the statistician knows only a set $\cB \subset \mathbb{R}^p$ such that $\cB \ni \bbeta^\star$, then combining Theorem~\ref{theo:lower_LASSO} and Theorem~\ref{theo:upper_LASSO}, the statistician can claim that the worst-case estimation error of $\hat\bbeta$ on $\cB$ \emph{is} given by $\sup(r^2(V_{J_\star}, V_{J_\star^c})(\bbeta^\star) : \bbeta^\star \in \cB)$ -- if there is only an upper bound and no matching lower bound on the estimation error, then the statistician can only claim that the worst-case estimation error of $\hat\bbeta$ on $\cB$ is \emph{at most} the established uniform upper bound. This precisely explains our claim: the FSD method not only provides a estimation error lower bound, but also serves as a compromise between the estimation error lower bound and the minimax lower bound
, by providing the worst-case estimation error for a class of problems. Admittedly, using the FSD method to establish lower bounds on the estimation error is considerably more difficult than using information-theoretic methods to establish minimax lower bounds. Nevertheless, we believe that this provides a potential tool for going beyond the minimax paradigm, and that such a tool should be particularly suitable for studying the refined statistical properties of estimators that are already widely used in important problems, rather than for examining the statistical properties of newly proposed estimators designed to address new statistical problems.

On top of that the {\it{instance-dependent}} convergence rate approach (as in Theorem~\ref{theo:upper_LASSO}) allows for a better choice of hyper-parameters, in theory first, such as choosing $\lambda$ minimizing the convergence rate $r^2(V_{J_\star}, V_{J_\star^c})$ and, secondly, with a data-driven choice of $\lambda$ -- based, for instance, via a Lepski's method adapted to the FSD bias/variance decomposition -- in order to achieve a better convergence rate than the one obtained from a minimax approach that suggest to take $\lambda\sim \sigma_\xi \sqrt{\log(ep)/N}$ whatever $\bbeta^\star$ and $\Sigma$ are. 

Finally, the minimax approach only tells about the worst case behavior of an estimator. It does not reveal its estimation mechanism in general unlike the highlights provided by the FSD approach as depicted by results such as Proposition~\ref{prop:effective-regularization-lasso} and its subsequent interpretation. Having such an understanding helps design improved estimators as in Remark~\ref{rem:improved_LASSO} or data-driven hyperparameters as in Section~\ref{sub:redefine_the_bias_variance_decomposition_via_fsd}.

\section{Summary and Outlook}\label{sec:summary}

In this paper, we use the FSD method to establish a sharp instance-optimal convergence rate for the estimation error of LASSO in the Gaussian linear regression problems for any given signal, diagonal covariance matrix and tuning parameter above $\sigma_\xi/\sqrt{N}$. We then discuss the connections of this result to the existing non-asymptotic theory of LASSO and to mean-field asymptotics theory, and provide new insights into LASSO brought by the FSD method --- effective risk and effective regularization. Furthermore, we put forward three insights of the FSD method for mathematical statistics:
\begin{enumerate}
    \item the low-dimensional estimation problem that is restricted to $V_{J_\star}$ bypasses several issues related to the computation of local complexity parameters previously encountered in the classical approach derived from the empirical processes' toolbox,
    \item the bias-variance decomposition of (R)ERM should be redefined through FSD for the design of data-driven choices of regularization parameter; and 
    \item  FSD is a fine-grained analysis of the generalization error that provides a deeper understanding of the estimation mechanism of high-dimensional estimators that overcomes the limitations of the worst case analysis of the minimax paradigm.
\end{enumerate}

\section*{Acknowledgments}

{\footnotesize 

This paper was completed while the second-named author was affiliated with CREST-ENSAE, Institut Polytechnique de Paris, following the defense of his doctoral dissertation. ZS thanks Matthieu Lerasle for his questions during the defense, which inspired the reflections on the bias-variance decomposition in this paper.

\textbf{LLM Usage.} The authors developed the proofs presented in this paper with the assistance of AI-based tools. The AI models used were ChatGPT 5.6 Pro Sol and Claude Opus 5.5. The authors take full responsibility for the correctness of the proofs.}

\appendix

\section{Proof of Theorem~\texorpdfstring{\ref{theo:upper_LASSO}}{upper bound}}\label{sec:proof_upper_bound}

We first prove Theorem~\ref{theo:upper_LASSO} for $\Sigma=I_p$. The modifications for diagonal $\Sigma$ are given in Section~\ref{sec:diagonal-upper}.

The first step of the feature space decomposition method  \cite{lecueGeometricalViewpointBenign2024,gavrilopoulosGeometricalAnalysisKernel2025,lecueSharpConvergenceRates2025,adamczakFeatureSpaceDecompositioninpreperation,shangFeatureSpaceDecomposition2026} is to decompose the prediction risk according to the decomposition of the feature space $\bR^p = V_J \oplus^\top V_{J^c}$:
\begin{equation}\label{eq:main_decomp_risk}
    \norm{\hat\bbeta - \bbeta^\star}_2^2\leq 2\norm{\hat\bbeta_{J_\star} - \bbeta_{J_\star}^\star}_2^2  + 4\norm{\hat\bbeta_{J_\star^c}}_2^2 + 4\norm{\bbeta_{J_\star^c}^\star}_2^2.
\end{equation}Then, we show that with $\bP^{\otimes N}$-high probability, there hold simultaneously, 
\begin{align*}
    \norm{\hat\bbeta_{J_\star} - \bbeta_{J_\star}^\star}_2^2\leq C r^2(V_{J_\star},V_{J_\star^c}), \mbox{ and }  \norm{\hat\bbeta_{J_\star^c}}_2^2\leq C r^2(V_{J_\star},V_{J_\star^c}).
\end{align*}We call the third term $\norm{\bbeta_{J_\star^c}^\star}_2$ appearing  \eqref{eq:main_decomp_risk} the {\it{alignment term}}. It is a deterministic quantity that measures the 'mass/energy' of the signal on $V_{J^c}$; that is where the estimator will not estimate the signal. There is nothing to do on this term and it appears as it is in the final convergence rate $r^2(V_{J_\star},V_{J_\star^c})$. This term is given by the problem and the only hope is that it is small, i.e. that the signal is well aligned with the sub-space $V_{J_\star}$ of the feature space.  In the case of the LASSO, $J_\star$ is selected such that the alignment term is minimized (this is in great contrast with ridge and many other spectral methods, see \cite{lecueSharpConvergenceRates2025}).

With a bit abuse of notation, we abbreviate $\bP^{\otimes N}$ by $\bP$ in the proof. For any non empty, closed and convex set $K$ we denote for all $x$,
\begin{equation*}
   \operatorname{Proj}_{K}(x)\in\argmin_{z\in K}\norm{z-x}_2 \mbox{ and } \operatorname{dist}(x, K)=\min_{z\in K}\norm{x-z}_2.
\end{equation*}

\subsection{Proof of Proposition~\ref{prop:effective-regularization-lasso}}

\begin{Proposition}[Effective regularization]\label{prop:effective-regularization-lasso}
For every \(J\subset[p]\), let \(V_J=\operatorname{span}(\boldsymbol e_j:j\in J)\). This proposition holds for any arbitrary positive diagonal $\Sigma$. For any $\lambda>0$, let
\begin{align}\label{eq:def_cA}
    \cA:\vmu\in\bR^N\mapsto \cA[\vmu]\in \argmin\bigg( \frac{1}{2N\lambda}\norm{\vmu - \bX_{J^c}\bbeta_{J^c}}_2^2 + \norm{\bbeta_{J^c}}_1:\, \bbeta_{J^c}\in V_{J^c} \bigg).
\end{align}
Then $\hat\bbeta = \hat\bbeta_J + \hat\bbeta_{J^c}$ where $\hat\bbeta_{J^c} = \cA[\vy-\bX_J\hat\bbeta_J]$ and
$
\widehat{\boldsymbol\beta}_J\in\operatorname*{argmin}_{\boldsymbol\beta_J\in V_J}\left(P_N\ell_{\boldsymbol\beta_J}^{(\lambda,J)}+\|\boldsymbol\beta_J\|_1\right),
$
for 
\begin{align}
P_N\ell_\bullet^{(\lambda,J)}:\boldsymbol\beta_J\in V_J&\mapsto\frac{1}{2N\lambda}\norm{\boldsymbol y-\mathbb X_J\boldsymbol\beta_J- \bX_{J^c}\cA[\vy - \bX_J\bbeta_J] }_2^2+\norm{\cA[\vy-\bX_J\bbeta_J]}_1\label{eq:def_effective_loss}\\
&=\max\left\{\langle\boldsymbol\lambda,\boldsymbol y-\mathbb X_J\boldsymbol\beta_J\rangle-\frac{N\lambda}{2}\|\boldsymbol\lambda\|_2^2:\, \vlambda\in\bR^N:\, \|\bX_{J^c}^\top\vlambda\|_\infty\leq1\right\}.\label{eq:dual_problem_effective_loss}
\end{align}
Moreover, the following facts hold.
\begin{enumerate}
    \item \label{item:effective_lambda_star} For every \(\boldsymbol\beta_J\in V_J\), \eqref{eq:dual_problem_effective_loss} has the unique solution
    \begin{align*}
    \boldsymbol\lambda^\star[\boldsymbol y-\mathbb X_J\boldsymbol\beta_J]=\operatorname{Proj}_{(\mathbb X_{J^c}B_1^{J^c})^\circ}\left(\frac{\boldsymbol y-\mathbb X_J\boldsymbol\beta_J}{N\lambda}\right),
    \end{align*}where $(\cdot)^\circ$ denotes the polar body, \cite[Definition 1.1.5]{artstein-avidanAsymptoticGeometricAnalysis2015}. Moreover, for any $\vmu_1,\vmu_2\in\bR^N$, $\|\vlambda^\star[\vmu_1] - \vlambda^\star[\vmu_2]\|_2\leq (1/(N\lambda))\|\vmu_1-\vmu_2\|_2$, and $\|\D \vlambda^\star[\cdot](\vmu)\|_{\op}\leq 1/(N\lambda)$ almost surely for any $\vmu\in\bR^N$, where for any $\vmu\in\bR^N$, $\D\vlambda^\star[\cdot](\vmu):\bR^N\to\bR^N$ is the Fr\'echet derivative of $\vlambda^\star[\cdot]$ evaluated at $\vmu$.
    \item For any $\lambda>0$ and any $\bbeta\in V_J$,
    \begin{equation}\label{eq:implicit_reg_effective_loss}
    P_N\ell_{\bbeta}^{(\lambda,J)} = \frac{1}{2N\lambda}\|\boldsymbol y-\mathbb X_J\boldsymbol\beta\|_2^2-\frac{1}{2N\lambda}\operatorname{dist}\left(\boldsymbol y-\mathbb X_J\boldsymbol\beta, N\lambda(\mathbb X_{J^c}B_1^{J^c})^\circ\right)^2.
    \end{equation}
    
    \item If \(\widehat{\boldsymbol\nu}\) attains the defining minimum of \eqref{eq:def_cA} for $\vmu = \boldsymbol y-\mathbb X_J\boldsymbol\beta_J$, then
    \begin{align}\label{eq:scaled-dual-solution}
    \boldsymbol y-\mathbb X_J\boldsymbol\beta_J-\mathbb X_{J^c}\widehat{\boldsymbol\nu}=N\lambda\boldsymbol\lambda^\star[\boldsymbol y-\mathbb X_J\boldsymbol\beta_J]=\operatorname{Proj}_{N\lambda(\mathbb X_{J^c}B_1^{J^c})^\circ}(\boldsymbol y-\mathbb X_J\boldsymbol\beta_J).
    \end{align}

    \item There holds
    \begin{align}
    \nabla P_N\ell_{\boldsymbol\beta_J}^{(\lambda,J)}=-\mathbb X_J^\top\boldsymbol\lambda^\star[\boldsymbol y-\mathbb X_J\boldsymbol\beta_J].
    \label{eq:effective-gradient}
    \end{align}
\end{enumerate}
\end{Proposition}

\begin{proof} The decomposition of the $\ell_1$-norm in the canonical basis  \(\|\boldsymbol\beta\|_1=\|\boldsymbol\beta_J\|_1+\|\boldsymbol\beta_{J^c}\|_1\) and division of the objective by \(\lambda\) when $\lambda>0$ prove the FSD of the LASSO as $\hat\bbeta = \hat\bbeta_J + \hat\bbeta_{J^c}$ for the definitions of $\hat\bbeta_J$ and $\hat\bbeta_{J^c}$ as above.

Let us now show that \eqref{eq:dual_problem_effective_loss} holds and prove {\it{item.1}} as well.  For every \(\boldsymbol\beta_J\in V_J\), \(\boldsymbol\nu\in V_{J^c}\), and \(\boldsymbol\lambda\in(\mathbb X_{J^c}B_1^{J^c})^\circ\), ie $\|\bX_{J^c}^\top\vlambda\|_\infty\leq1$,
\begin{align*}
&\langle\boldsymbol\lambda,\boldsymbol y-\mathbb X_J\boldsymbol\beta_J\rangle-\frac{N\lambda}{2}\|\boldsymbol\lambda\|_2^2 = \langle\boldsymbol\lambda,\boldsymbol y-\mathbb X_J\boldsymbol\beta_J-\mathbb X_{J^c}\boldsymbol\nu\rangle-\frac{N\lambda}{2}\|\boldsymbol\lambda\|_2^2+\langle\mathbb X_{J^c}^\top\boldsymbol\lambda,\boldsymbol\nu\rangle\\
&\leq\|\boldsymbol\lambda\|_2\|\boldsymbol y-\mathbb X_J\boldsymbol\beta_J-\mathbb X_{J^c}\boldsymbol\nu\|_2-\frac{N\lambda}{2}\|\boldsymbol\lambda\|_2^2+\|\mathbb X_{J^c}^\top\boldsymbol\lambda\|_\infty\|\boldsymbol\nu\|_1\\
&\leq\frac{1}{2N\lambda}\|\boldsymbol y-\mathbb X_J\boldsymbol\beta_J-\mathbb X_{J^c}\boldsymbol\nu\|_2^2+\|\mathbb X_{J^c}^\top\boldsymbol\lambda\|_\infty\|\boldsymbol\nu\|_1\leq\frac{1}{2N\lambda}\|\boldsymbol y-\mathbb X_J\boldsymbol\beta_J-\mathbb X_{J^c}\boldsymbol\nu\|_2^2+\|\boldsymbol\nu\|_1,
\end{align*}where we used that $\max_{t\geq0}(t \square - (N\lambda/2)t^2) = \square^2/(2N\lambda)$ for $\square = \|\boldsymbol y-\mathbb X_J\boldsymbol\beta_J-\mathbb X_{J^c}\boldsymbol\nu\|_2$.

By convexity and coercivity, the minimum in \eqref{eq:def_cA} is attained and by the KKT conditions (see, for instance, \cite[Theorem VI.2.2.1]{hiriart-urrutyConvexAnalysisMinimization1993}), every solution \(\widehat{\boldsymbol\nu}\) of \eqref{eq:def_cA} for $\vmu = \boldsymbol y-\mathbb X_J\boldsymbol\beta_J$ satisfies $\frac{1}{N\lambda}\mathbb X_{J^c}^\top(\boldsymbol y-\mathbb X_J\boldsymbol\beta_J-\mathbb X_{J^c}\widehat{\boldsymbol\nu})\in\partial\|\widehat{\boldsymbol\nu}\|_1.$
Consequently,
\begin{align*}
    \norm{\frac{1}{N\lambda}\mathbb X_{J^c}^\top(\boldsymbol y-\mathbb X_J\boldsymbol\beta_J-\mathbb X_{J^c}\widehat{\boldsymbol\nu})}_\infty \leq 1,\mbox{ and }\left\langle\frac{\mathbb X_{J^c}^\top(\boldsymbol y-\mathbb X_J\boldsymbol\beta_J-\mathbb X_{J^c}\widehat{\boldsymbol\nu})}{N\lambda},\widehat{\boldsymbol\nu}\right\rangle=\|\widehat{\boldsymbol\nu}\|_1.
\end{align*}Now combining the above displayed formula with $\vy - \bX_J\bbeta_J = (\vy - \bX_J\bbeta_J - \bX_{J^c}\hat\vnu) + \bX_{J^c}\hat\vnu$, one obtains
\begin{align*}
&\left\langle\frac{\boldsymbol y-\mathbb X_J\boldsymbol\beta_J-\mathbb X_{J^c}\widehat{\boldsymbol\nu}}{N\lambda},\boldsymbol y-\mathbb X_J\boldsymbol\beta_J\right\rangle-\frac{N\lambda}{2}\left\|\frac{\boldsymbol y-\mathbb X_J\boldsymbol\beta_J-\mathbb X_{J^c}\widehat{\boldsymbol\nu}}{N\lambda}\right\|_2^2\\
&= \frac{1}{N\lambda}\norm{\vy - \bX_J\bbeta_J - \bX_{J^c}\hat\vnu}_2^2 + \left\langle \frac{1}{N\lambda}\bX_{J^c}^\top\big(\vy - \bX_J\bbeta_J - \bX_{J^c}\hat\vnu\big),\hat\vnu \right\rangle - \frac{1}{2N\lambda}\|\vy - \bX_J\bbeta_J - \bX_{J^c}\hat\vnu\|_2^2\\
&=\frac{1}{2N\lambda}\|\boldsymbol y-\mathbb X_J\boldsymbol\beta_J-\mathbb X_{J^c}\widehat{\boldsymbol\nu}\|_2^2+\|\widehat{\boldsymbol\nu}\|_1.
\end{align*}This shows that \eqref{eq:dual_problem_effective_loss} holds and since the dual objective is strictly concave on its feasible set, the unique solution to \eqref{eq:dual_problem_effective_loss} is
\begin{align*}
\boldsymbol\lambda^\star[\boldsymbol y-\mathbb X_J\boldsymbol\beta_J]=\frac{\boldsymbol y-\mathbb X_J\boldsymbol\beta_J-\mathbb X_{J^c}\widehat{\boldsymbol\nu}}{N\lambda}.
\end{align*}
Completing the square gives
\begin{align}\label{eq:completing_squares_1}
&\langle\boldsymbol\lambda,\boldsymbol y-\mathbb X_J\boldsymbol\beta_J\rangle-\frac{N\lambda}{2}\|\boldsymbol\lambda\|_2^2=\frac{1}{2N\lambda}\|\boldsymbol y-\mathbb X_J\boldsymbol\beta_J\|_2^2-\frac{N\lambda}{2}\left\|\boldsymbol\lambda-\frac{\boldsymbol y-\mathbb X_J\boldsymbol\beta_J}{N\lambda}\right\|_2^2,
\end{align}
then taking argmax to \eqref{eq:completing_squares_1} over $\vlambda\in (\bX_{J^c}B_1^{J^c})^\circ$, gives
\begin{align*}
    &\argmax_{\vlambda\in (\bX_{J^c}B_1^{J^c})^\circ}\bigg( \langle\boldsymbol\lambda,\boldsymbol y-\mathbb X_J\boldsymbol\beta_J\rangle-\frac{N\lambda}{2}\|\boldsymbol\lambda\|_2^2\bigg) = \argmax_{\vlambda\in (\bX_{J^c}B_1^{J^c})^\circ} \bigg( \frac{1}{2N\lambda}\|\boldsymbol y-\mathbb X_J\boldsymbol\beta_J\|_2^2-\frac{N\lambda}{2}\left\|\boldsymbol\lambda-\frac{\boldsymbol y-\mathbb X_J\boldsymbol\beta_J}{N\lambda}\right\|_2^2\bigg)\\
    &=\argmin_{\vlambda\in (\bX_{J^c}B_1^{J^c})^\circ} \bigg( \left\|\boldsymbol\lambda-\frac{\boldsymbol y-\mathbb X_J\boldsymbol\beta_J}{N\lambda}\right\|_2^2\bigg) = \Proj_{(\bX_{J^c}B_1^{J^c})^\circ}\bigg(\frac{\vy-\bX_J\bbeta_J}{N\lambda}\bigg),
\end{align*}which implies that $\vlambda^\star[\vy - \bX_J\bbeta_J] = \Proj_{(\bX_{J^c}B_1^{J^c})^\circ}\big(\frac{\vy-\bX_J\bbeta_J}{N\lambda}\big)$ as announced in \textit{Item~1.}. Moreover, the projection onto a convex body is non-expansive, see, for instance, \cite[Chapter III, Proposition 3.1.3]{hiriart-urrutyConvexAnalysisMinimization1993}, and so $\vlambda^\star[\cdot]: \vmu\in (\bR^N,\ell_2)\to \Proj_{(\bX_{J^c}B_1^{J^c})^\circ}(\vmu/(N\lambda))\in (\bR^N,\ell_2)$ is $1/(N\lambda)$-Lipschitz, and hence $\|\D\vlambda^\star[\cdot](\vmu)\|_{\op}\leq\frac{1}{N\lambda}$ almost surely. 

{\textit{Item~2.}} follows from \eqref{eq:completing_squares_1} since $\boldsymbol\lambda\in(\mathbb X_{J^c}B_1^{J^c})^\circ$ iff $\|\bX_{J^c}^\top\vlambda\|_\infty\leq1$. To derive \textit{Item~3.}, we let $\hat\vnu$ be a solution to \eqref{eq:def_cA} for $\vmu = \boldsymbol y-\mathbb X_J\boldsymbol\beta_J$. Then, we have
\begin{align*}
\begin{aligned}
    &\boldsymbol y-\mathbb X_J\boldsymbol\beta_J-\mathbb X_{J^c}\widehat{\boldsymbol\nu}=N\lambda\boldsymbol\lambda^\star[\boldsymbol y-\mathbb X_J\boldsymbol\beta_J] = N\lambda\operatorname{Proj}_{(\mathbb X_{J^c}B_1^{J^c})^\circ}\left(\frac{\boldsymbol y-\mathbb X_J\boldsymbol\beta_J}{N\lambda}\right)=\operatorname*{argmin}_{\boldsymbol v\in N\lambda(\mathbb X_{J^c}B_1^{J^c})^\circ}\|\boldsymbol y-\mathbb X_J\boldsymbol\beta_J-\boldsymbol v\|_2^2\\
    &=\operatorname{Proj}_{N\lambda(\mathbb X_{J^c}B_1^{J^c})^\circ}(\boldsymbol y-\mathbb X_J\boldsymbol\beta_J).
\end{aligned}
\end{align*}

Finally, for {\textit{Item~4.}}, we remark that using \eqref{eq:dual_problem_effective_loss}, $\vbeta_J\to P_N\ell_{\vbeta_J}^{(\lambda,J)}$ is a maximal function of linear functions and that this maximum is unique at every $\vbeta_J$. It follows from Danskin's theorem that $\vbeta_J\to P_N\ell_{\vbeta_J}^{(\lambda,J)}$ is convex, differentiable and that its gradient is given by the derivative of the linear function $\bbeta_J \to \langle\boldsymbol\lambda,\boldsymbol y-\mathbb X_J\boldsymbol\beta_J\rangle-\frac{N\lambda}{2}\|\boldsymbol\lambda\|_2^2$ at the optimal value $\boldsymbol\lambda = \boldsymbol\lambda^\star[\boldsymbol y-\mathbb X_J\boldsymbol\beta_J]$ ie $\nabla P_N\ell_{\boldsymbol\beta_J}^{(\lambda,J)}=-\mathbb X_J^\top\boldsymbol\lambda^\star[\boldsymbol y-\mathbb X_J\boldsymbol\beta_J]$ (see \cite[Corollary VI.4.4.5]{hiriart-urrutyConvexAnalysisMinimization1993}).
\end{proof}

\subsection{Geometry of the effective goodness-of-fit term on \texorpdfstring{$V_{J_\star}$}{VJstar}}\label{sec:equiv_dist_polytope}

The goodness-of-fit term \eqref{eq:implicit_reg_LASSO} of the estimation part of the LASSO involves the distance from the residual $\vz=\vy-\bX_{J_\star}\bbeta$ to the dual polytope $N\lambda K^\circ=\big\{\vv\in\bR^N:\|\bX_{J_\star^c}^\top\vv\|_\infty\leq N\lambda\big\}$ where $K:=\bX_{J_\star^c}B_1^{J_\star^c}=\operatorname{conv}\big(\pm\bX_{(j)}:j\in J_\star^c\big)$
is the absolute convex hull of the $|J_\star^c|$ columns of $\bX_{J_\star^c}$. It is given by the function
\begin{equation}\label{eq:def_effective_fit}
    \cG_\lambda(\vz):=\frac12\|\vz\|_2^2-\frac12\operatorname{dist}\big(\vz,N\lambda K^\circ\big)^2,\qquad\vz\in\bR^N,
\end{equation}
so that \eqref{eq:implicit_reg_LASSO} equals $\frac{1}{N\lambda}\cG_\lambda(\vy-\bX_{J_\star}\bbeta)$. The aim of this section is to support the discussion following Proposition~\ref{prop:effective-regularization-lasso} regarding the shape of the goodness-of-fit term \eqref{eq:implicit_reg_LASSO}. To that end, we compare $\cG_\lambda$ with two Euclidean Huber functions with shape
\begin{equation}\label{eq:def_Huber}
    H_R(\vz):=\frac12\|\vz\|_2^2-\frac12\operatorname{dist}\big(\vz,RB_2^N\big)^2=\frac12\|\vz\|_2^2-\frac12\big(\|\vz\|_2-R\big)_+^2=\begin{cases}\frac12\|\vz\|_2^2&\mbox{if }\|\vz\|_2\leq R,\\ R\|\vz\|_2-\frac{R^2}{2}&\mbox{if }\|\vz\|_2>R,\end{cases}\qquad R>0,
\end{equation}
which are quadratic near the origin and linear in $\|\vz\|_2$ at infinity. Write $M:=|J_\star^c|$, $L:=\log(eM/N)$, and, for $\vw\in\bR^N$, let $\|\vw\|_K:=\min\{\|\vnu\|_1:\vnu\in V_{J_\star^c},\ \bX_{J_\star^c}\vnu=\vw\}$ denote the gauge of $K$, with the convention $\min\emptyset=+\infty$.

\begin{Theorem}\label{theo:equiv_dist_to_polytope}
Let $\lambda>0$.
\begin{enumerate}
    \item\label{item:dist_deterministic} \emph{(Deterministic description.)} Assume $\operatorname{rank}(\bX_{J_\star^c})=N$. Then for every $\vz\in\bR^N$,
    \begin{equation}\label{eq:dist_dual_formula}
        \operatorname{dist}\big(\vz,N\lambda K^\circ\big)=\max_{\vw\in B_2^N}\Big(\langle\vw,\vz\rangle-N\lambda\|\vw\|_K\Big),
    \end{equation}
    and, for every $\vz\neq\boldsymbol0$ and $\boldsymbol\theta:=\vz/\|\vz\|_2$,
    \begin{equation}\label{eq:dist_directional}
        \big(\|\vz\|_2-N\lambda\|\boldsymbol\theta\|_K\big)_+\leq\operatorname{dist}\big(\vz,N\lambda K^\circ\big)\leq\bigg(\|\vz\|_2-\frac{N\lambda}{\|\bX_{J_\star^c}^\top\boldsymbol\theta\|_\infty}\bigg)_+,\qquad\mbox{where}\qquad\|\boldsymbol\theta\|_K\,\|\bX_{J_\star^c}^\top\boldsymbol\theta\|_\infty\geq1.
    \end{equation}
    \item\label{item:dist_uniform} \emph{(Uniform two-sided bound.)} Grant the assumptions of Lemma~\ref{lem:geometric-event}, and let $c_{\mathrm G}$ be the constant of \eqref{eq:event-definitions}. Set
    \begin{equation*}
        R_-:=\frac{\sqrt N\lambda}{2}\qquad\mbox{and}\qquad R_+:=\frac{N\lambda}{c_{\mathrm G}\sqrt{L}}.
    \end{equation*}
    With probability at least $1-4\exp(-cN)$, simultaneously for every $\vz\in\bR^N$,
    \begin{equation}\label{eq:dist_uniform}
        \big(\|\vz\|_2-R_+\big)_+\leq\operatorname{dist}\big(\vz,N\lambda K^\circ\big)\leq\big(\|\vz\|_2-R_-\big)_+,\qquad H_{R_-}(\vz)\leq\cG_\lambda(\vz)\leq H_{R_+}(\vz).
    \end{equation}
    In particular, $\operatorname{dist}(\vz,N\lambda K^\circ)=0$ and $\cG_\lambda(\vz)=\frac12\|\vz\|_2^2$ whenever $\|\vz\|_2\leq R_-$; and $\frac12\|\vz\|_2\leq\operatorname{dist}(\vz,N\lambda K^\circ)\leq\|\vz\|_2$ and $\frac{R_-}{2}\|\vz\|_2\leq\cG_\lambda(\vz)\leq R_+\|\vz\|_2$ whenever $\|\vz\|_2\geq2R_+$.
    \item\label{item:dist_fixed} \emph{(Residuals independent of the free design.)} Assume $\Sigma=I_p$ and $M\geq N$. Let $\vz$ be a random vector of $\bR^N\setminus\{\boldsymbol0\}$ independent of $\bX_{J_\star^c}$, and let $u>0$. Set $R_u:=\frac{N\lambda}{\sqrt{2\log M+2u}}$. With probability (w.r.t. both $\bX_{J_\star^c}$ and $\vz$) at least $1-e^{-u}$,
    \begin{equation}\label{eq:dist_fixed}
        \operatorname{dist}\big(\vz,N\lambda K^\circ\big)\leq\big(\|\vz\|_2-R_u\big)_+\qquad\mbox{and}\qquad\cG_\lambda(\vz)\geq H_{R_u}(\vz).
    \end{equation}
    \item\label{item:dist_anisotropy} \emph{(The gap between $R_-$ and $R_+$ is intrinsic.)} Assume $\Sigma=I_p$ and $M\geq\max\{N,2\}$, fix $j_0\in J_\star^c$ and a deterministic $\boldsymbol\theta_0\in S^{N-1}$, and let $\boldsymbol\theta_{j_0}:=\bX_{(j_0)}/\|\bX_{(j_0)}\|_2$. With probability at least $1-e^{-N/5}-M^{-1}$, for every $t\geq0$,
    \begin{equation}\label{eq:dist_anisotropy}
        \operatorname{dist}\big(t\boldsymbol\theta_{j_0},N\lambda K^\circ\big)\geq t-2\sqrt N\lambda\qquad\mbox{and}\qquad\operatorname{dist}\big(t\boldsymbol\theta_0,N\lambda K^\circ\big)=0\ \mbox{ if }\ t\leq\frac{N\lambda}{2\sqrt{\log M}}.
    \end{equation}
    Consequently, if $N\geq64\log M$, then for every $t\in[4\sqrt N\lambda,\frac{N\lambda}{2\sqrt{\log M}}]$ the two 'residuals vectors' $t\boldsymbol\theta_{j_0}$ and $t\boldsymbol\theta_0$ have the same Euclidean norm, while $\operatorname{dist}(t\boldsymbol\theta_{j_0},N\lambda K^\circ)\geq\frac t2$ and $\operatorname{dist}(t\boldsymbol\theta_0,N\lambda K^\circ)=0$.
\end{enumerate}
\end{Theorem}

\begin{proof}
The proof proceeds by: (i) identifying the support function of $N\lambda K^\circ$ with $N\lambda\|\cdot\|_K$ through the bipolar theorem, which yields \eqref{eq:dist_dual_formula}, and bounding the distance along the ray through $\vz$, which yields \eqref{eq:dist_directional}; (ii) sandwiching $N\lambda K^\circ$ between two Euclidean balls on the event of Lemma~\ref{lem:geometric-event} and using the monotonicity of $C\mapsto\operatorname{dist}(\vz,C)$; and (iii) controlling $\|\bX_{J_\star^c}^\top\boldsymbol\theta\|_\infty$ and $\|\bX_{(j_0)}\|_2$ by Gaussian tail bounds in \eqref{eq:dist_directional}.

\textit{Proof of Item~\ref{item:dist_deterministic}.} Since $\operatorname{rank}(\bX_{J_\star^c})=N$, the affine set $\{\vnu\in V_{J_\star^c}:\bX_{J_\star^c}\vnu=\vw\}$ is non-empty and closed for every $\vw\in\bR^N$, and $\|\cdot\|_1$ is coercive; hence the minimum defining $\|\vw\|_K$ is attained and finite. Moreover, $K$ is a symmetric convex body, since it is compact, convex, symmetric, and contains a neighbourhood of $\boldsymbol0$ (it is the image of $B_1^{J_\star^c}$ under a surjective linear map). Consequently, $K^\circ$ is also a symmetric convex body, and by the bipolar theorem $(K^\circ)^\circ=K$, see \cite[Theorem 14.5]{rockafellarConvexAnalysis1970}. The support function of $K^\circ$ is the gauge of $(K^\circ)^\circ$; thus, for every $\vw\in\bR^N$,
\begin{equation}\label{eq:support_dual_polytope}
    h_{N\lambda K^\circ}(\vw):=\max_{\vv\in N\lambda K^\circ}\langle\vw,\vv\rangle=N\lambda\|\vw\|_K.
\end{equation}
Let $\vz\in\bR^N$. For every $\vw\in B_2^N$ and $\vv\in N\lambda K^\circ$, the Cauchy--Schwarz inequality and \eqref{eq:support_dual_polytope} give $\|\vz-\vv\|_2\geq\langle\vw,\vz-\vv\rangle\geq\langle\vw,\vz\rangle-N\lambda\|\vw\|_K$; minimizing over $\vv$ and maximizing over $\vw$ yields the inequality ``$\geq$'' in \eqref{eq:dist_dual_formula}. For the reverse inequality, if $\vz\in N\lambda K^\circ$, then the choice $\vw=\boldsymbol0$ gives the value $0=\operatorname{dist}(\vz,N\lambda K^\circ)$. Otherwise let $\vv^\star:=\Proj_{N\lambda K^\circ}(\vz)\neq\vz$ and $\vw^\star:=(\vz-\vv^\star)/\|\vz-\vv^\star\|_2\in S^{N-1}$. The variational characterization of the projection onto a closed convex set, see \cite[Chapter III, Theorem 3.1.1]{hiriart-urrutyConvexAnalysisMinimization1993}, gives $\langle\vz-\vv^\star,\vv-\vv^\star\rangle\leq0$ for every $\vv\in N\lambda K^\circ$, that is, $h_{N\lambda K^\circ}(\vw^\star)=\langle\vw^\star,\vv^\star\rangle$. Therefore, by \eqref{eq:support_dual_polytope},
\begin{equation*}
    \langle\vw^\star,\vz\rangle-N\lambda\|\vw^\star\|_K=\langle\vw^\star,\vz-\vv^\star\rangle=\|\vz-\vv^\star\|_2=\operatorname{dist}\big(\vz,N\lambda K^\circ\big),
\end{equation*}
which proves \eqref{eq:dist_dual_formula}; the maximum is attained since $\vw\mapsto\langle\vw,\vz\rangle-N\lambda\|\vw\|_K$ is continuous on the compact set $B_2^N$.

Let $\vz\neq\boldsymbol0$ and $\boldsymbol\theta=\vz/\|\vz\|_2$. The lower bound in \eqref{eq:dist_directional} follows from \eqref{eq:dist_dual_formula} with $\vw=\boldsymbol\theta$ and from $\operatorname{dist}(\vz,N\lambda K^\circ)\geq0$. For the upper bound, note that $\bX_{J_\star^c}^\top$ is injective, so $s_\theta:=\|\bX_{J_\star^c}^\top\boldsymbol\theta\|_\infty>0$, and $\vv_\theta:=\frac{N\lambda}{s_\theta}\boldsymbol\theta$ satisfies $\|\bX_{J_\star^c}^\top\vv_\theta\|_\infty=N\lambda$, i.e., $\vv_\theta\in N\lambda K^\circ$. If $\|\vz\|_2\geq\frac{N\lambda}{s_\theta}$, then $\operatorname{dist}(\vz,N\lambda K^\circ)\leq\|\vz-\vv_\theta\|_2=\|\vz\|_2-\frac{N\lambda}{s_\theta}$. If $\|\vz\|_2<\frac{N\lambda}{s_\theta}$, then $\vz$ lies on the segment $[\boldsymbol0,\vv_\theta]$, which is contained in the convex set $N\lambda K^\circ\ni\boldsymbol0$, so that $\operatorname{dist}(\vz,N\lambda K^\circ)=0$. Finally, let $\vnu\in V_{J_\star^c}$ satisfy $\bX_{J_\star^c}\vnu=\boldsymbol\theta$ and $\|\vnu\|_1=\|\boldsymbol\theta\|_K$. By the duality between $\ell_1$ and $\ell_\infty$,
\begin{equation*}
    1=\langle\boldsymbol\theta,\boldsymbol\theta\rangle=\langle\vnu,\bX_{J_\star^c}^\top\boldsymbol\theta\rangle\leq\|\vnu\|_1\|\bX_{J_\star^c}^\top\boldsymbol\theta\|_\infty=\|\boldsymbol\theta\|_K\,\|\bX_{J_\star^c}^\top\boldsymbol\theta\|_\infty.
\end{equation*}

\textit{Proof of Item~\ref{item:dist_uniform}.} By Lemma~\ref{lem:geometric-event}, $\bP(\Omega_{\mathrm{poly}}\cap\Omega_{\mathrm{col}})\geq\bP(\Omega)\geq1-4\exp(-cN)$; we work on $\Omega_{\mathrm{poly}}\cap\Omega_{\mathrm{col}}$. We first show
\begin{equation}\label{eq:sandwich_dual_polytope}
    R_-B_2^N\subseteq N\lambda K^\circ\subseteq R_+B_2^N.
\end{equation}
If $\|\vv\|_2\leq R_-$, then, by the Cauchy--Schwarz inequality and the definition of $\Omega_{\mathrm{col}}$, $\|\bX_{J_\star^c}^\top\vv\|_\infty=\max_{j\in J_\star^c}|\langle\bX_{(j)},\vv\rangle|\leq2\sqrt N\cdot\frac{\sqrt N\lambda}{2}=N\lambda$, which proves the first inclusion. Let $\vv\in N\lambda K^\circ\setminus\{\boldsymbol0\}$. On $\Omega_{\mathrm{poly}}$, $c_{\mathrm G}\sqrt L\,\vv/\|\vv\|_2\in K$, so there exists $\vnu\in B_1^{J_\star^c}$ with $\bX_{J_\star^c}\vnu=c_{\mathrm G}\sqrt L\,\vv/\|\vv\|_2$, and
\begin{equation*}
    c_{\mathrm G}\sqrt L\,\|\vv\|_2=\langle\bX_{J_\star^c}\vnu,\vv\rangle=\langle\vnu,\bX_{J_\star^c}^\top\vv\rangle\leq\|\vnu\|_1\|\bX_{J_\star^c}^\top\vv\|_\infty\leq N\lambda,
\end{equation*}
that is, $\|\vv\|_2\leq R_+$; this proves the second inclusion. For closed convex sets $C_1\subseteq C_2$ and $\vz\in\bR^N$, $\operatorname{dist}(\vz,C_2)\leq\operatorname{dist}(\vz,C_1)$. Together with \eqref{eq:sandwich_dual_polytope} and $\operatorname{dist}(\vz,RB_2^N)=(\|\vz\|_2-R)_+$, this yields the first chain of inequalities in \eqref{eq:dist_uniform}; the second one follows by the definitions \eqref{eq:def_effective_fit} and \eqref{eq:def_Huber}. If $\|\vz\|_2\leq R_-$, then $\vz\in N\lambda K^\circ$ by \eqref{eq:sandwich_dual_polytope}, so the distance vanishes. If $\|\vz\|_2\geq2R_+$, then $\|\vz\|_2-R_+\geq\frac12\|\vz\|_2$, while $\operatorname{dist}(\vz,N\lambda K^\circ)\leq\|\vz\|_2$ because $\boldsymbol0\in N\lambda K^\circ$. In this case, by \eqref{eq:def_Huber} and $R_-\leq R_+\leq\frac12\|\vz\|_2$,
\begin{equation*}
    \cG_\lambda(\vz)\leq H_{R_+}(\vz)\leq R_+\|\vz\|_2\qquad\mbox{and}\qquad\cG_\lambda(\vz)\geq H_{R_-}(\vz)=R_-\|\vz\|_2-\frac{R_-^2}{2}\geq\frac{R_-}{2}\|\vz\|_2.
\end{equation*}

\textit{Proof of Item~\ref{item:dist_fixed}.} Since $M\geq N$ and the entries of $\bX_{J_\star^c}$ are independent standard Gaussian variables, $\operatorname{rank}(\bX_{J_\star^c})=N$ almost surely, so \eqref{eq:dist_directional} applies almost surely. Let $\boldsymbol\theta=\vz/\|\vz\|_2$. Conditionally on $\vz$, the coordinates $\langle\bX_{(j)},\boldsymbol\theta\rangle$, $j\in J_\star^c$, are independent $\cN(0,1)$ variables, because the columns $\bX_{(j)}\sim\cN(\boldsymbol0,I_N)$ are independent of each other and of $\vz$, and $\|\boldsymbol\theta\|_2=1$. For $s\geq0$, the identity $\varphi(s+x)=\varphi(x)e^{-sx-s^2/2}\leq\varphi(x)e^{-s^2/2}$ for $x\geq0$, where $\varphi$ is the standard Gaussian density, gives $\Phi_G(s)=\int_0^\infty\varphi(s+x)\,\mathrm dx\leq\frac12e^{-s^2/2}$. Hence, by the union bound, with $s^2=2\log M+2u$,
\begin{equation*}
    \bP\Big\{\|\bX_{J_\star^c}^\top\boldsymbol\theta\|_\infty>s\,\Big|\,\vz\Big\}\leq2M\Phi_G(s)\leq Me^{-s^2/2}=e^{-u}.
\end{equation*}
On the complementary event, $\frac{N\lambda}{\|\bX_{J_\star^c}^\top\boldsymbol\theta\|_\infty}\geq R_u$, and the upper bound in \eqref{eq:dist_directional} gives $\operatorname{dist}(\vz,N\lambda K^\circ)\leq(\|\vz\|_2-R_u)_+$. The bound on $\cG_\lambda(\vz)$ follows from \eqref{eq:def_effective_fit} and \eqref{eq:def_Huber}. Integrating over $\vz$ proves \eqref{eq:dist_fixed}.

\textit{Proof of Item~\ref{item:dist_anisotropy}.} Let $\vg:=\bX_{(j_0)}\sim\cN(\boldsymbol0,I_N)$. By Markov's inequality applied to $e^{-\|\vg\|_2^2/2}$, and since $\bE e^{-\|\vg\|_2^2/2}=2^{-N/2}$,
\begin{equation*}
    \bP\Big\{\|\vg\|_2^2\leq\frac N4\Big\}\leq e^{N/8}\,\bE e^{-\|\vg\|_2^2/2}=\exp\Big(-N\Big(\frac{\log2}{2}-\frac18\Big)\Big)\leq e^{-N/5}.
\end{equation*}
On the event $\{\|\vg\|_2>\frac{\sqrt N}{2}\}$, $\boldsymbol\theta_{j_0}=\bX_{J_\star^c}\ve_{j_0}/\|\vg\|_2$ gives $\|\boldsymbol\theta_{j_0}\|_K\leq\|\ve_{j_0}\|_1/\|\vg\|_2<\frac{2}{\sqrt N}$, and the lower bound in \eqref{eq:dist_directional}, applied to $\vz=t\boldsymbol\theta_{j_0}$, gives $\operatorname{dist}(t\boldsymbol\theta_{j_0},N\lambda K^\circ)\geq t-N\lambda\|\boldsymbol\theta_{j_0}\|_K\geq t-2\sqrt N\lambda$ for every $t>0$; the case $t=0$ is trivial. Next, apply Item~\ref{item:dist_fixed} with $u=\log M$ to the deterministic direction $\boldsymbol\theta_0$: with probability at least $1-M^{-1}$, $\|\bX_{J_\star^c}^\top\boldsymbol\theta_0\|_\infty\leq2\sqrt{\log M}$, and the upper bound in \eqref{eq:dist_directional} gives $\operatorname{dist}(t\boldsymbol\theta_0,N\lambda K^\circ)=0$ for every $t\leq\frac{N\lambda}{2\sqrt{\log M}}$. A union bound proves \eqref{eq:dist_anisotropy}. If $N\geq64\log M$, then $4\sqrt N\lambda\leq\frac{N\lambda}{2\sqrt{\log M}}$, and for $t\geq4\sqrt N\lambda$, $t-2\sqrt N\lambda\geq\frac t2$.
\end{proof}

\begin{Remark}\label{rem:equiv_dist_polytope}
\begin{enumerate}
    \item Theorem~\ref{theo:equiv_dist_to_polytope} shows that the effective goodness-of-fit term $\frac{1}{N\lambda}\cG_\lambda(\vy-\bX_{J_\star}\bbeta)$ of \eqref{eq:implicit_reg_LASSO} is sandwiched between two Huber functions of the Euclidean norm of the residual: it is the squared loss $\frac{1}{2N\lambda}\|\vz\|_2^2$ for residuals inside $N\lambda K^\circ$, and it grows linearly in $\|\vz\|_2$, between $R_-\|\vz\|_2/2$ and $R_+\|\vz\|_2$, for residuals of norm at least $2R_+$. The interpolation is therefore between $\|\cdot\|_2^2$ and $\|\cdot\|_2$, in the spirit of the square-root LASSO. This is consistent with the rotation invariance of the law of $\bX_{J_\star^c}$, under which the law of $\cG_\lambda(\vz)$ depends on a deterministic $\vz$ only through $\|\vz\|_2$.
    \item The two radii in \eqref{eq:dist_uniform} come from the inradius $c_{\mathrm G}\sqrt L$ of the Gaussian polytope $K$, given by Gluskin's theorem \cite{gluskinExtremalPropertiesOrthogonal1989,litvakSmallestSingularValue2005,guedonGeometryPolytopesGenerated2022}, and from its circumradius $\max_{j}\|\bX_{(j)}\|_2\leq2\sqrt N$. Item~\ref{item:dist_anisotropy} shows that their ratio $\sqrt{N/L}$ cannot be removed from a bound that is uniform in $\vz\in\bR^N$ and depends on $\vz$ only through $\|\vz\|_2$: in the random directions of the vertices of $K$, $N\lambda K^\circ$ has radius at most $2\sqrt N\lambda$, whereas in a fixed direction its radius is at least $\frac{N\lambda}{2\sqrt{\log M}}$.
    \item For residuals independent of $\bX_{J_\star^c}$, such as $\vz=\bxi$, Items~\ref{item:dist_uniform} and~\ref{item:dist_fixed} give, with probability at least $1-4\exp(-cN)-e^{-u}$, $H_{R_u}(\vz)\leq\cG_\lambda(\vz)\leq H_{R_+}(\vz)$. If $M\geq N^{1+\delta}$ for some $\delta>0$ and $u\leq\log M$, then $\log M\leq\frac{1+\delta}{\delta}\log(M/N)\leq\frac{1+\delta}{\delta}L$, so that $R_+/R_u\leq\frac{2}{c_{\mathrm G}}\sqrt{\frac{1+\delta}{\delta}}$, and the two Huber functions have comparable radii. For $\|\vz\|_2\asymp\sigma_\xi\sqrt N$, this exhibits the phase transition anticipated in the less-regularized and well-regularized regimes of Proposition~\ref{prop:local-bernstein}: if $\lambda\leq\frac{c_{\mathrm G}}{4}\frac{\|\vz\|_2}{\sqrt N}\sqrt{\frac LN}$, then $\|\vz\|_2\geq2R_+$ and $\cG_\lambda(\vz)$ is linear in $\|\vz\|_2$ and so $\hat\bbeta_{J_\star}$ is like a square-root LASSO; if $\lambda\geq\frac{\|\vz\|_2}{\sqrt N}\sqrt{\frac{2\log M+2u}{N}}$, then $\|\vz\|_2\leq R_u$, so $\vz\in N\lambda K^\circ$ and $\cG_\lambda(\vz)=\frac12\|\vz\|_2^2$, i.e., the effective goodness-of-fit term coincides with the squared loss and $\hat\bbeta_{J_\star}$ is like the LASSO. This discussion holds only when the noise $\bxi$ is the dominant term in the residual $\vy - \bX_{J_\star}\hat\bbeta_{J_\star}$.
\end{enumerate}
\end{Remark}

\subsection{Restricted Isomorphic Property for the effective empirical risk} 

For some $\rho, r > 0$ to be determined later, we fix arbitrary $\bbeta_{J_\star}\in(B_1^{J_\star}(\bbeta_{J_\star}^\star;\rho) \cap S_2^{J_\star}(\bbeta_{J_\star}^\star;r))$. Define
\begin{align*}
    P\ell_\bullet^{(\lambda,J_\star)}:\bbeta_{J_\star}\in V_{J_\star} \mapsto \bE_{\bX_{J_\star},\bxi}\big[P_N\ell_{\bbeta_{J_\star}}^{(\lambda,J_\star)}\big|\bX_{J_\star^c}],\mbox{ and } P\cL_\bullet^{(\lambda,J_\star)}:\bbeta_{J_\star}\in V_{J_\star}\mapsto P\ell_{\bbeta_{J_\star}}^{(\lambda,J_\star)} - P\ell_{\bbeta_{J_\star}^\star}^{(\lambda,J_\star)}
\end{align*}be the (conditional) population risk and population excess risk functional. The following proposition establishes the local Bernstein's condition.

Let \(c_{\mathrm G},C_{\mathrm{ratio}}>0\) be absolute constants, let \(c_{\mathrm{CP}}>0\) be the constant in Lemma~\ref{lem:carl-pajor}, and let $C>0$ be a fixed absolute constant dominating the absolute constants denoted by $C$ in Propositions~\ref{prop:RIP} and~\ref{prop:RIP-small-lambda}. Assume $J_\star\neq\emptyset$, $N\leq c_m |J_\star^c|$, $|J_\star|\log\bigg(\frac{e|J_\star^c|}{N}\bigg)\leq c_\star N$, $\lambda\geq\frac{\sigma_\xi}{\sqrt{N}}$, for some absolute constants $0<c_m\leq C_{\mathrm{ratio}}^{-1}$ and $0<c_\star\leq\frac18$. Recall that $P_N\ell_{\boldsymbol\beta_{J_\star}}^{(\lambda,J_\star)}=\min_{\boldsymbol\nu\in\mathbb R^{J_\star^c}}\bigg\{\frac{\|\boldsymbol y-\mathbb X_{J_\star}\boldsymbol\beta_{J_\star}-\mathbb X_{J_\star^c}\boldsymbol\nu\|_2^2}{2N\lambda}+\|\boldsymbol\nu\|_1\bigg\}$, $P\ell_{\boldsymbol\beta_{J_\star}}^{(\lambda,J_\star)}=\mathbb E_{\mathbb X_{J_\star},\boldsymbol\xi}\big[P_N\ell_{\boldsymbol\beta_{J_\star}}^{(\lambda,J_\star)}\mid\mathbb X_{J_\star^c}\big]$, and $P\mathcal L_{\boldsymbol\beta_{J_\star}}^{(\lambda,J_\star)}=P\ell_{\boldsymbol\beta_{J_\star}}^{(\lambda,J_\star)}-P\ell_{\boldsymbol\beta_{J_\star}^\star}^{(\lambda,J_\star)}.$
Fix the absolute constants
\begin{align}
\begin{aligned}
&C_0=\max\bigg\{1,\frac1{c_{\mathrm G}}\bigg\},\, c_0=\min\bigg\{1,\frac{c_{\mathrm{CP}}}{4\sqrt{\pi e}}\bigg\},\, C_\beta=\frac{8C_0^3}{c_0^2}, \, c_{\mathrm{LB}}=\frac{c_0^2}{8C_0^2},\\
&C_\star=\max\bigg\{C_\beta,\frac{8\sqrt2}{c_{\mathrm{LB}}}\big(1+CC_0\big)\bigg\},\mbox{ and }c=\min\bigg\{c_{\mathrm G},C_{\mathrm{ratio}}^{-1},\frac{1}{(9C_\star)^2}\bigg\}.
\end{aligned}\label{eq:fixed-constants}
\end{align}

In the following Proposition~\ref{prop:local-bernstein}, we prove that the effective population excess risk $\bbeta \in V_{J_\star} \to P\cL_{\bbeta}^{(\lambda,J_\star)}$ satisfies the local Bernstein's property (\cite{chinotRobustStatisticalLearning2020,chinotRobustHighDimensional2021,lecueLearningLinearLoss2024})  --- a property that lays the foundation for our subsequent establishment of the Restricted Isomorphic Property.  Roughly speaking, Proposition~\ref{prop:local-bernstein} states that, depending on the relative size of $\lambda$ and $\sigma_\xi \sqrt{\frac{\log(\frac{e|J_\star^c|}{N})}{N}}$, the lower bound on the curvature of the effective population excess risk $r$-away from $\bbeta_{J_\star}^\star$ behaves like $r^2/\lambda$ (when $\lambda \gtrsim \sigma_\xi \sqrt{\frac{1}{N}\log(\frac{e|J_\star^c|}{N})}$), or like $\frac{r^2}{\sigma_\xi^2+r^2} \sqrt{\frac{N}{\log(\frac{e|J_\star^c|}{N})}}$ (when $\lambda \lesssim \sigma_\xi \sqrt{\frac{1}{N}\log(\frac{e|J_\star^c|}{N})}$). Accordingly, we call the former the well-regularized regime and the latter the less-regularized regime.

\begin{Proposition}[local Bernstein's property on $V_{J_\star}$]\label{prop:local-bernstein} Grant Assumption~\ref{assumption:LASSO_upper}.
There exist absolute constants $c_{\mathrm{LB}},C_\beta>0$ such that  with probability at least \(1-4\exp(-cN)\) with respect to \(\mathbb X_{J_\star^c}\), for any $r\geq C_\beta\|\bbeta_{J_\star^c}^\star\|_2$,
\begin{align*}
&\forall \bbeta_{J_\star}\in S_2^{J_\star}(\bbeta_{J_\star}^\star;r),\qquad P\mathcal L_{\boldsymbol\beta_{J_\star}}^{(\lambda,J_\star)}\geq c_{\mathrm{LB}}r^2\min\bigg\{\frac1\lambda,\frac{1}{\sqrt{\sigma_\xi^2+r^2}}\sqrt{\frac{N}{\log\big(\frac{e|J_\star^c|}{N}\big)}}\bigg\}.
\end{align*}
\end{Proposition}

Unlike the classical approach in statistical learning theory for verifying the local Bernstein's condition, here the verification is more difficult because we are dealing with a random loss function; indeed, $P\mathcal L_{\bullet}^{(\lambda,J_\star)}$ is an excess risk over $V_{J_\star}$ that depends on $\bX_{J_\star^c}$. This is where the Gaussian assumption on the noise and the design plays a crucial role (we refer the reader to \cite[Section 1.3.3]{shangFeatureSpaceDecomposition2026}). More specifically, the proof below views the effective population excess risk $P\cL_{\bbeta_{J_\star}}^{(\lambda,J_\star)}$ as the gap between an offset Gaussian mean width (up to a multiplicative factor), and a certain restricted projection of the 'not estimated part' of the signal $\bbeta_{J_\star^c}$. Here, an offset Gaussian mean width refers to a quantity of the form $\mathbb{E} \sup(\langle\vg,\vv\rangle-\varphi(\vv) : \vv \in K)$, where $\varphi$ is a fixed function, $K$ is a measurable set, and $\vg$ is a standard Gaussian random vector. This proof strategy is not confined to LASSO; only when examining the lower bound of this offset Gaussian mean width do we introduce standard AGA results on the lower bound of the Gaussian mean width of a Gaussian random (cross)-polytope and its $\ell_2^N$-diameter upper bound; see \cite{litvakSmallestSingularValue2005}. To the best of our knowledge, this approach of establishing the local Bernstein's property is novel.

\subsubsection{Proof of Proposition~\texorpdfstring{\ref{prop:local-bernstein}}{local Bernstein's property}}

\begin{Lemma}\label{lem:geometric-event}
Assume $\Sigma=I_p$. Define $\Omega=\Omega_{\mathrm{poly}}\cap\Omega_{\mathrm{col}}\cap\Omega_{\mathrm{signal}}$ where
\begin{align}
\begin{aligned}
&\Omega_{\mathrm{poly}}=\bigg\{\mathbb X_{J_\star^c}B_1^{J_\star^c}\supseteq c_{\mathrm G}\sqrt{\log\bigg(\frac{e|J_\star^c|}{N}\bigg)}B_2^N\bigg\},\qquad \Omega_{\mathrm{col}}=\bigg\{\max_{j\in J_\star^c}\|\mathbb X_{(j)}\|_2\leq2\sqrt N\bigg\},\\
&\mbox{ and  }\Omega_{\mathrm{signal}}=\bigg\{\|\mathbb X_{J_\star^c}\boldsymbol\beta_{J_\star^c}^\star\|_2\leq2\sqrt N\|\boldsymbol\beta_{J_\star^c}^\star\|_2\bigg\}.
\end{aligned}\label{eq:event-definitions}
\end{align}
Then \(\mathbb P(\Omega)\geq1-4\exp(-cN)\). On \(\Omega\),
\begin{align}
\begin{aligned}
&\diam\big((\bX_{J_\star^c}B_1^{J_\star^c})^\circ,\ell_2^N\big)=2\sup\big(\|\vlambda\|_2:\, \|\bX_{J_\star^c}^\top\vlambda\|_\infty\leq 1\big) \leq\frac{2C_0}{\sqrt{\log\big(\frac{e|J_\star^c|}{N}\big)}}
\end{aligned}\label{eq:dual-radius}
\end{align}
and, for \(\boldsymbol G\sim\mathcal N(0,I_N)\) independent of \(\mathbb X_{J_\star^c}\),
\begin{align}
\begin{aligned}
&c_0\sqrt{\frac{N}{\log\big(e|J_\star^c|/N\big)}}\leq\mathbb E_{\boldsymbol G}\sup_{\|\mathbb X_{J_\star^c}^\top\boldsymbol u\|_\infty\leq1}\langle\boldsymbol G,\boldsymbol u\rangle\leq C_0\sqrt{\frac{N}{\log\big(e|J_\star^c|/N\big)}}.
\end{aligned}\label{eq:dual-width}
\end{align}
\end{Lemma}

\begin{proof}[Proof of Lemma~\ref{lem:geometric-event}.]
Using the Borell-TIS inequality, we have for all $t>0$, with probability at least $1-(1+|J_\star^c|)\exp(-t)$, $\|\mathbb X_{J_\star^c}\boldsymbol\beta_{J_\star^c}^\star\|_2\leq \bE \|\mathbb X_{J_\star^c}\boldsymbol\beta_{J_\star^c}^\star\|_2 + \sqrt{t}$ and for all $j\in J_\star^c, \|\mathbb X_{(j)}\|_2\leq \bE \|\mathbb X_{(j)}\|_2  + \sqrt{t}$. Hence, when $N\geq 2\log(1+|J_\star^c|)$, (which is implied by the assumption $|J_\star|\log(e|J_\star^c|/N)\leq c_\star N$ where $c_\star\leq\frac{1}{8}$),  we get that $\bP[\Omega_{\mathrm{col}}\cap\Omega_{\mathrm{signal}}]\geq 1-\exp(-N/2)$. Next, one can check that $\mathbb X_{J_\star^c}B_1^{J_\star^c}\supseteq \square B_2^N$ for $\square=c_{\mathrm G}\sqrt{\log\big(e|J_\star^c|/N\big)}$ with probability at least $1-2\exp(-cN)$ from Gluskin's theorem, see, for instance,  \cite[Theorem 1]{guedonGeometryPolytopesGenerated2022}. This shows that $\bP(\Omega)\geq 1-4\exp(-cN)$. In the following, we work on $\Omega$. 

Notice that
\begin{align*}
&\operatorname{Vol}_N\big(\big\{\boldsymbol u:\|\mathbb X_{J_\star^c}^\top\boldsymbol u\|_\infty\leq1\big\}\big)^{\frac{1}{N}}=\frac{1}{2\sqrt N}\operatorname{Vol}_N\big(\big\{\boldsymbol u:\|\mathbb X_{J_\star^c}^\top\boldsymbol u\|_\infty\leq2\sqrt N\big\}\big)^{\frac{1}{N}}.
\end{align*}
By Carl-Gluksin-Pajor's theorem, see \cite[Theorem 4.6.8]{artstein-avidanAsymptoticGeometricAnalysis2021}, see also Lemma~\ref{lem:carl-pajor},
\begin{align}
\begin{aligned}
&\operatorname{Vol}_N\big(\big\{\boldsymbol u\in\mathbb R^N:\|\mathbb X_{J_\star^c}^\top\boldsymbol u\|_\infty\leq1\big\}\big)^{\frac{1}{N}}\geq\frac{c_{\mathrm{CP}}}{2\sqrt N\sqrt{\log\big(e|J_\star^c|/N\big)}}.
\end{aligned}\label{eq:dual-volume}
\end{align}
The set in \eqref{eq:dual-volume} is the intersection of the \(2|J_\star^c|\) closed half-spaces  $\{\bu:|\langle\mathbb X_{(j)},\boldsymbol u\rangle|\leq1\}, j\in J_\star^c$, and
\begin{align*}
&\big\{\boldsymbol u:\|\mathbb X_{J_\star^c}^\top\boldsymbol u\|_\infty\leq1\big\}\subseteq\frac{C_0}{\sqrt{\log(e|J_\star^c|/N)}}B_2^N.
\end{align*}

Moreover, since $\frac{\boldsymbol G}{\|\boldsymbol G\|_2}$ is uniform over $S^{N-1}$, $ \frac{\boldsymbol G}{\|\boldsymbol G\|_2}\perp\!\!\!\perp\|\boldsymbol G\|_2,$ $ \mathbb E\|\boldsymbol G\|_2\geq\sqrt{\frac N2}$, and $\operatorname{Vol}_N(B_2^N)=\frac{\pi^{\frac{N}{2}}}{\Gamma(\frac{N}{2}+1)}\leq\big(\frac{2\pi e}{N}\big)^{\frac{N}{2}}.$
Thus by Urysohn's inequality \cite[Corollary~1.4]{pisierVolumeConvexBodies1989}, see also Lemma~\ref{lem:urysohn} and \eqref{eq:dual-volume} yield
\begin{align*}
&\mathbb E_{\boldsymbol G}\sup_{\|\mathbb X_{J_\star^c}^\top\boldsymbol u\|_\infty\leq1}\langle\boldsymbol G,\boldsymbol u\rangle=\mathbb E\|\boldsymbol G\|_2\int_{S^{N-1}}\sup_{\|\mathbb X_{J_\star^c}^\top\boldsymbol u\|_\infty\leq1}\langle\boldsymbol\theta,\boldsymbol u\rangle\,\mathrm d\sigma(\boldsymbol\theta)\geq\frac{c_{\mathrm{CP}}}{4\sqrt{\pi e}}\sqrt{\frac{N}{\log\big(e|J_\star^c|/N\big)}}.
\end{align*}
Finally, \eqref{eq:dual-radius} and \(\mathbb E\|\boldsymbol G\|_2\leq\sqrt N\) imply
\begin{align*}
&\mathbb E_{\boldsymbol G}\sup_{\|\mathbb X_{J_\star^c}^\top\boldsymbol u\|_\infty\leq1}\langle\boldsymbol G,\boldsymbol u\rangle\leq\mathbb E\|\boldsymbol G\|_2\sup_{\|\mathbb X_{J_\star^c}^\top\boldsymbol u\|_\infty\leq1}\|\boldsymbol u\|_2\leq C_0\sqrt{\frac{N}{\log\big(e|J_\star^c|/N\big)}}.
\end{align*}
This proves \eqref{eq:dual-width}.
\end{proof}

Now we prove Proposition~\ref{prop:local-bernstein}.

\begin{proof}[Proof of Proposition~\ref{prop:local-bernstein}.]
We adopt a proof strategy similar to that of \cite[Lemma 19]{shangFeatureSpaceDecomposition2026} --- using the fact that the joint  distribution of $(\bX_{J_\star}, \bxi)$ is Gaussian, to reduce the proof of the local Bernstein's property to a one-dimensional setting.

Fix $r>0$ and let $\bbeta_{J_\star}\in S_2^{J_\star}(\bbeta_{J_\star}^\star;r)$. Then, $\bxi-\bX_{J_\star}(\bbeta_{J_\star}-\bbeta_{J_\star}^\star)$ follows a centered Gaussian distribution with covariance matrix $(\sigma_\xi^2 + r^2)I_N$. Let $\vG\sim\cN(\vzero,I_N)$, then $P\ell_{\bbeta_{J_\star}}^{(\lambda,J_\star)}$ has the following offset Gaussian mean width representation
\begin{align}\label{eq:conditonal_population_risk}
    P\ell_{\bbeta_{J_\star}}^{(\lambda,J_\star)} = \bE_\vG\max\bigg( \bigg\langle \vu,\bX_{J_\star^c}\bbeta_{J_\star^c}^\star + \sqrt{\sigma_\xi^2+r^2}\vG \bigg\rangle - \frac{N\lambda}{2}\|\vu\|_2^2 :\, \|\bX_{J_\star^c}^\top\vu\|_\infty\leq 1\bigg),
\end{align}and $P\ell_{\bbeta_{J_\star}^\star}^{(\lambda,J_\star)}$ is precisely equal to \eqref{eq:conditonal_population_risk} but with $r^2$ replaced by $0$. Notice for every feasible $\vu$, there holds 
\begin{align*}
    \langle\vu,\bX_{J_\star^c}\bbeta_{J_\star^c}^\star + \sigma_\xi\vG\rangle - \frac{N\lambda}{2}\|\vu\|_2^2 &= \frac{\sigma_\xi}{\sqrt{\sigma_\xi^2 + r^2}}\big( \langle\vu,\bX_{J_\star^c}\bbeta_{J_\star^c}^\star + \sqrt{\sigma_\xi^2+r^2}\vG\rangle - \frac{N\lambda}{2}\|\vu\|_2^2 \big) \\
    &+ \big(1-\frac{\sigma_\xi}{\sqrt{\sigma_\xi^2 + r^2}}\big)\big(\langle\vu,\bX_{J_\star^c}\bbeta_{J_\star^c}^\star\rangle - \frac{N\lambda}{2}\|\vu\|_2^2\big).
\end{align*}
By Jensen's inequality,
\begin{align*}
    &P\ell_{\bbeta_{J_\star}^\star}^{(\lambda,J_\star)} = \bE_\vG\max\bigg( \bigg\langle \vu,\bX_{J_\star^c}\bbeta_{J_\star^c}^\star + \sigma_\xi \vG \bigg\rangle - \frac{N\lambda}{2}\|\vu\|_2^2 :\, \|\bX_{J_\star^c}^\top\vu\|_\infty\leq 1\bigg)\\
    &\leq \frac{\sigma_\xi}{\sqrt{\sigma_\xi^2 + r^2}}\bE_\vG\max\bigg( \langle\vu,\bX_{J_\star^c}\bbeta_{J_\star^c}^\star + \sqrt{\sigma_\xi^2+r^2}\vG\rangle - \frac{N\lambda}{2}\|\vu\|_2^2 :\|\bX_{J_\star^c}^\top\vu\|_\infty\leq 1\bigg)\\
    &+\big(1-\frac{\sigma_\xi}{\sqrt{\sigma_\xi^2 + r^2}}\big)\max\bigg(\langle\vu,\bX_{J_\star^c}\bbeta_{J_\star^c}^\star\rangle - \frac{N\lambda}{2}\|\vu\|_2^2:\, \|\bX_{J_\star^c}^\top\vu\|_\infty\leq 1\bigg).
\end{align*}Notice that by \eqref{eq:conditonal_population_risk}, $\bE_\vG\max\big( \langle\vu,\bX_{J_\star^c}\bbeta_{J_\star^c}^\star + \sqrt{\sigma_\xi^2+r^2}\vG\rangle - \frac{N\lambda}{2}\|\vu\|_2^2 :\|\bX_{J_\star^c}^\top\vu\|_\infty\leq 1\big)=P\ell_{\bbeta_{J_\star}}^{(\lambda,J_\star)}$. Moreover, by definition, $P\cL_{\bbeta_{J_\star}}^{(\lambda,J_\star)} = P\ell_{\bbeta_{J_\star}}^{(\lambda,J_\star)}-P\ell_{\bbeta_{J_\star}^\star}^{(\lambda,J_\star)}$, we obtain
\begin{align}\label{eq:pre_lower_bound_local_Bernstein}
    &P\cL_{\bbeta_{J_\star}}^{(\lambda,J_\star)} \geq \bigg(1-\frac{\sigma_\xi}{\sqrt{\sigma_\xi^2 + r^2}}\bigg)P\ell_{\bbeta_{J_\star}^\star}^{(\lambda,J_\star)} - \bigg(1-\frac{\sigma_\xi}{\sqrt{\sigma_\xi^2 + r^2}}\bigg)\max\bigg(\langle\vu,\bX_{J_\star^c}\bbeta_{J_\star^c}^\star\rangle - \frac{N\lambda}{2}\|\vu\|_2^2:\, \|\bX_{J_\star^c}^\top\vu\|_\infty\leq 1\bigg).
\end{align}We bound the right-hand-side of \eqref{eq:pre_lower_bound_local_Bernstein} from below. For any $\vlambda\in\bR^N$, let \(\boldsymbol v_{\vlambda}\) be the unique element of \(\operatorname*{arg\,max}(\langle\vlambda,\boldsymbol z\rangle:\, \|\mathbb X_{J_\star^c}^\top\boldsymbol z\|_\infty\leq1)\) (which is well-defined almost surely). Then on $\Omega$, there holds
\begin{align}
\begin{aligned}
&\|\mathbb X_{J_\star^c}^\top\boldsymbol v_{\vlambda}\|_\infty\leq1,\qquad \boldsymbol v_{-\vlambda}=-\boldsymbol v_{\vlambda},\qquad \langle\vlambda,\boldsymbol v_{\vlambda}\rangle=\sup_{\|\mathbb X_{J_\star^c}^\top\boldsymbol u\|_\infty\leq1}\langle\vlambda,\boldsymbol u\rangle,\mbox{ and }\qquad \|\boldsymbol v_{\vlambda}\|_2\leq\frac{C_0}{\sqrt{\log\big(\frac{e|J_\star^c|}{N}\big)}}.
\end{aligned}\label{eq:v-properties}
\end{align}Let $\vlambda=\vG$, then \(\mathbb E_{\boldsymbol G}\|\boldsymbol v_{\boldsymbol G}\|_2\leq \frac{C_0}{\sqrt{\log(\frac{e|J_\star^c|}{N})}}\), \(\boldsymbol G\overset{d}{=}-\boldsymbol G\), and \eqref{eq:v-properties} imply \(\mathbb E_{\boldsymbol G}\boldsymbol v_{\boldsymbol G}=\mathbb E_{\boldsymbol G}\boldsymbol v_{-\boldsymbol G}=-\mathbb E_{\boldsymbol G}\boldsymbol v_{\boldsymbol G}=0\).
\begin{enumerate}
    \item We bound the first term of the right-hand-side of \eqref{eq:pre_lower_bound_local_Bernstein}, that is, $(1-\frac{\sigma_\xi}{\sqrt{\sigma_\xi^2+r^2}})P\ell_{\bbeta_{J_\star}^\star}^{(\lambda,J_\star)}$, from below. In this lower bound, we use the lower bound on the Gaussian mean width and the upper bound on the $\ell_2^N$-diameter of a Gaussian random polytope. For every \(\theta\in[0,1]\), \eqref{eq:v-properties} gives \(\|\mathbb X_{J_\star^c}^\top(\theta\boldsymbol v_{\boldsymbol G})\|_\infty\leq1\) and
    \begin{align*}
    \begin{aligned}
    &\max_{\|\mathbb X_{J_\star^c}^\top\boldsymbol u\|_\infty\leq1}\bigg\{\big\langle\boldsymbol u,\mathbb X_{J_\star^c}\boldsymbol\beta_{J_\star^c}^\star+\sqrt{\sigma_\xi^2+r^2}\boldsymbol G\big\rangle-\frac{N\lambda}{2}\|\boldsymbol u\|_2^2\bigg\} \geq\theta\langle\boldsymbol v_{\boldsymbol G},\mathbb X_{J_\star^c}\boldsymbol\beta_{J_\star^c}^\star\rangle+\theta\sqrt{\sigma_\xi^2+r^2}\langle\boldsymbol G,\boldsymbol v_{\boldsymbol G}\rangle-\frac{N\lambda\theta^2}{2}\|\boldsymbol v_{\boldsymbol G}\|_2^2.
    \end{aligned}
    \end{align*}
    Moreover, \(\mathbb E_{\boldsymbol G}\langle\boldsymbol v_{\boldsymbol G},\mathbb X_{J_\star^c}\boldsymbol\beta_{J_\star^c}^\star\rangle=\langle\mathbb E_{\boldsymbol G}\boldsymbol v_{\boldsymbol G},\mathbb X_{J_\star^c}\boldsymbol\beta_{J_\star^c}^\star\rangle=0\). Taking expectations and using \eqref{eq:dual-width} and \eqref{eq:v-properties} yield, for every \(\theta\in[0,1]\),
    \begin{align}
    \begin{aligned}
    &\mathbb E_{\boldsymbol G}\max_{\|\mathbb X_{J_\star^c}^\top\boldsymbol u\|_\infty\leq1}\bigg\{\big\langle\boldsymbol u,\mathbb X_{J_\star^c}\boldsymbol\beta_{J_\star^c}^\star+\sqrt{\sigma_\xi^2+r^2}\boldsymbol G\big\rangle-\frac{N\lambda}{2}\|\boldsymbol u\|_2^2\bigg\} \geq c_0\theta\sqrt{\sigma_\xi^2+r^2}\sqrt{\frac{N}{\log\big(\frac{e|J_\star^c|}{N}\big)}}-\frac{C_0^2N\lambda\theta^2}{2\log\big(\frac{e|J_\star^c|}{N}\big)}.
    \end{aligned}\label{eq:scalar-optimization}
    \end{align}
    \begin{enumerate}
    \item\label{item:scalar-interior} If \(\frac{c_0}{C_0^2}\frac{\sqrt{\sigma_\xi^2+r^2}}{\lambda}\sqrt{\frac{\log(\frac{e|J_\star^c|}{N})}{N}}\leq1\), substitute \(\theta=\frac{c_0}{C_0^2}\frac{\sqrt{\sigma_\xi^2+r^2}}{\lambda}\sqrt{\frac{\log(\frac{e|J_\star^c|}{N})}{N}}\) into \eqref{eq:scalar-optimization}; the right-hand side equals \(\frac{c_0^2}{2C_0^2}\frac{\sigma_\xi^2+r^2}{\lambda}\).
    \item\label{item:scalar-boundary} If \(\frac{c_0}{C_0^2}\frac{\sqrt{\sigma_\xi^2+r^2}}{\lambda}\sqrt{\frac{\log(\frac{e|J_\star^c|}{N})}{N}}>1\), substitute \(\theta=1\) into \eqref{eq:scalar-optimization}. The case assumption gives \(\frac{C_0^2N\lambda}{2\log(\frac{e|J_\star^c|}{N})}<\frac{c_0\sqrt{\sigma_\xi^2+r^2}\sqrt{\frac{N}{\log(\frac{e|J_\star^c|}{N})}}}{2}\). The right-hand side of \eqref{eq:scalar-optimization} is at least \(\frac{c_0\sqrt{\sigma_\xi^2+r^2}\sqrt{\frac{N}{\log(\frac{e|J_\star^c|}{N})}}}{2}\).
    \end{enumerate}
    Items~\ref{item:scalar-interior} and~\ref{item:scalar-boundary}, together with \(c_0\leq1\leq C_0\), provide a lower bound for the first term of \eqref{eq:pre_lower_bound_local_Bernstein}
    \begin{align}
    \begin{aligned}
    &\mathbb E_{\boldsymbol G}\max_{\|\mathbb X_{J_\star^c}^\top\boldsymbol u\|_\infty\leq1}\bigg\{\big\langle\boldsymbol u,\mathbb X_{J_\star^c}\boldsymbol\beta_{J_\star^c}^\star+\sqrt{\sigma_\xi^2+r^2}\boldsymbol G\big\rangle-\frac{N\lambda}{2}\|\boldsymbol u\|_2^2\bigg\} \geq\frac{c_0^2}{2C_0^2}\min\bigg\{\frac{\sigma_\xi^2+r^2}{\lambda},\sqrt{\sigma_\xi^2+r^2}\sqrt{\frac{N}{\log\big(\frac{e|J_\star^c|}{N}\big)}}\bigg\}.
    \end{aligned}\label{eq:gaussian-lower-bound}
    \end{align}

    \item We bound the second term of the right-hand-side of \eqref{eq:pre_lower_bound_local_Bernstein}, that is, $\big(1-\frac{\sigma_\xi}{\sqrt{\sigma_\xi^2 + r^2}}\big)\max\big(\langle\vu,\bX_{J_\star^c}\bbeta_{J_\star^c}^\star\rangle - \frac{N\lambda}{2}\|\vu\|_2^2:\, \|\bX_{J_\star^c}^\top\vu\|_\infty\leq 1\big)$, from below. In this upper bound, we use the upper bound on the $\ell_2^N$-diameter of a Gaussian random polytope. For every \(\boldsymbol u\in\mathbb R^N\),
    \begin{align}
    \begin{aligned}
    &\max_{\boldsymbol u\in\mathbb R^N}\bigg\{\langle\boldsymbol u,\mathbb X_{J_\star^c}\boldsymbol\beta_{J_\star^c}^\star\rangle-\frac{N\lambda}{2}\|\boldsymbol u\|_2^2\bigg\}=\frac{\|\mathbb X_{J_\star^c}\boldsymbol\beta_{J_\star^c}^\star\|_2^2}{2N\lambda},
    \end{aligned}\tag{\ref{eq:variational_representation_ell_2_square}}
    \end{align}
    and the maximizer in \eqref{eq:variational_representation_ell_2_square} is \(\frac{\mathbb X_{J_\star^c}\boldsymbol\beta_{J_\star^c}^\star}{N\lambda}\). For every \(\boldsymbol u\) satisfying \(\|\mathbb X_{J_\star^c}^\top\boldsymbol u\|_\infty\leq1\), $\langle\boldsymbol u,\mathbb X_{J_\star^c}\boldsymbol\beta_{J_\star^c}^\star\rangle-\frac{N\lambda}{2}\|\boldsymbol u\|_2^2\leq\langle\boldsymbol u,\mathbb X_{J_\star^c}\boldsymbol\beta_{J_\star^c}^\star\rangle\leq\|\boldsymbol u\|_2\|\mathbb X_{J_\star^c}\boldsymbol\beta_{J_\star^c}^\star\|_2.$ Combining \eqref{eq:dual-radius}, \eqref{eq:variational_representation_ell_2_square}, and \(\Omega_{\mathrm{signal}}\) gives
    \begin{align}
    \begin{aligned}
    &\max_{\|\mathbb X_{J_\star^c}^\top\boldsymbol u\|_\infty\leq1}\bigg\{\langle\boldsymbol u,\mathbb X_{J_\star^c}\boldsymbol\beta_{J_\star^c}^\star\rangle-\frac{N\lambda}{2}\|\boldsymbol u\|_2^2\bigg\}\\
    &\leq\min\bigg\{\frac{\|\mathbb X_{J_\star^c}\boldsymbol\beta_{J_\star^c}^\star\|_2^2}{2N\lambda},\sup_{\|\mathbb X_{J_\star^c}^\top\boldsymbol u\|_\infty\leq1}\|\boldsymbol u\|_2\|\mathbb X_{J_\star^c}\boldsymbol\beta_{J_\star^c}^\star\|_2\bigg\} \leq2C_0\min\bigg\{\frac{\|\boldsymbol\beta_{J_\star^c}^\star\|_2^2}{\lambda},\|\boldsymbol\beta_{J_\star^c}^\star\|_2\sqrt{\frac{N}{\log\big(\frac{e|J_\star^c|}{N}\big)}}\bigg\}.
    \end{aligned}\label{eq:signal-upper-bound}
    \end{align}
    The assumptions on \(r\) and \eqref{eq:fixed-constants} imply \(C_\beta\geq1\), \(\|\boldsymbol\beta_{J_\star^c}^\star\|_2\leq \frac{r}{C_\beta}\leq\frac{\sqrt{\sigma_\xi^2+r^2}}{C_\beta}\), and \(\frac{2C_0}{C_\beta}=\frac{c_0^2}{4C_0^2}\). Hence
    \begin{align}
    \begin{aligned}
    &2C_0\frac{\|\boldsymbol\beta_{J_\star^c}^\star\|_2^2}{\lambda}\leq\frac{c_0^2}{4C_0^2}\frac{\sigma_\xi^2+r^2}{\lambda},\mbox{ and }2C_0\|\boldsymbol\beta_{J_\star^c}^\star\|_2\sqrt{\frac{N}{\log\big(\frac{e|J_\star^c|}{N}\big)}}\leq\frac{c_0^2}{4C_0^2}\sqrt{\sigma_\xi^2+r^2}\sqrt{\frac{N}{\log\big(\frac{e|J_\star^c|}{N}\big)}}.
    \end{aligned}\label{eq:signal-absorption}
    \end{align}
\end{enumerate}
Finally,
\begin{align}
\begin{aligned}
&1-\frac{\sigma_\xi}{\sqrt{\sigma_\xi^2+r^2}}=\frac{r^2}{\sqrt{\sigma_\xi^2+r^2}\big(\sqrt{\sigma_\xi^2+r^2}+\sigma_\xi\big)}\geq\frac{r^2}{2(\sigma_\xi^2+r^2)},\mbox{ and }\\
&\bigg(1-\frac{\sigma_\xi}{\sqrt{\sigma_\xi^2+r^2}}\bigg)\min\bigg\{\frac{\sigma_\xi^2+r^2}{\lambda},\sqrt{\sigma_\xi^2+r^2}\sqrt{\frac{N}{\log\big(\frac{e|J_\star^c|}{N}\big)}}\bigg\} \geq\frac{r^2}{2}\min\bigg\{\frac1\lambda,\frac{1}{\sqrt{\sigma_\xi^2+r^2}}\sqrt{\frac{N}{\log\big(\frac{e|J_\star^c|}{N}\big)}}\bigg\}.
\end{aligned}\label{eq:curvature-identity}
\end{align}
Since \(r\) and \(\boldsymbol\beta_{J_\star}\) were arbitrary and \(\Omega\) is independent of both, together with \(\mathbb P(\Omega)\geq1-4\exp(-cN)\), this proves Proposition~\ref{prop:local-bernstein}.
\end{proof}

\subsection{Restricted Isomorphic Property}

Let $T=B_1^{J_\star}(\boldsymbol\beta_{J_\star}^\star;\rho)\cap B_2^{J_\star}(\boldsymbol\beta_{J_\star}^\star;r)$ for some $\rho,r>0$. Recall from Proposition~\ref{prop:effective-regularization-lasso} that for all $\boldsymbol{s}\in\bR^N$,
\begin{align}
&\boldsymbol\lambda^\star[\boldsymbol s]=\operatorname{Proj}_{(\mathbb X_{J_\star^c}B_1^{J_\star^c})^\circ}\!\left(\frac{\boldsymbol s}{N\lambda}\right).\label{eq:projection-representation-rip}
\end{align}
The main result of this subsection are Proposition~\ref{prop:RIP} and Proposition~\ref{prop:RIP-small-lambda}. These two propositions establish the Restricted Isomorphic Property (RIP) on $T$. Here, by RIP we do not mean that the empirical norm $\|\bX\cdot\|_2$ and the $\ell_2$ (or ellipsoidal) norm are isomorphic on a restricted subset; rather, we mean a more general equivalence, up to a multiplicative constant, between $P_N \cL_{\bbeta_{J_\star}}^{(\lambda,J_\star)}$ and $\|\langle X, \bbeta_{J_\star} - \bbeta_{J_\star}^\star \rangle\|_{L^2(\mathbb{P}_X)}^2$ (which equals $r^2$ under Assumption~\ref{assumption:LASSO_upper}) that holds uniformly when $\bbeta_{J_\star}$ belongs to a certain restricted subset. For our purposes, we only need to establish a one-sided RIP, namely a uniform upper bound on $(P - P_N)\cL_{\bbeta_{J_\star}}^{(\lambda,J_\star)}$ on a restricted subset.

\begin{Proposition}[Uniform convergence, when $\lambda$ is large]\label{prop:RIP}
Assume $\Sigma=I_p$. Assume that $r\le\rho$ and $r\ge C_\beta\|\boldsymbol\beta_{J_\star^c}^\star\|_2$, where $C_\beta\ge1$. Then $\forall \lambda>0$,
\begin{align*}
\bP\bigg\{\sup_{\boldsymbol\beta_{J_\star}\in T}\big|(P-P_N)\mathcal L_{\boldsymbol\beta_{J_\star}}^{(\lambda,J_\star)}\big|\le C\frac{\sqrt{\sigma_\xi^2+r^2}}{\lambda\sqrt N}r\sqrt{|J_\star|}\bigg\}\geq 1-C\exp(-c|J_\star|).
\end{align*}
\end{Proposition}
In terms of proof strategy, the proof of Proposition~\ref{prop:RIP} is via generic chaining, \cite{talagrandUpperLowerBounds2021}. However, since $(P-P_N)\cL_{\bbeta_{J_\star}}^{(\lambda,J_\star)}$ is by no means an empirical process --- it cannot be written as a sum of independent random variables --- the usual generic chaining for empirical processes cannot be applied. Here, we instead prove that its increments have uniformly bounded sub-Gaussian norms, and, noting that $\vzero \in T$, use the generic chaining upper bound for general sub-Gaussian processes, \cite[Theorem 2.7.13]{talagrandUpperLowerBounds2021}, to establish the RIP. Thus, the key to this proof lies in establishing the uniform boundedness of the sub-Gaussian norms of the increments of the stochastic process indexed by $T$.

To establish this property, for any $\bbeta_{J_\star} \in T$, we decompose $(P-P_N)\cL_{\bbeta_{J_\star}}^{(\lambda,J_\star)}$ into multiplier and quadratic components $(P-P_N)\cM_{\bbeta_{J_\star}}$ and $(P-P_N)\cQ_{\bbeta_{J_\star}}$, where $(P-P_N)\cM_{\bbeta_{J_\star}}$ is the process generated by the gradient of $P_N\ell_{\bbeta_{J_\star}^\star}^{(\lambda,J_\star)}$, and $(P-P_N)\cQ_{\bbeta_{J_\star}}$ is defined as $(P-P_N)\cL_{\bbeta_{J_\star}}^{(\lambda,J_\star)} - (P-P_N)\cM_{\bbeta_{J_\star}}$. However, unlike the classical multiplier and quadratic processes decomposition (e.g., \cite{lecueLearningSubgaussianClasses2016}), $(P-P_N)\cM_{\bbeta_{J_\star}}$ and $(P-P_N)\cQ_{\bbeta_{J_\star}}$ also cannot be written as sums of independent random variables. Nevertheless, $(P-P_N)\cM_{\bbeta_{J_\star}}$ can still be regarded as a multiplier process with multipliers independent of the design vectors, and hence we can use \cite[Lemma 2]{shangFeatureSpaceDecomposition2026}; for $(P-P_N)\cQ_{\bbeta_{J_\star}}$, however, a new treatment is needed: we write $(P-P_N)\cL_{\bbeta_{J_\star}}^{(\lambda,J_\star)}$ in a restricted variational representation (that is, \eqref{eq:dual-residual-representation} later), then use Danskin's theorem to establish properties of its gradient, and, together with the fundamental theorem of calculus and the definition of $(P-P_N)\cQ_{\bbeta_{J_\star}}$, relate the upper bound on the increments of $(P-P_N)\cQ_{\bbeta_{J_\star}}$ to the increments of this gradient.
\begin{proof}[Proof of Proposition~\ref{prop:RIP}.]
Fix $\lambda > 0$. Throughout the proof, we write $P\cdot=\mathbb E[\,\cdot\mid\mathbb X_{J_\star^c}]$. Let $\vzeta=\bX_{J_\star^c}\bbeta_{J_\star^c}^\star+\bxi$. We investigate the stochastic processes $(P_N\cM_{\bbeta_{J_\star}})_{\bbeta_{J_\star}\in T}$ and $(P_N\cQ_{\bbeta_{J_\star}})_{\bbeta_{J_\star}\in T}$, where
\begin{align}
&P_N\cM_{\bbeta_{J_\star}}:=\left\langle-\boldsymbol\lambda^\star[\vzeta],\mathbb X_{J_\star}(\boldsymbol\beta_{J_\star}-\boldsymbol\beta_{J_\star}^\star)\right\rangle,\quad \mbox{ and }\quad P_N\cQ_{\bbeta_{J_\star}}:=P_N\mathcal L_{\boldsymbol\beta_{J_\star}}^{(\lambda,J_\star)}-P_N\cM_{\bbeta_{J_\star}}.\label{eq:modified-process}
\end{align}
Independence and $\mathbb E\mathbb X_{J_\star}=\vzero$ imply $P(P_N\cM_{\bbeta_{J_\star}})=0$ and $P(P_N\cQ_{\bbeta_{J_\star}})=P\mathcal L_{\boldsymbol\beta_{J_\star}}^{(\lambda,J_\star)}$. We abbreviate $P(P_N\cM_{\bbeta_{J_\star}})$ by $P\cM_{\bbeta_{J_\star}}$ and similarly for $P(P_N\cQ_{\bbeta_{J_\star}})$. Consequently,
\begin{align}\label{eq:original-modified-decomposition}
&(P_N-P)\mathcal L_{\boldsymbol\beta_{J_\star}}^{(\lambda,J_\star)}=P_N\cM_{\bbeta_{J_\star}}+(P_N-P)\cQ_{\bbeta_{J_\star}}.
\end{align}
Taking absolute values gives
\begin{align*}
&\sup_{\boldsymbol\beta_{J_\star}\in T}\left|(P_N-P)\mathcal L_{\boldsymbol\beta_{J_\star}}^{(\lambda,J_\star)}\right|\le\sup_{\boldsymbol\beta_{J_\star}\in T}|P_N\cM_{\bbeta_{J_\star}}|+\sup_{\boldsymbol\beta_{J_\star}\in T}|(P_N-P)\cQ_{\bbeta_{J_\star}}|.
\end{align*}
Let
\begin{align*}
&\Omega_{\mathrm{res}}=\left\{\|\vzeta\|_2\le\sqrt{3}\sqrt{\sigma_\xi^2+\|\boldsymbol\beta_{J_\star^c}^\star\|_2^2}\sqrt N\right\}.
\end{align*}
Then $\mathbb P(\Omega_{\mathrm{res}}^c)\le\exp(-c|J_\star|)$. By generic chaining and the independence between $\vlambda^\star[\vzeta]$ and $\bX_{J_\star}$, see, for instance, \cite[Lemma 2]{shangFeatureSpaceDecomposition2026}, together with Proposition~\ref{prop:effective-regularization-lasso}, item~\emph{\ref{item:effective_lambda_star}} for the upper bound of $\|\vlambda^\star[\vzeta]\|_2$, we have
\begin{align}\label{eq:multiplier_process_upper}
&\bP\bigg\{\sup_{\boldsymbol\beta_{J_\star}\in T}|P_N\cM_{\bbeta_{J_\star}}|\le C\min\left\{\mathrm{diam}\!\left((\mathbb X_{J_\star^c}B_1^{J_\star^c})^\circ\right),\frac{\sqrt{3}\sqrt{\sigma_\xi^2+\|\boldsymbol\beta_{J_\star^c}^\star\|_2^2}}{\lambda\sqrt N}\right\}r\sqrt{|J_\star|}\bigg\}\geq 1-\exp(-c|J_\star|).
\end{align}
Let $\Omega_{\mathrm{lin}}$ denote the event in \eqref{eq:multiplier_process_upper}.

Let $\Omega_{\mathrm{op}}=\{\|\mathbb X_{J_\star}\|_{\mathrm{op}}\le3\sqrt N\}$. Then $\mathbb P(\Omega_{\mathrm{op}}^c)\le\exp(-c|J_\star|)$, and we write $P^{\Omega_{\mathrm{op}}}=\mathbb E[\,\cdot\mid\mathbb X_{J_\star^c},\Omega_{\mathrm{op}}]$. Let
\begin{align}\label{eq:dual-residual-representation}
    F:t\in[0,1] \mapsto \max_{\boldsymbol a\in(\mathbb X_{J_\star^c}B_1^{J_\star^c})^\circ}\left\{\left\langle\boldsymbol a,\vzeta-t\mathbb X_{J_\star}(\boldsymbol\beta_{J_\star}-\boldsymbol\beta_{J_\star}^\star)\right\rangle-\frac{N\lambda}{2}\|\boldsymbol a\|_2^2\right\}.
\end{align}
By \eqref{eq:dual_problem_effective_loss}, $P_N\mathcal L_{\boldsymbol\beta_{J_\star}}^{(\lambda,J_\star)}=F(1)-F(0)$, and the maximizer of $F(t)$ is unique and equal to $\boldsymbol\lambda^\star[\vzeta-t\mathbb X_{J_\star}(\boldsymbol\beta_{J_\star}-\boldsymbol\beta_{J_\star}^\star)]$. Thus, by Danskin's theorem, for almost every $t\in[0,1]$, $\frac{\mathrm{d}}{\mathrm{d}t}F(t)=-\langle\vlambda^\star[\vzeta-t\bX_{J_\star}(\bbeta_{J_\star}-\bbeta_{J_\star}^\star)],\bX_{J_\star}(\bbeta_{J_\star}-\bbeta_{J_\star}^\star)\rangle$.
The fundamental theorem of calculus and \eqref{eq:modified-process}  give $P_N\cL_{\bbeta_{J_\star}}^{(\lambda,J_\star)}=\int_0^1 \frac{\mathrm{d}}{\mathrm{d}t}F(t)\mathrm{d}t = P_N\cQ_{\bbeta_{J_\star}}+P_N\cM_{\bbeta_{J_\star}} = P_N\cQ_{\bbeta_{J_\star}} + \langle-\vlambda^\star[\vzeta],\bX_{J_\star}(\bbeta_{J_\star}-\bbeta_{J_\star}^\star)\rangle$, and  therefore
\begin{align}
&P_N\cQ_{\bbeta_{J_\star}}=\int_0^1\Big\langle\boldsymbol\lambda^\star[\vzeta]-\boldsymbol\lambda^\star[\vzeta-t\mathbb X_{J_\star}(\boldsymbol\beta_{J_\star}-\boldsymbol\beta_{J_\star}^\star)],\mathbb X_{J_\star}(\boldsymbol\beta_{J_\star}-\boldsymbol\beta_{J_\star}^\star)\Big\rangle\,\mathrm dt.\label{eq:original-smooth-convexity}
\end{align}
Recall that by Proposition~\ref{prop:effective-regularization-lasso}, $\vlambda^\star[\bullet] = \Proj_{(\bX_{J_\star^c}B_1^{J_\star^c})^\circ}(\frac{1}{N\lambda}\bullet)$. For every $t\in(0,1]$, the nonexpansiveness of the metric projection, the Cauchy--Schwarz inequality, and the Lipschitzness of $\vlambda^\star[\cdot]$ from Proposition~\ref{prop:effective-regularization-lasso} give
\begin{align*}
&0\le\frac{N\lambda}{t}\|\boldsymbol\lambda^\star[\vzeta]-\boldsymbol\lambda^\star[\vzeta-t\mathbb X_{J_\star}(\boldsymbol\beta_{J_\star}-\boldsymbol\beta_{J_\star}^\star)]\|_2^2\\
&\le\Big\langle\boldsymbol\lambda^\star[\vzeta]-\boldsymbol\lambda^\star[\vzeta-t\mathbb X_{J_\star}(\boldsymbol\beta_{J_\star}-\boldsymbol\beta_{J_\star}^\star)],\mathbb X_{J_\star}(\boldsymbol\beta_{J_\star}-\boldsymbol\beta_{J_\star}^\star)\Big\rangle\le\frac{t}{N\lambda}\|\mathbb X_{J_\star}(\boldsymbol\beta_{J_\star}-\boldsymbol\beta_{J_\star}^\star)\|_2^2.
\end{align*}
Integration in \eqref{eq:original-smooth-convexity} yields, almost surely, for every $\boldsymbol\beta_{J_\star}\in T$,
\begin{align*}
&0\le P_N\cQ_{\bbeta_{J_\star}}\le\frac{\|\mathbb X_{J_\star}(\boldsymbol\beta_{J_\star}-\boldsymbol\beta_{J_\star}^\star)\|_2^2}{2N\lambda}.
\end{align*}
This nonnegativity justifies our terminology ``quadratic component'' for $P_N\cQ_{\bbeta_{J_\star}}$ --- it is an analogue of the nonnegative stochastic process known as the quadratic process; however, the bound $P_N\cQ_{\bbeta_{J_\star}} \le \frac{\|\mathbb X_{J_\star}(\boldsymbol\beta_{J_\star} - \boldsymbol\beta_{J_\star}^\star)\|_2^2}{2N\lambda}$ is not sharp enough. We may see $P_N\cQ_{\bbeta_J^\star}\geq 0$ easily from the fact that $P_N\ell_\bullet^{(\lambda,J_\star)}$ is convex. We next examine, similarly, the uniform boundedness of the sub-Gaussian norms of the increments of $((P-P_N)\cQ_{\bbeta_{J_\star}})_{\bbeta_{J_\star} \in T}$. Our tool is the (restricted) Gaussian Lipschitz concentration inequality. To this end, we first examine the gradient of the increments. Recall from \eqref{eq:modified-process} that $P_N\cQ_{\bbeta_{J_\star}} = P_N\cL_{\bbeta_{J_\star}}^{(\lambda,J_\star)} - P_N\cM_{\bbeta_{J_\star}}$, then we compute the gradients of $P_N\cL_{\bbeta_{J_\star}}^{(\lambda,J_\star)}$ and that of $P_N\cM_{\bbeta_{J_\star}}^{(\lambda,J_\star)}$ to compute the gradient of $P_N\cQ_{\bbeta_{J_\star}}$.
\begin{enumerate}
    \item For any $\bbeta_{J_\star},\bbeta_{J_\star}'\in T$ and almost every $\boldsymbol\xi$, Danskin's theorem yields $\nabla_\bxi P_N\cL_{\bbeta_{J_\star}}^{(\lambda,J_\star)} = \vlambda^\star[\vzeta - \bX_{J_\star}(\bbeta_{J_\star}-\bbeta_{J_\star}^\star)]-\vlambda^\star[\bxi]$, and the chain rule gives $\nabla_\bxi P_N\cM_{\bbeta_{J_\star}} = -(\D\vlambda^\star[\cdot](\vzeta))\bX_{J_\star}(\bbeta_{J_\star}-\bbeta_{J_\star}^\star)$, hence
    \begin{align}\label{eq:noise-gradient} 
    \begin{aligned}
    &\nabla_{\boldsymbol\xi}(P_N\cQ_{\bbeta_{J_\star}}-P_N\cQ_{\bbeta'_{J_\star}})\\
    &=\boldsymbol\lambda^\star[\vzeta-\mathbb X_{J_\star}(\boldsymbol\beta_{J_\star}-\boldsymbol\beta_{J_\star}^\star)]-\boldsymbol\lambda^\star[\vzeta-\mathbb X_{J_\star}(\boldsymbol\beta'_{J_\star}-\boldsymbol\beta_{J_\star}^\star)] +\big(\D\boldsymbol\lambda^\star[\cdot](\vzeta)\big)\mathbb X_{J_\star}(\boldsymbol\beta_{J_\star}-\boldsymbol\beta'_{J_\star}).
    \end{aligned}
    \end{align}
    On $\Omega_{\mathrm{op}}$,
    \begin{align*}
    &\|\boldsymbol\lambda^\star[\vzeta-\mathbb X_{J_\star}(\boldsymbol\beta_{J_\star}-\boldsymbol\beta_{J_\star}^\star)]-\boldsymbol\lambda^\star[\vzeta-\mathbb X_{J_\star}(\boldsymbol\beta'_{J_\star}-\boldsymbol\beta_{J_\star}^\star)]\|_2\le\frac{\|\mathbb X_{J_\star}(\boldsymbol\beta_{J_\star}-\boldsymbol\beta'_{J_\star})\|_2}{N\lambda}\le\frac{3}{\lambda\sqrt N}\|\boldsymbol\beta_{J_\star}-\boldsymbol\beta'_{J_\star}\|_2,
    \end{align*}
    and $\|\big(\D\boldsymbol\lambda^\star[\cdot](\vzeta)\big)\mathbb X_{J_\star}(\boldsymbol\beta_{J_\star}-\boldsymbol\beta'_{J_\star})\|_2\le\frac{3}{\lambda\sqrt N}\|\boldsymbol\beta_{J_\star}-\boldsymbol\beta'_{J_\star}\|_2.$ Thus \eqref{eq:noise-gradient} implies
    \begin{align*}
    &\|\nabla_{\boldsymbol\xi}(P_N\cQ_{\bbeta_{J_\star}}-P_N\cQ_{\bbeta'_{J_\star}})\|_2\le\frac{6}{\lambda\sqrt N}\|\boldsymbol\beta_{J_\star}-\boldsymbol\beta'_{J_\star}\|_2.
    \end{align*}

    \item Similarly, $\nabla_{\bX_{J_\star}}P_N\cL_{\bbeta_{J_\star}}^{(\lambda,J_\star)} = -\vlambda^\star[\vzeta-\bX_{J_\star}(\bbeta_{J_\star}-\bbeta_{J_\star}^\star)]\otimes(\bbeta_{J_\star}-\bbeta_{J_\star}^\star)$, and $\nabla_{\bX_{J_\star}}P_N\cM_{\bbeta_{J_\star}} = -\vlambda^\star[\vzeta]\otimes(\bbeta_{J_\star}-\bbeta_{J_\star}^\star)$, hence
    \begin{align}\label{eq:design-gradient}
    \begin{aligned}
    &\nabla_{\mathbb X_{J_\star}}(P_N\cQ_{\bbeta_{J_\star}}-P_N\cQ_{\bbeta'_{J_\star}})\\
    &=\Big(\boldsymbol\lambda^\star[\vzeta]-\boldsymbol\lambda^\star[\vzeta-\mathbb X_{J_\star}(\boldsymbol\beta_{J_\star}-\boldsymbol\beta_{J_\star}^\star)]\Big)\otimes(\boldsymbol\beta_{J_\star}-\boldsymbol\beta'_{J_\star})\\
    &+\Big(\boldsymbol\lambda^\star[\vzeta-\mathbb X_{J_\star}(\boldsymbol\beta'_{J_\star}-\boldsymbol\beta_{J_\star}^\star)]-\boldsymbol\lambda^\star[\vzeta-\mathbb X_{J_\star}(\boldsymbol\beta_{J_\star}-\boldsymbol\beta_{J_\star}^\star)]\Big)\otimes(\boldsymbol\beta'_{J_\star}-\boldsymbol\beta_{J_\star}^\star).
    \end{aligned}
    \end{align}
    Since $\|\boldsymbol\beta_{J_\star}-\boldsymbol\beta_{J_\star}^\star\|_2,\|\boldsymbol\beta'_{J_\star}-\boldsymbol\beta_{J_\star}^\star\|_2\le r$, the Frobenius norms of the first and second terms on the right-hand side of \eqref{eq:design-gradient} are bounded, respectively, by $\frac{\|\mathbb X_{J_\star}(\boldsymbol\beta_{J_\star}-\boldsymbol\beta_{J_\star}^\star)\|_2}{N\lambda}\|\boldsymbol\beta_{J_\star}-\boldsymbol\beta'_{J_\star}\|_2\le\frac{3r}{\lambda\sqrt N}\|\boldsymbol\beta_{J_\star}-\boldsymbol\beta'_{J_\star}\|_2,$ and $\frac{\|\mathbb X_{J_\star}(\boldsymbol\beta_{J_\star}-\boldsymbol\beta'_{J_\star})\|_2}{N\lambda}\|\boldsymbol\beta'_{J_\star}-\boldsymbol\beta_{J_\star}^\star\|_2\le\frac{3r}{\lambda\sqrt N}\|\boldsymbol\beta_{J_\star}-\boldsymbol\beta'_{J_\star}\|_2.$ Therefore,
    \begin{align*}
    &\|\nabla_{\mathbb X_{J_\star}}(P_N\cQ_{\bbeta_{J_\star}}-P_N\cQ_{\bbeta'_{J_\star}})\|_{\mathrm F}\le\frac{6r}{\lambda\sqrt N}\|\boldsymbol\beta_{J_\star}-\boldsymbol\beta'_{J_\star}\|_2.
    \end{align*}
\end{enumerate}
Notice that the probability distribution of $\bX_{J_\star}$ condition on $\Omega_{\op}$ is the normalized Gaussian distribution restricted on the closed convex set $\{\vA:\|\vA\|_{\op}\leq 3\sqrt{N}\}$, hence by Lemma~\ref{lem:restricted_Gaussian_Lipschitz},
\begin{align*}
&\left\|(P_N\cQ_{\bbeta_{J_\star}}-P_N\cQ_{\bbeta'_{J_\star}})-P^{\Omega_{\mathrm{op}}}(P_N\cQ_{\bbeta_{J_\star}}-P_N\cQ_{\bbeta'_{J_\star}})\right\|_{\psi_2(P^{\Omega_{\op}})}\le C\frac{\sqrt{\sigma_\xi^2+r^2}}{\lambda\sqrt N}\|\boldsymbol\beta_{J_\star}-\boldsymbol\beta'_{J_\star}\|_2.
\end{align*}
Since $P_N\cQ_{\boldsymbol\beta_{J_\star}^\star}=0$, generic chaining gives, with $P^{\Omega_{\mathrm{op}}}$-probability at least $1-2\exp(-c|J_\star|)$,
\begin{align}\label{eq:conditional-chaining}
&\sup_{\boldsymbol\beta_{J_\star}\in T}\left|P_N\cQ_{\bbeta_{J_\star}}-P^{\Omega_{\mathrm{op}}}(P_N\cQ_{\bbeta_{J_\star}})\right|\le C\frac{\sqrt{\sigma_\xi^2+r^2}}{\lambda\sqrt N}r\sqrt{|J_\star|}.
\end{align}
We are left with comparing $P^{\Omega_{\op}}(P_N\cQ_{\bbeta_{J_\star}})$ with $P\cQ_{\bbeta_{J_\star}}$. The law of total expectation gives
\begin{align*}
    \forall\bbeta_{J_\star}\in T,\quad & P\cQ_{\bbeta_{J_\star}} = \bP(\Omega_{\op})P^{\Omega_{\op}}(P_N\cQ_{\bbeta_{J_\star}}) + P\big(P_N\cQ_{\bbeta_{J_\star}}\1_{\Omega_{\op}^c}\big),
\end{align*}hence $|P^{\Omega_{\op}}(P_N\cQ_{\bbeta_{J_\star}}) - P\cQ_{\bbeta_{J_\star}}| \leq \bP(\Omega_{\op}^c)P^{\Omega_{\op}}(P_N\cQ_{\bbeta_{J_\star}}) + P(P_N\cQ_{\bbeta_{J_\star}}\1_{\Omega_{\op}^c})$, and from $0\le P_N\cQ_{\bbeta_{J_\star}}\le \frac{r^2\|\mathbb X_{J_\star}\|_{\mathrm{op}}^2}{2N\lambda}$, we obtain
\begin{align*}
&\sup_{\boldsymbol\beta_{J_\star}\in T}\left|P^{\Omega_{\mathrm{op}}}(P_N\cQ_{\bbeta_{J_\star}})-P\cQ_{\bbeta_{J_\star}}\right|\le\frac{r^2}{2N\lambda}\mathbb E\!\left[(\|\mathbb X_{J_\star}\|_{\mathrm{op}}^2+9N)\mathbf 1_{\Omega_{\mathrm{op}}^c}\right].
\end{align*}
For every $t\ge3\sqrt N$, $\mathbb P\{\|\mathbb X_{J_\star}\|_{\mathrm{op}}>t\}\le\exp(-\frac{(t-\sqrt N-\sqrt{|J_\star|})^2}{2}).$ Let $A$ be a positive random variable. The identity $\bE[A^2\1_{\{Y>a\}}]=a^2\bP\{A>a\}+\int_a^\infty 2t\bP\{A>t\}\mathrm{d}t$ obtain from Fubini applied to $A=\|\bX_{J_\star}\|_{\op}$ and $a=3\sqrt{N}$ yields
\begin{align*}
&\mathbb E\!\left[(\|\mathbb X_{J_\star}\|_{\mathrm{op}}^2+9N)\mathbf 1_{\Omega_{\mathrm{op}}^c}\right]=18N\mathbb P\{\|\mathbb X_{J_\star}\|_{\mathrm{op}}>3\sqrt N\}+\int_{3\sqrt N}^{\infty}2t\,\mathbb P\{\|\mathbb X_{J_\star}\|_{\mathrm{op}}>t\}\,\mathrm dt\\
&\le18N\exp\left(-\frac{(2\sqrt N-\sqrt{|J_\star|})^2}{2}\right)+\int_{2\sqrt N-\sqrt{|J_\star|}}^\infty2(u+\sqrt N+\sqrt{|J_\star|})e^{-\frac{u^2}{2}}\,\mathrm du\\
&\le\left(18N+2+\frac{2(\sqrt N+\sqrt{|J_\star|})}{2\sqrt N-\sqrt{|J_\star|}}\right)\exp\left(-\frac{(2\sqrt N-\sqrt{|J_\star|})^2}{2}\right)\le22N\exp\left(-\frac{361N}{200}\right).
\end{align*}
Therefore, $\sup_{\boldsymbol\beta_{J_\star}\in T}|P^{\Omega_{\mathrm{op}}}(P_N\cQ_{\bbeta_{J_\star}})-P\cQ_{\bbeta_{J_\star}}|$ is at most $\frac{11r^2}{\lambda}\exp\left(-\frac{361N}{200}\right)$, which is absorbed into $C\frac{\sqrt{\sigma_\xi^2+r^2}}{\lambda\sqrt N}r\sqrt{|J_\star|}$. Together with \eqref{eq:conditional-chaining}, this gives an event $\Omega_{\mathrm{Br}}$ such that $\mathbb P(\Omega_{\mathrm{op}}\cap\Omega_{\mathrm{Br}}^c\mid\mathbb X_{J_\star^c})\le C\exp(-c|J_\star|)$ and, on $\Omega_{\mathrm{op}}\cap\Omega_{\mathrm{Br}}$,
\begin{align*}
&\sup_{\boldsymbol\beta_{J_\star}\in T}|(P_N-P)\cQ_{\bbeta_{J_\star}}|\le C\frac{\sqrt{\sigma_\xi^2+r^2}}{\lambda\sqrt N}r\sqrt{|J_\star|}.
\end{align*}
Finally, let
\begin{align*}
&\Omega=\Omega_{\mathrm{op}}\cap\Omega_{\mathrm{res}}\cap\Omega_{\mathrm{lin}}\cap\Omega_{\mathrm{Br}},\quad\mbox{ then }\quad \mathbb P(\Omega)\ge1-C\exp(-c|J_\star|).
\end{align*}
On $\Omega$, substitution of the bounds for $P_N\cM_{\bbeta_{J_\star}}$ and $(P_N-P)\cQ_{\bbeta_{J_\star}}$ into \eqref{eq:original-modified-decomposition}, together with $r\ge C_\beta\|\boldsymbol\beta_{J_\star^c}^\star\|_2$, proves the proposition.
\end{proof}

The RIP established in Proposition~\ref{prop:RIP} holds for every $\lambda > 0$, but it is sharp mainly in the regime $\lambda \ge \sigma_\xi \sqrt{\frac{1}{N}\log(\frac{e|J_\star^c|}{N})}$; in the regime $\frac{\sigma_\xi}{\sqrt N} \le \lambda \le \sigma_\xi \sqrt{\frac{1}{N}\log(\frac{e|J_\star^c|}{N})}$, the resulting upper bound is too large to combine with the local Bernstein's condition obtained in Proposition~\ref{prop:local-bernstein} to yield the desired RIP. Therefore, in the following Proposition~\ref{prop:RIP-small-lambda}, we establish a one-sided RIP that is effective in the regime $\frac{\sigma_\xi}{\sqrt N} \le \lambda \le \sigma_\xi \sqrt{\frac{1}{N}\log(\frac{e|J_\star^c|}{N})}$.

\begin{Proposition}[Uniform convergence, when $\lambda$ is small]\label{prop:RIP-small-lambda}Grant Assumption~\ref{assumption:LASSO_upper}.
Let  $0<r\le\rho$ and $\frac{\sigma_\xi}{\sqrt N}\le\lambda\le \sigma_\xi\sqrt{\frac1N\log(\frac{e|J_\star^c|}{N})}$. Conditionally on $\bX_{J_\star^c}$, with probability at least $1-C\exp(-c|J_\star|)$,
\begin{align}
&\sup_{\boldsymbol\beta_{J_\star}\in B_1^{J_\star}(\boldsymbol\beta_{J_\star}^\star;\rho)\cap S_2^{J_\star}(\boldsymbol\beta_{J_\star}^\star;r)}\left|(P-P_N)\mathcal L_{\boldsymbol\beta_{J_\star}}^{(\lambda,J_\star)}\right|\le C\,\mathrm{diam}\!\left((\mathbb X_{J_\star^c}B_1^{J_\star^c})^\circ\right)r\sqrt{|J_\star|}+C\frac{r^2\sqrt{|J_\star|}}{\sqrt{\sigma_\xi^2+r^2}+\sigma_\xi}.\label{eq:RIP-small-lambda-bound}
\end{align}
Consequently, with $\bP$-probability at least $1-C\exp(-c|J_\star|)$,
\begin{align*}
&\sup_{\boldsymbol\beta_{J_\star}\in B_1^{J_\star}(\boldsymbol\beta_{J_\star}^\star;\rho)\cap S_2^{J_\star}(\boldsymbol\beta_{J_\star}^\star;r)}\left|(P-P_N)\mathcal L_{\boldsymbol\beta_{J_\star}}^{(\lambda,J_\star)}\right|\le \frac{2C}{c_G}r\sqrt{\frac{|J_\star|}{\log(\frac{e|J_\star^c|}{N})}}+C\frac{r^2\sqrt{|J_\star|}}{\sqrt{\sigma_\xi^2+r^2}+\sigma_\xi}.
\end{align*}
In particular, if $r\ge C_\beta\|\boldsymbol\beta_{J_\star^c}^\star\|_2$,
\begin{equation}\label{eq:cond_r_rho_prop5}
r\ge\frac{8C}{c_Gc_{\mathrm{LB}}}\sqrt{\sigma_\xi^2+r^2}\sqrt{\frac{|J_\star|}{N}} \mbox{ and } |J_\star|\log\left(\frac{e|J_\star^c|}{N}\right)\le\frac{c_{\mathrm{LB}}^2}{16C^2}N,
\end{equation}
then, with probability at least $1-C\exp(-c|J_\star|)$,
\begin{align}
&\inf_{\boldsymbol\beta_{J_\star}\in B_1^{J_\star}(\boldsymbol\beta_{J_\star}^\star;\rho)\cap S_2^{J_\star}(\boldsymbol\beta_{J_\star}^\star;r)}P_N\mathcal L_{\boldsymbol\beta_{J_\star}}^{(\lambda,J_\star)}\ge\frac{c_{\mathrm{LB}}}{2}r^2\sqrt{\frac{N}{(\sigma_\xi^2+r^2)\log(\frac{e|J_\star^c|}{N})}}.\label{eq:RIP-small-lambda-lower}
\end{align}
\end{Proposition}

The proof of Proposition~\ref{prop:RIP-small-lambda} is similar to that of Proposition~\ref{prop:RIP}. The difference is that the index set of the stochastic process is now restricted to $B_1^{J_\star}(\bbeta_{J_\star}^\star;\rho) \cap S_2^{J_\star}(\bbeta_{J_\star}^\star;r)$, the set of almost sparse vectors; hence we can use spherical coordinates and the covariance properties of Gaussian random vectors to examine the derivative of the function $F$ defined in \eqref{eq:dual-residual-representation}, thereby avoiding the operator norm bound of $\D\vlambda^\star[\cdot]$. Compared with Proposition~\ref{prop:RIP}, the proof of Proposition~\ref{prop:RIP-small-lambda} includes an extra step: to use generic chaining \cite[Theorem 2.7.13]{talagrandUpperLowerBounds2021}, we need to verify the concentration property at a single point (since $\vzero$ is no longer in $B_1^{J_\star}(\bbeta_{J_\star}^\star;\rho) \cap S_2^{J_\star}(\bbeta_{J_\star}^\star;r)$).

\begin{proof}
Throughout the proof, we condition on $\mathbb X_{J_\star^c}$ and retain the notation $P=\mathbb E[\,\cdot\mid\mathbb X_{J_\star^c}]$. We also retain the processes $P_N\cQ_{\bbeta_{J_\star}}$ and $P_N\cM_{\bbeta_{J_\star}}$ from \eqref{eq:modified-process} and the decomposition \eqref{eq:original-modified-decomposition}.

The proof of Proposition~\ref{prop:RIP} gives
\begin{align}\tag{\ref{eq:multiplier_process_upper}}
\bP\bigg\{\sup_{\boldsymbol\beta_{J_\star}\in B_1^{J_\star}(\boldsymbol\beta_{J_\star}^\star;\rho)\cap S_2^{J_\star}(\boldsymbol\beta_{J_\star}^\star;r)}|P_N\cM_{\bbeta_{J_\star}}|\le C\,\mathrm{diam}\!\left((\mathbb X_{J_\star^c}B_1^{J_\star^c})^\circ\right)r\sqrt{|J_\star|}\bigg| \bX_{J_\star^c}\bigg\}\geq 1-C\exp(-c|J_\star|).
\end{align}

We next control $(P_N-P)\cQ_{\bbeta_{J_\star}}$. Fix $\boldsymbol\beta_{J_\star},\boldsymbol\beta'_{J_\star}\in B_1^{J_\star}(\boldsymbol\beta_{J_\star}^\star;\rho)\cap S_2^{J_\star}(\boldsymbol\beta_{J_\star}^\star;r)$. Then
\begin{align*}
&\operatorname{Cov}\!\begin{pmatrix}\mathbb X_{J_\star}(\boldsymbol\beta_{J_\star}-\boldsymbol\beta_{J_\star}^\star)\\ \mathbb X_{J_\star}(\boldsymbol\beta'_{J_\star}-\boldsymbol\beta_{J_\star}^\star)\end{pmatrix}=\begin{pmatrix}r^2I_N&\langle\boldsymbol\beta_{J_\star}-\boldsymbol\beta_{J_\star}^\star,\boldsymbol\beta'_{J_\star}-\boldsymbol\beta_{J_\star}^\star\rangle I_N\\ \langle\boldsymbol\beta_{J_\star}-\boldsymbol\beta_{J_\star}^\star,\boldsymbol\beta'_{J_\star}-\boldsymbol\beta_{J_\star}^\star\rangle I_N&r^2I_N\end{pmatrix},
\end{align*}
and hence for independent $\boldsymbol G_1,\boldsymbol G_2\sim\mathcal N(0,I_N)$,
\begin{align*}
&\left(\mathbb X_{J_\star}(\boldsymbol\beta_{J_\star}-\boldsymbol\beta_{J_\star}^\star),\mathbb X_{J_\star}(\boldsymbol\beta'_{J_\star}-\boldsymbol\beta_{J_\star}^\star)\right)\overset{\mathrm d}=r\left(\boldsymbol G_1,\frac{\langle\boldsymbol\beta_{J_\star}-\boldsymbol\beta_{J_\star}^\star,\boldsymbol\beta'_{J_\star}-\boldsymbol\beta_{J_\star}^\star\rangle}{r^2}\boldsymbol G_1+\sqrt{1-\frac{\langle\boldsymbol\beta_{J_\star}-\boldsymbol\beta_{J_\star}^\star,\boldsymbol\beta'_{J_\star}-\boldsymbol\beta_{J_\star}^\star\rangle^2}{r^4}}\,\boldsymbol G_2\right).
\end{align*}
By \eqref{eq:dual-residual-representation}, \eqref{eq:modified-process}, and Danskin's theorem, see, for instance, \cite[Theorem 4.13, pp. 273]{bonnansPerturbationAnalysisOptimization2000}, for $0\le t\le$, let
\begin{align*}
    G(t) &= \max_{\boldsymbol a\in(\mathbb X_{J_\star^c}B_1^{J_\star^c})^\circ}\left\{\left\langle\boldsymbol a,\vzeta-r(\cos(t)\boldsymbol G_1+\sin(t)\boldsymbol G_2)\right\rangle-\frac{N\lambda}{2}\|\boldsymbol a\|_2^2\right\} -\max_{\boldsymbol a\in(\mathbb X_{J_\star^c}B_1^{J_\star^c})^\circ}\left\{\langle\boldsymbol a,\vzeta\rangle-\frac{N\lambda}{2}\|\boldsymbol a\|_2^2\right\}\\
    &\qquad+\left\langle\boldsymbol\lambda^\star[\vzeta],r(\cos(t)\boldsymbol G_1+\sin(t)\boldsymbol G_2)\right\rangle.
\end{align*}Then together with \eqref{eq:original-modified-decomposition}, we know that $G(0)\overset{\mathrm{d}}{=}P_N\cQ_{\bbeta_{J_\star}}$ and $G(\theta)\overset{\mathrm{d}}{=}P_N\cQ_{\bbeta_{J_\star}'}$. Moreover,
\begin{align}\label{eq:spherical-path-derivative}
\begin{aligned}
&\frac{\mathrm d}{\mathrm dt}G(t) =r\Big\langle\boldsymbol\lambda^\star[\vzeta]-\boldsymbol\lambda^\star[\vzeta-r(\cos(t)\boldsymbol G_1+\sin(t)\boldsymbol G_2)],-\sin(t)\boldsymbol G_1+\cos(t)\boldsymbol G_2\Big\rangle.
\end{aligned}
\end{align}
For every $t$, $\cos(t)\boldsymbol G_1+\sin(t)\boldsymbol G_2$ and $-\sin(t)\boldsymbol G_1+\cos(t)\boldsymbol G_2$ are independent standard Gaussian vectors, while $\boldsymbol\lambda^\star[\vzeta],\boldsymbol\lambda^\star[\vzeta-r(\cos(t)\boldsymbol G_1+\sin(t)\boldsymbol G_2)]\in(\mathbb X_{J_\star^c}B_1^{J_\star^c})^\circ$. Hence, for every $q\ge2$, the $L_q(P)$-norm of the right-hand side of \eqref{eq:spherical-path-derivative} is at most $C\,\mathrm{diam}((\mathbb X_{J_\star^c}B_1^{J_\star^c})^\circ)r\sqrt q$. By elementary inequality $\arccos(u)\le\pi\sqrt{\frac{1-u}{2}}$ for $u\in[-1,1]$, together with the fact that $\|\bbeta_{J_\star}-\bbeta_{J_\star}'\|_2^2 = 2r^2(1-\frac{1}{r^2}(\langle\bbeta_{J_\star}-\bbeta_{J_\star}^\star,\bbeta_{J_\star}'-\bbeta_{J_\star}^\star\rangle))$, there holds
\begin{align*}
&r\arccos\left(\frac{\langle\boldsymbol\beta_{J_\star}-\boldsymbol\beta_{J_\star}^\star,\boldsymbol\beta'_{J_\star}-\boldsymbol\beta_{J_\star}^\star\rangle}{r^2}\right)\le\frac\pi2\|\boldsymbol\beta_{J_\star}-\boldsymbol\beta'_{J_\star}\|_2.
\end{align*}
The fundamental theorem of calculus applied to $G(t)$ over the interval $[0,\theta]$ therefore yields
\begin{align*}
&\left\|P_N\cQ_{\bbeta_{J_\star}}-P_N\cQ_{\bbeta'_{J_\star}}\right\|_{\psi_2(P)}\le C\,\mathrm{diam}\!\left((\mathbb X_{J_\star^c}B_1^{J_\star^c})^\circ\right)\,\|\boldsymbol\beta_{J_\star}-\boldsymbol\beta'_{J_\star}\|_2.
\end{align*}
The conditional distribution of $P_N\cQ_{\bbeta_{J_\star}}$ depends on $\boldsymbol\beta_{J_\star}$ only through $\|\boldsymbol\beta_{J_\star}-\boldsymbol\beta_{J_\star}^\star\|_2=r$. Thus $P\cQ_{\bbeta_{J_\star}}=P\cQ_{\bbeta'_{J_\star}}$, and the preceding inequality is equivalently
\begin{align*}
&\left\|\left[P_N\cQ_{\bbeta_{J_\star}}-P\cQ_{\bbeta_{J_\star}}\right]-\left[P_N\cQ_{\bbeta'_{J_\star}}-P\cQ_{\bbeta'_{J_\star}}\right]\right\|_{\psi_2(P)}\le C\,\mathrm{diam}\!\left((\mathbb X_{J_\star^c}B_1^{J_\star^c})^\circ\right)\,\|\boldsymbol\beta_{J_\star}-\boldsymbol\beta'_{J_\star}\|_2.
\end{align*}
To use generic chaining's upper bound for sub-Gaussian increment random process, that is, \cite[Theorem 2.7.13]{talagrandUpperLowerBounds2021}, it remains to control the process at one given point: fix any $\boldsymbol\beta_{J_\star}\in B_1^{J_\star}(\boldsymbol\beta_{J_\star}^\star;\rho)\cap S_2^{J_\star}(\boldsymbol\beta_{J_\star}^\star;r)$. For independent $\boldsymbol G_1,\boldsymbol G_2\sim\mathcal N(0,I_N)$, the vectors
\begin{align*}
&\frac{\sigma_\xi\boldsymbol G_1-r\boldsymbol G_2}{\sqrt{\sigma_\xi^2+r^2}}\qquad\mbox{and}\qquad\frac{r\boldsymbol G_1+\sigma_\xi\boldsymbol G_2}{\sqrt{\sigma_\xi^2+r^2}}
\end{align*}
are independent and standard Gaussian. From \eqref{eq:dual-residual-representation}, $P_N\mathcal L_{\boldsymbol\beta_{J_\star}}^{(\lambda,J_\star)}$ is distributed $H_1+H_2$, where
\begin{align}
H_1&=\max_{\boldsymbol a\in(\mathbb X_{J_\star^c}B_1^{J_\star^c})^\circ}\left\{\left\langle\boldsymbol a,\mathbb X_{J_\star^c}\boldsymbol\beta_{J_\star^c}^\star+\sqrt{\sigma_\xi^2+r^2}\,\boldsymbol G_1\right\rangle-\frac{N\lambda}{2}\|\boldsymbol a\|_2^2\right\}-\max_{\boldsymbol a\in(\mathbb X_{J_\star^c}B_1^{J_\star^c})^\circ}\left\{\left\langle\boldsymbol a,\mathbb X_{J_\star^c}\boldsymbol\beta_{J_\star^c}^\star+\sigma_\xi\boldsymbol G_1\right\rangle-\frac{N\lambda}{2}\|\boldsymbol a\|_2^2\right\},\mbox{ and}\label{eq:fixed-point-first-difference}\\
H_2&=\max_{\boldsymbol a\in(\mathbb X_{J_\star^c}B_1^{J_\star^c})^\circ}\left\{\left\langle\boldsymbol a,\mathbb X_{J_\star^c}\boldsymbol\beta_{J_\star^c}^\star+\sigma_\xi\boldsymbol G_1\right\rangle-\frac{N\lambda}{2}\|\boldsymbol a\|_2^2\right\}\notag\\
&-\max_{\boldsymbol a\in(\mathbb X_{J_\star^c}B_1^{J_\star^c})^\circ}\left\{\left\langle\boldsymbol a,\mathbb X_{J_\star^c}\boldsymbol\beta_{J_\star^c}^\star+\frac{\sigma_\xi^2}{\sqrt{\sigma_\xi^2+r^2}}\boldsymbol G_1+\frac{\sigma_\xi r}{\sqrt{\sigma_\xi^2+r^2}}\boldsymbol G_2\right\rangle-\frac{N\lambda}{2}\|\boldsymbol a\|_2^2\right\}.\label{eq:fixed-point-second-difference}
\end{align}
\begin{enumerate}
    \item We first deal with $|(P-P_N)H_2|$. First, we observe that $P H_2= 0$. Moreover, we observe that
    \begin{align*}
        H_2 = h_2(0) - h_2(\arctan(\frac{r}{\sigma_\xi})),\mbox{ where }h_2(t) := \max_{\va\in(\bX_{J_\star^c}B_1^{J_\star^c})^\circ}\bigg\{ \bigg\langle\va,\bX_{J_\star^c}\bbeta_{J_\star^c}^\star + \sigma_\xi\big(\cos(t)\vG_1+\sin(t)\vG_2\big) \bigg\rangle - \frac{N\lambda}{2}\|\va\|_2^2\bigg\},
    \end{align*}and where we have used the fact that $\cos(\arctan(\frac{r}{\sigma_\xi}))=\frac{\sigma_\xi}{\sqrt{\sigma_\xi^2+r^2}}$, and $\sin(\arctan(\frac{r}{\sigma_\xi}))=\frac{r}{\sqrt{\sigma_\xi^2+r^2}}$. Therefore, the Danskin's theorem gives
    \begin{align}\label{eq:fixed-point-rotation-derivative}
    &\frac{\mathrm d}{\mathrm dt}h_2(t) =\sigma_\xi\Big\langle\boldsymbol\lambda^\star[\mathbb X_{J_\star^c}\boldsymbol\beta_{J_\star^c}^\star+\sigma_\xi(\cos(t)\boldsymbol G_1+\sin(t)\boldsymbol G_2)],-\sin(t)\boldsymbol G_1+\cos(t)\boldsymbol G_2\Big\rangle.
    \end{align}
    Since $\sigma_\xi\arctan(\frac{r}{\sigma_\xi})\le r$, the fundamental theorem of calculus yields that $H_2=\int_0^1\frac{\mathrm{d}}{\mathrm{d}t}h_2(t)\mathrm{d}t$ and by the independence between $\cos(t)\boldsymbol G_1+\sin(t)\boldsymbol G_2$ and $-\sin(t)\boldsymbol G_1+\cos(t)\boldsymbol G_2$; Proposition~\ref{prop:effective-regularization-lasso} for the upper bound of $\|\vlambda^\star[\mathbb X_{J_\star^c}\boldsymbol\beta_{J_\star^c}^\star+\sigma_\xi(\cos(t)\boldsymbol G_1+\sin(t)\boldsymbol G_2)]\|_2$; together with the Lipschitz concentration property of $-\sin(t)\boldsymbol G_1+\cos(t)\boldsymbol G_2$, we obtain
    \begin{align*}
    &\left\|H_2\right\|_{\psi_2(P)}\le C\,\mathrm{diam}\!\left((\mathbb X_{J_\star^c}B_1^{J_\star^c})^\circ\right)r.
    \end{align*}
    
    \item We then deal with $H_1$. By Danskin's theorem,
    \begin{align}\label{eq:fixed-point-first-gradient}
    &\nabla_{\boldsymbol G_1}H_1=\sqrt{\sigma_\xi^2+r^2}\,\boldsymbol\lambda^\star[\mathbb X_{J_\star^c}\boldsymbol\beta_{J_\star^c}^\star+\sqrt{\sigma_\xi^2+r^2}\,\boldsymbol G_1]-\sigma_\xi\boldsymbol\lambda^\star[\mathbb X_{J_\star^c}\boldsymbol\beta_{J_\star^c}^\star+\sigma_\xi\boldsymbol G_1].
    \end{align}
    By \eqref{eq:projection-representation-rip}, the triangular inequality together with Proposition~\ref{prop:effective-regularization-lasso},
    \begin{align}\label{eq:fixed-point-projection-bound}
    \begin{aligned}
        &\Big\|\sqrt{\sigma_\xi^2+r^2}\,\boldsymbol\lambda^\star[\mathbb X_{J_\star^c}\boldsymbol\beta_{J_\star^c}^\star+\sqrt{\sigma_\xi^2+r^2}\,\boldsymbol G_1]-\sigma_\xi\boldsymbol\lambda^\star[\mathbb X_{J_\star^c}\boldsymbol\beta_{J_\star^c}^\star+\sigma_\xi\boldsymbol G_1]\Big\|_2\\
        &\le\left(\sqrt{\sigma_\xi^2+r^2}-\sigma_\xi\right)\|\boldsymbol\lambda^\star[\mathbb X_{J_\star^c}\boldsymbol\beta_{J_\star^c}^\star+\sqrt{\sigma_\xi^2+r^2}\,\boldsymbol G_1]\|_2 + \sigma_\xi\left\|\boldsymbol\lambda^\star[\mathbb X_{J_\star^c}\boldsymbol\beta_{J_\star^c}^\star+\sqrt{\sigma_\xi^2+r^2}\,\boldsymbol G_1]-\boldsymbol\lambda^\star[\mathbb X_{J_\star^c}\boldsymbol\beta_{J_\star^c}^\star+\sigma_\xi\boldsymbol G_1]\right\|_2\\
        &\le\left(\sqrt{\sigma_\xi^2+r^2}-\sigma_\xi\right)\left[\frac12\,\mathrm{diam}\!\left((\mathbb X_{J_\star^c}B_1^{J_\star^c})^\circ\right)+\frac{\sigma_\xi\|\boldsymbol G_1\|_2}{N\lambda}\right].
    \end{aligned}
    \end{align}
\end{enumerate}
By the $L_m$-moment version of the Gaussian Poincar\'e inequality, that is, $\|f-\bE f\|_{L^m}\lesssim \sqrt{m}\|\|\nabla f\|_2\|_{L^m}$, \eqref{eq:fixed-point-first-gradient}, \eqref{eq:fixed-point-projection-bound}, and $\|\|\boldsymbol G_1\|_2\|_{L_m}\le\sqrt N+\sqrt m$ imply that, for every $2\le m\le2N$,
\begin{align*}
&\left\| H_2 - P H_2\right\|_{L_m(P)}\le C\left(\sqrt{\sigma_\xi^2+r^2}-\sigma_\xi\right)\left[\mathrm{diam}\!\left((\mathbb X_{J_\star^c}B_1^{J_\star^c})^\circ\right)+\frac{\sigma_\xi(\sqrt N+\sqrt m)}{N\lambda}\right]\sqrt m.
\end{align*}
Here $\lambda\ge\frac{\sigma_\xi}{\sqrt N}$ and $q\le2N$ imply $\frac{\sigma_\xi(\sqrt N+\sqrt q)}{N\lambda}\le1+\sqrt2<\frac{5}{2}$, while $\sqrt{\sigma_\xi^2+r^2}-\sigma_\xi=\frac{r^2}{\sqrt{\sigma_\xi^2+r^2}+\sigma_\xi}\le r.$ It follows that
\begin{align*}
&\left\|H_2 - PH_2\right\|_{L_m(P)}\le C\left[\mathrm{diam}\!\left((\mathbb X_{J_\star^c}B_1^{J_\star^c})^\circ\right)r+\frac{r^2}{\sqrt{\sigma_\xi^2+r^2}+\sigma_\xi}\right]\sqrt m.
\end{align*}
Finally, conditionally on $(\mathbb X_{J_\star^c},\boldsymbol\xi)$, $P_N\cM_{\bbeta_{J_\star}}$ is centered Gaussian and
\begin{align}\tag{\ref{eq:multiplier_process_upper}}
&\forall m\geq 2,\qquad \|P_N\cM_{\bbeta_{J_\star}}\|_{L_m(P)}\le\frac12\,\mathrm{diam}\!\left((\mathbb X_{J_\star^c}B_1^{J_\star^c})^\circ\right)r\sqrt m.
\end{align}
Since $P_N\cQ_{\bbeta_{J_\star}}-P\cQ_{\bbeta_{J_\star}}=(P_N-P)\mathcal L_{\boldsymbol\beta_{J_\star}}^{(\lambda,J_\star)}-(P_N\cM_{\bbeta_{J_\star}}-P\cM_{\bbeta_{J_\star}})$ and $P\cM_{\bbeta_{J_\star}}=0$, the preceding estimates give,
\begin{align*}
&\forall 2\le q\le2N,\qquad \left\|P_N\cQ_{\bbeta_{J_\star}}-P\cQ_{\bbeta_{J_\star}}\right\|_{L_q(P)}\le C\left[\mathrm{diam}\!\left((\mathbb X_{J_\star^c}B_1^{J_\star^c})^\circ\right)r+\frac{r^2}{\sqrt{\sigma_\xi^2+r^2}+\sigma_\xi}\right]\sqrt q.
\end{align*}
Combining the preceding estimate at the fixed index with the sub-Gaussian increment estimate and $\gamma_2\!(B_1^{J_\star}(\boldsymbol\beta_{J_\star}^\star;\rho)\cap S_2^{J_\star}(\boldsymbol\beta_{J_\star}^\star;r),\|\cdot\|_2)\le\gamma_2(rB_2^{J_\star},\|\cdot\|_2)\le Cr\sqrt{|J_\star|},$
we apply the $L_m$ generic-chaining inequality with $m=2|J_\star|$, see, for instance, \cite[Theorem 3.2]{dirksenTailBoundsGeneric2015}. With probability at least $1-C\exp(-c|J_\star|)$,
\begin{align*}
&\sup_{\boldsymbol\beta_{J_\star}\in B_1^{J_\star}(\boldsymbol\beta_{J_\star}^\star;\rho)\cap S_2^{J_\star}(\boldsymbol\beta_{J_\star}^\star;r)}\left|(P_N-P)\cQ_{\bbeta_{J_\star}}\right|\le C\,\mathrm{diam}\!\left((\mathbb X_{J_\star^c}B_1^{J_\star^c})^\circ\right)r\sqrt{|J_\star|}+C\frac{r^2\sqrt{|J_\star|}}{\sqrt{\sigma_\xi^2+r^2}+\sigma_\xi}.
\end{align*}
Substitution of the preceding bounds for $P_N\cM_{\bbeta_{J_\star}}$ and $(P_N-P)\cQ_{\bbeta_{J_\star}}$ into \eqref{eq:original-modified-decomposition} proves \eqref{eq:RIP-small-lambda-bound}. Moreover, on $\Omega$ (from Lemma~\ref{lem:geometric-event}), we have
\begin{align*}
&(\mathbb X_{J_\star^c}B_1^{J_\star^c})^\circ\subset\frac{1}{c_G\sqrt{\log(\frac{e|J_\star^c|}{N})}}B_2^N,\qquad \mathrm{diam}\!\left((\mathbb X_{J_\star^c}B_1^{J_\star^c})^\circ\right)\le\frac{2}{c_G\sqrt{\log(\frac{e|J_\star^c|}{N})}},
\end{align*}
which proves the second assertion.

It remains to prove \eqref{eq:RIP-small-lambda-lower}. Since $\lambda\le\sigma_\xi\sqrt{\frac{\log(\frac{e|J_\star^c|}{N})}{N}}$,
\begin{align*}
&\frac1\lambda\ge\sqrt{\frac{N}{\sigma_\xi^2\log(\frac{e|J_\star^c|}{N})}}\ge\sqrt{\frac{N}{(\sigma_\xi^2+r^2)\log(\frac{e|J_\star^c|}{N})}}.
\end{align*}
Proposition~\ref{prop:local-bernstein} together with the assumptions on $r$ and $|J_\star|$ imply, respectively,
\begin{align*}
&\frac{2C}{c_G}r\sqrt{\frac{|J_\star|}{\log(\frac{e|J_\star^c|}{N})}}\le\frac{c_{\mathrm{LB}}}{4}r^2\sqrt{\frac{N}{(\sigma_\xi^2+r^2)\log(\frac{e|J_\star^c|}{N})}},\mbox{ and } C\frac{r^2\sqrt{|J_\star|}}{\sqrt{\sigma_\xi^2+r^2}+\sigma_\xi}\le\frac{c_{\mathrm{LB}}}{4}r^2\sqrt{\frac{N}{(\sigma_\xi^2+r^2)\log(\frac{e|J_\star^c|}{N})}}.
\end{align*}
Therefore, \eqref{eq:RIP-small-lambda-bound} and the preceding local Bernstein bound yield \eqref{eq:RIP-small-lambda-lower}. The probability estimate follows from the events established above and the event in Proposition~\ref{prop:local-bernstein}.
\end{proof}

\subsection{Upper bound for the estimation error of the estimation part}

Once the local Bernstein's property (Proposition~\ref{prop:local-bernstein}) and the uniform convergence results (Proposition~\ref{prop:RIP} and Proposition~\ref{prop:RIP-small-lambda}) are established, the remaining steps to obtain the estimation error of the estimation part are routine. Here, the difference from other FSD methods in establishing the estimation error of the estimation part is that we take $\frac{\rho_\star}{r_\star} = 2\sqrt{|J_\star|}$; thus there is no decomposition into almost sparse and well-spread parts for $T = B_1(\bbeta_{J_\star}^\star;\rho_\star) \cap B_2(\bbeta_{J_\star}^\star;r_\star)$ --- we only consider $S_2^{J_\star}(\bbeta_{J_\star}^\star;r_\star)$. By the local Bernstein's property and uniform convergence, we obtain an Isomorphic Property (IP) that holds on the entire $S_2^{J_\star}(\bbeta_{J_\star}^\star;r_\star)$, and consequently the effective empirical excess risk will dominate the excess regularization. Recall that $C_\beta$, $c_{\mathrm{LB}}$, $c_{\mathrm{G}}$ are absolute constants defined in \eqref{eq:fixed-constants}. Define
\begin{align}\label{eq:def_r_star}
    \rho_\star = 2r_\star\sqrt{|J_\star|},\mbox{ and } r_\star = \begin{cases}
        C_\star\bigg(\|\bbeta_{J_\star^c}^\star\|_2+\sigma_\xi\sqrt{\frac{|J_\star|\log(\frac{e|J_\star^c|}{N})}{N}}+\sqrt{\Var_{J_\star^c}}\bigg), & \mbox{ if } \frac{\sigma_\xi}{\sqrt{N}}\leq \lambda \leq \sigma_\xi\sqrt{\frac{\log(\frac{e|J_\star^c|}{N})}{N}},\\
        C_\star\bigg(\|\bbeta_{J_\star^c}^\star\|_2+\sigma_\xi\sqrt{\frac{|J_\star|}{N}}+\lambda\sqrt{|J_\star|}\bigg), &\mbox{ if }\lambda>\sigma_\xi\sqrt{\frac{\log(\frac{e|J_\star^c|}{N})}{N}}.
    \end{cases}
\end{align}
\begin{Proposition}[Estimation error]\label{prop:estimation_part_final}
    Grant the same assumption as in Theorem~\ref{theo:upper_LASSO}. If $r_\star\leq \sigma_\xi$, then
    \begin{align*}
        \bP\bigg\{ \|\hat\bbeta_{J_\star} - \bbeta_{J_\star}^\star\|_2^2 \leq r_\star^2 \bigg\}\geq 1-C\exp(-c|J_\star|).
    \end{align*}
\end{Proposition}
\begin{proof}[Proof of Proposition~\ref{prop:estimation_part_final}]
    We first perform a homogenization argument, which is a routine technique; see, e.g., \cite[Lemma 1]{shangFeatureSpaceDecomposition2026}. Note that we have chosen $\rho_\star = 2r_\star\sqrt{|J_\star|}$, so $T = B_1(\bbeta_{J_\star}^\star;\rho_\star)\cap B_2(\bbeta_{J_\star}^\star;r_\star) = B_2(\bbeta_{J_\star}^\star;r_\star)$. It suffices to prove that for every $\bbeta_{J_\star} \in S_2(\bbeta_{J_\star}^\star;r_\star)$, we have $P_N\cL_{\bbeta_{J_\star}}^{(\lambda,J_\star)} + (\|\bbeta_{J_\star}\|_1 - \|\bbeta_{J_\star}^\star\|_1) > 0$. We then separate into two cases.
    \begin{enumerate}
        \item When $\frac{\sigma_\xi}{\sqrt{N}}\leq \lambda \leq \sigma_\xi\sqrt{\frac{\log(\frac{e|J_\star^c|}{N})}{N}}$. Then Proposition~\ref{prop:RIP-small-lambda} gives that with probability at least $1-C\exp(-c|J_\star|)$, for any such $\bbeta_{J_\star}$, $P_N\cL_{\bbeta_{J_\star}}^{(\lambda,J_\star)}\geq \frac{c_{\mathrm{LB}}}{2\sqrt{2}}\frac{r_\star^2}{\sigma_\xi}\sqrt{\frac{N}{\log(\frac{e|J_\star^c|}{N})}}$. From H{\"o}lder's inequality,
        \begin{align*}
            \forall \bbeta_{J_\star}\in S_2(\bbeta_{J_\star}^\star;r_\star),\qquad \|\bbeta_{J_\star}\|_1 - \|\bbeta_{J_\star}^\star\|_1 \geq -\|\bbeta_{J_\star}-\bbeta_{J_\star}^\star\|_1\geq -r_\star\sqrt{|J_\star|}.
        \end{align*}Together with the definition of $r_\star$ in \eqref{eq:def_r_star},
        \begin{align*}
            \forall \bbeta_{J_\star}\in S_2(\bbeta_{J_\star}^\star;r_\star),\qquad &P_N\cL_{\bbeta_{J_\star}}^{(\lambda,J_\star)} + (\|\bbeta_{J_\star}\|_1 - \|\bbeta_{J_\star}^\star\|_1) > 3\sqrt{|J_\star|}r_\star>0.
        \end{align*}

        \item When $\lambda>\sigma_\xi\sqrt{\frac{\log(\frac{e|J_\star^c|}{N})}{N}}$. Then Proposition~\ref{prop:RIP} together with Proposition~\ref{prop:local-bernstein} give that with probability at least $1-C\exp(-c|J_\star|)$, $P_N\cL_{\bbeta_{J_\star}}^{(\lambda,J_\star)}\geq \frac{c_{\mathrm{LB}}}{2\lambda}r_\star^2$. Moreover, from $r_\star\geq \frac{8\sqrt{2}}{c_{\mathrm{LB}}}\lambda\sqrt{|J_\star|}$, we obtain that $\|\bbeta_{J_\star}\|_1-\|\bbeta_{J_\star}^\star\|_1\geq -\|\bbeta_{J_\star}-\bbeta_{J_\star}^\star\|_1\geq -\rho_\star \geq -\frac{c_{\mathrm{LB}}}{4\sqrt{2}\lambda}r_\star^2$. Together with the definition of $r_\star$ in \eqref{eq:def_r_star},
        \begin{align*}
            \forall \bbeta_{J_\star}\in S_2(\bbeta_{J_\star}^\star;r_\star),\qquad &P_N\cL_{\bbeta_{J_\star}}^{(\lambda,J_\star)} + (\|\bbeta_{J_\star}\|_1 - \|\bbeta_{J_\star}^\star\|_1) > \frac{c_{\mathrm{LB}}}{4}\frac{r_\star^2}{\lambda}>0.
        \end{align*}
    \end{enumerate}
\end{proof}

\subsection{Upper bound of the energy of noise absorption part.}

Grant Assumption~\ref{assumption:LASSO_upper}. In this section, our goal is to prove that, with high probability,
\begin{align*}
    \|\hat\bbeta_{J_\star^c}\|_2^2\leq C r^2(V_{J_\star},V_{J_\star^c}),\mbox{ where } \hat\bbeta_{J_\star^c} = \cA[\vy-\bX_{J_\star}\hat\bbeta_{J_\star}]=\cA[\bxi+\vzeta],\mbox{ and } \vzeta = \bX_{J_\star}(\bbeta_{J_\star}^\star-\hat\bbeta_{J_\star}) + \bX_{J_\star^c}\bbeta_{J_\star^c}^\star.
\end{align*}For any $\bbeta_{J_\star}\in V_{J_\star}$, define
\begin{align*}
    P_N\ell_{(\bbeta_{J_\star},\bullet)}:\bbeta_{J_\star^c}\in V_{J_\star^c}\mapsto P_N\ell_{\bbeta} = \frac{1}{2N\lambda}\norm{\vy-\bX_{J_\star}\bbeta_{J_\star} - \bX_{J_\star^c}\bbeta_{J_\star^c}}_2^2 + \norm{\bbeta_{J_\star^c}}_1 + \norm{\bbeta_{J_\star}}_1, \mbox{ where }\bbeta = (\bbeta_{J_\star},\bbeta_{J_\star^c})
\end{align*}then $\hat\bbeta_{J_\star^c}\in\argmin(P_N\ell_{\hat\bbeta_{J_\star},\bbeta_{J_\star^c}}: \bbeta_{J_\star^c}\in V_{J_\star^c})$.

\subsubsection{When \texorpdfstring{$\lambda\geq 2^8 \sigma_\xi\sqrt{\frac{\log(\frac{e|J_\star^c|}{N})}{N}}$}{lambda is large}}

\begin{Proposition}\label{prop:lambda_large_shift}
This proposition holds for every positive diagonal $\Sigma$.
Assume $N\le |J_\star^c|$, $0<r_\star\le\sigma_\xi$, $\log(\frac{e|J_\star^c|}{N})\le \frac{N}{2^{10}}$, and $\lambda\ge 2^{8}\sigma_\xi\sqrt{\frac{(\max_{j\in J_\star^c}\sigma_j)\log(\frac{e|J_\star^c|}{N})}{N}}$. Let $\boldsymbol\zeta = \bX_{J_\star^c}\bbeta_{J_\star^c}^\star + \bX_{J_\star}(\bbeta_{J_\star}^\star-\hat\bbeta_{J_\star})$, and set $\Omega_{\mathrm{esti}}=\{\|\boldsymbol\zeta\|_2\le\sqrt N r_\star\}$. Then
\begin{align*}
\|\Sigma_{J_\star^c}^{1/2}(\mathcal A[\boldsymbol\xi+\boldsymbol\zeta]-\mathcal A[\boldsymbol\xi])\|_2^2\le 5r_\star^2+2^{12}\frac{(\max_{j\in J_\star^c}\sigma_j)\log(\frac{e|J_\star^c|}{N})}{N\lambda^2}r_\star^4\le 6r_\star^2
\end{align*}
with probability at least $1-\mathbb P(\Omega_{\mathrm{esti}}^c)-3\exp(-\frac{N}{64})$.
\end{Proposition}

\begin{proof} The strategy of this proof is to mimic the proof strategy in \cite[pp. 52 - 54]{chafaiInteractionsCompressedSensing2012}, see also \cite[Section 6.2]{foucartMathematicalIntroductionCompressive2013}, namely, to use the Restricted Isomorphic Property. To this end, we first construct the random events needed for this proof.
Define
\begin{align*}
&\Omega_{\mathrm{RIP}}=\bigg\{\forall I\subset J_\star^c,\ |I|\le 2\bigg\lfloor\frac{N}{2^{9}\log(\frac{e|J_\star^c|}{N})}\bigg\rfloor,\forall\boldsymbol v\in V_I,\ \frac N3\|\Sigma_{J_\star^c}^{1/2}\boldsymbol v\|_2^2\le\|\mathbb X_{J_\star^c}\boldsymbol v\|_2^2\le2N\|\Sigma_{J_\star^c}^{1/2}\boldsymbol v\|_2^2\bigg\}.
\end{align*}Let $\Omega_\xi=\{\|\bxi\|_2\leq 2\sigma_\xi\sqrt{N}\}$ and $\Omega = \Omega_{\mathrm{RIP}}\cap \Omega_\xi \cap \Omega_{\mathrm{esti}}.$ It is well known that, see, for instance, \cite[Section 7.3, pp. 159]{vershyninHighDimensionalProbabilityIntroduction2018}, for every fixed $I\subset J_\star^c$ with $|I|\le N$ and every $t\ge0$,
\begin{align*}
&\mathbb P\bigl\{s_{\min}(\mathbb X_I\Sigma_I^{-1/2})\le\sqrt N-\sqrt{|I|}-t\bigr\}\le\exp(-\frac{t^2}{2}),\mbox{ and } \mathbb P\bigl\{\|\mathbb X_I\Sigma_I^{-1/2}\|_{\mathrm{op}}\ge\sqrt N+\sqrt{|I|}+t\bigr\}\le\exp(-\frac{t^2}{2}).
\end{align*}
Taking $t=\frac{\sqrt N}{3}$ for the smallest singular value, $t=\frac{\sqrt N}{3}$ for the largest singular value, and using union bound, and $\binom{|J_\star^c|}{k}\le(\frac{e|J_\star^c|}{k})^k$, we obtain $\mathbb P(\Omega_{\mathrm{RIP}}^c)\leq 2e^{-\frac{N}{64}}$. Moreover, $\mathbb P(\Omega_\xi^c)\le e^{-\frac{N}{2}}$, and hence $\mathbb P(\Omega)\ge 1-\mathbb P(\Omega_{\mathrm{esti}}^c)-3e^{-\frac{N}{64}}.$ In the following, we work on random event $\Omega$.

The property of $\Omega_{\mathrm{RIP}}$ that we use is as follows; its proof is routine and can be found in \cite[pp. 52 - 54]{chafaiInteractionsCompressedSensing2012}. On $\Omega_{\mathrm{RIP}}$, for every $I\subset J_\star^c$ with $|I|\le \lfloor\frac{N}{2^{9}\log(\frac{e|J_\star^c|}{N})}\rfloor$ and every $\boldsymbol v\in V_{J_\star^c}$,
\begin{align}
\|\Sigma_{J_\star^c}^{1/2}\boldsymbol v\|_2\le\frac{2}{\sqrt N}\|\mathbb X_{J_\star^c}\boldsymbol v\|_2+\frac{4}{\sqrt{\lfloor\frac{N}{2^{9}\log(\frac{e|J_\star^c|}{N})}\rfloor}}\|\Sigma_{I^c}^{1/2}\boldsymbol v_{I^c}\|_1.
\label{eq:lasso_rip_l1_l2}
\end{align}
Indeed, order $(\sqrt{\sigma_j}|v_j|)_{j\in I^c}$ decreasingly. Partition $I^c=\bigsqcup_{\ell\ge1}I_\ell$ into blocks of cardinality $\lfloor\frac{N}{2^{9}\log(\frac{e|J_\star^c|}{N})}\rfloor$. Then $\sum_{\ell\ge2}\|\Sigma_{I_\ell}^{1/2}\boldsymbol v_{I_\ell}\|_2\le\frac{\|\Sigma_{I^c}^{1/2}\boldsymbol v_{I^c}\|_1}{\sqrt{\lfloor\frac{N}{2^{9}\log(\frac{e|J_\star^c|}{N})}\rfloor}}$, while $\Omega_{\mathrm{RIP}}$ gives $\|\Sigma_{J_\star^c}^{1/2}\boldsymbol v\|_2\le\|\Sigma_{I\cup I_1}^{1/2}\boldsymbol v_{I\cup I_1}\|_2+\sum_{\ell\ge2}\|\Sigma_{I_\ell}^{1/2}\boldsymbol v_{I_\ell}\|_2 \le\frac{\sqrt3}{\sqrt N}\|\mathbb X_{J_\star^c}\boldsymbol v\|_2+(1+\sqrt6)\sum_{\ell\ge2}\|\Sigma_{I_\ell}^{1/2}\boldsymbol v_{I_\ell}\|_2,$ which proves \eqref{eq:lasso_rip_l1_l2}. 

To prove Proposition~\ref{prop:lambda_large_shift}, we set $\vv= \mathcal{A}[\boldsymbol{\xi}+\boldsymbol{\zeta}] - \mathcal{A}[\boldsymbol{\xi}]$. Thus we only need to construct an $I \subset J_\star^c$ with $|I|$ satisfying the RIP condition, in order to bound $\|\bX_{J_\star^c}(\mathcal{A}[\boldsymbol{\xi}+\boldsymbol{\zeta}] - \mathcal{A}[\boldsymbol{\xi}])\|_2$ and $\|(\mathcal{A}[\boldsymbol{\xi}+\boldsymbol{\zeta}] - \mathcal{A}[\boldsymbol{\xi}])_{I^c}\|_1$ from above.

\begin{enumerate}
    \item We first examine the upper bound of $\|\bX_{J_\star^c}(\mathcal{A}[\boldsymbol{\xi}+\boldsymbol{\zeta}] - \mathcal{A}[\boldsymbol{\xi}])\|_2$. This term is independent of the choice of $I$, and the tool for studying it is the KKT condition.  For every $\boldsymbol\mu\in\mathbb R^N$, the KKT condition is $\frac{\mathbb X_{J_\star^c}^\top(\boldsymbol\mu-\mathbb X_{J_\star^c}\mathcal A[\boldsymbol\mu])}{N\lambda}\in\partial^-\|\cdot\|_1(\mathcal A[\boldsymbol\mu]).$ Apply the KKT condition to $\boldsymbol\xi$ and $\boldsymbol\xi+\boldsymbol\zeta$. The subgradient inequality makes each of the two bracketed terms below nonnegative, and
    \begin{align}\label{eq:lasso_shift_basic}
    \begin{aligned}
        &0\le\lambda\bigg[\|\mathcal A[\boldsymbol\xi+\boldsymbol\zeta]\|_1-\|\mathcal A[\boldsymbol\xi]\|_1-\bigg\langle\frac{\mathbb X_{J_\star^c}^\top(\boldsymbol\xi-\mathbb X_{J_\star^c}\mathcal A[\boldsymbol\xi])}{N\lambda},\mathcal A[\boldsymbol\xi+\boldsymbol\zeta]-\mathcal A[\boldsymbol\xi]\bigg\rangle\bigg]\\
        &+\lambda\bigg[\|\mathcal A[\boldsymbol\xi]\|_1-\|\mathcal A[\boldsymbol\xi+\boldsymbol\zeta]\|_1-\bigg\langle\frac{\mathbb X_{J_\star^c}^\top(\boldsymbol\xi+\boldsymbol\zeta-\mathbb X_{J_\star^c}\mathcal A[\boldsymbol\xi+\boldsymbol\zeta])}{N\lambda},\mathcal A[\boldsymbol\xi]-\mathcal A[\boldsymbol\xi+\boldsymbol\zeta]\bigg\rangle\bigg]\\
        &=\frac{\langle\boldsymbol\zeta,\mathbb X_{J_\star^c}(\mathcal A[\boldsymbol\xi+\boldsymbol\zeta]-\mathcal A[\boldsymbol\xi])\rangle-\|\mathbb X_{J_\star^c}(\mathcal A[\boldsymbol\xi+\boldsymbol\zeta]-\mathcal A[\boldsymbol\xi])\|_2^2}{N}\\
        &\leq \frac{\|\boldsymbol\zeta\|_2\|\mathbb X_{J_\star^c}(\mathcal A[\boldsymbol\xi+\boldsymbol\zeta]-\mathcal A[\boldsymbol\xi])\|_2-\|\mathbb X_{J_\star^c}(\mathcal A[\boldsymbol\xi+\boldsymbol\zeta]-\mathcal A[\boldsymbol\xi])\|_2^2}{N} \leq \frac{\|\vzeta\|_2^2}{4N}\leq\frac{r_\star^2}{4},
    \end{aligned}
    \end{align}where we have use the elementary inequality $st-t^2\leq \frac{1}{4}s^2$ for any $s,t\geq 0$.
    In particular, $\|\mathbb X_{J_\star^c}(\mathcal A[\boldsymbol\xi+\boldsymbol\zeta]-\mathcal A[\boldsymbol\xi])\|_2^2\le\langle\boldsymbol\zeta,\mathbb X_{J_\star^c}(\mathcal A[\boldsymbol\xi+\boldsymbol\zeta]-\mathcal A[\boldsymbol\xi])\rangle$, and therefore
    \begin{align}
    \|\mathbb X_{J_\star^c}(\mathcal A[\boldsymbol\xi+\boldsymbol\zeta]-\mathcal A[\boldsymbol\xi])\|_2\le\|\boldsymbol\zeta\|_2\le\sqrt N r_\star.
    \label{eq:lasso_shift_prediction}
    \end{align}
    Before proceeding to the next item, we first record a stronger conclusion given by the KKT condition, which will be used in the next item. Fix the $j$-th coordinate with $j \in J_\star^c$, and let $\bX_{(j)}$ be the $j$-th column of $\bX$. From the KKT condition $\frac{\bX_{J_\star^c}^\top(\vmu - \bX_{J_\star^c}\cA[\vmu])}{N\lambda} \in \partial^-\|\cdot\|_1(\cA[\vmu])$ applied to $\vmu=\bxi$, we obtain the KKT condition on the $j$-th component: $\frac{1}{N\lambda}\langle \bX_{(j)}, \bxi - \bX_{J_\star^c}\cA[\bxi]\rangle\in\partial^-|\cdot|(\cA[\bxi]_j)$, and consequently,
    \begin{align}\label{eq:summand_positive}
        \forall j\in J_\star^c,\qquad \big| \cA[\bxi+\vzeta]_j \big| - \big| \cA[\bxi]_j \big| - \frac{1}{N\lambda}\big\langle \bX_{(j)}, \bxi - \bX_{J_\star^c}\cA[\bxi]\big\rangle \big( \cA[\bxi+\vzeta]_j - \cA[\bxi]_j \big)\geq 0.
    \end{align}

    \item We construct a coordinate subset $I$ satisfying the RIP property in order to examine the upper bound of $\|(\mathcal{A}[\boldsymbol{\xi}+\boldsymbol{\zeta}] - \mathcal{A}[\boldsymbol{\xi}])_{I^c}\|_1$. The coordinate subset we choose is
    \begin{align*}
        I = \bigg(j\in J_\star^c:\frac{|\langle\mathbb X_{(j)},\boldsymbol\xi-\mathbb X_{J_\star^c}\mathcal A[\boldsymbol\xi]\rangle|}{N\lambda}>\frac12\bigg).
    \end{align*}We first verify that $|I| < \big\lfloor\frac{N}{2^{9}\log(\frac{e|J_\star^c|}{N})}\big\rfloor$. By the definition of the solution of the convex optimization problem for $\mathcal{A}[\boldsymbol{\xi}]$, namely \eqref{eq:def_cA}, comparing the value of the objective function at the solution with that at $\mathbf{0}$, and then using $\boldsymbol{\xi} \in \Omega$, we obtain $\|\boldsymbol\xi-\mathbb X_{J_\star^c}\mathcal A[\boldsymbol\xi]\|_2\le\|\boldsymbol\xi\|_2\le2\sigma_\xi\sqrt N$. Now we assume the contrary, namely $|I| \ge \big\lfloor\frac{N}{2^{9}\log(\frac{e|J_\star^c|}{N})}\big\rfloor$. Then there exists $I' \subset I$ such that $|I'|=\lfloor\frac{N}{2^{9}\log(\frac{e|J_\star^c|}{N})}\rfloor$, and consequently, $\frac{N^2\lambda^2}{4}|I'|<\|\mathbb X_{I'}^\top(\boldsymbol\xi-\mathbb X_{J_\star^c}\mathcal A[\boldsymbol\xi])\|_2^2\le2N(\max_{j\in J_\star^c}\sigma_j)\|\boldsymbol\xi-\mathbb X_{J_\star^c}\mathcal A[\boldsymbol\xi]\|_2^2\le8N^2\sigma_\xi^2(\max_{j\in J_\star^c}\sigma_j)$, where we have used the definition on $I$, $\Omega_{\mathrm{RIP}}$ and the fact that $\|\boldsymbol\xi-\mathbb X_{J_\star^c}\mathcal A[\boldsymbol\xi]\|_2\le2\sigma_\xi\sqrt N$. That is, $\lambda^2|I'|<32\sigma_\xi^2(\max_{j\in J_\star^c}\sigma_j)$. On the other hand, by the assumption of Proposition~\ref{prop:lambda_large_shift}, $\log(\frac{e|J_\star^c|}{N})\le \frac{N}{2^{10}}$ gives $\lambda^2|I'|\ge\frac{\lambda^2N}{2^{10}\log(\frac{e|J_\star^c|}{N})}\ge64\sigma_\xi^2(\max_{j\in J_\star^c}\sigma_j)$, whereas the strict inequality just obtained gives $\lambda^2|I'|<32\sigma_\xi^2(\max_{j\in J_\star^c}\sigma_j)$, leading to a contradiction with $\lambda^2|I'|<32\sigma_\xi^2(\max_{j\in J_\star^c}\sigma_j)$. Therefore, $|I| < \big\lfloor\frac{N}{2^{9}\log(\frac{e|J_\star^c|}{N})}\big\rfloor$.

    Now we derive an upper bound for $\|(\mathcal{A}[\boldsymbol{\xi}+\boldsymbol{\zeta}] - \mathcal{A}[\boldsymbol{\xi}])_{I^c}\|_1$. We first make the following observation: for any $j \in I^c$, we must have $\cA[\bxi]_j = 0$; otherwise, $\partial^-|\cdot|(\cA[\bxi]_j) = \sign(\cA[\bxi]_j)$, whose absolute value is $1$, which cannot be no greater than $\frac{1}{2}$, hence contradicts with the KKT condition $\frac{1}{N\lambda}\langle \bX_{(j)}, \bxi - \bX_{J_\star^c}\cA[\bxi]\rangle\in\partial^-|\cdot|(\cA[\bxi]_j)$. Therefore by \eqref{eq:summand_positive} together with the definition of $j\in I^c$, there holds
    \begin{align*}
    &|\mathcal A[\boldsymbol\xi+\boldsymbol\zeta]_j| - |\cA[\bxi]|_j -\frac{\langle\mathbb X_{(j)},\boldsymbol\xi-\mathbb X_{J_\star^c}\mathcal A[\boldsymbol\xi]\rangle}{N\lambda}(\mathcal A[\boldsymbol\xi+\boldsymbol\zeta]_j-\cA[\bxi]_j)\ge\frac12\big|\mathcal A[\boldsymbol\xi+\boldsymbol\zeta]_j\big| = \frac12\big|\mathcal A[\boldsymbol\xi+\boldsymbol\zeta]_j-\mathcal A[\boldsymbol\xi]_j\big|.
    \end{align*}By \eqref{eq:lasso_shift_basic}, the fact that $\|\mathcal A[\boldsymbol\xi]\|_1-\|\mathcal A[\boldsymbol\xi+\boldsymbol\zeta]\|_1-\big\langle\frac{\mathbb X_{J_\star^c}^\top(\boldsymbol\xi+\boldsymbol\zeta-\mathbb X_{J_\star^c}\mathcal A[\boldsymbol\xi+\boldsymbol\zeta])}{N\lambda},\mathcal A[\boldsymbol\xi]-\mathcal A[\boldsymbol\xi+\boldsymbol\zeta]\big\rangle\geq 0$ together with \eqref{eq:summand_positive}, we obtain that
    \begin{align*}
        &\frac{r_\star^2}{4\lambda}\geq \|\mathcal A[\boldsymbol\xi+\boldsymbol\zeta]\|_1-\|\mathcal A[\boldsymbol\xi]\|_1-\bigg\langle\frac{\mathbb X_{J_\star^c}^\top(\boldsymbol\xi-\mathbb X_{J_\star^c}\mathcal A[\boldsymbol\xi])}{N\lambda},\mathcal A[\boldsymbol\xi+\boldsymbol\zeta]-\mathcal A[\boldsymbol\xi]\bigg\rangle\geq \frac{1}{2}\sum_{j\in I^c}\big|\mathcal A[\boldsymbol\xi+\boldsymbol\zeta]_j-\mathcal A[\boldsymbol\xi]_j\big|.
    \end{align*}
    That is,
    \begin{align}
    \bigg\|\bigl(\mathcal A[\boldsymbol\xi+\boldsymbol\zeta]-\mathcal A[\boldsymbol\xi]\bigr)_{j\in I^c}\bigg\|_1\le\frac{r_\star^2}{2\lambda}.
    \label{eq:lasso_shift_l1}
    \end{align}
\end{enumerate}
Using $\|\Sigma_{I^c}^{1/2}\boldsymbol v_{I^c}\|_1\leq\sqrt{\max_{j\in J_\star^c}\sigma_j}\|\boldsymbol v_{I^c}\|_1$, applying \eqref{eq:lasso_rip_l1_l2} with  $\boldsymbol v=\mathcal A[\boldsymbol\xi+\boldsymbol\zeta]-\mathcal A[\boldsymbol\xi]$, and then using \eqref{eq:lasso_shift_prediction} and \eqref{eq:lasso_shift_l1}, yields
\begin{align*}
&\|\Sigma_{J_\star^c}^{1/2}(\mathcal A[\boldsymbol\xi+\boldsymbol\zeta]-\mathcal A[\boldsymbol\xi])\|_2\le2r_\star+\frac{2r_\star^2\sqrt{\max_{j\in J_\star^c}\sigma_j}}{\lambda\sqrt{\lfloor\frac{N}{2^{9}\log(\frac{e|J_\star^c|}{N})}\rfloor}}\le2r_\star+2^{6}\frac{r_\star^2}{\lambda}\sqrt{\frac{(\max_{j\in J_\star^c}\sigma_j)\log(\frac{e|J_\star^c|}{N})}{N}}.
\end{align*}
The assumptions give $2^8\frac{r_\star^3}{\lambda}\sqrt{\frac{(\max_{j\in J_\star^c}\sigma_j)\log(\frac{e|J_\star^c|}{N})}{N}}\leq r_\star^2$, so the cross term in the square of the preceding bound is at most $r_\star^2$. Consequently,
\begin{align*}
\|\Sigma_{J_\star^c}^{1/2}(\mathcal A[\boldsymbol\xi+\boldsymbol\zeta]-\mathcal A[\boldsymbol\xi])\|_2^2\le5r_\star^2+2^{12}\frac{(\max_{j\in J_\star^c}\sigma_j)\log(\frac{e|J_\star^c|}{N})}{N\lambda^2}r_\star^4\le6r_\star^2.
\end{align*}
The last inequality uses $r_\star\le\sigma_\xi$ and $\lambda^2\ge\frac{2^{16}\sigma_\xi^2(\max_{j\in J_\star^c}\sigma_j)\log(\frac{e|J_\star^c|}{N})}{N}$.
\end{proof}


\subsubsection{The regime \texorpdfstring{$\frac{\sigma_\xi}{\sqrt N}\leq\lambda\leq 2^{-9}\sigma_\xi\sqrt{\frac{\log(\frac{e|J_\star^c|}{N})}{N}}$}{of small lambda}.}

\begin{Proposition}\label{prop:lambda_is_small_shift}
Grant Assumption~\ref{assumption:LASSO_upper}. Suppose $|J_\star^c|\geq2^{12}N$, $\log(\frac{e|J_\star^c|}{N})^2\leq2^{-8}N$, $|J_\star|\log(\frac{e|J_\star^c|}{N})\leq2^{-8}N$, and
\begin{align}
\frac{\sigma_\xi}{\sqrt N}\leq\lambda\le 2^{-9}\sigma_\xi\sqrt{\frac{\log(\frac{e|J_\star^c|}{N})}N}.
\label{eq:lambda-regime}
\end{align}We also assume $\|\bbeta_{J_\star^c}^\star\|_2^2\leq r_\star^2$. Then there is an event $\Omega$ such that
\begin{align*}
&\mathbb P(\Omega)\geq1-\mathbb P(\Omega_{J_\star}^c)-2^{9}\exp\bigg(-\frac{2^{-16}N}{\log(\frac{e|J_\star^c|}{N})^2}\bigg)-2^{9}\exp(-2^{-16}|J_\star|),
\end{align*}
and, on $\Omega$, if $\boldsymbol\zeta=\mathbb X_{J_\star}(\boldsymbol\beta_{J_\star}^\star-\widehat{\boldsymbol\beta}_{J_\star})+\mathbb X_{J_\star^c}\boldsymbol\beta_{J_\star^c}^\star$, then
\begin{align}
&\|\mathcal A[\boldsymbol\xi+\boldsymbol\zeta]\|_2^2=\|\widehat{\boldsymbol\beta}_{J_\star^c}[\boldsymbol y-\mathbb X_{J_\star}\widehat{\boldsymbol\beta}_{J_\star}]\|_2^2\leq R^2,\mbox{ where }R^2=2^{21}\bigg(\frac{\sigma_\xi^2}{\log(\frac{e|J_\star^c|}{N})}+\|\bbeta_{J_\star^c}^\star\|_2^2\bigg).
\label{eq:A-shift-bound}
\end{align}
\end{Proposition}
Starting from this proof, we will gradually begin to use CGMT as a proof tool (in the proof of Proposition~\ref{prop:lambda_is_small_shift}, Proposition~\ref{prop:null_signal_upper_bound}, and Proposition~\ref{prop:free-variance-lower-bound}). We have included in Section~\ref{sec:cgmt} some useful background knowledge on using CGMT as a proof method.

\begin{Lemma}\label{lem:soft-threshold}
Grant the assumptions of Proposition~\ref{prop:lambda_is_small_shift}. Let $\Phi_G(t)=\mathbb P\{G\geq t\}$ for $G\sim\mathcal N(0,1)$, let $z_{16N}=\Phi_G^{-1}(\frac{8N}{|J_\star^c|})$ and $\underline{z}_{16N}=z_{16N}-4/z_{16N}$, and let $(h_j)_{j\in J_\star^c}$ be independent standard Gaussian random variables. There is an event $\Omega_{\mathrm{soft}}^{\mathrm{aux}}$ such that
\begin{align*}
&\mathbb P(\Omega_{\mathrm{soft}}^{\mathrm{aux}})\geq1-5\exp\bigg(-\frac{N}{2\log(\frac{e|J_\star^c|}{N})}\bigg),
\end{align*}
on which the following conclusions hold simultaneously for every $\gamma\in[\frac{1}{\sqrt2},\sqrt{\frac{3}{2}}]$, $s\geq\sigma_\xi$, and $\tau\geq0$.
\begin{enumerate}
    \item 
    \begin{align}
    &\mbox{if }\quad \frac1N\sum_{j\in J_\star^c}h_j\bigg(S_\tau\bigg(h_j+\frac{\gamma\sqrt N}{s}\beta_j^\star\bigg)-\frac{\gamma\sqrt N}{s}\beta_j^\star\bigg)\leq2\quad\mbox{ then }\quad \tau\geq \underline{z}_{16N}\geq\frac12\sqrt{\log(\frac{e|J_\star^c|}{N})},
    \label{eq:soft-threshold-first-moment}
    \end{align}
    and
    \begin{align}
    &\frac1N\sum_{j\in J_\star^c}S_\tau\bigg(h_j+\frac{\gamma\sqrt N}{s}\beta_j^\star\bigg)^2\leq\frac{2^{14}}{\log(\frac{e|J_\star^c|}{N})}+\frac{2\gamma^2}{s^2}\|\boldsymbol\beta_{J_\star^c}^\star\|_2^2.
    \label{eq:soft-threshold-second-moment}
    \end{align}
    \item If $\tau\geq4\sqrt{\log(\frac{e|J_\star^c|}{N})}$, then $\frac1N\sum_{j\in J_\star^c}h_j\big(S_\tau\big(h_j+\frac{\gamma\sqrt N}{s}\beta_j^\star\big)-\frac{\gamma\sqrt N}{s}\beta_j^\star\big)\leq \frac{1}{4}$.
\end{enumerate}
\end{Lemma}

\begin{proof}
Let $\Omega_{\mathrm{soft}}^{\mathrm{aux}}=\Omega_{\mathrm{lin}}^{\mathrm{aux}}\cap\Omega_{\mathrm{count}}^{\mathrm{aux}}\cap\Omega_{\mathrm{energy}}^{\mathrm{aux}}\cap\Omega_{\mathrm{large}}^{\mathrm{aux}}$, where
\begin{align*}
&\Omega_{\mathrm{lin}}^{\mathrm{aux}}=\big\{|\langle\boldsymbol\beta_{J_\star^c}^\star,\boldsymbol h_{J_\star^c}\rangle|\leq\frac{\sigma_\xi\sqrt N}{8^{\frac{3}{2}}}\big\},\quad \Omega_{\mathrm{count}}^{\mathrm{aux}}=\bigg\{\sum_{\substack{j\in J_\star^c:\,\sqrt{\frac32}\sqrt N\frac{|\beta_j^\star|}{\sigma_\xi}\leq\frac{1}{z_{16N}}}}\1_{\{|h_j|>z_{16N}-\frac{2}{z_{16N}}\}}\geq8N\bigg\},\\
&\Omega_{\mathrm{energy}}^{\mathrm{aux}}=\bigg\{\sum_{j\in J_\star^c}(|h_j|-\underline{z}_{16N})_+^2\leq\frac{2^{13}N}{\log(\frac{e|J_\star^c|}{N})}\bigg\},\mbox{ and } \Omega_{\mathrm{large}}^{\mathrm{aux}}=\bigg\{\sum_{j\in J_\star^c}h_j^2\boldsymbol 1_{\{|h_j|>(4-2^{-8})\sqrt{\log(\frac{e|J_\star^c|}{N})}\}}\leq\frac N8\bigg\}.
\end{align*}We first prove that the random event $\Omega_{\mathrm{soft}}^{\mathrm{aux}}$ occurs with high probability. 
\begin{enumerate}
    \item It is easy to see that $\Omega_{\mathrm{lin}}^{\mathrm{aux}}$ occurs with probability at least $1 - 2\exp(-\frac{N}{1024})$, $\langle \bbeta_{J_\star^c}^\star, \vh_{J_\star^c} \rangle$ is a Gaussian random variable.

    \item $\Omega_{\mathrm{energy}}^{\mathrm{aux}}$ is also easy. Indeed, $F:\vh_{J_\star^c}\mapsto (\sum_{j\in J_\star^c}(|h_j|-\underline{z}_{16N})_+^2)^{\frac{1}{2}}$ is $1$-Lipschitz. Applying \eqref{eq:expectation_soft_max} with $a = \underline{z}_{16N}$ and \eqref{eq:upper_Phi_G_z_4z} with $z = z_{16N}$, together with the definition $z_{16N} = \Phi_G^{-1}(\frac{8N}{|J_\star^c|})$, we have $|J_\star^c| \mathbb{E} (|h_j| - \underline{z}_{16N})_+^2 \le |J_\star^c|\frac{4}{(\underline{z}_{16N})^2}e^4\Phi_G(z_{16N}) = \frac{32e^4N}{(\underline{z}_{16N})^2}$. The two bounds on $\lambda$ imply $\log(e|J_\star^c|/N)\geq2^{18}$. By \eqref{eq:estimate_z_q}, $z_{16N}^2\geq\log(e|J_\star^c|/N)-\log16$, so $z_{16N}^2\geq8$ and $(\underline{z}_{16N})^2\geq\frac12\log(e|J_\star^c|/N)$. Consequently, $|J_\star^c|\mathbb E(|h_j|-\underline{z}_{16N})_+^2\leq64e^4N/\log(e|J_\star^c|/N)$. The probability estimate of $\Omega_{\mathrm{energy}}^{\mathrm{aux}}$ follows from Gaussian Lipschitz concentration inequality.

    \item The probabilistic estimate of $\Omega_{\mathrm{large}}^{\mathrm{aux}}$ follows from standard Chernoff argument, see, for instance, \cite[Chapter 2]{vershyninHighDimensionalProbabilityIntroduction2018}. Indeed, $\bP(\big(\Omega_{\mathrm{large}}^{\mathrm{aux}}\big)^c)\leq \exp(-\frac{1}{64}N)\prod_{j\in J_\star^c}\bE\exp(\frac{1}{8}h_j^2\1_{\{|h_j|>a\}})$, where $a=(4-2^{-8})\sqrt{\log(\frac{e|J_\star^c|}{N})}$. Notice that $(h_j)_{j\in J_\star^c}$ are independent standard Gaussian variables. Let $g\sim\cN(0,1)$, then $\bE[\exp(\frac{1}{8}g^2)\1_{\{|g|>a\}}] = \frac{2}{\sqrt{2\pi}}\int_a^\infty \exp(-\frac{3}{8}u^2)\mathrm{d}u\leq \frac{8}{3a\sqrt{2\pi}}\exp(-\frac{3}{8}a^2)$. Together with $|J_\star^c|=N\exp(\log(e|\frac{J_\star^c}{N}|)-1)$ and $a^2 = (4-2^{-8})^2\log(\frac{e|J_\star^c|}{N})$, we obtain that $\bP(\big(\Omega_{\mathrm{large}}^{\mathrm{aux}}\big)^c)\leq \exp(-\frac{N}{128})$.

    \item Lastly we deal with $\Omega_{\mathrm{count}}^{\mathrm{aux}}$. We first show that there are at most (up to multiplicative constant) $N$ coordinates of $J_\star^c$ such that $\sqrt{\frac32}\sqrt N\frac{|\beta_j^\star|}{\sigma_\xi}>\frac{1}{z_{16N}}$. Indeed, denote by $B=(j\in J_\star^c: \sqrt{\frac32}\sqrt N\frac{|\beta_j^\star|}{\sigma_\xi}>\frac{1}{z_{16N}})$, then $\|\bbeta_{J_\star^c}^\star\|_2^2 \geq \sum_{j\in B}(\beta_j^\star)^2 > |B|\frac{2\sigma_\xi^2}{3Nz_{16N}^2}$. From the assumption, $\|\bbeta_{J_\star^c}^\star\|_2^2\leq r_\star^2\lesssim \sigma_\xi^2$, hence $|B|\lesssim N\log(e|J_\star^c|/N)\leq \frac{1}{2}|J_\star^c|$. By assumption of the proposition, $|B^c|\geq \frac{|J_\star^c|}{2}$. We only need to study $\bP\{|h_j|>z_{16N}-\frac{2}{z_{16N}}\}=2\Phi_G(z_{16N}-\frac{2}{z_{16N}})$. By \eqref{eq:lower_Phi_G_z_Jc}, $2\Phi_G(z_{16N}-\frac{2}{z_{16N}})\geq \frac{64N}{|J_\star^c|}$. Therefore, the binomial estimate yields that with probability at least $1-\exp(-N)$, $\Omega_{\mathrm{count}}^{\mathrm{aux}}$ happens.
\end{enumerate}
Recall that
\begin{align}\tag{\ref{eq:softmax_calculus}}
&\forall h,b\in\bR,\,\forall \tau\geq0,\quad \big(S_\tau(h+b)-b\big)h\geq|h|(|h|-\tau-|b|)_+-b^2-bh,\mbox{ and } |S_\tau(h+b)|^2\leq2b^2+2(|h|-\tau)_+^2.
\end{align}
On $\Omega_{\mathrm{count}}^{\mathrm{aux}}\cap\Omega_{\mathrm{lin}}^{\mathrm{aux}}$, if the condition of \eqref{eq:soft-threshold-first-moment} holds and $\tau<\underline{z}_{16N}$, then the first inequality of \eqref{eq:softmax_calculus} with $b = \frac{\gamma\sqrt{N}}{s}\beta_j^\star$, $h=h_j$ together with the assumption that $\gamma^2\leq \frac{3}{2}$, $\tau<\underline{z}_{16N}$ and $s\geq \sigma_\xi$ yields
\begin{align*}
&\frac1N\sum_{j\in J_\star^c}h_j\bigg(S_\tau\bigg(h_j+\frac{\gamma\sqrt N}{s}\beta_j^\star\bigg)-\frac{\gamma\sqrt N}{s}\beta_j^\star\bigg)\\
&\geq \frac{1}{N}\sum_{j\in J_\star^c}|h_j|\bigg(|h_j| - \tau - \frac{\gamma\sqrt{N}}{s}|\beta_j^\star|\bigg)_+ - \frac{\gamma^2}{s^2}\|\bbeta_{J_\star^c}^\star\|_2^2 - \frac{\gamma}{s\sqrt{N}}\langle\vh_{J_\star^c},\bbeta_{J_\star^c}^\star\rangle\\
&\geq \frac{8N}{N} \big(z_{16N} - \frac{2}{z_{16N}}\big) \bigg( z_{16N} - \frac{2}{z_{16N}} - \big( \underline{z}_{16N}\big) - \frac{1}{z_{16N}} \bigg)_+ - \frac{3}{2} - \frac{1}{8} >4-\frac32-\frac18>2,
\end{align*}leading to a contradiction. This proves the right-hand-side in \eqref{eq:soft-threshold-first-moment}; the second follows from the second inequality of \eqref{eq:softmax_calculus}. On $\Omega_{\mathrm{energy}}^{\mathrm{aux}}$,
\begin{align*}
&\frac1N\sum_{j\in J_\star^c}S_\tau\bigg(h_j+\frac{\gamma\sqrt N}{s}\beta_j^\star\bigg)^2\leq\frac{2\gamma^2}{s^2}\|\boldsymbol\beta_{J_\star^c}^\star\|_2^2+\frac2N\sum_{j\in J_\star^c}(|h_j|-\underline{z}_{16N})_+^2\leq\frac{2\gamma^2}{s^2}\|\boldsymbol\beta_{J_\star^c}^\star\|_2^2+\frac{2^{14}}{\log(\frac{e|J_\star^c|}{N})},
\end{align*}
which is \eqref{eq:soft-threshold-second-moment}. Finally, \eqref{eq:lambda-regime} and $\|\boldsymbol\beta_{J_\star^c}^\star\|_\infty\leq\lambda$ from the definition of $J_\star^c$ imply
\begin{align*}
&\frac{\gamma\sqrt N}{s}|\beta_j^\star|<2^{-8}\sqrt{\log(\frac{e|J_\star^c|}{N})}\qquad\forall j\in J_\star^c.
\end{align*}
If $\tau\geq4\sqrt{\log(\frac{e|J_\star^c|}{N})}$, a nonzero soft-thresholded coordinate therefore satisfies $|h_j|>(4-2^{-8})\sqrt{\log(\frac{e|J_\star^c|}{N})}$. On $\Omega_{\mathrm{large}}^{\mathrm{aux}}\cap\Omega_{\mathrm{lin}}^{\mathrm{aux}}$, the left-hand side of \eqref{eq:soft-threshold-first-moment} is at most $\frac{1}{8}+\frac{1}{8}=\frac{1}{4}$.
\end{proof}

\begin{proof}[Proof of Proposition~\ref{prop:lambda_is_small_shift}]
Let
\begin{align*}
\Omega_{J_\star} = \{\hat\bbeta_{J_\star}\in T\},\mbox{ where } T=B_2^{J_\star}(\boldsymbol\beta_{J_\star}^\star;r_\star)\cap B_1^{J_\star}(\boldsymbol\beta_{J_\star}^\star;\rho_\star).
\end{align*}
To avoid the issues arising from the dependence between $\hat\bbeta_{J_\star}$ and $\hat\bbeta_{J_\star^c}$, we prove a stronger statement: there exists a high-probability random event on which, for every $\bbeta_{J_\star} \in T$, we have $\|\hat\bbeta_{J_\star^c}[\vy - \bX_{J_\star} \bbeta_{J_\star}]\|_2^2 \le R^2$. Since the high-probability event $\Omega_{J_\star}$ ensures that $\hat\bbeta_{J_\star} \in T$ on it, on the intersection of the two events we have $\|\hat\bbeta_{J_\star^c}[\vy - \bX_{J_\star} \hat\bbeta_{J_\star}]\|_2^2 \le R^2$, which is precisely \eqref{eq:A-shift-bound} that we aim to prove.

We first define the high-probability random event under $\bP^{\otimes N}$. Let $\Omega_{\mathrm{PO}}=\Omega_\xi\cap\Omega_{\mathrm{op}}\cap\Omega_{\mathrm{tail}}$, where
\begin{align*}
&\Omega_\xi=\{\|\boldsymbol\xi\|_2\leq2\sigma_\xi\sqrt N\},\, \Omega_{\mathrm{op}}=\{\|\mathbb X_{J_\star}\|_{\mathrm{op}}\leq2\sqrt N\},\mbox{ and } \Omega_{\mathrm{tail}}=\{\|\mathbb X_{J_\star^c}\boldsymbol\beta_{J_\star^c}^\star\|_2\leq2\sqrt N\|\boldsymbol\beta_{J_\star^c}^\star\|_2\}.
\end{align*}
Next we define a high-probability random event under $\mathbb{P}_{\mathrm{aux}}$. Let $\boldsymbol g\in\mathbb R^N$ and $\boldsymbol h=(\boldsymbol h_{J_\star},\boldsymbol h_{J_\star^c},h_0)\in V_{J_\star}\times V_{J_\star^c}\times\bR$ be independent standard Gaussian vectors, write $\boldsymbol\omega=(\boldsymbol g,\boldsymbol h_{J_\star},\boldsymbol h_{J_\star^c},h_0)$, and set
\begin{align*}
&\Omega_{\mathrm{AO}}^{\mathrm{aux}}=\Omega_{\mathrm{soft}}^{\mathrm{aux}}\cap\bigg\{\frac12\leq\frac{\|\boldsymbol g\|_2^2}N\leq\frac32\bigg\}\cap\{\|\boldsymbol h_{J_\star}\|_2\leq2\sqrt{|J_\star|}\}\cap\{|h_0|\leq\frac{\sqrt N}{64}\},
\end{align*}where we recall that $\Omega_{\mathrm{soft}}^{\mathrm{aux}}$ is defined in Lemma~\ref{lem:soft-threshold}.
In the following, we work on $\Omega_{\mathrm{AO}}^{\mathrm{aux}}$ and prove that with $\bP^{\otimes N}$ high probability,
\begin{align*}
    \sup\bigg(\|\cA[\vy - \bX_{J_\star} \bbeta_{J_\star}]\|_2^2:\, \bbeta_{J_\star}\in T\bigg) \le R^2.
\end{align*}
Fix $\bbeta_{J_\star}\in T$. To apply CGMT, we first perform a compactification argument, which is routine; see Section~\ref{sec:cgmt}. We first make two observations.
\begin{enumerate}
    \item By the definition of $\cA[\vy-\bX_{J_\star}\bbeta_{J_\star}]$ from \eqref{eq:def_cA}, there holds $\frac{1}{2N\lambda}\|\vy-\bX_{J_\star}\bbeta_{J_\star} - \bX_{J_\star^c}\cA[\vy-\bX_{J_\star}\bbeta_{J_\star}]\|_2^2 + \|\cA[\vy-\bX_{J_\star}\bbeta_{J_\star}]\|_1\leq \frac{1}{2N\lambda}\|\vy-\bX_{J_\star}\bbeta_{J_\star}-\bX_{J_\star^c}\vzero\|_2^2 + \|\vzero\|_1$. Together with $\vy=\bX_{J_\star}\bbeta_{J_\star}^\star+\bX_{J_\star^c}\bbeta_{J_\star^c}^\star+\bxi$, there holds $\frac{1}{2N\lambda}\|\boldsymbol\xi+\mathbb X_{J_\star}(\boldsymbol\beta_{J_\star}^\star-\boldsymbol\beta_{J_\star})+\mathbb X_{J_\star^c}(\boldsymbol\beta_{J_\star^c}^\star-\widehat{\boldsymbol\beta}_{J_\star^c}[\boldsymbol y-\mathbb X_{J_\star}\boldsymbol\beta_{J_\star}])\|_2^2+\|\widehat{\boldsymbol\beta}_{J_\star^c}[\boldsymbol y-\mathbb X_{J_\star}\boldsymbol\beta_{J_\star}]\|_1\leq \frac{1}{2N\lambda}\|\bxi+\bX_{J_\star}(\bbeta_{J_\star}^\star - \bbeta_{J_\star}) + \bX_{J_\star^c}\bbeta_{J_\star^c}^\star\|_2^2.$ On $\Omega_{\mathrm{PO}}$, $\forall\bbeta_{J_\star}\in T$, $\|\boldsymbol\xi+\mathbb X_{J_\star}(\boldsymbol\beta_{J_\star}^\star-\boldsymbol\beta_{J_\star})+\mathbb X_{J_\star^c}\boldsymbol\beta_{J_\star^c}^\star\|_2\leq6\sigma_\xi\sqrt N.$ Therefore,
    \begin{align*}
        \mbox{ On } \Omega_{\mathrm{PO}},\quad \forall \bbeta_{J_\star}\in T,\qquad \norm{\hat\bbeta_{J_\star^c}[\vy-\bX_{J_\star}\bbeta_{J_\star}]}_2 \leq 18\sigma_\xi\sqrt N.
    \end{align*}
    
    \item From the last item, on $\Omega_{\mathrm{PO}}$, for every $\bbeta_{J_\star}\in T$, $\| \vy-\bX_{J_\star}\bbeta_{J_\star} -\bX_{J_\star^c}\cA[\vy-\bX_{J_\star}\bbeta_{J_\star}]\|_2 \leq \|\boldsymbol\xi+\mathbb X_{J_\star}(\boldsymbol\beta_{J_\star}^\star-\boldsymbol\beta_{J_\star})+\mathbb X_{J_\star^c}\boldsymbol\beta_{J_\star^c}^\star\|_2\leq6\sigma_\xi\sqrt N$. From \eqref{eq:variational_representation_ell_2_square}, $\frac{1}{2N\lambda}\|\vy-\bX_{J_\star}\bbeta_{J_\star}-\bX_{J_\star^c}\cA[\vy-\bX_{J_\star}\bbeta_{J_\star}]\|_2^2 = \max_{\vu\in\bR^N}\big( \langle\vu,\vy-\bX_{J_\star}\bbeta_{J_\star}-\bX_{J_\star^c}\cA[\vy-\bX_{J_\star}\bbeta_{J_\star}]\rangle - \frac{N\lambda}{2}\|\vu\|_2^2 \big)$, whose maximizer is given by $\vu^\circ = \frac{1}{N\lambda}(\vy-\bX_{J_\star}\bbeta_{J_\star}-\bX_{J_\star^c}\cA[\vy-\bX_{J_\star}\bbeta_{J_\star}]).$ Therefore,
    \begin{align*}
        \mbox{ On } \Omega_{\mathrm{PO}},\quad \forall \bbeta_{J_\star}\in T,\qquad \|\vu^\circ\|_2\leq \frac{6\sigma_\xi}{\lambda\sqrt{N}},\mbox{ where } \vu^\circ = \frac{1}{N\lambda}(\vy-\bX_{J_\star}\bbeta_{J_\star} - \bX_{J_\star^c}\cA[\vy-\bX_{J_\star}\bbeta_{J_\star}]).
    \end{align*}
    
\end{enumerate}Since $\hat\bbeta_{J_\star} \in T$ holds on $\Omega_{J_\star}$, combining this with \eqref{eq:variational_representation_LASSO}, we have the following fact on $\Omega_{J_\star}\cap\Omega_{\mathrm{PO}}$, that is, $\hat\bbeta = (\hat\bbeta_{J_\star}, \hat\bbeta_{J_\star^c}[\vy-\bX_{J_\star}\hat\bbeta_{J_\star}])\in\argmin(\Phi(\bbeta_{J_\star},\bbeta_{J_\star^c}):\bbeta_{J_\star}\in T,\,\|\bbeta_{J_\star^c}\|_2\leq 32\sigma_\xi\sqrt{N})$, where $\Phi(\bbeta_{J_\star},\bbeta_{J_\star^c})=\Phi$ and
\begin{align}
&\Phi=\min_{\substack{\boldsymbol\beta_{J_\star}\in T\\\|\boldsymbol\beta_{J_\star^c}\|_2\leq32\sigma_\xi\sqrt N}}\max_{\|\boldsymbol u\|_2\leq\frac{8\sigma_\xi}{\lambda\sqrt N}}\bigg\{ \langle\boldsymbol u,\boldsymbol\xi+\mathbb X_{J_\star}(\boldsymbol\beta_{J_\star}^\star-\boldsymbol\beta_{J_\star})+\mathbb X_{J_\star^c}(\boldsymbol\beta_{J_\star^c}^\star-\boldsymbol\beta_{J_\star^c})\rangle-\frac{N\lambda}{2}\|\boldsymbol u\|_2^2+\|\boldsymbol\beta_{J_\star}\|_1+\|\boldsymbol\beta_{J_\star^c}\|_1 \bigg\}.
\label{eq:primary-optimization}
\end{align}
Now, let us consider a restricted version of $\Phi$, that is,
\begin{align*}
&\Phi_R=\min_{\substack{\boldsymbol\beta_{J_\star}\in T\\ \|\bbeta_{J_\star^c}\|_2\geq R \\\|\boldsymbol\beta_{J_\star^c}\|_2\leq32\sigma_\xi\sqrt N}}\max_{\|\boldsymbol u\|_2\leq\frac{8\sigma_\xi}{\lambda\sqrt N}}\bigg\{ \langle\boldsymbol u,\boldsymbol\xi+\mathbb X_{J_\star}(\boldsymbol\beta_{J_\star}^\star-\boldsymbol\beta_{J_\star})+\mathbb X_{J_\star^c}(\boldsymbol\beta_{J_\star^c}^\star-\boldsymbol\beta_{J_\star^c})\rangle-\frac{N\lambda}{2}\|\boldsymbol u\|_2^2+\|\boldsymbol\beta_{J_\star}\|_1+\|\boldsymbol\beta_{J_\star^c}\|_1 \bigg\}.
\end{align*}Consider the random event $\{\|\hat\bbeta_{J_\star^c}[\vy-\bX_{J_\star}\hat\bbeta_{J_\star}]\|_2\geq R\}$, then for any $\vu$ such that $\|\vu\|_2\leq \frac{8\sigma_\xi}{\lambda\sqrt{N}}$, the triple $(\hat\bbeta_{J_\star},\hat\bbeta_{J_\star^c}[\vy-\bX_{J_\star}\hat\bbeta_{J_\star}],\vu)$ is feasible for $\Phi_R$, therefore, we necessarily have $\Phi_R=\Phi$. As a result,
\begin{align*}
    \big( \|\hat\bbeta_{J_\star^c}[\vy-\bX_{J_\star}\hat\bbeta_{J_\star}]\|_2 > R \big) \subset \big\{ \Phi_R = \Phi \big\},
\end{align*}and we only need to estimate the upper bound for $\bP\{\Phi_R=\Phi\}$. This is a comparison problem of convex-concave optimization problems on two compact sets; therefore, we will use CGMT, namely Lemma~\ref{lem:common-cgmt}, together with the comparison lemma~\ref{lem:radial-restriction-gap} in the Appendix. To apply Lemma~\ref{lem:radial-restriction-gap}, we first use CGMT, namely Lemma~\ref{lem:common-cgmt}, to reduce $\Phi$ to an auxiliary problem.

We define the map $(\bbeta_{J_\star}, \bbeta_{J_\star^c}) \in T \times B_2(\vzero;32\sigma_\xi\sqrt{N}) \mapsto \vw := (\bbeta_{J_\star}, \bbeta_{J_\star^c}, -\sigma_\xi) \in \mathbb{R}^{p+1}$, and denote its image by $S_\vw$; define $S_\vu = B_2(\vzero; \frac{8\sigma_\xi}{\sqrt N})$. Then by our definitions of $\vw$, $S_\vw$, and $S_\vu$ above, $\Phi$ is equivalent to $\min_{\vw \in S_\vw} \max_{\vu \in S_\vu} \{\langle \vu, \bG \vw \rangle - \frac{N\lambda}{2}\|\vu\|_2^2 + \|\bbeta_{J_\star}\|_1 + \|\bbeta_{J_\star^c}\|_1\}$, where $\bG = [\bX \,|\, -\frac{1}{\sigma_\xi}\bxi] \in \mathbb{R}^{N \times (p+1)}$ is a standard Gaussian random matrix. By Lemma~\ref{lem:common-cgmt}, together with the one-dimensional reduction trick in Section~\ref{sec:dimension_reduction_trick}, we obtain the auxiliary problem $\phi$ of $\Phi$. To formulate $\phi$, we first make some notation. For $\boldsymbol\beta_{J_\star}\in T$, $\boldsymbol\beta_{J_\star^c}\in\mathbb R^{J_\star^c}$, and any $\vu\in S_\vu$, write
\begin{align*}
&s_{\boldsymbol\beta_{J_\star}}(\boldsymbol\beta_{J_\star^c})=\big( \sigma_\xi^2+\|\boldsymbol\beta_{J_\star}-\boldsymbol\beta_{J_\star}^\star\|_2^2+\|\boldsymbol\beta_{J_\star^c}-\boldsymbol\beta_{J_\star^c}^\star\|_2^2\big)^{\frac{1}{2}},\mbox{ and } a = \lambda\sqrt{N}\|\vu\|_2.
\end{align*}By splitting the minimization over $\vw \in S_\vw$ into a two-step minimization over $\bbeta_{J_\star} \in T$ and $\bbeta_{J_\star^c} \in B_2(\vzero, 32\sigma_\xi\sqrt{N})$, we obtain
\begin{align}
&\phi=\min_{\boldsymbol\beta_{J_\star}\in T}\phi(\boldsymbol\beta_{J_\star}),\mbox{ where } \phi(\boldsymbol\beta_{J_\star})=\min_{\|\boldsymbol\beta_{J_\star^c}\|_2\leq32\sigma_\xi\sqrt N}\max_{0\leq a\leq8\sigma_\xi}\mathcal F_{\boldsymbol\beta_{J_\star}}(\boldsymbol\beta_{J_\star^c},a),
\label{eq:auxiliary-optimization}
\end{align}
and
\begin{align*}
&\mathcal F_{\boldsymbol\beta_{J_\star}}(\boldsymbol\beta_{J_\star^c},a)=\frac{a\|\boldsymbol g\|_2}{\lambda\sqrt N}s_{\boldsymbol\beta_{J_\star}}(\boldsymbol\beta_{J_\star^c})-\frac{a}{\lambda\sqrt N}\langle\boldsymbol h,\boldsymbol\beta-\boldsymbol\beta^\star\rangle-\frac{a\sigma_\xi}{\lambda\sqrt N}h_0-\frac{a^2}{2\lambda}+\|\boldsymbol\beta_{J_\star}\|_1+\|\boldsymbol\beta_{J_\star^c}\|_1.
\end{align*}Similarly, we define $\phi_R$ be the auxiliary problem of $\Phi_R$, and $\phi_R = \min_{\bbeta_{J_\star}\in T}\phi_R(\bbeta_{J_\star})$, where $\phi_R(\bbeta_{J_\star})$ is the same as $\phi(\bbeta_{J_\star})$, but with the feasible set replaced by $(R\leq \|\boldsymbol\beta_{J_\star^c}\|_2\leq32\sigma_\xi\sqrt N)$. It is trivial that for any $\bbeta_{J_\star}\in T$, $\phi_R(\bbeta_{J_\star})\geq\phi(\bbeta_{J_\star})$, but we shall prove in the following that this inequality is strict, that is, $\phi_R(\bbeta_{J_\star})>\phi(\bbeta_{J_\star})$, and consequently, $\phi_R>\phi$.

Let us temporarily assume that we already know the following facts hold:
\begin{enumerate}
    \item On $\Omega_{\mathrm{AO}}^{\mathrm{aux}}$,
    \begin{align}\label{eq:auxiliary-value-gap}
    &\phi_R-\phi\geq\frac{\sqrt N\sigma_\xi}{\log(\frac{e|J_\star^c|}{N})^{\frac{3}{2}}}.
    \end{align}
    \item $\phi$ is $\frac{10\sigma_\xi}{\sqrt{\log(\frac{e|J_\star^c|}{N})}}$ Lipschitz on $\Omega_{\mathrm{AO}}^{\mathrm{aux}}$.
\end{enumerate}
Then applying Lemma~\ref{lem:auxiliary-separation} to $\phi$ and $\phi_R$, with $E = \Omega_{\mathrm{AO}}^{\mathrm{aux}}$, we have
\begin{align*}
\begin{aligned}
&\mathbb P_{\mathrm{aux}}\{\phi\geq m+\eta\},\mathbb P_{\mathrm{aux}}\{\phi_R\leq m+2\eta\}\leq \mathsf p,\mbox{ where }
\mathsf p = 2^{6}\exp\bigg(-\frac{2^{-16}N}{\log(\frac{e|J_\star^c|}{N})^2}\bigg)+2^{6}\exp(-2^{-16}|J_\star|),
\end{aligned}
\end{align*}where $m$ is the median of $\tilde\phi$ defined in Lemma~\ref{lem:auxiliary-separation}, and $\eta=\frac{\sqrt N\sigma_\xi}{6\log(\frac{e|J_\star^c|}{N})^{\frac{3}{2}}}$.
Now, by Lemma~\ref{lem:common-cgmt} (upper tail) applied to $x=m+\eta$ and $\Phi=\Phi$, there holds $\bP\{\Phi>m+\eta\}\leq 2\mathsf p$; and by Lemma~\ref{lem:common-cgmt} again (lower tail), applied to $x=m+2\eta$, and $\Phi=\Phi_R$, there holds $\bP\{\Phi_R<m+2\eta\}\leq 2\mathsf p$. Therefore,
\begin{align*}
    \bP\left\{\norm{\hat\bbeta_{J_\star^c}[\vy-\bX_{J_\star}\hat\bbeta_{J_\star}}_2^2>R^2\right\}\leq \bP\{\Phi_R=\Phi\}\leq \bP\{\Phi>m+\eta\}+\bP\{\Phi_R<m+2\eta\}\leq 4\mathsf p.
\end{align*}
Therefore, it remains only to prove \eqref{eq:auxiliary-value-gap} and the Lipschitzness of $\phi$.

\paragraph{Proof of \eqref{eq:auxiliary-value-gap}.} Fix $\boldsymbol\beta_{J_\star}\in T$. Sion's theorem gives a saddle point $(\boldsymbol\beta_{J_\star^c}^\circ,a^\circ):=(\bbeta_{J_\star^c}(\bbeta_{J_\star}),a^\circ(\bbeta_{J_\star}))$. We first prove the following claim: on $\Omega_{\mathrm{AO}}^{\mathrm{aux}}$, every saddle point $(\boldsymbol\beta_{J_\star^c}^\circ,a^\circ)$ of $\phi(\boldsymbol\beta_{J_\star})$, uniformly in $\boldsymbol\beta_{J_\star}\in T$, satisfies
\begin{align}
&\|\boldsymbol\beta_{J_\star^c}^\circ\|_2^2\leq2^{17}\bigg(\frac{\sigma_\xi^2}{\log(\frac{e|J_\star^c|}{N})}+\|\bbeta_{J_\star^c}^\star\|_2^2\bigg)=\frac{R^2}{16},\mbox{ and } \frac{\lambda\sqrt N}{4\sqrt{\log(\frac{e|J_\star^c|}{N})}}<a^\circ\leq\frac{2\lambda\sqrt N}{\sqrt{\log(\frac{e|J_\star^c|}{N})}},
\label{eq:auxiliary-saddle-bounds}
\end{align}

\begin{enumerate}
    \item We first prove that 
    \begin{align}
    &\beta_j^\circ=\frac{s_{\boldsymbol\beta_{J_\star}}(\boldsymbol\beta_{J_\star^c}^\circ)}{\|\boldsymbol g\|_2}S_{\frac{\lambda\sqrt N}{a^\circ}}\bigg(h_j+\frac{\|\boldsymbol g\|_2}{s_{\boldsymbol\beta_{J_\star}}(\boldsymbol\beta_{J_\star^c}^\circ)}\beta_j^\star\bigg),\qquad j\in J_\star^c.
    \label{eq:coordinatewise-KKT}
    \end{align}To prove \eqref{eq:coordinatewise-KKT}, we need to solve the KKT equations at the saddle point. To this end, we first need to examine the normal cone of the feasible set; see, e.g., Section~\ref{sec:tools_convex_analysis}.
    \begin{enumerate}
        \item We first prove that $\|\bbeta_{J_\star^c}^\circ\|_2 < 32\sigma_\xi\sqrt{N}$. By the definition of $\phi(\bbeta_{J_\star})$, we know that $$\phi(\bbeta_{J_\star})\leq \max_{0\leq a\leq 8\sigma_\xi}\cF_{\bbeta_{J_\star}}(\vzero,a).$$ On $\Omega_{\mathrm{AO}}^{\mathrm{aux}}$, one may compute that $\|\vg\|_2\leq\sqrt{\frac{3}{2}N}$, $s_{\bbeta_{J_\star}}(\vzero)\leq \sqrt{3}\sigma_\xi$, $\|\vh_{J_\star}\|_2\leq 2\sqrt{|J_\star|}$, $\langle\vh_{J_\star},\bbeta_{J_\star^c}^\star\rangle\leq \frac{\sigma_\xi\sqrt N}{8^{\frac{3}{2}}}$, and $|h_0|\leq\frac{1}{64}\sqrt{N}$, hence $\cF_{\bbeta_{J_\star}}(\vzero,a)\leq \|\bbeta_{J_\star}\|_1+\frac{3a\sigma_\xi}{\lambda}-\frac{a^2}{2\lambda}$, and consequently, $\phi(\bbeta_{J_\star}) \leq \max_{0\leq a\leq 8\sigma_\xi}\cF_{\bbeta_{J_\star}}(\vzero,a)\leq \|\bbeta_{J_\star}\|_1+\frac{9\sigma_\xi^2}{2\lambda}$. On the other side, for $\bbeta_{J_\star^c}$ such that $\|\bbeta_{J_\star^c}\|_2 = 32\sigma_\xi\sqrt{N}$,  by the definition again, $\max_{0\leq a\leq8\sigma_\xi}\cF_{\bbeta_{J_\star}}(\bbeta_{J_\star^c},a)\geq \cF_{\bbeta_{J_\star}}(\bbeta_{J_\star^c},0)=\|\bbeta_{J_\star}\|_1+\|\bbeta_{J_\star^c}\|_1 = \|\bbeta_{J_\star}\|_1 + 32\sigma_\xi\sqrt{N}\geq \|\bbeta_{J_\star}\|_1+\frac{32\sigma_\xi^2}{\lambda}$, where we have used the assumption that $\lambda\geq\frac{\sigma_\xi}{\sqrt{N}}$. Suppose for contradiction that $\|\bbeta_{J_\star^c}^\circ\|_2 = 32\sigma_\xi\sqrt{N}$, then $\phi(\bbeta_{J_\star})\geq \|\bbeta_{J_\star}\|_1+\frac{32\sigma_\xi^2}{\lambda}$. However, we also have $\phi(\bbeta_{J_\star}) \leq \|\bbeta_{J_\star}\|_1+\frac{9\sigma_\xi^2}{2\lambda}$, leading to a contradiction. Therefore, $\|\bbeta_{J_\star^c}^\circ\|_2<32\sigma_\xi\sqrt{N}$.

        \item We next prove that $a^\circ > 0$. Suppose $a^\circ=0$, then by the definition of $\phi(\bbeta_{J_\star})$, $\min\cF_{\bbeta_{J_\star}}(\bbeta_{J_\star^c},0)=\|\bbeta\|_1=0$, indicating that $\bbeta_{J_\star^c}^\circ=\vzero$. Note that here $\max_{0\leq a\leq 8\sigma_\xi}\cF_{\bbeta_{J_\star}}(\bbeta_{J_\star^c},a)$ is a maximization problem, we must have $\partial_{a^+}\cF_{\bbeta_{J_\star}}(\vzero,0)\leq 0$, see Section~\ref{sec:tools_convex_analysis}. We compute that
        \begin{align}\label{eq:derivative_cF}
            \partial_a\cF_{\bbeta_{J_\star}}(\bbeta_{J_\star^c}^\circ,a) = \frac{\|\vg\|_2}{\lambda\sqrt{N}}s_{\bbeta_{J_\star}}(\bbeta_{J_\star^c}^\circ) - \frac{\langle\vh,\bbeta-\bbeta^\star\rangle}{\lambda\sqrt{N}} - \frac{\sigma_\xi h_0}{\lambda\sqrt{N}} - \frac{a}{\lambda}.
        \end{align}
        Hence $\partial_{a^+}\cF_{\bbeta_{J_\star}}(\vzero,0) = \frac{1}{\lambda}(\frac{\|\vg\|_2}{\sqrt{N}}s_{\bbeta_{J_\star}}(\vzero) - \frac{\langle\vh_{J_\star},\bbeta_{J_\star}-\bbeta_{J_\star}^\star\rangle}{\sqrt{N}} + \frac{\langle\vh_{J_\star^c},\bbeta_{J_\star^c}^\star\rangle}{\sqrt{N}} - \frac{\sigma_\xi}{\sqrt{N}}h_0)$. On $\Omega_{\mathrm{AO}}^{\mathrm{aux}}$, $\|\vg\|_2\geq\sqrt{\frac{N}{2}}$, $s_{\bbeta_{J_\star}}(\vzero)\geq\sigma_\xi$, $\frac{\langle\vh_{J_\star},\bbeta_{J_\star}-\bbeta_{J_\star}^\star\rangle}{\sqrt{N}} \leq 2\sigma_\xi\sqrt{\frac{|J_\star|}{N}}$, $\frac{\langle\vh_{J_\star^c},\bbeta_{J_\star^c}^\star\rangle}{\sqrt{N}} \leq \frac{\sigma_\xi}{8^{\frac{3}{2}}}$, and $\frac{\sigma_\xi}{\sqrt{N}}h_0\leq \frac{\sigma_\xi}{64}$, therefore $\partial_{a^+}\cF_{\bbeta_{J_\star}}(\vzero,0)>0$, leading to a contradiction. Consequently, $a^\circ>0$.  
    \end{enumerate}
    By the KKT condition of $\bbeta_{J_\star^c}^\circ$, with $N_{B_2(\vzero,32\sigma_\xi\sqrt{N})}(\bbeta_{J_\star^c}^\circ)=(\vzero)$, we obtain that $\vzero\in\partial_{\bbeta_{J_\star^c}}^-\cF_{\bbeta_{J_\star}}(\bbeta_{J_\star^c}^\circ,a^\circ)$. Moreover, $\nabla_{\bbeta_{J_\star^c}}s_{\bbeta_{J_\star}}(\bbeta_{J_\star^c}) = \frac{1}{s_{\bbeta_{J_\star}}(\bbeta_{J_\star^c})}(\bbeta_{J_\star^c}-\bbeta_{J_\star^c}^\star)$. Therefore, the KKT condition gives that for any $j\in J_\star^c$, $0\in \beta_j^\circ - \big(\beta_j^\star + \frac{s_{\bbeta_{J_\star}}(\bbeta_{J_\star^c}^\circ)}{\|\vg\|_2}h_j\big) + \frac{s_{\bbeta_{J_\star}}(\bbeta_{J_\star^c}^\circ)}{\|\vg\|_2}\frac{\lambda\sqrt{N}}{a^\circ}\partial^-|\cdot|(\beta_j^\circ)$. By \eqref{eq:softmax_equivalence} together with \eqref{eq:scaling_identity_softmax}, \eqref{eq:coordinatewise-KKT} follows.

    \item We prove that on $\Omega_{\mathrm{AO}}^{\mathrm{aux}}$, $\|\boldsymbol\beta_{J_\star^c}^\circ\|_2^2\leq2^{17}\frac{\sigma_\xi^2}{\log(\frac{e|J_\star^c|}{N})}=\frac{R^2}{16}$. From $a^\circ>0$, and from Section~\ref{sec:tools_convex_analysis}, we know that on $0 < a \le 8\sigma_\xi$ the derivative of $a\mapsto \cF_{\bbeta_{J_\star}}(\bbeta_{J_\star^c}^\circ,a)$ must be nonnegative, namely, $\partial_a\cF_{\bbeta_{J_\star}}(\bbeta_{J_\star^c}^\circ,a)\geq 0$. Splitting $\langle\vh,\bbeta-\bbeta^\star\rangle$ in \eqref{eq:derivative_cF} into $\langle\vh_{J_\star},\bbeta_{J_\star}-\bbeta_{J_\star}^\star\rangle + \langle\vh_{J_\star^c},\bbeta_{J_\star^c}-\bbeta_{J_\star^c}^\star\rangle$, rearranging terms, using $\partial_a\cF_{\bbeta_{J_\star}}(\bbeta_{J_\star^c}^\circ,a)\ge 0$, and multiplying $\frac{\|\vg\|_2}{s_{\bbeta_{J_\star}}(\bbeta_{J_\star^c}^\circ)\sqrt{N}}$ on both sides, together with the definition of $\Omega_{\mathrm{AO}}^{\mathrm{aux}}$, we obtain
    \begin{align}
    &\frac{\|\boldsymbol g\|_2\langle\boldsymbol h_{J_\star^c},\boldsymbol\beta_{J_\star^c}^\circ-\boldsymbol\beta_{J_\star^c}^\star\rangle}{N s_{\boldsymbol\beta_{J_\star}}(\boldsymbol\beta_{J_\star^c}^\circ)}\leq\frac{\|\boldsymbol g\|_2^2}{N}-\frac{a^\circ\|\boldsymbol g\|_2}{\sqrt N s_{\boldsymbol\beta_{J_\star}}(\boldsymbol\beta_{J_\star^c}^\circ)}-\frac{\|\boldsymbol g\|_2(\langle\boldsymbol h_{J_\star},\boldsymbol\beta_{J_\star}-\boldsymbol\beta_{J_\star}^\star\rangle+\sigma_\xi h_0)}{N s_{\boldsymbol\beta_{J_\star}}(\boldsymbol\beta_{J_\star^c}^\circ)}\leq2.
    \label{eq:soft-threshold-inner-product}
    \end{align}Notice that for $\gamma=\frac{\|\boldsymbol g\|_2}{\sqrt N}$, $s=s_{\boldsymbol\beta_{J_\star}}(\boldsymbol\beta_{J_\star^c}^\circ)$, and $\tau=\frac{\lambda\sqrt N}{a^\circ}$, by \eqref{eq:coordinatewise-KKT},
    \begin{align*}
       \forall j\in J_\star^c,\qquad S_\tau\big(h_j + \frac{\gamma\sqrt{N}}{s}\beta_j^\star\big) = \frac{\|\vg\|_2}{s_{\bbeta_{J_\star}}(\bbeta_{J_\star^c}^\circ)}\beta_j^\circ \mbox{ and }\, S_\tau\big( h_j + \frac{\gamma\sqrt{N}}{s}\beta_j^\star \big)-\frac{\gamma\sqrt{N}}{s}\beta_j^\star = \frac{\|\vg\|_2}{s_{\bbeta_{J_\star}}(\bbeta_{J_\star^c}^\circ)}\big(\beta_j^\circ - \beta_j^\star\big),
    \end{align*}
    which, together with \eqref{eq:soft-threshold-inner-product}, verifies the left-hand-side of \eqref{eq:soft-threshold-first-moment}. Applying Lemma~\ref{lem:soft-threshold} together with $\Omega_{\mathrm{AO}}^{\mathrm{aux}}$ and the definition of $s_{\bbeta_{J_\star}}(\bbeta_{J_\star^c}^\circ)$ yields that
    \begin{align*}
        \norm{\bbeta_{J_\star^c}^\circ}_2^2 \leq \frac{2^{14}Ns_{\bbeta_{J_\star}}(\bbeta_{J_\star^c}^\circ)^2}{\|\vg\|_2^2\log(\frac{e|J_\star^c|}{N})} + 2\|\bbeta_{J_\star^c}^\star\|_2^2 \leq \frac{2^{15}\big(\sigma_\xi^2+r_\star^2+2\|\bbeta_{J_\star^c}^\circ\|_2^2+2\|\bbeta_{J_\star^c}^\star\|_2^2\big)}{\log(\frac{e|J_\star^c|}{N})} + 2\|\bbeta_{J_\star^c}^\star\|_2^2,
    \end{align*}
    which indicates the first inequality of \eqref{eq:auxiliary-saddle-bounds}, thanks to the fact that $\log(\frac{e|J_\star^c|}{N})\geq 2^{18}$ (which is a consequence of the assumption of Proposition~\ref{prop:lambda_is_small_shift}).

    \item Now we estimate $a^\circ$, that is, we prove the second inequality of \eqref{eq:auxiliary-saddle-bounds}.
    Apply Lemma~\ref{lem:soft-threshold} with the same choices of  $\gamma$, $s$, and $\tau$. Equations \eqref{eq:coordinatewise-KKT}--\eqref{eq:soft-threshold-inner-product} give $\tau\geq\frac{\sqrt{\log(\frac{e|J_\star^c|}{N})}}{2}$, and hence by the relation $\tau = \frac{\lambda\sqrt{N}}{a^\circ}$, we know that $a^\circ\leq\frac{2\lambda\sqrt N}{\sqrt{\log(\frac{e|J_\star^c|}{N})}}<8\sigma_\xi$. Thus $0<a^\circ<8\sigma_\xi$, and by Section~\ref{sec:tools_convex_analysis}, the first-order condition implies that \eqref{eq:soft-threshold-inner-product} is an equality, that is,
    \begin{align}\tag{\ref{eq:soft-threshold-inner-product}}
        \frac{\|\boldsymbol g\|_2\langle\boldsymbol h_{J_\star^c},\boldsymbol\beta_{J_\star^c}^\circ-\boldsymbol\beta_{J_\star^c}^\star\rangle}{N s_{\boldsymbol\beta_{J_\star}}(\boldsymbol\beta_{J_\star^c}^\circ)} = \frac{\|\boldsymbol g\|_2^2}{N}-\frac{a^\circ\|\boldsymbol g\|_2}{\sqrt N s_{\boldsymbol\beta_{J_\star}}(\boldsymbol\beta_{J_\star^c}^\circ)}-\frac{\|\boldsymbol g\|_2(\langle\boldsymbol h_{J_\star},\boldsymbol\beta_{J_\star}-\boldsymbol\beta_{J_\star}^\star\rangle+\sigma_\xi h_0)}{N s_{\boldsymbol\beta_{J_\star}}(\boldsymbol\beta_{J_\star^c}^\circ)}.
    \end{align}
    Suppose for contradiction that $\tau\geq4\sqrt{\log(\frac{e|J_\star^c|}{N})}$. Plugging-in the solution of $\bbeta_{J_\star^c}^\circ$ from \eqref{eq:coordinatewise-KKT} into the left-hand-side of \eqref{eq:soft-threshold-inner-product}, we know that the left-hand-side of \eqref{eq:soft-threshold-inner-product} is precisely $\frac1N\sum_{j\in J_\star^c}h_j\big(S_\tau\big(h_j+\frac{\gamma\sqrt N}{s}\beta_j^\star\big)-\frac{\gamma\sqrt N}{s}\beta_j^\star\big)$, hence Lemma~\ref{lem:soft-threshold} bounds it from above by $\frac{1}{4}$; whereas \eqref{eq:soft-threshold-inner-product} and the definition of $\Omega_{\mathrm{AO}}^{\mathrm{aux}}$ give $\frac{\|\boldsymbol g\|_2\langle\boldsymbol h_{J_\star^c},\boldsymbol\beta_{J_\star^c}^\circ-\boldsymbol\beta_{J_\star^c}^\star\rangle}{N s_{\boldsymbol\beta_{J_\star}}(\boldsymbol\beta_{J_\star^c}^\circ)}\geq \frac12-2^{-10}-\sqrt{\frac32}\bigg(\frac{2\sqrt{N|J_\star|}\,r_\star}{N\sigma_\xi}+\frac1{64}\bigg)>\frac13.$ This contradiction proves $\tau<4\sqrt{\log(\frac{e|J_\star^c|}{N})}$ and the asserted two-sided bound for $a^\circ$.
\end{enumerate}

We prove \eqref{eq:auxiliary-value-gap}. Fix $\boldsymbol\beta_{J_\star}\in T$ and a saddle point $(\boldsymbol\beta_{J_\star^c}^\circ,a^\circ) = (\bbeta_{J_\star^c}(\bbeta_{J_\star}),a^\circ(\bbeta_{J_\star}))$. Along every segment joining $\boldsymbol\beta_{J_\star^c}^\circ$ to $RB_2^{J_\star^c}$, $s_{\boldsymbol\beta_{J_\star}}(\boldsymbol\beta_{J_\star^c})^2\leq2^{5}\sigma_\xi^2$. 
The function $\boldsymbol\beta_{J_\star^c}\mapsto\mathcal F_{\boldsymbol\beta_{J_\star}}(\boldsymbol\beta_{J_\star^c},a^\circ) - \|\bbeta_{J_\star^c}\|_1$ has Hessian
\begin{align*}
&\frac{a^\circ\|\boldsymbol g\|_2}{\lambda\sqrt N}\bigg(\frac{I_{J_\star^c}}{s_{\boldsymbol\beta_{J_\star}}(\boldsymbol\beta_{J_\star^c})}-\frac{(\boldsymbol\beta_{J_\star^c}-\boldsymbol\beta_{J_\star^c}^\star)\otimes(\boldsymbol\beta_{J_\star^c}-\boldsymbol\beta_{J_\star^c}^\star)}{s_{\boldsymbol\beta_{J_\star}}(\boldsymbol\beta_{J_\star^c})^3}\bigg).
\end{align*}
Consequently, for every $\boldsymbol v\in V_{J_\star^c}$, by using $s_{\bbeta_{J_\star}}(\bbeta_{J_\star^c})^2\geq \|\bbeta_{J_\star^c}-\bbeta_{J_\star^c}^\star\|_2^2$ and $\langle\bbeta_{J_\star^c}-\bbeta_{J_\star^c}^\star,\vv\rangle^2\leq\|\bbeta_{J_\star^c}-\bbeta_{J_\star^c}^\star\|_2^2\|\vv\|_2^2$,
\begin{align}
&\big\langle\boldsymbol v,\nabla^2_{\boldsymbol\beta_{J_\star^c}}(\mathcal F_{\boldsymbol\beta_{J_\star}}(\boldsymbol\beta_{J_\star^c},a^\circ)-\|\boldsymbol\beta_{J_\star^c}\|_1)\boldsymbol v\big\rangle\geq \frac{a^\circ\|\vg\|_2\sigma_\xi^2}{\lambda\sqrt{N}s_{\bbeta_{J_\star}}(\bbeta_{J_\star^c})^3}\|\vv\|_2^2\geq\frac{2^{-8}a^\circ}{\lambda\sigma_\xi}\|\boldsymbol v\|_2^2.
\label{eq:Hessian-lower-bound}
\end{align}
The term $\|\boldsymbol\beta_{J_\star^c}\|_1$ is convex. Lemma~\ref{lem:radial-restriction-gap}, \eqref{eq:auxiliary-saddle-bounds}, and \eqref{eq:Hessian-lower-bound} therefore give, uniformly in $\boldsymbol\beta_{J_\star}\in T$,
\begin{align*}
&\phi_R(\boldsymbol\beta_{J_\star})\geq\phi(\boldsymbol\beta_{J_\star})+\frac{2^{-9}a^\circ}{\lambda\sigma_\xi}\frac{9R^2}{16}\geq\phi(\boldsymbol\beta_{J_\star})+\frac{\sqrt N\sigma_\xi}{\log(\frac{e|J_\star^c|}{N})^{\frac{3}{2}}}.
\end{align*}
Taking the infimum over $\boldsymbol\beta_{J_\star}\in T$ proves \eqref{eq:auxiliary-value-gap}.
\paragraph{On the Lipschitz constant of \texorpdfstring{$\phi$}{phi}.}
At every saddle point of \eqref{eq:auxiliary-optimization}, $s_{\boldsymbol\beta_{J_\star}^\circ}(\boldsymbol\beta_{J_\star^c}^\circ)\leq3\sigma_\xi$. The saddle inequalities at two auxiliary Gaussian inputs, together with \eqref{eq:auxiliary-saddle-bounds}, now for any $\vomega=(\vg,\vh_{J_\star},\vh_{J_\star^c},h_0)$ and $\vomega'=(\vg',\vh_{J_\star}',\vh_{J_\star^c}',h_0')$, let $(\boldsymbol\beta_{J_\star}^\circ,\boldsymbol\beta_{J_\star^c}^\circ,a^\circ)$ and $(\boldsymbol\beta_{J_\star}^{\circ\prime},\boldsymbol\beta_{J_\star^c}^{\circ\prime},a^{\circ\prime})$ be corresponding global saddle points. By the saddle-point inequalities,
\begin{align*}
\begin{aligned}
&\phi(\vomega)-\phi(\vomega')
\leq\mathcal F_{\boldsymbol\beta_{J_\star}^{\circ\prime}}(\boldsymbol\beta_{J_\star^c}^{\circ\prime},a^\circ;\vomega)-\mathcal F_{\boldsymbol\beta_{J_\star}^{\circ\prime}}(\boldsymbol\beta_{J_\star^c}^{\circ\prime},a^\circ;\vomega')\\
&=\frac{a^\circ}{\lambda\sqrt N}\big(\|\vg\|_2-\|\vg'\|_2\big)s_{\boldsymbol\beta_{J_\star}^{\circ\prime}}(\boldsymbol\beta_{J_\star^c}^{\circ\prime})-\frac{a^\circ}{\lambda\sqrt N}\langle\vh_{J_\star}-\vh_{J_\star}',\boldsymbol\beta_{J_\star}^{\circ\prime}-\boldsymbol\beta_{J_\star}^\star\rangle-\frac{a^\circ}{\lambda\sqrt N}\langle\vh_{J_\star^c}-\vh_{J_\star^c}',\boldsymbol\beta_{J_\star^c}^{\circ\prime}-\boldsymbol\beta_{J_\star^c}^\star\rangle-\frac{a^\circ\sigma_\xi}{\lambda\sqrt N}(h_0-h_0').
\end{aligned}
\end{align*}
Moreover,
\begin{align*}
\begin{aligned}
&\big|\langle\vh_{J_\star}-\vh_{J_\star}',\boldsymbol\beta_{J_\star}^{\circ\prime}-\boldsymbol\beta_{J_\star}^\star\rangle+\langle\vh_{J_\star^c}-\vh_{J_\star^c}',\boldsymbol\beta_{J_\star^c}^{\circ\prime}-\boldsymbol\beta_{J_\star^c}^\star\rangle+\sigma_\xi(h_0-h_0')\big| \leq s_{\boldsymbol\beta_{J_\star}^{\circ\prime}}(\boldsymbol\beta_{J_\star^c}^{\circ\prime})\big(\|\vh_{J_\star}-\vh_{J_\star}'\|_2^2+\|\vh_{J_\star^c}-\vh_{J_\star^c}'\|_2^2+|h_0-h_0'|^2\big)^{\frac{1}{2}}.
\end{aligned}
\end{align*}
Together with $|\|\vg\|_2-\|\vg'\|_2|\leq\|\vg-\vg'\|_2$, $a^\circ\leq\frac{2\lambda\sqrt N}{\sqrt{\log(\frac{e|J_\star^c|}{N})}}$, and $s_{\boldsymbol\beta_{J_\star}^{\circ\prime}}(\boldsymbol\beta_{J_\star^c}^{\circ\prime})\leq3\sigma_\xi$, this gives
\begin{align*}
\begin{aligned}
\phi(\vomega)-\phi(\vomega')
&\leq\frac{6\sigma_\xi}{\sqrt{\log(\frac{e|J_\star^c|}{N})}}\bigg[\|\vg-\vg'\|_2+\big(\|\vh_{J_\star}-\vh_{J_\star}'\|_2^2+\|\vh_{J_\star^c}-\vh_{J_\star^c}'\|_2^2+|h_0-h_0'|^2\big)^{\frac{1}{2}}\bigg]\\
&\leq\frac{10\sigma_\xi}{\sqrt{\log(\frac{e|J_\star^c|}{N})}}\|\vomega-\vomega'\|_2.
\end{aligned}
\end{align*}
Exchanging $\vomega$ and $\vomega'$ completes the proof.
\end{proof}

\subsubsection{Upper bound for the energy}

\begin{Proposition}\label{prop:null_signal_upper_bound}
Grant Assumption~\ref{assumption:LASSO_upper}. Suppose $\lambda>0$, $|J_\star^c|\geq 100N$.
Let $q,d_q$, and $\mathrm{Var}_{J_\star^c}$ be given by \eqref{eq:common-q-calibration} and \eqref{eq:common-free-variance}.
Fix $c_0=2^{-57}$ and $C_1=2^{6}$. Then
\begin{align*}
\mathbb P_{\mathbb X_{J_\star^c},\boldsymbol\xi}\bigg\{ \|\widehat{\boldsymbol\beta}_{J_\star^c}[\boldsymbol\xi]\|_2^2>2^{12}\mathrm{Var}_{J_\star^c}\bigg\}\leq C_1e^{-c_0d_q}.
\end{align*}
\end{Proposition}

\begin{Lemma}\label{lemma:null_signal_events}
Let $\boldsymbol g\in\mathbb R^N$ and $\boldsymbol h\in\mathbb R^{J_\star^c}$ be independent standard Gaussian vectors, independent of $\boldsymbol\xi$. For $\boldsymbol h\in\mathbb R^{J_\star^c}$ and $u\geq0$, write
\begin{align*}
A_{\boldsymbol h}(u)=\sum_{j\in J_\star^c}|h_j|(|h_j|-u)_+,\qquad B_{\boldsymbol h}(u)=\sum_{j\in J_\star^c}(|h_j|-u)_+^2.
\end{align*}
Define $\Omega_\xi=\big\{\big|\frac{\|\boldsymbol\xi\|_2}{\sigma_\xi\sqrt N}-1\big|\leq2^{-23}\sqrt{\frac{d_q}{N}}\big\}$. For each $\bxi\in\Omega_\xi$, define $\Omega_{\mathrm{aux}}^{\mathrm{aux}}(\bxi)$ by the intersection of the following events
\begin{align*}
&\Omega_{\boldsymbol g}^{\mathrm{aux}}=\bigg\{\bigg|\frac{\|\boldsymbol g\|_2^2}{N}-1\bigg|\leq2^{-23}\sqrt{\frac{d_q}{N}},\ \frac{|\langle\boldsymbol\xi,\boldsymbol g\rangle|}{\|\boldsymbol\xi\|_2}\leq2^{-23}\sqrt{d_q}\bigg\},\\
&\Omega_{\underline{z}_q,1}^{\mathrm{aux}}=\bigg\{\bigg|\frac{A_{\boldsymbol h}(\underline{z}_q)}{2|J_\star^c|\Phi_G(\underline{z}_q)}-1\bigg|\leq2^{-23}\bigg\},\,\Omega_{\overline{z}_q,1}^{\mathrm{aux}}=\bigg\{\bigg|\frac{A_{\boldsymbol h}(\overline{z}_q)}{2|J_\star^c|\Phi_G(\overline{z}_q)}-1\bigg|\leq2^{-23}\bigg\},\mbox{ and }\\
&\Omega_{\underline{z}_q,2}^{\mathrm{aux}}=\bigg\{\bigg|\sqrt{\frac{B_{\boldsymbol h}(\underline{z}_q)}{|J_\star^c|\mathbb E(|G|-\underline{z}_q)_+^2}}-1\bigg|\leq2^{-23}\bigg\},\,\Omega_{\overline{z}_q,2}^{\mathrm{aux}}=\bigg\{\bigg|\sqrt{\frac{B_{\boldsymbol h}(\overline{z}_q)}{|J_\star^c|\mathbb E(|G|-\overline{z}_q)_+^2}}-1\bigg|\leq2^{-23}\bigg\}.
\end{align*}
If $C_1e^{-c_0d_q}<1$, then
\begin{align}
\mathbb P(\Omega_\xi^c)\leq2^{-4}C_1e^{-c_0d_q},\mbox{ and } \sup_{\boldsymbol\xi\in\Omega_\xi}\mathbb P_{\mathrm{aux}}(\big(\Omega_{\mathrm{aux}}^{\mathrm{aux}}(\bxi)\big)^c\mid\boldsymbol\xi)\leq2^{-4}C_1e^{-c_0d_q}.
\label{eq:omega_probability}
\end{align}
\end{Lemma}

\begin{proof}[Proof of Lemma~\ref{lemma:null_signal_events}]
The estimate for $\Omega_\xi$ and $\Omega_\vg^{\mathrm{aux}}$ are trivial, thus omitted.
\begin{enumerate}
    \item For $\Omega_{\underline{z}_q,1}^{\mathrm{aux}}$ and $\Omega_{\overline{z}_q,1}^{\mathrm{aux}}$. Let $t\in(\underline{z}_q,\overline{z}_q)$. \eqref{eq:exp_G_soft_max} gives $\mathbb E[|G|(|G|-t)_+]=2\Phi_G(t)$, hence $\bE[A_\vh(t)] = 2|J_\star^c|\Phi_G(t)$. Item~\emph{\ref{item:expectation_soft_max_moment}} in Section~\ref{sec:Gaussian_calculations} yields that $\mathbb E[|G|^k(|G|-t)_+^k]\leq2\Phi_G(t)2^kk!$, for any $k\geq2$. Therefore, Bernstein's inequality, see, for instance, \cite[Theorem 2.10]{boucheronbook}, gives that for any $x>0$, with probability at least $1-2\exp(-x)$, $|A_\vh(t)-\bE[A_\vh(t)]|\leq 8\sqrt{\bE[A_\vh(t)]x}+4x$. Take $x=2^{-52}d_q$. By item~\emph{\ref{item:Gaussian_local_comparison}} in Section~\ref{sec:Gaussian_calculations}, $\bE[A_\vh(t)]\geq e^{-\frac{1}{8}}q$. Together with $d_q\leq\frac{2q}{z_q^2}$ and $z_q>2.57$ from Lemma~\ref{lem:common-calibration}, we obtain that $\bE[A_\vh(t)]\gtrsim d_q$ and hence $8\sqrt{\bE[A_\vh(t)]x}+4x\lesssim 2^{-23}\bE[A_\vh(t)]$.

    \item For $\Omega_{\underline{z}_q,2}^{\mathrm{aux}}$ and $\Omega_{\overline{z}_q,2}^{\mathrm{aux}}$. Let $t\in(\underline{z}_q,\overline{z}_q)$. By Lemma~\ref{lem:common-calibration}, $z_q>2.57$, and hence $t>2.5$. The map $\boldsymbol h\mapsto B_{\boldsymbol h}(t)^{\frac{1}{2}}$ is $1$-Lipschitz, and Gaussian Poincar\'e gives $\big|\mathbb EB_{\boldsymbol h}(t)^{\frac{1}{2}}-[|J_\star^c|\mathbb E(|G|-t)_+^2]^{\frac{1}{2}}\big|\leq1$. Moreover, item~\emph{\ref{item:Gaussian_local_comparison}} in Section~\ref{sec:Gaussian_calculations} gives $2^{-4}d_q\leq|J_\star^c|\mathbb E(|G|-t)_+^2\leq2^4d_q$. Gaussian concentration therefore gives, for any $x>0$, with probability at least $1-2\exp(-x)$, $\big|B_{\boldsymbol h}(t)^{\frac{1}{2}}-\mathbb EB_{\boldsymbol h}(t)^{\frac{1}{2}}\big|\leq\sqrt{2x}.$ By triangular inequality,
    \begin{align*}
        &\big| B_\vh(t)^{\frac12} - \big( |J_\star^c|\bE[(|G|-t)_+^2] \big)^{\frac{1}{2}} \big| \leq \big| B_\vh(t)^{\frac{1}{2}} - \bE[B_\vh(t)^{\frac{1}{2}}]\big| + \big| \bE[B_\vh(t)^{\frac{1}{2}}] - \big( |J_\star^c|\bE[(|G|-t)_+^2] \big)^{\frac{1}{2}} \big| \leq 1+\sqrt{2x}.
    \end{align*} Taking $x=2^{-52}d_q$, dividing by $\big( |J_\star^c|\bE[(|G|-t)_+^2] \big)^{\frac{1}{2}}$, and using $|J_\star^c|\bE[(|G|-t)_+^2]\geq 2^{-4}d_q$ from \eqref{eq:Gaussian_local_comparison_soft_max} in Section~\ref{sec:Gaussian_calculations}, the proof is complete.
\end{enumerate}
\end{proof}

In the following, we fix $\bxi \in \Omega_\xi$, and work on $\Omega_{\mathrm{aux}}^{\mathrm{aux}}$. Let $H, K > 0$, and $I \subset [0, K]$ be an interval. Define the following function:
\begin{align}\label{eq:def_Psi}
\Psi:=\Psi_{\vg,\vh}:(t,s)\in I\times[0,H]\mapsto s\|\boldsymbol\xi+t\boldsymbol g\|_2-\frac{s^2}{2}-\frac{\|\boldsymbol\xi\|_2^2}{2}+t\begin{cases}-[B_{\boldsymbol h}(\frac{N\lambda}{s})]^{\frac{1}{2}}s,&s\|\boldsymbol h\|_\infty>N\lambda,\\N\lambda-s\|\boldsymbol h\|_\infty,&s\|\boldsymbol h\|_\infty\leq N\lambda.
\end{cases}
\end{align}
Let us first assume that the following two lemmas hold:
\begin{Lemma}\label{lemma:null_signal_scalarization}
For $K,H>0$ and a nonempty closed set $I\subset[0,K]$, define
\begin{align*}
\Phi_I^{K,H}&=\min_{\substack{\boldsymbol\beta_{J_\star^c}\in\mathbb R^{J_\star^c}\\\|\boldsymbol\beta_{J_\star^c}\|_2\in I}}\max_{\substack{\boldsymbol u\in\mathbb R^N\\\|\boldsymbol u\|_2\leq H}}\bigg\{\langle\boldsymbol u,\boldsymbol\xi-\mathbb X_{J_\star^c}\boldsymbol\beta_{J_\star^c}\rangle-\frac12\|\boldsymbol u\|_2^2+N\lambda\|\boldsymbol\beta_{J_\star^c}\|_1-\frac12\|\boldsymbol\xi\|_2^2\bigg\},\\
\phi_I^{K,H}&=\min_{\substack{\boldsymbol\beta_{J_\star^c}\in\mathbb R^{J_\star^c}\\\|\boldsymbol\beta_{J_\star^c}\|_2\in I}}\max_{\substack{\boldsymbol u\in\mathbb R^N\\\|\boldsymbol u\|_2\leq H}}\bigg\{\|\boldsymbol\beta_{J_\star^c}\|_2\langle\boldsymbol g,\boldsymbol u\rangle+\|\boldsymbol u\|_2\langle\boldsymbol h,\boldsymbol\beta_{J_\star^c}\rangle+\langle\boldsymbol u,\boldsymbol\xi\rangle-\frac12\|\boldsymbol u\|_2^2+N\lambda\|\boldsymbol\beta_{J_\star^c}\|_1-\frac12\|\boldsymbol\xi\|_2^2\bigg\}.
\end{align*}
Then $\phi_I^{K,H}=\min_{t\in I}\max_{0\leq s\leq H}\Psi(t,s).$ Moreover, for every $x\in\mathbb R$,
\begin{align}
\mathbb P_{\mathbb X_{J_\star^c}}\{\Phi_I^{K,H}<x\mid\boldsymbol\xi\}&\leq2\mathbb P_{\mathrm{aux}}\{\phi_I^{K,H}\leq x\mid\boldsymbol\xi\},\mbox{ and } \mathbb P_{\mathbb X_{J_\star^c}}\{\Phi_{[0,K]}^{K,H}>x\mid\boldsymbol\xi\}\leq2\mathbb P_{\mathrm{aux}}\{\phi_{[0,K]}^{K,H}\geq x\mid\boldsymbol\xi\}.
\label{eq:cgmt_comparisons}
\end{align}
Assume $C_1e^{-c_0d_q}<1$, $K\geq2^{7}\sqrt{\mathrm{Var}_{J_\star^c}}$, and $H\geq\frac{2^{5}N\lambda}{z_q},$. On $\Omega^{\mathrm{aux}}$, every saddle point $(t^\circ,s^\circ)$ of $\min_{0\leq t\leq K}\max_{0\leq s\leq H}\Psi(t,s)$ satisfies
\begin{align}
0\leq t^\circ\leq2^{5}\sqrt{\mathrm{Var}_{J_\star^c}},\qquad \frac{2^{-4}N\lambda}{z_q}\leq s^\circ\leq\frac{2^4N\lambda}{z_q}.
\label{eq:null_signal_saddle_bounds}
\end{align}
\end{Lemma}

\begin{Lemma}\label{lemma:null_signal_cost_separation}
Assume $C_1e^{-c_0d_q}<1$, $K\geq2^{7}\sqrt{\mathrm{Var}_{J_\star^c}}$, and $H\geq\frac{2^{5}N\lambda}{z_q}.$ There are conditional deterministic levels $\alpha_1(\boldsymbol\xi)<\gamma_1(\boldsymbol\xi)$ such that, for every $\boldsymbol\xi\in\Omega_\xi$, such that
\begin{align}
\mathbb P_{\mathrm{aux}}\{\phi_{[0,K]}^{K,H}\geq\alpha_1(\boldsymbol\xi)\mid\boldsymbol\xi\}+\mathbb P_{\mathrm{aux}}\big\{\phi_{[2^{6}\sqrt{\mathrm{Var}_{J_\star^c}},K]}^{K,H}\leq\gamma_1(\boldsymbol\xi)\mid\boldsymbol\xi\big\}\leq2^{-2}C_1e^{-c_0d_q}.
\label{eq:deterministic_levels}
\end{align}
\end{Lemma}

\begin{proof}[Proof of Proposition~\ref{prop:null_signal_upper_bound}]

Our goal is to prove an upper bound on the probability of the random event $\{\|\hat\bbeta_{J_\star^c}[\bxi]\|_2^2 > 2^{12} \Var_{J_\star^c}\}$. We shall embed this event into a random event concerning the saddle point of $\Phi_I^{K,H}$ defined in Lemma~\ref{lemma:null_signal_scalarization}. By \eqref{eq:variational_representation_LASSO}, we know that $(\hat\bbeta_{J_\star^c}[\bxi], \bxi - \bX_{J_\star^c}\hat\bbeta_{J_\star^c}[\bxi])$ is a saddle point of $\Phi_{\mathbb{R}_+}^{\infty,\infty}$. We first show that for sufficiently large $K$ and all sufficiently large $H_n$, $(\hat\bbeta_{J_\star^c}[\bxi], \bxi - \bX_{J_\star^c}\hat\bbeta_{J_\star^c}[\bxi])$ is also a saddle point of $\Phi_{[0,K]}^{K,H_n}$.

By the convex optimization definition of $\hat\bbeta_{J_\star^c}$, we have $\frac{1}{2N\lambda}\|\bxi - \bX_{J_\star^c}\hat\bbeta_{J_\star^c}[\bxi]\|_2^2 + \|\hat\bbeta_{J_\star^c}[\bxi]\|_1 \le \frac{1}{2N\lambda}\|\bxi - \bX_{J_\star^c}\vzero\|_2^2 + \|\vzero\|_1$, which yields $\|\hat\bbeta_{J_\star^c}[\bxi]\|_2 \le \frac{1}{2N\lambda}\|\bxi\|_2^2\leq \frac{\sigma_\xi^2}{\lambda}$ since we work on $\bxi\in\Omega_\xi$. Taking $K > \max\{1 + \frac{\sigma_\xi^2}{\lambda}, 2^{7}\sqrt{\Var_{J_\star^c}}\}$, we have $\|\hat\bbeta_{J_\star^c}[\bxi]\|_2 < K$. Set $I=(2^{6}\sqrt{\Var_{J_\star^c}},K)$. Consequently, $\{\|\hat\bbeta_{J_\star^c}[\bxi]\|_2^2 > 2^{12} \Var_{J_\star^c}\}\subset\{\|\hat\bbeta_{J_\star^c}[\bxi]\|_2\in I\}.$

We next localize the Lagrangian $\vu$. For every $n\in\mathbb N$, set $H_n=\frac{2^{5}N\lambda}{z_q}+n$. For any $(\bX_{J_\star^c},\bxi)$, take an integer $n_0=n_0(\bX_{J_\star^c},\bxi)$ such that $n_0>\big(\|\bxi\|_2+K\|\bX_{J_\star^c}\|_{\mathrm{op}}-\frac{2^{5}N\lambda}{z_q}\big)_+.$
Then, for every $n\geq n_0$, $H_n>\|\bxi\|_2+K\|\bX_{J_\star^c}\|_{\mathrm{op}}\geq\sup_{\|\bbeta_{J_\star^c}\|_2\leq K}\|\bxi-\bX_{J_\star^c}\bbeta_{J_\star^c}\|_2,$ and hence $\|\bxi-\bX_{J_\star^c}\hat\bbeta_{J_\star^c}[\bxi]\|_2<H_n$. Therefore, on $\{\|\hat\bbeta_{J_\star^c}[\bxi]\|_2\in I\}$, for every $n\geq n_0$, $\hat\bbeta_{J_\star^c}[\bxi]$ is the $V_{J_\star^c}$ component of a saddle point of both $\Phi_I^{K,H_n}$ and $\Phi_{[0,K]}^{K,H_n}$. Hence
\[
\{\|\hat\bbeta_{J_\star^c}[\bxi]\|_2\in I\}\subset\liminf_{n\to\infty}\{\Phi_{[0,K]}^{K,H_n}=\Phi_I^{K,H_n}\}=\bigcup_{m=1}^\infty\bigcap_{n\geq m}\{\Phi_{[0,K]}^{K,H_n}=\Phi_I^{K,H_n}\}.
\]

We next consider a further upper bound for these events. By $\alpha_1(\bxi)<\gamma_1(\bxi)$ provided by Lemma~\ref{lemma:null_signal_cost_separation}, for every $n\in\mathbb N$,
\[
\{\Phi_{[0,K]}^{K,H_n}=\Phi_I^{K,H_n}\}\subset\{\Phi_{[0,K]}^{K,H_n}>\alpha_1(\bxi)\}\cup\{\Phi_I^{K,H_n}\leq\gamma_1(\bxi)\}.
\]
Indeed, if there holds $\{\Phi_{[0,K]}^{K,H_n}\leq\alpha_1(\bxi)\}\cap\{\Phi_I^{K,H_n}>\gamma_1(\bxi)\}$, then from $\alpha_1(\bxi)<\gamma_1(\bxi)$, $\Phi_{[0,K]}^{K,H_n}\neq\Phi_I^{K,H_n}$. Therefore, by the preceding inclusion and the Fatou inequality,
\begin{align*}
\bP\big\{\|\hat\bbeta_{J_\star^c}[\bxi]\|_2^2>2^{12}\Var_{J_\star^c}\big|\bxi\big\}
&\leq\bP\bigg\{\liminf_{n\to\infty}\{\Phi_{[0,K]}^{K,H_n}=\Phi_I^{K,H_n}\}\bigg|\bxi\bigg\} \leq\liminf_{n\to\infty}\bP\big\{\Phi_{[0,K]}^{K,H_n}=\Phi_I^{K,H_n}\big|\bxi\big\}\\
&\leq\liminf_{n\to\infty}\bigg[\bP\big\{\Phi_{[0,K]}^{K,H_n}>\alpha_1(\bxi)\big|\bxi\big\}+\bP\big\{\Phi_I^{K,H_n}\leq\gamma_1(\bxi)\big|\bxi\big\}\bigg],
\end{align*}
and Proposition~\ref{prop:null_signal_upper_bound} follows by combining \eqref{eq:cgmt_comparisons}, \eqref{eq:deterministic_levels}, and Lemma~\ref{lemma:null_signal_events}. 
\end{proof}

\begin{proof}[Proof of Lemma~\ref{lemma:null_signal_scalarization}]
By conditioning on $\bxi$ and using CGMT, namely Lemma~\ref{lem:common-cgmt}, \eqref{eq:cgmt_comparisons} follows immediately. For the auxiliary problem, using the trick described in Section~\ref{sec:dimension_reduction_trick}, let $\|\vu\|_2 = s$. Then for any $\|\bbeta_{J_\star^c}\|_2 \in I$, the optimal value of the inner maximization problem is equivalent to $\max_{0 \le s \le H} ( s \langle \vh, \bbeta_{J_\star^c} \rangle + s \| \bxi + \|\bbeta_{J_\star^c}\|_2 \vg \|_2 - \frac{1}{2}s^2 + N\lambda \|\bbeta_{J_\star^c}\|_1 - \frac{1}{2}\|\bxi\|_2^2 )$. Next, considering the outer minimization problem, by Lemma~\ref{lemma:squared_hinge}, we have $\phi_I^{K,H} = \min_{t \in I} \max_{0 \le s \le H} \Psi(t, s)$.

We next estimate the saddle point $(t^\circ, s^\circ)$ of $\phi_I^{K,H} = \min_{t \in I} \max_{0 \le s \le H} \Psi(t, s)$. We first examine the estimate of $s^\circ$. The idea of the proof here is to use the $1$-strong concavity of $s \mapsto \min_{0 \le t \le K} \Psi(t, s)$ (which is easy to verify). Therefore, if we can construct $s^- < s^+$ such that $\partial_s \min_{0 \le t \le K} \Psi(t, s^-) > 0$ and $\partial_s \min_{0 \le t \le K} \Psi(t, s^+) < 0$, then we must have $s^- \le s^\circ \le s^+$. By Danskin's theorem, see, for instance, \cite[Theorem 4.13, pp. 273 - 274]{bonnansPerturbationAnalysisOptimization2000},
\begin{align*}
    \partial_s\min_{0\leq t\leq K}\Psi(t,s) = \partial_s\Psi(t(s),s) = \|\bxi+t(s)\vg\|_2 - s - t(s)\begin{cases}
        \frac{A_\vh(\frac{N\lambda}{s})}{B_\vh(\frac{N\lambda}{s})^{\frac{1}{2}}}, & s\|\vh\|_\infty>N\lambda,\\
        \|\vh\|_\infty, & \mbox{otherwise},
    \end{cases}
\end{align*}where $t(s) = \argmin(\Psi(t,s):0\leq t\leq K)$.
Next, we take $s^- = \frac{N\lambda}{\overline{z}_q}$ and $s^+ = \frac{N\lambda}{\underline{z}_q}$. Clearly $s^- < s^+$. We now verify that $s^- \|\vh\|_\infty > N\lambda$. On $\Omega_{\overline{z}_q,2}^{\mathrm{aux}}$, $B_\vh(\overline{z}_q)>0$, and therefore $\|\vh\|_\infty>\overline{z}_q=\frac{N\lambda}{s^-}$. Hence $s^-\|\vh\|_\infty>N\lambda$, and the same holds at $s^+$ since $s^-<s^+$. 

We are left to compute the partial derivatives at these two points. We first compute $t(s) = \arg\min(\Psi(t, s) : 0 \le t \le K)$. For $s\in(s^-,s^+)$, writing $u=\frac{N\lambda}{s}\in(\underline{z}_q,\overline{z}_q)$ and using $s\|\vh\|_\infty>N\lambda$, the $t$-dependent part of $\Psi(t,s)$ is $s\|\bxi+t\vg\|_2-tsB_\vh(u)^{\frac{1}{2}}$. Hence $\partial_t\Psi(t,s)=s\frac{\langle\bxi,\vg\rangle+t\|\vg\|_2^2}{\|\bxi+t\vg\|_2}-sB_\vh(u)^{\frac{1}{2}}$. By the first-order condition  $\frac{\langle\bxi,\vg\rangle+t(s)\|\vg\|_2^2}{\|\bxi+t(s)\vg\|_2}=B_\vh(u)^{\frac{1}{2}}$, which gives
\begin{align}\label{eq:solution_t_s}
t(s)=\frac{\|\bxi\|_2}{\|\vg\|_2^2}\bigg[B_\vh(u)^{\frac{1}{2}}\frac{\big(\|\vg\|_2^2-\frac{\langle\bxi,\vg\rangle^2}{\|\bxi\|_2^2}\big)^{\frac{1}{2}}}{\big(\|\vg\|_2^2-B_\vh(u)\big)^{\frac{1}{2}}}-\frac{\langle\bxi,\vg\rangle}{\|\bxi\|_2}\bigg].
\end{align}We then verify that $0 < t(s) < K$ and that $t(s)$ is unique. One may compute that
\begin{align}\label{eq:derivative_Hessian_Psi}
\partial_t^2\Psi(t,s)=s\frac{\|\bxi\|_2^2\|\vg\|_2^2-\langle\bxi,\vg\rangle^2}{\|\bxi+t\vg\|_2^3}>0.
\end{align}
Moreover, Lemma~\ref{lem:common-calibration} and Lemma~\ref{lemma:null_signal_events} give $0<B_\vh(u)<\|\vg\|_2^2$, hence $\partial_{t+}\Psi(0,s)<0$, where $\partial_{t+}$ denotes the right derivative, hence $t(s)>0$ and is unique. Moreover, at $u=\overline{z}_q$ and $u=\underline{z}_q$, respectively, the preceding formula for $t(s)$, together with Lemma~\ref{lem:common-calibration} and Lemma~\ref{lemma:null_signal_events}, yields $0<t(s^-),t(s^+)$. Moreover, on $\Omega^{\mathrm{aux}}$, $\|\bxi\|_2\leq2\sigma_\xi\sqrt N$, $\|\vg\|_2^2\geq\frac{3N}{4}$, $B_\vh(u)\leq2^5d_q$, $\|\vg\|_2^2-B_\vh(u)\geq \frac{N}{3}$, and $\frac{|\langle\bxi,\vg\rangle|}{\|\bxi\|_2}\leq2^{-23}\sqrt{d_q}$. Substituting these estimates into the formula for $t(s)$ gives $t(s^-),t(s^+)\leq2^5\sigma_\xi\sqrt{\frac{d_q}{N}}\leq2^5\sqrt{\Var_{J_\star^c}}\leq2^{5}\sqrt{\Var_{J_\star^c}}<K$, where we used $\Var_{J_\star^c}=\frac{\sigma_\xi^2d_q}{N-d_q}$ and $K>2^{7}\sqrt{\Var_{J_\star^c}}$. Hence $t(s^-)$ and $t(s^+)$, which are the unique unconstrained minimizers of $x\mapsto\Psi(x,s^-)$ and $x\mapsto\Psi(x,s^+)$, belong to $(0,K)$ and are therefore also their unique minimizers over $[0,K]$. Consequently, they are precisely the minimizers appearing in Danskin's formula above.

Substituting the preceding formula for $t(s)$ into $\partial_s\Psi(t(s),s)$, and recalling that $u=\frac{N\lambda}{s}$, we first obtain
\begin{align*}
\begin{aligned}
&\|\bxi+t(s)\vg\|_2^2
=\|\bxi\|_2^2+2t(s)\langle\bxi,\vg\rangle+t(s)^2\|\vg\|_2^2\\
&=\|\bxi\|_2^2+\frac{2\langle\bxi,\vg\rangle}{\|\vg\|_2^2}\bigg[\frac{B_\vh(u)^{\frac{1}{2}}\big(\|\bxi\|_2^2\|\vg\|_2^2-\langle\bxi,\vg\rangle^2\big)^{\frac{1}{2}}}{\big(\|\vg\|_2^2-B_\vh(u)\big)^{\frac{1}{2}}}-\langle\bxi,\vg\rangle\bigg]+\frac{1}{\|\vg\|_2^2}\bigg[\frac{B_\vh(u)^{\frac{1}{2}}\big(\|\bxi\|_2^2\|\vg\|_2^2-\langle\bxi,\vg\rangle^2\big)^{\frac{1}{2}}}{\big(\|\vg\|_2^2-B_\vh(u)\big)^{\frac{1}{2}}}-\langle\bxi,\vg\rangle\bigg]^2\\
&=\|\bxi\|_2^2+\frac{2\langle\bxi,\vg\rangle}{\|\vg\|_2^2}\bigg[\frac{B_\vh(u)^{\frac{1}{2}}\big(\|\bxi\|_2^2\|\vg\|_2^2-\langle\bxi,\vg\rangle^2\big)^{\frac{1}{2}}}{\big(\|\vg\|_2^2-B_\vh(u)\big)^{\frac{1}{2}}}-\langle\bxi,\vg\rangle\bigg]\\
&+\frac{1}{\|\vg\|_2^2}\bigg[\frac{B_\vh(u)\big(\|\bxi\|_2^2\|\vg\|_2^2-\langle\bxi,\vg\rangle^2\big)}{\|\vg\|_2^2-B_\vh(u)}-\frac{2B_\vh(u)^{\frac{1}{2}}\langle\bxi,\vg\rangle\big(\|\bxi\|_2^2\|\vg\|_2^2-\langle\bxi,\vg\rangle^2\big)^{\frac{1}{2}}}{\big(\|\vg\|_2^2-B_\vh(u)\big)^{\frac{1}{2}}}+\langle\bxi,\vg\rangle^2\bigg]\\
&=\|\bxi\|_2^2-\frac{\langle\bxi,\vg\rangle^2}{\|\vg\|_2^2}+\frac{B_\vh(u)}{\|\vg\|_2^2}\frac{\|\bxi\|_2^2\|\vg\|_2^2-\langle\bxi,\vg\rangle^2}{\|\vg\|_2^2-B_\vh(u)} =\frac{\|\bxi\|_2^2\big(\|\vg\|_2^2-\frac{\langle\bxi,\vg\rangle^2}{\|\bxi\|_2^2}\big)}{\|\vg\|_2^2-B_\vh(u)}.
\end{aligned}
\end{align*}
Using $u=\frac{N\lambda}{s}$, and \eqref{eq:solution_t_s}, there holds
\begin{align*}
\begin{aligned}
&\frac{u}{\|\bxi\|_2}\partial_s\Psi(t(s),s)
=\frac{u}{\|\bxi\|_2}\bigg(\|\bxi+t(s)\vg\|_2-s-t(s)\frac{A_\vh(u)}{B_\vh(u)^{\frac{1}{2}}}\bigg) =\frac{u\big(\|\vg\|_2^2-\frac{\langle\bxi,\vg\rangle^2}{\|\bxi\|_2^2}\big)^{\frac{1}{2}}}{\big(\|\vg\|_2^2-B_\vh(u)\big)^{\frac{1}{2}}}-\frac{N\lambda}{\|\bxi\|_2}-\frac{ut(s)A_\vh(u)}{\|\bxi\|_2B_\vh(u)^{\frac{1}{2}}}\\
&=\frac{u\big(\|\vg\|_2^2-\frac{\langle\bxi,\vg\rangle^2}{\|\bxi\|_2^2}\big)^{\frac{1}{2}}}{\big(\|\vg\|_2^2-B_\vh(u)\big)^{\frac{1}{2}}}-\frac{N\lambda}{\|\bxi\|_2} -\frac{uA_\vh(u)}{\|\vg\|_2^2B_\vh(u)^{\frac{1}{2}}}\bigg[B_\vh(u)^{\frac{1}{2}}\frac{\big(\|\vg\|_2^2-\frac{\langle\bxi,\vg\rangle^2}{\|\bxi\|_2^2}\big)^{\frac{1}{2}}}{\big(\|\vg\|_2^2-B_\vh(u)\big)^{\frac{1}{2}}}-\frac{\langle\bxi,\vg\rangle}{\|\bxi\|_2}\bigg]\\
&=\frac{u\big(\|\vg\|_2^2-\frac{\langle\bxi,\vg\rangle^2}{\|\bxi\|_2^2}\big)^{\frac{1}{2}}}{\big(\|\vg\|_2^2-B_\vh(u)\big)^{\frac{1}{2}}}-\frac{uA_\vh(u)}{\|\vg\|_2^2}\frac{\big(\|\vg\|_2^2-\frac{\langle\bxi,\vg\rangle^2}{\|\bxi\|_2^2}\big)^{\frac{1}{2}}}{\big(\|\vg\|_2^2-B_\vh(u)\big)^{\frac{1}{2}}}+\frac{uA_\vh(u)}{\|\vg\|_2^2B_\vh(u)^{\frac{1}{2}}}\frac{\langle\bxi,\vg\rangle}{\|\bxi\|_2}-\frac{N\lambda}{\|\bxi\|_2}\\
&=\frac{u}{\|\vg\|_2^2}\bigg[\frac{\big(\|\vg\|_2^2-\frac{\langle\bxi,\vg\rangle^2}{\|\bxi\|_2^2}\big)^{\frac{1}{2}}\big(\|\vg\|_2^2-A_\vh(u)\big)}{\big(\|\vg\|_2^2-B_\vh(u)\big)^{\frac{1}{2}}}+\frac{\langle\bxi,\vg\rangle}{\|\bxi\|_2}\frac{A_\vh(u)}{B_\vh(u)^{\frac{1}{2}}}\bigg]-\frac{N\lambda}{\|\bxi\|_2}.
\end{aligned}
\end{align*}
For $u\in(\underline{z}_q,\overline{z}_q)$, the events defining $\Omega^{\mathrm{aux}}$ give $\big|\frac{\|\vg\|_2^2}{N}-1\big|\leq2^{-23}\sqrt{\frac{d_q}{N}}$, $\frac{|\langle\bxi,\vg\rangle|}{\|\bxi\|_2}\leq2^{-23}\sqrt{d_q}$, 
\begin{align*}
    \big|\frac{A_\vh(u)}{2|J_\star^c|\Phi_G(u)}-1\big|\leq2^{-23},\mbox{ and }\big|(\frac{B_\vh(u)}{[|J_\star^c|\bE(|G|-u)_+^2]})^{\frac{1}{2}}-1\big|\leq2^{-23}.
\end{align*}
Moreover, item~\emph{\ref{item:Gaussian_local_comparison}} of Section~\ref{sec:Gaussian_calculations} and \eqref{eq:expectation_soft_max} give $2|J_\star^c|\Phi_G(u)\leq e^{\frac{1}{8}}q$, $|J_\star^c|\bE[(|G|-u)_+^2]\leq 2^4d_q$, and $|J_\star^c|\bE(|G|-u)_+^2\leq\frac{4|J_\star^c|\Phi_G(u)}{u^2}\leq\frac{2e^{\frac{1}{8}}N}{u^2}$. Hence
\begin{align*}
\begin{aligned}
B_\vh(u)&\leq(1+2^{-23})^2|J_\star^c|\bE(|G|-u)_+^2<\frac{3N}{8},\mbox{ and }\|\vg\|_2^2\geq N\bigg(1-2^{-23}\sqrt{\frac{d_q}{N}}\bigg)>\frac{3N}{4},
\end{aligned}
\end{align*}
and consequently $\|\vg\|_2^2-B_\vh(u)>\frac{3N}{8}>\frac{N}{3}$. In particular, $\frac{B_\vh(u)}{\|\vg\|_2^2}<\frac{1}{2}$ and $\frac{|J_\star^c|\bE(|G|-u)_+^2}{N}<\frac{1}{2}$.

We next compare
\begin{align*}
\begin{aligned}
&\frac{u}{\|\vg\|_2^2}\bigg[\frac{\big(\|\vg\|_2^2-\frac{\langle\bxi,\vg\rangle^2}{\|\bxi\|_2^2}\big)^{\frac{1}{2}}\big(\|\vg\|_2^2-A_\vh(u)\big)}{\big(\|\vg\|_2^2-B_\vh(u)\big)^{\frac{1}{2}}}+\frac{\langle\bxi,\vg\rangle}{\|\bxi\|_2}\frac{A_\vh(u)}{B_\vh(u)^{\frac{1}{2}}}\bigg]\\
&=\frac{u}{\|\vg\|_2^2}\bigg[\frac{\|\vg\|_2\big(1-\frac{\langle\bxi,\vg\rangle^2}{\|\bxi\|_2^2\|\vg\|_2^2}\big)^{\frac{1}{2}}\|\vg\|_2^2\big(1-\frac{A_\vh(u)}{\|\vg\|_2^2}\big)}{\|\vg\|_2\big(1-\frac{B_\vh(u)}{\|\vg\|_2^2}\big)^{\frac{1}{2}}}+\frac{\langle\bxi,\vg\rangle}{\|\bxi\|_2}\frac{A_\vh(u)}{B_\vh(u)^{\frac{1}{2}}}\bigg]\\
&=u\frac{\big(1-\frac{\langle\bxi,\vg\rangle^2}{\|\bxi\|_2^2\|\vg\|_2^2}\big)^{\frac{1}{2}}\big(1-\frac{A_\vh(u)}{\|\vg\|_2^2}\big)}{\big(1-\frac{B_\vh(u)}{\|\vg\|_2^2}\big)^{\frac{1}{2}}}+\frac{u\langle\bxi,\vg\rangle A_\vh(u)}{\|\bxi\|_2\|\vg\|_2^2B_\vh(u)^{\frac{1}{2}}}.
\end{aligned}
\end{align*}
with
\begin{align*}
\frac{u[N-2|J_\star^c|\Phi_G(u)]}{\sqrt N(N-|J_\star^c|\bE(|G|-u)_+^2)^{\frac{1}{2}}}=u\frac{1-\frac{2|J_\star^c|\Phi_G(u)}{N}}{(1-\frac{|J_\star^c|\bE(|G|-u)_+^2}{N})^{\frac{1}{2}}}.
\end{align*}
The preceding estimates imply
\begin{align*}
\begin{aligned}
\bigg|\frac{A_\vh(u)}{\|\vg\|_2^2}-\frac{2|J_\star^c|\Phi_G(u)}{N}\bigg|
&\leq\frac{|A_\vh(u)-2|J_\star^c|\Phi_G(u)|}{\|\vg\|_2^2}+2|J_\star^c|\Phi_G(u)\bigg|\frac1{\|\vg\|_2^2}-\frac1N\bigg|\\
&\leq\frac{4}{3N}2^{-23}2|J_\star^c|\Phi_G(u)+\frac{4}{3}\frac{2|J_\star^c|\Phi_G(u)}{N}\bigg|\frac{\|\vg\|_2^2}{N}-1\bigg|
\leq2^{-21}\frac{q}{N},
\end{aligned}
\end{align*}
and, since $\big|\frac{B_\vh(u)}{[|J_\star^c|\bE(|G|-u)_+^2]}-1\big|\leq2^{1-23}+2^{-46}$,
\begin{align*}
\begin{aligned}
&\bigg|\frac{B_\vh(u)}{\|\vg\|_2^2}-\frac{|J_\star^c|\bE(|G|-u)_+^2}{N}\bigg|
\leq\frac{|B_\vh(u)-|J_\star^c|\bE(|G|-u)_+^2|}{\|\vg\|_2^2}+|J_\star^c|\bE(|G|-u)_+^2\bigg|\frac1{\|\vg\|_2^2}-\frac1N\bigg|\\
&\leq\frac{4}{3N}(2^{1-23}+2^{-46})|J_\star^c|\bE(|G|-u)_+^2+\frac{4}{3}\frac{|J_\star^c|\bE(|G|-u)_+^2}{N}\bigg|\frac{\|\vg\|_2^2}{N}-1\bigg|\\
&\leq\frac{4}{3N}(2^{1-23}+2^{-46})2^4d_q+\frac{4}{3}\frac{2^4d_q}{N}2^{-23}\sqrt{\frac{d_q}{N}}
\leq2^{-17}\frac{d_q}{N}.
\end{aligned}
\end{align*}
Also, using $1-\sqrt{1-x}=\frac{x}{1+\sqrt{1-x}}\leq x$ for every $x\in[0,1]$,
\begin{align*}
0\leq1-\bigg(1-\frac{\langle\bxi,\vg\rangle^2}{\|\bxi\|_2^2\|\vg\|_2^2}\bigg)^{\frac{1}{2}}\leq\frac{\langle\bxi,\vg\rangle^2}{\|\bxi\|_2^2\|\vg\|_2^2}\leq2^{-45}\frac{d_q}{N}.
\end{align*}
Since $x\mapsto(1-x)^{-\frac{1}{2}}$ has derivative at most $2$ on $[0,\frac{1}{2}]$, the preceding three estimates give
\begin{align*}
\begin{aligned}
&\bigg|u\frac{\big(1-\frac{\langle\bxi,\vg\rangle^2}{\|\bxi\|_2^2\|\vg\|_2^2}\big)^{\frac{1}{2}}\big(1-\frac{A_\vh(u)}{\|\vg\|_2^2}\big)}{\big(1-\frac{B_\vh(u)}{\|\vg\|_2^2}\big)^{\frac{1}{2}}}-u\frac{1-\frac{2|J_\star^c|\Phi_G(u)}{N}}{(1-\frac{|J_\star^c|\bE(|G|-u)_+^2}{N})^{\frac{1}{2}}}\bigg|\\
&\leq2u\bigg[\bigg|1-\bigg(1-\frac{\langle\bxi,\vg\rangle^2}{\|\bxi\|_2^2\|\vg\|_2^2}\bigg)^{\frac{1}{2}}\bigg|+\bigg|\frac{A_\vh(u)}{\|\vg\|_2^2}-\frac{2|J_\star^c|\Phi_G(u)}{N}\bigg|+2\bigg|\frac{B_\vh(u)}{\|\vg\|_2^2}-\frac{|J_\star^c|\bE(|G|-u)_+^2}{N}\bigg|\bigg]\\
&\leq2u\bigg[2^{-45}\frac{d_q}{N}+2^{-21}\frac{q}{N}+2\cdot2^{-17}\frac{d_q}{N}\bigg]\leq2^{-14}u\bigg(\frac{q}{N}+\frac{d_q}{N}\bigg).
\end{aligned}
\end{align*}
Moreover, $B_\vh(u)^{\frac{1}{2}}\geq(1-2^{-23})[|J_\star^c|\bE(|G|-u)_+^2]^{\frac{1}{2}}>2^{-3}\sqrt{d_q}$, while $A_\vh(u)\leq(1+2^{-23})2|J_\star^c|\Phi_G(u)\leq2q$ and $\|\vg\|_2^2\geq\frac{3N}{4}$. Hence
\begin{align*}
\bigg|\frac{u\langle\bxi,\vg\rangle A_\vh(u)}{\|\bxi\|_2\|\vg\|_2^2B_\vh(u)^{\frac{1}{2}}}\bigg|\leq u\,2^{-23}\sqrt{d_q}\frac{2q}{(\frac{3N}{4})2^{-3}\sqrt{d_q}}\leq2^{-18}u\frac{q}{N}.
\end{align*}
Combining the preceding two displays, and using $u\leq\overline{z}_q<2z_q$ and $d_q\leq\frac{2q}{z_q^2}$, we obtain
\begin{align*}
\begin{aligned}
&\bigg|\frac{u}{\|\vg\|_2^2}\bigg[\frac{\big(\|\vg\|_2^2-\frac{\langle\bxi,\vg\rangle^2}{\|\bxi\|_2^2}\big)^{\frac{1}{2}}\big(\|\vg\|_2^2-A_\vh(u)\big)}{\big(\|\vg\|_2^2-B_\vh(u)\big)^{\frac{1}{2}}}+\frac{\langle\bxi,\vg\rangle}{\|\bxi\|_2}\frac{A_\vh(u)}{B_\vh(u)^{\frac{1}{2}}}\bigg]-\frac{u[N-2|J_\star^c|\Phi_G(u)]}{\sqrt N(N-|J_\star^c|\bE(|G|-u)_+^2)^{\frac{1}{2}}}\bigg|\\
&\leq2^{-13}u\bigg(\frac{q}{N}+\frac{d_q}{N}\bigg)\leq2^{-11}\bigg(\frac1{z_q}+\frac{qz_q}{N}\bigg).
\end{aligned}
\end{align*}

On the other hand, $\Omega_\xi$ gives $\big|\frac{\|\bxi\|_2}{\sigma_\xi\sqrt N}-1\big|\leq2^{-23}\sqrt{\frac{d_q}{N}}<\frac{1}{2}$. Hence
\begin{align*}
\begin{aligned}
\bigg|\frac{N\lambda}{\|\bxi\|_2}-\frac{\sqrt N\lambda}{\sigma_\xi}\bigg|
&=\frac{\sqrt N\lambda}{\sigma_\xi}\bigg|\frac{\sigma_\xi\sqrt N}{\|\bxi\|_2}-1\bigg|\leq2^{1-23}\sqrt{\frac{d_q}{N}}\frac{\sqrt N\lambda}{\sigma_\xi}\\
&=2^{1-23}\sqrt{\frac{d_q}{N}}\frac{z_q(1-\frac{q}{N})}{(1-\frac{d_q}{N})^{\frac{1}{2}}}\leq2^{-21}z_q\sqrt{\frac{d_q}{N}}\leq2^{-20}\bigg(\frac1{z_q}+\frac{qz_q}{N}\bigg),
\end{aligned}
\end{align*}
where we used \eqref{eq:common-q-calibration}, $\frac{d_q}{N}<0.31$, $d_q\leq\frac{2q}{z_q^2}$, and $z_q^{-1}+\frac{qz_q}{N}\geq2\sqrt{\frac{q}{N}}$.

At $s=s^-$, $u=\overline{z}_q$, so \eqref{eq:common-endpoint-calibration} and the preceding two estimates yield
\begin{align*}
\frac{\overline{z}_q}{\|\bxi\|_2}\partial_s\Psi(t(s^-),s^-)\geq\big(2^{-10}-2^{-11}-2^{-20}\big)\bigg(\frac1{z_q}+\frac{qz_q}{N}\bigg)>0.
\end{align*}
At $s=s^+$, $u=\underline{z}_q$, so \eqref{eq:common-endpoint-calibration} gives
\begin{align*}
\frac{\underline{z}_q}{\|\bxi\|_2}\partial_s\Psi(t(s^+),s^+)\leq-\big(2^{-10}-2^{-11}-2^{-20}\big)\bigg(\frac1{z_q}+\frac{qz_q}{N}\bigg)<0.
\end{align*}
Since $\frac{\overline{z}_q}{\|\bxi\|_2}>0$ and $\frac{\underline{z}_q}{\|\bxi\|_2}>0$, Danskin's formula therefore gives
\begin{align*}
\bigg[\partial_s\min_{0\leq x\leq K}\Psi(x,s)\bigg]_{s=s^-}=\partial_s\Psi(t(s^-),s^-)>0,\qquad \bigg[\partial_s\min_{0\leq x\leq K}\Psi(x,s)\bigg]_{s=s^+}=\partial_s\Psi(t(s^+),s^+)<0.
\end{align*}
The $1$-strong concavity of $s\mapsto\min_{0\leq x\leq K}\Psi(x,s)$ now yields $s^-\leq s^\circ\leq s^+$, or equivalently $\underline{z}_q\leq \frac{N\lambda}{s^\circ}\leq\overline{z}_q$. In particular, $0<s^\circ<H$ and $s^\circ\|\vh\|_\infty>N\lambda$. Since $t^\circ=t(s^\circ)$, applying the preceding formula for $t(s)$ uniformly for $\frac{N\lambda}{s}\in[\underline{z}_q,\overline{z}_q]$ and using Lemmas~\ref{lem:common-calibration} and~\ref{lemma:null_signal_events}, we obtain
\begin{align*}
0<t^\circ\leq2^{5}\sqrt{\Var_{J_\star^c}},\qquad \frac{2^{-4}N\lambda}{z_q}\leq s^\circ\leq\frac{2^4N\lambda}{z_q},
\end{align*}
which proves \eqref{eq:null_signal_saddle_bounds}.

\end{proof}

\begin{proof}[Proof of Lemma~\ref{lemma:null_signal_cost_separation}]
Set $r^-=2^{6}\sqrt{\Var_{J_\star^c}}$, $r^+=2^{7}\sqrt{\Var_{J_\star^c}}$, and $H^+=\frac{2^{5}N\lambda}{z_q}$.

Fix a saddle point $(t^\circ,s^\circ)$ of $\min_{0\leq t\leq K}\max_{0\leq s\leq H}\Psi(t,s)$. For $x\geq0$, from \eqref{eq:derivative_Hessian_Psi} we know that $\partial_x^2\Psi(x,s^\circ)=s^\circ\frac{\|\boldsymbol g\|_2^2\|\boldsymbol\xi\|_2^2-\langle\boldsymbol\xi,\boldsymbol g\rangle^2}{\|\boldsymbol\xi+x\boldsymbol g\|_2^3}$. Lemma~\ref{lem:common-calibration},~\eqref{eq:null_signal_saddle_bounds}, and the definition of $\Omega^{\mathrm{aux}}$ imply, for $0\leq x\leq r^-$, $\partial_x^2\Psi(x,s^\circ)\geq2^{-21}\frac{N^{\frac{3}{2}}\lambda}{\sigma_\xi z_q}$. Since \eqref{eq:null_signal_saddle_bounds} gives $0<t^\circ<r^-$ and $t^\circ$ minimizes $x\mapsto\Psi(x,s^\circ)$, $\partial_x\Psi(t^\circ,s^\circ)=0$. Therefore,
\begin{align*}
\forall t^\circ\leq y\leq r^-,\qquad \partial_y\Psi(y,s^\circ)=\int_{t^\circ}^y\partial_x^2\Psi(x,s^\circ)\,\mathrm dx\geq2^{-21}\frac{N^{\frac{3}{2}}\lambda}{\sigma_\xi z_q}(y-t^\circ),
\end{align*}
and integrating once more gives
\begin{align*}
\Psi(r^-,s^\circ)-\Psi(t^\circ,s^\circ)=\int_{t^\circ}^{r^-}\partial_y\Psi(y,s^\circ)\,\mathrm dy\geq2^{-22}\frac{N^{\frac{3}{2}}\lambda}{\sigma_\xi z_q}(r^--t^\circ)^2\geq2^{-12}\frac{N^{\frac{3}{2}}\lambda\Var_{J_\star^c}}{\sigma_\xi z_q}.
\end{align*}
Indeed, $t^\circ\leq2^{5}\sqrt{\Var_{J_\star^c}}$ gives $(r^--t^\circ)^2\geq2^{10}\Var_{J_\star^c}$. Since $x\mapsto\Psi(x,s^\circ)$ is convex and minimized at $t^\circ<r^-$, it is nondecreasing on $[r^-,K]$. Moreover, the saddle inequalities are $\Psi(t^\circ,s)\leq\Psi(t^\circ,s^\circ)\leq\Psi(x,s^\circ)$ for every $(x,s)\in[0,K]\times[0,H]$. Hence
\begin{align*}
\phi_{[r^-,K]}^{K,H}=\min_{r^-\leq x\leq K}\max_{0\leq s\leq H}\Psi(x,s)\geq\min_{r^-\leq x\leq K}\Psi(x,s^\circ)=\Psi(r^-,s^\circ),\qquad \phi_{[0,K]}^{K,H}=\Psi(t^\circ,s^\circ),
\end{align*}
and consequently
\begin{align}
\phi_{[r^-,K]}^{K,H}-\phi_{[0,K]}^{K,H}\geq2^{-12}\frac{N^{\frac{3}{2}}\lambda\Var_{J_\star^c}}{\sigma_\xi z_q}.
\label{eq:restricted_cost_gap}
\end{align}

We next record the Lipschitz constants of the two optimal values used below. Recalling the definition of $\Psi$,
\begin{align}
\Psi(t,s):=\Psi_{\vg,\vh}(t,s)=s\|\bxi+t\vg\|_2-\frac{s^2}{2}-\frac{\|\bxi\|_2^2}{2}+t\begin{cases}-sB_\vh(\frac{N\lambda}{s})^{\frac{1}{2}},&s\|\vh\|_\infty>N\lambda,\\N\lambda-s\|\vh\|_\infty,&s\|\vh\|_\infty\leq N\lambda,\end{cases}
\tag{\ref{eq:def_Psi}}
\end{align}
where its $\vh$-dependent term is equivalently $\min_{\|\bbeta_{J_\star^c}\|_2=t}\{s\langle\vh,\bbeta_{J_\star^c}\rangle+N\lambda\|\bbeta_{J_\star^c}\|_1\}$. Hence, for every $t,s\geq0$,
\begin{align*}
\big|s\|\bxi+t\vg\|_2-s\|\bxi+t\vg'\|_2\big|\leq st\|\vg-\vg'\|_2,
\end{align*}
while
\begin{align}\tag{\ref{eq:Lipschitz_s_mapsto_min}}
\bigg|\min_{\|\bbeta_{J_\star^c}\|_2=t}\{s\langle\vh,\bbeta_{J_\star^c}\rangle+N\lambda\|\bbeta_{J_\star^c}\|_1\}-\min_{\|\bbeta_{J_\star^c}\|_2=t}\{s\langle\vh',\bbeta_{J_\star^c}\rangle+N\lambda\|\bbeta_{J_\star^c}\|_1\}\bigg|\leq st\|\vh-\vh'\|_2.
\end{align}
Therefore,
\begin{align*}
|\Psi_{\vg,\vh}(t,s)-\Psi_{\vg',\vh'}(t,s)|\leq\sqrt2\,st\big(\|\vg-\vg'\|_2^2+\|\vh-\vh'\|_2^2\big)^{\frac{1}{2}}.
\end{align*}
By Lemma~\ref{lemma:Lipschitz_convex_concave_function}, there holds
\begin{align}
\max\bigg\{\operatorname{Lip}_{\boldsymbol g,\boldsymbol h}\big(\phi_{[0,r^+]}^{r^+,H^+}\big),\operatorname{Lip}_{\boldsymbol g,\boldsymbol h}\big(\phi_{(r^-)}^{r^+,H^+}\big)\bigg\}\leq\sqrt2\,r^+H^+\leq2^{13}\frac{N\lambda\sqrt{\Var_{J_\star^c}}}{z_q}.
\label{eq:capped_lipschitz}
\end{align}

Applying~\eqref{eq:restricted_cost_gap} with $K=r^+$ and $H=H^+$ gives $\phi_{[r^-,r^+]}^{r^+,H^+}-\phi_{[0,r^+]}^{r^+,H^+}\geq2^{-12}\frac{N^{\frac{3}{2}}\lambda\Var_{J_\star^c}}{\sigma_\xi z_q}$. Since $(r^-)\subset[r^-,r^+]$, minimization over the singleton can only increase the optimal value:
\begin{align*}
\phi_{(r^-)}^{r^+,H^+}=\max_{0\leq s\leq H^+}\Psi(r^-,s)\geq\min_{r^-\leq t\leq r^+}\max_{0\leq s\leq H^+}\Psi(t,s)=\phi_{[r^-,r^+]}^{r^+,H^+}.
\end{align*}
Consequently,
\begin{align}
\phi_{(r^-)}^{r^+,H^+}-\phi_{[0,r^+]}^{r^+,H^+}\geq2^{-12}\frac{N^{\frac{3}{2}}\lambda\Var_{J_\star^c}}{\sigma_\xi z_q}.
\label{eq:capped_cost_gap}
\end{align}

Gaussian isoperimetry inequality, see, for instance, \cite[Theorem 10.17]{boucheronbook}, and~\eqref{eq:capped_lipschitz} give, for every $x>0$,
\begin{align*}
&\max\bigg\{\mathbb P_{\mathrm{aux}}\bigg\{\bigg|\phi_{[0,r^+]}^{r^+,H^+}-\operatorname{Med}_{\boldsymbol g,\boldsymbol h}\big(\phi_{[0,r^+]}^{r^+,H^+}\mid\boldsymbol\xi\big)\bigg|>x\,\bigg|\,\boldsymbol\xi\bigg\},\mathbb P_{\mathrm{aux}}\bigg\{\bigg|\phi_{(r^-)}^{r^+,H^+}-\operatorname{Med}_{\boldsymbol g,\boldsymbol h}\big(\phi_{(r^-)}^{r^+,H^+}\mid\boldsymbol\xi\big)\bigg|>x\,\bigg|\,\boldsymbol\xi\bigg\}\bigg\}\\
&\leq2\exp\bigg(-\frac{x^2z_q^2}{2^{27}N^2\lambda^2\Var_{J_\star^c}}\bigg),
\end{align*}
where $\operatorname{Med}_{\boldsymbol g,\boldsymbol h}(\,\cdot\mid\boldsymbol\xi)$ denotes a conditional median. At $x=\frac{2^{-15}N^{\frac{3}{2}}\lambda\Var_{J_\star^c}}{\sigma_\xi z_q}$, the exponent is at least $2^{-57}d_q$. Hence, by the union bound, \eqref{eq:omega_probability}, and $C_1e^{-c_0d_q}<1$,
\begin{align*}
&\mathbb P_{\mathrm{aux}}\bigg\{\big(\Omega^{\mathrm{aux}}\big)^c\cup\bigg\{\bigg|\phi_{[0,r^+]}^{r^+,H^+}-\operatorname{Med}_{\boldsymbol g,\boldsymbol h}\big(\phi_{[0,r^+]}^{r^+,H^+}\mid\boldsymbol\xi\big)\bigg|>x\bigg\} \cup\bigg\{\bigg|\phi_{(r^-)}^{r^+,H^+}-\operatorname{Med}_{\boldsymbol g,\boldsymbol h}\big(\phi_{(r^-)}^{r^+,H^+}\mid\boldsymbol\xi\big)\bigg|>x\bigg\}\,\bigg|\,\boldsymbol\xi\bigg\}\\
&\leq\mathbb P_{\mathrm{aux}}(\big(\Omega^{\mathrm{aux}}\big)^c\mid\boldsymbol\xi)+4e^{-2^{-57}d_q}<1.
\end{align*}
Consequently, there exists $(\vg,\vh)\in\Omega^{\mathrm{aux}}$ for which both deviations in the preceding display are at most $x$. At this realization, \eqref{eq:capped_cost_gap} gives
\begin{align*}
&\operatorname{Med}_{\boldsymbol g,\boldsymbol h}\big(\phi_{(r^-)}^{r^+,H^+}\mid\boldsymbol\xi\big)\geq\phi_{(r^-)}^{r^+,H^+}-x\geq\phi_{[0,r^+]}^{r^+,H^+}+2^{-12}\frac{N^{\frac{3}{2}}\lambda\Var_{J_\star^c}}{\sigma_\xi z_q}-x \geq\operatorname{Med}_{\boldsymbol g,\boldsymbol h}\big(\phi_{[0,r^+]}^{r^+,H^+}\mid\boldsymbol\xi\big)+2^{-12}\frac{N^{\frac{3}{2}}\lambda\Var_{J_\star^c}}{\sigma_\xi z_q}-2x,
\end{align*}
and therefore
\begin{align}
\operatorname{Med}_{\boldsymbol g,\boldsymbol h}\big(\phi_{(r^-)}^{r^+,H^+}\mid\boldsymbol\xi\big)-\operatorname{Med}_{\boldsymbol g,\boldsymbol h}\big(\phi_{[0,r^+]}^{r^+,H^+}\mid\boldsymbol\xi\big)\geq2^{-13}\frac{N^{\frac{3}{2}}\lambda\Var_{J_\star^c}}{\sigma_\xi z_q}.
\label{eq:median_gap}
\end{align}
Define
\begin{align*}
\alpha_1(\boldsymbol\xi)=\frac34\operatorname{Med}_{\boldsymbol g,\boldsymbol h}\big(\phi_{[0,r^+]}^{r^+,H^+}\mid\boldsymbol\xi\big)+\frac14\operatorname{Med}_{\boldsymbol g,\boldsymbol h}\big(\phi_{(r^-)}^{r^+,H^+}\mid\boldsymbol\xi\big),\mbox{ and }\gamma_1(\boldsymbol\xi)=\frac14\operatorname{Med}_{\boldsymbol g,\boldsymbol h}\big(\phi_{[0,r^+]}^{r^+,H^+}\mid\boldsymbol\xi\big)+\frac34\operatorname{Med}_{\boldsymbol g,\boldsymbol h}\big(\phi_{(r^-)}^{r^+,H^+}\mid\boldsymbol\xi\big).
\end{align*}
Then~\eqref{eq:median_gap} gives $\alpha_1(\boldsymbol\xi)<\gamma_1(\boldsymbol\xi)$ and
\begin{align*}
&\alpha_1(\boldsymbol\xi)-\operatorname{Med}_{\boldsymbol g,\boldsymbol h}\big(\phi_{[0,r^+]}^{r^+,H^+}\mid\boldsymbol\xi\big)=\operatorname{Med}_{\boldsymbol g,\boldsymbol h}\big(\phi_{(r^-)}^{r^+,H^+}\mid\boldsymbol\xi\big)-\gamma_1(\boldsymbol\xi)\\
&=\frac14\bigg[\operatorname{Med}_{\boldsymbol g,\boldsymbol h}\big(\phi_{(r^-)}^{r^+,H^+}\mid\boldsymbol\xi\big)-\operatorname{Med}_{\boldsymbol g,\boldsymbol h}\big(\phi_{[0,r^+]}^{r^+,H^+}\mid\boldsymbol\xi\big)\bigg]\geq2^{-15}\frac{N^{\frac{3}{2}}\lambda\Var_{J_\star^c}}{\sigma_\xi z_q}.
\end{align*}
Gaussian concentration with $x=\frac{2^{-15}N^{\frac{3}{2}}\lambda\Var_{J_\star^c}}{\sigma_\xi z_q}$ therefore gives
\begin{align}
\mathbb P_{\mathrm{aux}}\big\{\phi_{[0,r^+]}^{r^+,H^+}\geq\alpha_1(\boldsymbol\xi)\mid\boldsymbol\xi\big\}+\mathbb P_{\mathrm{aux}}\big\{\phi_{(r^-)}^{r^+,H^+}\leq\gamma_1(\boldsymbol\xi)\mid\boldsymbol\xi\big\}\leq4e^{-2^{-57}d_q}\leq2^{-4}C_1e^{-c_0d_q}.
\label{eq:capped_levels}
\end{align}

For arbitrary $K\geq r^+$ and $H\geq H^+$,~\eqref{eq:null_signal_saddle_bounds} gives, on $\Omega^{\mathrm{aux}}$, every saddle point $(t^\circ,s^\circ)$ of the full scalar problem satisfies
\begin{align*}
(t^\circ,s^\circ)\in\bigg[0,2^{5}\sqrt{\Var_{J_\star^c}}\bigg]\times\bigg[\frac{2^{-4}N\lambda}{z_q},\frac{2^4N\lambda}{z_q}\bigg]\subset[0,r^+]\times[0,H^+]\subset[0,K]\times[0,H].
\end{align*}
Hence the same $(t^\circ,s^\circ)$ is a saddle point after restricting the full problem to $[0,r^+]\times[0,H^+]$, and therefore $\phi_{[0,K]}^{K,H}=\phi_{[0,r^+]}^{r^+,H^+}$. Moreover, $t\mapsto\max_{0\leq s\leq H}\Psi(t,s)$ is convex and is minimized at $t^\circ<r^-$, so it is nondecreasing on $[r^-,K]$. Consequently,
\begin{align*}
\phi_{[r^-,K]}^{K,H}=\max_{0\leq s\leq H}\Psi(r^-,s)\geq\max_{0\leq s\leq H^+}\Psi(r^-,s)=\phi_{(r^-)}^{r^+,H^+}.
\end{align*}
Thus, on $\Omega^{\mathrm{aux}}$, $\{\phi_{[0,K]}^{K,H}\geq\alpha_1(\boldsymbol\xi)\}\subset\{\phi_{[0,r^+]}^{r^+,H^+}\geq\alpha_1(\boldsymbol\xi)\}$ and $\{\phi_{[r^-,K]}^{K,H}\leq\gamma_1(\boldsymbol\xi)\}\subset\{\phi_{(r^-)}^{r^+,H^+}\leq\gamma_1(\boldsymbol\xi)\}$. Combining these inclusions with~\eqref{eq:capped_levels} and paying $\mathbb P_{\mathrm{aux}}(\big(\Omega^{\mathrm{aux}}\big)^c\mid\boldsymbol\xi)$ for each event gives
\begin{align*}
\mathbb P_{\mathrm{aux}}\big\{\phi_{[0,K]}^{K,H}\geq\alpha_1(\boldsymbol\xi)\mid\boldsymbol\xi\big\}+\mathbb P_{\mathrm{aux}}\big\{\phi_{[r^-,K]}^{K,H}\leq\gamma_1(\boldsymbol\xi)\mid\boldsymbol\xi\big\}\leq(2^{-4}+2^{-3})C_1e^{-c_0d_q}\leq2^{-2}C_1e^{-c_0d_q},
\end{align*}
which is~\eqref{eq:deterministic_levels}.
\end{proof}

\subsection{End of the proof of Theorem~\ref{theo:upper_LASSO}}

\begin{proof}
Fix $\lambda\geq\frac{\sigma_\xi}{\sqrt N}$ satisfying Assumption~\ref{assumption:LASSO_upper}.

Suppose first that $\frac{\sigma_\xi}{\sqrt N}\leq\lambda\leq2^{-9}\sigma_\xi\sqrt{\frac{\log(\frac{e|J_\star^c|}{N})}{N}}.$ Lemma~\ref{lem:estimate_free_variance}, applied with $\varepsilon=2^{-18}$, gives
\begin{align}
\mathrm{Var}_{J_\star^c}\asymp\frac{\sigma_\xi^2}{\log(\frac{e|J_\star^c|}{N})}.
\label{eq:small-lambda-free-variance-comparison}
\end{align}
Moreover, the additional condition in Assumption~\ref{assumption:LASSO_upper} yields
\begin{align*}
\sigma_\xi^2\frac{|J_\star|\log(\frac{e|J_\star^c|}{N})}{N}+\lambda^2|J_\star|\lesssim c\frac{\sigma_\xi^2}{\log(\frac{e|J_\star^c|}{N})}\lesssim\mathrm{Var}_{J_\star^c}.
\end{align*}
In particular, the first term in the preceding display is of smaller order than $\mathrm{Var}_{J_\star^c}$ after the absolute constant in Assumption~\ref{assumption:LASSO_upper} is fixed sufficiently small. Since \eqref{eq:def_r_star} is an equality and $C_\star>1$, it also gives $\|\boldsymbol\beta_{J_\star^c}^\star\|_2^2\leq r_\star^2.$
Together with \eqref{eq:small-lambda-free-variance-comparison}, the definition of $r_\star$, and the choice of the absolute constants in \eqref{eq:fixed-constants}, this verifies the remaining assumptions of Proposition~\ref{prop:lambda_is_small_shift}. Let
\begin{align*}
\Omega_{\mathrm{small}}=\left\{\|\widehat{\boldsymbol\beta}_{J_\star}-\boldsymbol\beta_{J_\star}^\star\|_2^2\leq r_\star^2\right\}\cap\left\{\|\widehat{\boldsymbol\beta}_{J_\star^c}\|_2^2\leq2^{21}\bigg( \frac{\sigma_\xi^2}{\log(\frac{e|J_\star^c|}{N})} + \|\bbeta_{J_\star^c}^\star\|_2^2\bigg)\right\}.
\end{align*}
Propositions~\ref{prop:estimation_part_final} and~\ref{prop:lambda_is_small_shift} imply
\begin{align*}
\mathbb P^{\otimes N}(\Omega_{\mathrm{small}})\geq1-C\exp(-c|J_\star|)-C\exp\left(-\frac{cN}{\log^2(\frac{e|J_\star^c|}{N})}\right).
\end{align*}
Since $|J_\star|\log^2(\frac{e|J_\star^c|}{N})\leq cN$, the right-hand side is bounded from below by $1-C\exp(-c|J_\star|)$. On $\Omega_{\mathrm{small}}$, the feature-space decomposition at the beginning of this section, we obtain that
\begin{align*}
\|\widehat{\boldsymbol\beta}-\boldsymbol\beta^\star\|_2^2
\leq2r_\star^2+4\|\boldsymbol\beta_{J_\star^c}^\star\|_2^2+2^{23}\frac{\sigma_\xi^2}{\log(\frac{e|J_\star^c|}{N})}\leq Cr^2(V_{J_\star},V_{J_\star^c}).
\end{align*}
Suppose now that $\lambda\geq2^{8}\sigma_\xi\sqrt{\frac{\log(\frac{e|J_\star^c|}{N})}{N}}.$
The definition of $r_\star$ gives
\begin{align*}
r_\star^2\leq C\left(\|\boldsymbol\beta_{J_\star^c}^\star\|_2^2+\sigma_\xi^2\frac{|J_\star|}{N}+\lambda^2|J_\star|\right)\leq Cr^2(V_{J_\star},V_{J_\star^c}).
\end{align*}
The choice of the absolute constants in \eqref{eq:fixed-constants} and Assumption~\ref{assumption:LASSO_upper} ensures that $3r_\star\leq\sigma_\xi$. Moreover, the bound $\lambda^2|J_\star|\leq r^2(V_{J_\star},V_{J_\star^c})$, together with $|J_\star|\geq1$, verifies the logarithmic assumption in Proposition~\ref{prop:lambda_large_shift}. Recall that $\boldsymbol\zeta=\mathbb X_{J_\star}(\boldsymbol\beta_{J_\star}^\star-\widehat{\boldsymbol\beta}_{J_\star})+\mathbb X_{J_\star^c}\boldsymbol\beta_{J_\star^c}^\star.$ On the intersection of the events $\left\{\|\widehat{\boldsymbol\beta}_{J_\star}-\boldsymbol\beta_{J_\star}^\star\|_2\leq r_\star\right\}$, $\left\{\|\mathbb X_{J_\star}\|_{\mathrm{op}}\leq2\sqrt N\right\}$, and $\left\{\|\mathbb X_{J_\star^c}\boldsymbol\beta_{J_\star^c}^\star\|_2\leq2\sqrt N\|\boldsymbol\beta_{J_\star^c}^\star\|_2\right\}$, the equality \eqref{eq:def_r_star} and $C_\star\geq2$ imply $\|\boldsymbol\zeta\|_2\leq3\sqrt N r_\star.$ Thus Proposition~\ref{prop:lambda_large_shift}, applied with $3r_\star$, Proposition~\ref{prop:null_signal_upper_bound}, Proposition~\ref{prop:estimation_part_final}, and the standard Gaussian norm estimates imply that, with probability at least
\begin{align*}
1-C\exp(-c|J_\star|)-C\exp(-cd_q)\geq1-C\exp\big(-c\min\{|J_\star|,d_q\}\big),
\end{align*}
there hold simultaneously
\begin{align*}
\|\widehat{\boldsymbol\beta}_{J_\star}-\boldsymbol\beta_{J_\star}^\star\|_2^2\leq r_\star^2,\qquad
\|\mathcal A[\boldsymbol\xi+\boldsymbol\zeta]-\mathcal A[\boldsymbol\xi]\|_2^2\leq54r_\star^2,\quad\mbox{ and }\quad
\|\mathcal A[\boldsymbol\xi]\|_2^2\leq2^{12}\mathrm{Var}_{J_\star^c}.
\end{align*}
On this event,
\begin{align*}
\|\widehat{\boldsymbol\beta}_{J_\star^c}\|_2^2
=\|\mathcal A[\boldsymbol\xi+\boldsymbol\zeta]\|_2^2 \leq2\|\mathcal A[\boldsymbol\xi+\boldsymbol\zeta]-\mathcal A[\boldsymbol\xi]\|_2^2+2\|\mathcal A[\boldsymbol\xi]\|_2^2 \leq108r_\star^2+2^{13}\mathrm{Var}_{J_\star^c}.
\end{align*}
Therefore $\|\widehat{\boldsymbol\beta}-\boldsymbol\beta^\star\|_2^2\leq Cr^2(V_{J_\star},V_{J_\star^c}).$
\end{proof}

\section{Proof of Theorem~\texorpdfstring{\ref{theo:lower_LASSO}}{lower bound}}\label{sec:proof_lower_bound}

We first prove the theorem for $\Sigma=I_p$. The modifications for diagonal $\Sigma$ are given in Section~\ref{sec:diagonal-lower}.

\subsection{Lower bound for bias term}
\begin{Proposition}\label{prop:lower_bound_bias}
    Grant the assumptions of Theorem~\ref{theo:lower_LASSO}.
    Then
    \begin{align*}
    \mathbb P^{\otimes N}\left\{\|\widehat{\boldsymbol\beta}-\boldsymbol\beta^\star\|_2^2\geq\frac1{60}\big(\lambda^2|J_\star| + \|\bbeta_{J_\star^c}^\star\|_2^2\big)\right\}\geq1-4\exp\left(-\frac{|J_\star|}{16}\right).
    \end{align*}
\end{Proposition}
\subsubsection{Stochastic Argument}
We shall repeatedly use the following Gaussian estimates. If \(\mathbb G\in\mathbb R^{N\times k}\) has independent standard Gaussian entries and \(t\geq0\), then \(\mathbb P\{\|\mathbb G\|_{\mathrm{op}}>\sqrt N+\sqrt k+t\}\leq e^{-\frac{t^2}{2}}\). Indeed, \(\mathbb E\|\mathbb G\|_{\mathrm{op}}\leq\sqrt N+\sqrt k\), and the operator norm is \(1\)-Lipschitz with respect to the Frobenius norm. Likewise, if \(\boldsymbol g\sim N(0,I_N)\), then \(\mathbb P\{\|\boldsymbol g\|_2>\sqrt N+t\}\leq e^{-\frac{t^2}{2}}\). The lower \(\chi^2\)-tail inequality gives \(\mathbb P\{\|\boldsymbol g\|_2^2<N-2\sqrt{Nt}\}\leq e^{-t}\) for every \(t\geq0\).

\begin{Lemma}\label{lem:general-position-support}
This lemma holds for every positive diagonal $\Sigma$. Define
\begin{align}
\Omega_{\mathrm{gp}}=\bigcap_{\substack{I\subset[p]\\|I|=N+1}}\bigcap_{\boldsymbol\varepsilon\in\{-1,1\}^{I}}\left\{(\varepsilon_j\mathbb X\boldsymbol e_j:j\in I)\text{ is affinely independent}\right\}.
\label{eq:omega-gp}
\end{align}
Then \(\mathbb P(\Omega_{\mathrm{gp}})=1\). On \(\Omega_{\mathrm{gp}}\), if $\hat\bbeta$ is unique, then \(|\operatorname{supp}(\widehat{\boldsymbol\beta})|\leq N\).
\end{Lemma}

\begin{proof}
Fix \(I=(i_0,\ldots,i_N)\subset[p]\) and \(\boldsymbol\varepsilon\in\{-1,1\}^I\). Since \((\varepsilon_{i_0}\mathbb X\boldsymbol e_{i_0},\ldots,\varepsilon_{i_N}\mathbb X\boldsymbol e_{i_N})\) has the same distribution as \((\mathbb X\boldsymbol e_{i_0},\ldots,\mathbb X\boldsymbol e_{i_N})\),
\begin{align*}
&\mathbb P\left\{\mathbb X\boldsymbol e_{i_k}\in\operatorname{aff}(\mathbb X\boldsymbol e_{i_0},\ldots,\mathbb X\boldsymbol e_{i_{k-1}})\mid\mathbb X\boldsymbol e_{i_0},\ldots,\mathbb X\boldsymbol e_{i_{k-1}}\right\}=0,\qquad 1\leq k\leq N,\\
&\mathbb P\left\{\dim\operatorname{aff}(\varepsilon_{i_0}\mathbb X\boldsymbol e_{i_0},\ldots,\varepsilon_{i_N}\mathbb X\boldsymbol e_{i_N})<N\right\}\leq\sum_{k=1}^N\mathbb P\left\{\mathbb X\boldsymbol e_{i_k}\in\operatorname{aff}(\mathbb X\boldsymbol e_{i_0},\ldots,\mathbb X\boldsymbol e_{i_{k-1}})\right\}=0.
\end{align*}
The two intersections in \eqref{eq:omega-gp} are finite, and hence \(\mathbb P(\Omega_{\mathrm{gp}})=1\).

By the KKT condition, see, for instance, \cite[Theorem VI.2.2.1]{hiriart-urrutyConvexAnalysisMinimization1993},
\begin{align}
\frac1N\mathbb X^\top(\boldsymbol y-\mathbb X\widehat{\boldsymbol\beta})\in\lambda\partial\|\widehat{\boldsymbol\beta}\|_1.
\label{eq:full-kkt}
\end{align}
We distinguish two cases.
\begin{enumerate}
\item Suppose that \(\boldsymbol y\neq\mathbb X\widehat{\boldsymbol\beta}\). For every \(j\in\operatorname{supp}(\widehat{\boldsymbol\beta})\), $\left\langle\operatorname{sign}(\widehat\beta_j)\mathbb X\boldsymbol e_j,\boldsymbol y-\mathbb X\widehat{\boldsymbol\beta}\right\rangle=N\lambda.$
Consequently,
\begin{align*}
\mathrm{aff}(\operatorname{sign}(\widehat\beta_j)\mathbb X\boldsymbol e_j:j\in\operatorname{supp}(\widehat{\boldsymbol\beta}))\subset\left(\boldsymbol x\in\mathbb R^N:\langle\boldsymbol x,\boldsymbol y-\mathbb X\widehat{\boldsymbol\beta}\rangle=N\lambda\right).
\end{align*}
Notice that when $\vy\neq\bX\hat\bbeta$, $(\boldsymbol x\in\mathbb R^N:\langle\boldsymbol x,\boldsymbol y-\mathbb X\widehat{\boldsymbol\beta}\rangle=N\lambda)$ has dimension $N-1$. Suppose $|\supp(\hat\bbeta)|\geq N+1$, then for $I\subset\supp(\hat\bbeta)$, $|I|= N+1$ and $\varepsilon_j=\sign(\hat\beta_j)$ for $j\in I$, $(\varepsilon_j\bX\ve_j)_{j\in I}$ belongs to a $(N-1)$-dimensional space, hence is affinely dependent. This can happen only outside $\Omega_{\mathrm{gp}}$.
\item Suppose that \(\boldsymbol y=\mathbb X\widehat{\boldsymbol\beta}\). If \(\lambda>0\), then the KKT condition gives \(\widehat{\boldsymbol\beta}=\boldsymbol0\). For \(\lambda=0\), when \(\widehat{\boldsymbol\beta}\) is the unique minimum-\(\ell_1^p\)-norm interpolant, \cite[Theorem 3.1]{foucartMathematicalIntroductionCompressive2013} gives \(|\operatorname{supp}(\widehat{\boldsymbol\beta})|\leq N\).
\end{enumerate}
\end{proof}

Set \(\boldsymbol\zeta=\mathbb X_{J_\star^c}\boldsymbol\beta_{J_\star^c}^\star+\boldsymbol\xi\), and define
\begin{align}
\Omega_{J_\star,\mathrm{sign}}=\left\{\sum_{j\in J_\star}\mathbf 1_{\left\{\operatorname{sign}(\beta_j^\star)\left\langle\mathbb X\boldsymbol e_j,\operatorname{Proj}_{N\lambda(\mathbb X_{J_\star^c}B_1^{J_\star^c})^\circ}(\boldsymbol\zeta)\right\rangle\leq0\right\}}\geq\frac{8|J_\star|}{25}\right\}.
\label{eq:omega-J-sign}
\end{align} On \((\operatorname{Proj}_{N\lambda(\mathbb X_{J_\star^c}B_1^{J_\star^c})^\circ}(\boldsymbol\zeta)=0)\), the sum in \eqref{eq:omega-J-sign} equals \(|J_\star|\). 
On \((\operatorname{Proj}_{N\lambda(\mathbb X_{J_\star^c}B_1^{J_\star^c})^\circ}(\boldsymbol\zeta)\neq0)\),
\begin{align*}
&\left.\left(\frac{\operatorname{sign}(\beta_j^\star)\left\langle\mathbb X\boldsymbol e_j,\operatorname{Proj}_{N\lambda(\mathbb X_{J_\star^c}B_1^{J_\star^c})^\circ}(\boldsymbol\zeta)\right\rangle}{\left\|\operatorname{Proj}_{N\lambda(\mathbb X_{J_\star^c}B_1^{J_\star^c})^\circ}(\boldsymbol\zeta)\right\|_2}\right)_{j\in J_\star}\right|(\mathbb X_{J_\star^c},\boldsymbol\xi)\sim N(0,I_{|J_\star|}).
\end{align*}
Hence by Hoeffding's inequality, see, for instance, \cite[Theorem 2.2.6]{vershyninHighDimensionalProbabilityIntroduction2018}, $\mathbb P\left(\Omega_{J_\star,\mathrm{sign}}^c\mid\mathbb X_{J_\star^c},\boldsymbol\xi\right)\leq\exp\left(-\frac{|J_\star|}{16}\right).$
Define
\begin{align*}
\Omega_{J_\star,\mathrm{op}}=\left\{\|\mathbb X_{J_\star}\|_{\mathrm{op}}^2\leq\frac{13N}{10}\right\},\qquad \Omega_{T,\mathrm{op}}=\left\{\|\mathbb X_T\|_{\mathrm{op}}^2\leq\frac{21N}{5}\right\}.
\end{align*}
The Gaussian operator-norm estimate and Assumption~\ref{assumption:LASSO_lower} give \(\mathbb P(\Omega_{J_\star,\mathrm{op}}^c)\leq e^{-|J_\star|/16}\) and \(\mathbb P(\Omega_{T,\mathrm{op}}^c)\leq e^{-|J_\star|/16}\).
Let \(T\subset J_\star^c\) be the lexicographically first subset satisfying
    \begin{align}
    |T|=N,\qquad \|\boldsymbol\beta_T^\star\|_2^2=\max_{\substack{I\subset J_\star^c\\|I|=N}}\|\boldsymbol\beta_I^\star\|_2^2.
    \label{eq:def-T}
    \end{align}
\begin{enumerate}
\item\label{item:T-sign-large} If \(\|\boldsymbol\beta_T^\star\|_2^2\geq\lambda^2|J_\star|\), define
\begin{align}
\Omega_{T,\mathrm{sign}}=\left\{\sum_{j\in T}(\beta_j^\star)^2\mathbf 1_{\left\{\operatorname{sign}(\beta_j^\star)\left\langle\mathbb X\boldsymbol e_j,\operatorname{Proj}_{N\lambda(\mathbb X_{T^c}B_1^{T^c})^\circ}\left(\mathbb X_{T^c}\boldsymbol\beta_{T^c}^\star+\boldsymbol\xi\right)\right\rangle\leq0\right\}}\geq\frac8{25}\|\boldsymbol\beta_T^\star\|_2^2\right\}.
\label{eq:omega-T-sign}
\end{align}On \((\operatorname{Proj}_{N\lambda(\mathbb X_{T^c}B_1^{T^c})^\circ}(\mathbb X_{T^c}\boldsymbol\beta_{T^c}^\star+\boldsymbol\xi)=0)\), the sum in \eqref{eq:omega-T-sign} equals \(\|\boldsymbol\beta_T^\star\|_2^2\). 
On \((\operatorname{Proj}_{N\lambda(\mathbb X_{T^c}B_1^{T^c})^\circ}(\mathbb X_{T^c}\boldsymbol\beta_{T^c}^\star+\boldsymbol\xi)\neq0)\),
\begin{align*}
&\left.\left(\mathbf 1_{\left\{\operatorname{sign}(\beta_j^\star)\left\langle\mathbb X\boldsymbol e_j,\operatorname{Proj}_{N\lambda(\mathbb X_{T^c}B_1^{T^c})^\circ}\left(\mathbb X_{T^c}\boldsymbol\beta_{T^c}^\star+\boldsymbol\xi\right)\right\rangle\leq0\right\}}\right)_{\substack{j\in T\\\beta_j^\star\neq0}}\right|(\mathbb X_{T^c},\boldsymbol\xi) \sim\operatorname{Bernoulli}(\frac{1}{2})^{\otimes|(j\in T:\beta_j^\star\neq0)|}.
\end{align*}
Since \(T\subset J_\star^c\), \((\beta_j^\star)^2\leq\lambda^2\) for every \(j\in T\). Hence by Hoeffding's inequality again,
\begin{align*}
&\mathbb P\left(\Omega_{T,\mathrm{sign}}^c\mid\mathbb X_{T^c},\boldsymbol\xi\right)\leq\exp\left(-\frac{81\|\boldsymbol\beta_T^\star\|_2^4}{1250\sum_{j\in T}(\beta_j^\star)^4}\right)\leq\exp\left(-\frac{81\|\boldsymbol\beta_T^\star\|_2^2}{1250\lambda^2}\right)\leq e^{-|J_\star|/16}.
\end{align*}
Consequently, \(\mathbb P(\Omega_{T,\mathrm{sign}}^c)\leq e^{-|J_\star|/16}\).

\item\label{item:T-sign-small} If \(\|\boldsymbol\beta_T^\star\|_2^2<\lambda^2|J_\star|\), set \(\Omega_{T,\mathrm{sign}}\) equal to the whole sample space.
\end{enumerate}
Then
\begin{align*}
\mathbb P(\Omega_{\mathrm{signal}}^c)\leq e^{-\frac{|J_\star|}{16}}+3e^{-\frac{|J_\star|}{16}},\mbox{ where }\Omega_{\mathrm{signal}}=\Omega_{\mathrm{gp}}\cap\Omega_{J_\star,\mathrm{op}}\cap\Omega_{J_\star,\mathrm{sign}}\cap\Omega_{T,\mathrm{op}}\cap\Omega_{T,\mathrm{sign}}.
\end{align*}

\subsubsection{Deterministic Argument}

The following KKT identities and Lemma~\ref{lem:coordinatewise-lower-bound} also apply when the design has a positive diagonal covariance matrix \(\Sigma=\operatorname{diag}(\sigma_1,\ldots,\sigma_p)\). Here \(\Sigma_J\) denotes its restriction to \(V_J\), and \(\Sigma_J^{-1/2}\) is defined on \(V_J\). We assume \(\lambda>0\), as required by the definition of \(\boldsymbol\lambda^\star[\cdot]\).

Since \(\widehat{\boldsymbol\beta}_J\) minimizes \(\boldsymbol\beta_J\mapsto P_N\ell_{\boldsymbol\beta_J}^{(\lambda,J)}+\|\boldsymbol\beta_J\|_1\), \cite[Theorem VI.2.2.1]{hiriart-urrutyConvexAnalysisMinimization1993} and \eqref{eq:effective-gradient} give
\begin{align}
\mathbb X_J^\top\boldsymbol\lambda^\star\bigg[\mathbb X_{J^c}\boldsymbol\beta_{J^c}^\star+\boldsymbol\xi+\mathbb X_J(\boldsymbol\beta_J^\star-\widehat{\boldsymbol\beta}_J)\bigg]\in\partial\|\widehat{\boldsymbol\beta}_J\|_1.\label{eq:effective-kkt}
\end{align}
By \cite[Proposition III.3.1.3]{hiriart-urrutyConvexAnalysisMinimization1993}, whenever \(\|\mathbb X_J\Sigma_J^{-1/2}\|_{\mathrm{op}}^2\leq\frac{21N}{5}\),
\begin{align}
&\lambda^2\bigg\|\Sigma_J^{-1/2}\mathbb X_J^\top\bigg(\boldsymbol\lambda^\star\bigg[\mathbb X_{J^c}\boldsymbol\beta_{J^c}^\star+\boldsymbol\xi+\mathbb X_J(\boldsymbol\beta_J^\star-\widehat{\boldsymbol\beta}_J)\bigg]-\boldsymbol\lambda^\star\bigg[\mathbb X_{J^c}\boldsymbol\beta_{J^c}^\star+\boldsymbol\xi\bigg]\bigg)\bigg\|_2^2\notag\\
&\leq\frac{\|\mathbb X_J\Sigma_J^{-1/2}\|_{\mathrm{op}}^2}{N^2}\big\|\mathbb X_J(\boldsymbol\beta_J^\star-\widehat{\boldsymbol\beta}_J)\big\|_2^2\leq\bigg(\frac{21}{5}\bigg)^2\big\|\Sigma_J^{1/2}(\widehat{\boldsymbol\beta}_J-\boldsymbol\beta_J^\star)\big\|_2^2.\label{eq:dual-displacement}
\end{align}

\begin{Lemma}\label{lem:coordinatewise-lower-bound}
For every \(J\subset[p]\) and every positive diagonal \(\Sigma\),
\begin{align}
&\sum_{j\in J}\min\bigg\{\frac{\lambda^2}{\sigma_j},\sigma_j(\beta_j^\star)^2\bigg\}\mathbf1_{\big\{\operatorname{sign}(\beta_j^\star)\langle\mathbb X\boldsymbol e_j,\boldsymbol\lambda^\star[\mathbb X_{J^c}\boldsymbol\beta_{J^c}^\star+\boldsymbol\xi]\rangle\leq0\big\}}\notag\\
&\leq\big\|\Sigma_J^{1/2}(\widehat{\boldsymbol\beta}_J-\boldsymbol\beta_J^\star)\big\|_2^2+\lambda^2\bigg\|\Sigma_J^{-1/2}\mathbb X_J^\top\bigg(\boldsymbol\lambda^\star\bigg[\mathbb X_{J^c}\boldsymbol\beta_{J^c}^\star+\boldsymbol\xi+\mathbb X_J(\boldsymbol\beta_J^\star-\widehat{\boldsymbol\beta}_J)\bigg]-\boldsymbol\lambda^\star\bigg[\mathbb X_{J^c}\boldsymbol\beta_{J^c}^\star+\boldsymbol\xi\bigg]\bigg)\bigg\|_2^2.\label{eq:coordinatewise-lower-bound}
\end{align}
\end{Lemma}

\begin{proof}
Fix \(j\in J\) for which the indicator in \eqref{eq:coordinatewise-lower-bound} equals one.
\begin{enumerate}
\item If \(|\widehat\beta_j-\beta_j^\star|\geq|\beta_j^\star|\), then \(\min\{\frac{\lambda^2}{\sigma_j},\sigma_j(\beta_j^\star)^2\}\leq\sigma_j(\widehat\beta_j-\beta_j^\star)^2\).
\item If \(|\widehat\beta_j-\beta_j^\star|<|\beta_j^\star|\), then \(\widehat\beta_j\neq0\) and \(\operatorname{sign}(\widehat\beta_j)=\operatorname{sign}(\beta_j^\star)\). The \(j\)-th coordinate of \eqref{eq:effective-kkt} gives
\begin{align*}
\operatorname{sign}(\beta_j^\star)\bigg\langle\mathbb X\boldsymbol e_j,\boldsymbol\lambda^\star\bigg[\mathbb X_{J^c}\boldsymbol\beta_{J^c}^\star+\boldsymbol\xi+\mathbb X_J(\boldsymbol\beta_J^\star-\widehat{\boldsymbol\beta}_J)\bigg]\bigg\rangle=1.
\end{align*}
Since the indicator equals one,
\begin{align*}
\operatorname{sign}(\beta_j^\star)\bigg\langle\mathbb X\boldsymbol e_j,\boldsymbol\lambda^\star\bigg[\mathbb X_{J^c}\boldsymbol\beta_{J^c}^\star+\boldsymbol\xi+\mathbb X_J(\boldsymbol\beta_J^\star-\widehat{\boldsymbol\beta}_J)\bigg]-\boldsymbol\lambda^\star\bigg[\mathbb X_{J^c}\boldsymbol\beta_{J^c}^\star+\boldsymbol\xi\bigg]\bigg\rangle\geq1.
\end{align*}
Consequently,
\begin{align*}
\min\bigg\{\frac{\lambda^2}{\sigma_j},\sigma_j(\beta_j^\star)^2\bigg\}\leq\frac{\lambda^2}{\sigma_j}\leq\frac{\lambda^2}{\sigma_j}\bigg\langle\mathbb X\boldsymbol e_j,\boldsymbol\lambda^\star\bigg[\mathbb X_{J^c}\boldsymbol\beta_{J^c}^\star+\boldsymbol\xi+\mathbb X_J(\boldsymbol\beta_J^\star-\widehat{\boldsymbol\beta}_J)\bigg]-\boldsymbol\lambda^\star\bigg[\mathbb X_{J^c}\boldsymbol\beta_{J^c}^\star+\boldsymbol\xi\bigg]\bigg\rangle^2.
\end{align*}
\end{enumerate}
Summing gives \eqref{eq:coordinatewise-lower-bound}.
\end{proof}
We now return to the isotropic assumptions of Proposition~\ref{prop:lower_bound_bias}.
\paragraph{Lower bound for \texorpdfstring{$\|\widehat{\boldsymbol\beta}_{J_\star}-\boldsymbol\beta_{J_\star}^\star\|_2^2$}{VJ}.}
On \(\Omega_{J_\star,\mathrm{sign}}\cap\Omega_{J_\star,\mathrm{op}}\), \eqref{eq:scaled-dual-solution}, Lemma~\ref{lem:coordinatewise-lower-bound}, and \eqref{eq:dual-displacement} with \(J=J_\star\) give
\begin{align*}
&\frac{8\lambda^2|J_\star|}{25}\leq\|\widehat{\boldsymbol\beta}_{J_\star}-\boldsymbol\beta_{J_\star}^\star\|_2^2+\lambda^2\left\|\mathbb X_{J_\star}^\top\left(\boldsymbol\lambda^\star\left[\boldsymbol\zeta+\mathbb X_{J_\star}(\boldsymbol\beta_{J_\star}^\star-\widehat{\boldsymbol\beta}_{J_\star})\right]-\boldsymbol\lambda^\star[\boldsymbol\zeta]\right)\right\|_2^2\leq\left[1+\left(\frac{13}{10}\right)^2\right]\|\widehat{\boldsymbol\beta}_{J_\star}-\boldsymbol\beta_{J_\star}^\star\|_2^2,
\end{align*}hence
\begin{align}
    \|\widehat{\boldsymbol\beta}_{J_\star}-\boldsymbol\beta_{J_\star}^\star\|_2^2\geq\frac{8\lambda^2|J_\star|}{25[1+(13/10)^2]}=\frac{32\lambda^2|J_\star|}{269}.
\label{eq:J-signal-lower}
\end{align}
\paragraph{Lower bound for\texorpdfstring{$\|\widehat{\boldsymbol\beta}_{J_\star^c}-\boldsymbol\beta_{J_\star^c}^\star\|_2^2$}{VJc}.}
By Lemma~\ref{lem:general-position-support} and \eqref{eq:def-T}, $\sum_{j\in\operatorname{supp}(\widehat{\boldsymbol\beta})\cap J_\star^c}(\beta_j^\star)^2\leq\|\boldsymbol\beta_T^\star\|_2^2.$
Therefore,
\begin{align}
\begin{aligned}
    &\|\widehat{\boldsymbol\beta}_{J_\star^c}-\boldsymbol\beta_{J_\star^c}^\star\|_2^2=\sum_{j\in\operatorname{supp}(\widehat{\boldsymbol\beta})\cap J_\star^c}(\widehat\beta_j-\beta_j^\star)^2+\sum_{j\in J_\star^c\setminus\operatorname{supp}(\widehat{\boldsymbol\beta})}(\beta_j^\star)^2\\
&\geq\sum_{j\in J_\star^c\setminus\operatorname{supp}(\widehat{\boldsymbol\beta})}(\beta_j^\star)^2=\|\boldsymbol\beta_{J_\star^c}^\star\|_2^2-\sum_{j\in\operatorname{supp}(\widehat{\boldsymbol\beta})\cap J_\star^c}(\beta_j^\star)^2 \geq\|\boldsymbol\beta_{J_\star^c}^\star\|_2^2-\|\boldsymbol\beta_T^\star\|_2^2=\|\boldsymbol\beta_{J_\star^c\setminus T}^\star\|_2^2.
\end{aligned}
\label{eq:Jc-tail-lower}
\end{align}

\begin{enumerate}
\item\label{item:deterministic-T-large} Suppose that \(\|\boldsymbol\beta_T^\star\|_2^2\geq\lambda^2|J_\star|\). Since \(T\subset J_\star^c\), \(\min\{\lambda^2,(\beta_j^\star)^2\}=(\beta_j^\star)^2\) for every \(j\in T\). Hence, on \(\Omega_{T,\mathrm{sign}}\cap\Omega_{T,\mathrm{op}}\), \eqref{eq:scaled-dual-solution}, Lemma~\ref{lem:coordinatewise-lower-bound}, and \eqref{eq:dual-displacement} with \(J=T\) give
\begin{align*}
&\frac8{25}\|\boldsymbol\beta_T^\star\|_2^2\leq\|\widehat{\boldsymbol\beta}_T-\boldsymbol\beta_T^\star\|_2^2 +\lambda^2\left\|\mathbb X_T^\top\left(\boldsymbol\lambda^\star\left[\mathbb X_{T^c}\boldsymbol\beta_{T^c}^\star+\boldsymbol\xi+\mathbb X_T(\boldsymbol\beta_T^\star-\widehat{\boldsymbol\beta}_T)\right]-\boldsymbol\lambda^\star\left[\mathbb X_{T^c}\boldsymbol\beta_{T^c}^\star+\boldsymbol\xi\right]\right)\right\|_2^2\\
&\leq\left[1+\left(\frac{21}{5}\right)^2\right]\|\widehat{\boldsymbol\beta}_T-\boldsymbol\beta_T^\star\|_2^2,
\end{align*}hence
\begin{align}
    \|\widehat{\boldsymbol\beta}_T-\boldsymbol\beta_T^\star\|_2^2\geq\frac{8\|\boldsymbol\beta_T^\star\|_2^2}{25[1+(21/5)^2]}=\frac{4\|\boldsymbol\beta_T^\star\|_2^2}{233}.
\label{eq:T-signal-lower}
\end{align}
Since \(T\subset J_\star^c\), \eqref{eq:T-signal-lower} gives
\begin{align*}
&\|\widehat{\boldsymbol\beta}_{J_\star^c}-\boldsymbol\beta_{J_\star^c}^\star\|_2^2=\|\widehat{\boldsymbol\beta}_T-\boldsymbol\beta_T^\star\|_2^2+\|\widehat{\boldsymbol\beta}_{J_\star^c\setminus T}-\boldsymbol\beta_{J_\star^c\setminus T}^\star\|_2^2\geq\|\widehat{\boldsymbol\beta}_T-\boldsymbol\beta_T^\star\|_2^2\geq\frac4{233}\|\boldsymbol\beta_T^\star\|_2^2.
\end{align*}
Combining this inequality with \eqref{eq:Jc-tail-lower} yields
\begin{align*}
&\|\widehat{\boldsymbol\beta}_{J_\star^c}-\boldsymbol\beta_{J_\star^c}^\star\|_2^2\geq\max\left\{\|\boldsymbol\beta_{J_\star^c\setminus T}^\star\|_2^2,\frac4{233}\|\boldsymbol\beta_T^\star\|_2^2\right\}.
\end{align*}Moreover, $\|\boldsymbol\beta_{J_\star^c}^\star\|_2^2=\|\boldsymbol\beta_{J_\star^c\setminus T}^\star\|_2^2+\frac{233}{4}\frac{4\|\boldsymbol\beta_T^\star\|_2^2}{233} \leq\frac{237}{4}\max\left\{\|\boldsymbol\beta_{J_\star^c\setminus T}^\star\|_2^2,\frac4{233}\|\boldsymbol\beta_T^\star\|_2^2\right\}$, hence
\begin{align}
&\|\widehat{\boldsymbol\beta}_{J_\star^c}-\boldsymbol\beta_{J_\star^c}^\star\|_2^2\geq\frac4{237}\|\boldsymbol\beta_{J_\star^c}^\star\|_2^2.
\label{eq:Jc-large-lower}
\end{align}

\item\label{item:deterministic-T-small} Suppose that \(\|\boldsymbol\beta_T^\star\|_2^2<\lambda^2|J_\star|\). Then
\begin{align*}
\lambda^2|J_\star| + \|\bbeta_{J_\star^c}^\star\|_2^2  \leq2\lambda^2|J_\star|+\|\boldsymbol\beta_{J_\star^c\setminus T}^\star\|_2^2.
\end{align*}
Consequently, \eqref{eq:J-signal-lower} and \eqref{eq:Jc-tail-lower} give
\begin{align*}
\|\widehat{\boldsymbol\beta}-\boldsymbol\beta^\star\|_2^2\geq\frac{32\lambda^2|J_\star|}{269}+\|\boldsymbol\beta_{J_\star^c\setminus T}^\star\|_2^2\geq\frac{16}{269}\left(\lambda^2|J_\star|+\|\boldsymbol\beta_{J_\star^c}^\star\|_2^2\right).
\end{align*}
\end{enumerate}
Thus,
\begin{align}
\bP\bigg\{\|\widehat{\boldsymbol\beta}-\boldsymbol\beta^\star\|_2^2\geq\frac1{60}\left(\lambda^2 |J_\star|+\|\boldsymbol\beta_{J_\star^c}^\star\|_2^2\right)\bigg\}\geq 1-4\exp(-\frac{|J_\star|}{16}).
\label{eq:base-lower-bound}
\end{align}

\subsection{Lower bound for variance terms on \texorpdfstring{$V_{J_\star}$}{VJ}}

We prove the following proposition, which serves as the variance lower bound on $V_{J_\star}$.
\begin{Proposition}\label{prop:active-sign-dropout}
Grant the assumptions of Theorem~\ref{theo:lower_LASSO}. Then for any $\lambda\geq0,$
\begin{align}
&\mathbb P\!\left\{\|\hat{\boldsymbol\beta}_{J_\star}[\vy]-\boldsymbol\beta_{J_\star}^\star\|_2^2\ge\frac{|J_\star|}{16}\min\bigg\{\frac{\sigma_\xi^2}{N},\lambda^2\bigg\}\right\}\ge 1-e^{-0.2|J_\star|}.\label{eq:conditional-tail-lower-bound}
\end{align}
\end{Proposition}

\beginproof The first step of the proof is the following lemma~\ref{lem:many-correct-signs}. It shows that on $V_{J_\star}$, if the estimation error of $\hat\bbeta_{J_\star}[\cdot]$ is small, then the signs of most coordinates of $\hat\bbeta_{J_\star}$ almost match those of $\bbeta_{J_\star}^\star$. That is, small estimation error implies sign matching for most coordinates.
\begin{Lemma}\label{lem:many-correct-signs}
Let $\lambda\geq 0$ be any real number. For every $\boldsymbol v\in\mathbb R^N$, there holds
\begin{align}
&\left\{\|\hat{\boldsymbol\beta}_{J_\star}[\boldsymbol v]-\boldsymbol\beta_{J_\star}^\star\|_2^2<\frac{|J_\star|}{16}\min\{\frac{\sigma_\xi^2}{N},\lambda^2\}\right\}\subseteq\bigcup\left\{\|\hat{\boldsymbol\beta}_I[\boldsymbol v]-\boldsymbol\beta_I^\star\|_2^2<\frac{|J_\star|}{16}\min\{\frac{\sigma_\xi^2}{N},\lambda^2\},\ \operatorname{sign}(\hat{\boldsymbol\beta}_I[\boldsymbol v])=\operatorname{sign}(\boldsymbol\beta_I^\star)\right\},\label{eq:many-correct-signs}
\end{align}where the union is over all $I\subset J_\star$ such that $|I|=|J_\star|-\lfloor\frac{|J_\star|}{16}\rfloor$.
\end{Lemma}
\begin{proof}
If $\lambda=0$, then the event on the left-hand side of \eqref{eq:many-correct-signs} is empty. Suppose that $\lambda>0$. For every $j\in J_\star$, if $\operatorname{sign}(\hat\beta_j[\boldsymbol v])\ne\operatorname{sign}(\beta_j^\star)$, then $\hat\beta_j[\boldsymbol v]\beta_j^\star\leq 0$, and consequently $(\hat\beta_j[\boldsymbol v]-\beta_j^\star)^2\ge(\beta_j^\star)^2>\lambda^2\ge\min\{\frac{\sigma_\xi^2}{N},\lambda^2\}$. Therefore, $\|\hat{\boldsymbol\beta}_{J_\star}[\boldsymbol v]-\boldsymbol\beta_{J_\star}^\star\|_2^2<\frac{|J_\star|}{16}\min\{\frac{\sigma_\xi^2}{N},\lambda^2\}$ implies
\begin{align*}
&\left|(j\in J_\star:\operatorname{sign}(\hat\beta_j[\boldsymbol v])\ne\operatorname{sign}(\beta_j^\star))\right|<\frac{|J_\star|}{16}.
\end{align*}
At least $|J_\star|-\lfloor\frac{|J_\star|}{16}\rfloor$ coordinates consequently have the correct sign. Let $I$ be any such set of coordinates. Then
\begin{align*}
&\|\hat{\boldsymbol\beta}_I[\boldsymbol v]-\boldsymbol\beta_I^\star\|_2^2\le\|\hat{\boldsymbol\beta}_{J_\star}[\boldsymbol v]-\boldsymbol\beta_{J_\star}^\star\|_2^2,
\end{align*}
which proves \eqref{eq:many-correct-signs}.
\end{proof}
The claim of Proposition~\ref{prop:active-sign-dropout} is immediate when $\lambda=0$. Henceforth, suppose that $\lambda>0$. Therefore, to prove the high-probability lower bound for the estimation error, it suffices to show that for any set $I$ described by the proof of Lemma~\ref{lem:many-correct-signs}, the probability of the random event $\{\|\hat{\boldsymbol\beta}_I[\boldsymbol v]-\boldsymbol\beta_I^\star\|_2^2<\frac{|J_\star|}{16}\min\{\frac{\sigma_\xi^2}{N},\lambda^2\},\ \operatorname{sign}(\hat{\boldsymbol\beta}_I[\boldsymbol v])=\operatorname{sign}(\boldsymbol\beta_I^\star)\}$ is so small that, after taking a union bound, the probability of the stronger union obtained in the proof of Lemma~\ref{lem:many-correct-signs} is negligibly small. Hence, we study the random event $\{\|\hat{\boldsymbol\beta}_I[\boldsymbol v]-\boldsymbol\beta_I^\star\|_2^2<\frac{|J_\star|}{16}\min\{\frac{\sigma_\xi^2}{N},\lambda^2\},\ \operatorname{sign}(\hat{\boldsymbol\beta}_I[\boldsymbol v])=\operatorname{sign}(\boldsymbol\beta_I^\star)\}$. We need the following lemma~\ref{lem:common-offset}.

\begin{Lemma}\label{lem:common-offset}
Fix $\mathbb X\in\Omega_{\mathrm{gp}}$, $\boldsymbol u\in\mathbb R^N$, a nonempty $I\subseteq[p]$, and $\boldsymbol s\in\{-1,1\}^{I}$. For any $\boldsymbol t,\boldsymbol t'\in\mathbb R^I$ such that
\begin{align}
&\operatorname{sign}(\hat{\boldsymbol\beta}_I[\mathbb X\boldsymbol\beta^\star+\boldsymbol u+\mathbb X_I\boldsymbol t])=\operatorname{sign}(\hat{\boldsymbol\beta}_I[\mathbb X\boldsymbol\beta^\star+\boldsymbol u+\mathbb X_I\boldsymbol t'])=\boldsymbol s.\label{eq:common-sign-condition}
\end{align}
one has
\begin{align*}
&\hat{\boldsymbol\beta}_I[\mathbb X\boldsymbol\beta^\star+\boldsymbol u+\mathbb X_I\boldsymbol t]-\boldsymbol t=\hat{\boldsymbol\beta}_I[\mathbb X\boldsymbol\beta^\star+\boldsymbol u+\mathbb X_I\boldsymbol t']-\boldsymbol t'.
\end{align*}
\end{Lemma}
\begin{proof}
We first construct a convex optimization problem with a unique solution, such that its solution is exactly $(\vt - \hat\bbeta_I[\bX\bbeta^\star+\vu+\bX_I\vt], \hat\bbeta_{I^c}[\bX\bbeta^\star+\vu+\bX_I\vt])$.
We prove the following equation:
\begin{align}
&\left(\boldsymbol t-\hat{\boldsymbol\beta}_I[\mathbb X\boldsymbol\beta^\star+\boldsymbol u+\mathbb X_I\boldsymbol t],\hat{\boldsymbol\beta}_{I^c}[\mathbb X\boldsymbol\beta^\star+\boldsymbol u+\mathbb X_I\boldsymbol t]\right)\in\operatorname*{argmin}_{\substack{\boldsymbol\delta\in\mathbb R^I\\\boldsymbol\nu\in\mathbb R^{I^c}}}\left(\frac{1}{2N}\|\mathbb X\boldsymbol\beta^\star+\boldsymbol u+\mathbb X_I\boldsymbol\delta-\mathbb X_{I^c}\boldsymbol\nu\|_2^2-\lambda\langle\boldsymbol s,\boldsymbol\delta\rangle+\lambda\|\boldsymbol\nu\|_1\right).\label{eq:transformed-program}
\end{align}Indeed, setting the gradients of this objective function with respect to $\vdelta,\vnu$ to $\vzero$, we obtain 
\begin{align}
&\mathbb X_I^\top(\mathbb X\boldsymbol\beta^\star+\boldsymbol u+\mathbb X_I\boldsymbol\delta-\mathbb X_{I^c}\boldsymbol\nu)=N\lambda\boldsymbol s,\mbox{ and } \mathbb X_{I^c}^\top(\mathbb X\boldsymbol\beta^\star+\boldsymbol u+\mathbb X_I\boldsymbol\delta-\mathbb X_{I^c}\boldsymbol\nu)\in N\lambda\partial^-\|\cdot\|_1(\boldsymbol\nu).\label{eq:transformed-kkt}
\end{align}
It suffices to verify that for $\vdelta = \vt - \hat\bbeta_I[\bX\bbeta^\star+\vu+\bX_I\vt]$, $\vnu = \hat\bbeta_{I^c}[\bX\bbeta^\star+\vu+\bX_I\vt]$, the above two conditions are satisfied. This only requires us to use the KKT conditions provided by the definition of $\hat\bbeta$, together with the condition $\sign(\hat\bbeta_I[\bX\bbeta^\star+\vu+\bX_I\vt])=\vs$ provided by \eqref{eq:common-sign-condition}, that is,
\begin{align*}
\begin{aligned}
    &\mathbb X_I^\top\!\left(\mathbb X\boldsymbol\beta^\star+\boldsymbol u+\mathbb X_I\boldsymbol t-\mathbb X_{I^c}\hat{\boldsymbol\beta}[\mathbb X\boldsymbol\beta^\star+\boldsymbol u+\mathbb X_I\boldsymbol t]\right)=N\lambda\boldsymbol s,\\
    &\mathbb X_{I^c}^\top\!\left(\mathbb X\boldsymbol\beta^\star+\boldsymbol u+\mathbb X_I\boldsymbol t-\mathbb X_{I^c}\hat{\boldsymbol\beta}[\mathbb X\boldsymbol\beta^\star+\boldsymbol u+\mathbb X_I\boldsymbol t]\right)\in N\lambda\partial^-\|\cdot\|_1(\hat{\boldsymbol\beta}_{I^c}[\mathbb X\boldsymbol\beta^\star+\boldsymbol u+\mathbb X_I\boldsymbol t]).
\end{aligned}
\end{align*}
This exactly proves equation \eqref{eq:transformed-program}. Replacing $\boldsymbol t$ with $\boldsymbol t'$ in the KKT condition, together with \eqref{eq:common-sign-condition}, shows that the tuple $(\boldsymbol t'-\hat{\boldsymbol\beta}_I[\mathbb X\boldsymbol\beta^\star+\boldsymbol u+\mathbb X_I\boldsymbol t'],\hat{\boldsymbol\beta}_{I^c}[\mathbb X\boldsymbol\beta^\star+\boldsymbol u+\mathbb X_I\boldsymbol t'])$ also minimizes the criterion in \eqref{eq:transformed-program}. 

Next, we prove the uniqueness of its solution. Before that, we point out that once uniqueness is established, since both $(\boldsymbol t'-\hat{\boldsymbol\beta}_I[\mathbb X\boldsymbol\beta^\star+\boldsymbol u+\mathbb X_I\boldsymbol t'],\hat{\boldsymbol\beta}_{I^c}[\mathbb X\boldsymbol\beta^\star+\boldsymbol u+\mathbb X_I\boldsymbol t'])$ and $(\boldsymbol t-\hat{\boldsymbol\beta}_I[\mathbb X\boldsymbol\beta^\star+\boldsymbol u+\mathbb X_I\boldsymbol t],\hat{\boldsymbol\beta}_{I^c}[\mathbb X\boldsymbol\beta^\star+\boldsymbol u+\mathbb X_I\boldsymbol t])$ are solutions to \eqref{eq:transformed-program}, it necessarily follows that $\boldsymbol t-\hat{\boldsymbol\beta}_I[\mathbb X\boldsymbol\beta^\star+\boldsymbol u+\mathbb X_I\boldsymbol t] = \boldsymbol t'-\hat{\boldsymbol\beta}_I[\mathbb X\boldsymbol\beta^\star+\boldsymbol u+\mathbb X_I\boldsymbol t']$, which is precisely the target equation. 

Let $(\boldsymbol\delta,\boldsymbol\nu)$ and $(\boldsymbol\delta',\boldsymbol\nu')$ be two minimizers in \eqref{eq:transformed-program}. Strict convexity of the squared Euclidean norm gives
\begin{align}
&\mathbb X\boldsymbol\beta^\star+\boldsymbol u+\mathbb X_I\boldsymbol\delta-\mathbb X_{I^c}\boldsymbol\nu=\mathbb X\boldsymbol\beta^\star+\boldsymbol u+\mathbb X_I\boldsymbol\delta'-\mathbb X_{I^c}\boldsymbol\nu'.\label{eq:common-residual}
\end{align}Indeed, if \eqref{eq:common-residual} does not hold, then $(\frac{1}{2}\vdelta+\frac{1}{2}\vdelta',\frac{1}{2}\vnu+\frac{1}{2}\vnu')$, by the convexity of the objective function, must produce a smaller value of the objective function, contradicting with the assumption that $(\boldsymbol\delta,\boldsymbol\nu)$ and $(\boldsymbol\delta',\boldsymbol\nu')$ are minimizers.

Applying \eqref{eq:transformed-kkt} to both minimizers gives
\begin{align*}
&I\cup\operatorname{supp}(\boldsymbol\nu)\cup\operatorname{supp}(\boldsymbol\nu')\subseteq\left(j\in[p]:\left|\left\langle\mathbb X\boldsymbol e_j,\mathbb X\boldsymbol\beta^\star+\boldsymbol u+\mathbb X_I\boldsymbol\delta-\mathbb X_{I^c}\boldsymbol\nu\right\rangle\right|=N\lambda\right).
\end{align*}
Since $I\ne\varnothing$ and $\lambda>0$, the first relation in \eqref{eq:transformed-kkt} implies that $\mathbb X\boldsymbol\beta^\star+\boldsymbol u+\mathbb X_I\boldsymbol\delta-\mathbb X_{I^c}\boldsymbol\nu\ne\boldsymbol 0$. By Lemma \ref{lem:general-position-support}, the columns indexed by the set on the right-hand side of the inclusion above are linearly independent. By \eqref{eq:common-residual},
\begin{align*}
&\boldsymbol 0=\mathbb X_I(\boldsymbol\delta-\boldsymbol\delta')-\mathbb X_{I^c}(\boldsymbol\nu-\boldsymbol\nu')=\sum_{j\in I}(\delta_j-\delta_j')\mathbb X\boldsymbol e_j-\sum_{j\in\operatorname{supp}(\boldsymbol\nu)\cup\operatorname{supp}(\boldsymbol\nu')}(\nu_j-\nu_j')\mathbb X\boldsymbol e_j.
\end{align*}
The linear independence of the columns indexed by $I\cup\operatorname{supp}(\boldsymbol\nu)\cup\operatorname{supp}(\boldsymbol\nu')$ therefore gives $\delta_j=\delta_j',$ for any $j\in I,$ and $\nu_j=\nu_j',$ for any $j\in\operatorname{supp}(\boldsymbol\nu)\cup\operatorname{supp}(\boldsymbol\nu').$
For every $j\in I^c\backslash(\operatorname{supp}(\boldsymbol\nu)\cup\operatorname{supp}(\boldsymbol\nu'))$, one has $\nu_j=\nu_j'=0$. Hence, $\boldsymbol\delta=\boldsymbol\delta'$ and $\boldsymbol\nu=\boldsymbol\nu'$.
\end{proof}
On $\Omega_{\mathrm{gp}}\cap\{\|\mathbb X_{J_\star}\|_{\mathrm{op}}\le\sqrt{2N}\}$, we apply Lemma~\ref{lem:common-offset} to study the random event $$\{\|\hat{\boldsymbol\beta}_I[\mathbb X\boldsymbol\beta^\star+\boldsymbol\xi]-\boldsymbol\beta_I^\star\|_2^2<\frac{|J_\star|}{16}\min\{\frac{\sigma_\xi^2}{N},\lambda^2\},\ \operatorname{sign}(\hat{\boldsymbol\beta}_I[\mathbb X\boldsymbol\beta^\star+\boldsymbol\xi])=\operatorname{sign}(\boldsymbol\beta_I^\star)\}.$$ Take $I$ in Lemma~\ref{lem:common-offset} to be the set on the right-hand side of Lemma~\ref{lem:many-correct-signs}, and take $\vs = \sign(\bbeta_I^\star)$. Let $R^2 = \frac{|J_\star|}{16}\min\{\frac{\sigma_\xi^2}{N},\lambda^2\}$. Let $\vt = (\bX_I^\top\bX_I)^{-1}\bX_I^\top\bxi$ and $\vu = (I_N-\operatorname{Proj}_{\mathbb XV_I})\boldsymbol\xi$; then $\hat\bbeta_I[\bX\bbeta^\star+\bxi] = \hat\bbeta_I[\bX\bbeta^\star+\vu + \bX_I \vt]$, where almost surely,
\begin{align*}
    &\operatorname{Proj}_{\mathbb XV_I}=\mathbb X_I(\mathbb X_I^\top\mathbb X_I)^{-1}\mathbb X_I^\top,\mbox{ and hence } \boldsymbol\xi=(I_N-\operatorname{Proj}_{\mathbb XV_I})\boldsymbol\xi+\mathbb X_I(\mathbb X_I^\top\mathbb X_I)^{-1}\mathbb X_I^\top\boldsymbol\xi = \vu+\bX_I\vt.
\end{align*}Since $\bxi$ is assumed to be $\cN(\vzero,\sigma_\xi^2 I_N)$, the conditional distribution of $\vt$ satisfies
\begin{align}
&\mathrm{Law}\!\left(\vt\middle|\,\mathbb X,\vu\right)=\cN\!\left(\boldsymbol 0,\sigma_\xi^2(\mathbb X_I^\top\mathbb X_I)^{-1}\right).\label{eq:conditional-gaussian-coordinates}
\end{align}
We first prove the following inclusion: there exists some $\va\in\bR^I$ such that
\begin{align}\label{eq:correct-sign-ball-containment}
\begin{aligned}
    &\left(\boldsymbol t\in\mathbb R^I:\|\hat{\boldsymbol\beta}_I[\mathbb X\boldsymbol\beta^\star+\vu+\mathbb X_I\boldsymbol t]-\boldsymbol\beta_I^\star\|_2<R,\ \operatorname{sign}(\hat{\boldsymbol\beta}_I[\mathbb X\boldsymbol\beta^\star+\vu+\mathbb X_I\boldsymbol t])=\operatorname{sign}(\boldsymbol\beta_I^\star)\right) \subseteq(\boldsymbol t\in\mathbb R^I:\|\boldsymbol t-\boldsymbol a\|_2<R).
\end{aligned}
\end{align}
If the set on the left-hand side of \eqref{eq:correct-sign-ball-containment} is empty, then the inclusion is trivial; otherwise, for any $\vt,\vt'$ from the set on the left hand side of \eqref{eq:correct-sign-ball-containment}, possibly $\vt=\vt'$, Lemma~\ref{lem:common-offset} implies that $\hat\bbeta_I[\bX\bbeta^\star+\vu+\bX_I\vt] - \bbeta_I^\star = \vt - (\vt' + \bbeta_I^\star - \hat\bbeta_I[\bX\bbeta^\star+\vu+\bX_I\vt'])$.
Therefore, for $\boldsymbol a=\boldsymbol t'+\boldsymbol\beta_I^\star-\hat{\boldsymbol\beta}_I[\mathbb X\boldsymbol\beta^\star+\vu+\mathbb X_I\boldsymbol t']$, the inclusion holds.

Combining \eqref{eq:correct-sign-ball-containment} and \eqref{eq:conditional-gaussian-coordinates},
\begin{align}\label{eq:conditional-density-bound}
\begin{aligned}
    &\mathbb P_{\boldsymbol\xi}\!\left\{\|\hat{\boldsymbol\beta}_I[\mathbb X\boldsymbol\beta^\star+\boldsymbol\xi]-\boldsymbol\beta_I^\star\|_2<R,\ \operatorname{sign}(\hat{\boldsymbol\beta}_I[\mathbb X\boldsymbol\beta^\star+\boldsymbol\xi])=\operatorname{sign}(\boldsymbol\beta_I^\star)\,\middle|\,\mathbb X,\vu\right\}\\
    &\le\sup_{\boldsymbol a\in\mathbb R^I}\int_{(\boldsymbol t\in\mathbb R^I:\|\boldsymbol t-\boldsymbol a\|_2<R)}\frac{\det(\mathbb X_I^\top\mathbb X_I)^{\frac{1}{2}}}{(2\pi\sigma_\xi^2)^{\frac{|I|}{2}}}\cdot\exp\!\left(-\frac{\boldsymbol t^\top\mathbb X_I^\top\mathbb X_I\boldsymbol t}{2\sigma_\xi^2}\right)\,\mathrm d\boldsymbol t.
\end{aligned}
\end{align}
Since every singular value of $\mathbb X_I$ is at most $\|\mathbb X_I\|_{\mathrm{op}}\le\|\mathbb X_{J_\star}\|_{\mathrm{op}}\le\sqrt{2N}$, there holds $\det(\mathbb X_I^\top\mathbb X_I)^{\frac{1}{2}}\le(2N)^{\frac{|I|}{2}}$. Moreover, $(2\pi)^{\frac{|I|}{2}}=\int_{\mathbb R^I}e^{-\frac{\|\boldsymbol u\|_2^2}{2}}\,\mathrm d\boldsymbol u\ge\int_{(\boldsymbol u\in\mathbb R^I:\|\boldsymbol u\|_2\le\sqrt{|I|})}e^{-\frac{\|\boldsymbol u\|_2^2}{2}}\,\mathrm d\boldsymbol u\ge e^{-\frac{|I|}{2}}|I|^{\frac{|I|}{2}}\operatorname{vol}(\boldsymbol u\in\mathbb R^I:\|\boldsymbol u\|_2\le1)$, see, for instance, \cite[Exercise 4.40, pp. 97]{aubrunAliceBobMeet2017}.
For every $\boldsymbol a\in\mathbb R^I$, the exponential factor in \eqref{eq:conditional-density-bound} is at most one, and the translation invariance of the Lebesgue measure gives $\operatorname{vol}(\boldsymbol t\in\mathbb R^I:\|\boldsymbol t-\boldsymbol a\|_2<R)=R^{|I|}\operatorname{vol}(\boldsymbol u\in\mathbb R^I:\|\boldsymbol u\|_2\le1).$
Consequently, the right-hand side of \eqref{eq:conditional-density-bound} is bounded, uniformly in $\boldsymbol a$ and $\vu$, by
\begin{align*}
&\frac{(2N)^{\frac{|I|}{2}}R^{|I|}}{(2\pi\sigma_\xi^2)^{\frac{|I|}{2}}}\operatorname{vol}(\boldsymbol u\in\mathbb R^I:\|\boldsymbol u\|_2\le1)\le\left(\frac{2eNR^2}{\sigma_\xi^2|I|}\right)^{\frac{|I|}{2}}.
\end{align*}
Taking the conditional expectation over $\vu$ therefore gives
\begin{align*}
&\mathbb P_{\boldsymbol\xi}\!\left\{\|\hat{\boldsymbol\beta}_I[\mathbb X\boldsymbol\beta^\star+\boldsymbol\xi]-\boldsymbol\beta_I^\star\|_2<R,\ \operatorname{sign}(\hat{\boldsymbol\beta}_I[\mathbb X\boldsymbol\beta^\star+\boldsymbol\xi])=\operatorname{sign}(\boldsymbol\beta_I^\star)\,\middle|\,\mathbb X\right\}\le\left(\frac{2eNR^2}{\sigma_\xi^2|I|}\right)^{\frac{|I|}{2}}.
\end{align*}
Apply Lemma \ref{lem:many-correct-signs} with $\boldsymbol v=\vy = \bX\bbeta^\star+\bxi$, and use the union bound with $R^2=\frac{|J_\star|}{16}\min\{\frac{\sigma_\xi^2}{N},\lambda^2\}$. We obtain
\begin{align*}
&\mathbb P_{\boldsymbol\xi}\!\left\{\|\hat{\boldsymbol\beta}_{J_\star}[\mathbb X\boldsymbol\beta^\star+\boldsymbol\xi]-\boldsymbol\beta_{J_\star}^\star\|_2^2<\frac{|J_\star|}{16}\min\{\frac{\sigma_\xi^2}{N},\lambda^2\}\,\middle|\,\mathbb X\right\}\\
&\le\binom{|J_\star|}{\lfloor\frac{|J_\star|}{16}\rfloor}\left(\frac{2e|J_\star|}{16(|J_\star|-\lfloor\frac{|J_\star|}{16}\rfloor)}\min\left\{1,\frac{N\lambda^2}{\sigma_\xi^2}\right\}\right)^{\frac{|J_\star|-\lfloor\frac{|J_\star|}{16}\rfloor}{2}}.
\end{align*}
The binomial coefficient satisfies $\binom{|J_\star|}{\lfloor\frac{|J_\star|}{16}\rfloor}\le\exp\!\left(\frac{|J_\star|}{16}\log(16e)\right)$. Since $|J_\star|-\lfloor\frac{|J_\star|}{16}\rfloor\ge\frac{15|J_\star|}{16}$, $\frac{2e}{15}<1$, and $\min\{1,\frac{N\lambda^2}{\sigma_\xi^2}\}\le1$,
\begin{align*}
&\mathbb P_{\boldsymbol\xi}\!\left\{\|\hat{\boldsymbol\beta}_{J_\star}[\mathbb X\boldsymbol\beta^\star+\boldsymbol\xi]-\boldsymbol\beta_{J_\star}^\star\|_2^2<\frac{|J_\star|}{16}\min\left\{\frac{\sigma_\xi^2}{N},\lambda^2\right\}\,\middle|\,\mathbb X\right\}\\
&\le\exp\!\left(|J_\star|\left[\frac{1}{16}\log(16e)+\frac{15}{32}\log\!\left(\frac{2e}{15}\right)\right]\right)\le e^{-0.23|J_\star|}.
\end{align*}
The Gaussian operator-norm deviation recalled above and $\frac{|J_\star|}{N}\le c<\frac{1}{100}$ give
\begin{align*}
&\mathbb P_{\mathbb X}\!\left\{\|\mathbb X_{J_\star}\|_{\mathrm{op}}>\sqrt{2N}\right\}\le\exp\!\left(-\frac{(\sqrt{2}-1-\sqrt{\frac{|J_\star|}{N}})^2N}{2}\right)\le e^{-N/25}\le e^{-4|J_\star|}.
\end{align*}
Combining the last two displays with $e^{-0.23|J_\star|}+e^{-4|J_\star|}\le e^{-0.2|J_\star|}$ proves Proposition~\ref{prop:active-sign-dropout}.

\subsection{Lower bound in \texorpdfstring{$\mathrm{Var}_{J_\star^c}$}{Var Jc}}

\begin{Proposition}
\label{prop:free-variance-lower-bound}
Grant Assumption~\ref{assumption:LASSO_lower}. Then there exist absolute constant $c_{\mathrm{est}}=2^{-30},c_{\mathrm{prob}}=2^{-35}$, and $C_{\mathrm{prob}}=2^{11}$, such that
\begin{align*}
    \bP\big\{\|\hat\bbeta[\vy]-\bbeta^\star\|_2^2\geq c_{\mathrm{est}}\Var_{J_\star^c}\big\}\geq 1-C_{\mathrm{prob}}\exp\big (-c_{\mathrm{prob}}d_q\big).
\end{align*}
\end{Proposition}

\begin{proof}
This proof again relies on CGMT; to this end, we first write down the truncated primary problem and its auxiliary problem for $\hat\bbeta[\vy]$. We first define two quantities that will be used throughout the entire proof. Let
\begin{align*}
    r^- = 2^{-15}\sqrt{\Var_{J_\star^c}},\mbox{ and }r^+ = 2^{-14}\sqrt{\Var_{J_\star^c}}.
\end{align*}For each $\bbeta\in\bR^p$, let $\vw=\bbeta-\bbeta^\star$, and let
\begin{align}
\begin{aligned}
&\Phi_{t,H}=\min_{\|\boldsymbol w\|_2\leq t}\max_{\|\boldsymbol u\|_2\leq H}\bigg\{\langle\boldsymbol u,\boldsymbol\xi-\mathbb X\boldsymbol w\rangle-\frac{N\lambda}{2}\|\boldsymbol u\|_2^2+\|\boldsymbol\beta^\star+\boldsymbol w\|_1-\|\boldsymbol\beta^\star\|_1\bigg\}.
\end{aligned}
\label{eq:primary-optimization_lower}
\end{align}The relationship between the objective function of $\Phi_{t,H}$ and the objective function of LASSO $\hat\bbeta[\vy]$ is given in Section~\ref{sec:variational_representation_LASSO}.
To control the probability of the random event $\{\|\hat\bbeta[\vy]-\bbeta^\star\|_2^2 \le 2^{-30} \Var_{J_\star^c}\} = \{\|\hat\bbeta[\vy]-\bbeta^\star\|_2 \le r^-\}$ from above by $2^{11} \exp(-2^{-35} d_q)$, we embed this random event into a random event determined by the optimal solutions of two primary problems. Here we rely on the following observation: if $\hat\vw = \hat\bbeta[\vy]-\bbeta^\star$ satisfies $\|\hat\vw\|_2 \le r^-$, then we must have $\Phi_{r^-,\infty} = \Phi_{r^+,\infty}$. Indeed, as a minimizer on the full-space feasible set of $\Phi_{\infty,\infty}$, if $\hat\vw = \hat\bbeta - \bbeta^\star$ satisfies $\|\hat\vw\|_2 \le r^-$, then its optimal value on the larger feasible set $B_2(\vzero, r^+)$ must be the same as this global minimum.

To apply CGMT, namely Lemma~\ref{lem:common-cgmt}, we perform a compactification argument; see Section~\ref{sec:cgmt} for an explanation of this argument. Take $H_n\uparrow \infty$ be any sequence of integers (which will be determined later). Then for any $(\mathbb X, \bxi)$, there exists $n_0(\mathbb X,\bxi)\in\bN_+$ such that $H_n>\frac{\|\boldsymbol\xi\|_2+r^+\|\mathbb X\|_{\mathrm{op}}}{N\lambda}$ for any $n\geq n_0(\mathbb X,\bxi)$. Then for every $\vw\in\bR^p$, $\|\vw\|_2\leq t$ for some $t>0$, \eqref{eq:variational_representation_ell_2_square} gives that the maximizer $\vu(\vw)$ of $\max(\langle\vu,\bxi-\mathbb X\vw\rangle-\frac{N\lambda}{2}\|\vu\|_2^2:\vu\in\bR^N)$ by $\vu(\vw):=\frac{1}{N\lambda}(\bxi-\mathbb X\vw)$, hence $\|\vu(\vw)\|_2<H_n$. Since the choice of $\vw$ is arbitrary over $\vw\in B_2(\vzero;t)$, together with \eqref{eq:primary-optimization_lower}, gives $\Phi_{t,H_n}=\Phi_{t,\infty}$ for each $t\in(r^-,r^+)$, that is, $u(\vw)$, the global maximizer, is already contained in $B_2(\vzero;H_n)$. Here, $B_2(\vzero,t)\times B_2(\vzero,H_n)$ is compact, hence satisfies the assumption of CGMT. Thus $\{\|\hat\bbeta[\vy]-\bbeta^\star\|_2 \le r^-\}\subseteq\liminf_{n\to\infty} \{ \Phi_{r^+,H_n}=\Phi_{r^-,H_n}\}$, and by Fatou's lemma,
\begin{align*}
&\mathbb P_{\mathbb X,\boldsymbol\xi}\{\{\|\hat\bbeta[\vy]-\bbeta^\star\|_2^2 \le 2^{-30} \Var_{J_\star^c}\}\}\leq\mathbb P_{\mathbb X,\boldsymbol\xi}\big\{\liminf_{n\to\infty}\mathcal \{ \Phi_{r^+,H_n}=\Phi_{r^-,H_n}\}\big\}\leq\liminf_{n\to\infty}\mathbb P_{\mathbb X,\boldsymbol\xi}\{\Phi_{r^+,H_n}=\Phi_{r^-,H_n}\}.
\end{align*}Thus we only need to bound the probability $\mathbb P_{\mathbb X,\boldsymbol\xi}\{\Phi_{r^+,H_n}=\Phi_{r^-,H_n}\}$ from above. Note that we have not yet used CGMT. In order to apply CGMT, in the following proof we condition on $\bxi$ and apply CGMT only to $\mathbb X$. It is easy to see that $\{ \Phi_{r^+,H_n}=\Phi_{r^-,H_n}\}\subset \{ \Phi_{r^-,H_n}\leq x \} \cup \{\Phi_{r^+,H_n}\geq x\}$ for any $x>0$. Below, we only need to choose a fixed $x$ (once conditioned on $\bxi$) such that both random events on the right-hand side have small probabilities.

\paragraph{Stochastic Arguments.} We first present the stochastic arguments needed for this proof. Let $\vh\in V_{J_\star^c}$ be a standard Gaussian random vector and $G$ be a standard Gaussian random variable. Let \begin{align*}
&\boldsymbol\Delta_{\rho,u}(\boldsymbol h)=\big(S_u(\rho\beta_j^\star-h_j)-\rho\beta_j^\star\big)_{j\in J_\star^c},\qquad v_{\rho,u}=\Phi_G^{-1}\bigg(\frac1{2|J_\star^c|}\sum_{j\in J_\star^c}\mathbb P\{|G+\rho\beta_j^\star|>u\}\bigg),\mbox{ and}\\ 
&I_{\rho,u}(\boldsymbol h)=(j\in J_\star^c:|\rho\beta_j^\star-h_j|>u).
\end{align*} The following lemma collects the stochastic arguments needed for this proof. 
\begin{Lemma}\label{lem:stochastic_argument_lower_bound_free_subspace} 
Grant the same assumptions as in Proposition~\ref{prop:free-variance-lower-bound}. Define the primary probability space random event
\begin{align*}
    &\Omega_\xi=\bigg\{\boldsymbol\xi:\sigma_\xi\sqrt N\bigg[1-2^{-13}\bigg(\frac qN+z_q^{-2}\bigg)\bigg]\leq\|\boldsymbol\xi\|_2\leq2\sigma_\xi\sqrt N\bigg\}.
\end{align*}Then $\bP^{\otimes N}(\Omega_\xi^c)\leq 2\exp(-2^{-27}d_q)$. Fix $\bxi\in\Omega_\xi$, we define the following auxiliary probability space random events.
Define $\Omega_{\mathrm{prof},1}^{\mathrm{aux}},\Omega_{\mathrm{prof},2}^{\mathrm{aux}}\subset\mathbb R^p$ by \begin{align*} &\Omega_{\mathrm{prof},1}^{\mathrm{aux}}=\bigcap_{\substack{\rho,u\geq0:\ v_{\rho,u}\leq\bar z_q}}\bigg\{\boldsymbol h:\|\boldsymbol\Delta_{\rho,u}(\boldsymbol h)\|_2^2\geq\frac{d_q}{256}\bigg\},\\ 
&\Omega_{\mathrm{prof},2}^{\mathrm{aux}}=\bigcap_{\substack{\rho,u\geq0:\ v_{\rho,u}>\bar z_q}}\bigg\{\boldsymbol h:\big[-\langle\boldsymbol h_{J_\star^c},\boldsymbol\Delta_{\rho,u}(\boldsymbol h)\rangle\big]_+\leq\frac{8q}{\sqrt{d_q}}\|\boldsymbol\Delta_{\rho,u}(\boldsymbol h)\|_2,\ |I_{\rho,u}(\boldsymbol h)|\leq\lfloor2q\rfloor\bigg\}\\ 
&\qquad\qquad\cap\bigcap_{\substack{A\subset J_\star^c:\ |A|\leq\lfloor2q\rfloor}}\bigg\{\boldsymbol h:\|\boldsymbol h_A\|_2\leq4z_q\sqrt q,\, \frac{|\langle\boldsymbol h_{J_\star^c\setminus A},\boldsymbol\beta_{J_\star^c\setminus A}^\star\rangle|}{\|\boldsymbol\beta_{J_\star^c\setminus A}^\star\|_2}\leq4z_q\sqrt q\bigg\}, \end{align*} 
where we let $\frac{|\langle\boldsymbol h_{J_\star^c\setminus A},\boldsymbol\beta_{J_\star^c\setminus A}^\star\rangle|}{\|\boldsymbol\beta_{J_\star^c\setminus A}^\star\|_2}=0$, when its denominator vanishes. Then
\begin{align} \begin{aligned} 
&\mathbb P_{\mathrm{aux}}\big((\Omega_{\mathrm{prof},1}^{\mathrm{aux}}\cap\Omega_{\mathrm{prof},2}^{\mathrm{aux}})^c\big)\leq2^7\exp(-2^{-20}d_q). 
\end{aligned} \label{eq:profile-event-bound} \end{align} 
Define 
\begin{align*} 
&\Omega_g^{\mathrm{aux}}(\bxi)=\bigg\{\vg:\, \frac N2\leq\|\boldsymbol g\|_2^2\leq2N,\ |\langle\boldsymbol\xi,\boldsymbol g\rangle|\leq2^{-7}\|\boldsymbol\xi\|_2\sqrt{d_q},\ \langle\boldsymbol\xi,\boldsymbol g\rangle\geq-2^{-15}\sqrt{\mathrm{Var}_{J_\star^c}}\|\boldsymbol g\|_2^2\bigg\},\\ 
&\Omega_h^{\mathrm{aux}}=\big\{\boldsymbol h:\|\boldsymbol h_{J_\star}\|_2\leq2\sqrt N,\ \|\boldsymbol h\|_2\leq2\sqrt p\big\},\, \Omega_{\mathrm{tie}}^{\mathrm{aux}}=\bigcap_{1\leq i<j\leq p}\{|h_i|\neq|h_j|\},\\ 
&\Omega_0^{\mathrm{aux}}(\bxi)=\big\{(\boldsymbol g,\boldsymbol h):\|\boldsymbol\xi+r^+\boldsymbol g\|_2+\langle\boldsymbol h,\boldsymbol w_0\rangle>0\big\},\mbox{ where } \boldsymbol w_0\in\mathop{\arg\min}_{\|\boldsymbol w\|_2\leq r^+}\|\boldsymbol\beta^\star+\boldsymbol w\|_1,\\ 
&\Omega_{\mathrm{aux}}^{\mathrm{aux}}(\boldsymbol\xi)=\big\{(\boldsymbol g,\boldsymbol h)\in\mathbb R^N\times\bR^p:\vg\in\Omega_g^{\mathrm{aux}}(\bxi),\ \boldsymbol h\in\Omega_h^{\mathrm{aux}}\cap\Omega_{\mathrm{tie}}^{\mathrm{aux}}\cap\Omega_{\mathrm{prof},1}^{\mathrm{aux}}\cap\Omega_{\mathrm{prof},2}^{\mathrm{aux}},\ (\boldsymbol g,\boldsymbol h)\in\Omega_0^{\mathrm{aux}}(\bxi)\big\}.
\end{align*} 
Then \begin{align} \begin{aligned} 
\sup_{\boldsymbol\xi\in\Omega_\xi}\mathbb P_{\mathrm{aux}}(\Omega_{\boldsymbol g,\boldsymbol h}^{\mathrm{aux}}(\boldsymbol\xi)^c)\leq2^{8}\exp\big (-2^{-35}d_q\big). 
\end{aligned} \label{eq:stochastic-event-bounds} 
\end{align}
\end{Lemma}

\paragraph{Choice of level \texorpdfstring{$x_H(\bxi)$}{x}.}
To make a good choice of $x = x_H(\bxi)$, we first introduce some definitions. For $\boldsymbol w\in\mathbb R^p$, $s\geq0$, and $(\boldsymbol g,\boldsymbol h)\in\mathbb R^{N+p}$, write
\begin{align*}
&\Psi(\vw,s):=\Psi_{\vg,\vh}(\boldsymbol w,s)=s\|\boldsymbol\xi+\|\boldsymbol w\|_2\boldsymbol g\|_2-\frac{N\lambda s^2}{2}+s\langle\boldsymbol h,\boldsymbol w\rangle+\|\boldsymbol\beta^\star+\boldsymbol w\|_1-\|\boldsymbol\beta^\star\|_1,\mbox{ and}\\
&\tilde\Psi(\vw,s):=\widetilde\Psi_{\vg,\vh}(\boldsymbol w,s)=s\|\boldsymbol\xi+r^+\boldsymbol g\|_2-\frac{N\lambda s^2}{2}+s\langle\boldsymbol h,\boldsymbol w\rangle+\|\boldsymbol\beta^\star+\boldsymbol w\|_1-\|\boldsymbol\beta^\star\|_1.
\end{align*}
Namely, we only replace $\|\vw\|_2$ by $r^+$, which greatly simplifies the analysis of $\tilde\Psi$. Define
\begin{align}
\begin{aligned}
&\phi_{t,H}=\min_{\|\boldsymbol w\|_2\leq t}\max_{0\leq s\leq H}\Psi_{\boldsymbol w,s}(\boldsymbol g,\boldsymbol h),\mbox{ and } \widetilde\varphi_H(\boldsymbol g,\boldsymbol h)=\min_{\|\boldsymbol w\|_2\leq r^+}\max_{0\leq s\leq H}\widetilde\Psi_{\boldsymbol w,s}(\boldsymbol g,\boldsymbol h).
\end{aligned}
\label{eq:auxiliary-optimization_lower}
\end{align}
The definition and existence of $x_H(\bxi)$ are given by the following lemma.
\begin{Lemma}\label{lem:def_x_H}
    Grant the same assumptions as Proposition~\ref{prop:free-variance-lower-bound}. For $x\in\mathbb R$, define
    \begin{align*}
    &A_x^{\mathrm{aux}}=\big\{(\boldsymbol g,\boldsymbol h):\phi_{r^-,H}\leq x\big\},\, B_x^{\mathrm{aux}}=\big\{(\boldsymbol g,\boldsymbol h):\widetilde\varphi_H(\boldsymbol g,\boldsymbol h)\geq x\big\}\cap\Omega_{\boldsymbol g,\boldsymbol h}^{\mathrm{aux}}(\boldsymbol\xi).
    \end{align*}Then for any $H>\frac{2\sigma_\xi\sqrt N+r^+(\sqrt{2N}+2\sqrt p)}{N\lambda}$,
    \begin{align}
    \begin{aligned}
    \forall x\in\bR,\qquad &\mathbb P_{\mathrm{aux}}\big\{\widetilde\varphi_H(\boldsymbol g,\boldsymbol h)=x,\Omega_{\boldsymbol g,\boldsymbol h}^{\mathrm{aux}}(\boldsymbol\xi)\mid\boldsymbol\xi\big\}=0.
    \end{aligned}
    \label{eq:no-conditional-atom}
    \end{align}Therefore, there exists $x_H(\bxi)$ such that\begin{align*}
    &x_H(\boldsymbol\xi)=\inf\big(x:\mathbb P_{\mathrm{aux}}(B_x^{\mathrm{aux}}\mid\boldsymbol\xi)\leq\Phi_G\big(2^{-8}\sqrt{d_q}\big)\big).
    \end{align*}
\end{Lemma}
With the definition of $x_H(\bxi)$ in place, we shall prove the following lemma.
\begin{Lemma}\label{lem:auxiliary_lower_bound}
    Grant the same assumptions as Proposition~\ref{prop:free-variance-lower-bound}. Then for any $H$ as defined in Lemma~\ref{lem:def_x_H}, there holds
    \begin{align}
    \begin{aligned}
    &\mathbb P_{\mathrm{aux}}\big\{\phi_{r^-,H}\leq x_H(\boldsymbol\xi)\mid\boldsymbol\xi\big\}+\mathbb P_{\mathrm{aux}}\big\{\widetilde\varphi_H(\boldsymbol g,\boldsymbol h)\geq x_H(\boldsymbol\xi)\mid\boldsymbol\xi\big\}\leq2^{9}\exp\big(-2^{-35}d_q\big),
    \end{aligned}
    \label{eq:conditional-auxiliary-tails}
    \end{align}and on $\Omega^{\mathrm{aux}}$,
    \begin{align}
    \begin{aligned}
    &\phi_{r^+,H}\leq\widetilde\varphi_H(\boldsymbol g,\boldsymbol h).
    \end{aligned}
    \label{eq:large-ball-auxiliary-bound}
    \end{align}
\end{Lemma}
With Lemma~\ref{lem:auxiliary_lower_bound}, we only need to choose $(H_n)_{n \in \mathbb{N}_+}$ and using CGMT, that is, Lemma~\ref{lem:common-cgmt}, to examine the upper bound of $\mathbb{P}_{\mathbb X}\{\Phi_{r^-,H_n} \le x_{H_n}(\bxi) \mid \bxi\} + \mathbb{P}_{\mathbb X}\{\Phi_{r^+,H_n} \ge x_{H_n}(\bxi) \mid \bxi\}$. Here we choose $H_n$ in the following way:
\begin{align*}
&H_n\geq\frac{2\sigma_\xi\sqrt N+r^+(\sqrt{2N}+2\sqrt p)+n}{N\lambda}.
\end{align*}
Applying both conclusions of Lemma~\ref{lem:common-cgmt} to $x_{H_n}(\boldsymbol\xi)\pm\varepsilon$ and then letting $\varepsilon\downarrow0$ yields
\begin{align*}
&\mathbb P_{\mathbb X}\big\{\Phi_{r^-,H_n}\leq x_{H_n}(\boldsymbol\xi)\mid\boldsymbol\xi\big\}\leq2\mathbb P_{\mathrm{aux}}\big\{\phi_{r^-,H_n}\leq x_{H_n}(\boldsymbol\xi)\mid\boldsymbol\xi\big\},\mbox{ and }\\
&\mathbb P_{\mathbb X}\big\{\Phi_{r^+,H_n}\geq x_{H_n}(\boldsymbol\xi)\mid\boldsymbol\xi\big\}\leq2\mathbb P_{\mathrm{aux}}\big\{\phi_{r^+,H_n}\geq x_{H_n}(\boldsymbol\xi)\mid\boldsymbol\xi\big\}.
\end{align*} From Lemma~\ref{lem:auxiliary_lower_bound},
\begin{align}
\begin{aligned}
&\mathbb P_{\mathbb X}\big\{\Phi_{r^-,H_n}\leq x_{H_n}(\boldsymbol\xi)\mid\boldsymbol\xi\big\}+\mathbb P_{\mathbb X}\big\{\Phi_{r^+,H_n}\geq x_{H_n}(\boldsymbol\xi)\mid\boldsymbol\xi\big\} \leq3\cdot2^{9}\exp\big(-2^{-35}d_q\big).
\end{aligned}
\label{eq:conditional-primary-tails}
\end{align}
Integrating \eqref{eq:conditional-primary-tails} over $\bxi\in\Omega_\xi$, using \eqref{eq:stochastic-event-bounds} together with taking the liminf, we finish the proof of Proposition~\ref{prop:free-variance-lower-bound}.
\end{proof}
We are left with the verification of Lemma~\ref{lem:stochastic_argument_lower_bound_free_subspace}, Lemma~\ref{lem:def_x_H}, and Lemma~\ref{lem:auxiliary_lower_bound}. Before that, we first introduce some technical lemmas.

\paragraph{Proof of Lemma~\ref{lem:stochastic_argument_lower_bound_free_subspace}.}

\begin{enumerate}
\item Recall that $\boldsymbol\Delta_{\rho,u}(\boldsymbol h)=\big(S_u(\rho\beta_j^\star-h_j)-\rho\beta_j^\star\big)_{j\in J_\star^c}$, and $ v_{\rho,u}=\Phi_G^{-1}\big(\frac1{2|J_\star^c|}\sum_{j\in J_\star^c}\mathbb P\{|G+\rho\beta_j^\star|>u\}\big)$. We first prove
\begin{align*}
    \bP_{\mathrm{aux}}\big(\Omega_{\mathrm{prof},1}^{\mathrm{aux}}\big)\geq 1-2^6\exp(-2^{-20}d_q),\mbox{ where }\Omega_{\mathrm{prof},1}^{\mathrm{aux}}=\bigcap_{\substack{\rho,u\geq0:\ v_{\rho,u}\leq\bar z_q}}\bigg\{\boldsymbol h:\|\boldsymbol\Delta_{\rho,u}(\boldsymbol h)\|_2^2\geq\frac{d_q}{256}\bigg\}.
\end{align*}
Our proof strategy is to apply Lemma~\ref{lem:probabilistic_lemma}; to this end, we first need to construct the family of Borel subsets in Lemma~\ref{lem:probabilistic_lemma}. We fix $\rho, u \ge 0$ and $j \in J_\star^c$. We first prove the following displayed formula
\begin{align*}
    &C_j:=C(\rho,u,j):=\big(y\in\bR:\,  \big| S_u\big( \rho\beta_j^\star - y \big)-\rho\beta_j^\star \big| \geq \frac{1}{6z_q} \big)=C_1(\rho,u,j)\sqcup C_2(\rho,u,j)\sqcup C_3(\rho,u,j)\\
    &:=\big(y:\,  y<\rho\beta_j^\star - u,\, |y+u|\geq \frac{1}{6z_q} \big)\sqcup \big(y:\, y>\rho\beta_j^\star + u,\, |y-u|\geq \frac{1}{6z_q} \big)\sqcup \big(y:\, |y-\rho\beta_j^\star| \leq u,\, |\rho\beta_j^\star|\geq \frac{1}{6z_q} \big).
\end{align*}Indeed, if $y<\rho\beta_j^\star-u$, then $S_u(\rho\beta_j^\star-y)-\rho\beta_j^\star = -y-u$; if $|y-\rho\beta_j^\star|\leq u$, then $S_u(\rho\beta_j^\star-y)-\rho\beta_j^\star=-\rho\beta_j^\star$; and if $y>\rho\beta_j^\star+u$, then $S_u(\rho\beta_j^\star-y)-\rho\beta_j^\star=u-y$. Combining these pairwise disjoint events, the above displayed formula is verified.

Now let $\cC_1=\big( (C(\rho,u,m))_{m\in J_\star^c} :\, \rho,u\geq 0\big).$ This is precisely the $\mathcal{C}$ we choose in Lemma~\ref{lem:probabilistic_lemma}. Next we verify the assumption of Lemma~\ref{lem:probabilistic_lemma}, that is,
\begin{align*}
    \sup_{(y_j,\widetilde y_j)_{j\in J_\star^c}\in\mathbb R^{2|J_\star^c|}}\bigg|\bigg(\big(\boldsymbol1_{\{y_j\in C_j\}},\boldsymbol1_{\{\widetilde y_j\in C_j\}}\big)_{j\in J_\star^c}:C\in\mathcal C_1\bigg)\bigg|.
\end{align*}
Fix $(y_j,\widetilde y_j)_{j\in J_\star^c}\in\mathbb R^{2|J_\star^c|}$. For each $j\in J_\star^c$, the pair $(\boldsymbol1_{\{y_j\in C_j\}},\boldsymbol1_{\{\widetilde y_j\in C_j\}})$ is a function of $(\rho,u)$. Indeed, it is determined by the signs of the twelve affine functions obtained, for $r=y_j,\widetilde y_j$, from
\begin{align*}
&r-\rho\beta_j^\star\pm u,\qquad r\pm u-(6z_q)^{-1},\qquad r\pm u+(6z_q)^{-1}.
\end{align*}
Together with the signs of $\rho\beta_j^\star\pm(6z_q)^{-1}$, these signs determine the pair $(\boldsymbol1_{\{y_j\in C_j\}},\boldsymbol1_{\{\widetilde y_j\in C_j\}})$. Summarizing, the vector $(\boldsymbol1_{\{y_j\in C_j\}},\boldsymbol1_{\{\widetilde y_j\in C_j\}})_{j\in J_\star^c}$ is determined by by at most $14|J_\star^c|$ affine lines in the $(\rho,u)$-plane. Therefore, to obtain this cardinality upper bound, we only need to estimate the total number of 2-faces (open polygonal regions enclosed by the lines), 1-faces (segments or rays obtained by cutting a line by other lines), and 0-faces (intersection points of several lines) of the polyhedral partition of the $(\rho, u)$-plane induced by these $14|J_\star^c|$ lines --- since on the interior of each 2-face, 1-face, and 0-face, the indicators $(\mathbf{1}_{\{y_j \in C\}}, \mathbf{1}_{\{\tilde y_j \in C_j\}})_{j \in J_\star^c}$ remain constant. An arrangement of $L$ affine lines in $\mathbb R^2$ has at most $1+L+\binom{L}{2}$ two-dimensional faces, at most $L^2$ one-dimensional faces, and at most $\binom{L}{2}$ zero-dimensional faces, hence at most $1+2L^2$ relatively open faces. Since $|J_\star^c|\geq4q\geq480$,
\begin{align*}
&\sup_{(y_j,\widetilde y_j)_{j\in J_\star^c}\in\mathbb R^{2|J_\star^c|}}\bigg|\bigg(\big(\boldsymbol1_{\{y_j\in C_j\}},\boldsymbol1_{\{\widetilde y_j\in C_j\}}\big)_{j\in J_\star^c}:C\in\mathcal C_1\bigg)\bigg|\leq1+392|J_\star^c|^2\leq\bigg(\frac{e|J_\star^c|}{8}\bigg)^{16}.
\end{align*}
The last step before applying Lemma~\ref{lem:probabilistic_lemma} is to examine the lower bound on the probability $\sum_{j \in J_\star^c} \mathbb{P}\{h_j \in C_j\}$, that is, $\sum_{j\in J_\star^c}\bP\{|S_u(\rho\beta_j^\star-h_j)-\rho\beta_j^\star|\geq\frac{1}{6z_q}\}$.
By Lemma~\ref{lem:v_rho_u}, whenever $v_{\rho,u}\leq\bar z_q$,
\begin{align*}&\sum_{j\in J_\star^c}\mathbb P\bigg\{|S_u(\rho\beta_j^\star-G)-\rho\beta_j^\star|\geq\frac1{6z_q}\bigg\}>\frac34q.\end{align*}
Applying the first conclusion of Lemma~\ref{lem:probabilistic_lemma} with $V=16$, $x=\frac{3q}{4}$, $Z_j=h_j$, and $M=|J_\star^c|$, and using \eqref{eq:common-calibration-identities} together with $d_q>2^{35}\log(2^{11})$, gives
\begin{align*}\begin{aligned}&\mathbb P_{\mathrm{aux}}\bigg\{\exists\rho,u\geq0:v_{\rho,u}\leq\bar z_q,\ \sum_{j\in J_\star^c}\boldsymbol1_{\{|S_u(\rho\beta_j^\star-h_j)-\rho\beta_j^\star|\geq(6z_q)^{-1}\}}<\frac23\sum_{j\in J_\star^c}\mathbb P\bigg\{|S_u(\rho\beta_j^\star-G)-\rho\beta_j^\star|\geq\frac1{6z_q}\bigg\}\bigg\}\\&\leq2\bigg(\frac{e|J_\star^c|}{8}\bigg)^{16}\exp\bigg(-\frac{q}{2304}\bigg)\leq2^6\exp(-2^{-20}d_q).\end{aligned}\end{align*}
Consequently, with probability at least $1-2^6\exp(-2^{-20}d_q)$, simultaneously for every $\rho,u\geq0$ satisfying $v_{\rho,u}\leq\bar z_q$,
\begin{align*}&\sum_{j\in J_\star^c}\boldsymbol1_{\{|S_u(\rho\beta_j^\star-h_j)-\rho\beta_j^\star|\geq(6z_q)^{-1}\}}\geq\frac23\sum_{j\in J_\star^c}\mathbb P\bigg\{|S_u(\rho\beta_j^\star-G)-\rho\beta_j^\star|\geq\frac1{6z_q}\bigg\}>\frac q2.\end{align*}
Therefore
\begin{align*}&\|\boldsymbol\Delta_{\rho,u}(\boldsymbol h)\|_2^2>\frac q2\frac1{36z_q^2}=\frac{q}{72z_q^2}\geq\frac{d_q}{144}>\frac{d_q}{256},\end{align*}
where the first inequality in the last comparison follows from \eqref{eq:common-calibration-identities}. This proves the defining inequality of $\Omega_{\mathrm{prof},1}^{\mathrm{aux}}$. Taking $\rho=0$ and $u=\bar z_q$ gives $\sum_{j\in J_\star^c}(|h_j|-\bar z_q)_+^2\geq \frac{d_q}{256}$.

\item\label{item:profile-prof2} Recall that, for every $\rho,u\geq0$, $I_{\rho,u}(\boldsymbol h)=(j\in J_\star^c:|\rho\beta_j^\star-h_j|>u)$. In this item, we prove that with probability at least $1-2^5\exp(-2^{-20}d_q)$, simultaneously for any $\rho,u\geq 0$, $v_{\rho,u}>\bar z_q$, $|I_{\rho,u}(\vh)|\leq \lfloor 2q\rfloor$. For $v_{\rho,u}>\bar z_q$,
\begin{align*}&\sum_{j\in J_\star^c}\mathbb P\{|G+\rho\beta_j^\star|>u\}=2|J_\star^c|\Phi_G(v_{\rho,u})<2|J_\star^c|\Phi_G(z_q)=q.\end{align*}
Consider
\begin{align*}&\mathcal C_2=\bigg(\bigg(\bigg(y\in\mathbb R:|\rho\beta_j^\star-y|>u\bigg)\bigg)_{j\in J_\star^c}:\rho,u\geq0\bigg).\end{align*}
Fix $(y_j,\widetilde y_j)_{j\in J_\star^c}\in\mathbb R^{2|J_\star^c|}$. For every $j\in J_\star^c$, the pair $\big(\boldsymbol1_{\{|\rho\beta_j^\star-y_j|>u\}},\boldsymbol1_{\{|\rho\beta_j^\star-\widetilde y_j|>u\}}\big)$ is determined by the signs of the four affine functions $\rho\beta_j^\star-y_j-u$, $\rho\beta_j^\star-y_j+u$, $\rho\beta_j^\star-\widetilde y_j-u$, and $\rho\beta_j^\star-\widetilde y_j+u$. Hence the vector $\big(\boldsymbol1_{\{|\rho\beta_j^\star-y_j|>u\}},\boldsymbol1_{\{|\rho\beta_j^\star-\widetilde y_j|>u\}}\big)_{j\in J_\star^c}$ is constant on each relatively open face generated by at most $4|J_\star^c|$ affine lines. By the preceding face-counting argument and $|J_\star^c|\geq480$,
\begin{align*}&\sup_{(y_j,\widetilde y_j)_{j\in J_\star^c}\in\mathbb R^{2|J_\star^c|}}\bigg|\bigg(\big(\boldsymbol1_{\{|\rho\beta_j^\star-y_j|>u\}},\boldsymbol1_{\{|\rho\beta_j^\star-\widetilde y_j|>u\}}\big)_{j\in J_\star^c}:\rho,u\geq0\bigg)\bigg|\leq1+32|J_\star^c|^2\leq\bigg(\frac{e|J_\star^c|}{2}\bigg)^4.\end{align*}
Applying the second conclusion of Lemma~\ref{lem:probabilistic_lemma} with $M=|J_\star^c|$, $V=4$, $x=q$, and $Z_j=h_j$, and using \eqref{eq:common-calibration-identities} together with $d_q>2^{35}\log(2^{11})$, gives
\begin{align*}&\mathbb P_{\mathrm{aux}}\bigg\{\exists\rho,u\geq0:v_{\rho,u}>\bar z_q,\ |I_{\rho,u}(\boldsymbol h)|>2q\bigg\}\leq3\bigg(\frac{e|J_\star^c|}{2}\bigg)^4e^{-\frac{q}{196}}\leq2^5\exp(-2^{-20}d_q).\end{align*}
Thus, with probability at least $1-2^5\exp(-2^{-20}d_q)$, simultaneously for every $\rho,u\geq0$ satisfying $v_{\rho,u}>\bar z_q$, $|I_{\rho,u}(\boldsymbol h)|\leq\lfloor2q\rfloor$.

\item We prove the remaining claims in $\Omega_{\mathrm{prof},2}^{\mathrm{aux}}$. For every fixed $A\subset J_\star^c$ with $|A|\leq\lfloor2q\rfloor$, if $\|\boldsymbol\beta_{J_\star^c\setminus A}^\star\|_2=0$, then $\langle\boldsymbol h_{J_\star^c\setminus A},\boldsymbol\beta_{J_\star^c\setminus A}^\star\rangle=0$; otherwise, $\frac{\langle\boldsymbol h_{J_\star^c\setminus A},\boldsymbol\beta_{J_\star^c\setminus A}^\star\rangle}{\|\boldsymbol\beta_{J_\star^c\setminus A}^\star\|_2}\sim\cN(0,1)$. Moreover, $\mathbb E\exp(t\|\boldsymbol h_A\|_2^2)=(1-2t)^{-\frac{|A|}{2}}$ for $0<t<\frac{1}{2}$, see, for instance, \cite[Section 2.4, pp. 28]{boucheronbook}. Taking $t=\frac{1}{4}$ and using the one-dimensional Gaussian tail gives
\begin{align*}&\mathbb P\{\|\boldsymbol h_A\|_2>4z_q\sqrt q\}\leq e^{-4z_q^2q+q\log2},\qquad \mathbb P\bigg\{\frac{|\langle\boldsymbol h_{J_\star^c\setminus A},\boldsymbol\beta_{J_\star^c\setminus A}^\star\rangle|}{\|\boldsymbol\beta_{J_\star^c\setminus A}^\star\|_2}>4z_q\sqrt q\bigg\}\leq2e^{-8qz_q^2},\end{align*}
where the second probability is zero when $\|\boldsymbol\beta_{J_\star^c\setminus A}^\star\|_2=0$. Since $\sum_{k\leq\lfloor2q\rfloor}\binom{|J_\star^c|}{k}\leq\exp\big(2q\log\frac{e|J_\star^c|}{q}\big)\leq e^{2qz_q^2}$ a union bound gives
\begin{align*}\begin{aligned}&\mathbb P\bigg\{\exists A\subset J_\star^c,\ |A|\leq\lfloor2q\rfloor:\|\boldsymbol h_A\|_2>4z_q\sqrt q\bigg\}\\&+\mathbb P\bigg\{\exists A\subset J_\star^c,\ |A|\leq\lfloor2q\rfloor:\frac{|\langle\boldsymbol h_{J_\star^c\setminus A},\boldsymbol\beta_{J_\star^c\setminus A}^\star\rangle|}{\|\boldsymbol\beta_{J_\star^c\setminus A}^\star\|_2}>4z_q\sqrt q\bigg\}\leq4e^{-2^{-4}qz_q^2}.\end{aligned}\end{align*}
Fix $\rho,u\geq0$ with $v_{\rho,u}>\bar z_q$. Since $|I_{\rho,u}(\boldsymbol h)|\leq\lfloor2q\rfloor$, the preceding two simultaneous estimates apply with $A=I_{\rho,u}(\boldsymbol h)$. For every $j\in J_\star^c\setminus I_{\rho,u}(\boldsymbol h)$, $|\rho\beta_j^\star-h_j|\leq u$, hence $S_u(\rho\beta_j^\star-h_j)=0$ and $\Delta_{\rho,u,j}(\boldsymbol h)=-\rho\beta_j^\star$. Decomposing the inner product over $I_{\rho,u}(\boldsymbol h)$ and $J_\star^c\setminus I_{\rho,u}(\boldsymbol h)$, using $(a+b)_+\leq|a|+|b|$, the Cauchy--Schwarz inequality, and \eqref{eq:common-calibration-identities}, gives
\begin{align*}\begin{aligned}&\big(-\langle\boldsymbol h_{J_\star^c},\boldsymbol\Delta_{\rho,u}(\boldsymbol h)\rangle\big)_+\leq\|\boldsymbol h_{I_{\rho,u}(\boldsymbol h)}\|_2\|\boldsymbol\Delta_{\rho,u,I_{\rho,u}(\boldsymbol h)}(\boldsymbol h)\|_2+\rho\big|\langle\boldsymbol h_{J_\star^c\setminus I_{\rho,u}(\boldsymbol h)},\boldsymbol\beta^\star_{J_\star^c\setminus I_{\rho,u}(\boldsymbol h)}\rangle\big|\\&\leq4z_q\sqrt q\big(\|\boldsymbol\Delta_{\rho,u,I_{\rho,u}(\boldsymbol h)}(\boldsymbol h)\|_2+\|\boldsymbol\Delta_{\rho,u,J_\star^c\setminus I_{\rho,u}(\boldsymbol h)}(\boldsymbol h)\|_2\big)\\&\leq4\sqrt2z_q\sqrt q\|\boldsymbol\Delta_{\rho,u}(\boldsymbol h)\|_2\leq\frac{8q}{\sqrt{d_q}}\|\boldsymbol\Delta_{\rho,u}(\boldsymbol h)\|_2.\end{aligned}\end{align*}
Thus, with probability at least $1-2^5\exp(-2^{-20}d_q)-4e^{-2^{-4}qz_q^2}$, all defining inequalities of $\Omega_{\mathrm{prof},2}^{\mathrm{aux}}$ hold. Finally, by \eqref{eq:common-calibration-identities} and $d_q>2^{35}\log(2^{11})$,
\begin{align*}&2^6e^{-2^{-20}d_q}+2^5e^{-2^{-20}d_q}+4e^{-2^{-4}qz_q^2}\leq2^7e^{-2^{-20}d_q}.\end{align*}
The union bound proves \eqref{eq:profile-event-bound}.
\end{enumerate}
\endproof

We shall also use that on $\Omega^{\mathrm{aux}}$, \begin{align}
\begin{aligned}
&0\leq\frac{\langle\boldsymbol\xi+r^+\boldsymbol g,\boldsymbol g\rangle}{\|\boldsymbol\xi+r^+\boldsymbol g\|_2}\leq2^{-6}\sqrt{d_q},\mbox{ and } \|\boldsymbol\xi+r^+\boldsymbol g\|_2-r^+\bigg(\frac{8q}{\sqrt{d_q}}+\|\boldsymbol h_{J_\star}\|_2\bigg)>\frac{N\lambda}{\bar z_q}.
\end{aligned}
\label{eq:auxiliary-hypotheses}
\end{align}
Indeed, on $\Omega_g^{\mathrm{aux}}$, $\langle\boldsymbol\xi+r^+\boldsymbol g,\boldsymbol g\rangle\geq \frac{r^+}{2\|\boldsymbol g\|_2^2}\geq0$, and $ \|\boldsymbol\xi+r^+\boldsymbol g\|_2\geq\|\boldsymbol\xi\|_2.$ Consequently, by \eqref{eq:common-calibration-identities}, $\|\boldsymbol g\|_2^2\leq2N$, and $\boldsymbol\xi\in\Omega_\xi$,
\begin{align*}
&\frac{\langle\boldsymbol\xi+r^+\boldsymbol g,\boldsymbol g\rangle}{\|\boldsymbol\xi+r^+\boldsymbol g\|_2}\leq2^{-7}\sqrt{d_q}+\frac{2Nr^+}{\|\boldsymbol\xi\|_2}\leq2^{-6}\sqrt{d_q}.
\end{align*}
Furthermore,
\begin{align}
\begin{aligned}
&\frac{r^+}{\sigma_\xi\sqrt N}\bigg(\frac{8q}{\sqrt{d_q}}+2\sqrt N\bigg)\leq3\cdot2^{-12}\bigg(\frac qN+z_q^{-2}\bigg).
\end{aligned}
\label{eq:penalty-radius-bound}
\end{align}
The lower bound defining $\Omega_\xi$, \eqref{eq:penalty-radius-bound}, and \eqref{eq:common-radial-calibration} imply
\begin{align*}
&\|\boldsymbol\xi+r^+\boldsymbol g\|_2-r^+\bigg(\frac{8q}{\sqrt{d_q}}+\|\boldsymbol h_{J_\star}\|_2\bigg)\geq\sigma_\xi\sqrt N-\big(2^{-13}+3\cdot2^{-12}\big)\sigma_\xi\sqrt N\bigg(\frac qN+z_q^{-2}\bigg)>\frac{N\lambda}{\bar z_q},
\end{align*}
which proves \eqref{eq:auxiliary-hypotheses}.

\paragraph{Saddle point analysis.}

The following Lemma is the saddle point analysis of $\min_{\|\vw\|_2\leq r^+}\max_{0\leq s\leq H}\tilde\Psi$.
\begin{Lemma}\label{lem:saddle_point_analysis_tilde_Psi}
Grant the same assumptions as Proposition~\ref{prop:free-variance-lower-bound}. Fix $H>\frac{2\sigma_\xi\sqrt N+r^+(\sqrt{2N}+2\sqrt p)}{N\lambda}$. On the event $(\vg,\vh)\in \Omega^{\mathrm{aux}}$ defined in Lemma~\ref{lem:stochastic_argument_lower_bound_free_subspace}, let $\cS_H(\vg,\vh)$ be the set of all the saddle points of $(\vw,s)\in B_2(\vzero;r^+)\times [0,H]\mapsto\tilde\Psi_{\vg,\vh}(\vw,s)$, and let $(\boldsymbol w_H,s_H)\in\cS_H(\vg,\vh)$. Then there exists $\gamma_H>0$ such that the following conclusions hold.
\begin{enumerate}
\item The maximizing coordinate satisfies
\begin{align}
\begin{aligned}
&0<s_H<H,\qquad N\lambda s_H=\|\boldsymbol\xi+r^+\boldsymbol g\|_2+\langle\boldsymbol h,\boldsymbol w_H\rangle.
\end{aligned}
\tag{\ref{eq:scalar-first-order}}
\end{align}

\item There holds $\|\boldsymbol w_H\|_2=r^+$, and
\begin{align}
\begin{aligned}
&\vw_H=\frac{s_H}{\gamma_H}\bigg(S_{\frac{1}{s_H}}\bigg(\frac{\gamma_H}{s_H}\bbeta^\star-\vh\bigg)-\frac{\gamma_H}{s_H}\bbeta^\star\bigg).
\end{aligned}
\tag{\ref{eq:solution_S_Nlambda_s_H}}
\end{align}

\item The multiplier satisfies
\begin{align}
\begin{aligned}
&\gamma_Hr^+-s_H\frac{\langle\boldsymbol\xi+r^+\boldsymbol g,\boldsymbol g\rangle}{\|\boldsymbol\xi+r^+\boldsymbol g\|_2}\geq s_H(2^{-4}-2^{-6})\sqrt{d_q}\geq2^{-5}s_H\sqrt{d_q}.
\end{aligned}
\tag{\ref{eq:saddle-structure}}
\end{align}

\item The following pointwise gap and auxiliary-value gap hold:
\begin{align}
\begin{aligned}
&\forall\,\|\boldsymbol w\|_2\leq r^-,\qquad \Psi(\boldsymbol w,s_H)\geq\widetilde\phi_H(\boldsymbol g,\boldsymbol h)+2^{-6}s_Hr^+\sqrt{d_q},\mbox{ and }\phi_{r^-,H}\geq\widetilde\phi_H(\boldsymbol g,\boldsymbol h)+2^{-6}s_Hr^+\sqrt{d_q}.
\end{aligned}
\tag{\ref{eq:small-ball-auxiliary-gap}}
\end{align}
\end{enumerate}
\end{Lemma}

\beginproof
Fix $(\boldsymbol\xi,\boldsymbol g,\boldsymbol h)$. Let $\mathcal S_H(\boldsymbol g,\boldsymbol h)$ be the set of saddle points $(\boldsymbol w_H,s_H)\in B_2^p(\boldsymbol 0,r^+)\times[0,H]$ of $\tilde\Psi(\vw,s)$, that is, for every $\boldsymbol w$ in this ball and every $s\in[0,H]$, $\widetilde\Psi(\boldsymbol w_H,s)\leq\widetilde\Psi(\boldsymbol w_H,s_H)\leq\widetilde\Psi(\boldsymbol w,s_H).$ Lemma~\ref{lem:common-saddle} gives $\mathcal S_H(\boldsymbol g,\boldsymbol h)\neq\varnothing$. Fix $(\boldsymbol w_H,s_H)\in\mathcal S_H(\boldsymbol g,\boldsymbol h)$. Notice that $s\in[0,H]\mapsto\tilde\Psi(\vw,s) = s\|\boldsymbol\xi+r^+\boldsymbol g\|_2-\frac{N\lambda s^2}{2}+s\langle\boldsymbol h,\boldsymbol w\rangle+\|\boldsymbol\beta^\star+\boldsymbol w\|_1-\|\boldsymbol\beta^\star\|_1$ is maximized at $s=s_H$ due to the fact that $\max_{0\leq s\leq H}(s(\|\boldsymbol\xi+r^+\boldsymbol g\|_2+\langle\boldsymbol h,\boldsymbol w\rangle)-\frac{N\lambda s^2}{2})$ is a projection onto interval $[0,H]$, where
\begin{align*}
&s_H=\operatorname{Proj}_{[0,H]}\bigg(\frac{\|\boldsymbol\xi+r^+\boldsymbol g\|_2+\langle\boldsymbol h,\boldsymbol w_H\rangle}{N\lambda}\bigg),
\end{align*}
and the KKT condition implies that there exists some $\gamma_H\geq 0$ such that
\begin{align}
\begin{aligned}
&\boldsymbol 0\in s_H\boldsymbol h+\partial\|\boldsymbol\beta^\star+\boldsymbol w_H\|_1+\gamma_H\boldsymbol w_H,\, \|\boldsymbol w_H\|_2^2\leq(r^+)^2,\mbox{ and }\gamma_H\big(\|\boldsymbol w_H\|_2^2-(r^+)^2\big)=0,
\end{aligned}
\label{eq:auxiliary-kkt}
\end{align}where the condition $\gamma_H\big(\|\boldsymbol w_H\|_2^2-(r^+)^2\big)=0$ is from the normal cone computation, or from the so-called complementary slackness condition.
The assumption on $H$ and the projection formula imply $s_H<H$.

\textbf{We first solve $s_H=s_H(\vw_H)$}, that is, we prove that
\begin{align}
\begin{aligned}
&N\lambda s_H=\|\boldsymbol\xi+r^+\boldsymbol g\|_2+\langle\boldsymbol h,\boldsymbol w_H\rangle.
\end{aligned}
\label{eq:scalar-first-order}
\end{align}
Since we already know that $s_H = \operatorname{Proj}_{[0,H]}\big(\frac{\|\boldsymbol\xi+r^+\boldsymbol g\|_2+\langle\boldsymbol h,\boldsymbol w_H\rangle}{N\lambda}\big)$, and $s_H<H$, it only remains to prove that $\|\boldsymbol\xi+r^+\boldsymbol g\|_2+\langle\boldsymbol h,\boldsymbol w_H\rangle>0$.

The minimizer $\boldsymbol w_0\in\mathop{\arg\min}_{\|\boldsymbol w\|_2\leq r^+}\|\boldsymbol\beta^\star+\boldsymbol w\|_1$ defined in Lemma~\ref{lem:stochastic_argument_lower_bound_free_subspace} is unique. For $s>0$, choose
\begin{align*}
&\boldsymbol w_s\in\mathop{\arg\min}\big(m_s(\vw):\, \|\boldsymbol w\|_2\leq r^+\big),\, \mbox{where }m_s(\vw) = s\langle\boldsymbol h,\boldsymbol w\rangle+\|\boldsymbol\beta^\star+\boldsymbol w\|_1-\|\boldsymbol\beta^\star\|_1.
\end{align*}
Compactness and uniqueness of $\boldsymbol w_0$ give $\boldsymbol w_s\to\boldsymbol w_0$ as $s\downarrow0$. Furthermore, $$\langle\boldsymbol h,\boldsymbol w_s\rangle\leq\frac{1}{s}\big[\min_{\|\vw\|_2\leq r^+}m_s(\vw)-\min_{\|\vw\|_2\leq r^+}m_0(\vw)\big]\leq\langle\boldsymbol h,\boldsymbol w_0\rangle.$$ Indeed, since $\boldsymbol w_0$ minimizes $\boldsymbol w\mapsto\|\boldsymbol\beta^\star+\boldsymbol w\|_1$ on $B_2^p(\boldsymbol 0,r^+)$, one has $\|\bbeta^\star+\vw_s\|_1\geq\|\bbeta^\star+\vw_0\|_1$, hence $m_s(\vw_s)\geq s\langle\vh,\vw_s\rangle+\|\bbeta^\star+\vw_0\|_1-\|\bbeta^\star\|_1=s\langle\vh,\vw_s\rangle+m_0(\vw_0)$, and consequently, $\langle\boldsymbol h,\boldsymbol w_s\rangle\leq\frac{1}{s}\big[m_s(\vw_s)-m_0(\vw_0)\big]=\frac{1}{s}\big[\min_{\|\vw\|_2\leq r^+}m_s(\vw)-\min_{\|\vw\|_2\leq r^+}m_0(\vw)\big]$. Conversely, by the definition of $\vw_s$, $m_s(\vw_s)\leq m_s(\vw_0)=s\langle\vh,\vw_0\rangle+\|\bbeta^\star+\vw_0\|_1-\|\bbeta^\star\|_1=s\langle\vh,\vw_0\rangle+m_0(\vw_0)$, and consequently, $\frac{1}{s}\big[\min_{\|\vw\|_2\leq r^+}m_s(\vw)-\min_{\|\vw\|_2\leq r^+}m_0(\vw)\big]\leq\langle\vh,\vw_0\rangle$. Letting $s\downarrow0$, and using $\vw_s\to\vw_0$, gives 
\begin{align*}
    \frac{\mathrm d}{\mathrm{ d}s} \bigg( s\mapsto \min_{\|\vw\|_2\leq r^+} m_s(\vw) \bigg)\bigg|_{s=0^+} = \lim_{s\to0^+}\frac{1}{s}\big(\min_{\|\vw\|_2\leq r^+}m_s(\vw)-\min_{\|\vw\|_2\leq r^+}m_0(\vw)\big)=\langle\vh,\vw_0\rangle.
\end{align*}
Together with the definition of $\Omega_0^{\mathrm{aux}}$, there holds
\begin{align*}
    &\frac{\mathrm d}{\mathrm d s}\big[\min_{\vw\in B_2^p(\boldsymbol 0,r^+)}\widetilde\Psi(\vw,s)\big]_{s=0^+} = \frac{\mathrm{d}}{\mathrm{d}s^+}\big( s\|\bxi+r^+\vg\|_2 - \frac{N\lambda s^2}{2} \big)+ \frac{\mathrm{d}}{\mathrm{d}s^+} \big(\min_{\|\vw\|_2\leq r^+} m_s(\vw)\big)  =\|\boldsymbol\xi+r^+\boldsymbol g\|_2+\langle\boldsymbol h,\boldsymbol w_0\rangle>0
\end{align*}
Since the saddle inequalities give $\widetilde\varphi_H(\boldsymbol g,\boldsymbol h)=\max_{0\leq s\leq H}\min_{\vw\in B_2^p(\boldsymbol 0,r^+)}\widetilde\Psi(\vw,s)$, if $s_H=0$ were a maximizer, the first-order normal-cone condition for the maximization over $[0,H]$ would imply $\frac{\mathrm d}{\mathrm ds}_+\big[\min_{\vw\in B_2^p(\boldsymbol 0,r^+)}\widetilde\Psi(\vw,s)\big]_{s=0}\in N_{[0,H]}(0)=(-\infty,0]$ (see, for instance, Section~\ref{sec:tools_convex_analysis}), contradicting the strict positivity above. Hence  $s_H>0$.

\textbf{We next show that $\gamma_H>0$.} 

Suppose for contradiction that $\gamma_H=0$, then \eqref{eq:auxiliary-kkt} gives that for any $j\in[p]$, $-s_Hh_j\in\partial^-|\cdot|(\beta_j^\star+w_{H,j})$. Consequently, $|h_j|\leq\frac{1}{s_H}$, and in particular,
\begin{align}\label{eq:h_j_smaller_strict}
    \forall j\in [p],\qquad\mbox{if } |h_j| < \frac{1}{s_H},\qquad\mbox{ then necessarily } \qquad w_{H_j} = -\beta_j^\star.
\end{align}Moreover, \eqref{eq:auxiliary-kkt} can be rewritten coordinatewise as 
\begin{align*}
    0\in\frac{\gamma_H}{s_H}(\beta_j^\star+w_{H,j})-\big(\frac{\gamma_H}{s_H}\beta_j^\star-h_j\big)+\frac{1}{s_H}\partial^-|\cdot|(\beta_j^\star+w_{H,j}).
\end{align*}
By \eqref{eq:softmax_equivalence}, $\frac{\gamma_H}{s_H}(\beta_j^\star+w_{H,j})=S_{\frac{1}{s_H}}\big(\frac{\gamma_H}{s_H}\beta_j^\star-h_j\big)$, and hence substracting $\frac{\gamma_H}{s_H}\beta_j^\star$ on both sides yields
\begin{align}\label{eq:solution_S_Nlambda_s_H}
    S_{\frac{1}{s_H}}\big(\frac{\gamma_H}{s_H}\beta_j^\star-h_j\big)-\frac{\gamma_H}{s_H}\beta_j^\star=\frac{\gamma_H}{s_H}w_{H,j}.
\end{align}
Taking $\rho=\frac{\gamma_H}{s_H}$ and $u=\frac{1}{s_H}$ in the definition of $\boldsymbol\Delta_{\rho,u}(\boldsymbol h)$ therefore gives
\begin{align*}
    \boldsymbol\Delta_{\frac{\gamma_H}{s_H},\frac{1}{s_H}}(\boldsymbol h) = \boldsymbol\Delta_{0,\frac{1}{s_H}}(\boldsymbol h) = \frac{\gamma_H}{s_H}\vw_{H}=\vzero,
\end{align*}where we have used the assumption that $\gamma_H=0$.
We now prove that $\frac{1}{s_H}>\bar z_q$. Suppose to the contrary that $\frac{1}{s_H}\leq\bar z_q$. Since $\gamma_H=0$, by the definition of $v_{\rho,u}$ and the identity $\mathbb P\{|G|>u\}=2\Phi_G(u)$ for $u\geq0$, one has $v_{0,\frac{1}{s_H}}=\Phi_G^{-1}\big(\frac{1}{2|J_\star^c|}\sum_{j\in J_\star^c}\mathbb P\{|G|>\frac{1}{s_H}\}\big)=\Phi_G^{-1}\big(\Phi_G(\frac{1}{s_H})\big)=\frac{1}{s_H}\leq\bar z_q$. 
By the definition of $\Omega_{\mathrm{prof},1}^{\mathrm{aux}}$ applied at $(\rho,u)=(0,\frac{1}{s_H})$ then gives $\|\boldsymbol\Delta_{0,\frac{1}{s_H}}(\boldsymbol h)\|_2^2\geq\frac{d_q}{256}>0$, a contradiction. Therefore, necessarily $\frac{1}{s_H}>\bar z_q$.

Now consider $A=(j\in J_\star^c:|h_j|=\frac{1}{s_H})$, then almost surely, $|A|\leq 1\leq \lfloor 2q\rfloor$ (otherwise there exists $j\neq k\in J_\star^c$ such that $|h_j|=|h_k|$, which happens only on a null set); and we know that for any $j\in J_\star^c\backslash A$, $|h_j|<\frac{1}{s_H}$, and by \eqref{eq:h_j_smaller_strict}, $\boldsymbol w_{H,J_\star^c\setminus A}=-\boldsymbol\beta^\star_{J_\star^c\setminus A}$. 
As a consequence,
\begin{align*}
&-\langle\boldsymbol h_{J_\star^c},\boldsymbol w_{H,J_\star^c}\rangle\leq\|\boldsymbol h_A\|_2\|\boldsymbol w_{H,A}\|_2+\frac{|\langle\boldsymbol h_{J_\star^c\setminus A},\boldsymbol\beta^\star_{J_\star^c\setminus A}\rangle|}{\|\boldsymbol\beta^\star_{J_\star^c\setminus A}\|_2}\|\boldsymbol w_{H,J_\star^c\setminus A}\|_2.
\end{align*}
Using $(a_1b_1+a_2b_2)^2\leq(a_1^2+a_2^2)(b_1^2+b_2^2)$ gives
\begin{align*}
&-\langle\boldsymbol h_{J_\star^c},\boldsymbol w_{H,J_\star^c}\rangle\leq\bigg[\|\boldsymbol h_A\|_2^2+\frac{|\langle\boldsymbol h_{J_\star^c\setminus A},\boldsymbol\beta^\star_{J_\star^c\setminus A}\rangle|^2}{\|\boldsymbol\beta^\star_{J_\star^c\setminus A}\|_2^2}\bigg]^{\frac{1}{2}}\|\boldsymbol w_{H,J_\star^c}\|_2.
\end{align*}
By the definition of $\Omega_{\mathrm{prof},2}^{\mathrm{aux}}$ applied to this $A=(j\in J_\star^c:|h_j|=\frac{1}{s_H})$, there holds $\|\boldsymbol h_A\|_2\leq4z_q\sqrt q$, and $\frac{|\langle\boldsymbol h_{J_\star^c\setminus A},\boldsymbol\beta^\star_{J_\star^c\setminus A}\rangle|}{\|\boldsymbol\beta^\star_{J_\star^c\setminus A}\|_2}\leq4z_q\sqrt q$. Hence $-\langle\boldsymbol h_{J_\star^c},\boldsymbol w_{H,J_\star^c}\rangle\leq4\sqrt2z_q\sqrt q\,\|\boldsymbol w_{H,J_\star^c}\|_2$. By \eqref{eq:common-calibration-identities} from Lemma~\ref{lem:common-calibration}, $d_q\leq\frac{2q}{z_q^2}$, and hence $z_q\sqrt q\leq\frac{\sqrt2q}{\sqrt{d_q}}$. Together with $\|\boldsymbol w_{H,J_\star^c}\|_2\leq\|\boldsymbol w_H\|_2\leq r^+$, this yields $-\langle\boldsymbol h_{J_\star^c},\boldsymbol w_{H,J_\star^c}\rangle\leq\frac{8q}{\sqrt{d_q}}r^+\leq\frac{8q}{\sqrt{d_q}}r^+$, provided that $\gamma_H=0$.
Recall that
\begin{align}
\begin{aligned}
&N\lambda s_H=\|\boldsymbol\xi+r^+\boldsymbol g\|_2+\langle\boldsymbol h,\boldsymbol w_H\rangle.
\end{aligned}
\tag{\ref{eq:scalar-first-order}}
\end{align}
Together with $\|\boldsymbol w_{H,J_\star}\|_2\leq r^+$ and \eqref{eq:auxiliary-hypotheses}, this yields
\begin{align*}
\begin{aligned}
&N\lambda s_H\geq \|\bxi+r^+\vg\|_2 - \|\vh_{J_\star}\|_2\|\vw_{H,J_\star}\|_2 - \langle\vh_{J_\star^c}, \vw_{H,J_\star^c}\rangle \geq\|\boldsymbol\xi+r^+\boldsymbol g\|_2-r^+\bigg[\frac{8q}{\sqrt{d_q}}+\|\boldsymbol h_{J_\star}\|_2\bigg]>\frac{N\lambda}{\bar z_q},
\end{aligned}
\end{align*}
which contradicts with $\frac{1}{s_H}>\bar z_q$. Thus $\gamma_H>0$, and \eqref{eq:auxiliary-kkt} gives $$\|\boldsymbol w_H\|_2=r^+.$$

\textbf{Saddle Point Analysis.}
An important conclusion is that 
\begin{align}\label{eq:equation_at_saddle_point}
    \tilde\Psi(\vw_H,s_H) = \Psi(\vw_H,s_H) = \tilde\varphi_H(\vg,\vh)
\end{align}
Since $s_H\gamma_H>0$, then \eqref{eq:solution_S_Nlambda_s_H} may be rewritten as
\begin{align}\tag{\ref{eq:solution_S_Nlambda_s_H}}
&w_{H,j}=\frac{s_H}{\gamma_H}\bigg(S_{\frac{1}{s_H}}\big(\frac{\gamma_H\beta_j^\star}{s_H}-h_j\big)-\frac{\gamma_H\beta_j^\star}{s_H}\bigg),
\end{align}
and therefore, since $\|\boldsymbol w_H\|_2=r^+$, squaring both sides of \eqref{eq:solution_S_Nlambda_s_H}, summing over $1\leq j\leq p$, and taking the square root give
\begin{align}
\begin{aligned}
&\frac{\gamma_Hr^+}{s_H}=\bigg(\sum_{j=1}^p\bigg(S_{\frac{1}{s_H}}\bigg(\frac{\gamma_H\beta_j^\star}{s_H}-h_j\bigg)-\frac{\gamma_H\beta_j^\star}{s_H}\bigg)^2\bigg)^{\frac{1}{2}}.
\end{aligned}
\label{eq:soft-threshold-kkt}
\end{align}
\textbf{We first prove that $v_{\frac{\gamma_H}{s_H},\frac{1}{s_H}}\leq\bar z_q$.} Suppose to the contrary that $v_{\frac{\gamma_H}{s_H},\frac{1}{s_H}}>\bar z_q$. By the definition of $\Omega_{\mathrm{prof},2}^{\mathrm{aux}}$ with $\rho=\frac{\gamma_H}{s_H}$ and $u=\frac{1}{s_H}$ gives
\begin{align*}
&\bigg(-\bigg\langle\boldsymbol h_{J_\star^c},\boldsymbol\Delta_{\frac{\gamma_H}{s_H},\frac{1}{s_H}}(\boldsymbol h)\bigg\rangle\bigg)_+\leq\frac{8q}{\sqrt{d_q}}\bigg\|\boldsymbol\Delta_{\frac{\gamma_H}{s_H},\frac{1}{s_H}}(\boldsymbol h)\bigg\|_2.
\end{align*}
Together with \eqref{eq:scalar-first-order} and \eqref{eq:auxiliary-hypotheses}, the bound established above gives $N\lambda s_H\geq\|\boldsymbol\xi+r^+\boldsymbol g\|_2-r^+\bigg(\frac{8q}{\sqrt{d_q}}+\|\boldsymbol h_{J_\star}\|_2\bigg)>\frac{N\lambda}{\bar z_q}$. Therefore, $s_H>\frac{1}{\bar z_q}$. For every $\theta\in\mathbb R$ and $u\geq0$, $\mathbb P\big\{|G+\theta|>u\big\}\geq2\Phi_G(u)$. Hence, by the definition of $v_{\rho,u}$, $\Phi_G(v_{\rho,u})=\frac{1}{2|J_\star^c|}\sum_{j\in J_\star^c}\mathbb P\big\{|G-\rho\beta_j^\star|>u\big\}\geq\Phi_G(u)$. Since $\Phi_G$ is strictly decreasing, $v_{\rho,u}\leq u$. Taking $\rho=\frac{\gamma_H}{s_H}$ and $u=\frac{1}{s_H}$ therefore gives $v_{\frac{\gamma_H}{s_H},\frac{1}{s_H}}\leq\frac{1}{s_H}<\bar z_q$, contradicting the assumption $v_{\frac{\gamma_H}{s_H},\frac{1}{s_H}}>\bar z_q$. Therefore, $v_{\frac{\gamma_H}{s_H},\frac{1}{s_H}}\leq\bar z_q$.

Since we have verified that $v_{\frac{\gamma_H}{s_H},\frac{1}{s_H}}\leq\bar z_q$, the defining inequality of $\Omega_{\mathrm{prof},1}^{\mathrm{aux}}$ with $\rho=\frac{\gamma_H}{s_H}$ and $u=\frac{1}{s_H}$ gives $\|\boldsymbol\Delta_{\frac{\gamma_H}{s_H},\frac{1}{s_H}}(\boldsymbol h)\|_2^2\geq\frac{d_q}{256}$. Since all summands in \eqref{eq:soft-threshold-kkt} are nonnegative,
\begin{align}
\begin{aligned}
&\frac{\gamma_Hr^+}{s_H}=\bigg(\sum_{j=1}^p\bigg(S_{\frac{1}{s_H}}\bigg(\frac{\gamma_H\beta_j^\star}{s_H}-h_j\bigg)-\frac{\gamma_H\beta_j^\star}{s_H}\bigg)^2\bigg)^{\frac{1}{2}}\\
&\geq\bigg(\sum_{j\in J_\star^c}\bigg(S_{\frac{1}{s_H}}\bigg(\frac{\gamma_H\beta_j^\star}{s_H}-h_j\bigg)-\frac{\gamma_H\beta_j^\star}{s_H}\bigg)^2\bigg)^{\frac{1}{2}}
=\bigg\|\boldsymbol\Delta_{\frac{\gamma_H}{s_H},\frac{1}{s_H}}(\boldsymbol h)\bigg\|_2\geq2^{-4}\sqrt{d_q}.
\end{aligned}
\label{eq:multiplier-lower-bound}
\end{align} The first inequality in \eqref{eq:auxiliary-hypotheses} now yields that $\exists 0<s_H<H$ and $\exists \vw$, such that $\|\boldsymbol w_H\|_2=r^+$, such that
\begin{align}
\begin{aligned}
& \boldsymbol 0\in s_H\boldsymbol h+\partial\|\boldsymbol\beta^\star+\boldsymbol w_H\|_1+\gamma_H\boldsymbol w_H,\mbox{ and } \gamma_Hr^+-s_H\frac{\langle\boldsymbol\xi+r^+\boldsymbol g,\boldsymbol g\rangle}{\|\boldsymbol\xi+r^+\boldsymbol g\|_2}\geq s_H(2^{-4}-2^{-6})\sqrt{d_q}\geq2^{-5}s_H\sqrt{d_q}.
\end{aligned}
\label{eq:saddle-structure}
\end{align}
Therefore, from the KKT condition $\boldsymbol 0\in s_H\boldsymbol h+\partial\|\boldsymbol\beta^\star+\boldsymbol w_H\|_1+\gamma_H\boldsymbol w_H$, there holds 
\begin{align*}
    \forall\|\vw\|_2\leq r^+,\qquad \|\bbeta^\star+\vw\|_1-\|\bbeta^\star-\vw_H\|_1\geq -\langle s_H\vh + \gamma_H\vw_H,\vw-\vw_H\rangle.
\end{align*}
Moreover, from $\|\vw_H\|_2 = r^+$, there holds $\gamma_H(\|\vw_H\|_2^2 - \langle\vw,\vw_H\rangle)\geq\gamma_H((r^+)^2-r^+\|\vw\|_2)=\gamma_Hr^+(r^+-\|\vw\|_2)$; and consequently
\begin{align*}
    \forall\|\vw\|_2\leq r^+,\qquad s_H\langle\vh,\vw-\vw_H\rangle+\|\bbeta^\star+\vw\|_1-\|\bbeta^\star+\vw_H\|_1\geq \gamma_Hr^+(r^+-\|\vw\|_2).
\end{align*}
By \eqref{eq:equation_at_saddle_point}, $\Psi(\vw_H,s_H)=\tilde\varphi_H(\vg,\vh)$. hence $\Psi(\vw,s_H)-\tilde\varphi(\vg,\vh) = \Psi(\vw,s_H)-\Psi(\vw_H,s_H)$. Therefore, $\forall\,\|\boldsymbol w\|_2\leq r^+$, together with the definition of $\Psi(\vw,s_H)$ and that of $\Psi(\vw_H,s_H)$, there holds
\begin{align*}
&\Psi(\boldsymbol w,s_H)-\widetilde\varphi_H(\boldsymbol g,\boldsymbol h)=s_H\big(\|\boldsymbol\xi+\|\boldsymbol w\|_2\boldsymbol g\|_2-\|\boldsymbol\xi+r^+\boldsymbol g\|_2\big)+s_H\langle\boldsymbol h,\boldsymbol w-\boldsymbol w_H\rangle+\|\boldsymbol\beta^\star+\boldsymbol w\|_1-\|\boldsymbol\beta^\star+\boldsymbol w_H\|_1.
\end{align*}
The map $t\mapsto\|\boldsymbol\xi+t\boldsymbol g\|_2$ is convex and $\frac{\mathrm{d}}{\mathrm{d}t}\|\bxi+t\vg\|_2\big|_{t=r^+} = \frac{\langle\bxi+r^+\vg,\vg\rangle}{\|\bxi+r^+\vg\|_2}$, hence
\begin{align*}
&\forall\|\vw\|_2\leq r^+,\qquad\|\boldsymbol\xi+\|\boldsymbol w\|_2\boldsymbol g\|_2-\|\boldsymbol\xi+r^+\boldsymbol g\|_2\geq\frac{\langle\boldsymbol\xi+r^+\boldsymbol g,\boldsymbol g\rangle}{\|\boldsymbol\xi+r^+\boldsymbol g\|_2}\big(\|\boldsymbol w\|_2-r^+\big).
\end{align*}
Together with the preceding inequality, this gives
\begin{align*}
\begin{aligned}
&\forall\|\vw\|_2\leq r^+,\qquad\Psi(\boldsymbol w,s_H)-\widetilde\varphi_H(\boldsymbol g,\boldsymbol h)\geq\bigg(\gamma_Hr^+-s_H\frac{\langle\boldsymbol\xi+r^+\boldsymbol g,\boldsymbol g\rangle}{\|\boldsymbol\xi+r^+\boldsymbol g\|_2}\bigg)\big(r^+-\|\boldsymbol w\|_2\big).
\end{aligned}
\end{align*}
Moreover, from $r^-<r^+$,  and \eqref{eq:saddle-structure}, there holds
\begin{align}
\begin{aligned}
&\forall\|\vw\|_2\leq r^-,\qquad \Psi(\boldsymbol w,s_H)\geq\widetilde\varphi_H(\boldsymbol g,\boldsymbol h)+2^{-6}s_Hr^+\sqrt{d_q}.
\end{aligned}
\label{eq:pointwise-auxiliary-gap}
\end{align}
Maximizing over $s$ and then minimizing over $\|\boldsymbol w\|_2\leq r^-$ gives
\begin{align}
\begin{aligned}
&\phi_{r^-,H}\geq\widetilde\varphi_H(\boldsymbol g,\boldsymbol h)+2^{-6}s_Hr^+\sqrt{d_q}.
\end{aligned}
\label{eq:small-ball-auxiliary-gap}
\end{align}
\endproof

\paragraph{Proof of Lemma~\ref{lem:def_x_H}.}
By Section~\ref{sec:cgmt}, and \eqref{eq:variational_representation_LASSO}, it is easy to see that $\phi_{t,H}$ is the auxiliary problem of $\Phi_{t,H}$. Notice that
\begin{align}
\begin{aligned}
&\forall (\vw,s)\in B_2(\vzero;t)\times[0,H],\qquad |\Psi_{\boldsymbol g',\boldsymbol h'}(\boldsymbol w,s)-\Psi_{\boldsymbol g,\boldsymbol h}(\boldsymbol w,s)|\leq\sqrt2s\|\boldsymbol w\|_2\|(\boldsymbol g'-\boldsymbol g,\boldsymbol h'-\boldsymbol h)\|_2.
\end{aligned}
\label{eq:auxiliary-lipschitz}
\end{align}
Thus $(\boldsymbol g,\boldsymbol h)\mapsto\Psi_{\boldsymbol g,\boldsymbol h}(\boldsymbol w,s)$ is $\sqrt2s\|\boldsymbol w\|_2$-Lipschitz.

We first verify their separation. Fix $(\boldsymbol g',\boldsymbol h')\in B_x^{\mathrm{aux}}$ and $(\boldsymbol w_H',s_H')\in\mathcal S_H(\boldsymbol g',\boldsymbol h')$. From the definition of $B_x^{\mathrm{aux}}$, there holds $\tilde\varphi_H(\vg',\vh')\geq x$, and from Lemma~\ref{lem:saddle_point_analysis_tilde_Psi}, $s_H'>0$. Fix $(\boldsymbol g,\boldsymbol h)\in A_x^{\mathrm{aux}}$ and
\begin{align*}
&\boldsymbol w\in\mathop{\arg\min}_{\|\boldsymbol z\|_2\leq r^-}\max_{0\leq s\leq H}\Psi_{\boldsymbol g,\boldsymbol h}(\boldsymbol z,s).
\end{align*}Then $\Psi_{\vg,\vh}(\vw,s_H')\leq \max_{0\leq s\leq H}\Psi_{\vh,\vh}(\vw,s) = \phi_{r^-,H}(\vg,\vh) \leq x$, where we have used the assumption that $(\vg,\vh)\in A_x^{\mathrm{aux}}$.
Moreover, Lemma~\ref{lem:saddle_point_analysis_tilde_Psi} yields that $\Psi(\boldsymbol w,s_H')\geq\widetilde\varphi_H(\boldsymbol g',\boldsymbol h')+2^{-6}s_H'r^+\sqrt{d_q}$. Together with the Lipschitzness, that is, \eqref{eq:auxiliary-lipschitz},
\begin{align*}
&x+2^{-6}s_H'r^+\sqrt{d_q}\leq\Psi_{\boldsymbol g',\boldsymbol h'}(\boldsymbol w,s_H') \leq\Psi_{\boldsymbol g,\boldsymbol h}(\boldsymbol w,s_H')+\sqrt2s_H'\|\boldsymbol w\|_2\|(\boldsymbol g'-\boldsymbol g,\boldsymbol h'-\boldsymbol h)\|_2 \leq x+\sqrt2s_H'r^+\|(\boldsymbol g'-\boldsymbol g,\boldsymbol h'-\boldsymbol h)\|_2.
\end{align*}
Cancellation of $s_H'r^+>0$ yields $\|(\boldsymbol g'-\boldsymbol g,\boldsymbol h'-\boldsymbol h)\|_2\geq\frac{2^{-6}\sqrt{d_q}}{\sqrt2}\geq2^{-7}\sqrt{d_q}$. Taking the infimum over $A_x^{\mathrm{aux}}\times B_x^{\mathrm{aux}}$ proves
\begin{align}
\begin{aligned}
&\operatorname{dist}_{\mathbb R^{N+p}}(A_x^{\mathrm{aux}},B_x^{\mathrm{aux}})\geq2^{-7}\sqrt{d_q}.
\end{aligned}
\label{eq:gaussian-set-separation}
\end{align}

It remains to prove \eqref{eq:no-conditional-atom}. Fix $\boldsymbol h$ and $x\in\mathbb R$. Suppose that $\boldsymbol g_1,\boldsymbol g_2$ satisfy $(\boldsymbol g_i,\boldsymbol h)\in\Omega_{\boldsymbol g,\boldsymbol h}^{\mathrm{aux}}(\boldsymbol\xi)$ and $\widetilde\varphi_H(\boldsymbol g_i,\boldsymbol h)=x$ for $i\in(1,2)$, and set $a_i=\|\boldsymbol\xi+r^+\boldsymbol g_i\|_2$. Suppose that $a_1<a_2$, and take $(\boldsymbol w_{H,1},s_{H,1})\in S_H(\boldsymbol g_1,\boldsymbol h)$. Lemma~\ref{lem:saddle_point_analysis_tilde_Psi} gives $s_{H,1}>0$. Since the dependence of $\widetilde\Psi$ on $\boldsymbol g$ is only through $s\|\boldsymbol\xi+r^+\boldsymbol g\|_2$, for every $\|\boldsymbol w\|_2\leq r^+$ and $0\leq s\leq H$, $\widetilde\Psi_{\boldsymbol g_2,\boldsymbol h}(\boldsymbol w,s)=\widetilde\Psi_{\boldsymbol g_1,\boldsymbol h}(\boldsymbol w,s)+s(a_2-a_1)$. The saddle inequalities therefore give $$\widetilde\varphi_H(\boldsymbol g_2,\boldsymbol h)\geq\min_{\|\boldsymbol w\|_2\leq r^+}\widetilde\Psi_{\boldsymbol g_2,\boldsymbol h}(\boldsymbol w,s_{H,1})=\min_{\|\boldsymbol w\|_2\leq r^+}\widetilde\Psi_{\boldsymbol g_1,\boldsymbol h}(\boldsymbol w,s_{H,1})+s_{H,1}(a_2-a_1)=\widetilde\phi_H(\boldsymbol g_1,\boldsymbol h)+s_{H,1}(a_2-a_1)>x,$$ a contradiction. Thus, for fixed $\boldsymbol h$ and $x$, the set $\big(\|\boldsymbol\xi+r^+\boldsymbol g\|_2:(\boldsymbol g,\boldsymbol h)\in\Omega_{\boldsymbol g,\boldsymbol h}^{\mathrm{aux}}(\boldsymbol\xi),\ \widetilde\phi_H(\boldsymbol g,\boldsymbol h)=x\big)$ contains at most one number. Each set $\big(\boldsymbol g:\|\boldsymbol\xi+r^+\boldsymbol g\|_2=a\big)$ is a Euclidean sphere of $\gamma_N$-measure zero. Fubini's theorem therefore gives, for every $x\in\mathbb R$, $\mathbb P_{\mathrm{aux}}\big\{\widetilde\phi_H(\boldsymbol g,\boldsymbol h)=x,\Omega_{\boldsymbol g,\boldsymbol h}^{\mathrm{aux}}(\boldsymbol\xi)\mid\boldsymbol\xi\big\}=0$. By \eqref{eq:stochastic-event-bounds} and the strict assumption at the beginning of the proof, $\mathbb P_{\mathrm{aux}}\big(\Omega_{\boldsymbol g,\boldsymbol h}^{\mathrm{aux}}(\boldsymbol\xi)\mid\boldsymbol\xi\big)>\Phi_G\big(2^{-8}\sqrt{d_q}\big)$, whereas $\mathbb P_{\mathrm{aux}}\big(B_x^{\mathrm{aux}}\mid\boldsymbol\xi\big)\downarrow0$ as $x\uparrow\infty$. Together with the absence of conditional atoms proved above, this proves the assertion and the existence of $x_H(\boldsymbol\xi)$.

\paragraph{Proof of Lemma~\ref{lem:auxiliary_lower_bound}.}
On $\Omega_g^{\mathrm{aux}}$, $\|\boldsymbol\xi+r^+\boldsymbol g\|_2^2-\|\boldsymbol\xi\|_2^2=2r^+\langle\boldsymbol\xi,\boldsymbol g\rangle+(r^+)^2\|\boldsymbol g\|_2^2\geq0$, hence $\|\boldsymbol\xi\|_2\leq\|\boldsymbol\xi+r^+\boldsymbol g\|_2$, and the convexity of $t\mapsto\|\boldsymbol\xi+t\boldsymbol g\|_2$ gives, for every $0\leq t\leq r^+$, there holds $\|\boldsymbol\xi+t\boldsymbol g\|_2\leq\bigg(1-\frac{t}{r^+}\bigg)\|\boldsymbol\xi\|_2+\frac{t}{r^+}\|\boldsymbol\xi+r^+\boldsymbol g\|_2\leq\|\boldsymbol\xi+r^+\boldsymbol g\|_2.$ In particular, for every $\|\boldsymbol w\|_2\leq r^+$ and $0\leq s\leq H$, taking $t=\|\boldsymbol w\|_2$ gives $\Psi(\boldsymbol w,s)\leq\widetilde\Psi(\boldsymbol w,s)$. Recalling the definition of $\widetilde\varphi_H(\boldsymbol g,\boldsymbol h)$, we obtain
\begin{align}
\begin{aligned}
&\phi_{r^+,H}=\min_{\|\boldsymbol w\|_2\leq r^+}\max_{0\leq s\leq H}\Psi(\boldsymbol w,s)\leq\min_{\|\boldsymbol w\|_2\leq r^+}\max_{0\leq s\leq H}\widetilde\Psi(\boldsymbol w,s)=\widetilde\varphi_H(\boldsymbol g,\boldsymbol h).
\end{aligned}
\tag{\ref{eq:large-ball-auxiliary-bound}}
\end{align}

Without loss of generality, we assume $1-C_{\mathrm{prob}}\exp(-c_{\mathrm{prob}}d_q)>0.$ We use the estimates in \eqref{eq:common-calibration-identities} and \eqref{eq:common-radial-calibration} throughout. From $1-C_{\mathrm{prob}}\exp(-c_{\mathrm{prob}}d_q)>0$, \eqref{eq:common-calibration-identities} implies $q\geq120$. Hence, $\mathbb P_{\mathrm{aux}}(\Omega_{\boldsymbol g,\boldsymbol h}^{\mathrm{aux}}(\boldsymbol\xi)\mid\boldsymbol\xi)>\Phi_G(2^{-8}\sqrt{d_q})$ for every $\boldsymbol\xi\in\Omega_\xi$. Apply Lemma~\ref{lem:common-gaussian-separation} with $U=\phi_{r^-,H}$, $V=\widetilde\varphi_H$, $E=\Omega_{\boldsymbol g,\boldsymbol h}^{\mathrm{aux}}(\boldsymbol\xi)$, and $a=2^{-8}\sqrt{d_q}$. Equations \eqref{eq:gaussian-set-separation} and \eqref{eq:no-conditional-atom} verify its hypotheses, and consequently,
\begin{align*}
&\mathbb P_{\mathrm{aux}}\big\{\phi_{r^-,H}\leq x_H(\boldsymbol\xi)\mid\boldsymbol\xi\big\}+\mathbb P_{\mathrm{aux}}\big\{\widetilde\varphi_H(\boldsymbol g,\boldsymbol h)\geq x_H(\boldsymbol\xi)\mid\boldsymbol\xi\big\}\leq2\Phi_G(2^{-8}\sqrt{d_q})+\mathbb P_{\mathrm{aux}}\big(\Omega_{\boldsymbol g,\boldsymbol h}^{\mathrm{aux}}(\boldsymbol\xi)^c\mid\boldsymbol\xi\big).
\end{align*}
The estimate $\Phi_G(t)\leq \frac{1}{2}e^{-\frac{t^2}{2}}$ and \eqref{eq:stochastic-event-bounds} then give
\begin{align}
\begin{aligned}
&\mathbb P_{\mathrm{aux}}\big\{\phi_{r^-,H}\leq x_H(\boldsymbol\xi)\mid\boldsymbol\xi\big\}+\mathbb P_{\mathrm{aux}}\big\{\widetilde\varphi_H(\boldsymbol g,\boldsymbol h)\geq x_H(\boldsymbol\xi)\mid\boldsymbol\xi\big\}\leq2^{9}\exp\bigg(-2^{-35}d_q\bigg).
\end{aligned}
\tag{\ref{eq:conditional-auxiliary-tails}}
\end{align}

\subsection{End of the proof of Theorem~\ref{theo:lower_LASSO}}

\begin{proof}
Fix $\lambda>0$ satisfying Assumption~\ref{assumption:LASSO_lower}, and define
\begin{align*}
&\Omega_{\mathrm{bias}}=\left\{\|\widehat{\boldsymbol\beta}-\boldsymbol\beta^\star\|_2^2\geq\frac1{60}\left(\lambda^2|J_\star|+\|\boldsymbol\beta_{J_\star^c}^\star\|_2^2\right)\right\},\, \Omega_{\mathrm{var},J_\star}=\left\{\|\widehat{\boldsymbol\beta}_{J_\star}-\boldsymbol\beta_{J_\star}^\star\|_2^2\geq\frac{|J_\star|}{16}\min\left\{\frac{\sigma_\xi^2}{N},\lambda^2\right\}\right\},\\
&\mbox{ and }\Omega_{\mathrm{var},J_\star^c}=\left\{\|\widehat{\boldsymbol\beta}-\boldsymbol\beta^\star\|_2^2\geq2^{-30}\mathrm{Var}_{J_\star^c}\right\}.
\end{align*}
Propositions~\ref{prop:lower_bound_bias}, \ref{prop:active-sign-dropout}, and~\ref{prop:free-variance-lower-bound} give
\begin{align*}
&\mathbb P^{\otimes N}(\Omega_{\mathrm{bias}}^c)\leq4\exp\left(-\frac{|J_\star|}{16}\right),\qquad
\mathbb P^{\otimes N}(\Omega_{\mathrm{var},J_\star}^c)\leq\exp(-0.2|J_\star|),\mbox{ and }\mathbb P^{\otimes N}(\Omega_{\mathrm{var},J_\star^c}^c)\leq2^{11}\exp(-2^{-35}d_q).
\end{align*}
On $\Omega_{\mathrm{var},J_\star}$, $\|\widehat{\boldsymbol\beta}-\boldsymbol\beta^\star\|_2^2\geq\frac1{16}\frac{\sigma_\xi^2|J_\star|}{N}.$
If $d_q\geq|J_\star|$, then on $\Omega_{\mathrm{bias}}\cap\Omega_{\mathrm{var},J_\star}\cap\Omega_{\mathrm{var},J_\star^c}$,
\begin{align*}
\|\widehat{\boldsymbol\beta}-\boldsymbol\beta^\star\|_2^2
&\geq2^{-31}\left(\frac{\sigma_\xi^2|J_\star|}{N}+\lambda^2|J_\star|+\|\boldsymbol\beta_{J_\star^c}^\star\|_2^2+\mathrm{Var}_{J_\star^c}\right)=2^{-31}r^2(V_{J_\star},V_{J_\star^c}).
\end{align*}
If $d_q<|J_\star|$, Lemma~\ref{lem:common-calibration} gives $d_q/N<0.31$, and hence $\mathrm{Var}_{J_\star^c}=\frac{\sigma_\xi^2d_q}{N-d_q}\leq\frac{100}{69}\frac{\sigma_\xi^2|J_\star|}{N}.$ Consequently, on $\Omega_{\mathrm{bias}}\cap\Omega_{\mathrm{var},J_\star}$, $\|\widehat{\boldsymbol\beta}-\boldsymbol\beta^\star\|_2^2\geq2^{-31}r^2(V_{J_\star},V_{J_\star^c})$. 

\end{proof}

\section{Technical Background}

\subsection{Tools from Asymptotic Geometric Analysis}

\begin{Lemma}[\cite{carlGelfandNumbersOperators1988},\cite{gluskinExtremalPropertiesOrthogonal1989}]\label{lem:carl-pajor}
Let $\mathrm{Vol}_N(\cdot)$ be the $N$-dimensional Lebesgue measure volume. If \(\{\mathbb X_{\{j\}}/(4\sqrt N):j\in J_\star^c\}\) spans \(\mathbb R^N\) and \(\max_{j\in J_\star^c}\|\mathbb X_{\{j\}}\|_2\leq4\sqrt N\), then
\begin{align*}
&\operatorname{Vol}_N\bigg(\bigg\{\boldsymbol u\in\mathbb R^N:\big\|\mathbb X_{J_\star^c}^\top\boldsymbol u\big\|_\infty\leq4\sqrt N\bigg\}\bigg)^{1/N}\geq\frac{c_{\mathrm{CP}}}{\sqrt{\log\big(\frac{e|J_\star^c|}{N}\big)}}.
\end{align*}
\end{Lemma}

\begin{Lemma}[\cite{artstein-avidanAsymptoticGeometricAnalysis2015}, Theorem 1.5.11]\label{lem:urysohn}
For every convex body \(K\subset\mathbb R^N\), where \(\sigma\) is the normalized Haar probability measure on \(S^{N-1}\),
\begin{align*}
&\int_{S^{N-1}}\sup_{\boldsymbol u\in K}\langle\boldsymbol\theta,\boldsymbol u\rangle\,\mathrm d\sigma(\boldsymbol\theta)\geq\bigg(\frac{\operatorname{Vol}_N(K)}{\operatorname{Vol}_N(B_2^N)}\bigg)^{1/N}.
\end{align*}
\end{Lemma}

\begin{Lemma}
\label{lem:common-gaussian-isoperimetry}
Let \(\gamma_n\) be the standard Gaussian measure on \(\mathbb R^n\). If \(A,B\subset\mathbb R^n\) are Borel and \(\operatorname{dist}(A,B) = \inf(\|\va-\vb\|_2:\va\in A,\vb\in B)\geq t\), then \(\Phi_G^{-1}(\gamma_n(A))+\Phi_G^{-1}(\gamma_n(B))\geq t\).
\end{Lemma}
\begin{proof}
    For any $0<u<t$, from $\dist(A,B)\geq t$, there holds $B\subset A_u^c$, where $A_u=\{\vx:\exists\va\in A,\mbox{ s.t., }\|\vx-\va\|_2\leq u\}$. By the Sudakov--Tsirelson--Borell inequality, see, for instance, \cite[Theorem 3.1.1]{artstein-avidanAsymptoticGeometricAnalysis2021}, $\gamma_n(A_u)\geq \Phi_G\bigl(\Phi_G^{-1}(\gamma_n(A))-u\bigr)$. Therefore,
    \[
    \gamma_n(B)\leq 1-\gamma_n(A_u)\leq 1-\Phi_G\bigl(\Phi_G^{-1}(\gamma_n(A))-u\bigr)=\Phi_G\bigl(-\Phi_G^{-1}(\gamma_n(A))+u\bigr).
    \]
    Since $\Phi_G$ is strictly decreasing, applying $\Phi_G^{-1}$ to both sides yields
    \[
    \Phi_G^{-1}(\gamma_n(B))\geq-\Phi_G^{-1}(\gamma_n(A))+u.
    \]
    Hence $\Phi_G^{-1}(\gamma_n(A))+\Phi_G^{-1}(\gamma_n(B))\geq u$. Sending $u\uparrow t$ concludes the proof.
\end{proof}

\begin{Lemma}
\label{lem:common-gaussian-separation}
Let \(U,V:\mathbb R^n\to\mathbb R\) be Borel functions, let \(E\subset\mathbb R^n\) be a Borel set, and let \(a>0\). Assume
\begin{align*}
&\operatorname{dist}(\{U\leq x\},\{V\geq x\}\cap E)\geq2a,\quad \gamma_n\{\{V=x\}\cap E\}=0,\quad\forall x\in\mathbb R,\mbox{ and assume } \gamma_n(E)>\Phi_G(a).
\end{align*}
Let
\begin{align*}
&x_E=\inf\{x:\gamma_n\{\{V\geq x\}\cap E\}\leq\Phi_G(a)\}
\end{align*}
Then
\begin{align*}
&\gamma_n\{U\leq x_E\}+\gamma_n\{V\geq x_E\}\leq2\Phi_G(a)+\gamma_n(E^c).
\end{align*}
\end{Lemma}

\begin{proof}
Notice that $x\mapsto \gamma_n\{\{V\geq x\}\cap E\}$ is decreasing. From $\gamma_n\{\{V=x\}\cap E\}=0$ for any $x\in\bR$ and the continuity of $\gamma_n$, the function $x\mapsto \gamma_n\{\{V\geq x\}\cap E\}$ is continuous. Since $\lim_{x\to-\infty}\gamma_n\{\{V\geq x\}\cap E\} = \gamma_n(E)>\Phi_G(a)$, and $\lim_{x\to+\infty}\gamma_n\{\{V\geq x\}\cap E\} = 0<\Phi_G(a)$, there must hold $\gamma_n\{\{V\geq x_E\}\cap E\}=\Phi_G(a)$.

Applying Lemma~\ref{lem:common-gaussian-isoperimetry} with $A = \{U\leq x_E\}$, $B=\{V\geq x_E\}\cap E$, and $t=2a$, we obtain $\Phi_G^{-1}(\gamma_n\{U\leq x_E\}) + \Phi_G^{-1}(\Phi_G(a)) = \Phi_G^{-1}(\gamma_n\{U\leq x_E\}) + \Phi_G^{-1}(\gamma_n\{\{V\geq x_E\}\cap E\})\geq 2a$. Hence $\Phi_G^{-1}(\gamma_n\{U\leq x_E\})\geq a$. Applying $\Phi_G$ on both sides together with the fact that $\Phi_G$ is decreasing, we obtain that $\gamma_n\{U\leq x_E\}\leq \Phi_G(a)$. Now since $\{V\geq x_E\}\subset (\{V\geq x_E\}\cap E)\cup E^c$, $\gamma_n\{V\geq x_E\}\leq \Phi_G(a)+\gamma_n(E^c)$.
\end{proof}

\begin{Lemma}\label{lem:radial-restriction-gap}
Let $K$ be a compact convex subset of a Euclidean space, let $A$ be a compact interval, and let $F:K\times A\to\mathbb R$ be convex--concave. Suppose $(\boldsymbol w^\circ,a^\circ)$ is a saddle point, $\|\boldsymbol w^\circ\|_2\leq r<R$, and $\boldsymbol w\mapsto F(\boldsymbol w,a^\circ)-\kappa\|\boldsymbol w\|_2^2/2$ is convex on every segment joining $\boldsymbol w^\circ$ to $K\cap RB_2$. Then
\begin{align*}
&\min_{\substack{\boldsymbol w\in K\\\|\boldsymbol w\|_2\geq R}}\max_{a\in A}F(\boldsymbol w,a)\geq\min_{\boldsymbol w\in K}\max_{a\in A}F(\boldsymbol w,a)+\frac\kappa2(R-r)^2.
\end{align*}
\end{Lemma}

\begin{proof}
Fix $\boldsymbol w\in K$ with $\|\boldsymbol w\|_2\geq R$. There exists $0<\alpha\leq1$ such that $\overline{\boldsymbol w}=(1-\alpha)\boldsymbol w^\circ+\alpha\boldsymbol w$ satisfies $\|\overline{\boldsymbol w}\|_2=R$. Since $\boldsymbol w^\circ$ minimizes $F(\cdot,a^\circ)$ on $K$, convexity gives
\begin{align*}
F(\overline{\boldsymbol w},a^\circ)\leq(1-\alpha)F(\boldsymbol w^\circ,a^\circ)+\alpha F(\boldsymbol w,a^\circ)\leq F(\boldsymbol w,a^\circ).
\end{align*}
The assumed strong convexity on the segment from $\boldsymbol w^\circ$ to $\overline{\boldsymbol w}$ and minimality at $\boldsymbol w^\circ$ imply $F(\overline{\boldsymbol w},a^\circ)\geq F(\boldsymbol w^\circ,a^\circ)+\kappa\|\overline{\boldsymbol w}-\boldsymbol w^\circ\|_2^2/2$. For example, apply strong convexity at $(1-t)\boldsymbol w^\circ+t\overline{\boldsymbol w}$, use minimality, divide by $t>0$, and let $t\downarrow0$. Thus
\begin{align*}
\max_{a\in A}F(\boldsymbol w,a)\geq F(\boldsymbol w^\circ,a^\circ)+\frac\kappa2(R-r)^2.
\end{align*}
Taking the minimum over the restricted set and using the saddle-point identity proves the result.
\end{proof}

\begin{Lemma}\label{lem:auxiliary-separation}
Let $\phi$ and $\phi_{S^c}$ be two functions on a Gaussian Euclidean space and let $E$ be a nonempty event, and let $L,\eta>0$. Suppose $\phi$ is $L$-Lipschitz on $E$ and $\phi_{S^c}-\phi\geq6\eta$ on $E$. Let
\begin{align*}
&\widetilde\phi(\boldsymbol\omega)=\inf_{\boldsymbol\omega'\in E}\{\phi(\boldsymbol\omega')+L\|\boldsymbol\omega-\boldsymbol\omega'\|_2\},
\end{align*}
and let $m$ be a median of $\widetilde\phi$ under the standard Gaussian probability measure. Then
\begin{align*}
&\mathbb P\{\phi\geq m+\eta\}\leq\mathbb P(E^c)+\exp\bigg(-\frac{\eta^2}{2L^2}\bigg),\qquad \mathbb P\{\phi_{S^c}\leq m+2\eta\}\leq\mathbb P(E^c)+\exp\bigg(-\frac{\eta^2}{2L^2}\bigg).
\end{align*}
\end{Lemma}

\begin{proof}
We first verify that $\tilde\phi$ is $L$-Lipschitz. In fact, such a $\tilde\phi$ is called the \emph{Kirszbraun-McShane extension} of $\phi$, as in \cite[Exercise 3.6.14, pp. 227]{gineMathematicalFoundationsInfiniteDimensional2015}. Take arbitrary $\vomega_1, \vomega_2$ from the ambient space, and arbitrary $\vomega' \in E$. We first prove that $\tilde\phi(\vomega_1) - \tilde\phi(\vomega_2) \le L\|\vomega_1 - \vomega_2\|_2$. By the triangle inequality, we have $\phi(\vomega') + L\|\vomega_1 - \vomega'\|_2 \le \phi(\vomega') + L\|\vomega_2 - \vomega'\|_2 + L\|\vomega_1 - \vomega_2\|_2$. Taking the infimum over $\vomega' \in E$ and combining with the definition of $\tilde\phi$, we obtain $\tilde\phi(\vomega_1) - \tilde\phi(\vomega_2) \le L\|\vomega_1 - \vomega_2\|_2$. Interchanging $\vomega_1$ and $\vomega_2$ proves the $L$-Lipschitzness. 

Note that on $E$, $\tilde\phi = \phi$, so we immediately have $\{\phi \ge m + \eta\} \subset E^c \cup \{\tilde\phi \ge m + \eta\}$. Applying the Gaussian Lipschitz concentration inequality, see, for instance, \cite[Section 3.1]{artstein-avidanAsymptoticGeometricAnalysis2015}, to the $L$-Lipschitz function $\tilde\phi$, we obtain $\mathbb{P}\{\phi \ge m + \eta\} \le \mathbb{P}(E^c) + \exp(-\frac{1}{2L^2}\eta^2)$. Moreover, since we have assumed that $\phi_{S^c} \ge \phi + 6\eta$ on $E$, and again using the fact that $\tilde\phi = \phi$ on $E$, on the set $E \cap \{\tilde\phi > m - \eta\}$ we must have $\phi_{S^c} \ge \phi + 6\eta = \tilde\phi + 6\eta > m + 2\eta$. Therefore its contrapositive states that $\{\phi_{S^c} \le m + 2\eta\} \subset E^c \cup \{\tilde\phi \le m - \eta\}$. Taking Gaussian measure on both sides again and applying Gaussian Lipschitz concentration, Lemma~\ref{lem:auxiliary-separation} follows.
\end{proof}

\begin{Lemma}\label{lem:restricted_Gaussian_Lipschitz}
Let $m,n\geq1$, let $\gamma_m$ be the standard Gaussian measure on $\mathbb R^m$, and let $K\subset\mathbb R^m$ be a closed convex set such that $\gamma_m(K)>0$. Let $\boldsymbol Z$ have the normalized restriction of $\gamma_m$ to $K$, that is, for every Borel set $A\subset\mathbb R^m$, $\mathbb P\{\boldsymbol Z\in A\}=\gamma_m(A\cap K)/\gamma_m(K)$, and let $\boldsymbol G\sim\mathcal N(0,I_n)$ be independent of $\boldsymbol Z$. If a locally Lipschitz function $f:K\times\mathbb R^n\to\mathbb R$ satisfies $\|\nabla_{\boldsymbol z}f\|_2^2+\|\nabla_{\boldsymbol g}f\|_2^2\leq L^2$ almost everywhere, then there exists an absolute constant $C>0$ such that
\[
\left\|f(\boldsymbol Z,\boldsymbol G)-\mathbb Ef(\boldsymbol Z,\boldsymbol G)\right\|_{\psi_2}\leq CL.
\]
\end{Lemma}

\begin{proof}
Set $C_0=K\times\mathbb R^n\subset\mathbb R^{m+n}$. By independence, the distribution of $(\boldsymbol Z,\boldsymbol G)$ is precisely the normalized restriction of the standard Gaussian measure $\gamma_{m+n}$ to the closed convex set $C_0$.

We first recall the Gaussian isoperimetric inequality of Bakry--Ledoux~\cite[Corollary 2.2 and Theorem 3.2]{bakryLevyGromovsIsoperimetricInequality1996}, see also \cite[Theorem 5.2.15]{vershyninHighDimensionalProbabilityIntroduction2018}: if a probability measure $\mu$ on $\mathbb R^d$ has a smooth density $Z^{-1}e^{-V}$ and $\nabla^2V\succeq I_d$, then there exists an absolute constant $C>0$ such that every $L$-Lipschitz function $h$ satisfies $\|h-\int h\,\mathrm d\mu\|_{\psi_2}\leq CL$.

Take a nonnegative $C^\infty$ function $\rho$ satisfying $\int_{\mathbb R^{m+n}}\rho=1$ and $\operatorname{supp}\rho\subset B_2^{m+n}$, and set $\rho_k(x)=k^{2(m+n)}\rho(k^2x)$. Let $d_{C_0}$ denote the Euclidean distance to $C_0$, and set $\varphi_k=k(d_{C_0}^2*\rho_k)$. Since $C_0$ is closed and convex, $d_{C_0}^2$ is convex, and hence $\varphi_k$ is smooth and convex. Define the probability measure
\[
\mathrm d\mu_k(x)=\frac{\exp(-\|x\|_2^2/2-\varphi_k(x))}{\int_{\mathbb R^{m+n}}\exp(-\|y\|_2^2/2-\varphi_k(y))\,\mathrm dy}\,\mathrm dx.
\]
The Hessian of the corresponding potential is bounded from below by $I_{m+n}$, so the Bakry--Ledoux conclusion applies to $\mu_k$.

If \(x\in C_0\), then \(d_{C_0}(x-y)\leq\|y\|_2\). Using \(\rho_k(y)=k^{2(m+n)}\rho(k^2y)\) and the change of variables \(z=k^2y\), for which \(y=k^{-2}z\) and \(\mathrm dy=k^{-2(m+n)}\,\mathrm dz\), we obtain
\[
\begin{aligned}
\varphi_k(x)&=k\int d_{C_0}(x-y)^2\rho_k(y)\,\mathrm dy\leq k\int\|y\|_2^2k^{2(m+n)}\rho(k^2y)\,\mathrm dy\\
&=k\int\|k^{-2}z\|_2^2k^{2(m+n)}\rho(z)k^{-2(m+n)}\,\mathrm dz=k^{-3}\int\|z\|_2^2\rho(z)\,\mathrm dz\longrightarrow0.
\end{aligned}
\]
If \(x\notin C_0\), then \(d_{C_0}(x)>0\). Since \(\operatorname{supp}\rho\subseteq B_2^{m+n}\), the definition of \(\rho_k\) gives \(\operatorname{supp}\rho_k\subseteq k^{-2}B_2^{m+n}\). Choose \(k\) sufficiently large that \(k^{-2}\leq d_{C_0}(x)/2\). Then every \(y\in\operatorname{supp}\rho_k\) satisfies \(\|y\|_2\leq k^{-2}\), and the \(1\)-Lipschitz property of \(d_{C_0}\) gives
\[
\begin{aligned}
d_{C_0}(x-y)&\geq d_{C_0}(x)-\|y\|_2\geq d_{C_0}(x)-k^{-2}\geq\frac{d_{C_0}(x)}{2},\\
\varphi_k(x)&=k\int d_{C_0}(x-y)^2\rho_k(y)\,\mathrm dy\geq\frac{kd_{C_0}(x)^2}{4}\int\rho_k(y)\,\mathrm dy=\frac{kd_{C_0}(x)^2}{4}\longrightarrow\infty.
\end{aligned}
\]
Thus \(e^{-\varphi_k}\to\mathbf1_{C_0}\) pointwise and \(0\leq e^{-\varphi_k}\leq1\). Dominated convergence yields \(\int e^{-\varphi_k}\,\mathrm d\gamma_{m+n}\to\gamma_{m+n}(C_0)\). Consequently,
\[
\frac{\mathrm d\mu_k}{\mathrm d\gamma_{m+n}}=\frac{e^{-\varphi_k}}{\int e^{-\varphi_k}\,\mathrm d\gamma_{m+n}}\longrightarrow\frac{\mathbf1_{C_0}}{\gamma_{m+n}(C_0)}\quad\text{pointwise},\qquad\int\left|\frac{\mathrm d\mu_k}{\mathrm d\gamma_{m+n}}-\frac{\mathbf1_{C_0}}{\gamma_{m+n}(C_0)}\right|\,\mathrm d\gamma_{m+n}\longrightarrow0,
\]
where the second conclusion follows from Scheff\'e's lemma. Hence \(\mu_k\) converges in total variation, and therefore weakly, to \(\gamma_{m+n}^{C_0}\).

By the assumption and convexity of \(C_0\), \(f\) is \(L\)-Lipschitz on \(C_0\). Take a McShane extension \(\widetilde f:\mathbb R^{m+n}\to\mathbb R\) defined in Lemma~\ref{lem:auxiliary-separation} preserving the Lipschitz constant. Since \(|\widetilde f(x)|\leq|\widetilde f(0)|+L\|x\|_2\), while \(\int e^{-\varphi_k}\,\mathrm d\gamma_{m+n}\geq\gamma_{m+n}(C_0)/2\) for all sufficiently large \(k\), we have \(\bigg|\frac{\widetilde f(x)e^{-\varphi_k(x)}}{\int e^{-\varphi_k}\,\mathrm d\gamma_{m+n}}\bigg|\leq2\left(|\widetilde f(0)|+L\|x\|_2\right)\gamma_{m+n}(C_0)^{-1}\). The right-hand side is \(\gamma_{m+n}\)-integrable, so dominated convergence gives
\[
\int\widetilde f\,\mathrm d\mu_k=\int\widetilde f(x)\frac{e^{-\varphi_k(x)}}{\int e^{-\varphi_k}\,\mathrm d\gamma_{m+n}}\,\mathrm d\gamma_{m+n}(x)\longrightarrow\frac{\int_{C_0}\widetilde f\,\mathrm d\gamma_{m+n}}{\gamma_{m+n}(C_0)}=\int_{C_0}f\,\mathrm d\gamma_{m+n}^{C_0}.
\]
The Bakry--Ledoux inequality gives \(\|\widetilde f-\int\widetilde f\,\mathrm d\mu_k\|_{\psi_2(\mu_k)}\leq CL\). Hence, for every \(a>CL\),
\[
\int\exp\left(\frac{\left|\widetilde f(x)-\int\widetilde f\,\mathrm d\mu_k\right|^2}{a^2}\right)\frac{e^{-\varphi_k(x)}}{\int e^{-\varphi_k}\,\mathrm d\gamma_{m+n}}\,\mathrm d\gamma_{m+n}(x)\leq2.
\]
The convergence of the centering constants and of the densities shows that the integrand converges pointwise to \(\exp\left(\frac{1}{a^2}\left|\widetilde f(x)-\int_{C_0}f\,\mathrm d\gamma_{m+n}^{C_0}\right|^2\right)\mathbf1_{C_0}(x)\gamma_{m+n}(C_0)^{-1}\). Fatou's lemma and \(\widetilde f=f\) on \(C_0\) therefore give
\[
\int_{C_0}\exp\left(\frac{\left|f(x)-\int_{C_0}f\,\mathrm d\gamma_{m+n}^{C_0}\right|^2}{a^2}\right)\mathrm d\gamma_{m+n}^{C_0}(x)\leq\liminf_{k\to\infty}\int\exp\left(\frac{\left|\widetilde f(x)-\int\widetilde f\,\mathrm d\mu_k\right|^2}{a^2}\right)\mathrm d\mu_k(x)\leq2.
\]
Since this holds for every \(a>CL\), \(\|f-\int_{C_0}f\,\mathrm d\gamma_{m+n}^{C_0}\|_{\psi_2(\gamma_{m+n}^{C_0})}\leq CL\). Since \((\boldsymbol Z,\boldsymbol G)\sim\gamma_{m+n}^{C_0}\), the conclusion follows.
\end{proof}

The proof is based on the treatment of normalized restricted measures by Bobkov--Nayar--Tetali, \cite{Bobkov2017}; in particular, they point out that normalized Gaussian restrictions to convex sets retain the sub-Gaussian constant and trace this fact back to Bakry--Ledoux in \cite{bakryLevyGromovsIsoperimetricInequality1996}.

\subsection{Convex Gaussian Min-Max Theorem}\label{sec:cgmt}

The Gaussian Convex Min-max Theorem (CGMT) is a tool for handling tail probabilities of the optimal value of random Gaussian convex-concave optimization problems involving a Gaussian bilinear form. Its power lies in its ability to handle tail probabilities of arbitrary scales, and in converting the problem into an auxiliary problem, through estimates of the saddle point of the auxiliary problem (via dimensionality reduction, see, for instance, Section~\ref{sec:dimension_reduction_trick}; or other specific optimization techniques, see, for instance, Section~\ref{sec:tools_convex_analysis}), to obtain estimates of the saddle point of the primary problem.

\begin{Lemma}[\cite{thrampoulidisPreciseErrorAnalysis2016}, Theorem 6.1, pp.~25--26]
\label{lem:common-cgmt}
Let \(S_{\boldsymbol w}\subset\mathbb R^m\) and \(S_{\boldsymbol u}\subset\mathbb R^n\) be compact, let \(\psi:S_{\boldsymbol w}\times S_{\boldsymbol u}\to\mathbb R\) be continuous, let $\vG\in\bR^{n\times m},\vg\in\bR^n,\vh\in\bR^m$ be standard Gaussian random matrix and vectors, and define
\begin{align*}
&\Phi(\boldsymbol G)=\min_{\boldsymbol w\in S_{\boldsymbol w}}\max_{\boldsymbol u\in S_{\boldsymbol u}}\{\langle\boldsymbol u,\boldsymbol G\boldsymbol w\rangle+\psi(\boldsymbol w,\boldsymbol u)\},\\
&\phi(\boldsymbol g,\boldsymbol h)=\min_{\boldsymbol w\in S_{\boldsymbol w}}\max_{\boldsymbol u\in S_{\boldsymbol u}}\{\|\boldsymbol w\|_2\langle\boldsymbol g,\boldsymbol u\rangle+\|\boldsymbol u\|_2\langle\boldsymbol h,\boldsymbol w\rangle+\psi(\boldsymbol w,\boldsymbol u)\}.
\end{align*}
For every \(x\in\mathbb R\), the following lower tail comparison holds $\mathbb P\{\Phi(\boldsymbol G)<x\}\leq2\mathbb P\{\phi(\boldsymbol g,\boldsymbol h)\leq x\}.$ If \(S_{\boldsymbol w}\) and \(S_{\boldsymbol u}\) are convex and \(\psi\) is convex in \(\boldsymbol w\) and concave in \(\boldsymbol u\), then the following upper tail comparison holds
\begin{align*}
&\mathbb P\{\Phi(\boldsymbol G)>x\}\leq2\mathbb P\{\phi(\boldsymbol g,\boldsymbol h)\geq x\}.
\end{align*}
Moreover, let \(S\) be relatively open in \(S_{\boldsymbol w}\). If deterministic \(\overline\phi,\overline\phi_{S^c}\in\mathbb R\), $0<\mathsf p<1$ and \(\eta>0\) satisfy
\begin{align*}
&\overline\phi_{S^c}\geq\overline\phi+3\eta,\qquad \mathbb P\{\phi(\vg,\vh)<\overline\phi+\eta\}\geq1-\mathsf p,\qquad \mathbb P\{\phi_{S^c}(\vg,\vh)>\overline\phi_{S^c}-\eta\}\geq1-\mathsf p,
\end{align*}
then every prescribed primary minimizer belongs to \(S\) with probability at least \(1-4\mathsf p\).
\end{Lemma}
Since CGMT involves two probability measures and two probability spaces, in the proofs of this paper we append the superscript $\mathrm{aux}$ to random events to indicate that the event is in the auxiliary probability space, and use $\mathbb{P}_{\mathrm{aux}}$ to denote the auxiliary probability measure, i.e., the probability measure of the auxiliary random variables such as $\vh, \vg$, etc.

\subsubsection{Variational representation of LASSO}\label{sec:variational_representation_LASSO}
We shall repeatedly use the following variational representation of $\|\cdot\|_2^2$.
\begin{align}\label{eq:variational_representation_ell_2_square}
&\forall a>0,\qquad \frac{\|\boldsymbol z\|_2^2}{2a}=\max_{\boldsymbol u\in\mathbb R^n}\bigg\{\langle\boldsymbol u,\boldsymbol z\rangle-\frac a2\|\boldsymbol u\|_2^2\bigg\},\mbox{ where the maximum is achieved at } \boldsymbol u=\frac{\boldsymbol z}{a}.
\end{align}
This variational representation of $\|\cdot\|_2^2$ is a common trick in the CGMT method, see, for instance, \cite[Section 3.2.2]{thrampoulidisRegularizedLinearRegression2015}, whose role is to convert the RERM of the minimization of the (regularized) squared risk into a min-max structure. For instance, for LASSO, with $a=N\lambda$ and $\vz=\vy-\bX\bbeta$, the min-max version of LASSO is given by
\begin{align}\label{eq:variational_representation_LASSO}
    \hat\bbeta \in \argmin_{\bbeta\in\bR^p}\max_{\vu\in\bR^N}\bigg( \big\langle \vu,\bxi + \bX(\bbeta^\star-\bbeta)\rangle - \frac{N\lambda}{2}\|\vu\|_2^2 + \|\bbeta\|_1 \bigg).
\end{align}We call $\vu$ the Lagrangian.
\subsubsection{Dimension reduction trick}\label{sec:dimension_reduction_trick}

In the CGMT method, another common trick is the one-dimensional reduction of the Lagrangian $\vu$. That is, fixing $\bbeta \in \mathbb{R}^p$, consider the inner maximization problem of the concave function $\vu \mapsto \langle \vu, \bX(\bbeta^\star - \bbeta) \rangle + \langle \vu, \bxi \rangle - \frac{N\lambda}{2}\|\vu\|_2^2 + \|\bbeta\|_1$ in \eqref{eq:variational_representation_LASSO}. One observation is that this function depends on the direction of $\vu$ only through $\langle \vu, \bX(\bbeta^\star - \bbeta) \rangle + \langle \vu, \bxi \rangle$. Therefore, if we set $r = \|\vu\|_2$, then the maximum of this inner concave function is exactly equal to
\begin{align*}
    \max\bigg( r \|\bxi + \bX(\bbeta^\star - \bbeta)\|_2 - \frac{N\lambda}{2}r^2 + \|\bbeta\|_1 : r\geq 0 \bigg).
\end{align*}Therefore, the optimization problem for $\hat\bbeta$ becomes
\begin{align*}
    \hat\bbeta\in\argmin_{\bbeta\in\bR^p}\max_{r\geq 0}\bigg( r\norm{\bxi+\bX(\bbeta^\star-\bbeta)}_2 - \frac{N\lambda}{2}r^2 + \norm{\bbeta}_1 \bigg).
\end{align*}
In particular, if $(\bbeta^\circ, r^\circ)$ is a saddle point of this convex-concave optimization problem, then by \eqref{eq:variational_representation_ell_2_square}, $r^\circ = \|\vu^\circ\|_2 = \frac{1}{N\lambda}\|\vy - \bX\bbeta^\circ\|_2$. The advantage of such a transformation is that the feasible set of the inner maximization problem is one-dimensional, and hence easier to examine. This technique has been widely adopted in the literature on the CGMT method; see, e.g., \cite{wangTightBoundsMinimum2022}.

\subsubsection{Compactification argument}

The CGMT requires both optimization domains to be compact. We therefore use the following compactification argument whenever one of the optimization domains is non-compact. Let $\mathcal S_{\boldsymbol w}\subset\mathbb R^m$ and $\mathcal S_{\boldsymbol u}\subset\mathbb R^n$ be closed (but not necessarily bounded) convex sets, and let $\psi:\mathcal S_{\boldsymbol w}\times\mathcal S_{\boldsymbol u}\to\mathbb R$ be continuous, convex in $\boldsymbol w$, and concave in $\boldsymbol u$. For $K,H>0$, define
\[
\Phi_{K,H}(\boldsymbol G)=\min_{\boldsymbol w\in\mathcal S_{\boldsymbol w}\cap KB_2^m}\max_{\boldsymbol u\in\mathcal S_{\boldsymbol u}\cap HB_2^n}\bigl\{\langle\boldsymbol u,\boldsymbol G\boldsymbol w\rangle+\psi(\boldsymbol w,\boldsymbol u)\bigr\},
\]
and
\[
\phi_{K,H}(\boldsymbol g,\boldsymbol h)=\min_{\boldsymbol w\in\mathcal S_{\boldsymbol w}\cap KB_2^m}\max_{\boldsymbol u\in\mathcal S_{\boldsymbol u}\cap HB_2^n}\bigl\{\|\boldsymbol w\|_2\langle\boldsymbol g,\boldsymbol u\rangle+\|\boldsymbol u\|_2\langle\boldsymbol h,\boldsymbol w\rangle+\psi(\boldsymbol w,\boldsymbol u)\bigr\}.
\]
For every finite $K,H$, the Theorem~\ref{lem:common-cgmt} applies to $\Phi_{K,H}$ and $\phi_{K,H}$.

We also define
\begin{align*}
    \Phi_{K,\infty}(\boldsymbol G)&=\min_{\boldsymbol w\in\mathcal S_{\boldsymbol w}\cap KB_2^m}\max_{\boldsymbol u\in\mathcal S_{\boldsymbol u}}\{\langle\boldsymbol u,\boldsymbol G\boldsymbol w\rangle+\psi(\boldsymbol w,\boldsymbol u)\},\mbox{ and } \Phi_{\infty,H}(\boldsymbol G)=\min_{\boldsymbol w\in\mathcal S_{\boldsymbol w}}\max_{\boldsymbol u\in\mathcal S_{\boldsymbol u}\cap HB_2^n}\{\langle\boldsymbol u,\boldsymbol G\boldsymbol w\rangle+\psi(\boldsymbol w,\boldsymbol u)\}.
\end{align*}
To remove the compactification, let $(K_\ell)_{\ell\geq1}$ and $(H_\ell)_{\ell\geq1}$ be deterministic sequences such that $K_\ell,H_\ell\uparrow\infty$. We try to verify that $\forall\omega$, there exist $K_0(\omega),H_0(\omega)<\infty$ such that $\Phi_{K,H}(\omega)=\Phi(\omega)$ for all $K\geq K_0(\omega)$ and $H\geq H_0(\omega)$. Hence there exists $\ell_0(\omega)\in\mathbb N$ such that $\Phi_{K_\ell,H_\ell}(\omega)=\Phi(\omega)$ for all $\ell\geq\ell_0(\omega)$. In particular, for any Borel set $A\subset\mathbb R$, if $E=\{\Phi\in A\}$ and $E_\ell=\{\Phi_{K_\ell,H_\ell}\in A\}$, then $\forall\omega\in E$, $\omega\in\bigcap_{\ell\geq\ell_0(\omega)}E_\ell$. Therefore $E\subset\bigcup_{k\geq1}\bigcap_{\ell\geq k}E_\ell=\liminf_{\ell\to\infty}E_\ell$, and, by Fatou's lemma, $\mathbb P(E)\leq\mathbb P(\liminf_{\ell\to\infty}E_\ell)\leq\liminf_{\ell\to\infty}\mathbb P(E_\ell)$. We refer to this passage as the compactification argument.

\subsection{The technical route of using CGMT to study the estimation error}

Let $\hat\bbeta\in\bR^p$ be a random vector that can written as the solution of a convex-concave optimization problem
\begin{align*}
    \hat\bbeta\in\argmin_{\bbeta\in\bR^p}\max_{\vu\in\bR^N}L(\bbeta,\vu),
\end{align*}where $L:\bR^p\times\bR^N\to\bR$ satisfies the assumption of CGMT (up to compactification). Suppose that our goal is to prove that $\widehat{\boldsymbol\beta}\in S$ with high $\mathbb P^{\otimes N}$-probability, where $S$ is a Borel subset of $\bR^p$. Denote
\begin{align*}
    \Phi = \min_{\bbeta\in\bR^p}\max_{\vu\in\bR^N} L(\bbeta,\vu),\mbox{ and } \Phi_{S^c} = \min_{\bbeta\in\bR^p\backslash S}\max_{\vu\in\bR^N} L(\bbeta,\vu)
\end{align*}
by the optimal value of this primary problem and its restriction onto $S^c$. Since restricting the feasible set can only increase the minimum, $\Phi_{S^c}\geq\Phi$. Moreover, on $\{\widehat{\boldsymbol\beta}\in S^c\}$, the minimizer $\widehat{\boldsymbol\beta}$ remains feasible for the restricted problem, and hence $\Phi_{S^c}=\Phi$. Therefore, $\{\widehat{\boldsymbol\beta}\in S^c\}\subset\{\Phi_{S^c}=\Phi\}.$ Thus, we need to compare the optimal values of two convex--concave optimization problems that differ only in their primal feasible sets, while their objective functions are identical and determined by the same realization of random variables. To this end, we seek a fixed real number $\bar\Phi$ such that, with high $\mathbb P^{\otimes N}$-probability, $\Phi<\bar\Phi<\Phi_{S^c}.$ The CGMT transforms these two primary problems into the corresponding auxiliary problems, whose optimal values are denoted by $\phi$ and $\phi_{S^c}$, respectively. Namely,
\begin{align*}
\phi=\inf_{\boldsymbol\beta\in K}\sup_{\boldsymbol u\in U}F(\boldsymbol\beta,\boldsymbol u),\mbox{ and }\phi_{S^c}=\inf_{\boldsymbol\beta\in K\cap S^c}\sup_{\boldsymbol u\in U}F(\boldsymbol\beta,\boldsymbol u),
\end{align*}where $F:\bR^p\times\bR^N\to\bR$ is the objective function of the auxiliary problem.
It therefore suffices to construct a high-probability event $\Omega_{\mathrm{aux}}$ in the auxiliary probability space and, through the analysis of the auxiliary problems, prove
\[
\exists 0<\delta_1,\delta_2<\frac{1}{4},\quad \mbox{ such that }\quad \mathbb P_{\mathrm{aux}}\bigl\{\phi<\bar\Phi\bigr\}\geq1-\delta_1,\qquad \mathbb P_{\mathrm{aux}}\bigl\{\phi_{S^c}>\bar\Phi\bigr\}\geq1-\delta_2.
\]
The corresponding CGMT tail comparisons then yield
\[
\mathbb P^{\otimes N}\bigl\{\widehat{\boldsymbol\beta}\in S^c\bigr\}\leq\mathbb P^{\otimes N}\bigl\{\Phi_{S^c}=\Phi\bigr\}\leq\mathbb P^{\otimes N}\bigl\{\Phi\geq\bar\Phi\bigr\}+\mathbb P^{\otimes N}\bigl\{\Phi_{S^c}\leq\bar\Phi\bigr\}\leq2\mathbb P_{\mathrm{aux}}\bigl\{\phi\geq\bar\Phi\bigr\}+2\mathbb P_{\mathrm{aux}}\bigl\{\phi_{S^c}\leq\bar\Phi\bigr\}\leq2\delta_1+2\delta_2.
\]

Let $(\boldsymbol\beta^\circ,\boldsymbol u^\circ)$ be a saddle point of $\phi$, where $\boldsymbol u^\circ$ is the optimal Lagrangian variable. Then, $\phi=F(\boldsymbol\beta^\circ,\boldsymbol u^\circ)$ and $\phi_{S^c}\geq\inf_{\boldsymbol\beta\in K\cap S^c}F(\boldsymbol\beta,\boldsymbol u^\circ)$. Hence, to prove the separation of the optimal values of the two auxiliary problems by the fixed real number $\bar\Phi$, it suffices to prove
\begin{align*}\mathbb P_{\mathrm{aux}}\bigl\{F(\boldsymbol\beta^\circ,\boldsymbol u^\circ)<\bar\Phi\bigr\}\geq1-\delta_1,\mbox{ and }\mathbb P_{\mathrm{aux}}\left\{\inf_{\boldsymbol\beta\in K\cap S^c}F(\boldsymbol\beta,\boldsymbol u^\circ)>\bar\Phi\right\}\geq1-\delta_2.\end{align*}
To analyze $F(\boldsymbol\beta^\circ,\boldsymbol u^\circ)$ and $F(\boldsymbol\beta,\boldsymbol u^\circ)$, we estimate the saddle point $(\boldsymbol\beta^\circ,\boldsymbol u^\circ)$.

\subsubsection{Saddle point analysis}
The analysis starts from the KKT conditions of $\phi$, that is,
\begin{align*}
0\in\partial_{\boldsymbol\beta}F(\boldsymbol\beta^\circ,\boldsymbol u^\circ)+N_K(\boldsymbol\beta^\circ),\mbox{ and }0\in-\partial_{\boldsymbol u}F(\boldsymbol\beta^\circ,\boldsymbol u^\circ)+N_U(\boldsymbol u^\circ).
\end{align*}
\paragraph{Boundary conditions.}
To solve the saddle point, one first determines the corresponding boundary conditions, that is, to compute $N_K(\bbeta^\circ)$ and $N_U(\vu^\circ)$. This is particularly important when the feasible set of the primal variable is a Euclidean ball, or when a one-dimensional reduction places the Lagrangian variable in a closed interval of the real line as in Section~\ref{sec:dimension_reduction_trick}. In the former case, whether $\boldsymbol\beta^\circ$ lies in the interior or on the boundary determines whether a nonnegative normal-cone multiplier is present in the primal KKT condition. In the latter case, the left endpoint, the interior, and the right endpoint correspond to different one-sided first-order conditions. Thus, the boundary configuration of the saddle point is determined before the KKT system is solved.
\paragraph{System of (nonlinear)-equations.}
Once the normal-cone conditions have been identified, we first evaluate the Lagrangian at $\boldsymbol u^\circ$ and solve the primal minimization problem
\begin{align*}
\boldsymbol\beta^\circ\in\mathop{\rm argmin}_{\boldsymbol\beta\in K}F(\boldsymbol\beta,\boldsymbol u^\circ).
\end{align*}
For auxiliary problems involving an $\ell_1$ regularization term, the coordinatewise KKT condition can usually be solved through the soft-thresholding identity \eqref{eq:softmax_equivalence}. Consequently, the primal KKT condition expresses $\boldsymbol\beta^\circ$ through a soft-thresholding formula that still contains the unknown optimal Lagrangian $\boldsymbol u^\circ$ and, when the primal constraint is active, the corresponding normal-cone multiplier. Thus, at this stage the primal minimizer is solved conditionally on the optimal Lagrangian.

We next estimate the optimal Lagrangian $\boldsymbol u^\circ$. Once the boundary analysis shows that $N_U(\vu^\circ)=\{\vzero\}$, its KKT condition reduces to
\begin{align*}
\partial_{\boldsymbol u}F(\boldsymbol\beta^\circ,\boldsymbol u^\circ)=0.
\end{align*}
Substituting the preceding soft-thresholding expression for $\boldsymbol\beta^\circ$ into this equation reduces the coupled KKT system of equations to an equation involving the optimal Lagrangian and the remaining scalar quantities.

\subsection{Tools from Convex Analysis}\label{sec:tools_convex_analysis}

For a closed convex set $K\subset\bR^n$, and some $\vx\in K$, we denote $N_K(\vx)=\{ \vs: \sup(\langle\vs,\vy\rangle:\vy\in K)\leq\langle\vs,\vx\rangle \}$ by its normal cone evaluated at $\vx$, see, for instance, \cite[Chapter III, pp. 136, Definition 5.2.3]{hiriart-urrutyConvexAnalysisMinimization1993}. For instance, when $K=rB_2^n$ for some $r>0$, then $N_{rB_2^n}(\vx)=\{\tau\vx:\tau\geq 0\}$, if $\|\vx\|_2=r$, and $N_{rB_2^n}(\vx)=\{\vzero\}$, if $\|x\|_2<r$.
We write $\partial^-f(\vx) =\{\vg\in\bR^n: f(\vz)-f(\vx)\geq\langle \vg,\vz-\vx\rangle,\,\forall \vz\in\bR^n\}$ by the sub-differential of a convex function $f$ evaluated at $\vx$, where $\vg\in\partial^-f(\vx)$ is called a sub-gradient of $f$ evaluated at $\vx$; and write $\nabla f(\vx)$ by the gradient of $f$ evaluated at $\vx$. We write $\nabla^2 f(\vx):\bR^n\to\bR^n$ by the Hessian of $f$ evaluated at $\vx$. For a convex minimization problem $\inf(f(\vx):\vx\in K)$, its KKT condition refers to the following equivalence: $\vx^\star\in K$ is the solution of $\inf(f(\vx):\vx\in K)$, if and only if $\vzero\in\partial^-f(\vx^\star)+N_K(\vx^\star)$, see \cite[Chapter VII, Theorem 1.1.1]{hiriart-urrutyConvexAnalysisMinimization1993}. We give a useful fact: when $K=[a,b]\subset\bR$ for some $a<b\in\bR$. Then $N_K(a)=\{x:x\leq 0\}$ and $N_K(b)=\{x:x\geq 0\}$. In particular, if $f$ is differentiable, then from the KKT condition $0\in f'(x^\star)+N_K(x^\star)$, there must hold $f'(a)\geq 0$ and $f'(b)\leq 0$. Let $(\cX,\|\cdot\|_\cX)$, $(\cY,\|\cdot\|_\cY)$ be two normed spaces, and $T:\cX\to\cY$, we let $\D T[\cdot]:\vx\in\cX\mapsto \D T[\cdot](\vx)\in L(\cX;\cY)$ be its Fr\'echet derivative, and in particular, when evaluated at each $\vx\in\cX$, $\D T[\cdot](\vx)$ is a bounded linear operator from $\cX$ to $\cY$.

\begin{Lemma}
\label{lem:common-saddle}
Let \(D\) be a nonempty compact convex subset of a Euclidean space and let \(Q:D\times[0,H]\to\mathbb R\) be continuous, convex in its first argument, and concave in its second argument. Then
\begin{align*}
&\min_{\boldsymbol w\in D}\max_{0\leq s\leq H}Q(\boldsymbol w,s)=\max_{0\leq s\leq H}\min_{\boldsymbol w\in D}Q(\boldsymbol w,s),
\end{align*}
and the saddle-point set is nonempty. If \(D=rB_2\), \((\boldsymbol w_H,s_H)\) is a saddle point, and \(Q\) is differentiable in \(s\), then \(0<s_H<H\) implies \(\partial_sQ(\boldsymbol w_H,s_H)=0\), and for any such $s_H$, there exists \(\gamma_H\geq0\) such that
\begin{align*}
&\boldsymbol 0\in\partial_{\boldsymbol w}Q(\boldsymbol w_H,s_H)+\gamma_H\boldsymbol w_H,\qquad \gamma_H(\|\boldsymbol w_H\|_2^2-r^2)=0.
\end{align*}
\end{Lemma}
\begin{proof}
    By Sion's theorem, see, for instance, \cite[Theorem VII. 4.3.1, pp. 344]{hiriart-urrutyConvexAnalysisMinimization1993}, the saddle-point set is non-empty, and by \cite[Section 4.3, VII, pp. 330]{hiriart-urrutyConvexAnalysisMinimization1993}, there holds $\min_{\boldsymbol w\in D}\max_{0\leq s\leq H}Q(\boldsymbol w,s)=\max_{0\leq s\leq H}\min_{\boldsymbol w\in D}Q(\boldsymbol w,s)$. Moreover, $\partial_sQ(\boldsymbol w_H,s_H)=0$ is the first order condition. When $D=rB_2$, by the KKT condition, see, for instance, \cite[Theorem VII  1.1.1, pp. 293 - 294]{hiriart-urrutyConvexAnalysisMinimization1993}, together with the fact that the normal cone of $rB_2$ evaluated at $\vw_H$ is given by $N_{rB_2}(\vw_H)=\{\gamma\vw_H:\gamma\geq 0\}$, there holds $\boldsymbol 0\in\partial_{\boldsymbol w}Q(\boldsymbol w_H,s_H)+\gamma_H\boldsymbol w_H$. The equation $\gamma_H(\|\boldsymbol w_H\|_2^2-r^2)=0$ is the complementary slackness condition.
\end{proof}

\begin{Lemma}\label{lemma:squared_hinge}
    Let $s,t,\lambda>0$, let $\vh\in V_{J^c}$, let $N\in\bN_+$. Then
    \begin{align*}
        \min_{\|\bbeta_{J^c}\|_2=t}\big( s\langle\vh,\bbeta_{J^c}\rangle + N\lambda\|\bbeta_{J^c}\|_1 \big) = t\begin{cases}
            N\lambda - s\|\vh\|_\infty, & \mbox{ if } s\|\vh\|_\infty\leq N\lambda,\\
            -\big( \sum_{j\in J^c} ( s|h_j| - N\lambda )_+^2 \big)^{\frac{1}{2}}, &\mbox{otherwise}.
        \end{cases}
    \end{align*}Moreover, $s\mapsto \min_{\|\bbeta_{J^c}\|_2=t}\big( s\langle\vh,\bbeta_{J^c}\rangle + N\lambda\|\bbeta_{J^c}\|_1 \big)$ is concave, and satisfies that
    \begin{align}\label{eq:Lipschitz_s_mapsto_min}
        \forall \vh,\vh'\in V_{J^c},\qquad \bigg|\min_{\|\bbeta_{J_\star^c}\|_2=t}\{s\langle\vh,\bbeta_{J_\star^c}\rangle+N\lambda\|\bbeta_{J_\star^c}\|_1\}-\min_{\|\bbeta_{J_\star^c}\|_2=t}\{s\langle\vh',\bbeta_{J_\star^c}\rangle+N\lambda\|\bbeta_{J_\star^c}\|_1\}\bigg|\leq st\|\vh-\vh'\|_2.
    \end{align}
\end{Lemma}
\begin{proof}
    Let $\vv\in V_{J^c}$, such that $v_j\geq 0$ for any $j\in J^c$, and $\|\vv\|_2=1$. Let $|\beta_j|=tv_j$ for any $j\in J^c$, then the problem is reduced to $t\min(\sum_{j\in J^c}(N\lambda - s|h_j|)v_j:\, v_j\geq 0,\|\vv\|_2=1)$. We separate into two cases:
    \begin{enumerate}
        \item when $s\|\vh\|_\infty\leq N\lambda$, then $N\lambda-s|h_j|\geq 0$ for any $j\in J^c$. One may check that the minimum is given by $t(N\lambda-s\|\vh\|_\infty)$;
        \item when $s\|\vh\|_\infty>N\lambda$, then $\sum_{j\in J^c}(N\lambda - s|h_j|)v_j\geq -\sum_{j\in J^c}(s|h_j|-N\lambda)_+v_j\geq -( \sum_{j\in J^c}(s|h_j|-N\lambda)_+^2 )^{1/2}$.
    \end{enumerate}The concavity of $s\mapsto \min_{\|\bbeta_{J^c}\|_2=t}\big( s\langle\vh,\bbeta_{J^c}\rangle + N\lambda\|\bbeta_{J^c}\|_1 \big)$ follows directly since it is a minimum over affine function of $s$. Last, the Lipschitzness follows from the following Lemma~\ref{lemma:Lipschitz_convex_concave_function}.
\end{proof}

\begin{Lemma}\label{lemma:Lipschitz_convex_concave_function}
Let $(\Theta,d)$ be a metric space, let $A,B$ be nonempty sets, and let $\Psi:\Theta\times A\times B\to\mathbb R$ satisfy
\[
|\Psi(\theta,a,b)-\Psi(\theta',a,b)|\leq Ld(\theta,\theta'),\qquad \theta,\theta'\in\Theta,\quad (a,b)\in A\times B.
\]
Then $\theta\mapsto\inf_{a\in A}\sup_{b\in B}\Psi(\theta,a,b)$ is $L$-Lipschitz.
\end{Lemma}

\begin{proof}
For every $\theta,\theta'\in\Theta$,
\[
\inf_{a\in A}\sup_{b\in B}\Psi(\theta,a,b)\leq\inf_{a\in A}\sup_{b\in B}\big\{\Psi(\theta',a,b)+Ld(\theta,\theta')\big\}=\inf_{a\in A}\sup_{b\in B}\Psi(\theta',a,b)+Ld(\theta,\theta').
\]
Exchanging $\theta$ and $\theta'$ proves the claim.
\end{proof}

The following lemma may be found in \cite{bonnansPerturbationAnalysisOptimization2000}, Theorem~4.13, pp.~273--274; or \cite{hiriart-urrutyConvexAnalysisMinimization1993}, Corollary~VI.4.4.5, p.~269.
\begin{Lemma}[Danksin's theorem]
Let \(E\) be a finite-dimensional Euclidean space, let \(K\) be a nonempty compact set, and let \(f:E\times K\to\mathbb R\) be continuous. Assume that, for every \(\boldsymbol a\in K\), \(f(\cdot,\boldsymbol a)\) is differentiable and that \((\boldsymbol u,\boldsymbol a)\mapsto\nabla_{\boldsymbol u}f(\boldsymbol u,\boldsymbol a)\) is continuous. Then, for every \(\boldsymbol u,\boldsymbol v\in E\),
\[
\mathbf D\Big(\boldsymbol u\mapsto\max_{\boldsymbol a\in K}f(\boldsymbol u,\boldsymbol a)\Big)(\boldsymbol u)[\boldsymbol v]=\max_{\boldsymbol a\in\operatorname*{arg\,max}_{\boldsymbol b\in K}f(\boldsymbol u,\boldsymbol b)}\big\langle\nabla_{\boldsymbol u}f(\boldsymbol u,\boldsymbol a),\boldsymbol v\big\rangle,
\]
and
\[
\mathbf D\Big(\boldsymbol u\mapsto\min_{\boldsymbol a\in K}f(\boldsymbol u,\boldsymbol a)\Big)(\boldsymbol u)[\boldsymbol v]=\min_{\boldsymbol a\in\operatorname*{arg\,min}_{\boldsymbol b\in K}f(\boldsymbol u,\boldsymbol b)}\big\langle\nabla_{\boldsymbol u}f(\boldsymbol u,\boldsymbol a),\boldsymbol v\big\rangle.
\]
If the maximizer is unique, then
\[
\nabla_{\boldsymbol u}\max_{\boldsymbol a\in K}f(\boldsymbol u,\boldsymbol a)=\nabla_{\boldsymbol u}f(\boldsymbol u,\boldsymbol a),\qquad \boldsymbol a\in\operatorname*{arg\,max}_{\boldsymbol b\in K}f(\boldsymbol u,\boldsymbol b),
\]
and if the minimizer is unique, then
\[
\nabla_{\boldsymbol u}\min_{\boldsymbol a\in K}f(\boldsymbol u,\boldsymbol a)=\nabla_{\boldsymbol u}f(\boldsymbol u,\boldsymbol a),\qquad \boldsymbol a\in\operatorname*{arg\,min}_{\boldsymbol b\in K}f(\boldsymbol u,\boldsymbol b).
\]
The same conclusions hold when \(E\) is a finite-dimensional matrix space equipped with the Frobenius inner product.
\end{Lemma}
The following lemma characterizes the non-expansiveness of the projection operator onto a convex body.
\begin{Lemma}[\cite{hiriart-urrutyConvexAnalysisMinimization1993}, Proposition~III.3.1.3, p.~118]
Let \(K\subset\mathbb R^N\) be a nonempty closed convex set. Then
\[
\forall\boldsymbol\mu_1,\boldsymbol\mu_2\in\mathbb R^N,\quad\big\|\operatorname{Proj}_K(\boldsymbol\mu_1)-\operatorname{Proj}_K(\boldsymbol\mu_2)\big\|_2^2\leq\big\langle\operatorname{Proj}_K(\boldsymbol\mu_1)-\operatorname{Proj}_K(\boldsymbol\mu_2),\boldsymbol\mu_1-\boldsymbol\mu_2\big\rangle.
\]
\end{Lemma}

\subsection{Soft-max and Gaussian variables calculations}

\subsubsection{Gaussian variables calculations}\label{sec:Gaussian_calculations}

For any $u\in\bR$, let $\Phi_G(u)=\bP\{G\geq u\}$. Then $\Phi_G(u)=\int_u^\infty \frac{1}{\sqrt{2\pi}}\exp(-\frac{1}{2}x^2)\mathrm{d}x$ and $\Phi_G'(u) = -\frac{1}{\sqrt{2\pi}}\exp(-\frac{1}{2}u^2)$.
For any $x>0$ and $t\in\bR$, there holds the well-known Mills ratio estimate
\begin{align}\label{eq:Gaussian_calculus}
    x\leq\frac{\frac{1}{\sqrt{2\pi}}\exp(-\frac{1}{2}x^2)}{\Phi_G(x)} \leq x+\frac{1}{x},\quad \frac{t}{t+x}\exp\big(-tx-\frac{1}{2}x^2\big) \leq \frac{\Phi_G(t+x)}{\Phi_G(t)} \leq \exp\big(-tx-\frac{1}{2}x^2\big).
\end{align}The first inequality may be found in \cite[pp. 164, Equation 3.7.2 and Lemma 3.7.4]{artstein-avidanAsymptoticGeometricAnalysis2021}; while the second is from the following observation: $\Phi_G'(u) = -\frac{1}{\sqrt{2\pi}}\exp(-\frac{1}{2}u^2)$, and from the fundamental theorem of calculus,
\begin{align*}
    \log\bigg(\frac{\Phi_G(t+x)}{\Phi_G(t)}\bigg) = -\int_t^{t+x}\frac{\frac{1}{\sqrt{2\pi}}\exp(-\frac{1}{2}u^2)}{\Phi_G(u)}\mathrm{d}u \leq -\int_t^{t+x}u \mathrm{d}u = -tx - \frac{1}{2}x^2.
\end{align*}Similarly, the lower side of the second inequality is verified.

We show some frequently used facts of \eqref{eq:Gaussian_calculus}.
\begin{enumerate}
    \item Let $x = \frac{2}{z_q}$, $t = z_q - \frac{2}{z_q}$, where $z_q = \Phi_G^{-1}(\frac{q}{2|J^c|})$, \eqref{eq:Gaussian_calculus} yields
    \begin{align}\label{eq:lower_Phi_G_z_Jc}
        \Phi_G\bigg(z_q - \frac{2}{z_q}\bigg) \geq \exp\bigg( 2 - \frac{2}{z_q^2} \bigg)\Phi_G(z_q) = \exp\bigg( 2 - \frac{2}{z_q^2} \bigg) \frac{q}{2|J^c|}.
    \end{align}
    On the other side, if $z^2 \geq 8$, then take $t = z-\frac{4}{z}$ and $x=\frac{4}{z}$ from \eqref{eq:Gaussian_calculus}, then $\frac{\Phi_G(z)}{\Phi_G(z-\frac{4}{z})} \geq (1-\frac{4}{z^2})\exp(-4+\frac{8}{z^2})\geq e^{-4}$, therefore,
    \begin{align}\label{eq:upper_Phi_G_z_4z}
        \forall z^2\geq 8,\qquad \Phi_G\big(z-\frac{4}{z}\big) \leq e^4\Phi_G(z).
    \end{align}
    
    \item \label{item:expectation_soft_max_square} For any $a>0$, $\bE(|G|-a)_+^2 = 2\int_0^\infty u\bP\{|G|-a>u\}\mathrm{d}u = 4\int_0^\infty u\Phi_G(a+u)\mathrm{d}u$. By \eqref{eq:Gaussian_calculus}, $\Phi_G(a+u)\leq \Phi_G(a)\exp(-au - \frac{1}{2}u^2)\leq \Phi_G(a)\exp(-au)$. Therefore,
    \begin{align}\label{eq:expectation_soft_max}
        \bE\big( |G|-a \big)_+^2\leq 4\Phi_G(a)\int_0^\infty u\exp(-au)\mathrm{d}u = \frac{4}{a^2}\Phi_G(a).
    \end{align}
    Moreover, if $a\geq2.5$, then for every $0\leq u\leq(2a)^{-1}$, \eqref{eq:Gaussian_calculus} gives $\Phi_G(a+u)\geq\frac14\Phi_G(a)$, and hence
    \begin{align}\label{eq:lower_expectation_soft_max}
        \bE\big(|G|-a\big)_+^2\geq4\int_0^{1/(2a)}u\Phi_G(a+u)\mathrm{d}u\geq\frac{\Phi_G(a)}{8a^2}.
    \end{align}
    The same argument gives $\bE(|G|-a)_+\leq2\Phi_G(a)/a$. Another fact about $\bE(|G|-a)_+^2$ that is used is its derivative: for any $a\geq0$,
    \begin{align}\label{eq:identity_derivative_exp_positive}
        \begin{aligned}
            &\frac{\mathrm{d}}{\mathrm{d}a}\bE\big(|G|-a\big)_+^2 = 2\frac{\mathrm{d}}{\mathrm{d}a}\int_a^\infty (x-a)^2\frac{\exp(-\frac{1}{2}x^2)}{\sqrt{2\pi}}\mathrm{d}x = -4\int_a^\infty (x-a)\frac{\exp(-\frac{1}{2}x^2)}{\sqrt{2\pi}}\mathrm{d}x\\
            &= -4\big( \frac{1}{\sqrt{2\pi}}\exp(-\frac{1}{2}a^2) - a\Phi_G(a) \big).
        \end{aligned}
    \end{align}
    By integration by parts applied to $\int_a^\infty \frac{x^2}{\sqrt{2\pi}}\exp(-\frac{1}{2}x^2)\mathrm{d}x = \Phi_G(a) + a\frac{\exp(-\frac{1}{2}a^2)}{\sqrt{2\pi}}$, the following identity holds:
    \begin{align*}
        \bE\big(|G|-a\big)_+^2 = 2(1+a^2)\Phi_G(a) - 2a\frac{\exp(-\frac{1}{2}a^2)}{\sqrt{2\pi}}.
    \end{align*}
    Consequently,
    \begin{align}\label{eq:identity_exp_square_soft_max}
        2\Phi_G(a) - \bE[(|G|-a)_+^2] = 2a\bigg( \frac{1}{\sqrt{2\pi}}\exp(-\frac{1}{2}a^2) - a\Phi_G(a) \bigg)>0,
    \end{align}where the last inequality follows from \eqref{eq:Gaussian_calculus}.
    Moreover, letting $a=z_q$, by $q=2|J_\star^c|\Phi_G(z_q)$, then combining \eqref{eq:expectation_soft_max} and \eqref{eq:lower_expectation_soft_max}, $\frac{d_q}{q}\asymp \frac{1}{z_q^2}$, where we recall that $d_q=|J_\star^c|\bE[(|G|-z_q)_+^2]$.

    \item \label{item:derivatives} Writing $\log\big(\frac{u\big(N-2|J^c|\Phi_G(u)\big)}{\sqrt{N}\big( N - |J^c|\bE[(|G|-u)_+^2] \big)^{\frac{1}{2}}}\big)$ as the sum of logarithms, together with the fact that $\Phi_G'(u)=-\frac{1}{\sqrt{2\pi}}\exp(-\frac{1}{2}u^2)$, we obtain
    \begin{align*}
        &\frac{\mathrm{d}}{\mathrm{d}u}\log\bigg(\frac{u\big(N-2|J^c|\Phi_G(u)\big)}{\sqrt{N}\big( N - |J^c|\bE[(|G|-u)_+^2] \big)^{\frac{1}{2}}}\bigg) = \frac{1}{u} + \frac{\frac{2|J^c|}{\sqrt{2\pi}}\exp(-\frac{1}{2}u^2)}{N-2|J^c|\Phi_G(u)} - \frac{2|J^c|\big( \frac{1}{\sqrt{2\pi}}\exp(-\frac{1}{2}u^2) - u\Phi_G(u) \big)}{ N-|J^c|\bE[(|G|-u)_+^2] }.
    \end{align*}
    Combining \eqref{eq:identity_exp_square_soft_max},
    \begin{align*}
        &\frac{\frac{2|J^c|}{\sqrt{2\pi}}\exp(-\frac{1}{2}u^2)}{N-2|J^c|\Phi_G(u)} - \frac{2|J^c|\big( \frac{1}{\sqrt{2\pi}}\exp(-\frac{1}{2}u^2) - u\Phi_G(u) \big)}{ N-|J^c|\bE[(|G|-u)_+^2] }\\
        &= \frac{2|J^c|u\Phi_G(u)}{N-2|J^c|\Phi_G(u)} + \frac{\big( 2|J^c|\Phi_G(u) - |J^c|\bE[(|G|-u)_+^2] \big)^2}{ u\big( N - 2|J^c|\Phi_G(u) \big)\big( N - |J^c|\bE[(|G|-u)_+^2] \big) }.
    \end{align*}
    Moreover, by \eqref{eq:identity_exp_square_soft_max} again,
    \begin{align*}
        &\frac{\mathrm{d}}{\mathrm{d}u}\bigg(\frac{u\big(N-2|J^c|\Phi_G(u)\big)}{ \big( N - |J^c|\bE[(|G|-u)_+^2] \big)^{\frac{1}{2}}}\bigg)\\
        &= \frac{N-2|J^c|\Phi_G(u) + \frac{2|J^c|}{\sqrt{2\pi}}u\exp(-\frac{1}{2}u^2)}{  \big( N - |J^c|\bE[(|G|-u)_+^2] \big)^{\frac{1}{2}} } - \frac{ 2|J^c|u\big( N-2|J^c|\Phi_G(u) \big)\big( \frac{1}{\sqrt{2\pi}}\exp(-\frac{1}{2}u^2) - u\Phi_G(u) \big) }{ \big( N - |J^c|\bE[(|G|-u)_+^2] \big)^{\frac{3}{2}}}\\
        &= \frac{ N+2|J^c|(u^2-1)\Phi_G(u) }{  \big( N - |J^c|\bE[(|G|-u)_+^2] \big)^{\frac{1}{2}} } + \frac{ \big( 2|J^c|\Phi_G(u) - |J^c|\bE[(|G|-u)_+^2] \big)^2 }{ \big( N - |J^c|\bE[(|G|-u)_+^2] \big)^{\frac{3}{2}}}.
    \end{align*}

    \item For any $q$ such that $z_q = \Phi_G^{-1}(\frac{q}{2|J^c|})>2.57$,
    \begin{align}\label{eq:estimate_z_q}
        \log\big( \frac{e|J^c|}{q} \big) \leq z_q^2 \leq 2\log\big(\frac{|J^c|}{q}\big).
    \end{align}
    To see this, let $t=0$ and $x=z_q$ in \eqref{eq:Gaussian_calculus}, together with $\Phi_G(0)=\frac{1}{2}$, we obtain that $\frac{q}{2|J^c|}=\Phi_G(z_q)\leq \frac{1}{2}\exp(-\frac{1}{2}z_q^2)$, that is, $z_q^2\leq 2\log(\frac{|J^c|}{q})$. On the other side, from $\frac{\frac{1}{\sqrt{2\pi}}\exp(-\frac{1}{2}x^2)}{\Phi_G(x)} \leq x+\frac{1}{x}$, we obtain that $2\Phi_G(x) \geq \frac{2}{\sqrt{2\pi}(x+\frac{1}{x})}\exp(-\frac{1}{2}x^2)$. When $x\geq 2.57$, there holds $\frac{2}{\sqrt{2\pi}(x+\frac{1}{x})}\exp(-\frac{1}{2}x^2)\geq \exp(1-x^2)$, hence $\frac{q}{|J^c|} = 2\Phi_G(z_q)\geq \exp(1-z_q^2)$, that is, $z_q^2\geq \log(\frac{e|J^c|}{q})$.

    \item \label{item:expectation_soft_max_moment} By integration by parts, together with the fact that $\int_a^\infty x\frac{\exp(-\frac{1}{2}x^2)}{\sqrt{2\pi}}\mathrm{d}x = \frac{1}{\sqrt{2\pi}}\exp(-\frac{1}{2}a^2)$, and the fact that $\int_a^\infty x^2 \frac{\exp(-\frac{1}{2}x^2)}{\sqrt{2\pi}}\mathrm{d}x = a\frac{\exp(-\frac{1}{2}a^2)}{\sqrt{2\pi}} + \Phi_G(a)$, which follows again from integration by parts,
    \begin{align}\label{eq:exp_G_soft_max}
        \begin{aligned}
            &\bE[|G|(|G|-a)_+] = 2\int_a^\infty x(x-a)\frac{\exp(-\frac{1}{2}x^2)}{\sqrt{2\pi}}\mathrm{d}x = 2\bigg( \int_a^\infty x^2 \frac{\exp(-\frac{1}{2}x^2)}{\sqrt{2\pi}}\mathrm{d}x - a\int_a^\infty x\frac{\exp(-\frac{1}{2}x^2)}{\sqrt{2\pi}}\mathrm{d}x \bigg)\\
            &= 2\Phi_G(a).
        \end{aligned}
    \end{align}
    Moreover, by \eqref{eq:Gaussian_calculus}, for any $u\geq 0$,
    \begin{align*}
        \bP\big\{|G|>a+u\big||G|>a\big\} = \frac{\Phi_G(a+u)}{\Phi_G(a)} \leq \exp\big(-\frac{1}{2}u^2-au\big)\leq \exp\big(-\frac{1}{2}u(a+u)\big).
    \end{align*}
    Therefore, $\bP\big\{|G|(|G|-a)_+>t\big||G|>a\big\}\leq\exp(-t/2)$ for every $t\geq0$, and for any $k\geq1$,
    \begin{align*}
        \bE\big[ |G|^k(|G|-a)_+^k \big| |G|>a \big]\leq k\int_0^\infty t^{k-1}\exp(-\frac{1}{2}t)\mathrm{d}t = 2^k k!,
    \end{align*}
    and hence $\bE\big[ |G|^k(|G|-a)_+^k \big]\leq 2\Phi_G(a)2^kk!$. In particular, for every integer $k\geq2$,
    \begin{align}\label{eq:Bernstein_moment_soft_max}
        \bE\big|\,|G|(|G|-a)_+-2\Phi_G(a)\big|^k\leq 2^k\bE\big[|G|^k(|G|-a)_+^k\big]\leq\frac{k!}{2}\big(64\Phi_G(a)\big)4^{k-2}.
    \end{align}

    \item \label{item:Gaussian_local_comparison} For any $z>2.57$ and $u$ such that $|u-z|\leq(2^8z)^{-1}$, applying \eqref{eq:Gaussian_calculus} separately when $u\geq z$ and $u\leq z$ gives
    \begin{align}\label{eq:Gaussian_local_comparison}
        e^{-1/8}\Phi_G(z)\leq\Phi_G(u)\leq e^{1/8}\Phi_G(z).
    \end{align}
    Moreover, by \eqref{eq:identity_derivative_exp_positive}, \eqref{eq:lower_expectation_soft_max}, and $\bE(|G|-a)_+\leq2\Phi_G(a)/a$, for every $a\geq2.5$,
    \begin{align*}
        0\leq-\frac{\mathrm{d}}{\mathrm{d}a}\log\bE(|G|-a)_+^2=\frac{2\bE(|G|-a)_+}{\bE(|G|-a)_+^2}\leq32a.
    \end{align*}
    Since every point between $u$ and $z$ is larger than $2.5$, the fundamental theorem of calculus yields
    \begin{align*}
        \left|\log\frac{\bE(|G|-u)_+^2}{\bE(|G|-z)_+^2}\right|\leq32\bigg(z+\frac{1}{2^8z}\bigg)|u-z|<\frac14.
    \end{align*}
    Consequently,
    \begin{align}\label{eq:Gaussian_local_comparison_soft_max}
        2^{-4}\bE(|G|-z)_+^2\leq\bE(|G|-u)_+^2\leq2^4\bE(|G|-z)_+^2.
    \end{align}

    \item\label{item:Gaussian_profile_comparison} For any $z>2.57$, set $\bar z=z+(2^8z)^{-1}$ and $a=(6z)^{-1}$. Then
    \begin{align}\label{eq:Gaussian_profile_comparison}
    \sup_{0\leq v\leq\bar z}\frac{a\varphi_G((v-a)_+)}{\Phi_G(v)}<\frac{23}{100},\qquad \frac{\Phi_G(\bar z)}{\Phi_G(z)}>\frac{99}{100}.
    \end{align}
    Indeed, write $r(v)=\varphi_G(v)/\Phi_G(v)$. By \eqref{eq:Gaussian_calculus}, $r(v)\geq v$, and hence $r'(v)=r(v)(r(v)-v)\geq0$. If $0\leq v\leq a$, then $a<1$ and \eqref{eq:Gaussian_calculus} at $1$ gives
    \begin{align*}
    \frac{a\varphi_G(0)}{\Phi_G(v)}\leq\frac{a\varphi_G(0)}{\Phi_G(1)}\leq2ae^{1/2}\leq\frac{e^{1/2}}{3\cdot2.57}<\frac{23}{100}.
    \end{align*}
    If $a<v\leq\bar z$, monotonicity of $r$ and \eqref{eq:Gaussian_calculus} yield
    \begin{align*}
    \frac{a\varphi_G(v-a)}{\Phi_G(v)}=ae^{av-a^2/2}r(v)\leq ae^{a\bar z}r(\bar z)\leq ae^{a\bar z}(\bar z+\bar z^{-1})=\frac{e^{1/6+1/(1536z^2)}}{6}\bigg(1+\frac1{256z^2}+\frac1{z^2+1/256}\bigg)<\frac{23}{100},
    \end{align*}
    where the last expression is decreasing for $z>0$ and is smaller than $23/100$ at $z=2.57$. Finally, applying the lower bound in \eqref{eq:Gaussian_calculus} with $t=z$ and $x=(2^8z)^{-1}$ gives
    \begin{align*}
    \frac{\Phi_G(\bar z)}{\Phi_G(z)}\geq\frac{z}{\bar z}\exp\bigg(-\frac1{256}-\frac1{2\cdot256^2z^2}\bigg)=\frac{\exp\{-1/256-1/(2\cdot256^2z^2)\}}{1+1/(256z^2)}>\frac{99}{100},
    \end{align*}
    where the last expression is increasing in $z$ and is larger than $99/100$ at $z=2.57$.
\end{enumerate}

\subsubsection{Softmax calculations}
For \(u\geq0\), write \(S_u(x)=\operatorname{sign}(x)(|x|-u)_+\) be the soft-max function at level $u$. For any $h,b\in\bR$, $\tau\geq 0$, the following inequalities hold
\begin{align}\label{eq:softmax_calculus}
    \big( S_\tau(h+b)-b \big)h \geq |h|\big( |h| - \tau - |b| \big)_+ - b^2 - bh,\quad\mbox{ and } \quad\big|S_\tau(h+b)\big|^2 \leq 2b^2 + 2\big(|h|-\tau\big)_+^2.
\end{align}
For the first inequality, it suffices to prove $hS_\tau(h+b)\geq |h|(|h|-\tau-|b|)_+-b^2.$ If $h(h+b)\geq0$, then $hS_\tau(h+b)=|h|(|h+b|-\tau)_+\geq |h|(|h|-|b|-\tau)_+,$ where we used $|h+b|\geq |h|-|b|$. If $h(h+b)<0$, then $|b|=|h|+|h+b|$, so $(|h|-\tau-|b|)_+=0$, while $hS_\tau(h+b)\geq-|h||h+b|\geq-b^2.$
This proves the first inequality after subtracting $bh$ from both sides.

For the second inequality, $|S_\tau(h+b)|=(|h+b|-\tau)_+\leq(|h|-\tau)_++|b|.$
Hence, by $(x+y)^2\leq2x^2+2y^2$, $|S_\tau(h+b)|^2\leq2b^2+2(|h|-\tau)_+^2.$

We also use the following equivalence in this work
\begin{align}\label{eq:softmax_equivalence}
    & 0\in x-y+\eta\partial^-|\cdot|(x) \quad\Leftrightarrow\quad x=S_\eta(y),\qquad \forall \eta\geq0.
\end{align}To see this, we separate into three cases: (i) when $x>0$, then from the left-hand-side, $x=y-\eta$, and by $x>0$, we have $y>\eta$. Hence $S_\eta(y)=y-\eta=x$; (ii) when $x<0$, from the left-hand-side, $x=y+\eta$, and from $x<0$, we know $y<-\eta$, hence $S_\eta(y)=y+\eta=x$; (iii) $x=0$, then $0\in -y+\eta[-1,1]$ is equivalent to $|y|\leq \eta$, which is further equivalent to $S_\eta(y)=0=x$.

Moreover, the following scaling identity holds:
\begin{align}\label{eq:scaling_identity_softmax}
    \forall c>0,\,\eta\geq 0,\,y\in\bR,\qquad S_{c\eta}(cy) = \sign(cy)\big(|cy|-c\eta\big)_+ = \sign(y)c\big(|y|-\eta\big)_+ = cS_\eta(y).
\end{align}

When $x\neq\pm u$, it is easy to check that $S_u'(x) = \1_{\{|x|>u\}}$. Since $\{\pm u\}$ is a null set, there holds $\bE[S_u'(g)]=\bP\{|g|>u\}=2\Phi_G(u)$, where $g\sim\cN(0,1)$.

\begin{Lemma}\label{lem:v_rho_u}
Let $G\sim\mathcal N(0,1)$. For every $\rho,u\geq0$, recall that
\begin{align*}
v_{\rho,u}=\Phi_G^{-1}\left(\frac{1}{2|J_\star^c|}\sum_{j\in J_\star^c}\mathbb P\bigl\{|G+\rho\beta_j^\star|>u\bigr\}\right).
\end{align*}
Then $u\geq r_{\rho,u}$. Moreover, if $v_{\rho,u}\leq\bar z_q$, then
\begin{align*}
\sum_{j\in J_\star^c}\mathbb P\left\{\left|S_u(\rho\beta_j^\star-G)-\rho\beta_j^\star\right|\geq\frac1{6z_q}\right\}>\frac34q.
\end{align*}
\end{Lemma}

\begin{proof}We first prove $u\geq v_{\rho,u}$. Since for any $u>0$, $\bP\{|G+\rho\beta_j^\star|>u\}\geq \bP\{|G|>u\}=2\Phi_G(u)$, $\frac{1}{2|J^c|}\sum_{j\in J^c}\bP\{|G+\rho\beta_j^\star|>u\}\geq\Phi_G(u)$. By definition, the left-hand-side is precisely $v_{\rho,u}=\Phi_G^{-1}(\frac{1}{2|J_\star^c|}\sum_{j\in J_\star^c}\mathbb P\{|G+\rho\beta_j^\star|>u\})$, hence $\Phi_G(v_{\rho,u})\geq\Phi_G(u)$, and by the monotonicity of $\Phi_G$, there holds $v_{\rho,u}\leq u$.

Fix $\rho,u\geq0$. For every $j\in J_\star^c$, on $\{y\in\bR:\,|\rho\beta_j^\star-y|>u\}$, either $y<\rho\beta_j^\star-u$ or $y>\rho\beta_j^\star+u$. In the first case, $S_u(\rho\beta_j^\star-y)-\rho\beta_j^\star=-y-u$, while in the second case, $S_u(\rho\beta_j^\star-y)-\rho\beta_j^\star=u-y$. Hence
\begin{align*}
&\big\{y:\left|S_u(\rho\beta_j^\star-y)-\rho\beta_j^\star\right|<\frac1{6z_q},\mbox{ and } |\rho\beta_j^\star-y|>u\big\}\\
&\subset \big\{y:y<\rho\beta_j^\star-u,\ |y+u|<\frac1{6z_q}\big\}\cup \big\{y:y>\rho\beta_j^\star+u,\ |y-u|<\frac1{6z_q}\big\}.
\end{align*}Consequently,
\begin{align*}
&\mathbb P\big\{\big|S_u(\rho\beta_j^\star-G)-\rho\beta_j^\star\big|\geq\frac1{6z_q}\big\} \geq \bP\big\{ \big|\rho\beta_j^\star-G\big|>u, \, \big| S_u(\rho\beta_j^\star-G)-\rho\beta_j^\star \big|\geq \frac{1}{6z_q} \big\}\\
&=\bP\{|\rho\beta_j^\star-G|>u\} - \bP\big\{\left|S_u(\rho\beta_j^\star-G)-\rho\beta_j^\star\right|<\frac1{6z_q},\, |\rho\beta_j^\star-G|>u\big\}\\
&\geq \bP\{|\rho\beta_j^\star-G|>u\} - \bP\bigg\{\big\{G<\rho\beta_j^\star-u,\ |G+u|<\frac1{6z_q}\big\}\cup \big\{G>\rho\beta_j^\star+u,\ |G-u|<\frac1{6z_q}\big\}\bigg\}.
\end{align*}The interval $\big\{y:y<\rho\beta_j^\star-u,\ |y+u|<\frac1{6z_q}\big\}\cup \big\{y:y>\rho\beta_j^\star+u,\ |y-u|<\frac1{6z_q}\big\}$ has length (under the Lebesgue measure) at most $\frac{2}{6z_q}$. Moreover, on the first interval $|y|\geq u-|y+u|>u-(6z_q)^{-1}$, and on the second interval $|y|\geq u-|y-u|>u-(6z_q)^{-1}$; thus every point in their union has absolute value at least $(u-(6z_q)^{-1})_+$. Hence
\begin{align*}
    &\bP\bigg\{\big\{G<\rho\beta_j^\star-u,\ |G+u|<\frac1{6z_q}\big\}\cup \big\{G>\rho\beta_j^\star+u,\ |G-u|<\frac1{6z_q}\big\}\bigg\} \leq \frac{2}{6z_q}\varphi_G((u-(6z_q)^{-1})).
\end{align*}
Together with the fact that $\mathbb P\{|G+\rho\beta_j^\star|>u\}\geq2\Phi_G(u)$; and by the definition of $v_{\rho,u}$, this yields $u\geq v_{\rho,u}$, there holds
\begin{align*}
&\sum_{j\in J_\star^c}\mathbb P\left\{\left|S_u(\rho\beta_j^\star-G)-\rho\beta_j^\star\right|\geq\frac1{6z_q}\right\}
\geq\sum_{j\in J_\star^c}\mathbb P\{|G+\rho\beta_j^\star|>u\}-\frac{2|J_\star^c|}{6z_q}\varphi_G\left((u-(6z_q)^{-1})_+\right)\\
&\geq2|J_\star^c|\Phi_G(v_{\rho,u})-\frac{2|J_\star^c|}{6z_q}\varphi_G\left((v_{\rho,u}-(6z_q)^{-1})_+\right) = 2|J_\star^c|\Phi_G(v_{\rho,u})\left(1-\frac{(6z_q)^{-1}\varphi_G((v_{\rho,u}-(6z_q)^{-1})_+)}{\Phi_G(v_{\rho,u})}\right).
\end{align*}
Item~\ref{item:Gaussian_profile_comparison} of Section~\ref{sec:Gaussian_calculations} gives, for $v_{\rho,u}\leq\bar z_q$,
\begin{align*}
\frac{(6z_q)^{-1}\varphi_G((v_{\rho,u}-(6z_q)^{-1})_+)}{\Phi_G(v_{\rho,u})}<\frac{23}{100},\qquad \frac{\Phi_G(\bar z_q)}{\Phi_G(z_q)}>\frac{99}{100}.
\end{align*}
Since $\Phi_G$ is decreasing and $2|J_\star^c|\Phi_G(z_q)=q$,
\begin{align*}
\sum_{j\in J_\star^c}\mathbb P\left\{\left|S_u(\rho\beta_j^\star-G)-\rho\beta_j^\star\right|\geq\frac1{6z_q}\right\}>\frac{77}{100}\frac{99}{100}q>\frac34q.
\end{align*}
\end{proof}

\subsection{A probabilistic lemma}

\begin{Lemma}\label{lem:probabilistic_lemma}
Let \(M\in\mathbb N\), \(x\geq90\), and \(V\in\{1,\ldots,2M\}\). Let \(Z_1,\ldots,Z_M\) be independent real random variables, let \(\widetilde Z_1,\ldots,\widetilde Z_M\) be independent copies, and let \(\mathcal C\) be a family of \(M\)-tuples \(C=(C_1,\ldots,C_M)\) of Borel subsets of \(\mathbb R\). Assume that, for every \((z_j,\widetilde z_j)_{j=1}^M\in\mathbb R^{2M}\), the following holds:
\begin{align*}&\left|\left\{\left(\mathbf 1_{\{z_j\in C_j\}},\mathbf 1_{\{\widetilde z_j\in C_j\}}\right)_{j=1}^M:C\in\mathcal C\right\}\right|\leq\left(\frac{2eM}{V}\right)^V.\end{align*}
Then
\begin{align*}&\mathbb P\left\{\exists C\in\mathcal C:\sum_{j=1}^M\mathbb P\{Z_j\in C_j\}\geq x,\ \sum_{j=1}^M\mathbf 1_{\{Z_j\in C_j\}}<\frac23\sum_{j=1}^M\mathbb P\{Z_j\in C_j\}\right\}\leq2\left(\frac{2eM}{V}\right)^V e^{-x/1728},\\
&\mathbb P\left\{\exists C\in\mathcal C:\sum_{j=1}^M\mathbb P\{Z_j\in C_j\}\leq x,\ \sum_{j=1}^M\mathbf 1_{\{Z_j\in C_j\}}>2x\right\}\leq3\left(\frac{2eM}{V}\right)^V e^{-x/196}.\end{align*}
\end{Lemma}

\begin{proof}Without loss of generality, we assume that such a set $C \in \mathcal{C}$ exists (otherwise the probability is $0$). We prove only the first conclusion, since the second one is similar.
Fix \(C\in\mathcal C\) such that \(\sum_{j=1}^M\mathbb P\{Z_j\in C_j\}\geq x\). Independence gives \(\mathbb E\sum_{j=1}^M\mathbf 1_{\{\widetilde Z_j\in C_j\}}=\sum_{j=1}^M\mathbb P\{Z_j\in C_j\}\) and \(\operatorname{Var}(\sum_{j=1}^M\mathbf 1_{\{\widetilde Z_j\in C_j\}})\leq\sum_{j=1}^M\mathbb P\{Z_j\in C_j\}\). Hence Chebyshev's inequality yields
\begin{align*}&\mathbb P_{\widetilde Z}\left\{\frac56\sum_{j=1}^M\mathbb P\{Z_j\in C_j\}\leq\sum_{j=1}^M\mathbf 1_{\{\widetilde Z_j\in C_j\}}\leq\frac76\sum_{j=1}^M\mathbb P\{Z_j\in C_j\}\right\}\geq1-\frac{36}{\sum_{j=1}^M\mathbb P\{Z_j\in C_j\}}\geq\frac35>\frac12.\end{align*}
Consequently, for every \((z_j)_{j=1}^M\in\mathbb R^M\),
\begin{align*}&\mathbf 1_{\left\{\exists C\in\mathcal C:\,\sum_{j=1}^M\mathbb P\{Z_j\in C_j\}\geq x,\ \sum_{j=1}^M\mathbf 1_{\{z_j\in C_j\}}<\frac23\sum_{j=1}^M\mathbb P\{Z_j\in C_j\}\right\}}\\
&\leq2\mathbb P_{\widetilde Z}\Bigg\{\exists C\in\mathcal C:\sum_{j=1}^M\mathbb P\{Z_j\in C_j\}\geq x,\ \sum_{j=1}^M\mathbf 1_{\{z_j\in C_j\}}<\frac23\sum_{j=1}^M\mathbb P\{Z_j\in C_j\},\mbox{ and}\\
&\qquad\frac56\sum_{j=1}^M\mathbb P\{Z_j\in C_j\}\leq\sum_{j=1}^M\mathbf 1_{\{\widetilde Z_j\in C_j\}}\leq\frac76\sum_{j=1}^M\mathbb P\{Z_j\in C_j\}\Bigg\}.\end{align*}
Taking probability with respect to $\bP_Z$ on both sides, and notice that for each $C\in\cC$ such that the right-hand-side random event holds, we have: $\sum_{j=1}^M(\mathbf 1_{\{\widetilde Z_j\in C_j\}}-\mathbf 1_{\{Z_j\in C_j\}})>\frac16\sum_{j=1}^M\mathbb P\{Z_j\in C_j\}\geq\frac{x}{6}$, and $\sum_{j=1}^M\big|\mathbf 1_{\{\widetilde Z_j\in C_j\}}-\mathbf 1_{\{Z_j\in C_j\}}\big|\leq\sum_{j=1}^M\mathbf 1_{\{\widetilde Z_j\in C_j\}}+\sum_{j=1}^M\mathbf 1_{\{Z_j\in C_j\}}<\frac{11}{6}\sum_{j=1}^M\mathbb P\{Z_j\in C_j\}$. Therefore, 
\begin{align*}
    \sum_{j=1}^M\left(\mathbf 1_{\{\widetilde Z_j\in C_j\}}-\mathbf 1_{\{Z_j\in C_j\}}\right)>\frac1{12}\sum_{j=1}^M\left|\mathbf 1_{\{\widetilde Z_j\in C_j\}}-\mathbf 1_{\{Z_j\in C_j\}}\right|.
\end{align*}
Moreover,
\begin{align*}&\sum_{j=1}^M\left|\mathbf 1_{\{\widetilde Z_j\in C_j\}}-\mathbf 1_{\{Z_j\in C_j\}}\right|\geq\left|\sum_{j=1}^M\left(\mathbf 1_{\{\widetilde Z_j\in C_j\}}-\mathbf 1_{\{Z_j\in C_j\}}\right)\right|>\frac{x}{6}.\end{align*}

Let \((\varepsilon_j)_{j=1}^M\) be independent Rademacher variables independent of \((Z_j,\widetilde Z_j)_{j=1}^M\). Then
\begin{align*}
&\bP\left\{\exists C\in\mathcal C:\,\sum_{j=1}^M\mathbb P\{Z_j\in C_j\}\geq x,\ \sum_{j=1}^M\mathbf 1_{\{z_j\in C_j\}}<\frac23\sum_{j=1}^M\mathbb P\{Z_j\in C_j\}\right\}\\
&\leq 2\mathbb P_{Z,\widetilde Z}\left\{\exists C\in\mathcal C:\sum_{j=1}^M\left(\mathbf 1_{\{\widetilde Z_j\in C_j\}}-\mathbf 1_{\{Z_j\in C_j\}}\right)>\frac1{12}\sum_{j=1}^M\left|\mathbf 1_{\{\widetilde Z_j\in C_j\}}-\mathbf 1_{\{Z_j\in C_j\}}\right|,\ \sum_{j=1}^M\left|\mathbf 1_{\{\widetilde Z_j\in C_j\}}-\mathbf 1_{\{Z_j\in C_j\}}\right|>\frac{x}{6}\right\}\\
&=2\mathbb E_{Z,\widetilde Z}\mathbb P_{\boldsymbol\varepsilon}\left\{\exists C\in\mathcal C:\sum_{j=1}^M\varepsilon_j\left(\mathbf 1_{\{\widetilde Z_j\in C_j\}}-\mathbf 1_{\{Z_j\in C_j\}}\right)>\frac1{12}\sum_{j=1}^M\left|\mathbf 1_{\{\widetilde Z_j\in C_j\}}-\mathbf 1_{\{Z_j\in C_j\}}\right|,\ \sum_{j=1}^M\left|\mathbf 1_{\{\widetilde Z_j\in C_j\}}-\mathbf 1_{\{Z_j\in C_j\}}\right|>\frac{x}{6}\right\}.\end{align*}
Fix \((z_j,\widetilde z_j)_{j=1}^M\in\mathbb R^{2M}\) and \(C\in\mathcal C\). Since
\(\sum_{j=1}^M(\mathbf 1_{\{\widetilde z_j\in C_j\}}-\mathbf 1_{\{z_j\in C_j\}})^2=\sum_{j=1}^M|\mathbf 1_{\{\widetilde z_j\in C_j\}}-\mathbf 1_{\{z_j\in C_j\}}|\), Hoeffding's inequality gives
\begin{align*}&\mathbb P_{\boldsymbol\varepsilon}\left\{\sum_{j=1}^M\varepsilon_j\left(\mathbf 1_{\{\widetilde z_j\in C_j\}}-\mathbf 1_{\{z_j\in C_j\}}\right)>\frac1{12}\sum_{j=1}^M\left|\mathbf 1_{\{\widetilde z_j\in C_j\}}-\mathbf 1_{\{z_j\in C_j\}}\right|,\ \sum_{j=1}^M\left|\mathbf 1_{\{\widetilde z_j\in C_j\}}-\mathbf 1_{\{z_j\in C_j\}}\right|>\frac{x}{6}\right\}\\
&=\mathbf 1_{\left\{\sum_{j=1}^M|\mathbf 1_{\{\widetilde z_j\in C_j\}}-\mathbf 1_{\{z_j\in C_j\}}|>x/6\right\}}\mathbb P_{\boldsymbol\varepsilon}\left\{\sum_{j=1}^M\varepsilon_j\left(\mathbf 1_{\{\widetilde z_j\in C_j\}}-\mathbf 1_{\{z_j\in C_j\}}\right)>\frac1{12}\sum_{j=1}^M\left|\mathbf 1_{\{\widetilde z_j\in C_j\}}-\mathbf 1_{\{z_j\in C_j\}}\right|\right\}\\
&\leq\mathbf 1_{\left\{\sum_{j=1}^M|\mathbf 1_{\{\widetilde z_j\in C_j\}}-\mathbf 1_{\{z_j\in C_j\}}|>x/6\right\}}\exp\left(-\frac1{288}\sum_{j=1}^M\left|\mathbf 1_{\{\widetilde z_j\in C_j\}}-\mathbf 1_{\{z_j\in C_j\}}\right|\right)\leq e^{-x/1728}.\end{align*}
\end{proof}

\subsubsection{Estimates of \texorpdfstring{$\Var_{J_\star^c}$}{free variance}}

\begin{Lemma}
\label{lem:common-calibration}
Assume $|J_\star^c|>N$ and $\sigma_j>0$ for every $j\in J_\star^c$. For every $\lambda>0$, equation~\eqref{eq:common-q-calibration} has a unique solution $q\in(0,N)$. For this solution,
\begin{align*}
0<d_q<q<N,\qquad \mathrm{Var}_{J_\star^c}=\frac{\sigma_\xi^2d_q}{N-d_q},\qquad d_q=\frac{N\mathrm{Var}_{J_\star^c}}{\sigma_\xi^2+\mathrm{Var}_{J_\star^c}},\qquad \frac{N\lambda}{z_q}=\frac{\sigma_\xi(N-q)}{\sqrt{N-d_q}}.
\end{align*}
For fixed $J_\star^c$ and fixed $(\sigma_j)_{j\in J_\star^c}$, $\lambda\mapsto z_{q(\lambda)}$ is strictly increasing and $\lambda\mapsto\mathrm{Var}_{J_\star^c}$ is strictly decreasing on $(0,\infty)$. At $\lambda=0$, the continuous extension is $q(0)=N$ and $\mathrm{Var}_{J_\star^c}=\sigma_\xi^2d_N/(N-d_N)$.
For the remaining numerical estimates assume $\Sigma_{J_\star^c}=I_{V_{J_\star^c}}$ and $|J_\star^c|\geq100N$. With \(\bar z_q=z_q+\frac{1}{2^8z_q}\) and \(\underline z_q=z_q-\frac{1}{2^8z_q}\),
\begin{align}
&z_q>2.57,\qquad d_q=\frac{N\mathrm{Var}_{J_\star^c}}{\sigma_\xi^2+\mathrm{Var}_{J_\star^c}}\leq\frac{2q}{z_q^2}\leq\frac{2q}{\log(e|J_\star^c|/q)},\qquad \frac{d_q}{N}\leq\frac{2}{z_q^2}<0.31,\qquad \mathrm{Var}_{J_\star^c}\leq\frac{3\sigma_\xi^2q}{N\log(e|J_\star^c|/q)}.
\label{eq:common-calibration-identities}
\end{align}
Moreover,
\begin{align}
&\frac{\bar z_q(N-2|J_\star^c|\Phi_G(\bar z_q))}{\sqrt N(N-|J_\star^c|\mathbb E(|G|-\bar z_q)_+^2)^{1/2}}-\frac{\sqrt N\lambda}{\sigma_\xi}\geq2^{-10}\bigg(\frac1{z_q}+\frac{qz_q}{N}\bigg),\mbox{ and }\nonumber\\
&\frac{\sqrt N\lambda}{\sigma_\xi}-\frac{\underline z_q(N-2|J_\star^c|\Phi_G(\underline z_q))}{\sqrt N(N-|J_\star^c|\mathbb E(|G|-\underline z_q)_+^2)^{1/2}}\geq2^{-10}\bigg(\frac1{z_q}+\frac{qz_q}{N}\bigg).
\label{eq:common-endpoint-calibration}
\end{align}
Moreover,
\begin{align}
&\frac{N\lambda}{z_q}=\frac{\sigma_\xi(N-q)}{(N-d_q)^{1/2}},\mbox{ and } \sigma_\xi\sqrt N-\frac{N\lambda}{\bar z_q}\geq2^{-10}\sigma_\xi\sqrt N\bigg(\frac qN+\frac1{z_q^2}\bigg)\geq2^{-11}\sigma_\xi\sqrt{d_q}.
\label{eq:common-radial-calibration}
\end{align}
\end{Lemma}

\begin{proof}
For every $a>0$, integration by parts gives $\mathbb E[|G|(|G|-a)_+]=2\Phi_G(a)$. Consequently,
\begin{align*}
0<\mathbb E(|G|-a)_+^2=2\Phi_G(a)-a\mathbb E(|G|-a)_+<2\Phi_G(a).
\end{align*}
Differentiation under the integral is justified on each compact subinterval of $(0,\infty)$ by Gaussian integrability. Therefore,
\begin{align*}
u\frac{\mathrm d}{\mathrm du}\sum_{j\in J_\star^c}\mathbb E\bigg(|G|-\frac{u}{\sqrt{\sigma_j}}\bigg)_+^2=2\sum_{j\in J_\star^c}\bigg[\mathbb E\bigg(|G|-\frac{u}{\sqrt{\sigma_j}}\bigg)_+^2-2\Phi_G\bigg(\frac{u}{\sqrt{\sigma_j}}\bigg)\bigg].
\end{align*}
The function $u\mapsto\sum_{j\in J_\star^c}2\Phi_G(u/\sqrt{\sigma_j})$ is continuous and strictly decreasing from $|J_\star^c|$ to $0$. This proves that every $z_r$ is well defined, and that $u>z_N$ is equivalent to $\sum_{j\in J_\star^c}2\Phi_G(u/\sqrt{\sigma_j})<N$. On this interval the preceding identities yield
\begin{align}
&\frac{\mathrm d}{\mathrm du}\bigg[\frac{u\bigg(N-\sum_{j\in J_\star^c}2\Phi_G\big(\frac{u}{\sqrt{\sigma_j}}\big)\bigg)}{\sqrt N\bigg(N-\sum_{j\in J_\star^c}\mathbb E\big(|G|-\frac{u}{\sqrt{\sigma_j}}\big)_+^2\bigg)^{1/2}}\bigg]\notag\\
&=\frac{\bigg(N-\sum_{j\in J_\star^c}2\Phi_G\big(\frac{u}{\sqrt{\sigma_j}}\big)\bigg)^2}{\sqrt N\bigg(N-\sum_{j\in J_\star^c}\mathbb E\big(|G|-\frac{u}{\sqrt{\sigma_j}}\big)_+^2\bigg)^{3/2}}+\frac{u\sum_{j\in J_\star^c}\frac{2}{\sqrt{2\pi\sigma_j}}\exp\big(-\frac{u^2}{2\sigma_j}\big)}{\sqrt N\bigg(N-\sum_{j\in J_\star^c}\mathbb E\big(|G|-\frac{u}{\sqrt{\sigma_j}}\big)_+^2\bigg)^{1/2}}>0.\label{eq:monotonocity_calibration_threshold}
\end{align}
The expression being differentiated tends to $0$ as $u\downarrow z_N$, since $d_N<N$, and tends to $\infty$ as $u\to\infty$. It therefore takes the value $\sqrt N\lambda/\sigma_\xi$ exactly once. The bijection $q\mapsto z_q$ between $(0,N)$ and $(z_N,\infty)$ proves existence and uniqueness. The inequalities $0<d_q<q<N$ follow from the scalar inequality above, and the two identities involving $\mathrm{Var}_{J_\star^c}$ follow by rearrangement. Strict monotonicity in $\lambda$ follows from the positive derivative and the strict decrease of $u\mapsto\sum_{j\in J_\star^c}\mathbb E(|G|-u/\sqrt{\sigma_j})_+^2$. Continuity at $\lambda=0$ gives $q(0)=N$.

For the remaining estimates we work under the stated isotropic restriction.
Since \(q<N\) and \(|J_\star^c|\geq100N\), \(z_q>\Phi_G^{-1}(1/200)>2.57\). By \eqref{eq:expectation_soft_max} and \eqref{eq:estimate_z_q}, $d_q\leq\frac{4|J_\star^c|\Phi_G(z_q)}{z_q^2}=\frac{2q}{z_q^2}\leq\frac{2q}{\log(e|J_\star^c|/q)}$, and $\frac{d_q}{N}\leq\frac{2}{z_q^2}<0.31$. Moreover, \(\Var_{J_\star^c}=\sigma_\xi^2d_q/(N-d_q)\) gives \(N\Var_{J_\star^c}=d_q(\sigma_\xi^2+\Var_{J_\star^c})\), and hence $d_q=\frac{N\Var_{J_\star^c}}{\sigma_\xi^2+\Var_{J_\star^c}}$, and $\Var_{J_\star^c}\leq\frac{3\sigma_\xi^2q}{N\log(e|J_\star^c|/q)}$.

For \(2.5<\underline z_q\leq u\leq\bar z_q\), we know $u^2>1$, and $\sqrt{N}(N-|J_\star^c|\bE(|G|-u)_+^2)^{1/2}\leq N$; hence item~\emph{\ref{item:derivatives}} of Section~\ref{sec:Gaussian_calculations} gives
\begin{align*}
&\frac{\mathrm d}{\mathrm du}\bigg(\frac{u\bigl(N-2|J_\star^c|\Phi_G(u)\bigr)}{\sqrt N\bigl(N-|J_\star^c|\bE(|G|-u)_+^2\bigr)^{1/2}}\bigg)\geq \frac{ N+2|J_\star^c|(u^2-1)\Phi_G(u) }{ \sqrt N \big( N - |J_\star^c|\bE[(|G|-u)_+^2] \big)^{\frac{1}{2}} } \geq 1+\frac{2|J_\star^c|(u^2-1)\Phi_G(u)}{N}.
\end{align*}
By \(\underline z_q = z_q - \frac{1}{2^8z_q}\leq u\leq\bar z_q\) and \(z_q>2.57\), we have \(u^2-1\ \geq (z_q - \frac{1}{2^8z_q})^2-1\geq \frac{1}{2}z_q^2\). Moreover, item~\emph{\ref{item:Gaussian_local_comparison}} of Section~\ref{sec:Gaussian_calculations} gives \(\Phi_G(u)\geq e^{-1/8}\Phi_G(z_q)\geq\Phi_G(z_q)/2=q/(4|J_\star^c|)\). Consequently,
\begin{align*}
&\frac{\mathrm d}{\mathrm du}\bigg(\frac{u\bigl(N-2|J_\star^c|\Phi_G(u)\bigr)}{\sqrt N\bigl(N-|J_\star^c|\bE(|G|-u)_+^2\bigr)^{1/2}}\bigg)\geq 1+\frac{qz_q^2}{4N}.
\end{align*}

We next prove \eqref{eq:common-endpoint-calibration}. By the fundamental theorem of calculus,
\begin{align*}
&\frac{\bar z_q(N-2|J_\star^c|\Phi_G(\bar z_q))}{\sqrt N\bigl(N-|J_\star^c|\bE(|G|-\bar z_q)_+^2\bigr)^{1/2}}-\frac{z_q\bigl(N-2|J_\star^c|\Phi_G(z_q)\bigr)}{\sqrt N\bigl(N-|J_\star^c|\bE(|G|-z_q)_+^2\bigr)^{1/2}}\\
&=\int_{z_q}^{\bar z_q}\frac{\mathrm d}{\mathrm du}\bigg(\frac{u\bigl(N-2|J_\star^c|\Phi_G(u)\bigr)}{\sqrt N\bigl(N-|J_\star^c|\bE(|G|-u)_+^2\bigr)^{1/2}}\bigg)\,\mathrm du\geq\frac1{2^8z_q}\bigg(1+\frac{qz_q^2}{4N}\bigg).
\end{align*}
Similarly,
\begin{align*}
&\frac{z_q\bigl(N-2|J_\star^c|\Phi_G(z_q)\bigr)}{\sqrt N\bigl(N-|J_\star^c|\bE(|G|-z_q)_+^2\bigr)^{1/2}}-\frac{\underline z_q(N-2|J_\star^c|\Phi_G(\underline z_q))}{\sqrt N\bigl(N-|J_\star^c|\bE(|G|-\underline z_q)_+^2\bigr)^{1/2}}\\
&=\int_{\underline z_q}^{z_q}\frac{\mathrm d}{\mathrm du}\bigg(\frac{u\bigl(N-2|J_\star^c|\Phi_G(u)\bigr)}{\sqrt N\bigl(N-|J_\star^c|\bE(|G|-u)_+^2\bigr)^{1/2}}\bigg)\,\mathrm du\geq\frac1{2^8z_q}\bigg(1+\frac{qz_q^2}{4N}\bigg).
\end{align*}
Since \eqref{eq:common-q-calibration} gives
\[
\frac{z_q\bigl(N-2|J_\star^c|\Phi_G(z_q)\bigr)}{\sqrt N\bigl(N-|J_\star^c|\bE(|G|-z_q)_+^2\bigr)^{1/2}}=\frac{z_q(N-q)}{\sqrt N(N-d_q)^{1/2}}=\frac{\sqrt N\lambda}{\sigma_\xi},
\]
\eqref{eq:common-endpoint-calibration} follows.

Finally, \eqref{eq:common-q-calibration} gives \(\frac{N\lambda}{z_q}=\frac{\sigma_\xi(N-q)}{(N-d_q)^{1/2}}\). Since \(d_q\leq2q/z_q^2<q/2\), there holds \(1-\frac{1-q/N}{(1-d_q/N)^{1/2}}>\frac{q}{2N}\), while \(1-z_q/\bar z_q=(2^8z_q^2+1)^{-1}\geq(257z_q^2)^{-1}\) and \(z_q/\bar z_q\geq1/2\). Therefore
\[
\sigma_\xi\sqrt N-\frac{N\lambda}{\bar z_q}=\sigma_\xi\sqrt N\bigg(1-\frac{z_q}{\bar z_q}+\frac{z_q}{\bar z_q}\bigg(1-\frac{1-q/N}{(1-d_q/N)^{1/2}}\bigg)\bigg)\geq2^{-10}\sigma_\xi\sqrt N\bigg(\frac qN+\frac1{z_q^2}\bigg).
\]
Moreover, \(N(q/N+z_q^{-2})^2\geq4q/z_q^2\geq2d_q\), and thus \(\sigma_\xi\sqrt N-\frac{N\lambda}{\bar z_q}\geq2^{-11}\sigma_\xi\sqrt{d_q}\).
\end{proof}

\begin{Lemma}\label{lem:estimate_free_variance}
    Assume $\Sigma_{J_\star^c}=I_{V_{J_\star^c}}$.
    Let $W:\bR\to\bR$ be the Lambert $W$ function, that is, the function that satisfies $W(x)\exp(W(x))=x$ for any $x\geq -\frac{1}{e}$. Then for any $0<\varepsilon<2$ such that $\varepsilon\log(e|J^c|/N)\geq 1$ and $\frac{\sigma_\xi}{\sqrt{N}}\leq \lambda \leq \sigma_\xi\sqrt{\frac{\varepsilon\log(e|J^c|/N)}{N}}$,
        \begin{align*}
            \Var_{J_\star^c} \asymp_\varepsilon \frac{\sigma_\xi^2}{\log(e|J^c|/N)}.
        \end{align*}
\end{Lemma}

\begin{proof}[Proof of Lemma~\ref{lem:estimate_free_variance}]
Fix \(0<\varepsilon<2\), and choose \(c_\varepsilon=\min\{\frac{1-\frac{1}{2}\varepsilon}{8}, \frac{5}{12\sqrt{2\pi e}}\sqrt{\frac{1-\frac{1}{2}\varepsilon}{\varepsilon}}\}\). Suppose first that \(q\leq c_\varepsilon N\). By \eqref{eq:common-radial-calibration}, $\frac{\sqrt{N}\lambda}{\sigma_\xi}=\frac{z_q(1-q/N)}{\sqrt{1-d_q/N}} \geq z_q(1-q/N)$; by the upper bound on \(\lambda\), \(z_q(1-q/N)\leq\sqrt N\lambda/\sigma_\xi\leq\sqrt{\varepsilon\log(e|J^c|/N)}\), and therefore \(z_q\leq\sqrt{\varepsilon\log(e|J^c|/N)}/(1-c_\varepsilon)\). Since \(z_q>2.57\), we have \(z_q+z_q^{-1}\leq6z_q/5\). By the lower side of \eqref{eq:Gaussian_calculus},
\begin{align*}
&c_\varepsilon\geq \frac qN=\frac{2|J^c|}{N}\Phi_G(z_q)\geq\frac{5(1-c_\varepsilon)}{3e\sqrt{2\pi\varepsilon\log(e|J^c|/N)}}\exp\bigg\{\bigg(1-\frac{\varepsilon}{2(1-c_\varepsilon)^2}\bigg)\log\bigg(\frac{e|J^c|}{N}\bigg)\bigg\}.
\end{align*}
Since $c_\varepsilon\leq \frac{1}{8}$, and $(1-c_\varepsilon)^{-2}\leq 1+4c_\varepsilon$, we have $1-\frac{\varepsilon}{2(1-c_\varepsilon)^2}\geq 1-\frac{1}{2}\varepsilon - 2\varepsilon c_\varepsilon \geq \frac{1}{2}(1-\frac{1}{2}\varepsilon)$. Together with $1-c_\varepsilon\geq\frac{1}{2}$ and elementary inequality $\frac{1}{\sqrt{x}}\exp(\frac{1}{2}(1-\frac{1}{2}\varepsilon)x)\geq\sqrt{e(1-\frac{1}{2}\varepsilon)}$, we obtain
\begin{align*}
    c_\varepsilon \geq \frac{q}{N} \geq \frac{5}{6e\sqrt{2\pi\varepsilon\log(e|J_\star^c|/N)}}\exp\Big(\frac{1-\frac{1}{2}\varepsilon}{2}\log(e|J_\star^c|/N)\Big) \geq \frac{5}{6\sqrt{2\pi e}}\sqrt{\frac{1-\varepsilon/2}{\varepsilon}} \geq 2c_\varepsilon,
\end{align*}leading to a contradiction. Therefore $c_\varepsilon N < q < N$.

It follows from \eqref{eq:estimate_z_q} that \(\log(e|J^c|/N)\leq z_q^2\leq2\log(|J^c|/q)\lesssim_\varepsilon\log(e|J^c|/N)\). Combining this with \(c_\varepsilon N<q<N\) and \(q/(16z_q^2)\leq d_q\leq2q/z_q^2\) gives \(d_q/N\asymp_\varepsilon1/\log(e|J^c|/N)\). Finally, \(d_q/N<0.31\) by \eqref{eq:common-calibration-identities}, and therefore \(\Var_{J_\star^c}/\sigma_\xi^2=d_q/(N-d_q)\asymp d_q/N\), proving the first assertion.
\end{proof}

\begin{Lemma}\label{lem:monotonicity_total_fsd_variance}
Grant the same assumption as in Proposition~\ref{prop:monotonicity_bias_variance}, for any $\lambda_0\leq\lambda_1\leq\lambda_2$, such that $|J_\star^c(\lambda_1)|\ge 100N$,
\begin{align*}
    \mathrm{VAR}^{\mathrm{FSD}}(\lambda_2) \leq \mathrm{VAR}^{\mathrm{FSD}}(\lambda_1);
\end{align*} if and only if for any $\lambda>\lambda_0$ such that $\lambda = \sigma_j|\beta_j^\star|$ for some $j\in[p]$, there holds
    \begin{align*}
        d_{q(\lambda)} - d_{q(\lambda_-)} \leq \Big|\Big\{j\in[p]:\, \sigma_j|\beta_j^\star|=\lambda\Big\}\Big|\Big(1-\frac{d_{q(\lambda)}}{N}\Big)\Big(1-\frac{d_{q(\lambda_-)}}{N}\Big).
    \end{align*}
\end{Lemma}

\begin{proof}
We prove the equivalence for $\frac{\sigma_\xi^2|J_\star(\lambda)|}{N}+\Var_{J_\star(\lambda)^c}$.
\begin{enumerate}
\item\label{item:fsd-variance-fixed-set} The goal of this step is to prove monotonicity whenever $J_\star$ is constant for both $\lambda_1$ and $\lambda_2$. For $\lambda_0\leq\lambda_1\leq\lambda_2$ with $\lambda_1>0$ and $J_\star(\lambda_1)=J_\star(\lambda_2)$, Lemma~\ref{lem:common-calibration} gives $\lambda\mapsto \Var_{J_\star^c}(\lambda)$ is strictly decreasing, while $\sigma_\xi^2\frac{J_\star(\lambda_1)}{N}=\sigma_\xi^2\frac{J_\star(\lambda_2)}{N}$.  Consequently,
\begin{align}\begin{aligned}&\mathrm{VAR}^{\mathrm{FSD}}(\lambda_2)-\mathrm{VAR}^{\mathrm{FSD}}(\lambda_1)=\frac{\sigma_\xi^2N\big(d_{q(\lambda_2)}-d_{q(\lambda_1)}\big)}{(N-d_{q(\lambda_2)})(N-d_{q(\lambda_1)})}\leq0.\end{aligned}\label{eq:fsd-variance-fixed-set}\end{align}

\item\label{item:fsd-variance-left-limit} For every $\lambda\in(\lambda_0,\infty)\cap\{\sigma_j|\beta_j^\star|:j\in[p]\}$, $J_\star(\lambda)\subseteq J_\star(\lambda-)$ and $J_\star(\lambda-)\setminus J_\star(\lambda)=\{j\in[p]:\sigma_j|\beta_j^\star|=\lambda\}$. Hence
\begin{align}\begin{aligned}&\mathrm{VAR}^{\mathrm{FSD}}(\lambda)-\mathrm{VAR}^{\mathrm{FSD}}(\lambda-)=-\frac{\sigma_\xi^2}{N}|J_\star(\lambda-)\setminus J_\star(\lambda)|+\frac{\sigma_\xi^2N\big(d_{q(\lambda)}-d_{q(\lambda-)}\big)}{(N-d_{q(\lambda)})(N-d_{q(\lambda-)})}.\end{aligned}\label{eq:fsd-variance-left-difference}\end{align}
Both factors in the denominator are positive by Lemma~\ref{lem:common-calibration}. Multiplying \eqref{eq:fsd-variance-left-difference} by $\frac{N}{\sigma_\xi^2}$ shows that $\mathrm{VAR}^{\mathrm{FSD}}(\lambda)\leq\mathrm{VAR}^{\mathrm{FSD}}(\lambda-)$ is equivalent to
\begin{align}\begin{aligned}&d_{q(\lambda)}-d_{q(\lambda-)}\leq|J_\star(\lambda-)\setminus J_\star(\lambda)|\bigg(1-\frac{d_{q(\lambda)}}N\bigg)\bigg(1-\frac{d_{q(\lambda-)}}N\bigg).\end{aligned}\label{eq:fsd-variance-diagonal-criterion}\end{align}
If $\mathrm{VAR}^{\mathrm{FSD}}$ is non-increasing on $[\lambda_0,\infty)$, then $\mathrm{VAR}^{\mathrm{FSD}}(\lambda)\leq\mathrm{VAR}^{\mathrm{FSD}}(t)$ for every $\lambda_0\leq t<\lambda$. Passing to $t\uparrow\lambda$ proves \eqref{eq:fsd-variance-diagonal-criterion}.

Conversely, assume \eqref{eq:fsd-variance-diagonal-criterion} for all the stated values of $\lambda$. For $\lambda_0\leq\lambda_1<\lambda_2$, we perform the following induction. Suppose the following holds: for any $\lambda_0\leq \lambda_1<\lambda_2$, if $|(\lambda_1,\lambda_2]\cap\{\sigma_j|\beta_j^\star|:j\in[p]\}|\leq k-1$, and $J_\star^c(\lambda_1)\geq 100N$, then $\mathrm{VAR}^{\mathrm{FSD}}(\lambda_2)\leq\mathrm{VAR}^{\mathrm{FSD}}(\lambda_1)$. 
\begin{enumerate}
\item\label{item:fsd-variance-induction-base} If $k=1$, \eqref{eq:fsd-variance-fixed-set} and the continuity established in item~\ref{item:fsd-variance-fixed-set} give the assertion.

\item\label{item:fsd-variance-induction-step} Suppose the induction hypothesis holds for $k-1$, then for any $\lambda_1<\lambda\leq\lambda_2$ such that there exists $j\in[p]$ such that $\lambda=\sigma_j|\beta_j^\star|$, we have $|(\lambda,\lambda_2]\cap\{\sigma_j|\beta_j^\star|:j\in[p]\}|\leq k-1$. By the hypothesis, $\mathrm{VAR}^{\mathrm{FSD}}(\lambda_2)\leq\mathrm{VAR}^{\mathrm{FSD}}(\lambda)$. Morevoer, for any $\lambda_1\le t < \lambda$, $|(\lambda_1,t]\cap\{\sigma_j|\beta_j^\star|:j\in[p]\}|\leq k-1$. By the hypothesis, $\mathrm{VAR}^{\mathrm{FSD}}(t)\leq\mathrm{VAR}^{\mathrm{FSD}}(\lambda_1)$. Send $t\uparrow\lambda$, $\mathrm{VAR}^{\mathrm{FSD}}(\lambda_-)\leq \mathrm{VAR}^{\mathrm{FSD}}(\lambda_1)$. Combining \eqref{eq:fsd-variance-left-difference} and \eqref{eq:fsd-variance-diagonal-criterion}, $\mathrm{VAR}^{\mathrm{FSD}}(\lambda)\leq \mathrm{VAR}^{\mathrm{FSD}}(\lambda_-)$. Therefore $\mathrm{VAR}^{\mathrm{FSD}}(\lambda_2)\leq \mathrm{VAR}^{\mathrm{FSD}}(\lambda_1)$.
\end{enumerate}
\end{enumerate}
\end{proof}

\section{Diagonal covariance in the well-regularized regime}\label{sec:large-diagonal-proof}

In this section, we give the modifications of Sections~\ref{sec:proof_upper_bound} and~\ref{sec:proof_lower_bound} for $\Sigma=\operatorname{diag}(\sigma_1,\ldots,\sigma_p)$ with $\sigma_j>0$. Fix an absolute constant $C\ge8$ for which Proposition~\ref{prop:RIP} below holds with deviation constant $C$ and probability at least $1-C\exp(-\frac{|J_\star|}{C})$, and for which Lemma~\ref{lem:large-cross-process} below holds, and set $\varepsilon=2^{-14}$, $C_\star=256(1+C)$, and $A=2^{34}(1+C)^4$. In the diagonal cases of Assumptions~\ref{assumption:LASSO_upper} and~\ref{assumption:LASSO_lower}, take $c=A^{-32}$ and recall that by assumption,
\begin{align*} |J_\star^c|>N,\qquad 1\le|J_\star|\le cN,\qquad r^2(V_{J_\star},V_{J_\star^c})\le c\sigma_\xi^2,\quad\mbox{and}\quad \lambda\ge\frac{8\sigma_\xi z_N^{\mathrm{all}}}{\sqrt N}. \end{align*}
The individual propositions below state their sufficient assumptions. The calibration replacements used throughout are collected in Section~\ref{sec:diagonal-calibration}.
In this section, we only present the differences from the proofs in the isotropic case, namely those in Sections~\ref{sec:proof_upper_bound} and~\ref{sec:proof_lower_bound}.

\subsection{Upper bound}\label{sec:diagonal-upper}
In the case $\Sigma = \diag(\sigma_j : j \in [p])$, the proof of the upper bound is similar to that in the isotropic design case. We use $\|\Sigma^{1/2}(\hat\bbeta-\bbeta^\star)\|_2^2 \le \|\Sigma_{J_\star}^{1/2}(\hat\bbeta_{J_\star} - \bbeta_{J_\star}^\star)\|_2^2 + 2\|\Sigma_{J_\star^c}^{1/2}\hat\bbeta_{J_\star^c}\|_2^2 + 2\|\Sigma_{J_\star^c}^{1/2}\bbeta_{J_\star^c}^\star\|_2^2$. For the high-probability upper bound on $\|\Sigma_{J_\star}^{1/2}(\hat\bbeta_{J_\star} - \bbeta_{J_\star}^\star)\|_2^2$, we establish its effective-regularization property via Proposition~\ref{prop:effective-regularization-lasso}, namely by identifying $\hat\bbeta_{J_\star}$ as an RERM of an effective risk, and obtaining $\|\Sigma_{J_\star^c}^{1/2}\hat\bbeta_{J_\star^c}\|_2^2 = \|\Sigma_{J_\star^c}^{1/2}\cA[\bxi+\vzeta]\|_2^2$, where $\vzeta = \bX_{J_\star}(\bbeta_{J_\star}^\star-\hat\bbeta_{J_\star}) + \bX_{J_\star^c}\bbeta_{J_\star^c}^\star$. Here, for the free part of the estimator, namely $\hat\bbeta_{J_\star^c}$, we directly use Proposition~\ref{prop:lambda_large_shift} to obtain the overall upper bound, rather than, as in the isotropic case, establishing an upper bound on the discrepancy $\|\Sigma_{J_\star^c}^{1/2}(\cA[\bxi+\vzeta]-\cA[\bxi])\|_2^2$ via RIP; and for the high-probability upper bound on $\|\Sigma^{1/2}(\hat\bbeta_{J_\star} - \bbeta_{J_\star}^\star)\|_2^2$, we establish the local Bernstein's property and thereby the one-sided RIP.

For the process of establishing the local Bernstein's property and proving the one-sided RIP, the main difference lies in the establishment of the local Bernstein's property — the latter is almost a verbatim copy of the proof in the isotropic case. Here, since in the anisotropic case we are not aware of any substitute for Gluskin's theorem, after obtaining \eqref{eq:pre_lower_bound_local_Bernstein}, we cannot establish, as in \eqref{eq:v-properties}, the high-probability upper bound $\frac{C_0}{\sqrt{\log(\frac{e|J_\star^c|}{N})}}$ on $\|v_{\vlambda}\|_2$. Instead, here we directly examine the Hessian of $\vlambda^\star[\cdot]$.

\begin{Propositionextend}{prop:local-bernstein}
Grant the diagonal, well-regularized case of Assumption~\ref{assumption:LASSO_upper}. {Fix $C_\beta=16$ and $c_{\mathrm{LB}}=\frac{1}{32}$.} For $C_\beta\|{\Sigma_{J_\star^c}^{1/2}}\boldsymbol\beta_{J_\star^c}^\star\|_2\le r\le\sigma_\xi$, with probability at least $1-2\exp(-\frac{N}{8})$ with respect to $\mathbb X_{J_\star^c}$,
\begin{align*} \forall\boldsymbol\beta_{J_\star}\in{\boldsymbol\beta_{J_\star}^\star+\Sigma_{J_\star}^{-1/2}S_2^{J_\star}(\boldsymbol0;r)},\qquad P\mathcal L_{\boldsymbol\beta_{J_\star}}^{(\lambda,J_\star)}\ge c_{\mathrm{LB}}\frac{r^2}{\lambda}. \end{align*}
\end{Propositionextend}
\begin{proof}[Proof of Proposition~\ref{prop:local-bernstein}.]
Replace the geometric event of Lemma~\ref{lem:geometric-event} by
\begin{align*} \Omega_{\mathrm{act}}=\Omega_{\mathrm{gp}}\cap\bigg\{\forall\boldsymbol z\in2\sigma_\xi\sqrt N B_2^N,\ |\operatorname{supp}(\mathcal A[\boldsymbol z])|\le\frac{5N}{16}\bigg\},\quad\mbox{and}\quad \Omega_{\mathrm{signal}}=\big\{\|\mathbb X_{J_\star^c}\boldsymbol\beta_{J_\star^c}^\star\|_2\le2\sqrt N\|\Sigma_{J_\star^c}^{1/2}\boldsymbol\beta_{J_\star^c}^\star\|_2\big\}. \end{align*}
Here $\Omega_{\mathrm{gp}}$ is the general-position event for $\mathbb X_{J_\star^c}$. Lemma~\ref{lem:diagonal-active-support} gives $\mathbb P(\Omega_{\mathrm{act}}^c)\le\exp(-\frac{N}{8})$; the Gaussian norm bound gives $\mathbb P(\Omega_{\mathrm{signal}}^c)\le\exp(-\frac{N}{2})$. Consequently,
\begin{align*} \Omega=\Omega_{\mathrm{act}}\cap\Omega_{\mathrm{signal}},\qquad\mbox{and}\qquad \mathbb P(\Omega)\ge1-2\exp\bigg(-\frac{N}{8}\bigg). \end{align*}
Work on $\Omega$. Fix $r>0$ and $\boldsymbol\beta_{J_\star}\in{\boldsymbol\beta_{J_\star}^\star+\Sigma_{J_\star}^{-1/2}S_2^{J_\star}(\boldsymbol0;r)}$. Conditionally on $\mathbb X_{J_\star^c}$, $\boldsymbol\xi-\mathbb X_{J_\star}(\boldsymbol\beta_{J_\star}-\boldsymbol\beta_{J_\star}^\star)$ follows the centered Gaussian distribution with covariance matrix
\begin{align*} {\big(\sigma_\xi^2+\|\Sigma_{J_\star}^{1/2}(\boldsymbol\beta_{J_\star}-\boldsymbol\beta_{J_\star}^\star)\|_2^2\big)I_N}=(\sigma_\xi^2+r^2)I_N. \end{align*}
Let $\boldsymbol G\sim\mathcal N(\boldsymbol0,I_N)$. Then $P\ell_{\boldsymbol\beta_{J_\star}}^{(\lambda,J_\star)}$ has the following offset Gaussian mean width representation:
\begin{align*} P\ell_{\boldsymbol\beta_{J_\star}}^{(\lambda,J_\star)}=\mathbb E_{\boldsymbol G}\max_{\|\mathbb X_{J_\star^c}^{\top}\boldsymbol u\|_\infty\le1}\bigg\{\big\langle\boldsymbol u,\mathbb X_{J_\star^c}\boldsymbol\beta_{J_\star^c}^\star+\sqrt{\sigma_\xi^2+r^2}\boldsymbol G\big\rangle-\frac{N\lambda}{2}\|\boldsymbol u\|_2^2\bigg\}. \tag{\ref{eq:conditonal_population_risk}} \end{align*}
The expression for $P\ell_{\boldsymbol\beta_{J_\star}^\star}^{(\lambda,J_\star)}$ is obtained by replacing $r^2$ by zero. For every feasible $\boldsymbol u$,
\begin{align*} \langle\boldsymbol u,\mathbb X_{J_\star^c}\boldsymbol\beta_{J_\star^c}^\star+\sigma_\xi\boldsymbol G\rangle-\frac{N\lambda}{2}\|\boldsymbol u\|_2^2&=\frac{\sigma_\xi}{\sqrt{\sigma_\xi^2+r^2}}\bigg(\langle\boldsymbol u,\mathbb X_{J_\star^c}\boldsymbol\beta_{J_\star^c}^\star+\sqrt{\sigma_\xi^2+r^2}\boldsymbol G\rangle-\frac{N\lambda}{2}\|\boldsymbol u\|_2^2\bigg)\\ &\quad+\bigg(1-\frac{\sigma_\xi}{\sqrt{\sigma_\xi^2+r^2}}\bigg)\bigg(\langle\boldsymbol u,\mathbb X_{J_\star^c}\boldsymbol\beta_{J_\star^c}^\star\rangle-\frac{N\lambda}{2}\|\boldsymbol u\|_2^2\bigg). \end{align*}
Taking the maximum and then expectation, and using $P\ell_{\boldsymbol\beta_{J_\star}}^{(\lambda,J_\star)}\ge P\ell_{\boldsymbol\beta_{J_\star}^\star}^{(\lambda,J_\star)}$ by Gaussian convolution and Jensen's inequality, gives the same reduction as~\eqref{eq:pre_lower_bound_local_Bernstein}:
\begin{align*} \tag{\ref{eq:pre_lower_bound_local_Bernstein}} P\mathcal L_{\boldsymbol\beta_{J_\star}}^{(\lambda,J_\star)}\ge\bigg(1-\frac{\sigma_\xi}{\sqrt{\sigma_\xi^2+r^2}}\bigg)\bigg[P\ell_{\boldsymbol\beta_{J_\star}^\star}^{(\lambda,J_\star)}-\max_{\|\mathbb X_{J_\star^c}^{\top}\boldsymbol u\|_\infty\le1}\bigg\{\langle\boldsymbol u,\mathbb X_{J_\star^c}\boldsymbol\beta_{J_\star^c}^\star\rangle-\frac{N\lambda}{2}\|\boldsymbol u\|_2^2\bigg\}\bigg]. \end{align*}

For any $I\subset J_\star^c$, $|I|\leq N$ and any $\vveps\in\{-1,1\}^I$, let
\begin{align*}
    &\cC_{I,\vveps} := \bigg\{\vz\in\bR^N:\, \varepsilon_j\Big[(\bX_I^\top\bX_I)^{-1}(\bX_I^\top\vz-N\lambda\vveps)\Big]_j>0,\, \forall\, j\in I,\\
    &\qquad\qquad\mbox{ and }\, \Big|\Big\langle \bX\ve_j, \vz - \bX_I(\bX_I^\top\bX_I)^{-1}(\bX_I\vz-N\lambda\vveps) \Big\rangle\Big|<N\lambda,\,\forall\, j\in J_\star^c\backslash I\bigg\}.
\end{align*}
Fix $\vz\in\cC_{I,\vveps}$, then $\supp(\cA[\vz])=I$ and $\sign(\cA[\vz])_I=\vveps$. The KKT identity~\eqref{eq:def_cA} gives
\begin{align*} \mathcal A[\boldsymbol z]_I=(\mathbb X_I^{\top}\mathbb X_I)^{-1}\big(\mathbb X_I^{\top}\boldsymbol z-N\lambda\operatorname{sign}(\mathcal A[\boldsymbol z]_I)\big). \end{align*}Recall that from Proposition~\ref{prop:effective-regularization-lasso}, $\vlambda^\star[\vz] = \frac{1}{N\lambda}(\vz-\bX_{J_\star^c}\cA[\vz])$. Moreover, on $\cC_{I,\vveps}$, $\cA[\vz]_{I^c}=\vzero$, hence $\vlambda^\star[\vz] = \frac{1}{N\lambda}\Big(\vz-\bX_I(\bX_I^\top\bX_I)^{-1}(\bX_I^\top\vz-N\lambda\vveps)\Big)$. Therefore, for any sufficiently small $\vh$ such that $\vz+\vh\in\cC_{I,\vveps}$, $\vlambda^\star[\vz+\vh]-\vlambda^\star[\vz]=\frac{1}{N\lambda}\Big( I_N - \bX_I(\bX_I^\top\bX_I)^{-1}\bX_I^\top\Big)\vh$.
That is,
\begin{align*} \mathbf D\boldsymbol\lambda^\star[\boldsymbol z]=\frac{I_N-\operatorname{Proj}_{\mathbb XV_I}}{N\lambda},\qquad\mbox{hence}\qquad \mathrm{Tr}\mathbf D\boldsymbol\lambda^\star[\boldsymbol z]=\frac{N-|I|}{N\lambda}.\end{align*}
In particular, if $\|\vz\|_2\leq 2\sigma_\xi\sqrt{N}$, then $\Omega_{\mathrm{act}}$ gives $|\supp(\cA[\vz])|\leq \frac{5}{16}N$. In the following, we take $\vz=\bX_{J_\star^c}\bbeta_{J_\star^c}^\star+\sqrt{t}\vG$ for $0\leq t\leq \sigma_\xi^2$. Indeed, since $C_\beta=16$, on $\Omega_{\mathrm{signal}}$, for $0\le t\le\sigma_\xi^2$ and $\|\boldsymbol G\|_2\le\frac{5}{4}\sqrt N$,
\begin{align*} \|\mathbb X_{J_\star^c}\boldsymbol\beta_{J_\star^c}^\star+\sqrt t\boldsymbol G\|_2\le2\sqrt N\|\Sigma_{J_\star^c}^{1/2}\boldsymbol\beta_{J_\star^c}^\star\|_2+\frac{5}{4}\sigma_\xi\sqrt N\le\bigg(\frac{1}{8}+\frac{5}{4}\bigg)\sigma_\xi\sqrt N<2\sigma_\xi\sqrt N. \end{align*}
The Gaussian norm bound gives $\mathbb P\{\|\boldsymbol G\|_2\le\frac{5}{4}\sqrt N\}\ge1-\exp(-\frac{N}{32})\ge\frac{1}{2}$; $N\ge32$ follows from $1\le|J_\star|$ and the small-risk assumption. Notice that
\begin{align*}
    &\bE_\vG\max_{\|\bX_{J_\star^c}^\top\vu\|_\infty\leq 1}\Big\{\Big\langle\vu,\bX_{J_\star^c}\bbeta_{J_\star^c}^\star+\sqrt{t}\vG\Big\rangle - \frac{N\lambda}{2}\|\vu\|_2^2\Big\}\bigg|_{t=\sigma_\xi^2} = P\ell_{\bbeta_{J_\star}^\star}^{(\lambda,J_\star)},\mbox{ and }\\
    &\bE_\vG\max_{\|\bX_{J_\star^c}^\top\vu\|_\infty\leq 1}\Big\{\Big\langle\vu,\bX_{J_\star^c}\bbeta_{J_\star^c}^\star+\sqrt{t}\vG\Big\rangle - \frac{N\lambda}{2}\|\vu\|_2^2\Big\}\bigg|_{t=0} = \max_{\|\bX_{J_\star^c}^\top\vu\|_\infty\leq 1}\Big\{\Big\langle\vu,\bX_{J_\star^c}\bbeta_{J_\star^c}^\star\Big\rangle - \frac{N\lambda}{2}\|\vu\|_2^2\Big\}.
\end{align*}Fix $\bX_{J_\star}$. Let $F(\vz)=\max\{ \langle\vu,\vz\rangle - \frac{N\lambda}{2}\|\vu\|_2^2 :\, \|\bX_{J_\star^c}^\top\vu\|_\infty\leq 1\}$, and  let $\vx=\bX_{J_\star^c}\bbeta_{J_\star^c}$. Let $H(t)=\bE_\vG F(\vx+\sqrt{t}\vG)$. Then by the chain rule, $H'(t) = \bE_\vG\langle\nabla F(\vx+\sqrt{t}\vG),\frac{1}{\sqrt{2t}}\vG\rangle=\frac{1}{2\sqrt{t}}\sum_{i=1}^N\bE_\vG[G_i\lambda_i^\star[\vx+\sqrt{t}\vG]]$. By Gaussian integral by parts, $\bE_\vG[G_i\lambda_i^\star[\vx+\sqrt{t}\vG]]=\sqrt{t}\bE_\vG[\D\vlambda^\star[\cdot](\vx+\sqrt{t}\vG)]_{ii}$. Consequently, $H'(t)=\frac{1}{2}\bE_\vG\Tr(\D\vlambda^\star[\cdot](\vx+\sqrt{t}\vG))$. Applying the fundamental theorem of calculus,
\begin{align*} &P\ell_{\boldsymbol\beta_{J_\star}^\star}^{(\lambda,J_\star)}-\max_{\|\mathbb X_{J_\star^c}^{\top}\boldsymbol u\|_\infty\le1}\bigg\{\langle\boldsymbol u,\mathbb X_{J_\star^c}\boldsymbol\beta_{J_\star^c}^\star\rangle-\frac{N\lambda}{2}\|\boldsymbol u\|_2^2\bigg\}  = \int_0^{\sigma_\xi^2} \frac{\mathrm{d}}{\mathrm{d}t}  \bE_\vG\max_{\|\bX_{J_\star^c}^\top\vu\|_\infty\leq 1}\Big\{\Big\langle\vu,\bX_{J_\star^c}\bbeta_{J_\star^c}^\star+\sqrt{t}\vG\Big\rangle - \frac{N\lambda}{2}\|\vu\|_2^2\Big\} \\
&= \frac{1}{2}\int_0^{\sigma_\xi^2} \bE_\vG\Tr\D\vlambda^\star[\cdot]\Big(\bX_{J_\star^c}\bbeta_{J_\star^c}^\star+\sqrt{t}\vG\Big)\mathrm{d}t =\frac{1}{2}\int_0^{\sigma_\xi^2}\mathbb E_{\boldsymbol G}\mathrm{Tr}\mathbf D\boldsymbol\lambda^\star[\mathbb X_{J_\star^c}\boldsymbol\beta_{J_\star^c}^\star+\sqrt t\boldsymbol G]\,\mathrm dt\ge\frac{11\sigma_\xi^2}{64\lambda}. \end{align*}
Finally, using~\eqref{eq:curvature-identity}. Since $r\le\sigma_\xi$,
\begin{align*} P\mathcal L_{\boldsymbol\beta_{J_\star}}^{(\lambda,J_\star)}\ge\bigg(1-\frac{\sigma_\xi}{\sqrt{\sigma_\xi^2+r^2}}\bigg)\frac{11\sigma_\xi^2}{64\lambda}\ge\frac{11r^2\sigma_\xi^2}{128\lambda(\sigma_\xi^2+r^2)}\ge\frac{11r^2}{256\lambda}\ge c_{\mathrm{LB}}\frac{r^2}{\lambda}. \end{align*}
Since $r$ and $\boldsymbol\beta_{J_\star}$ were arbitrary and $\Omega$ is independent of both, this proves the proposition.
\end{proof}

The proof of the one-sided RIP, as we mentioned earlier, is almost identical to that in the isotropic case.
\begin{Propositionextend}{prop:RIP}
{Assume $\Sigma=\operatorname{diag}(\sigma_1,\ldots,\sigma_p)$ with $\sigma_j>0$, and $1\le |J_\star|\le\frac{N}{100}$. For each fixed $\lambda>0$, let
\begin{align*} T=B_1^{J_\star}(\boldsymbol\beta_{J_\star}^\star;\rho)\cap\big(\boldsymbol\beta_{J_\star}^\star+\Sigma_{J_\star}^{-1/2}B_2^{J_\star}(\boldsymbol0;r)\big). \end{align*}}
Assume that $r\le\rho$ and $r\ge C_\beta\|{\Sigma_{J_\star^c}^{1/2}}\boldsymbol\beta_{J_\star^c}^\star\|_2$, where $C_\beta\ge1$. Then
\begin{align*} \bP\bigg\{\sup_{\boldsymbol\beta_{J_\star}\in T}\big|(P-P_N)\mathcal L_{\boldsymbol\beta_{J_\star}}^{(\lambda,J_\star)}\big|\le C\frac{\sqrt{\sigma_\xi^2+r^2}}{\lambda\sqrt N}r\sqrt{|J_\star|}\bigg\}\geq 1-C\exp(-\frac{|J_\star|}{C}). \end{align*}

\end{Propositionextend}

\begin{proof}[Proof of Proposition~\ref{prop:RIP}.]
Fix $\lambda > 0$. {Write $P\cdot=\mathbb E[\,\cdot\mid\mathbb X_{J_\star^c}]$. Conditional estimates below are taken with respect to the estimation design and the noise; probabilities involving $\Omega_{\mathrm{res}}$ are joint over the free design and the noise.} Let $\vzeta=\bX_{J_\star^c}\bbeta_{J_\star^c}^\star+\bxi$. We investigate the stochastic processes $(P_N\cM_{\bbeta_{J_\star}})_{\bbeta_{J_\star}\in T}$ and $(P_N\cQ_{\bbeta_{J_\star}})_{\bbeta_{J_\star}\in T}$, where
\begin{align*} &P_N\cM_{\bbeta_{J_\star}}:=\left\langle-\boldsymbol\lambda^\star[\vzeta],\mathbb X_{J_\star}(\boldsymbol\beta_{J_\star}-\boldsymbol\beta_{J_\star}^\star)\right\rangle,\quad \mbox{ and }\quad P_N\cQ_{\bbeta_{J_\star}}:=P_N\mathcal L_{\boldsymbol\beta_{J_\star}}^{(\lambda,J_\star)}-P_N\cM_{\bbeta_{J_\star}}.\tag{\ref{eq:modified-process}} \end{align*}
Independence and $\mathbb E\mathbb X_{J_\star}=\vzero$ imply $P(P_N\cM_{\bbeta_{J_\star}})=0$ and $P(P_N\cQ_{\bbeta_{J_\star}})=P\mathcal L_{\boldsymbol\beta_{J_\star}}^{(\lambda,J_\star)}$. We abbreviate $P(P_N\cM_{\bbeta_{J_\star}})$ by $P\cM_{\bbeta_{J_\star}}$ and similarly for $P(P_N\cQ_{\bbeta_{J_\star}})$. Consequently,
\begin{align*} \tag{\ref{eq:original-modified-decomposition}} &(P_N-P)\mathcal L_{\boldsymbol\beta_{J_\star}}^{(\lambda,J_\star)}=P_N\cM_{\bbeta_{J_\star}}+(P_N-P)\cQ_{\bbeta_{J_\star}}. \end{align*}
Taking absolute values gives
\begin{align*} &\sup_{\boldsymbol\beta_{J_\star}\in T}\left|(P_N-P)\mathcal L_{\boldsymbol\beta_{J_\star}}^{(\lambda,J_\star)}\right|\le\sup_{\boldsymbol\beta_{J_\star}\in T}|P_N\cM_{\bbeta_{J_\star}}|+\sup_{\boldsymbol\beta_{J_\star}\in T}|(P_N-P)\cQ_{\bbeta_{J_\star}}|. \end{align*}
Let
\begin{align*} &\Omega_{\mathrm{res}}=\left\{\|\vzeta\|_2\le\sqrt{3}\sqrt{\sigma_\xi^2+\|{\Sigma_{J_\star^c}^{1/2}}\boldsymbol\beta_{J_\star^c}^\star\|_2^2}\sqrt N\right\}. \end{align*}
Then $\mathbb P(\Omega_{\mathrm{res}}^c)\le\exp(-\frac{|J_\star|}{C})$. {Conditioning first on $(\mathbb X_{J_\star^c},\boldsymbol\xi)$, the process has Gaussian increments with metric $\|\boldsymbol\lambda^\star[\boldsymbol\zeta]\|_2\|\Sigma_{J_\star}^{1/2}(\boldsymbol\beta_{J_\star}-\boldsymbol\beta'_{J_\star})\|_2$. Its index set is carried by whitening into a subset of $B_2^{J_\star}(\boldsymbol0;r)$.} By generic chaining and the independence between $\vlambda^\star[\vzeta]$ and $\bX_{J_\star}$, see, for instance, \cite[Lemma 2]{shangFeatureSpaceDecomposition2026}, together with Proposition~\ref{prop:effective-regularization-lasso}, item~\emph{\ref{item:effective_lambda_star}} for the upper bound of $\|\vlambda^\star[\vzeta]\|_2$, we have
\begin{align*} \tag{\ref{eq:multiplier_process_upper}} &\bP\bigg\{\sup_{\boldsymbol\beta_{J_\star}\in T}|P_N\cM_{\bbeta_{J_\star}}|\le C\min\left\{\mathrm{diam}\!\left((\mathbb X_{J_\star^c}B_1^{J_\star^c})^\circ\right),\frac{\sqrt{3}\sqrt{\sigma_\xi^2+\|{\Sigma_{J_\star^c}^{1/2}}\boldsymbol\beta_{J_\star^c}^\star\|_2^2}}{\lambda\sqrt N}\right\}r\sqrt{|J_\star|}\bigg\}\geq 1-{C}\exp(-\frac{|J_\star|}{C}). \end{align*}
Let $\Omega_{\mathrm{lin}}$ denote the event in \eqref{eq:multiplier_process_upper}.

{For the quadratic component, condition only on $\mathbb X_{J_\star^c}$. }Let $\Omega_{\mathrm{op}}=\{\|\mathbb X_{J_\star}{\Sigma_{J_\star}^{-1/2}}\|_{\mathrm{op}}\le3\sqrt N\}$. Then $\mathbb P(\Omega_{\mathrm{op}}^c)\le\exp(-\frac{|J_\star|}{C})$, and we write $P^{\Omega_{\mathrm{op}}}=\mathbb E[\,\cdot\mid\mathbb X_{J_\star^c},\Omega_{\mathrm{op}}]$. Let
\begin{align*} \tag{\ref{eq:dual-residual-representation}} F:t\in[0,1] \mapsto \max_{\boldsymbol a\in(\mathbb X_{J_\star^c}B_1^{J_\star^c})^\circ}\left\{\left\langle\boldsymbol a,\vzeta-t\mathbb X_{J_\star}(\boldsymbol\beta_{J_\star}-\boldsymbol\beta_{J_\star}^\star)\right\rangle-\frac{N\lambda}{2}\|\boldsymbol a\|_2^2\right\}. \end{align*}
By \eqref{eq:dual_problem_effective_loss}, $P_N\mathcal L_{\boldsymbol\beta_{J_\star}}^{(\lambda,J_\star)}=F(1)-F(0)$, and the maximizer of $F(t)$ is unique and equal to $\boldsymbol\lambda^\star[\vzeta-t\mathbb X_{J_\star}(\boldsymbol\beta_{J_\star}-\boldsymbol\beta_{J_\star}^\star)]$. Thus, by Danskin's theorem, for almost every $t\in[0,1]$, $\frac{\mathrm{d}}{\mathrm{d}t}F(t)=-\langle\vlambda^\star[\vzeta-t\bX_{J_\star}(\bbeta_{J_\star}-\bbeta_{J_\star}^\star)],\bX_{J_\star}(\bbeta_{J_\star}-\bbeta_{J_\star}^\star)\rangle$.
The fundamental theorem of calculus and \eqref{eq:modified-process}  give $P_N\cL_{\bbeta_{J_\star}}^{(\lambda,J_\star)}=\int_0^1 \frac{\mathrm{d}}{\mathrm{d}t}F(t)\mathrm{d}t = P_N\cQ_{\bbeta_{J_\star}}+P_N\cM_{\bbeta_{J_\star}} = P_N\cQ_{\bbeta_{J_\star}} + \langle-\vlambda^\star[\vzeta],\bX_{J_\star}(\bbeta_{J_\star}-\bbeta_{J_\star}^\star)\rangle$, and  therefore
\begin{align*} &P_N\cQ_{\bbeta_{J_\star}}=\int_0^1\Big\langle\boldsymbol\lambda^\star[\vzeta]-\boldsymbol\lambda^\star[\vzeta-t\mathbb X_{J_\star}(\boldsymbol\beta_{J_\star}-\boldsymbol\beta_{J_\star}^\star)],\mathbb X_{J_\star}(\boldsymbol\beta_{J_\star}-\boldsymbol\beta_{J_\star}^\star)\Big\rangle\,\mathrm dt.\tag{\ref{eq:original-smooth-convexity}} \end{align*}
Recall that by Proposition~\ref{prop:effective-regularization-lasso}, $\vlambda^\star[\bullet] = \Proj_{(\bX_{J_\star^c}B_1^{J_\star^c})^\circ}(\frac{1}{N\lambda}\bullet)$. For every $t\in(0,1]$, the nonexpansiveness of the metric projection, the Cauchy--Schwarz inequality, and the Lipschitzness of $\vlambda^\star[\cdot]$ in Proposition~\ref{prop:effective-regularization-lasso} give
\begin{align*} &0\le\frac{N\lambda}{t}\|\boldsymbol\lambda^\star[\vzeta]-\boldsymbol\lambda^\star[\vzeta-t\mathbb X_{J_\star}(\boldsymbol\beta_{J_\star}-\boldsymbol\beta_{J_\star}^\star)]\|_2^2\\ &\le\Big\langle\boldsymbol\lambda^\star[\vzeta]-\boldsymbol\lambda^\star[\vzeta-t\mathbb X_{J_\star}(\boldsymbol\beta_{J_\star}-\boldsymbol\beta_{J_\star}^\star)],\mathbb X_{J_\star}(\boldsymbol\beta_{J_\star}-\boldsymbol\beta_{J_\star}^\star)\Big\rangle\le\frac{t}{N\lambda}\|\mathbb X_{J_\star}(\boldsymbol\beta_{J_\star}-\boldsymbol\beta_{J_\star}^\star)\|_2^2. \end{align*}
Integration in \eqref{eq:original-smooth-convexity} yields, almost surely, for every $\boldsymbol\beta_{J_\star}\in T$,
\begin{align*} &0\le P_N\cQ_{\bbeta_{J_\star}}\le\frac{\|\mathbb X_{J_\star}(\boldsymbol\beta_{J_\star}-\boldsymbol\beta_{J_\star}^\star)\|_2^2}{2N\lambda}. \end{align*}
This nonnegativity justifies our terminology ``quadratic component'' for $P_N\cQ_{\bbeta_{J_\star}}$ --- it is an analogue of the nonnegative stochastic process known as the quadratic process; however, the bound $P_N\cQ_{\bbeta_{J_\star}} \le \frac{\|\mathbb X_{J_\star}(\boldsymbol\beta_{J_\star} - \boldsymbol\beta_{J_\star}^\star)\|_2^2}{2N\lambda}$ is not sharp enough. We next examine, similarly, the uniform boundedness of the sub-Gaussian norms of the increments of $((P-P_N)\cQ_{\bbeta_{J_\star}})_{\bbeta_{J_\star} \in T}$. Our tool is the (restricted) Gaussian Lipschitz concentration inequality. To this end, we first examine the gradient of the increments. Recall from \eqref{eq:modified-process} that $P_N\cQ_{\bbeta_{J_\star}} = P_N\cL_{\bbeta_{J_\star}}^{(\lambda,J_\star)} - P_N\cM_{\bbeta_{J_\star}}$, then we compute the gradients of $P_N\cL_{\bbeta_{J_\star}}^{(\lambda,J_\star)}$ and that of $P_N\cM_{\bbeta_{J_\star}}^{(\lambda,J_\star)}$ to compute the gradient of $P_N\cQ_{\bbeta_{J_\star}}$.
\begin{enumerate}
    \item For any $\bbeta_{J_\star},\bbeta_{J_\star}'\in T$ and almost every $\boldsymbol\xi$, Danskin's theorem yields $\nabla_\bxi P_N\cL_{\bbeta_{J_\star}}^{(\lambda,J_\star)} = \vlambda^\star[\vzeta - \bX_{J_\star}(\bbeta_{J_\star}-\bbeta_{J_\star}^\star)]-\vlambda^\star[{\vzeta}]$, and the chain rule gives $\nabla_\bxi P_N\cM_{\bbeta_{J_\star}} = -(\D\vlambda^\star[\cdot](\vzeta))\bX_{J_\star}(\bbeta_{J_\star}-\bbeta_{J_\star}^\star)$, hence
    \begin{align*} \tag{\ref{eq:noise-gradient}} \begin{aligned} &\nabla_{\boldsymbol\xi}(P_N\cQ_{\bbeta_{J_\star}}-P_N\cQ_{\bbeta'_{J_\star}})\\ &=\boldsymbol\lambda^\star[\vzeta-\mathbb X_{J_\star}(\boldsymbol\beta_{J_\star}-\boldsymbol\beta_{J_\star}^\star)]-\boldsymbol\lambda^\star[\vzeta-\mathbb X_{J_\star}(\boldsymbol\beta'_{J_\star}-\boldsymbol\beta_{J_\star}^\star)] +\big(\D\boldsymbol\lambda^\star[\cdot](\vzeta)\big)\mathbb X_{J_\star}(\boldsymbol\beta_{J_\star}-\boldsymbol\beta'_{J_\star}). \end{aligned} \end{align*}
    On $\Omega_{\mathrm{op}}$,
    \begin{align*} &\|\boldsymbol\lambda^\star[\vzeta-\mathbb X_{J_\star}(\boldsymbol\beta_{J_\star}-\boldsymbol\beta_{J_\star}^\star)]-\boldsymbol\lambda^\star[\vzeta-\mathbb X_{J_\star}(\boldsymbol\beta'_{J_\star}-\boldsymbol\beta_{J_\star}^\star)]\|_2\le\frac{\|\mathbb X_{J_\star}(\boldsymbol\beta_{J_\star}-\boldsymbol\beta'_{J_\star})\|_2}{N\lambda}\le\frac{3}{\lambda\sqrt N}\|{\Sigma_{J_\star}^{1/2}(}\boldsymbol\beta_{J_\star}-\boldsymbol\beta'_{J_\star}{)}\|_2, \end{align*}
    and $\|\big(\D\boldsymbol\lambda^\star[\cdot](\vzeta)\big)\mathbb X_{J_\star}(\boldsymbol\beta_{J_\star}-\boldsymbol\beta'_{J_\star})\|_2\le\frac{3}{\lambda\sqrt N}\|{\Sigma_{J_\star}^{1/2}(}\boldsymbol\beta_{J_\star}-\boldsymbol\beta'_{J_\star}{)}\|_2.$ Thus \eqref{eq:noise-gradient} implies
    \begin{align*} &\|\nabla_{\boldsymbol\xi}(P_N\cQ_{\bbeta_{J_\star}}-P_N\cQ_{\bbeta'_{J_\star}})\|_2\le\frac{6}{\lambda\sqrt N}\|{\Sigma_{J_\star}^{1/2}(}\boldsymbol\beta_{J_\star}-\boldsymbol\beta'_{J_\star}{)}\|_2. \end{align*}

    \item {Take derivatives with respect to $\mathbb X_{J_\star}\Sigma_{J_\star}^{-1/2}$, whose entries are independent standard Gaussian variables. }Similarly, 
\begin{align*} \nabla_{\bX_{J_\star}{\Sigma_{J_\star}^{-1/2}}}P_N\cL_{\bbeta_{J_\star}}^{(\lambda,J_\star)} = -\vlambda^\star[\vzeta-\bX_{J_\star}(\bbeta_{J_\star}-\bbeta_{J_\star}^\star)]\otimes{\Sigma_{J_\star}^{1/2}}(\bbeta_{J_\star}-\bbeta_{J_\star}^\star),\quad\mbox{and}\\ \nabla_{\bX_{J_\star}{\Sigma_{J_\star}^{-1/2}}}P_N\cM_{\bbeta_{J_\star}} = -\vlambda^\star[\vzeta]\otimes{\Sigma_{J_\star}^{1/2}}(\bbeta_{J_\star}-\bbeta_{J_\star}^\star), \end{align*}
hence
    \begin{align*} \tag{\ref{eq:design-gradient}} \begin{aligned} &\nabla_{\mathbb X_{J_\star}{\Sigma_{J_\star}^{-1/2}}}(P_N\cQ_{\bbeta_{J_\star}}-P_N\cQ_{\bbeta'_{J_\star}})=\Big(\boldsymbol\lambda^\star[\vzeta]-\boldsymbol\lambda^\star[\vzeta-\mathbb X_{J_\star}(\boldsymbol\beta_{J_\star}-\boldsymbol\beta_{J_\star}^\star)]\Big)\otimes{\Sigma_{J_\star}^{1/2}}(\boldsymbol\beta_{J_\star}-\boldsymbol\beta'_{J_\star})\\ &+\Big(\boldsymbol\lambda^\star[\vzeta-\mathbb X_{J_\star}(\boldsymbol\beta'_{J_\star}-\boldsymbol\beta_{J_\star}^\star)]-\boldsymbol\lambda^\star[\vzeta-\mathbb X_{J_\star}(\boldsymbol\beta_{J_\star}-\boldsymbol\beta_{J_\star}^\star)]\Big)\otimes{\Sigma_{J_\star}^{1/2}}(\boldsymbol\beta'_{J_\star}-\boldsymbol\beta_{J_\star}^\star). \end{aligned} \end{align*}
    Since $\|{\Sigma_{J_\star}^{1/2}(}\boldsymbol\beta_{J_\star}-\boldsymbol\beta_{J_\star}^\star{)}\|_2,\|{\Sigma_{J_\star}^{1/2}(}\boldsymbol\beta'_{J_\star}-\boldsymbol\beta_{J_\star}^\star{)}\|_2\le r$, the Frobenius norms of the first and second terms on the right-hand side of \eqref{eq:design-gradient} are bounded, respectively, by $\frac{\|\mathbb X_{J_\star}(\boldsymbol\beta_{J_\star}-\boldsymbol\beta_{J_\star}^\star)\|_2}{N\lambda}\|{\Sigma_{J_\star}^{1/2}(}\boldsymbol\beta_{J_\star}-\boldsymbol\beta'_{J_\star}{)}\|_2\le\frac{3r}{\lambda\sqrt N}\|{\Sigma_{J_\star}^{1/2}(}\boldsymbol\beta_{J_\star}-\boldsymbol\beta'_{J_\star}{)}\|_2,$ and $\frac{\|\mathbb X_{J_\star}(\boldsymbol\beta_{J_\star}-\boldsymbol\beta'_{J_\star})\|_2}{N\lambda}\|{\Sigma_{J_\star}^{1/2}(}\boldsymbol\beta'_{J_\star}-\boldsymbol\beta_{J_\star}^\star{)}\|_2\le\frac{3r}{\lambda\sqrt N}\|{\Sigma_{J_\star}^{1/2}(}\boldsymbol\beta_{J_\star}-\boldsymbol\beta'_{J_\star}{)}\|_2.$ Therefore,
    \begin{align*} &\|\nabla_{\mathbb X_{J_\star}{\Sigma_{J_\star}^{-1/2}}}(P_N\cQ_{\bbeta_{J_\star}}-P_N\cQ_{\bbeta'_{J_\star}})\|_{\mathrm F}\le\frac{6r}{\lambda\sqrt N}\|{\Sigma_{J_\star}^{1/2}(}\boldsymbol\beta_{J_\star}-\boldsymbol\beta'_{J_\star}{)}\|_2. \end{align*}
\end{enumerate}
Notice that the probability distribution of $\bX_{J_\star}{\Sigma_{J_\star}^{-1/2}}$ condition on $\Omega_{\op}$ is the normalized Gaussian distribution restricted on the closed convex set $\{\vA:\|\vA\|_{\op}\leq 3\sqrt{N}\}$, hence by Lemma~\ref{lem:restricted_Gaussian_Lipschitz},
\begin{align*} &\left\|(P_N\cQ_{\bbeta_{J_\star}}-P_N\cQ_{\bbeta'_{J_\star}})-P^{\Omega_{\mathrm{op}}}(P_N\cQ_{\bbeta_{J_\star}}-P_N\cQ_{\bbeta'_{J_\star}})\right\|_{\psi_2(P^{\Omega_{\op}})}\le C\frac{\sqrt{\sigma_\xi^2+r^2}}{\lambda\sqrt N}\|{\Sigma_{J_\star}^{1/2}(}\boldsymbol\beta_{J_\star}-\boldsymbol\beta'_{J_\star}{)}\|_2. \end{align*}
Since $P_N\cQ_{\boldsymbol\beta_{J_\star}^\star}=0$, {and $\Sigma_{J_\star}^{1/2}(T-\boldsymbol\beta_{J_\star}^\star)\subseteq B_2^{J_\star}(\boldsymbol0;r)$, }generic chaining gives, with $P^{\Omega_{\mathrm{op}}}$-probability at least $1-2\exp(-\frac{|J_\star|}{C})$,
\begin{align*} \tag{\ref{eq:conditional-chaining}} &\sup_{\boldsymbol\beta_{J_\star}\in T}\left|P_N\cQ_{\bbeta_{J_\star}}-P^{\Omega_{\mathrm{op}}}(P_N\cQ_{\bbeta_{J_\star}})\right|\le C\frac{\sqrt{\sigma_\xi^2+r^2}}{\lambda\sqrt N}r\sqrt{|J_\star|}. \end{align*}
We are left with comparing $P^{\Omega_{\op}}(P_N\cQ_{\bbeta_{J_\star}})$ with $P\cQ_{\bbeta_{J_\star}}$. The law of total expectation gives
\begin{align*} \forall\bbeta_{J_\star}\in T,\quad & P\cQ_{\bbeta_{J_\star}} = \bP(\Omega_{\op})P^{\Omega_{\op}}(P_N\cQ_{\bbeta_{J_\star}}) + P\big(P_N\cQ_{\bbeta_{J_\star}}\1_{\Omega_{\op}^c}\big), \end{align*}hence $|P^{\Omega_{\op}}(P_N\cQ_{\bbeta_{J_\star}}) - P\cQ_{\bbeta_{J_\star}}| \leq \bP(\Omega_{\op}^c)P^{\Omega_{\op}}(P_N\cQ_{\bbeta_{J_\star}}) + P(P_N\cQ_{\bbeta_{J_\star}}\1_{\Omega_{\op}^c})$, and from $0\le P_N\cQ_{\bbeta_{J_\star}}\le \frac{r^2\|\mathbb X_{J_\star}{\Sigma_{J_\star}^{-1/2}}\|_{\mathrm{op}}^2}{2N\lambda}$, we obtain
\begin{align*} &\sup_{\boldsymbol\beta_{J_\star}\in T}\left|P^{\Omega_{\mathrm{op}}}(P_N\cQ_{\bbeta_{J_\star}})-P\cQ_{\bbeta_{J_\star}}\right|\le\frac{r^2}{2N\lambda}\mathbb E\!\left[(\|\mathbb X_{J_\star}{\Sigma_{J_\star}^{-1/2}}\|_{\mathrm{op}}^2+9N)\mathbf 1_{\Omega_{\mathrm{op}}^c}\right]. \end{align*}
For every $t\ge3\sqrt N$, $\mathbb P\{\|\mathbb X_{J_\star}{\Sigma_{J_\star}^{-1/2}}\|_{\mathrm{op}}>t\}\le\exp(-\frac{(t-\sqrt N-\sqrt{|J_\star|})^2}{2}).$ Fubini's theorem {giv}es
\begin{align*} &\mathbb E\!\left[(\|\mathbb X_{J_\star}{\Sigma_{J_\star}^{-1/2}}\|_{\mathrm{op}}^2+9N)\mathbf 1_{\Omega_{\mathrm{op}}^c}\right]=18N\mathbb P\{\|\mathbb X_{J_\star}{\Sigma_{J_\star}^{-1/2}}\|_{\mathrm{op}}>3\sqrt N\}+\int_{3\sqrt N}^{\infty}2t\,\mathbb P\{\|\mathbb X_{J_\star}{\Sigma_{J_\star}^{-1/2}}\|_{\mathrm{op}}>t\}\,\mathrm dt\\ &\le18N\exp\left(-\frac{(2\sqrt N-\sqrt{|J_\star|})^2}{2}\right)+\int_{2\sqrt N-\sqrt{|J_\star|}}^\infty2(u+\sqrt N+\sqrt{|J_\star|})e^{-\frac{u^2}{2}}\,\mathrm du\\ &\le\left(18N+2+\frac{2(\sqrt N+\sqrt{|J_\star|})}{2\sqrt N-\sqrt{|J_\star|}}\right)\exp\left(-\frac{(2\sqrt N-\sqrt{|J_\star|})^2}{2}\right)\le22N\exp\left(-\frac{361N}{200}\right). \end{align*}
Therefore, $\sup_{\boldsymbol\beta_{J_\star}\in T}|P^{\Omega_{\mathrm{op}}}(P_N\cQ_{\bbeta_{J_\star}})-P\cQ_{\bbeta_{J_\star}}|$ is at most $\frac{11r^2}{\lambda}\exp\left(-\frac{361N}{200}\right)$, which is absorbed into $C\frac{\sqrt{\sigma_\xi^2+r^2}}{\lambda\sqrt N}r\sqrt{|J_\star|}$. Together with \eqref{eq:conditional-chaining}, this gives an event $\Omega_{\mathrm{Br}}$ such that $\mathbb P(\Omega_{\mathrm{op}}\cap\Omega_{\mathrm{Br}}^c\mid\mathbb X_{J_\star^c})\le C\exp(-\frac{|J_\star|}{C})$ and, on $\Omega_{\mathrm{op}}\cap\Omega_{\mathrm{Br}}$,
\begin{align*} &\sup_{\boldsymbol\beta_{J_\star}\in T}|(P_N-P)\cQ_{\bbeta_{J_\star}}|\le C\frac{\sqrt{\sigma_\xi^2+r^2}}{\lambda\sqrt N}r\sqrt{|J_\star|}. \end{align*}
Finally, let
\begin{align*} &\Omega=\Omega_{\mathrm{op}}\cap\Omega_{\mathrm{res}}\cap\Omega_{\mathrm{lin}}\cap\Omega_{\mathrm{Br}},\quad\mbox{ then }\quad \mathbb P(\Omega)\ge1-C\exp(-\frac{|J_\star|}{C}). \end{align*}
On $\Omega$, substitution of the bounds for $P_N\cM_{\bbeta_{J_\star}}$ and $(P_N-P)\cQ_{\bbeta_{J_\star}}$ into \eqref{eq:original-modified-decomposition}, together with $r\ge C_\beta\|{\Sigma_{J_\star^c}^{1/2}}\boldsymbol\beta_{J_\star^c}^\star\|_2$, proves the proposition.
\end{proof}

For the free part of the estimator $\hat\bbeta$, we replace Proposition~\ref{prop:lambda_large_shift} with the following extended version.
\begin{Propositionextend}{prop:lambda_large_shift}
Assume $|J_\star^c|>N$, $r_\star\ge0$, {$\lambda\ge\frac{2\sigma_\xi z_N^{\mathrm{all}}}{\sqrt N}$, and $r_\star^2+\mathrm{Var}_{J_\star^c}\le2^{-16}\sigma_\xi^2$}. Let $\boldsymbol\zeta=\mathbb X_{J_\star^c}\boldsymbol\beta_{J_\star^c}^\star+\mathbb X_{J_\star}(\boldsymbol\beta_{J_\star}^\star-\widehat{\boldsymbol\beta}_{J_\star})$, and retain $\Omega_{\mathrm{esti}}=\{\|\boldsymbol\zeta\|_2\le\sqrt N r_\star\}$. With probability at least {$1-\mathbb P(\Omega_{\mathrm{esti}}^c)-12\exp(-\frac{d_q}{4})$},
{\begin{align*} \|\Sigma_{J_\star^c}^{1/2}\mathcal A[\boldsymbol y-\mathbb X_{J_\star}\widehat{\boldsymbol\beta}_{J_\star}]\|_2^2\le 2^{12}(r_\star^2+\mathrm{Var}_{J_\star^c}). \end{align*}}
\end{Propositionextend}

\begin{proof}
Let $t=64\sqrt{r_\star^2 + \Var_{J_\star^c}}$. The overall strategy of the proof is divided into the following two steps. 
\begin{enumerate}
    \item Prove that with high probability,
    \begin{align*}
        \forall \bbeta_{J_\star^c}\in V_{J_\star^c}\mbox{ such that }\|\Sigma_{J_\star^c}^{1/2}\bbeta_{J_\star^c}\|_2\geq t,\qquad \frac{1}{2N\lambda}\|\bxi - \bX_{J_\star^c}\bbeta_{J_\star^c}\|_2^2 + \|\bbeta_{J_\star^c}\|_1 - \frac{1}{2N\lambda}\|\bxi\|_2^2 > \frac{r_\star^2+\Var_{J_\star^c}}{\lambda}.
    \end{align*}
    \item Prove that with high probability,
    \begin{align*}
        \frac{1}{2N\lambda}\|\bxi - \bX_{J_\star^c}\cA[\bxi+\vzeta]\|_2^2 + \|\cA[\bxi+\vzeta]\|_1 - \frac{1}{2N\lambda}\|\bxi\|_2^2 \leq  \frac{r_\star^2}{2\lambda}.
    \end{align*}
\end{enumerate}
Combining the above, we obtain that $\|\Sigma_{J_\star^c}^{1/2}\cA[\bxi+\vzeta]\|_2<t=64\sqrt{r_\star^2+\Var_{J_\star^c}}$.

Item~\emph{2.} is easy to prove; it only requires us to use the KKT condition of $\cA[\bxi+\vzeta]\in\argmin\{ \frac{1}{2N\lambda}\|\bxi+\vzeta-\bX_{J_\star^c}\bbeta_{J_\star^c}\|_2^2 + \|\bbeta_{J_\star^c}\|_1 :\, \bbeta_{J_\star^c}\in V_{J_\star^c}\}$. Indeed, this gives $\frac{1}{N\lambda}(\bX_{J_\star^c}^\top(\bxi+\vzeta - \bX_{J_\star^c}\cA[\bxi+\vzeta]))\in\partial^-\|\cdot\|_1(\cA[\bxi+\vzeta])$. By the definition of sub-differential,
\begin{align*}
    \frac{1}{N\lambda}\Big\langle \bX_{J_\star^c}^\top(\bxi+\vzeta - \bX_{J_\star^c}\cA[\bxi+\vzeta]), \cA[\bxi+\vzeta] \Big\rangle \geq \|\cA[\bxi+\vzeta]\|_1.
\end{align*}Therefore, on $\Omega_{\mathrm{esti}}$,
\begin{align*}
    &\frac{1}{2N\lambda}\|\bxi - \bX_{J_\star^c}\cA[\bxi+\vzeta]\|_2^2 + \|\cA[\bxi+\vzeta]\|_1 - \frac{1}{2N\lambda}\|\bxi\|_2^2 \\
    &\leq \frac{1}{2N\lambda}\|\bxi - \bX_{J_\star^c}\cA[\bxi+\vzeta]\|_2^2  - \frac{1}{2N\lambda}\|\bxi\|_2^2 + \frac{1}{N\lambda}\Big\langle \bxi+\vzeta - \bX_{J_\star^c}\cA[\bxi+\vzeta], \bX_{J_\star^c}\cA[\bxi+\vzeta] \Big\rangle\\
    &=\frac{1}{N\lambda}\Big\langle \vzeta,\bX_{J_\star^c}\cA[\bxi+\vzeta] \Big\rangle - \frac{1}{2N\lambda}\|\bX_{J_\star^c}\cA[\bxi+\vzeta]\|_2^2 \leq \frac{1}{2N\lambda}\|\vzeta\|_2^2 \leq \frac{r_\star^2}{2\lambda}.
\end{align*}We now prove item~\emph{1.}. By homogeneity argument, we only need to prove $\bP^{\otimes N}(\Omega_{\mathrm{cost}})\geq 1-12\exp(-\frac{1}{4}d_q)$, where
\begin{align*} \Omega_{\mathrm{cost}}=\bigg\{\inf_{\|\Sigma_{J_\star^c}^{1/2}\boldsymbol\beta_{J_\star^c}\|_2=t}\bigg[\frac12\|\boldsymbol\xi-\mathbb X_{J_\star^c}\boldsymbol\beta_{J_\star^c}\|_2^2+N\lambda\|\boldsymbol\beta_{J_\star^c}\|_1-\frac12\|\boldsymbol\xi\|_2^2\bigg]>N(r_\star^2+\mathrm{Var}_{J_\star^c})\bigg\}. \end{align*}
We again use CGMT and follow the proof of Proposition~\ref{prop:null_signal_upper_bound}. Thus we replace Lemma~\ref{lemma:null_signal_events}, Lemma~\ref{lemma:null_signal_scalarization} and Lemma~\ref{lemma:null_signal_cost_separation}. Let $\vg, \vh$ be random vectors in the auxiliary probability space. We replace the random events in Lemma~\ref{lemma:null_signal_events} with
\begin{align*} &\Omega_\xi=\big\{\big|\|\boldsymbol\xi\|_2-\sigma_\xi\sqrt N\big|\le\sigma_\xi\sqrt{d_q}\big\},\qquad\mbox{and}\qquad \Omega_{\boldsymbol g}^{\mathrm{aux}}=\bigg\{\frac N2\le\|\boldsymbol g\|_2^2\le2N,\ \mbox{and}\ |\langle\boldsymbol\xi,\boldsymbol g\rangle|\le\|\boldsymbol\xi\|_2\sqrt{d_q}\bigg\},\\ &\Omega_{z_q,2}^{\mathrm{aux}}=\bigg\{\bigg[\sum_{j\in J_\star^c}\bigg(|h_j|-\frac{z_q}{\sqrt{\sigma_j}}\bigg)_+^2\bigg]^{\frac12}\le2\sqrt{d_q}\bigg\},\qquad\mbox{and}\qquad \Omega^{\mathrm{aux}}=\Omega_{\boldsymbol g}^{\mathrm{aux}}\cap\Omega_{z_q,2}^{\mathrm{aux}}. \end{align*}Using the same proof method as in Lemma~\ref{lemma:null_signal_events}, the following conclusion readily follows:
\begin{align*} {\mathbb P(\Omega_\xi^c)\le2e^{-\frac{d_q}{4}},\qquad\mbox{and}\qquad \sup_{\boldsymbol\xi\in\Omega_\xi}\mathbb P_{\mathrm{aux}}((\Omega^{\mathrm{aux}})^c\mid\boldsymbol\xi)\le2e^{-\frac N{16}}+3e^{-\frac{d_q}{2}}\le5e^{-\frac{d_q}{2}}.} \tag{\ref{eq:omega_probability}} \end{align*}
Here $d_q\le\frac{N\mathrm{Var}_{J_\star^c}}{\sigma_\xi^2}\le2^{-16}N$. In particular, $\mathbb P(\Omega_\xi\cap\Omega^{\mathrm{aux}})\ge1-7e^{-\frac{d_q}{4}}$.

Here, we do not need the full version of Lemma~\ref{lemma:null_signal_scalarization}, but only the conclusion of CGMT, namely, for $K,H>0$ and a nonempty closed set $I\subset[0,K]$, define
\begin{align*}
\Phi_I^{K,H}&=\min_{\substack{\boldsymbol\beta_{J_\star^c}\in\mathbb R^{J_\star^c}\\\|\Sigma_{J_\star^c}^{1/2}\boldsymbol\beta_{J_\star^c}\|_2\in I}}\max_{\substack{\boldsymbol u\in\mathbb R^N\\\|\boldsymbol u\|_2\leq H}}\bigg\{\langle\boldsymbol u,\boldsymbol\xi-\mathbb X_{J_\star^c}\boldsymbol\beta_{J_\star^c}\rangle-\frac12\|\boldsymbol u\|_2^2+N\lambda\|\boldsymbol\beta_{J_\star^c}\|_1-\frac12\|\boldsymbol\xi\|_2^2\bigg\},\\
\phi_I^{K,H}&=\min_{\substack{\boldsymbol\beta_{J_\star^c}\in\mathbb R^{J_\star^c}\\\|\Sigma_{J_\star^c}^{1/2}\boldsymbol\beta_{J_\star^c}\|_2\in I}}\max_{\substack{\boldsymbol u\in\mathbb R^N\\\|\boldsymbol u\|_2\leq H}}\bigg\{\|\Sigma_{J_\star^c}^{1/2}\boldsymbol\beta_{J_\star^c}\|_2\langle\boldsymbol g,\boldsymbol u\rangle+\|\boldsymbol u\|_2\langle\boldsymbol h,\Sigma_{J_\star^c}^{1/2}\boldsymbol\beta_{J_\star^c}\rangle+\langle\boldsymbol u,\boldsymbol\xi\rangle-\frac12\|\boldsymbol u\|_2^2+N\lambda\|\boldsymbol\beta_{J_\star^c}\|_1-\frac12\|\boldsymbol\xi\|_2^2\bigg\}.
\end{align*}By CGMT, that is, Lemma~\ref{lem:common-cgmt},
\begin{align}
\forall x\in\bR,\qquad \mathbb P_{\mathbb X_{J_\star^c}}\{\Phi_I^{K,H}<x\mid\boldsymbol\xi\}&\leq2\mathbb P_{\mathrm{aux}}\{\phi_I^{K,H}\leq x\mid\boldsymbol\xi\}.
\tag{\ref{eq:cgmt_comparisons}}
\end{align}Note that $\Phi_{\{t\}}^{t,\frac{N\lambda}{z_q}}$ is precisely obtained by adding the additional dual variable constraint $\|\vu\|_2 \le \frac{N\lambda}{z_q}$ to the feasible set of the convex optimization problem in $\Omega_{\mathrm{cost}}$ — that is, the original optimization problem corresponds to the case $\vu \in \mathbb{R}^N$. Hence (noting that $\vu$ corresponds to a maximization problem), we have $\Phi_{\{t\}}^{t,\frac{N\lambda}{z_q}} \le \inf_{\|\Sigma_{J_\star^c}^{1/2}\boldsymbol\beta_{J_\star^c}\|_2=t}\bigg[\frac12\|\boldsymbol\xi-\mathbb X_{J_\star^c}\boldsymbol\beta_{J_\star^c}\|_2^2+N\lambda\|\boldsymbol\beta_{J_\star^c}\|_1-\frac12\|\boldsymbol\xi\|_2^2\bigg]>N(r_\star^2+\mathrm{Var}_{J_\star^c})$. That is, we have the following set inclusion relation
\begin{align*}
    \Omega_{\mathrm{cost}}^c \subset \Big\{\Phi_{\{t\}}^{t,\frac{N\lambda}{z_q}} < 2N(r_\star^2 + \Var_{J_\star^c}) \Big\}.
\end{align*}By \eqref{eq:cgmt_comparisons}, this reduces to a problem concerning the auxiliary problem, namely it suffices to prove that on $\Omega^{\mathrm{aux}}$, $\phi_{\{t\}}^{N,\frac{N\lambda}{z_q}} \geq 125N(r_\star^2+\Var_{J_\star^c})$.

Fix $\bbeta_{J_\star^c}\in V_{J_\star^c}$ such that $\|\Sigma_{J_\star^c}^{1/2}\bbeta_{J_\star^c}^\star\|_2=t$, then $64\sigma_\xi\sqrt{\frac{d_q}{N}}\le t\le\frac{\sigma_\xi}{4}$. Work on $\Omega_\xi\cap\Omega^{\mathrm{aux}}$. Take $\vu = \frac{N\lambda}{z_q}\frac{\bxi+t\vg}{\|\bxi+t\vg\|_2}$. Then $t\langle\vg,\vu\rangle+\langle\bxi,\vu\rangle = \frac{N\lambda}{z_q}\|\bxi+t\vg\|_2$, and
\begin{align*}
    &\phi_{\{t\}}^{t,\frac{N\lambda}{z_q}}\geq \frac{N\lambda}{z_q}\|\bxi+t\vg\|_2 -  \frac{1}{2}(\frac{N\lambda}{z_q})^2 - \frac{1}{2}\|\bxi\|_2^2 \\
    &+ \inf\Big( \frac{N\lambda}{z_q}\langle\vh,\Sigma_{J_\star^c}^{1/2}\bbeta_{J_\star^c}\rangle + N\lambda\|\bbeta_{J_\star^c}\|_1 :\, \bbeta_{J_\star^c}\in V_{J_\star^c}\mbox{ such that } \|\Sigma_{J_\star^c}^{1/2}\bbeta_{J_\star^c}^\star\|_2=t\Big).
\end{align*}By the calibration equation, that is, \eqref{eq:common-q-calibration}, $z_q = \frac{\lambda\sqrt{N}}{\sigma_\xi}\frac{\sqrt{N(N-d_q)}}{N-q}$. Since $d_q\leq q<N$, $(N-q)^2\leq N(N-q)$, and consequently, $z_q\ge\frac{\lambda\sqrt N}{\sigma_\xi}$. From the assumption $\lambda\geq\frac{2\sigma_\xi z_N^{\mathrm{all}}}{\sqrt{N}}$, we further obtain $z_q\ge2z_N^{\mathrm{all}}$. Now, applying Lemma~\ref{lem:large-spectral-control} to $I=J_\star^c$, and $u=z_q$, we obtain
\begin{align*} 
q^2 = \Big(\sum_{j\in J_\star^c}2\Phi_G\big(\frac{z_q}{\sqrt{\sigma_j}}\big)\Big)^2\le 2N\sum_{j\in J_\star^c}\bE\big(|G|-\frac{z_q}{\sqrt{\sigma_j}}\big)_+^2 = 2Nd_q. \end{align*}Moreover, from $d_q\geq 0$, $\frac{N-q}{\sqrt{N-d_q}}\geq \frac{N-q}{\sqrt{N}}=\sqrt{N}-\frac{q}{\sqrt{N}}$. From $q^2\leq 2Nd_q$, there holds $\frac{q}{\sqrt{N}}\leq \sqrt{2d_q}$. Consequently,
\begin{align*}
    \sigma_\xi\sqrt N-\sigma_\xi\sqrt{2d_q}\le\frac{N\lambda}{z_q}=\frac{\sigma_\xi(N-q)}{\sqrt{N-d_q}}\le\sigma_\xi\sqrt N,
\end{align*}where the right-hand-side follows from $d_q\leq q<N$, hence $(N-q)^2\leq N(N-d_q)$.
Applying $a-b=\frac{a^2-b^2}{a+b}$ to $a=\|\bxi+t\vg\|_2$ and $b=\|\bxi\|_2$, together with $\|\bxi+t\vg\|_2\leq \sqrt{2\|\bxi\|_2^2 + 2t^2\|\vg\|_2^2}$ and $\|\bxi+t\vg\|_2+\|\bxi\|_2\geq \|\bxi\|_2$ yield
{\begin{align*} \|\boldsymbol\xi+t\boldsymbol g\|_2-\|\boldsymbol\xi\|_2&\ge\frac{t^2\|\boldsymbol g\|_2^2}{\sqrt{\|\boldsymbol\xi\|_2^2+t^2\|\boldsymbol g\|_2^2}+\|\boldsymbol\xi\|_2}-\frac{2t|\langle\boldsymbol\xi,\boldsymbol g\rangle|}{\|\boldsymbol\xi\|_2}\ge\frac{t^2\sqrt N}{8\sigma_\xi}-2t\sqrt{d_q},\\ \frac12\bigg(\frac{N\lambda}{z_q}-\|\boldsymbol\xi\|_2\bigg)^2&\le\frac{\sigma_\xi^2}{2}(\sqrt{2d_q}+\sqrt{d_q})^2\le3\sigma_\xi^2d_q,\end{align*}}
where we have used the definition of $\Omega_\xi$, and of $\Omega_\vg^{\mathrm{aux}}$. Notice that for any $j\in J_\star^c$, $\frac{z_q}{\sqrt{\sigma_j}}+h_j\sign(\beta_j)\geq -(|h_j|-\frac{z_q}{\sqrt{\sigma_j}})_+$, hence
\begin{align*}
    &N\lambda\|\bbeta_{J_\star^c}\|_1 + \frac{N\lambda}{z_q}\langle\vh,\Sigma_{J_\star^c}^{1/2}\bbeta_{J_\star^c}\rangle = \frac{N\lambda}{z_q}\sum_{j\in J_\star^c}\sqrt{\sigma_j}|\beta_j|\Big(\frac{z_q}{\sqrt{\sigma_j}} + h_j\sign(\beta_j)\Big)\geq -\frac{N\lambda}{z_q}\sum_{j\in J_\star^c}\sqrt{\sigma_j}|\beta_j|\Big((|h_j|-\frac{z_q}{\sqrt{\sigma_j}})_+\Big),
\end{align*}which is further bounded from below by $-\frac{N\lambda}{z_q}\|\Sigma_{J_\star^c}^{1/2}\bbeta_{J_\star^c}\|_2\Big(\sum_{j\in J_\star^c}(|h_j|-\frac{z_q}{\sqrt{\sigma_j}})_+^2\Big)^{\frac{1}{2}}$ due to Cauchy-Schwartz. 

In summary,
\begin{align*} \phi_{\{t\}}^{t,\frac{N\lambda}{z_q}}&\ge\frac{N\lambda}{z_q}\big(\|\boldsymbol\xi+t\boldsymbol g\|_2-\|\boldsymbol\xi\|_2-2t\sqrt{d_q}\big)-\frac12\bigg(\frac{N\lambda}{z_q}-\|\boldsymbol\xi\|_2\bigg)^2\\ &\ge\frac{N\lambda}{z_q}\bigg(\frac{t^2\sqrt N}{8\sigma_\xi}-4t\sqrt{d_q}\bigg)-3\sigma_\xi^2d_q\ge\frac{Nt^2}{32}-3\sigma_\xi^2d_q\ge125N(r_\star^2+\mathrm{Var}_{J_\star^c}). \end{align*}

\end{proof}

\subsection{Lower bound}\label{sec:diagonal-lower}
Retain $C$, $\varepsilon$, and $A$ fixed at the beginning of this section. Write $z_{\mathrm{all}}=z_{q^{\mathrm{all}}}^{\mathrm{all}}$, $d_{\mathrm{all}}=d_{q^{\mathrm{all}}}^{\mathrm{all}}$, and $\kappa_\pm=1\pm\frac{\varepsilon r(V_{J_\star},V_{J_\star^c})}{\sigma_\xi}$. For any $u\geq 0$, recall that $S_u(\cdot)=\sign(\cdot)(|\cdot|-u)_+$. Replace the definition of $\boldsymbol\Delta_{\rho,u}(\boldsymbol h)$ preceding Lemma~\ref{lem:stochastic_argument_lower_bound_free_subspace} by the following all-coordinate formula. Gaussian symmetry permits replacing $\boldsymbol h$ by $-\boldsymbol h$ in~\eqref{eq:auxiliary-optimization_lower}. For $\rho,u\ge0$ and $\boldsymbol h\sim\mathcal N(\boldsymbol0,I_p)$,
\begin{align*} \Delta_{\rho,u,j}=S_{\frac{u}{\sqrt{\sigma_j}}}(h_j+\rho\sqrt{\sigma_j}\beta_j^\star)-\rho\sqrt{\sigma_j}\beta_j^\star,\qquad j\in[p]. \end{align*}

Replace the auxiliary objectives in~\eqref{eq:auxiliary-optimization_lower} and the proof of Proposition~\ref{prop:null_signal_upper_bound} as follows. Use $\boldsymbol w=\boldsymbol\beta-\boldsymbol\beta^\star$ in \eqref{eq:variational_representation_LASSO} and Section~\ref{sec:dimension_reduction_trick}, retain the normalization $\frac{1}{N\lambda}$, set $a=N\lambda\|\boldsymbol u\|_2$, and subtract $\|\boldsymbol\beta^\star\|_1$. Applying CGMT, that is, Lemma~\ref{lem:common-cgmt}, to $\mathbb X\Sigma^{-1/2}$ gives the auxiliary objective below, followed by its convex approximation:
\begin{align*} \Psi(\boldsymbol w,a)=\frac{a}{N\lambda}\big\|\|\Sigma^{1/2}\boldsymbol w\|_2\boldsymbol g-\boldsymbol\xi\big\|_2-\frac{a^2}{2N\lambda} -\frac{a}{N\lambda}\langle\boldsymbol h,\Sigma^{1/2}\boldsymbol w\rangle+\|\boldsymbol\beta^\star+\boldsymbol w\|_1-\|\boldsymbol\beta^\star\|_1,\quad\mbox{and}\\ \tilde\Psi(\boldsymbol w,a)=\frac{a}{N\lambda}\sqrt{\|\boldsymbol\xi\|_2^2+\|\boldsymbol g\|_2^2\|\Sigma^{1/2}\boldsymbol w\|_2^2}-\frac{a^2}{2N\lambda} -\frac{a}{N\lambda}\langle\boldsymbol h,\Sigma^{1/2}\boldsymbol w\rangle+\|\boldsymbol\beta^\star+\boldsymbol w\|_1-\|\boldsymbol\beta^\star\|_1,\qquad a\ge0. \end{align*}
The identity $\sqrt x-\sqrt y=\frac{x-y}{\sqrt x+\sqrt y}$ yields
\begin{align*} |\Psi(\boldsymbol w,a)-\tilde\Psi(\boldsymbol w,a)|\le\frac{2a|\langle\boldsymbol\xi,\boldsymbol g\rangle|}{N\lambda\|\boldsymbol\xi\|_2}\|\Sigma^{1/2}\boldsymbol w\|_2. \end{align*}

\begin{Lemmaextend}{lem:saddle_point_analysis_tilde_Psi}\label{lem:diagonal-auxiliary-saddle}
Let $A$ be as in Lemma~\ref{lem:large-soft-events}. Assume $\lambda\geq\frac{8\sigma_\xi z_N^{\mathrm{all}}}{\sqrt N}$ and $r^2(V_{J_\star},V_{J_\star^c})\leq A^{-16}\sigma_\xi^2$. The assertion also applies when $\boldsymbol\beta^\star=\boldsymbol0$ and $J_\star=\varnothing$. Define
\begin{align*} \Omega_\xi=\bigg\{\bigg|\frac{\|\boldsymbol\xi\|_2^2}{N\sigma_\xi^2}-1\bigg|\le\frac{\varepsilon r(V_{J_\star},V_{J_\star^c})}{2^{12}\sigma_\xi}\bigg\}. \end{align*}
For independent $\boldsymbol g\sim\mathcal N(\boldsymbol0,I_N)$ and $\boldsymbol h\sim\mathcal N(\boldsymbol0,I_p)$, independent of $\boldsymbol\xi$, with $\Omega_{\mathrm{soft}}^{\mathrm{aux}}$ as in Lemma~\ref{lem:large-soft-events}, define
\begin{align*} \Omega_{\mathrm g}^{\mathrm{aux}}=\bigg\{\bigg|\frac{\|\boldsymbol g\|_2^2}{N}-1\bigg|\le\frac{\varepsilon r(V_{J_\star},V_{J_\star^c})}{2^{12}\sigma_\xi}\bigg\}\cap\bigg\{\frac{|\langle\boldsymbol\xi,\boldsymbol g\rangle|}{\|\boldsymbol\xi\|_2}\le\frac{\sqrt N r(V_{J_\star},V_{J_\star^c})}{2^{20}A^6\sigma_\xi}\bigg\},\quad\mbox{and}\quad \Omega^{\mathrm{aux}}=\Omega_{\mathrm{soft}}^{\mathrm{aux}}\cap\Omega_{\mathrm g}^{\mathrm{aux}}. \end{align*}
Then
\begin{align*} \mathbb P(\Omega_\xi^c)+\sup_{\boldsymbol\xi\in\Omega_\xi}\mathbb P_{\mathrm{aux}}((\Omega^{\mathrm{aux}})^c\mid\boldsymbol\xi)\leq2^{10}A\exp\bigg(-\frac{Nr^2(V_{J_\star},V_{J_\star^c})}{2^{45}A^{12}\sigma_\xi^2}\bigg). \end{align*}
On $\Omega_\xi\cap\Omega^{\mathrm{aux}}$, $\tilde\Psi$ has a saddle point $(\boldsymbol w^\circ,a^\circ)$ satisfying
\begin{align*} \frac{\sigma_\xi\sqrt N}{2}\leq a^\circ\leq2\sigma_\xi\sqrt N,\qquad\mbox{and}\qquad \frac{r(V_{J_\star},V_{J_\star^c})}{2A}\leq\|\Sigma^{1/2}\boldsymbol w^\circ\|_2\leq8Ar(V_{J_\star},V_{J_\star^c}). \end{align*}
For $\|\Sigma^{1/2}\boldsymbol w\|_2\leq\frac{\sigma_\xi}{2}$,
\begin{align*} \tilde\Psi(\boldsymbol w,a^\circ)\geq\tilde\Psi(\boldsymbol w^\circ,a^\circ)+\frac{1}{32\lambda}\|\Sigma^{1/2}(\boldsymbol w-\boldsymbol w^\circ)\|_2^2. \end{align*}
\end{Lemmaextend}

Replace Proposition~\ref{prop:free-variance-lower-bound} by the following statement.
\begin{Propositionextend}{prop:free-variance-lower-bound}
Retain the constant $A$ of Lemma~\ref{lem:large-soft-events}. If $\lambda\ge\frac{8\sigma_\xi z_N^{\mathrm{all}}}{\sqrt N}$ and $r^2(V_{J_\star},V_{J_\star^c})\le A^{-16}\sigma_\xi^2$, then, with probability at least $1-A^6\exp\big(-\frac{Nr^2(V_{J_\star},V_{J_\star^c})}{A^{32}\sigma_\xi^2}\big)$, every LASSO minimizer satisfies
\begin{align*} \|\Sigma^{1/2}(\widehat{\boldsymbol\beta}-\boldsymbol\beta^\star)\|_2^2\ge\frac{r^2(V_{J_\star},V_{J_\star^c})}{16A^2}. \end{align*}
\end{Propositionextend}
\begin{proof}
Replace Lemma~\ref{lem:stochastic_argument_lower_bound_free_subspace}, including~\eqref{eq:profile-event-bound} and~\eqref{eq:stochastic-event-bounds}, by Lemmas~\ref{lem:large-soft-events} and~\ref{lem:diagonal-auxiliary-saddle}. Replace Lemma~\ref{lem:saddle_point_analysis_tilde_Psi} by Lemma~\ref{lem:diagonal-auxiliary-saddle}. Replace Lemma~\ref{lem:def_x_H} by the conditional-median calculation below, and replace Lemma~\ref{lem:auxiliary_lower_bound} by the two auxiliary bounds below. Retain Lemma~\ref{lem:common-cgmt}, Lemma~\ref{lem:common-gaussian-isoperimetry}, and Lemma~\ref{lemma:Lipschitz_convex_concave_function}.
Retain $\Omega_\xi$, $\Omega^{\mathrm{aux}}=\Omega_{\mathrm{soft}}^{\mathrm{aux}}\cap\Omega_{\mathrm g}^{\mathrm{aux}}$, $\Psi$, and $\tilde\Psi$ from Lemma~\ref{lem:diagonal-auxiliary-saddle}. Conditional on $\boldsymbol\xi\in\Omega_\xi$, put
\begin{align*} \widetilde\phi=\min_{\|\Sigma^{1/2}\boldsymbol w\|_2\le16Ar(V_{J_\star},V_{J_\star^c})}\max_{0\le a\le4\sigma_\xi\sqrt N}\tilde\Psi(\boldsymbol w,a). \tag{\ref{eq:auxiliary-optimization_lower}} \end{align*}
By the calculation in~\eqref{eq:auxiliary-lipschitz} and Lemma~\ref{lemma:Lipschitz_convex_concave_function}, $\operatorname{Lip}(\widetilde\phi)\le\frac{2^8A\sigma_\xi r(V_{J_\star},V_{J_\star^c})}{\lambda\sqrt N}$. Let $m$ be a conditional median of $\widetilde\phi$, and define
\begin{align*} \Omega_{\mathrm{med}}^{\mathrm{aux}}=\bigg\{|\widetilde\phi-m|\le\frac{r^2(V_{J_\star},V_{J_\star^c})}{2^{14}A^2\lambda}\bigg\},\qquad\mbox{and}\qquad \Omega_{\mathrm{lower}}^{\mathrm{aux}}=\Omega_{\mathrm{soft}}^{\mathrm{aux}}\cap\Omega_{\mathrm g}^{\mathrm{aux}}\cap\Omega_{\mathrm{med}}^{\mathrm{aux}}. \end{align*}
Thus $\Omega=\Omega_\xi\cap\Omega_{\mathrm{soft}}^{\mathrm{aux}}\cap\Omega_{\mathrm g}^{\mathrm{aux}}\cap\Omega_{\mathrm{med}}^{\mathrm{aux}}$. Lemma~\ref{lem:common-gaussian-isoperimetry} and Lemma~\ref{lem:diagonal-auxiliary-saddle} give
\begin{align*} \mathbb P(\Omega_\xi^c)+\sup_{\boldsymbol\xi\in\Omega_\xi}\mathbb P_{\mathrm{aux}}((\Omega_{\mathrm{lower}}^{\mathrm{aux}})^c\mid\boldsymbol\xi)\le2^{10}A\exp\bigg(-\frac{Nr^2(V_{J_\star},V_{J_\star^c})}{2^{45}A^{12}\sigma_\xi^2}\bigg) +2\exp\bigg(-\frac{Nr^2(V_{J_\star},V_{J_\star^c})}{2^{50}A^6\sigma_\xi^2}\bigg). \end{align*}

In~\eqref{eq:primary-optimization_lower}, use $\boldsymbol w=\boldsymbol\beta-\boldsymbol\beta^\star$ and $\mathbb X\boldsymbol w=(\mathbb X\Sigma^{-1/2})(\Sigma^{1/2}\boldsymbol w)$. With the normalization $\frac{1}{N\lambda}$ and the dual rescaling of Lemma~\ref{lem:diagonal-auxiliary-saddle}, take the two primal domains
\begin{align*} \bigg\{\boldsymbol w:\|\Sigma^{1/2}\boldsymbol w\|_2\le\frac{r(V_{J_\star},V_{J_\star^c})}{4A}\bigg\},\quad\mbox{and}\quad \bigg\{\boldsymbol w:\|\Sigma^{1/2}\boldsymbol w\|_2\le16Ar(V_{J_\star},V_{J_\star^c})\bigg\}. \end{align*}
Write $\Phi_{\frac{r(V_{J_\star},V_{J_\star^c})}{4A},H}$ and $\Phi_{16Ar(V_{J_\star},V_{J_\star^c}),H}$ for the resulting primary values, and use $\phi$ for their auxiliary counterparts. These replace~\eqref{eq:primary-optimization_lower} and~\eqref{eq:auxiliary-optimization_lower}; their auxiliary objective is $\Psi$.

Fix $\boldsymbol\xi\in\Omega_\xi$ and work on $\Omega_{\mathrm{lower}}^{\mathrm{aux}}$. The saddle point from Lemma~\ref{lem:diagonal-auxiliary-saddle} satisfies $\|\Sigma^{1/2}\boldsymbol w^\circ\|_2\ge\frac{r(V_{J_\star},V_{J_\star^c})}{2A}$. Its Hessian lower bound, as in~\eqref{eq:restricted_cost_gap}, therefore gives
\begin{align*} \inf_{\|\Sigma^{1/2}\boldsymbol w\|_2\le\frac{r(V_{J_\star},V_{J_\star^c})}{4A}}\big(\tilde\Psi(\boldsymbol w,a^\circ)-\widetilde\phi\big)\ge\frac{r^2(V_{J_\star},V_{J_\star^c})}{2^{10}A^2\lambda}. \tag{\ref{eq:pointwise-auxiliary-gap}} \end{align*}
The perturbation bound in the same lemma, evaluated at $a^\circ$ on the smaller ball and at $\boldsymbol w^\circ$ on the larger ball, yields, for every $H\ge4\sigma_\xi\sqrt N$,
\begin{align*} \phi_{\frac{r(V_{J_\star},V_{J_\star^c})}{4A},H}\ge\widetilde\phi+\frac{3r^2(V_{J_\star},V_{J_\star^c})}{2^{12}A^2\lambda}\ge m+\frac{11r^2(V_{J_\star},V_{J_\star^c})}{2^{14}A^2\lambda},\quad\mbox{and}\tag{\ref{eq:small-ball-auxiliary-gap}}\\ \phi_{16Ar(V_{J_\star},V_{J_\star^c}),H}\le\widetilde\phi+\frac{r^2(V_{J_\star},V_{J_\star^c})}{2^{12}A^2\lambda}\le m+\frac{5r^2(V_{J_\star},V_{J_\star^c})}{2^{14}A^2\lambda}. \tag{\ref{eq:large-ball-auxiliary-bound}} \end{align*}
Here the scalar maximizer at $\boldsymbol w^\circ$ is the positive part below, by~\eqref{eq:variational_representation_ell_2_square}:
\begin{align*} \bigg(\big\|\|\Sigma^{1/2}\boldsymbol w^\circ\|_2\boldsymbol g-\boldsymbol\xi\big\|_2-\langle\boldsymbol h,\Sigma^{1/2}\boldsymbol w^\circ\rangle\bigg)_+\le a^\circ+\frac{2|\langle\boldsymbol\xi,\boldsymbol g\rangle|}{\|\boldsymbol\xi\|_2}\|\Sigma^{1/2}\boldsymbol w^\circ\|_2\le4\sigma_\xi\sqrt N. \end{align*} These inequalities replace~\eqref{eq:small-ball-auxiliary-gap} and~\eqref{eq:large-ball-auxiliary-bound}, respectively. In place of~\eqref{eq:conditional-auxiliary-tails}, they imply
\begin{align*} \mathbb P_{\mathrm{aux}}\bigg\{\phi_{\frac{r(V_{J_\star},V_{J_\star^c})}{4A},H}\le m+\frac{r^2(V_{J_\star},V_{J_\star^c})}{2^{11}A^2\lambda}\,\bigg|\,\boldsymbol\xi\bigg\} +\mathbb P_{\mathrm{aux}}\bigg\{\phi_{16Ar(V_{J_\star},V_{J_\star^c}),H}\ge m+\frac{r^2(V_{J_\star},V_{J_\star^c})}{2^{11}A^2\lambda}\,\bigg|\,\boldsymbol\xi\bigg\}\\ \le 2\mathbb P_{\mathrm{aux}}((\Omega^{\mathrm{aux}})^c\mid\boldsymbol\xi)+4\exp\bigg(-\frac{Nr^2(V_{J_\star},V_{J_\star^c})}{2^{50}A^6\sigma_\xi^2}\bigg). \tag{\ref{eq:conditional-auxiliary-tails}} \end{align*}

Both primary domains are convex and compact. Applying CGMT, that is, Lemma~\ref{lem:common-cgmt}, at the intermediate level $m+\frac{r^2(V_{J_\star},V_{J_\star^c})}{2^{11}A^2\lambda}$ gives the following replacement for~\eqref{eq:conditional-primary-tails}:
\begin{align*} &\mathbb P\big\{\Phi_{\frac{r(V_{J_\star},V_{J_\star^c})}{4A},H}=\Phi_{16Ar(V_{J_\star},V_{J_\star^c}),H}\mid\boldsymbol\xi\big\}\le4\mathbb P_{\mathrm{aux}}((\Omega^{\mathrm{aux}})^c\mid\boldsymbol\xi)+8\exp\bigg(-\frac{Nr^2(V_{J_\star},V_{J_\star^c})}{2^{50}A^6\sigma_\xi^2}\bigg). \tag{\ref{eq:conditional-primary-tails}} \end{align*}
For every fixed primary sample, the bound
\begin{align*} \sup_{\|\Sigma^{1/2}\boldsymbol w\|_2\le16Ar(V_{J_\star},V_{J_\star^c})}\|\boldsymbol\xi-\mathbb X\boldsymbol w\|_2\le\|\boldsymbol\xi\|_2+16Ar(V_{J_\star},V_{J_\star^c})\|\mathbb X\Sigma^{-1/2}\|_{\mathrm{op}} \end{align*}
justifies letting $H\to\infty$, as in Section~\ref{sec:cgmt}. Existence of a global LASSO minimizer in the smaller ball implies equality of the two primary values for all sufficiently large $H$. Fatou's lemma and integration over $\Omega_\xi$ give the stated probability, since $A\ge2^{34}$.
\end{proof}

\begin{proof}[Proof of Lemma~\ref{lem:diagonal-auxiliary-saddle}]
Replace the saddle-point calculation in~\eqref{eq:soft-threshold-kkt}--\eqref{eq:saddle-structure} by the following calculation for $\tilde\Psi$.

\textbf{Stochastic Argument.}
For $\Omega=\Omega_\xi\cap\Omega_{\mathrm{soft}}^{\mathrm{aux}}\cap\Omega_{\mathrm g}^{\mathrm{aux}}$ defined in the statement, Gaussian concentration and Lemma~\ref{lem:large-soft-events} give
\begin{align*} \mathbb P(\Omega_\xi^c)+\sup_{\boldsymbol\xi\in\Omega_\xi}\mathbb P_{\mathrm{aux}}\big((\Omega^{\mathrm{aux}})^c\mid\boldsymbol\xi\big)\le2^{10}A\exp\bigg(-\frac{N r^2(V_{J_\star},V_{J_\star^c})}{2^{45}A^{12}\sigma_\xi^2}\bigg). \end{align*}

\textbf{Deterministic Argument.}
Work on $\Omega=\Omega_\xi\cap\Omega_{\mathrm{soft}}^{\mathrm{aux}}\cap\Omega_{\mathrm g}^{\mathrm{aux}}$. For
\begin{align*} \frac{N\lambda}{z_{\mathrm{all}}\kappa_+}\le a\le\frac{N\lambda}{z_{\mathrm{all}}\kappa_-}, \end{align*}
Lemma~\ref{lem:large-soft-events} implies
\begin{align*} \bigg(\sum_{j=1}^p\bigg(|h_j|-\frac{N\lambda}{a\sqrt{\sigma_j}}\bigg)_+^2\bigg)^{\frac{1}{2}}\le\frac{A\sqrt N r(V_{J_\star},V_{J_\star^c})}{\sigma_\xi}<\|\boldsymbol g\|_2. \end{align*}
Thus $\tilde\Psi(\cdot,a)$ is coercive and strictly convex. Its minimizer $\boldsymbol w^\circ(a)$ satisfies, by the coordinatewise KKT calculation in \eqref{eq:soft-threshold-kkt},
\begin{align*} \Sigma^{1/2}\boldsymbol w^\circ(a)=\frac{\boldsymbol\Delta_{\rho,\frac{N\lambda}{a}}}{\rho},\quad\mbox{and}\quad \rho=\frac{\|\boldsymbol g\|_2^2}{\sqrt{\|\boldsymbol\xi\|_2^2+\|\boldsymbol g\|_2^2\|\Sigma^{1/2}\boldsymbol w^\circ(a)\|_2^2}}. \tag{\ref{eq:solution_S_Nlambda_s_H}} \end{align*}
The inequalities $|S_t(h+b)-b|\le|h|+t$ and $|S_t(h+b)-b|\le(|h|-t)_++2|b|$, applied on $J_\star$ and $J_\star^c$, respectively, give
\begin{align*} \|\Sigma^{1/2}\boldsymbol w^\circ(a)\|_2\le\frac{2Ar(V_{J_\star},V_{J_\star^c})}{\sigma_\xi}\sqrt{2\sigma_\xi^2+\|\Sigma^{1/2}\boldsymbol w^\circ(a)\|_2^2}+2r(V_{J_\star},V_{J_\star^c}),\\ \bigg(1-\frac{2Ar(V_{J_\star},V_{J_\star^c})}{\sigma_\xi}\bigg)\|\Sigma^{1/2}\boldsymbol w^\circ(a)\|_2\le(2\sqrt2A+2)r(V_{J_\star},V_{J_\star^c}),\quad\mbox{and}\quad \|\Sigma^{1/2}\boldsymbol w^\circ(a)\|_2\le8Ar(V_{J_\star},V_{J_\star^c}). \end{align*}
Consequently $\frac{\sqrt N}{2\sigma_\xi}\le\rho\le\frac{2\sqrt N}{\sigma_\xi}$, and Lemma~\ref{lem:large-soft-events} gives
\begin{align*} \frac{r(V_{J_\star},V_{J_\star^c})}{2A}\le\|\Sigma^{1/2}\boldsymbol w^\circ(a)\|_2\le8Ar(V_{J_\star},V_{J_\star^c}),\quad\mbox{and}\quad \bigg|\frac{\sigma_\xi\rho}{\sqrt N}-1\bigg|\le\frac{\varepsilon r(V_{J_\star},V_{J_\star^c})}{128\sigma_\xi}. \end{align*}
The last inequality follows from the formula for $\rho$, $\Omega_\xi\cap\Omega_{\mathrm g}^{\mathrm{aux}}$, and $\frac{64A^2r^2(V_{J_\star},V_{J_\star^c})}{\sigma_\xi^2}\le\frac{\varepsilon r(V_{J_\star},V_{J_\star^c})}{2^{10}\sigma_\xi}$.

The map $a\mapsto\min_{\boldsymbol w}\tilde\Psi(\boldsymbol w,a)$ is strictly concave. Apply Danskin's theorem as in the proof of Lemma~\ref{lem:saddle_point_analysis_tilde_Psi}:
\begin{align*} \rho\frac{\mathrm d}{\mathrm da}\min_{\boldsymbol w}\tilde\Psi(\boldsymbol w,a)=\frac{\|\boldsymbol g\|_2^2-\langle\boldsymbol h,\boldsymbol\Delta_{\rho,\frac{N\lambda}{a}}\rangle-\rho a}{N\lambda}. \end{align*}
The calibration identity gives
\begin{align*} \frac{\sqrt N}{\sigma_\xi}\frac{N\lambda}{z_{\mathrm{all}}}=\frac{N-q^{\mathrm{all}}}{\sqrt{1-\frac{d_{\mathrm{all}}}{N}}},\quad\mbox{and}\quad 0\le\frac{N-q^{\mathrm{all}}}{\sqrt{1-\frac{d_{\mathrm{all}}}{N}}}-(N-q^{\mathrm{all}})\le2d_{\mathrm{all}}\le\frac{50N r^2(V_{J_\star},V_{J_\star^c})}{\sigma_\xi^2}\le\frac{\varepsilon Nr(V_{J_\star},V_{J_\star^c})}{64\sigma_\xi}. \end{align*}
Combining the endpoint correlation bounds of Lemma~\ref{lem:large-soft-events} with 
\begin{align*} \sum_j2\Phi_G\big(\frac{z_{\mathrm{all}}\kappa_+}{\sqrt{\sigma_j}}\big)\le q^{\mathrm{all}}\le\sum_j2\Phi_G\big(\frac{z_{\mathrm{all}}\kappa_-}{\sqrt{\sigma_j}}\big) \end{align*}
yields
\begin{align*} \rho\frac{\mathrm d}{\mathrm da}\min_{\boldsymbol w}\tilde\Psi(\boldsymbol w,a)\bigg|_{a=\frac{N\lambda}{z_{\mathrm{all}}\kappa_+}}\ge\frac{1}{N\lambda}\bigg(\frac{\varepsilon Nr(V_{J_\star},V_{J_\star^c})}{4\sigma_\xi}-\frac{3\varepsilon Nr(V_{J_\star},V_{J_\star^c})}{64\sigma_\xi}-2d_{\mathrm{all}}\bigg)>0,\quad\mbox{and}\\ \rho\frac{\mathrm d}{\mathrm da}\min_{\boldsymbol w}\tilde\Psi(\boldsymbol w,a)\bigg|_{a=\frac{N\lambda}{z_{\mathrm{all}}\kappa_-}}\le\frac{1}{N\lambda}\bigg(-\frac{\varepsilon Nr(V_{J_\star},V_{J_\star^c})}{4\sigma_\xi}+\frac{3\varepsilon Nr(V_{J_\star},V_{J_\star^c})}{64\sigma_\xi}\bigg)<0. \end{align*}
Coercivity is uniform on the displayed interval, so Danskin's theorem applies on a common compact sublevel set. The interior maximizer $a^\circ$ and $\boldsymbol w^\circ=\boldsymbol w^\circ(a^\circ)$ form a saddle point, with
\begin{align*} \frac{\sigma_\xi\sqrt N}{2}\le a^\circ\le2\sigma_\xi\sqrt N,\quad\mbox{and}\quad \frac{r(V_{J_\star},V_{J_\star^c})}{2A}\le\|\Sigma^{1/2}\boldsymbol w^\circ\|_2\le8Ar(V_{J_\star},V_{J_\star^c}). \end{align*}

Replace the second-derivative estimate used in~\eqref{eq:restricted_cost_gap} by
\begin{align*} \frac{a^\circ\|\boldsymbol g\|_2^2\|\boldsymbol\xi\|_2^2}{N\lambda\big(\|\boldsymbol\xi\|_2^2+\|\boldsymbol g\|_2^2\|\Sigma^{1/2}\boldsymbol w\|_2^2\big)^{\frac{3}{2}}}\geq\frac{18}{125\lambda}>\frac{1}{16\lambda},\qquad \|\Sigma^{1/2}\boldsymbol w\|_2\leq\frac{\sigma_\xi}{2}. \end{align*}
The subgradient inequality at $\boldsymbol w^\circ$, as in Lemma~\ref{lem:radial-restriction-gap}, gives the asserted quadratic bound.
\end{proof}

\begin{Lemma}\label{lem:large-cross-process}
For an absolute constant $C\geq1$ fixed at the beginning of this section, for $u\geq0$ and $x\geq1$, with probability at least $1-Ce^{-x}$,
\begin{align*} \sup_{0\leq\rho\leq\frac{2\sqrt N}{\sigma_\xi}}\big|(1-\mathbb E)\langle\boldsymbol h_{J_\star^c},\boldsymbol\Delta_{\rho,u,J_\star^c}-\boldsymbol\Delta_{0,u,J_\star^c}\rangle\big|\leq\frac{C\sqrt{Nx}}{\sigma_\xi}\|\Sigma_{J_\star^c}^{1/2}\boldsymbol\beta_{J_\star^c}^\star\|_2. \end{align*}
Moreover, for every $\rho\geq0$,
\begin{align*} 0\leq\mathbb E\langle\boldsymbol h_{J_\star^c},\boldsymbol\Delta_{\rho,u,J_\star^c}-\boldsymbol\Delta_{0,u,J_\star^c}\rangle\leq C\rho^2\|\Sigma_{J_\star^c}^{1/2}\boldsymbol\beta_{J_\star^c}^\star\|_2^2. \end{align*}
For every $I\subset[p]$,
\begin{align*} \mathbb P\bigg\{\bigg|\sum_{j\in I}h_j\Delta_{0,u,j}-\sum_{j\in I}2\Phi_G\bigg(\frac{u}{\sqrt{\sigma_j}}\bigg)\bigg|>C\bigg[\sqrt{x\sum_{j\in I}2\Phi_G\bigg(\frac{u}{\sqrt{\sigma_j}}\bigg)}+x\bigg]\bigg\}\leq2e^{-x}. \end{align*}
\end{Lemma}
\begin{proof}
Replace the common-threshold application of~\eqref{eq:Bernstein_moment_soft_max} in Lemma~\ref{lem:stochastic_argument_lower_bound_free_subspace} by the following coordinatewise calculation. The Gaussian integration by parts in~\eqref{eq:exp_G_soft_max}, with $S_t^{\prime}=\mathbf1_{\{|\cdot|>t\}}$, gives
\begin{align*} \mathbb E[G(S_t(G+b)-b)]=\mathbb E[S_t'(G+b)]=\Phi_G(t-b)+\Phi_G(t+b),\\ 0\leq\Phi_G(t-b)+\Phi_G(t+b)-2\Phi_G(t)\leq\frac{b^2}{\sqrt{2\pi e}},\quad\mbox{and}\quad \sup_{v\in\mathbb R}|\Phi_G''(v)|=\frac{1}{\sqrt{2\pi e}}. \end{align*}
Taking $b=\rho\sqrt{\sigma_j}\beta_j^\star$ and $t=\frac{u}{\sqrt{\sigma_j}}$ proves the mean bound. Since $S_t$ is $1$-Lipschitz, independence gives, for $\rho,\rho'\geq0$,
\begin{align*} |h_j(\Delta_{\rho,u,j}-\Delta_{\rho',u,j})|\leq2|\rho-\rho'|\sqrt{\sigma_j}|\beta_j^\star|\,|h_j|,\quad\mbox{and}\quad \big\|(1-\mathbb E)\langle\boldsymbol h_{J_\star^c},\boldsymbol\Delta_{\rho,u,J_\star^c}-\boldsymbol\Delta_{\rho',u,J_\star^c}\rangle\big\|_{\psi_2}\leq C|\rho-\rho'|\|\Sigma_{J_\star^c}^{1/2}\boldsymbol\beta_{J_\star^c}^\star\|_2. \end{align*}
For the dyadic grids $\{\frac{2m\sqrt N}{2^k\sigma_\xi}:0\leq m\leq2^k\}$, $k\geq0$, attach each point to its left parent, and attach the initial grid to zero. The union bound with deviations $x+2k+4$ gives
\begin{align*} \mathbb P\bigg\{&\sup_{0\leq\rho\leq\frac{2\sqrt N}{\sigma_\xi}}\big|(1-\mathbb E)\langle\boldsymbol h_{J_\star^c},\boldsymbol\Delta_{\rho,u,J_\star^c}-\boldsymbol\Delta_{0,u,J_\star^c}\rangle\big|\\ &>\frac{C\sqrt N}{\sigma_\xi}\|\Sigma_{J_\star^c}^{1/2}\boldsymbol\beta_{J_\star^c}^\star\|_2\sum_{k\geq0}2^{-k}\sqrt{x+2k+4}\bigg\}\leq\sum_{k\geq0}2(2^k+1)e^{-x-2k-4}\leq Ce^{-x}. \end{align*}
Continuity in $\rho$ and $\sum_{k\geq0}2^{-k}\sqrt{x+2k+4}\leq C\sqrt x$ give the first assertion. Finally,~\eqref{eq:exp_G_soft_max} and~\eqref{eq:Bernstein_moment_soft_max}, with $a=\frac{u}{\sqrt{\sigma_j}}$, give
\begin{align*} \mathbb E[h_j\Delta_{0,u,j}]=2\Phi_G\bigg(\frac{u}{\sqrt{\sigma_j}}\bigg),\quad\mbox{and}\quad \mathbb E\big|h_j\Delta_{0,u,j}-\mathbb E[h_j\Delta_{0,u,j}]\big|^k\leq\frac{k!}{2}\,64\Phi_G\bigg(\frac{u}{\sqrt{\sigma_j}}\bigg)4^{k-2},\qquad k\geq2. \end{align*}
Bernstein's inequality proves the last assertion.
\end{proof}

\begin{Lemmaextend}{lem:stochastic_argument_lower_bound_free_subspace}\label{lem:large-soft-events}
Let $C$ and $A$ be the constants fixed at the beginning of this section. Assume $\lambda\ge\frac{8\sigma_\xi z_N^{\mathrm{all}}}{\sqrt N}$ and $r^2(V_{J_\star},V_{J_\star^c})\le2^{-40}\varepsilon^2A^{-4}\sigma_\xi^2$. There is an event $\Omega_{\mathrm{soft}}^{\mathrm{aux}}$ of probability at least $1-A\exp\big(-\frac{Nr^2(V_{J_\star},V_{J_\star^c})}{A\sigma_\xi^2}\big)$ on which, simultaneously for $\frac{\sqrt N}{2\sigma_\xi}\le\rho\le\frac{2\sqrt N}{\sigma_\xi}$ and $z_{\mathrm{all}}\kappa_-\le u\le z_{\mathrm{all}}\kappa_+$,
\begin{align*} \frac{\sqrt N r(V_{J_\star},V_{J_\star^c})}{A\sigma_\xi}\le\|\boldsymbol\Delta_{\rho,u}\|_2\le\frac{A\sqrt N r(V_{J_\star},V_{J_\star^c})}{\sigma_\xi}. \end{align*}
On the same event,
\begin{align*} \|\boldsymbol h_{J_\star}\|_2+z_{\mathrm{all}}\kappa_+\bigg(\sum_{j\in J_\star}\frac{1}{\sigma_j}\bigg)^{\frac{1}{2}} +\bigg(\sum_{j=1}^p\bigg(|h_j|-\frac{z_{\mathrm{all}}\kappa_-}{\sqrt{\sigma_j}}\bigg)_+^2\bigg)^{\frac{1}{2}}\le\frac{A\sqrt N r(V_{J_\star},V_{J_\star^c})}{\sigma_\xi}. \end{align*}
At each of the two endpoints for $u$, simultaneously for $0\le\rho\le\frac{2\sqrt N}{\sigma_\xi}$,
\begin{align*} \bigg|\langle\boldsymbol h,\boldsymbol\Delta_{\rho,u}\rangle-\sum_{j=1}^p2\Phi_G\bigg(\frac{u}{\sqrt{\sigma_j}}\bigg)\bigg|\le\frac{\varepsilon Nr(V_{J_\star},V_{J_\star^c})}{128\sigma_\xi}. \end{align*}
\end{Lemmaextend}
\begin{proof}
Replace the stochastic norm and correlation estimates in Lemma~\ref{lem:stochastic_argument_lower_bound_free_subspace} by the following events. If $\frac{Nr^2(V_{J_\star},V_{J_\star^c})}{A\sigma_\xi^2}\le\log A$, take $\Omega_{\mathrm{soft}}^{\mathrm{aux}}=\varnothing$; its asserted probability lower bound is nonpositive. Henceforth, $\frac{Nr^2(V_{J_\star},V_{J_\star^c})}{A\sigma_\xi^2}>\log A$. Every infimum below is over the rectangle in the statement.

\textbf{Stochastic Argument.} Lemma~\ref{lem:large-spectral-control} gives, at both endpoints,
\begin{align*} \bigg|\bigg(\sum_{j=1}^p\mathbb E\bigg(|G|-\frac{u}{\sqrt{\sigma_j}}\bigg)_+^2\bigg)^{\frac{1}{2}}-\sqrt{d_{\mathrm{all}}}\bigg|\le\frac{3\varepsilon\sqrt N r(V_{J_\star},V_{J_\star^c})}{\sigma_\xi}. \end{align*}
The identity $|G|=(|G|-t)_++\min(|G|,t)$ yields $1\le2t^2+2\mathbb E(|G|-t)_+^2$. Summing on $J_\star$, and using~\eqref{eq:common-free-variance} and the all-coordinate calibration comparisons in Section~\ref{sec:diagonal-calibration}, gives
\begin{align*} \frac{Nr^2(V_{J_\star},V_{J_\star^c})}{\sigma_\xi^2}\le3\sum_{j=1}^p\min\bigg\{\frac{N\sigma_j(\beta_j^\star)^2}{\sigma_\xi^2},\frac{N\lambda^2}{\sigma_\xi^2\sigma_j}\bigg\}+2\sum_{j=1}^p\mathbb E\bigg(|G|-\frac{\lambda\sqrt N}{\sigma_\xi\sqrt{\sigma_j}}\bigg)_+^2,\quad\mbox{and}\\ \sum_{j=1}^p\mathbb E\bigg(|G|-\frac{\lambda\sqrt N}{\sigma_\xi\sqrt{\sigma_j}}\bigg)_+^2\le25d_{\mathrm{all}}. \end{align*}
By the three cases in~\eqref{eq:softmax_equivalence}, $h_j\beta_j^\star\le0$ implies, throughout the rectangle,
\begin{align*} |\Delta_{\rho,u,j}|\ge\frac{1}{2}\min\bigg\{\frac{\sqrt{N\sigma_j}|\beta_j^\star|}{\sigma_\xi},\frac{\lambda\sqrt N}{\sigma_\xi\sqrt{\sigma_j}}\bigg\},\quad\mbox{and}\quad |\Delta_{\rho,u,j}|\ge\bigg(|h_j|-\frac{z_{\mathrm{all}}\kappa_+}{\sqrt{\sigma_j}}\bigg)_+. \end{align*}
Conditionally on $(|h_j|)_{j=1}^p$, each indicator $\mathbf1_{\{h_j\beta_j^\star\le0\}}$ has mean at least $\frac{1}{2}$. For nonnegative $a_j$ with $\sum_j a_j^2>0$, the inequality $\sqrt{\sum_j a_j^2\mathbf1_{\{h_j\beta_j^\star\le0\}}}\ge\frac{\sum_j a_j^2\mathbf1_{\{h_j\beta_j^\star\le0\}}}{\sqrt{\sum_j a_j^2}}$ therefore bounds the conditional expected selected norm below by $\frac{1}{2}\sqrt{\sum_j a_j^2}$. Applying this to the two coordinatewise bounds gives
\begin{align*} \mathbb E\inf_{\rho,u}\|\boldsymbol\Delta_{\rho,u}\|_2\ge\frac{1}{4}\bigg(\sum_{j=1}^p\min\bigg\{\frac{N\sigma_j(\beta_j^\star)^2}{\sigma_\xi^2},\frac{N\lambda^2}{\sigma_\xi^2\sigma_j}\bigg\}\bigg)^{\frac{1}{2}},\quad\mbox{and}\\ \mathbb E\inf_{\rho,u}\|\boldsymbol\Delta_{\rho,u}\|_2\ge\frac{1}{2}\bigg(\sum_{j=1}^p\mathbb E\bigg(|G|-\frac{z_{\mathrm{all}}\kappa_+}{\sqrt{\sigma_j}}\bigg)_+^2\bigg)^{\frac{1}{2}}-\frac{1}{2}. \end{align*}
Gaussian Poincare, as used in Lemma~\ref{lemma:null_signal_events}, gives, for each fixed $u\geq0$,
\begin{align*} 0\leq\bigg(\sum_{j=1}^p\mathbb E\bigg(|G|-\frac{u}{\sqrt{\sigma_j}}\bigg)_+^2\bigg)^{\frac{1}{2}}-\mathbb E\bigg(\sum_{j=1}^p\bigg(|h_j|-\frac{u}{\sqrt{\sigma_j}}\bigg)_+^2\bigg)^{\frac{1}{2}}\leq1. \end{align*}
Combining the preceding bounds,
\begin{align*} \frac{\sqrt N r(V_{J_\star},V_{J_\star^c})}{\sigma_\xi}\leq28\mathbb E\inf_{\rho,u}\|\boldsymbol\Delta_{\rho,u}\|_2+10+\frac{30\varepsilon\sqrt N r(V_{J_\star},V_{J_\star^c})}{\sigma_\xi},\quad\mbox{and}\\ \mathbb E\inf_{\rho,u}\|\boldsymbol\Delta_{\rho,u}\|_2\geq\frac{1-30\varepsilon}{28}\frac{\sqrt N r(V_{J_\star},V_{J_\star^c})}{\sigma_\xi}-\frac{5}{14}\geq\frac{\sqrt N r(V_{J_\star},V_{J_\star^c})}{32\sigma_\xi}-\frac{1}{2}. \end{align*}

Define
\begin{align*} \Omega_-^{\mathrm{aux}}=\bigg\{\inf_{\rho,u}\|\boldsymbol\Delta_{\rho,u}\|_2\ge\frac{\sqrt N r(V_{J_\star},V_{J_\star^c})}{A\sigma_\xi}\bigg\},\quad\mbox{and}\\ \Omega_+^{\mathrm{aux}}=\bigg\{\|\boldsymbol h_{J_\star}\|_2\le\frac{2\sqrt N r(V_{J_\star},V_{J_\star^c})}{\sigma_\xi}\bigg\}\cap \bigg\{\bigg(\sum_{j=1}^p\bigg(|h_j|-\frac{z_{\mathrm{all}}\kappa_-}{\sqrt{\sigma_j}}\bigg)_+^2\bigg)^{\frac{1}{2}}\le\frac{6\sqrt N r(V_{J_\star},V_{J_\star^c})}{\sigma_\xi}\bigg\}. \end{align*}
The infimum defining $\Omega_-^{\mathrm{aux}}$ is $1$-Lipschitz. Gaussian concentration gives $\mathbb P((\Omega_-^{\mathrm{aux}})^c)\le\exp\big(-\frac{Nr^2(V_{J_\star},V_{J_\star^c})}{2^{15}\sigma_\xi^2}\big)$. The two norms defining $\Omega_+^{\mathrm{aux}}$ are also $1$-Lipschitz, with expectations at most $\frac{\sqrt N r(V_{J_\star},V_{J_\star^c})}{\sigma_\xi}$ and $\frac{(5+3\varepsilon)\sqrt N r(V_{J_\star},V_{J_\star^c})}{\sigma_\xi}$, respectively. Concentration with deviation $\frac{\sqrt{2N}r(V_{J_\star},V_{J_\star^c})}{\sqrt A\sigma_\xi}$ therefore gives $\mathbb P((\Omega_+^{\mathrm{aux}})^c)\le2\exp\big(-\frac{Nr^2(V_{J_\star},V_{J_\star^c})}{A\sigma_\xi^2}\big)$.

Let $\Omega_{\mathrm{corr}}^{\mathrm{aux}}$ be the intersection, over the two endpoints for $u$, of
\begin{align*} \bigg\{\sup_{0\le\rho\le\frac{2\sqrt N}{\sigma_\xi}}\bigg|\langle\boldsymbol h,\boldsymbol\Delta_{\rho,u}\rangle-\sum_{j=1}^p2\Phi_G\bigg(\frac{u}{\sqrt{\sigma_j}}\bigg)\bigg|\le\frac{\varepsilon Nr(V_{J_\star},V_{J_\star^c})}{128\sigma_\xi}\bigg\}. \end{align*}
On $\Omega_+^{\mathrm{aux}}$, the estimation subspace contribution to the difference from the unshifted correlation is, by~\eqref{eq:softmax_equivalence}, at most
\begin{align*} 3\|\boldsymbol h_{J_\star}\|_2^2+u^2\sum_{j\in J_\star}\frac{1}{\sigma_j}\le\frac{28Nr^2(V_{J_\star},V_{J_\star^c})}{\sigma_\xi^2}. \end{align*}
Apply Lemma~\ref{lem:large-cross-process} at the two endpoints with $x=\frac{Nr^2(V_{J_\star},V_{J_\star^c})}{A\sigma_\xi^2}$; apply its last assertion with $I=[p]$. Lemma~\ref{lem:large-spectral-control} gives $\sum_{j=1}^p2\Phi_G(\frac{u}{\sqrt{\sigma_j}})\le\frac{8Nr(V_{J_\star},V_{J_\star^c})}{\sigma_\xi}$ there. Outside a set of probability at most $(4+2C)\exp\big(-\frac{Nr^2(V_{J_\star},V_{J_\star^c})}{A\sigma_\xi^2}\big)$, the correlation error on $\Omega_+^{\mathrm{aux}}$, divided by $\frac{Nr(V_{J_\star},V_{J_\star^c})}{\sigma_\xi}$, is at most
\begin{align*} (40+6C)\frac{r(V_{J_\star},V_{J_\star^c})}{\sigma_\xi}+C\sqrt{\frac{8}{A}}\sqrt{\frac{r(V_{J_\star},V_{J_\star^c})}{\sigma_\xi}}\le\frac{\varepsilon}{128}. \end{align*}
The last inequality uses $\frac{r(V_{J_\star},V_{J_\star^c})}{\sigma_\xi}\le2^{-20}\varepsilon A^{-2}$. Consequently,
\begin{align*} \Omega_{\mathrm{soft}}^{\mathrm{aux}}=\Omega_-^{\mathrm{aux}}\cap\Omega_+^{\mathrm{aux}}\cap\Omega_{\mathrm{corr}}^{\mathrm{aux}},\quad\mbox{and}\quad \mathbb P((\Omega_{\mathrm{soft}}^{\mathrm{aux}})^c)\le(7+2C)\exp\bigg(-\frac{Nr^2(V_{J_\star},V_{J_\star^c})}{A\sigma_\xi^2}\bigg), \end{align*}

where $7+2C\le A$.

\textbf{Deterministic Argument.} On $\Omega_{\mathrm{soft}}^{\mathrm{aux}}$, the auxiliary sum in the statement is at most $\frac{12\sqrt N r(V_{J_\star},V_{J_\star^c})}{\sigma_\xi}$, since $z_{\mathrm{all}}\kappa_+\le\frac{4\lambda\sqrt N}{\sigma_\xi}$. The inequalities used in~\eqref{eq:softmax_calculus} and~\eqref{eq:softmax_equivalence} give
\begin{align*} \|\boldsymbol\Delta_{\rho,u}\|_2\le\|\boldsymbol h_{J_\star}\|_2+u\bigg(\sum_{j\in J_\star}\frac{1}{\sigma_j}\bigg)^{\frac{1}{2}}+\bigg(\sum_{j\in J_\star^c}\bigg(|h_j|-\frac{u}{\sqrt{\sigma_j}}\bigg)_+^2\bigg)^{\frac{1}{2}} +2\rho\|\Sigma_{J_\star^c}^{1/2}\boldsymbol\beta_{J_\star^c}^\star\|_2\le\frac{16\sqrt N r(V_{J_\star},V_{J_\star^c})}{\sigma_\xi}. \end{align*}
The remaining assertions are the definitions of $\Omega_-^{\mathrm{aux}}$ and $\Omega_{\mathrm{corr}}^{\mathrm{aux}}$.
\end{proof}

\subsection{Technical lemmas in the anisotropic case}\label{sec:diagonal-calibration}


\begin{Lemma}\label{lem:large-spectral-control}
For every $t\geq0$, $[2\Phi_G(t)]^2\leq2\mathbb E(|G|-t)_+^2\,2\Phi_G(\frac{t}{2})$. Moreover, for $I\subset[p]$ and $u\geq2z_N^{\mathrm{all}}$,
\begin{align*} \bigg[\sum_{j\in I}2\Phi_G\bigg(\frac{u}{\sqrt{\sigma_j}}\bigg)\bigg]^2\leq2N\sum_{j\in I}\mathbb E\bigg(|G|-\frac{u}{\sqrt{\sigma_j}}\bigg)_+^2,\quad\mbox{and}\quad \bigg|\frac{\mathrm d}{\mathrm d\log u}\bigg[\sum_{j\in I}\mathbb E\bigg(|G|-\frac{u}{\sqrt{\sigma_j}}\bigg)_+^2\bigg]^{\frac{1}{2}}\bigg|\leq\sqrt{2N}. \end{align*}
\end{Lemma}
\begin{proof}
Replace the common-threshold comparison~\eqref{eq:Gaussian_local_comparison_soft_max} used in Lemma~\ref{lemma:null_signal_events} by the coordinatewise inequalities below. By~\eqref{eq:exp_G_soft_max} and Cauchy--Schwarz, $[2\Phi_G(t)]^2=\big[\mathbb E|G|(|G|-t)_+\big]^2\leq\mathbb E(|G|-t)_+^2\,\mathbb E[G^2\mathbf1_{\{|G|>t\}}]$. Since $\bE[G^2\1_{\{|G|>t\}}]=\frac{2}{\sqrt{2\pi}}\int_t^\infty x^2\exp(-\frac{1}{2}x^2)\mathrm{d}x=\frac{16}{\sqrt{2\pi}}\int_{\frac{t}{2}}^\infty v^2\exp(-2v^2)\mathrm{d}v$ through change-of-variable $x=2v$; and since $v^2\exp(-2v^2)=(v^2\exp(-\frac{3}{2}v^2))\exp(-\frac{1}{2}v^2)\leq \frac{2}{3e}\exp(-\frac{1}{2}v^2)$, there holds $\mathbb E[G^2\mathbf1_{\{|G|>t\}}]=\frac{16}{\sqrt{2\pi}}\int_{\frac{t}{2}}^{\infty}v^2e^{-2v^2}\,\mathrm dv\leq\frac{16}{3e}\,2\Phi_G\big(\frac{t}{2}\big)<4\Phi_G\big(\frac{t}{2}\big)$. From the assumption $u\geq 2z_N^{\mathrm{all}}$, we obtain $\sum_{j\in I}2\Phi_G(\frac{u}{2\sqrt{\sigma_j}})\leq N$, hence
\begin{align*} \bigg[\sum_{j\in I}2\Phi_G\bigg(\frac{u}{\sqrt{\sigma_j}}\bigg)\bigg]^2\leq2\sum_{j\in I}2\Phi_G\bigg(\frac{u}{2\sqrt{\sigma_j}}\bigg)\sum_{j\in I}\mathbb E\bigg(|G|-\frac{u}{\sqrt{\sigma_j}}\bigg)_+^2\leq2N\sum_{j\in I}\mathbb E\bigg(|G|-\frac{u}{\sqrt{\sigma_j}}\bigg)_+^2. \end{align*}Let $f(a)=\bE[(|G|-a)_+^2]$. Then \eqref{eq:identity_derivative_exp_positive} gives $f'(a)=-4(\frac{1}{\sqrt{2\pi}}\exp(-\frac{1}{2}a^2)-a\Phi_G(a))$; while \eqref{eq:identity_exp_square_soft_max} gives $2\Phi_G(a)-f(a)=2a(\frac{1}{\sqrt{2\pi}}\exp(-\frac{1}{2}a^2)-a\Phi_G(a))$. Therefore, $af'(a)=-2(2\Phi_G(a)-f(a))$. Taking $a=a_j=\frac{u}{\sqrt{\sigma_j}}$.From $u\frac{\mathrm d}{\mathrm{d}u}f(a_j)=a_jf'(a_j)$, we obtain that $u\frac{\mathrm{d}}{\mathrm{d}u}\bE(|G|-\frac{u}{\sqrt{\sigma_j}})_+^2 = 2\bE(|G|-\frac{u}{\sqrt{\sigma_j}})_+^2 - 4\Phi_G(\frac{u}{\sqrt{\sigma_j}})$. Taking sum over $j\in I$,
\begin{align*} u\frac{\mathrm d}{\mathrm du}\sum_{j\in I}\mathbb E\bigg(|G|-\frac{u}{\sqrt{\sigma_j}}\bigg)_+^2=2\sum_{j\in I}\mathbb E\bigg(|G|-\frac{u}{\sqrt{\sigma_j}}\bigg)_+^2-2\sum_{j\in I}2\Phi_G\bigg(\frac{u}{\sqrt{\sigma_j}}\bigg). \end{align*}
Let $F(u) = \sum_{j\in I}\bE(|G|-\frac{u}{\sqrt{\sigma_j}})_+^2$, and let $Q(u) = \sum_{j\in I}2\Phi_G(\frac{u}{\sqrt{\sigma_j}})$. The aforementioned equality is $uF'(u) = 2F(u)-2Q(u)$. Since $\frac{\mathrm{d}}{\mathrm{d}\log(u)}=u\frac{\mathrm{d}}{\mathrm{d}u}$, we have $\frac{\mathrm{d}}{\mathrm{d}\log(u)}\sqrt{F(u)} = \frac{uF'(u)}{2\sqrt{F(u)}} = \frac{F(u)-Q(u)}{\sqrt{F(u)}}$. By \eqref{eq:identity_exp_square_soft_max} again, $F(u)\leq Q(u)$. Hence $|\frac{\mathrm{d}}{\mathrm{d}\log(u)}\sqrt{F(u)}| = \leq \frac{Q(u)}{\sqrt{F(u)}}$. From $\big[\sum_{j\in I}2\Phi_G\big(\frac{u}{\sqrt{\sigma_j}}\big)\big]^2\leq2N\sum_{j\in I}\mathbb E\big(|G|-\frac{u}{\sqrt{\sigma_j}}\big)_+^2$, what we just proved, there holds $Q(u)\leq \sqrt{2NF(u)}$, and consequently, the last claim of this lemma is verified.
\end{proof}

We collect some facts that are used throughout the proof. From the definition of $z_q$ in \eqref{eq:common-q-calibration}, $\frac{\sigma_\xi z_q}{\lambda\sqrt{N}}=\frac{\sqrt{1-\frac{d_q}{N}}}{1-\frac{q}{N}}$. From $z_q\geq \frac{\lambda\sqrt{N}}{\sigma_\xi}$ and $\lambda\geq 2\frac{\sigma_\xi z_N^{\mathrm{all}}}{\sqrt{N}}$ we obtain $z_q\geq 2z_N^{\mathrm{all}}$. Applying Lemma~\ref{lem:large-spectral-control} with $u=z_q$, we obtain $q^2 = ( \sum_{j\in J_\star^c} 2\Phi_G(\frac{z_q}{\sqrt{\sigma_j}}) )^2\leq 2N\sum_{j\in J_\star^c}\bE(|G|-\frac{z_q}{\sqrt{\sigma_j}})_+^2=2Nd_q$. Moreover, $\log(\frac{\sigma_\xi z_q}{\lambda\sqrt{N}}) = \log(\frac{\sqrt{1-\frac{d_q}{N}}}{1-\frac{q}{N}}) \leq -\log(1-\frac{q}{N}) \leq 2\frac{q}{N}$. From Lemma~\ref{lem:large-spectral-control}, we integral from $\frac{\lambda\sqrt{N}}{\sigma_\xi}$ to $z_q$, we obtain $(\sum_{j\in J_\star^c}\bE(|G|-\frac{\lambda\sqrt{N}}{\sigma_\xi\sqrt{\sigma_j}})_+^2)^{\frac{1}{2}}\leq \sqrt{d_q} + \sqrt{2N}\log(\frac{\sigma_\xi z_q}{\lambda\sqrt{N}})\leq \sqrt{d_q} + \sqrt{2N}\frac{2q}{N} \leq 5\sqrt{d_q}$. In summary,
\begin{align}\label{eq:consequence_large-spectral-control}
\begin{aligned}
    &z_q\geq\frac{\lambda\sqrt N}{\sigma_\xi},\qquad q\leq\sqrt{2Nd_q},\qquad \log\bigg(\frac{\sigma_\xi z_q}{\lambda\sqrt N}\bigg)=\log\bigg(\frac{\sqrt{1-\frac{d_q}{N}}}{1-\frac{q}{N}}\bigg)\leq\frac{2q}{N},\quad\mbox{and}
    \\ &\bigg[\sum_{j\in J_\star^c}\mathbb E\bigg(|G|-\frac{\lambda\sqrt N}{\sigma_\xi\sqrt{\sigma_j}}\bigg)_+^2\bigg]^{\frac{1}{2}}\leq\sqrt{d_q}+\sqrt{2N}\log\bigg(\frac{\sigma_\xi z_q}{\lambda\sqrt N}\bigg)\leq5\sqrt{d_q}.
\end{aligned}
\end{align} From $z_q\geq \frac{\lambda\sqrt{N}}{\sigma_\xi}$ together with monotonicity of $\eta\mapsto \bE(|G|-\eta)_+^2$, $d_q = \sum_{j\in J_\star^c}\bE(|G|-\frac{z_q}{\sqrt{\sigma_j}})_+^2 \leq \sum_{j\in J_\star^c}\bE(|G|-\frac{\lambda\sqrt{N}}{\sigma_\xi\sqrt{\sigma_j}})_+^2$. From \eqref{eq:consequence_large-spectral-control}, together with $\frac{d_q}{N}\leq 2^{-16}\leq \frac{1}{128}$, there hold
\begin{align}\label{eq:estimate_d_q_Var_anisotropic}
    d_q\leq\sum_{j\in J_\star^c}\mathbb E\bigg(|G|-\frac{\lambda\sqrt N}{\sigma_\xi\sqrt{\sigma_j}}\bigg)_+^2\leq25d_q,\qquad\mbox{and}\qquad \frac{\sigma_\xi^2d_q}{N}\leq\mathrm{Var}_{J_\star^c}\leq\frac{128\sigma_\xi^2d_q}{127N}.
\end{align}

Apply Lemma~\ref{lem:common-calibration} with $J_\star^c$ replaced by $[p]$ to define $q^{\mathrm{all}}\in(0,N)$ by
\begin{align*} q^{\mathrm{all}}=\sum_{j=1}^p2\Phi_G\bigg(\frac{z_{q^{\mathrm{all}}}^{\mathrm{all}}}{\sqrt{\sigma_j}}\bigg),\qquad d_{q^{\mathrm{all}}}^{\mathrm{all}}=\sum_{j=1}^p\mathbb E\bigg(|G|-\frac{z_{q^{\mathrm{all}}}^{\mathrm{all}}}{\sqrt{\sigma_j}}\bigg)_+^2,\quad\mbox{and}\quad \frac{z_{q^{\mathrm{all}}}^{\mathrm{all}}(N-q^{\mathrm{all}})}{\sqrt N(N-d_{q^{\mathrm{all}}}^{\mathrm{all}})^{\frac{1}{2}}}=\frac{\sqrt N\lambda}{\sigma_\xi}. \end{align*}
The preceding comparison and~\eqref{eq:def_r_star_without_constant} give
\begin{align*} d_{q^{\mathrm{all}}}^{\mathrm{all}}\leq\sum_{j=1}^p\mathbb E\bigg(|G|-\frac{\lambda\sqrt N}{\sigma_\xi\sqrt{\sigma_j}}\bigg)_+^2\leq|J_\star|+25d_q\leq\frac{25Nr^2(V_{J_\star},V_{J_\star^c})}{\sigma_\xi^2}, \end{align*}where we have used the fact that for any $j\in J_\star$, $\mathbb E\big(|G|-\frac{\lambda\sqrt N}{\sigma_\xi\sqrt{\sigma_j}}\big)_+^2\leq\bE|G|^2=1$.
The same calculation with $J_\star^c$ replaced by $[p]$ gives
\begin{align*} q^{\mathrm{all}}\leq\sqrt{2Nd_{q^{\mathrm{all}}}^{\mathrm{all}}},\qquad 1\leq\frac{\sigma_\xi z_{q^{\mathrm{all}}}^{\mathrm{all}}}{\lambda\sqrt N}\leq2,\qquad \log\bigg(\frac{\sigma_\xi z_{q^{\mathrm{all}}}^{\mathrm{all}}}{\lambda\sqrt N}\bigg)\leq\frac{2q^{\mathrm{all}}}{N},\quad\mbox{and}\quad d_{q^{\mathrm{all}}}^{\mathrm{all}}\leq\sum_{j=1}^p\mathbb E\bigg(|G|-\frac{\lambda\sqrt N}{\sigma_\xi\sqrt{\sigma_j}}\bigg)_+^2\leq25d_{q^{\mathrm{all}}}^{\mathrm{all}}, \end{align*} where we have used $\frac{\sigma_\xi z_{q^{\mathrm{all}}}^{\mathrm{all}}}{\lambda\sqrt{N}}=\frac{\sqrt{1-\frac{d_{q^{\mathrm{all}}}^{\mathrm{all}}}{N}}}{1-\frac{q^{\mathrm{all}}}{N}}\leq\frac{1}{1-\frac{q^{\mathrm{all}}}{N}} \leq 2$.

Replace Lemma~\ref{lem:geometric-event} in the well-regularized regime by Lemmas~\ref{lem:diagonal-offset-comparison} and~\ref{lem:diagonal-active-support} below. The constants fixed at the outset satisfy $A\ge2^{34}$ and $r^2(V_{J_\star},V_{J_\star^c})\le A^{-16}\sigma_\xi^2$. Lemma~\ref{lem:large-spectral-control} and~\eqref{eq:common-free-variance} give
\begin{align*} \sum_{j\in J_\star^c}\mathbb E\bigg(|G|-\frac{\lambda\sqrt N}{\sigma_\xi\sqrt{\sigma_j}}\bigg)_+^2\le\frac{25N}{\sigma_\xi^2}r^2(V_{J_\star},V_{J_\star^c})\le25A^{-16}N. \end{align*}

\begin{Lemma}\label{lem:diagonal-offset-comparison}
For every positive diagonal $\Sigma_{J_\star^c}$ and every $\lambda>0$,
\begin{align*} \mathbb E\sup_{\substack{\boldsymbol v\in V_{J_\star^c}\\\|\boldsymbol v\|_2=1}}\bigg(\|\mathbb X_{J_\star^c}\Sigma_{J_\star^c}^{-1/2}\boldsymbol v\|_2-\frac{\lambda\sqrt N}{4\sigma_\xi}\|\Sigma_{J_\star^c}^{-1/2}\boldsymbol v\|_1\bigg)\le\sqrt N+\bigg[\sum_{j\in J_\star^c}\mathbb E\bigg(|G|-\frac{\lambda\sqrt N}{4\sigma_\xi\sqrt{\sigma_j}}\bigg)_+^2\bigg]^{\frac{1}{2}}. \end{align*}
The random supremum on the left is $1$-Lipschitz in the independent standard Gaussian entries of $\mathbb X_{J_\star^c}\Sigma_{J_\star^c}^{-1/2}$.
\end{Lemma}
\begin{proof}
Replace the Gaussian-width estimate~\eqref{eq:dual-width} in the proof of Lemma~\ref{lem:geometric-event} by the following Gaussian comparison with deterministic offsets. For unit vectors $\boldsymbol u,\boldsymbol u'\in\mathbb R^N$ and $\boldsymbol v,\boldsymbol v'\in V_{J_\star^c}$, the centered Gaussian processes $\langle\boldsymbol u,\mathbb X_{J_\star^c}\Sigma_{J_\star^c}^{-1/2}\boldsymbol v\rangle$ and $\langle\boldsymbol g,\boldsymbol u\rangle+\langle\boldsymbol h_{J_\star^c},\boldsymbol v\rangle$ have increment variances
\begin{align*} 2-2\langle\boldsymbol u,\boldsymbol u'\rangle\langle\boldsymbol v,\boldsymbol v'\rangle\le4-2\langle\boldsymbol u,\boldsymbol u'\rangle-2\langle\boldsymbol v,\boldsymbol v'\rangle. \end{align*}
Take $\boldsymbol g$ and $\boldsymbol h_{J_\star^c}$ independent standard Gaussian vectors, independent of the design. For a finite net $(\boldsymbol u_i,\boldsymbol v_i)_i$ and $\eta>0$, let $p_i(t)$ be the normalized exponential weights in the following expression. The interpolation argument in \cite[Theorem 7.2.11 and its proof, pp.~165--166]{vershyninHighDimensionalProbabilityIntroduction2018} gives
\begin{align*} &\frac{\mathrm d}{\mathrm dt}\mathbb E\frac{1}{\eta}\log\sum_i\exp\!\bigg(\eta\bigg[\sqrt{1-t}\langle\boldsymbol u_i,\mathbb X_{J_\star^c}\Sigma_{J_\star^c}^{-1/2}\boldsymbol v_i\rangle+\sqrt t\big(\langle\boldsymbol g,\boldsymbol u_i\rangle+\langle\boldsymbol h_{J_\star^c},\boldsymbol v_i\rangle\big)-\frac{\lambda\sqrt N}{4\sigma_\xi}\|\Sigma_{J_\star^c}^{-1/2}\boldsymbol v_i\|_1\bigg]\bigg)\\ &\qquad=\frac{\eta}{2}\mathbb E\sum_{i,k}p_i(t)p_k(t)(1-\langle\boldsymbol u_i,\boldsymbol u_k\rangle)(1-\langle\boldsymbol v_i,\boldsymbol v_k\rangle)\ge0. \end{align*}
Letting $\eta$ tend to infinity, and then increasing the finite nets, gives
\begin{align*} \mathbb E\sup_{\|\boldsymbol v\|_2=1}\bigg(\|\mathbb X_{J_\star^c}\Sigma_{J_\star^c}^{-1/2}\boldsymbol v\|_2-\frac{\lambda\sqrt N}{4\sigma_\xi}\|\Sigma_{J_\star^c}^{-1/2}\boldsymbol v\|_1\bigg)\le\mathbb E\|\boldsymbol g\|_2+\mathbb E\sup_{\|\boldsymbol v\|_2=1}\bigg(\langle\boldsymbol h_{J_\star^c},\boldsymbol v\rangle-\frac{\lambda\sqrt N}{4\sigma_\xi}\|\Sigma_{J_\star^c}^{-1/2}\boldsymbol v\|_1\bigg). \end{align*}
The last supremum is at most $[\sum_{j\in J_\star^c}(|h_j|-\frac{\lambda\sqrt N}{4\sigma_\xi\sqrt{\sigma_j}})_+^2]^{\frac{1}{2}}$. Jensen's inequality proves the expectation bound. The difference of the two suprema at two matrices is bounded by the Frobenius norm of their difference, which proves the Lipschitz assertion.
\end{proof}

\begin{Lemmaextend}{lem:geometric-event}\label{lem:diagonal-active-support}
Assume $\lambda\ge\frac{8\sigma_\xi z_N^{\mathrm{all}}}{\sqrt N}$ and $r^2(V_{J_\star},V_{J_\star^c})\le A^{-16}\sigma_\xi^2$. There exists an event $\Omega_{\mathrm{act}}$, depending only on $\mathbb X_{J_\star^c}$, of probability at least $1-e^{-\frac{N}{8}}$, such that
\begin{align*} \forall\boldsymbol z\in2\sigma_\xi\sqrt N B_2^N,\qquad|\operatorname{supp}(\mathcal A[\boldsymbol z])|\le\frac{5N}{16}. \end{align*}
\end{Lemmaextend}
\begin{proof}
Replace the polar-radius event~\eqref{eq:dual-radius} in Lemma~\ref{lem:geometric-event} by the active-support event below.

\textbf{Stochastic Argument.} Equations~\eqref{eq:Gaussian_calculus} and~\eqref{eq:lower_expectation_soft_max} give $\mathbb E(|G|-16)_+^2>2^{-210}$. Consequently $25A^{-16}\le\frac{\mathbb E(|G|-16)_+^2}{64}$, and monotonicity gives
\begin{align*} \bigg|\bigg\{j\in J_\star^c:\frac{\lambda\sqrt N}{\sigma_\xi\sqrt{\sigma_j}}\le16\bigg\}\bigg|\le\frac{N}{64}. \end{align*}
On coordinates with $\frac{\lambda\sqrt N}{4\sigma_\xi\sqrt{\sigma_j}}>3$, equation~\eqref{eq:expectation_soft_max} bounds the second moment by $\frac{2}{9}\,2\Phi_G(\frac{\lambda\sqrt N}{4\sigma_\xi\sqrt{\sigma_j}})$. The remaining coordinates belong to the preceding set. Since $\frac{\lambda\sqrt N}{4\sigma_\xi}\ge2z_N^{\mathrm{all}}$, we obtain
\begin{align*} \sum_{j\in J_\star^c}\mathbb E\bigg(|G|-\frac{\lambda\sqrt N}{4\sigma_\xi\sqrt{\sigma_j}}\bigg)_+^2\le\frac{N}{64}+\frac{2}{9}\sum_{j\in J_\star^c}2\Phi_G\bigg(\frac{\lambda\sqrt N}{4\sigma_\xi\sqrt{\sigma_j}}\bigg)\le\frac{N}{4}. \end{align*}
Let $\Omega_{\mathrm{gp}}$ be the general-position event constructed in the proof of Lemma~\ref{lem:general-position-support}. Its probability remains one because the finitely many exceptional sets are polynomial zero sets and the positive diagonal Gaussian design has a density. Define $\Omega_{\mathrm{act}}$ as the intersection of $\Omega_{\mathrm{gp}}$ and
\begin{align*} \bigg\{\sup_{\substack{\boldsymbol v\in V_{J_\star^c}\\\|\boldsymbol v\|_2=1}}\bigg(\|\mathbb X_{J_\star^c}\Sigma_{J_\star^c}^{-1/2}\boldsymbol v\|_2-\frac{\lambda\sqrt N}{4\sigma_\xi}\|\Sigma_{J_\star^c}^{-1/2}\boldsymbol v\|_1\bigg)\le2\sqrt N\bigg\}. \end{align*}
Lemma~\ref{lem:diagonal-offset-comparison} bounds the expectation of this supremum by $\frac{3\sqrt N}{2}$. Gaussian concentration for its $1$-Lipschitz dependence on the whitened free matrix gives $\mathbb P(\Omega_{\mathrm{act}}^c)\le e^{-\frac{N}{8}}$.

\textbf{Deterministic Argument.} The dual norm identity gives, on $\Omega_{\mathrm{act}}$,
\begin{align*} \sup_{\|\boldsymbol u\|_2=1}\bigg[\sum_{j\in J_\star^c}\bigg(\frac{|\langle\mathbb X\boldsymbol e_j,\boldsymbol u\rangle|}{\sqrt{\sigma_j}}-\frac{\lambda\sqrt N}{4\sigma_\xi\sqrt{\sigma_j}}\bigg)_+^2\bigg]^{\frac{1}{2}}\le2\sqrt N. \end{align*}
Indeed, the left side equals the maximum of zero and the supremum defining $\Omega_{\mathrm{act}}$. Comparing the objective in~\eqref{eq:def_cA} at $\mathcal A[\boldsymbol z]$ and at zero gives $\|\boldsymbol z-\mathbb X_{J_\star^c}\mathcal A[\boldsymbol z]\|_2\le\|\boldsymbol z\|_2$. If this residual vanishes, KKT and $\lambda>0$ give $\mathcal A[\boldsymbol z]=0$. Otherwise use its unit direction in the preceding inequality. On every active coordinate the KKT equality implies
\begin{align*} \frac{|\langle\mathbb X\boldsymbol e_j,\boldsymbol z-\mathbb X_{J_\star^c}\mathcal A[\boldsymbol z]\rangle|}{\sqrt{\sigma_j}\|\boldsymbol z-\mathbb X_{J_\star^c}\mathcal A[\boldsymbol z]\|_2}=\frac{N\lambda}{\sqrt{\sigma_j}\|\boldsymbol z-\mathbb X_{J_\star^c}\mathcal A[\boldsymbol z]\|_2}\ge\frac{\lambda\sqrt N}{2\sigma_\xi\sqrt{\sigma_j}}. \end{align*}
Therefore $\sum_{j\in\operatorname{supp}(\mathcal A[\boldsymbol z])}\frac{N\lambda^2}{\sigma_\xi^2\sigma_j}\le64N$. At most $\frac{N}{4}$ active coordinates have $\frac{\lambda\sqrt N}{\sigma_\xi\sqrt{\sigma_j}}>16$, and at most $\frac{N}{64}\le \frac{N}{16}$ have the opposite inequality. This proves the assertion simultaneously for all the stated inputs.
\end{proof}

\subsection{Examples}\label{sec:examples}

\subsubsection{Example~\texorpdfstring{\ref{example:polynomial_decay}}{polynomial decay}}

We first estimate $u\in\bR_+\mapsto \sum_{j=1}^p 2\Phi_G(\frac{u}{\sqrt{\sigma_j}})$.
For every $u>0$, Tonelli's theorem gives
\begin{align*} \int_0^\infty2\Phi_G(ux^{\frac{\alpha}{2}})\,\mathrm dx=\mathbb E\int_0^\infty\mathbf1_{\{x\leq(\frac{|G|}{u})^{\frac{2}{\alpha}}\}}\,\mathrm dx=u^{-\frac{2}{\alpha}}\mathbb E|G|^{\frac{2}{\alpha}}. \end{align*}
Let $f_u:x\in\bR_+\mapsto 2\Phi_G(c u x^{\frac{\alpha}{2}})$, then $f_u$ is decreasing and is non-negative. Hence $\sum_{j=1}^p f_u(j)\leq \sum_{j=1}^\infty f_u(j)\leq\int_0^\infty f_u(x)\mathrm{d}x = \int_0^\infty 2\Phi_G(cux^{\frac{\alpha}{2}})\mathrm{d}x =  (cu)^{-\frac2\alpha}\bE|G|^{\frac2\alpha}$. Combining with the assumption that $\sigma_j\asymp j^{-\alpha}$, we obtain that
\begin{align*}
\sum_{j=1}^p2\Phi_G\bigg(\frac{u}{\sqrt{\sigma_j}}\bigg) \lesssim \sum_{j=1}^p f_u(j) \lesssim\min\big\{p,u^{-\frac{2}{\alpha}}\big\},\qquad u>0. \end{align*}We next show that this estimate can be reversed. Indeed, for any $j\leq u^{-\frac{2}{\alpha}}$, from $\sigma_j\asymp j^{-\alpha}$, we obtain that $uj^{\frac{\alpha}{2}}\lesssim 1$. Therefore, for any $j\leq \min\{p,\lfloor u^{-\frac{2}{\alpha}}\rfloor\}$, there holds $\frac{u}{\sqrt{\sigma_j}}\lesssim C$, and consequently,
\begin{align*} \sum_{j=1}^p2\Phi_G\bigg(\frac{u}{\sqrt{\sigma_j}}\bigg)\gtrsim\sum_{j=1}^{\min\{p,\lfloor u^{-\frac{2}{\alpha}}\rfloor\}}1\asymp\min\big\{p,u^{-\frac{2}{\alpha}}\big\}. \end{align*}
By the definition of $z_N^{\mathrm{all}}$,
\begin{align*} N=\sum_{j=1}^p2\Phi_G\bigg(\frac{z_N^{\mathrm{all}}}{\sqrt{\sigma_j}}\bigg)\lesssim(z_N^{\mathrm{all}})^{-\frac{2}{\alpha}},\qquad\mbox{and}\qquad z_N^{\mathrm{all}}\lesssim N^{-\frac{\alpha}{2}}\to0. \end{align*}
\paragraph{Convergence rate of estimation error.}
Consequently, the two bounds above apply at $u=z_N^{\mathrm{all}}$ for sufficiently large $N$ and yield $N\asymp\min\big\{p,(z_N^{\mathrm{all}})^{-\frac{2}{\alpha}}\big\}. $
Since $\frac{p}{N}\to\infty$, the second term realizes the minimum for sufficiently large $N$. Hence $z_N^{\mathrm{all}}\asymp N^{-\frac{\alpha}{2}}$ and $\|\Sigma\|_{\mathrm{op}}\|\Sigma^{-1}\|_{\mathrm{op}}\asymp p^\alpha\to\infty.$ 
Choose
\begin{align*} \lambda\asymp\frac{\sigma_\xi}{\sqrt N}N^{-\frac{\alpha}{2(\alpha s+1)}},\qquad\mbox{so that}\qquad \frac{\lambda\sqrt N}{\sigma_\xi z_N^{\mathrm{all}}}\asymp N^{\frac{\alpha^2s}{2(\alpha s+1)}}\to\infty. \end{align*}
The definition $J_\star=\{j\in[p]:\sigma_j|\beta_j^\star|>\lambda\}$ implies
\begin{align*} \lambda^2|J_\star|^{\alpha(s+1)+1}\lesssim\lambda^2\sum_{j\in J_\star}\sigma_j^{-(s+1)}\le\|\Sigma^{\frac{1-s}{2}}\boldsymbol\beta^\star\|_2^2,\qquad\mbox{and}\qquad |J_\star|\asymp N^{\frac{1}{\alpha s+1}}=o(N). \end{align*}
For every $k>0$, split the coordinate sum at $\lfloor k\rfloor$. The source condition yields
\begin{align*} \lambda^2\sum_{j\in J_\star}\frac{1}{\sigma_j}+\|\Sigma_{J_\star^c}^{1/2}\boldsymbol\beta_{J_\star^c}^\star\|_2^2&=\sum_{j=1}^p\min\bigg\{\frac{\lambda^2}{\sigma_j},\sigma_j(\beta_j^\star)^2\bigg\}\lesssim\inf_{k>0}\bigg\{\lambda^2k^{\alpha+1}+\frac{\|\Sigma^{\frac{1-s}{2}}\boldsymbol\beta^\star\|_2^2}{k^{\alpha s}}\bigg\}\\ &\lesssim\|\Sigma^{\frac{1-s}{2}}\boldsymbol\beta^\star\|_2^{\frac{2(\alpha+1)}{\alpha(s+1)+1}}\lambda^{\frac{2\alpha s}{\alpha(s+1)+1}}\lesssim\sigma_\xi^2N^{-\frac{\alpha s}{\alpha s+1}}. \end{align*}
Apply $\sum_{j=1}^p2\Phi_G\big(\frac{u}{\sqrt{\sigma_j}}\big) \asymp\min\big\{p,u^{-\frac{2}{\alpha}}\big\}$ to $u=\frac{\lambda\sqrt{N}}{\sigma_\xi}$ yields
\begin{align*} \sum_{j=1}^p2\Phi_G\bigg(\frac{\lambda\sqrt N}{\sigma_\xi\sqrt{\sigma_j}}\bigg)\asymp N^{\frac{1}{\alpha s+1}}\asymp\sum_{j=1}^p\mathbb E\bigg(|G|-\frac{\lambda\sqrt N}{\sigma_\xi\sqrt{\sigma_j}}\bigg)_+^2=o(N),\end{align*} where we have used the fact that there exists a constant $c$ such that $2c\Phi_G(2t)\leq\bE[(|G|-t)_+^2]\leq 2\Phi_G(t)$ for any $t\geq 0$\footnote{Indeed, for $0\leq t\leq 1$, by monotonicity, $\bE[(|G|-t)_+^2]\geq \bE[(|G|-1)_+^2]$, while $2\Phi_G(2t)\leq 1$; for $t\geq 1$, $\{|G|\geq 2t\}$ implies $(|G|-t)^2\geq t^2\geq 1$ and hence $\bE[(|G|-t)_+^2]\geq 2\Phi_G(2t)$. This proves the lower estimate. The upper estimate is proved in \eqref{eq:identity_exp_square_soft_max}}.

By \eqref{eq:estimate_d_q_Var_anisotropic}, $d_q\asymp \sum_{j\in J_\star^c}\bE(|G|-\frac{\lambda\sqrt{N}}{\sigma_\xi\sqrt{\sigma_j}})_+^2\asymp N^{\frac{1}{\alpha s+1}}$. Therefore
\begin{align*} d_q=\sum_{j\in J_\star^c}\mathbb E\bigg(|G|-\frac{z_q}{\sqrt{\sigma_j}}\bigg)_+^2\asymp N^{\frac{1}{\alpha s+1}},\qquad\mbox{and}\qquad \mathrm{Var}_{J_\star^c}=\frac{\sigma_\xi^2d_q}{N-d_q}\asymp\sigma_\xi^2N^{-\frac{\alpha s}{\alpha s+1}}. \end{align*}
Together with~\eqref{eq:def_r_star_without_constant}, these estimates give
\begin{align*} r^2(V_{J_\star},V_{J_\star^c})\asymp\sigma_\xi^2N^{-\frac{\alpha s}{\alpha s+1}}. \end{align*}

\paragraph{Monotonicity of \texorpdfstring{$\mathrm{VAR}^{\mathrm{FSD}}$}{FSD variance}}
By Proposition~\ref{prop:monotonicity_bias_variance}, it suffices to verify $|J_\star(\lambda_0)^c|>N$ and, for every $\lambda\in(\lambda_0,\infty)\cap\{\sigma_j|\beta_j^\star|:j\in[p]\}$,
\begin{align}\begin{aligned}&d_{q(\lambda)}-d_{q(\lambda-)}\leq|J_\star(\lambda-)\setminus J_\star(\lambda)|\bigg(1-\frac{d_{q(\lambda)}}N\bigg)\bigg(1-\frac{d_{q(\lambda-)}}N\bigg).\end{aligned}\tag{\ref{eq:fsd-variance-diagonal-criterion}}\end{align}
All implicit constants may depend on $\alpha$, $s$, the spectral comparison constants, the constants in $\lambda_0\asymp\frac{\sigma_\xi}{\sqrt N}$, and $\frac{\|\Sigma^{\frac{1-s}{2}}\boldsymbol\beta^\star\|_2}{\sigma_\xi}$, and are independent of $N,p$.
\begin{enumerate}
\item\label{item:fsd-example-uniform-estimates} The source-condition estimate established earlier in this subsection gives $|J_\star(\lambda_0)|\lesssim(\frac{\|\Sigma^{\frac{1-s}{2}}\boldsymbol\beta^\star\|_2}{\lambda_0})^{\frac{2}{\alpha(s+1)+1}}\lesssim N^{\frac{1}{\alpha(s+1)+1}}=o(N)$. Since $J_\star(\lambda)\subseteq J_\star(\lambda_0)$ for $\lambda\geq\lambda_0$, $|J_\star(\lambda)^c|\geq p-|J_\star(\lambda_0)|>N$ uniformly in $\lambda\geq\lambda_0$ for all sufficiently large $N$ and $p$.

For each $\lambda\geq\lambda_0$, Lemma~\ref{lem:common-calibration}, applied with the set $J_\star(\lambda)^c$, gives $0<d_{q(\lambda)}<q(\lambda)<N$ and $z_{q(\lambda)}\geq\frac{\sqrt N\lambda}{\sigma_\xi}\geq\frac{\sqrt N\lambda_0}{\sigma_\xi}$. The polynomial-tail estimate established earlier in this subsection gives
\begin{align}\begin{aligned}&\sum_{j=1}^p2\Phi_G\bigg(\frac{\sqrt N\lambda_0}{\sigma_\xi\sqrt{\sigma_j}}\bigg)\lesssim\bigg(\frac{\sqrt N\lambda_0}{\sigma_\xi}\bigg)^{-\frac{2}{\alpha}}\asymp1.\end{aligned}\label{eq:fsd-example-full-tail}\end{align}
Thus $q(\lambda)=\sum_{j\in J_\star(\lambda)^c}2\Phi_G(\frac{z_{q(\lambda)}}{\sqrt{\sigma_j}})\leq\sum_{j=1}^p2\Phi_G(\frac{\sqrt N\lambda_0}{\sigma_\xi\sqrt{\sigma_j}})$, and \eqref{eq:fsd-example-full-tail} yields $\sup_{\lambda\geq\lambda_0}d_{q(\lambda)}\leq\sup_{\lambda\geq\lambda_0}q(\lambda)\lesssim1$. For every $\lambda>\lambda_0$, $J_\star(t)=J_\star(\lambda-)$ for $t<\lambda$ sufficiently close to $\lambda$. Continuity for this fixed set in Lemma~\ref{lem:common-calibration}, followed by $t\uparrow\lambda$, gives $z_{q(\lambda-)}\geq\frac{\sqrt N\lambda_0}{\sigma_\xi}$ and $d_{q(\lambda-)}\leq q(\lambda-)\lesssim1$, uniformly in $\lambda>\lambda_0$.

\item\label{item:fsd-example-threshold-comparison} To control the difference on the left-hand side of \eqref{eq:fsd-variance-diagonal-criterion}, we prove $z_{q(\lambda-)}\leq z_{q(\lambda)}$. Lemma~\ref{lem:common-calibration} gives monotonicity for a fixed free-coordinate set; here $J_\star(\lambda)^c=J_\star(\lambda-)^c\cup(J_\star(\lambda-)\setminus J_\star(\lambda))$, with disjoint union. By item~\ref{item:fsd-example-uniform-estimates}, \eqref{eq:fsd-example-full-tail}, and $\bE(|G|-a)_+^2<2\Phi_G(a)$ for $a>0$, $\sum_{j\in J_\star(\lambda)^c}\bE(|G|-\frac{z_{q(\lambda-)}}{\sqrt{\sigma_j}})_+^2<\sum_{j\in J_\star(\lambda)^c}2\Phi_G(\frac{z_{q(\lambda-)}}{\sqrt{\sigma_j}})\lesssim1<N$ for all sufficiently large $N$. Together with $d_{q(\lambda-)}\leq q(\lambda-)\lesssim1$, this gives
\begin{align}\begin{aligned}&0\leq\frac{\sum_{j\in J_\star(\lambda-)\setminus J_\star(\lambda)}\bE\big(|G|-\frac{z_{q(\lambda-)}}{\sqrt{\sigma_j}}\big)_+^2}{N-d_{q(\lambda-)}}\leq\frac{\sum_{j\in J_\star(\lambda-)\setminus J_\star(\lambda)}2\Phi_G\big(\frac{z_{q(\lambda-)}}{\sqrt{\sigma_j}}\big)}{N-q(\lambda-)}<1.\end{aligned}\label{eq:fsd-example-relative-increments}\end{align}
Since $1-x\leq\sqrt{1-y}$ for $0\leq y\leq x<1$, \eqref{eq:fsd-example-relative-increments}, the disjoint union, and $\frac{z_{q(\lambda-)}(N-q(\lambda-))}{\sqrt N(N-d_{q(\lambda-)})^{\frac12}}=\frac{\sqrt N\lambda}{\sigma_\xi}$ imply
\begin{align*}&\frac{z_{q(\lambda-)}\bigg(N-\sum_{j\in J_\star(\lambda)^c}2\Phi_G\big(\frac{z_{q(\lambda-)}}{\sqrt{\sigma_j}}\big)\bigg)}{\sqrt N\bigg(N-\sum_{j\in J_\star(\lambda)^c}\bE\big(|G|-\frac{z_{q(\lambda-)}}{\sqrt{\sigma_j}}\big)_+^2\bigg)^{\frac12}}\leq\frac{\sqrt N\lambda}{\sigma_\xi}=\frac{z_{q(\lambda)}(N-q(\lambda))}{\sqrt N(N-d_{q(\lambda)})^{\frac12}}.\end{align*}
Since also $\sum_{j\in J_\star(\lambda)^c}2\Phi_G(\frac{z_{q(\lambda)}}{\sqrt{\sigma_j}})=q(\lambda)<N$, strict monotonicity in \eqref{eq:monotonocity_calibration_threshold}, applied with the same expanded set $J_\star(\lambda)^c$ at both thresholds, gives $z_{q(\lambda-)}\leq z_{q(\lambda)}$.

\item\label{item:fsd-example-criterion} We now verify \eqref{eq:fsd-variance-diagonal-criterion}. By item~\ref{item:fsd-example-threshold-comparison}, $\bE(|G|-\frac{z_{q(\lambda)}}{\sqrt{\sigma_j}})_+^2\leq\bE(|G|-\frac{z_{q(\lambda-)}}{\sqrt{\sigma_j}})_+^2$ for every $j\in J_\star(\lambda-)^c$. Hence the disjoint union in item~\ref{item:fsd-example-threshold-comparison}, $z_{q(\lambda)}\geq\frac{\sqrt N\lambda_0}{\sigma_\xi}$, and $\sigma_j\leq\|\Sigma\|_{\mathrm{op}}$ yield
\begin{align}\begin{aligned}
&d_{q(\lambda)}-d_{q(\lambda-)} = \sum_{j\in J_\star(\lambda_-)^c}\Big( \bE\Big(|G|-\frac{z_{q(\lambda)}}{\sqrt{\sigma_j}}\Big)_+^2 - \bE\Big(|G|-\frac{z_{q(\lambda_-)}}{\sqrt{\sigma_j}}\Big)_+^2 \Big) + \sum_{j\in J_\star(\lambda_-)\backslash J_\star(\lambda)}\bE\Big(|G|-\frac{z_{q(\lambda)}}{\sqrt{\sigma_j}}\Big)_+^2 \\
&\leq\sum_{j\in J_\star(\lambda-)\setminus J_\star(\lambda)}\bE\bigg(|G|-\frac{z_{q(\lambda)}}{\sqrt{\sigma_j}}\bigg)_+^2\leq|J_\star(\lambda-)\setminus J_\star(\lambda)|\,\bE\bigg(|G|-\frac{\sqrt N\lambda_0}{\sigma_\xi\sqrt{\|\Sigma\|_{\mathrm{op}}}}\bigg)_+^2.
\end{aligned}\label{eq:fsd-example-d-increment}\end{align}
Since $\|\Sigma\|_{\mathrm{op}}\asymp1$, $\frac{\sqrt N\lambda_0}{\sigma_\xi}\asymp1$, and $a\mapsto\bE(|G|-a)_+^2$ is continuous and strictly decreasing from $1$ on $[0,\infty)$, $1-\bE(|G|-\frac{\sqrt N\lambda_0}{\sigma_\xi\sqrt{\|\Sigma\|_{\mathrm{op}}}})_+^2\asymp1$. By item~\ref{item:fsd-example-uniform-estimates}, $1-(1-\frac{d_{q(\lambda)}}{N})(1-\frac{d_{q(\lambda-)}}{N})\leq\frac{d_{q(\lambda)}+d_{q(\lambda-)}}{N}\lesssim N^{-1}$. Therefore, uniformly in the stated values of $\lambda$, for all sufficiently large $N$,
\begin{align*}&\bE\bigg(|G|-\frac{\sqrt N\lambda_0}{\sigma_\xi\sqrt{\|\Sigma\|_{\mathrm{op}}}}\bigg)_+^2\leq\bigg(1-\frac{d_{q(\lambda)}}N\bigg)\bigg(1-\frac{d_{q(\lambda-)}}N\bigg).\end{align*}
Together with \eqref{eq:fsd-example-d-increment}, this proves \eqref{eq:fsd-variance-diagonal-criterion} for every $\lambda\in(\lambda_0,\infty)\cap\{\sigma_j|\beta_j^\star|:j\in[p]\}$. Proposition~\ref{prop:monotonicity_bias_variance} therefore implies that $\mathrm{VAR}^{\mathrm{FSD}}$ is non-increasing on $[\lambda_0,\infty)$.
\end{enumerate}

\subsubsection{Example~\texorpdfstring{\ref{example:infinite}}{infinite-dimensional LASSO}}

We use the notation of Example~\ref{example:infinite}. The argument proceeds in four steps: (i) the condition $X\in c_0$ is translated into summability of Gaussian tails, which makes $J_\star^\infty$ finite and the infinite-dimensional quantities in \eqref{eq:r_infty_infinite_example} well defined; (ii) the calibration equation is shown to be stable as $p\to\infty$, which gives $r_p^2\to r_\infty^2$; (iii) Assumption~\ref{assumption:LASSO_upper} is verified for large $p$, so that Theorem~\ref{theo:upper_LASSO} applies; and (iv) the bound is transferred to a minimizer of the LASSO over $\ell_1$ by weak-$*$ compactness.

\paragraph{Step (i): summability.} For $u>0$, the events $\{\sqrt{\sigma_j}|G_j|>u\}$, $j\geq1$, are independent and have probabilities $2\Phi_G(\frac{u}{\sqrt{\sigma_j}})$. Since $\sqrt{\sigma_j}G_j\to0$ almost surely, only finitely many of them occur almost surely, and the second Borel--Cantelli lemma gives
\begin{align}\label{eq:infinite_example_summability}
    \sum_{j\geq1}2\Phi_G\bigg(\frac{u}{\sqrt{\sigma_j}}\bigg)<\infty\qquad\mbox{for every }u>0.
\end{align}
In particular $\Phi_G(\frac{u}{\sqrt{\sigma_j}})\to0$, that is, $\sigma_j\to0$. Since $\sum_{j\geq1}\sigma_j(\beta_j^\star)^2<\infty$, we have $\sigma_j|\beta_j^\star|=\sqrt{\sigma_j}\,(\sigma_j(\beta_j^\star)^2)^{\frac12}\to0$; hence $J_\star^\infty$ is finite, and for every $p\geq p_0:=\max J_\star^\infty$ there hold $J_\star^{(p)}=J_\star^\infty$ and $J_\star^{(p),c}=J_\star^{\infty,c}\cap[p]$. Moreover, $\sigma_{\xi,p}^2\downarrow\sigma_\xi^2$ as the tail of a convergent series.

For $p\in\bN_+\cup\{\infty\}$ and $u>0$, let
\begin{align*}
    Q_p(u):=\sum_{j\in J_\star^{\infty,c},\,j\leq p}2\Phi_G\bigg(\frac{u}{\sqrt{\sigma_j}}\bigg),\qquad D_p(u):=\sum_{j\in J_\star^{\infty,c},\,j\leq p}\mathbb E\bigg(|G|-\frac{u}{\sqrt{\sigma_j}}\bigg)_+^2,
\end{align*}
which are finite by \eqref{eq:infinite_example_summability}. By the proof of Lemma~\ref{lem:common-calibration}, $0<\mathbb E(|G|-a)_+^2<2\Phi_G(a)$ for $a>0$; hence $D_\infty(u)<Q_\infty(u)$ for every $u>0$. Both $a\mapsto2\Phi_G(a)$ and $a\mapsto\mathbb E(|G|-a)_+^2$ are decreasing, so that for every $u_0>0$,
\begin{align}\label{eq:infinite_example_uniform}
    \sup_{u\geq u_0}\big|D_\infty(u)-D_p(u)\big|\leq\sup_{u\geq u_0}\big|Q_\infty(u)-Q_p(u)\big|\leq\sum_{j>p}2\Phi_G\bigg(\frac{u_0}{\sqrt{\sigma_j}}\bigg)\xrightarrow[p\to\infty]{}0.
\end{align}
Consequently $Q_\infty$ and $D_\infty$ are continuous on $(0,\infty)$, as locally uniform limits of continuous functions, and they are strictly decreasing, as sums of strictly decreasing functions. Since $J_\star^{\infty,c}$ is infinite, monotone convergence gives $Q_\infty(u)\to\infty$ as $u\downarrow0$, and dominated convergence (with dominating sequence given by the terms of $Q_\infty(1)$) gives $Q_\infty(u)\to0$ as $u\to\infty$. Therefore $z_r^\infty$, the unique solution of $Q_\infty(z_r^\infty)=r$, is well defined for every $r>0$, and $d_r^\infty=D_\infty(z_r^\infty)$. The same argument applied to the sums over $\bN_+$ shows that $z_N^{\mathrm{all},\infty}$ is well defined. Finally, $Q_p\leq Q_\infty$ and the strict decrease of $Q_p$ give $z_N^{(p)}\leq z_N^\infty$ whenever $|J_\star^{(p),c}|>N$, and likewise $z_N^{\mathrm{all},(p)}\leq z_N^{\mathrm{all},\infty}$.

\paragraph{Step (ii): stability of the calibration.} For $p\in\bN_+\cup\{\infty\}$ and $u>z_N^\infty$, let
\begin{align*}
    F_p(u):=\frac{u\big(N-Q_p(u)\big)}{\sqrt N\big(N-D_p(u)\big)^{\frac12}}.
\end{align*}
For $u>z_N^\infty$, we have $D_p(u)\leq D_\infty(u)<Q_\infty(u)<N$ and $Q_p(u)\leq Q_\infty(u)<N$, so $F_p$ is well defined. Let $p$ be such that $|J_\star^{(p),c}|>N$. Since $z_N^{(p)}\leq z_N^\infty$, \eqref{eq:monotonocity_calibration_threshold} applies to $F_p$ on $(z_N^\infty,\infty)$; dropping its second (non-negative) term and using $N-Q_p\geq N-Q_\infty>0$ and $N-D_p\leq N$, we obtain $F_p'(u)\geq(1-\frac{Q_\infty(u)}{N})^2$. Hence, for $z_N^\infty<a<b$,
\begin{align}\label{eq:infinite_example_increment}
    F_p(b)-F_p(a)\geq\kappa(a,b):=\int_a^b\bigg(1-\frac{Q_\infty(u)}{N}\bigg)^2\mathrm du>0.
\end{align}
On $[a,b]$, the denominators satisfy $N-D_\infty(a)\leq N-D_p(u)\leq N$, with $N-D_\infty(a)>0$. Together with \eqref{eq:infinite_example_uniform}, this gives $F_p\to F_\infty$ uniformly on $[a,b]$. Passing to the limit in \eqref{eq:infinite_example_increment}, $F_\infty$ is continuous and strictly increasing on $(z_N^\infty,\infty)$. As $u\downarrow z_N^\infty$, $N-Q_\infty(u)\to0$ while $N-D_\infty(u)\to N-D_\infty(z_N^\infty)>0$, so $F_\infty(u)\to0$; as $u\to\infty$, $F_\infty(u)\geq\frac{u(N-Q_\infty(u))}{N}\to\infty$. Hence there is a unique $u_\infty>z_N^\infty$ with $F_\infty(u_\infty)=\frac{\sqrt N\lambda}{\sigma_\xi}$; we set $q_\infty:=Q_\infty(u_\infty)\in(0,N)$, so that $u_\infty=z_{q_\infty}^\infty$ and $d_{q_\infty}^\infty=D_\infty(u_\infty)\in(0,q_\infty)$.

Let $u_p:=z_{q_p}^{(p)}$. By Lemma~\ref{lem:common-calibration} applied with $J_\star^{(p),c}$ and $\sigma_{\xi,p}$, $u_p$ is the unique solution in $(z_N^{(p)},\infty)$ of $F_p(u_p)=\frac{\sqrt N\lambda}{\sigma_{\xi,p}}$. Fix $z_N^\infty<a<u_\infty<b$. Then $F_\infty(a)<\frac{\sqrt N\lambda}{\sigma_\xi}<F_\infty(b)$. Since $F_p(a)\to F_\infty(a)$, $F_p(b)\to F_\infty(b)$ and $\sigma_{\xi,p}\to\sigma_\xi$, for every sufficiently large $p$ we have $F_p(a)<\frac{\sqrt N\lambda}{\sigma_{\xi,p}}<F_p(b)$, and the strict monotonicity of $F_p$ gives $u_p\in(a,b)$. As $a$ and $b$ are arbitrary, $u_p\to u_\infty$. By \eqref{eq:infinite_example_uniform} and the continuity of $Q_\infty$ and $D_\infty$,
\begin{align*}
    |q_p-q_\infty|\leq\sup_{u\geq a}|Q_p(u)-Q_\infty(u)|+|Q_\infty(u_p)-Q_\infty(u_\infty)|\to0,\qquad |d_{q_p}^{(p)}-d_{q_\infty}^\infty|\leq\sup_{u\geq a}|D_p(u)-D_\infty(u)|+|D_\infty(u_p)-D_\infty(u_\infty)|\to0.
\end{align*}
Since $d_{q_\infty}^\infty<N$, it follows that $\frac{\sigma_{\xi,p}^2d_{q_p}^{(p)}}{N-d_{q_p}^{(p)}}\to\frac{\sigma_\xi^2d_{q_\infty}^\infty}{N-d_{q_\infty}^\infty}$. For $p\geq p_0$, the first two terms of \eqref{eq:r_p_infinite_example} equal $\sigma_{\xi,p}^2\frac{|J_\star^\infty|}{N}$ and $\lambda^2\sum_{j\in J_\star^\infty}\sigma_j^{-1}$, and the third term equals $\sum_{j\in J_\star^{\infty,c},\,j\leq p}\sigma_j(\beta_j^\star)^2$, which increases to $\sum_{j\in J_\star^{\infty,c}}\sigma_j(\beta_j^\star)^2$. Comparing with \eqref{eq:r_infty_infinite_example}, $r_p^2\to r_\infty^2$.

\paragraph{Step (iii): application of Theorem~\ref{theo:upper_LASSO}.} For each $p$, the response satisfies $Y=\langle X_{(p)},\bbeta_{(p)}^\star\rangle+\xi^{(p)}$, where $\xi^{(p)}\sim\cN(0,\sigma_{\xi,p}^2)$ is independent of $X_{(p)}\sim\cN(\boldsymbol0,\Sigma_{(p)})$, so the $p$-dimensional problem is an instance of the Gaussian linear regression model of the Introduction, with covariance $\Sigma_{(p)}$, signal $\bbeta_{(p)}^\star$ and noise level $\sigma_{\xi,p}$. We check item 2 of Assumption~\ref{assumption:LASSO_upper} for $p$ large. For $p\geq p_0$, $1\leq|J_\star^{(p)}|=|J_\star^\infty|\leq cN$, and $|J_\star^{(p),c}|=p-|J_\star^\infty|>N$ as soon as $p>N+|J_\star^\infty|$. By Step (i), $\frac{8\sigma_{\xi,p}z_N^{\mathrm{all},(p)}}{\sqrt N}\leq\frac{8\sigma_{\xi,p}z_N^{\mathrm{all},\infty}}{\sqrt N}\to\frac{8\sigma_\xi z_N^{\mathrm{all},\infty}}{\sqrt N}<\lambda$. By Step (ii), $r_p^2\to r_\infty^2<c\sigma_\xi^2\leq c\sigma_{\xi,p}^2$, so that $r_p^2\leq c\sigma_{\xi,p}^2$ for $p$ large. Hence there is $p_1$ such that, for every $p\geq p_1$, Theorem~\ref{theo:upper_LASSO} applies to the $p$-dimensional problem: the event
\begin{align*}
    \Omega_p:=\Big\{\|\Sigma_{(p)}^{1/2}(\hat\bbeta_{(p)}-\bbeta_{(p)}^\star)\|_2^2\leq Cr_p^2\Big\}\quad\mbox{satisfies}\quad\bP^{\otimes N}(\Omega_p)\geq1-C\exp\big(-c\min\{|J_\star^\infty|,d_{q_p}^{(p)}\}\big).
\end{align*}

\paragraph{Step (iv): passage to the LASSO over $\ell_1$.} Almost surely $X_1,\dots,X_N\in c_0$ and each series $Y_i=\sum_{j\geq1}x_{ij}\beta_j^\star+\xi_i$ converges, being a series of independent centered Gaussian variables with summable variances. Let $L(\bbeta):=\frac{1}{2N\lambda}\sum_{i=1}^N(Y_i-\langle X_i,\bbeta\rangle)^2+\|\bbeta\|_1$ for $\bbeta\in\ell_1$, and identify $\bR^p$ with the sequences supported in $[p]$; then $\langle X_i,\bbeta\rangle=\langle X_{i,(p)},\bbeta_{(p)}\rangle$ for such $\bbeta$, so $\hat\bbeta_{(p)}$ minimizes $L$ over $\bR^p$. We record three facts.
\begin{enumerate}
    \item Since $|\langle X_i,\bbeta-\bbeta'\rangle|\leq\|X_i\|_\infty\|\bbeta-\bbeta'\|_1$, $L$ is continuous on $(\ell_1,\|\cdot\|_1)$. As finitely supported sequences are dense in $\ell_1$, $\lim_{p\to\infty}L(\hat\bbeta_{(p)})=\inf_{\ell_1}L$, the sequence $(L(\hat\bbeta_{(p)}))_p$ being non-increasing.
    \item By $L(\hat\bbeta_{(p)})\leq L(\boldsymbol0)$, $\|\hat\bbeta_{(p)}\|_1\leq\frac{1}{2N\lambda}\sum_{i=1}^NY_i^2$ for every $p$. Since $\ell_1=c_0^*$ and $c_0$ is separable, the Banach--Alaoglu theorem implies that every subsequence of $(\hat\bbeta_{(p)})_p$ has a further weak-$*$ convergent subsequence.
    \item Since $X_i\in c_0$, $\bbeta\mapsto\langle X_i,\bbeta\rangle$ is weak-$*$ continuous, and $\|\cdot\|_1$, being a dual norm, is weak-$*$ lower semicontinuous; hence $L$ is weak-$*$ lower semicontinuous. Moreover, weak-$*$ convergence implies coordinatewise convergence (pair with $\ve_j\in c_0$), so Fatou's lemma for series shows that $\bbeta\mapsto\sum_{j\geq1}\sigma_j(\beta_j-\beta_j^\star)^2$ is weak-$*$ lower semicontinuous.
\end{enumerate}
Let $\Omega_\infty:=\limsup_{p}\Omega_p$. By the reverse Fatou lemma and $d_{q_p}^{(p)}\to d_{q_\infty}^\infty$ from Step (ii),
\begin{align*}
    \bP^{\otimes N}(\Omega_\infty)\geq\limsup_{p\to\infty}\bP^{\otimes N}(\Omega_p)\geq1-C\exp\big(-c\min\{|J_\star^\infty|,d_{q_\infty}^\infty\}\big).
\end{align*}
Fix an outcome in $\Omega_\infty$ for which the almost sure statements above hold. There is a sequence $p_k\to\infty$ such that $\Omega_{p_k}$ holds for every $k$, and by the second fact we may assume that $\hat\bbeta_{(p_k)}\to\hat\bbeta$ weak-$*$ for some $\hat\bbeta\in\ell_1$. By the third and the first facts, $L(\hat\bbeta)\leq\liminf_kL(\hat\bbeta_{(p_k)})=\inf_{\ell_1}L$, so $\hat\bbeta$ minimizes $L$ over $\ell_1$. On $\Omega_{p_k}$,
\begin{align*}
    \sum_{j\geq1}\sigma_j\big((\hat\bbeta_{(p_k)})_j-\beta_j^\star\big)^2=\|\Sigma_{(p_k)}^{1/2}(\hat\bbeta_{(p_k)}-\bbeta_{(p_k)}^\star)\|_2^2+\sum_{j>p_k}\sigma_j(\beta_j^\star)^2\leq Cr_{p_k}^2+\sum_{j>p_k}\sigma_j(\beta_j^\star)^2.
\end{align*}
The third fact, $r_{p_k}^2\to r_\infty^2$ and $\sum_{j>p_k}\sigma_j(\beta_j^\star)^2\to0$ yield $\sum_{j\geq1}\sigma_j(\hat\beta_j-\beta_j^\star)^2\leq Cr_\infty^2$, which is the claim of Example~\ref{example:infinite}.

\bibliographystyle{alpha}
\bibliography{biblio,other}
\end{document}